\documentclass[type=drfinal, dr=rernat, 
blackrule,
longdoc, 12pt, book, BCOR10cm]{tudthesis}

\usepackage{enumerate}

\usepackage{amsmath}

\usepackage{amsthm}  

\usepackage[all,dvips]{xy}

\usepackage[algo2e, boxed]{algorithm2e}

\usepackage{pstricks}

\usepackage{wasysym} 
\usepackage{xspace}

\usepackage{nicefrac}

\usepackage{verbatim}

\theoremstyle{plain}

\newtheorem{theorem}{Theorem}[section]
\newtheorem{lemma}[theorem]{Lemma}
\newtheorem{proposition}[theorem]{Proposition}
\newtheorem{corollary}[theorem]{Corollary}

\theoremstyle{definition}
\newtheorem{definition}[theorem]{Definition}
\newtheorem{assumption}[theorem]{Assumption}
\newtheorem{example}[theorem]{Example}

\theoremstyle{remark}
\newtheorem{remark}[theorem]{Remark}

\newcommand*{\exptower}[1]{\ensuremath{\exp_{#1}}}

\newcommand*{\RamQ}[1]{\ensuremath{\mathrm{Ram}^{#1}}}
\newcommand*{\tci}[0]{\rm{tci}\xspace}

\newcommand*{\N}{\ensuremath{\mathbb{N}}}  

\newcommand*{\Push}[1]{\ensuremath{\mathrm{push}_{#1}}}
\newcommand*{\Pop}[1]{\ensuremath{\mathrm{pop}_{#1}}}
\newcommand*{\Id}[0]{\ensuremath{\mathrm{id}}}
\newcommand*{\Clone}[1]{\ensuremath{{\mathrm{clone}_{#1}}}}
\newcommand*{\Collapse}[0]{\ensuremath{\mathrm{collapse}}}
\newcommand*{\TOP}[1]{\ensuremath{\mathrm{top}_{#1}}}
\newcommand*{\Sym}[0]{\ensuremath{\mathrm{Sym}}}
\newcommand*{\Lvl}[0]{\ensuremath{\mathrm{CLvl}}}
\newcommand*{\Lnk}[0]{\ensuremath{\mathrm{CLnk}}}
\newcommand*{\Op}[0]{\mathrm{OP}}
\newcommand*{\op}[0]{\mathrm{op}}

\newcommand*{\Stacks}[0]{\mathrm{Stacks}}

\newcommand*{\prefixeq}[0]{\ensuremath{\mathop{\trianglelefteq}}}
\newcommand*{\notprefixeq}[0]{\ensuremath{\mathop{\not\trianglelefteq}}}

\newcommand*{\Reach}{\ensuremath{\mathrm{REACH}}}

\newcommand*{\exLoops}[0]{\ensuremath{\exists}\mathrm{Loops}}
\newcommand*{\exHighLoops}[0]{\ensuremath{\exists}\mathrm{HLoops}}
\newcommand*{\exLowLoops}[0]{\ensuremath{\exists}\mathrm{LLoops}}
\newcommand*{\exReturns}[0]{\ensuremath{\exists}\mathrm{Returns}}
\newcommand*{\exOneLoops}[0]{{\ensuremath{\exists1}}\text{-}{\mathrm{Loops}}}

\newcommand*{\ReturnFunc}[1]{\ensuremath{\mathrm{\#Ret}^{#1}}}
\newcommand*{\LoopFunc}[1]{\ensuremath{\mathrm{\#Loop}^{#1}}}
\newcommand*{\HighLoopFunc}[1]{\ensuremath{\mathrm{\#HLoop}^{#1}}}
\newcommand*{\LowLoopFunc}[1]{\ensuremath{\mathrm{\#LLoop}^{#1}}}

\newcommand*{\Decode}[0]{\mathrm{Dec}}
\newcommand*{\Encode}[0]{\mathrm{Enc}}

\newcommand*{\fReturns}[1]{\ensuremath{\mathrm{Rt}}}

\newcommand*{\BiHalfgrid}[0]{\ensuremath{\mathfrak{H}}}

\newcommand*{\arity}[0]{\ensuremath{\mathrm{ar}}}
\newcommand*{\Var}[0]{\ensuremath{\mathrm{Var}}}
\newcommand*{\Freevar}[0]{\ensuremath{\mathrm{Free}}}

\newcommand*{\FO}[1]{\ensuremath{\mathrm{FO}_{#1}}}
\newcommand*{\MSO}{\ensuremath{\mathrm{MSO}}\xspace}

\newcommand*{\lfp}{\ensuremath{\mathrm{lfp}}}

\newcommand*{\False}[0]{\ensuremath{\mathrm{False}}}
\newcommand*{\True}[0]{\ensuremath{\mathrm{True}}}

\newcommand*{\ReachguessSym}[0]{\ensuremath{f_\mathrm{CSR}}}

\newcommand*{\coloneqq}[0]{\ensuremath{\mathrel{\mathop:}=}} 
\newcommand*{\trans}[1]{\ensuremath{\mathrel{{\vdash^{#1}}}}}

\newcommand*{\invtrans}[1]{\ensuremath{\mathrel{{\dashv^{#1}}}}}

\newcommand*{\domain}[0]{\ensuremath{\mathrm{dom}}}

\newcommand*{\length}[0]{\mathrm{ln}}
\newcommand*{\Pot}[1]{\ensuremath{2^{#1}}} 

\newcommand*{\treeLR}[3]{\ensuremath{{#2}\leftarrow{#1} \rightarrow{#3}}}
\newcommand*{\treeL}[2]{\ensuremath{{#2}\leftarrow {#1}}}
\newcommand*{\treeR}[2]{\ensuremath{{#1} \rightarrow{#2}}}

\newcommand*{\LeftTree}[1]{\ensuremath{LT({#1})}}

\newcommand*{\LeftStack}[0]{\ensuremath{\mathrm{LStck}}}

\newcommand*{\EncTrees}[0]{\ensuremath{\mathbb{T}_{\mathrm{Enc}}}}

\newcommand*{\Trees}[1]{\ensuremath{\mathrm{Tree}_{#1}}}
\newcommand*{\infTrees}[1]{\ensuremath{\mathrm{Tree}^\omega_{#1}}}
\newcommand*{\allTrees}[1]{\ensuremath{\mathrm{Tree}^{\leq \omega}_{#1}}}
\newcommand*{\fullTrees}[1]{
  \ensuremath{\mathrm{B}}\text{-}\ensuremath{\mathrm{Tree}^{\omega}_{#1}}} 
\newcommand*{\inducedTreeof}[2]{\ensuremath{({#2})_{#1}}}
\newcommand*{\depth}[1]{\ensuremath{\mathrm{dp}({#1})}}
\newcommand*{\initialTreeq}[0]{\ensuremath{\mathrel{\preceq}}}
\newcommand*{\prune}[0]{\ensuremath{\mathrm{prune}}}
\newcommand*{\extract}[0]{\ensuremath{\mathrm{extract}}}
\newcommand*{\Closure}[0]{\mathrm{CL}}

\newcommand*{\pretree}[0]{\mathrm{extree}}
\newcommand*{\Set}[0]{\mathrm{Set}}
\newcommand{\In}{\mathrm{In}}
\newcommand*{\reduction}[0]{\mathrm{red}}

\newcommand{\SMALL}[0]{small\xspace}
\newcommand{\SMALLRel}[0]{\ensuremath{\mathrm{Small}}}

\newcommand{\OmegaExp}[0]{\ensuremath{\mathrm{Ext}}}

\newcommand*{\NPT}{\ensuremath{\mathrm{NPT}}}
\newcommand*{\CPS}[0]{\ensuremath{\mathrm{CPS}}\xspace}
\newcommand*{\CPG}[0]{\ensuremath{\mathrm{CPG}}\xspace}
\newcommand*{\Graph}[0]{\mathrm{Grph}}

\newcommand*{\Milestones}[0]{\ensuremath{\mathrm{MS}}}
\newcommand*{\genMilestones}[0]{\ensuremath{\mathrm{GMS}}}

\newcommand*{\OccNPT}[0]{\ensuremath{\Xi}}

\newcommand*{\FuncBoundReturnLength}[2]{\ensuremath{ \mathrm{BRL}_{#2}^{#1}}}  
\newcommand*{\FuncBoundLoopLength}[2]{\ensuremath{ \mathrm{BLL}_{#2}^{#1}}}  
\newcommand*{\FuncBoundHighLoopLength}[2]{\ensuremath{ \mathrm{BHLL}_{#2}^{#1}}}  
\newcommand*{\FuncBoundLowLoopLength}[2]{\ensuremath{ \mathrm{BLLL}_{#2}^{#1}}}  

\newcommand*{\height}[0]{\ensuremath{\mathrm{hgt}}}

\newcommand*{\TMCells}[0]{\ensuremath{\mathrm{Cells}}}
\newcommand*{\TMPos}[0]{\ensuremath{\mathrm{Pos}}}
\newcommand*{\TMFrom}[0]{\ensuremath{\mathrm{From}}}
\newcommand*{\TMLeft}[0]{\ensuremath{\mathrm{Left}}}
\newcommand*{\TMRight}[0]{\ensuremath{\mathrm{Right}}}
\newcommand*{\TMBefore}[0]{\ensuremath{\mathrm{Before}}}
\newcommand*{\TMState}[0]{\ensuremath{\mathrm{State}}}
\newcommand*{\TMTape}[0]{\ensuremath{\mathrm{Tape}}}
\newcommand*{\TMUpdate}[0]{\ensuremath{\mathrm{Update}}}
\newcommand*{\TMInit}[0]{\ensuremath{\mathrm{Init}}}

\newcommand*{\TMHalting}[0]{\ensuremath{\mathrm{Halting}}}

\newcommand*{\Fraisse}[0]{Fra\"iss\'e\xspace}

\newcommand*{\Runs}[1]{\mathrm{Runs}({#1})}

\newcommand*{\ReturnSimulator}[3]{\ensuremath{\mathrm{Rt}^{#3}_{#2}({#1})}}

\newcommand*{\SimToReturnSym}[1]{\ensuremath{\mathrm{sTr}_{#1}}} 
 
\newcommand*{\SimToLoopSym}[1]{\ensuremath{\mathrm{sTl}_{#1}}}

\newcommand*{\HONPT}[0]{\ensuremath{\mathrm{NPT}}\xspace}
\newcommand*{\width}[0]{\mathrm{wdt}}

\newcommand*{\RelAnc}[2]{\ensuremath{{\mathrm{RA}_{#1}({#2})}}}
\newcommand*{\AncestorClass}[4]{\mathfrak{rA}_{#1, #2, #3, #4}}
\newcommand*{\RelAncequiv}[4]{\mathrel{{_{#4}^{#1}}{\equiv}_{#2}^{#3}}}

\newcommand*{\Lin}[3]{\ensuremath\mathfrak{Lin}_{#1}^{#2;#3}}
\newcommand*{\Typ}[3]{\ensuremath\mathrm{Type}_{#1}^{#2;#3}}
\newcommand*{\wordequiv}[2]{\equiv_{#1}^{#2}}
\newcommand*{\stackequiv}[3]{\mathrel{{_{#3}}{\equiv}_{#1}^{#2}}}

\newcommand*{\stackequivTyp}[3]{{\stackequiv{#1}{#2}{#3}}\text{-}\mathrm{Type}}

\newcommand*{\BoundHeight}[0]{\mathrm{BH}} 
\newcommand*{\BoundWidth}[0]{\mathrm{BW}}
\newcommand*{\BoundRunLength}[0]{\mathrm{BL}}

\newcommand*{\BoundHeightOnestepConstruction}[0]{\mathrm{f}}
\newcommand*{\BoundHeightOnestepConstructionSimultanious}[0]{\mathrm{BH_{1}}}

\newcommand*{\FuncBoundTopWord}[0]{\mathrm{BTW}}
\newcommand*{\FuncBoundWidthWord}[0]{\mathrm{BWW}}
\newcommand*{\ConstBoundHeightWord}[0]{\mathrm{B_{hgt}}}

\begin{document}

\thesistitle{First-Order Model Checking on Generalisations of Pushdown
  Graphs}{} 
\author{Dipl.- Math. Alexander Kartzow}
\birthplace{Gie\ss en}
\date{Juli 2011}
\referee{Prof. Dr. Martin Otto}{Prof. Damian
  Niwi\'nski}[Prof. Dr. Stephan Kreutzer]
\department{Fachbereich Mathematik}
\group{Arbeitsgruppe Logik}
\dateofexam{09.12.2010}{11.05.2011}

\makethesistitle

\chapter*{Acknowledgement}
I am deeply grateful to my supervisor Martin Otto for his support. 
Beside his mathematical advice, I especially appreciated his lessons
in mathematical writing and his efforts for improving my English. 
I thank my referees Damian Niwi\'nski and Stephan Kreutzer for their
valuable comments on this work. 
Furthermore, I thank Achim Blumensath and Dietrich Kuske for many
helpful comments and the opportunity to discuss some of my ideas. 
I am grateful to Alex Kreuzer and my wife Franziska for spell checking
parts of this thesis. 
Finally, I thank my wife and my family for the moral support and the
DFG for the financial support during the last years. 

\chapter*{German Summary / Zusammenfassung}

In dieser Arbeit untersuchen wir das Model-Checking-Problem f\"ur
Pushdown-Graphen. Ein Model-Checking-Algorithmus  
 f\"ur eine Logik
$\mathcal{L}$ und eine Klasse von Strukturen $\mathcal{C}$  ist ein
Algorithmus, der bei Eingabe eines Paares $(\mathfrak{A},\varphi)$ mit
$\mathfrak{A}\in\mathcal{C}$ und $\varphi\in\mathcal{L}$ entscheidet,
ob die Struktur $\mathfrak{A}$ die Formel $\varphi$ erf\"ullt.

In dieser Arbeit konzentrieren wir uns gr\"o\ss tenteils auf die
Entwicklung von Model-Checking-Algorithmen f\"ur die Logik erster Stufe
(im folgenden FO abgek\"urzt)
und ihrer Erweiterung um Erreichbarkeitspr\"adikate auf Klassen
verallgemeinerter Pushdown-Graphen. 

Ein Pushdown-Graph ist der Konfigurationsgraph eines
Kellerautomaten. Kellerautomaten, die auch Pushdown-Systeme genannt
werden, sind endliche Automaten erweitert 
um die Speicherstruktur eines Stacks. 
Ein klassisches Resultat von Muller und Schupp \cite{MullerS85} beweist
die Entscheidbarkeit des Model-Checking-Problems f\"ur die monadische
Logik zweiter Stufe (im folgenden MSO abgek\"urzt) auf der Klasse der
Pushdown-Graphen.  
Insbesondere gibt es also auch einen Model-Checking-Algorithmus f\"ur
die Logik erster Stufe auf der Klasse der Pushdown-Graphen. 

In den letzten Jahren haben Verallgemeinerungen der Pushdown-Graphen
gro\ss es Interesse im Bereich
der automatischen Verifikation von funktionalen
Programmiersprachen erlangt. 
Pushdown-Graphen wurden im wesentlichen auf zwei Arten erweitert. 

Die erste Erweiterung f\"uhrt zum Konzept eines Pushdown-Systems
h\"oherer Ordnung. Hierbei wird der Stack eines Kellerautomaten
ersetzt durch eine Struktur ineinander geschachtelter Stacks.
Die Verschachtelungstiefe dieser Stacks wird dabei als die Stufe des
Systems bezeichnet. Ein Pushdown-System der Stufe $2$ hat also einen
Stack aus Stacks, ein System der Stufe $3$ einen Stack aus Stacks aus
Stacks und analog f\"ur jede Stufe $n\in\N$. 
Auf jeder Stufe $i\leq n$ dieser Schachtelung gibt es entsprechende
Stack-Operationen um den obersten Eintrag des Stufe $i$ Stacks zu
manipulieren. Mit diesem Ansatz wurden zwei Hierarchien
verallgemeinerter Pushdown-Graphen definiert. Die Hierarchie der
``Higher-Order-Pushdown-Graphen'' und die der
``Collapsible-Pushdown-Graphen''. Die beiden Klassen unterscheiden
sich in den verwendeten 
Stack-Operationen.  Die Pushdown-Systeme, die
Collapsible-Pushdown-Graphen erzeugen, erweitern die Pushdown-Systeme,
die Higher-Order-Pushdown-Graphen erzeugen, um eine neue Operation, die
``Collapse'' genannt 
wird. Trotz der \"ahnlichen Definition dieser beiden Hierarchien von
Graphen haben die Hierarchien sehr unterschiedliche
modelltheoretische Eigenschaften. 

Die Hierarchie der Higher-Order-Pushdown-Graphen f\"allt mit der
Caucal-Hierarchie zusammen. Diese Klasse von 
Graphen ist definiert durch iteriertes Anwenden von MSO-Interpretationen
und Abwicklungen beginnend von der Klasse der endlichen Graphen. 
Da sowohl Abwicklungen als auch MSO-Interpretationen die
Entscheidbarkeit von monadischer Logik zweiter Stufe erhalten, ist
MSO-Model-Checking  auf der Klasse der Higher-Order-Pushdown-Graphen
entscheidbar. 

Die Klasse der Collapsible-Pushdown-Graphen hat dagegen ganz andere
modelltheoretische Eigenschaften. Schon auf der zweiten Stufe dieser
Hierarchie gibt es Graphen mit unentscheidbarer  MSO-Theorie. 
Hingegen ist der modale $\mu$-Kalk\"ul auf der Klasse der 
Collapsible-Pushdown-Graphen entscheidbar. 
Dieses unterschiedliche Verhalten in Bezug auf MSO und $\mu$-Kalk\"ul
tritt nur bei wenigen nat\"urlichen Strukturklassen auf. 

Eine weitere Klasse mit dieser Eigenschaft erhalten wir durch die
zweite Verallgemeinerung von Pushdown-Graphen. 
Abwicklungen von Pushdown-Graphen haben sich in der
Software-Verifikation als n\"utzliche Abstraktion von
Programmabl\"aufen herausgestellt. Hierbei wird auf dem Stack vor
allem der Aufruf von Funktionen und die R\"uckkehr zum aufrufenden
Programm verwaltet. Viele interessante Eigenschaften von Programmen
lassen sich so durch MSO-Model-Checking auf der Abwicklung eines
Pushdown-Graphen \"uberpr\"ufen und nachweisen. 
Allerdings ist es in diesem Modell nicht m\"oglich, den Zustand des
Programms vor einem Funktionsaufruf mit dem Zustand am Ende dieser
Funktion zu vergleichen, denn in monadischer Logik zweiter Stufe kann
man bei 
unbeschr\"ankt verschachteltem Aufruf von Funktionen die
zusammenhgeh\"orenden Positionen von Funktionsaufruf und Funktionsende
nicht definieren. 

Um dieses Problem zu umgehen haben Alur et al.\ \cite{Alur06languagesof} die
Klasse der ``Nested-Pushdown-Trees'' eingef\"uhrt (Warnung: wir
bezeichnen diese 
bewusst nicht als ``Nested-Pushdown-B\"aume'', weil es keine B\"aume
sind). Ein Nested-Pushdown-Tree ist die 
Abwicklung eines Pushdown-Graphen mit einer zus\"atzlichen Relation
$\hookrightarrow$. Diese verbindet eine Push-Operation des
Kellerautomaten mit der dazugeh\"origen Pop-Operation. Wenn man einen
Pushdown-Graphen also als abstraktes Modell des Programmablaufs eines
Computerprogramms sieht, wird der Funktionsaufruf \"uber
$\hookrightarrow$ mit dem Ende der aufgerufenen Funktion verbunden. 
Mit diesem Modell kann man also die oben erw\"ahnten Nachteile der
Pushdown-Graphen \"uberwinden. Alur et al.\ konnten zeigen, dass f\"ur die
Klasse der Nested-Pushdown-Trees das $\mu$-Kalk\"ul-Model-Checking
entscheidbar ist. Jedoch gibt es einen Nested-Pushdown-Tree mit
unentscheidbarer MSO-Theorie.  

Da die monadische Logik zweiter Stufe f\"ur
Collapsible-Pushdown-Graphen und f\"ur Nested-Pushdown-Trees
unentscheidbar ist, stellt sich 
die nat\"urliche Frage, welche  Fragmente der
monadischen Logik zweiter Stufe auf diesen Klassen entscheidbar sind. 

In unserer Arbeit geben wir daf\"ur die folgenden partiellen Antworten.
\begin{enumerate}
\item Auf der zweiten Stufe der Hierarchie der
  Collapsible-Pushdown-Graphen ist das FO-Model-Checking-Problem
  entscheidbar.  
  Genauer ist die Erweiterung von FO um regul\"are
  Erreichbarkeitspr\"adikate und Ramsey-Quantoren entscheidbar. 
  Wir beweisen dies, indem wir eine baumautomatische Repr\"asentation
  (vgl. Punkt 4)
  f\"ur jeden Collapsible-Pushdown-Graphen der zweiten Stufe erzeugen.
\item Das FO-Model-Checking-Problem auf der Klasse der
  Nested-Pushdown-Trees ist in zweifach exponentiellem Platz
  entscheidbar.  
  Zus\"atzlich  kann jeder Nested-Pushdown-Tree durch eine
  FO-Interpretation aus einem Collapsible-Pushdown-Graphen der Stufe
  2 erzeugt werden. Mithilfe dieser Interpretation k\"onnen wir auch
  die Theorie der Logik erster Stufe erweitert um das
  Erreichbarkeitspr\"adikat f\"ur jeden Nested-Pushdown-Tree
  entscheiden. 
\end{enumerate}
Neben diesen Resultaten \"uber bekannte Erweiterungen von
Pushdown-Graphen beinhaltet diese Arbeit auch die folgenden
Ergebnisse.  
\begin{enumerate}
\item[3.] Durch die Kombination der Idee der geschachtelten Stacks mit der
  Definition der Nested-Pushdown-Trees definieren wir eine neue
  Hierarchie der Nested-Pushdown-Trees h\"oherer Ordnung. 
  Ein Nested-Pushdown-Tree der Stufe $l$ ist die Abwicklung eines
  Pushdown-Graphen der 
  Stufe $l$ erweitert um eine neue Relation $\hookrightarrow$, die
  zusammengeh\"orende Push- und Pop-Operationen verbindet. 
  Wir beweisen, dass diese neue Hierarchie eng verwandt mit den
  Hierarchien der Higher-Order-Pushdown-Graphen und der
  Collapsible-Pushdown-Graphen ist. Alle Abwicklungen von 
  Higher-Order-Pushdown-Graphen sind in der neuen Hierarchie
  enthalten. Au\ss erdem lassen 
  sich alle Higher-Order-Nested-Pushdown-Trees durch
  FO-Interpretationen aus der Klasse der Collapsible-Pushdown-Graphen
  erzeugen. Durch diese Interpretation kann man
  Higher-Order-Nested-Pushdown-Trees der Stufe $l$ als besonders
  einfache Collapsible-Pushdown-Graphen der Stufe $l+1$ betrachten. 
  Wir zeigen dann, dass f\"ur die zweite Stufe dieser neuen Hierarchie ein
  FO-Model-Checking-Algorithmus existiert.
\item[4.] Wir zeigen in dieser Arbeit auch, dass die Erweiterung 
  der Logik erster Stufe um Ramsey-Quantoren  auf
  baumautomatischen Strukturen entscheidbar ist. 
  Baumautomatische Strukturen sind Strukturen, die sich durch endliche 
  Baumautomaten repr\"asentieren lassen. Ein Ramsey-Quantor ist von der
  Gestalt $\RamQ{n}\bar x (\varphi(\bar x))$. Eine solche Formel wird
  von einer Struktur $\mathfrak{A}$ erf\"ullt, wenn es eine unendliche
  Teilmenge $M\subseteq \mathfrak{A}$ gibt, so dass jedes  $n$-Tupel
  aus $M$, von dem je zwei Elemente paarweise verschieden sind,  die
  Formel $\varphi$ erf\"ullt.   
  Unser Beweis, der in Zusammenarbeit mit Dietrich Kuske entstand,
  verallgemeinert ein analoges Resultat f\"ur die Klasse 
  der wortautomatischen Strukturen. 
\end{enumerate}


\tableofcontents

\chapter{Introduction}
In this thesis, we investigate the first-order model checking problem for
generalisations of pushdown graphs. Our work is a contribution to the
classification of all graphs that have decidable first-order
theories. The classes of graphs that we study are collapsible pushdown
graphs and nested pushdown trees. These classes of graphs have the
following interesting model-theoretic properties. The
monadic second-order theory of a graph from these classes is not
decidable in general, while its modal $\mu$-calculus theory is
always decidable. 
Most other classes of graphs do not share these properties. In
most cases, a  natural class of graphs will either have
decidable monadic second-order and modal $\mu$-calculus theories or 
undecidable monadic second-order and modal $\mu$-calculus theories. 
We start by briefly recalling the history of generalisations of
pushdown graphs. These classes of graphs arise naturally in the field
of software verification for higher-order functional programmes. 

\section{Verification and Model Checking}
\label{sec:Verification}

Verification of hard- and software is concerned with the problem of
proving that a certain piece of hard- or software fulfils the task
for which it was designed. Since computer systems are more and more used
in safety critical areas, failure of a system can have severe
consequences. Thus,  verification of these systems is very 
important. The most successful approach to verification is the model
checking paradigm introduced by Clarke and Emerson \cite{ClarkeE81}. 
In model checking, one derives an abstract structure $\mathfrak{A}$
as a model of 
some piece of hard- or software and one specifies the requirements of
the system in a formula $\varphi$ from some logic $\mathcal{L}$. The
problem whether the system is correct then reduces to the problem
whether the abstract model $\mathfrak{A}$ of the system satisfies the
formula.  This is called a model 
checking problem.
If the model satisfies the formula, we write
$\mathfrak{A}\models\varphi$. In this terminology, the $\mathcal{L}$
model checking problem on 
some class $\mathcal{C}$ of structures asks 
on input a structure $\mathfrak{A}\in\mathcal{C}$ and a
formula $\varphi\in \mathcal{L}$ whether
\mbox{$\mathfrak{A}\models\varphi$}. Since the 1980's, 
model 
checking on 
finite structures has been developed and is nowadays used for
real-world hardware verification problems. For hardware, it is
sufficient to consider finite structures. Each piece of hardware has
a finite amount of storage capacity whence it can always be modelled
as a finite state system. On the other hand, software verification
requires the use of infinite models as abstractions because the
storage capacity of the
underlying hardware is a priori unbounded. Hence, software
verification naturally leads to 
model checking problems on infinite structures. Of course, model
checking on infinite structures is only possible for certain classes
of structures. Since we expect an algorithm to process the structures
involved as input, we need finite descriptions of these infinite
structures. Hence, model checking on infinite structures is
only interesting for classes of finitely representable structures. 
A further restriction is imposed by the question of decidability of
the model checking problem. A very expressive logic on a large class
of finitely represented structures will result in an undecidable model
checking problem (the halting problem can be formulated as a special
version of model checking on structures representing Turing
machines). 
Thus, there is a tradeoff between the choice of the class
$\mathcal{C}$ and the logic $\mathcal{L}$. It is important to
identify those 
pairs $(\mathcal{C}, \mathcal{L})$ for which a model checking
algorithm exists, 
i.e., for which the $\mathcal{L}$ model checking on $\mathcal{C}$ is
decidable.  

Various techniques have been developed to finitely represent infinite
structures. According to B\'ar\'any et al.\ \cite{BGR10},  these may be
classified into the following approaches.
\begin{itemize}
\item Algebraic representations: a structure is described as the
  least solution of some recursive equation in some appropriate
  algebra of structures. An example of this class are 
  vertex replacement equational graphs \cite{Courcelle90}.
\item Transformational or logical representations: the structure is
  described as 
  the result of applying finitely many transformations to some finite
  structure. A transformation in this sense is, e.g., the tree-unfolding,
  the Muchnik-Iteration, or some logical interpretation (see
  \cite{BlumensathCL07} for a survey).
\item Internal representations: an isomorphic copy of the structure is
  explicitly described 
  using transducers or rewriting techniques. In most cases a set of
  words or trees is used as the universe of the structure. 
  The relations are then represented by rewriting rules or by
  transducers that process tuples of elements from this set.
  Rewriting rules often appear in the disguise of transitions of some
  computational model. In this case the universe consists of
  configurations of some computational model. There is an edge from
  one configuration to another configuration if one step of the computation
  leads form the first to the second configuration.
\end{itemize}
There is no clear separation between the approaches because there are many
classes of structures that may be represented using techniques from
different approaches. 

In this thesis we will only deal with structures that have internal
representations. 
We investigate 
configuration graphs of different types of automata. 
The universe of such a graph consists of 
the set of configurations of an automaton and the relations are given
by the transitions from one configuration to another. 
Automata that may be used for this approach are, e.g., Turing
machines, finite automata, pushdown systems or collapsible pushdown
systems. In this thesis we study configuration graphs of collapsible
pushdown systems. We will introduce these 
systems later in detail. 
A pushdown system can be seen as a finite automaton equipped with a stack. 
A collapsible pushdown system uses a nested stack, i.e., a stack of stacks
of stacks of \dots of stacks instead of the ordinary stack. 
On each stack level, the collapsible pushdown system can
manipulate the topmost entry of its stack.

Another concept that plays a major role within
this thesis is the concept of a tree generated by some pushdown
system. 
This tree is obtained by applying a graph
unfolding to the configuration graph. 
This can be seen as a transformational representation of the graph
that starts from the underlying configuration graph. On the other
hand, it can also be seen as an internal representation:
the nodes of a graph are represented by the set of runs of the given
automaton and the relations of the structure are defined by rewriting
rules 
that transform a run of length $n$ into a run of length $n+1$ that
extends the first run. If a graph is the configuration graph of some
automaton, we will refer to the unfolding of this graph as the tree
generated by this automaton.  
This notion becomes important when we discuss  nested
pushdown trees. These are trees generated by pushdown systems expanded
by a so-called jump relation. We will present this concept at the end
of the next section.  

The second form of internal representation for infinite structures
that we will use are tree-automatic structures. 
A structure is tree-automatic if it can be
represented as a regular set of trees such that for each relation
there is a finite tree-automaton that accepts those tuples of trees
from the universe that form a tuple of the relation. We provide a more
detailed introduction to tree-automatic structures as well as some
notes concerning the history of tree-automatic structures in Section
\ref{Chapter_AutomaticStructures}. The class of tree-automatic
structures is a nice class because first-order model checking is
decidable on this class: there are automata constructions that
correspond to negation, conjunction and existential
quantification. Thus, for any 
tree-automatic structure and any first-order formula, one can
construct a tree-automaton that accepts an input (representing a tuple of
parameters from the structure) if and only if the structure satisfies
the formula (where the free variables of $\varphi$ are assigned to the
parameters represented by the input).

\section{Collapsible Pushdown Graphs and Nested Pushdown Trees}
\label{sec:IntroCPGNPT}

The history of software verification is closely connected to two
important results on model checking. 
In 1969, Rabin \cite{Rab69} proved the decidability of  monadic second-order
logic (\MSO) on the infinite binary tree. In terms of model checking, his
result states that the \MSO model checking is decidable for the class
that only consists of one structure, namely, the full binary tree. 
Sixteen years later, Muller and Schupp \cite{MullerS85} showed the
decidability of the \MSO model checking on pushdown
graphs. This was a very important step towards
automated software verification 
because  pushdown graphs proved to be very
suitable for modelling procedural programmes with calls of first-order
recursive procedures. The function calls and returns are modelled using
the stack. At a function call, the state of the programme is pushed onto
the stack and at a return the old context is restored using a
pop operation. 

Collapsible pushdown systems can be seen as the result of the search
for a similar result for higher-order functional programming
languages. Already in the 1970's Maslov was the first to consider
so-called higher-order pushdown systems as accepting devices for word
languages. A higher-order pushdown system is a generalisation of a
pushdown system where one replaces the stack by a nested stack of
stacks of stacks of \ldots stacks. For each stack level the
higher-order pushdown system can use a push and a pop operation. 
In the last years, these automata have become an important topic of
interest because of two results.
\begin{enumerate}
\item Carayol and W\"ohrle \cite{cawo03} showed that the class of
  graphs generated by $\varepsilon$-contractions of configuration
  graphs of higher-order pushdown systems coincide with the class of
  graphs in the Caucal hierarchy. Caucal \cite{Caucal02} defined this
  class as follows. The initial level in the hierarchy contains all
  finite graphs. A graph in the next level is obtained by applying an
  unfolding and an \MSO-interpretation to a graph in the previous
  level. Since both operations preserve the \MSO decidability, \MSO
  model checking on higher-order pushdown graphs is decidable. In
  fact, the Caucal hierarchy is one of the largest classes where the
  \MSO model checking is known to be decidable. 
\item Knapik et al.\ \cite{KNU02} studied higher-order pushdown systems
  as generators of trees. They proved that the class of
  trees generated by higher-order pushdown systems coincides with the
  class of trees generated by \emph{safe} higher-order recursion
  schemes (safe higher-order functional programmes). Safety is a
  rather syntactic condition on the types of in- 
  and outputs to functions that are used in a recursion scheme. 
\end{enumerate}

The second result initiated a lot of study on the question whether
there 
is some computational model whose generated trees form exactly the
class of trees generated by arbitrary higher-order recursion schemes
and whether the trees generated by safe recursion schemes form a
proper subclass of the class of trees generated by arbitrary recursion
schemes.   
For instance, Aehlig et al.\ \cite{RR-04-23} showed that
safety is no restriction for string languages defined by level $2$
recursion schemes. 
Hague et al.\ \cite{Hague2008}
introduced collapsible pushdown systems. The concept of a collapsible
pushdown system 
is a stronger variant of the concept of a
higher-order pushdown systems. They showed that these
are as expressive as arbitrary higher-order recursion schemes, i.e.,
a tree is generated by a level $n$ recursion scheme if and only if it
is generated by some level $n$ collapsible pushdown system. 
Furthermore, they showed the decidability of  modal
$\mu$-calculus ($L\mu$) model checking on collapsible pushdown
graphs. 
Recently, Kobayashi \cite{DBLP:conf/popl/Kobayashi09} designed an
$L\mu$ model 
checker for higher-order recursion schemes and successfully applied
this model checker to the verification of higher-order functional
programmes. 
Even though the connection to higher-order recursion schemes turns
collapsible pushdown systems into a very interesting class of
structure for model checking purposes, there are few things known about
the structure of the trees and graphs generated by these systems. 
For example, it is conjectured -- but not proved -- that the
class of trees generated by collapsible pushdown systems properly
extends the class of trees generated by higher-order pushdown
systems. The same conjecture in terms of recursion schemes says that
there is a tree generated by some  unsafe higher-order recursion
scheme that is not generated by any  safe higher-order recursion
scheme.  

Concerning model checking, Hague et al.\ proved another
interesting fact about the class of graphs 
generated by collapsible pushdown systems: they presented a
collapsible pushdown system of level $2$ that has undecidable \MSO
theory. In terms of model checking, this is a proof of the fact that
\MSO model checking on the class of collapsible pushdown graphs is
undecidable. 

From a theoretical point of view, this turns collapsible pushdown
graphs into an interesting class of graphs. Besides the class of
nested pushdown trees it is the only
known natural class of graphs that has decidable $L\mu$ model checking
but undecidable \MSO model checking. Thus, a better understanding of
this class of 
graphs  may also give insight into the difference between
$L\mu$ 
and \MSO. In fact, this thesis tries to identify larger fragments of
\MSO that are still decidable on collapsible pushdown graphs. The most
prominent fragment of \MSO is, of course, first-order logic
(\FO{}). The author was the first to investigate first-order model
checking on collapsible pushdown graphs. 
In STACS'10 \cite{Kartzow10}, we proved that the model checking problem
for the extension of \FO{} by reachability predicates on the class of
collapsible pushdown graphs of level $2$ is decidable.
In this thesis, we present a slightly 
extended version of this result: if one enriches the  graphs by
$L\mu$-definable predicates, model checking is still decidable.  
Furthermore, we may also enrich \FO{} by Ramsey quantifiers.
Very recently, Broadbent
\cite{Broadbent2010Mail} matched our result with a tight upper
bound:  the first-order model
checking on level $3$ collapsible pushdown
graphs is undecidable. Moreover, Broadbent presented a fixed formula 
$\varphi\in\FO{}$ such that the question whether $\varphi$ is
satisfied by some level $3$ collapsible pushdown graph is
undecidable. Furthermore, he provided an example of a level $3$
collapsible pushdown graph with undecidable $\FO{}$-theory. 

We now turn to the history of nested pushdown trees. 
Alur et al.\ \cite{Alur06languagesof} introduced the concept of
so-called jump edges in order to overcome the following weakness of
model checking on pushdown graphs. Recall that pushdown
systems are useful abstractions of programmes which call first-order
recursive functions. Function calls and returns are handled by using
push and pop operations. But interesting properties of some programme
may include statements about the situation before a function call
happens in comparison to the situation at the end of this function,
i.e., at the return of this function. 
Unfortunately, even strong logics
like \MSO cannot express such properties. They cannot ``find'' the
exact corresponding pop operation for a given push operation in
general. As soon
as  a potentially unbounded nesting of function calls may occur, 
 \MSO like many 
other logics cannot keep track of the number of nestings in the call and
return structure. But this would be necessary
for identifying the pop operation that corresponds to a given push. 

Alur et al.\ wanted to make this correspondence of push and pop
operations explicit. Thus, a nested pushdown tree is defined to be the
unfolding of a pushdown graph enriched by jump edges that connect each
push operation with the corresponding pop operation. 
Unfortunately, this expansion of trees generated by pushdown systems
leads to undecidability of the \MSO model
checking \cite{Alur06languagesof}. Anyhow, Alur et al.\ were able to prove
that $L\mu$ model checking is still decidable on nested pushdown
trees. Thus, the class of nested pushdown trees is the second natural
class of structures with undecidable \MSO but decidable $L\mu$ model
checking. We were able to provide an elementary \FO{} model checking
algorithm for nested pushdown trees. This result was first presented
in MFCS'09 \cite{Kartzow09}. 

The similar behaviour of the class of nested pushdown trees and
collapsible pushdown graphs with respect to model checking has an easy
explanation. 
 Nested pushdown trees are
first-order interpretable in collapsible pushdown graphs of level
$2$. Furthermore, the interpretation is quite simple and
uniform.

\section{Goal and Outline of this Thesis}

This thesis is concerned with various model checking problems. 
Our most important results provide model checking algorithms for 
first-order logic (and slight extensions) 
on various classes of structures. 
The main focus is on structures defined
 by higher-order (collapsible) pushdown
systems. On the one hand, we study the hierarchy of collapsible
pushdown graphs that was introduced by Hague et al.\ \cite{Hague2008}. On
the other hand, we study a new hierarchy of \emph{higher-order nested
  pushdown trees}. This hierarchy is the class obtained by the
straightforward generalisation of the concept of a nested pushdown
tree to trees generated by higher-order pushdown systems. We consider
the expansions of these trees by jump-edges that connect corresponding
push and pop transitions (at the highest level of the underlying
higher-order pushdown system). 
This new hierarchy forms a class of graphs that contains the class of
trees generated by higher-order pushdown systems and that is contained (via
uniform \FO{}-interpretations) in the class of collapsible pushdown
graphs. Thus, we hope that the study of this new hierarchy can
reveal some insights into the differences between these two
hierarchies. 

In this thesis, we obtain the following  results on the model checking
problems for these hierarchies.
\begin{enumerate}
\item The second level of the collapsible pushdown hierarchy is
  tree-automatic and its \FO{}($\Reach$) model checking is decidable.
\item First-order model checking on nested pushdown trees is in
  $2$-EXPSPACE.
\item First-order model checking on level $2$ nested pushdown trees is
  decidable. 
\end{enumerate}
In order to prove these claims, we develop various new techniques. 

All of these proofs
rely on a structural analysis of runs of
higher-order collapsible pushdown systems. This analysis provides a
characterisation of the reachability of one configuration from
another. 

The second ingredient for our first result is 
a clever encoding of configurations in trees which turns the set of
reachable configurations into a regular set of trees. 

The other two
results use a new application of 
Ehrenfeucht-\Fraisse games to the model checking problem. 
We analyse strategies in the
Ehrenfeucht-\Fraisse game 
that are subject to certain restrictions. The existence of such restricted
winning strategies on a class of structures can be used to provide a
model checking algorithm on this class.
The basic idea is as follows: assume 
that Duplicator 
has a strategy that only requires to consider finitely many elements
in a structure. Model checking  on this structure 
can then be reduced  to model checking on a finite
substructure, namely, on the substructure induced by those elements
that are relevant for Duplicator's strategy. 

Using our analysis of runs of collapsible pushdown systems, we show
that there are such restricted strategies on the first two levels of
the nested pushdown hierarchy. 

Motivated by the tree-automaticity of level $2$ collapsible pushdown
graphs, we also study the model checking problem on the class of all
tree-automatic structures. We provide an extension of the known
first-order model checking algorithm to Ramsey quantifiers. 
These are also called Magidor-Malitz quantifiers because these
generalised quantifiers were first introduced by Magidor and Malitz  
\cite{MagidorM77}. 

The proof of this result is given by an explicit automata-construction
that corresponds to this quantifier. 
For the string-automatic structures, such a proof was given by 
Rubin \cite{Rubin2008} using 
the concept of word-combs. 
Rubin then proved that,
on string-automatic structures, 
each set witnessing a Ramsey quantifier contains a
word-comb. Using the theory of $\omega$-string-automata, he then uses
word-combs to design a finite string-automaton corresponding to the
Ramsey quantifier. 
In joint work with Dietrich Kuske, 
we extended this result to the tree-case. We define the concept of a
tree-comb  
and use $\omega$-tree-automata in order to provide a
finite tree-automata construction that corresponds to the Ramsey
quantifier on a tree-automatic structure. We stress that our result is
a nontrivial adaption of Rubin's work. The technical difference
between the string and the tree case is based on the fact that strings
have a uniquely defined length, while the lengths of paths in a
tree are not necessarily uniform.  

\paragraph{Outline of this Thesis}

In Chapter \ref{Chapter_Basics}, we first review all basic concepts
that are necessary for understanding this thesis. Namely, we review
different logics, logical
interpretations  and the concepts of trees and words. 
We also revisit the theory of Ehrenfeucht-\Fraisse games and develop a
new model checking approach based on the analysis of restricted
strategies in these games. 
After these
preliminaries, we introduce our objects 
of study. In Section \ref{Chapter_InfiniteStructures}, we introduce
higher-order pushdown systems, collapsible pushdown systems and nested
pushdown trees. After this, we provide some  
technical results on runs of collapsible pushdown systems
in Section \ref{ChapterLoops}. 
These results concern the existence and computability of certain runs
of level $2$ collapsible pushdown systems. 
The technical lemmas provided in this section play a crucial role in
proving our results concerning level $2$ (collapsible) pushdown systems. 
In Section \ref{Chapter_AutomaticStructures}, we
review the basic concepts and results on tree-automatic
structures. 
At the beginning of Chapter \ref{Chapter_MainResults} we briefly
present our 
main results in the following order.
\begin{enumerate}
\item The second level of the collapsible pushdown hierarchy is
  tree-automatic and its \FO{}+$\Reach$ theory is decidable.
\item First-order model checking on nested pushdown trees is in
  $2$-EXPSPACE.
\item First-order model checking on level $2$ nested pushdown trees is
  decidable. 
\item The model checking problem for \FO{} extended by Ramsey
  quantifiers  on 
  tree-automatic structures is decidable. 
\end{enumerate}
For each of these results
there is one section in Chapter \ref{Chapter_MainResults}  providing
the details of the proof 
and some discussion on related topics. 
Note that we postpone the formal definition of the hierarchy of
higher-order nested pushdown trees to Section
\ref{Chapter HONPT}. 
In that section, we relate this new hierarchy to the hierarchy of
higher-order pushdown graphs and to the hierarchy of collapsible
pushdown graphs.  
Finally, Chapter \ref{Chapter_Conclusions} contains concluding remarks
and some open problems.


\chapter{Basic Definitions and Technical Results}
\label{Chapter_Basics}
In the first part of this  chapter, we review different kinds of
logics and logical interpretations that will play a role in this
thesis. Most of this part is 
assumed to be known to the reader and is merely stated for fixing
notation. An exception to this rule is the part on
Ehrenfeucht-\Fraisse games. 
First, we briefly recall the definition and some well-known facts about 
Ehrenfeucht-\Fraisse games. Afterwards, we introduce a new 
application of these games to first-order model checking problems. 
We develop an approach for model checking via the analysis of
restricted strategies in 
the Ehrenfeucht-\Fraisse game played on two identical copies of a fixed
structure. If Duplicator has  winning strategies that satisfy certain
restrictions  on each structure of some class $\mathcal{C}$, then we
can turn these strategies into an $\FO{}$ model checking algorithm on
$\mathcal{C}$.  

In Section \ref{Sec:StrucandInt} we review the notions of
grids and trees. 
Grids only play a minor role for our results. We use a certain grid-like
structure as a counterexample in an undecidability proof. In contrast,
trees play  a crucial role for our first two main results.

Section \ref{Chapter_InfiniteStructures}
is an introduction to  collapsible pushdown graphs
and nested pushdown trees. The main focus of this thesis is on model
checking algorithms for the classes of these graphs. 
As a preparation for the development of these algorithms, we present
the most 
important tool for our results in Section \ref{ChapterLoops}. In that
section we give
a detailed analysis of the 
structure of runs of collapsible pushdown graphs of level $2$. 
Finally, in Section \ref{Chapter_AutomaticStructures} we  recall the
necessary notions 
concerning tree-automatic structures. 
Note that tree-automatic structures play
two different roles in this thesis:
our first main result studies
the class of tree-automatic structures on its own. 
We provide a model checking algorithm for 
first-order logic extended by Ramsey quantifiers
(also called Magidor-Malitz quantifiers) on this class. 
Our algorithm extends the known first-order model checking algorithm
on tree-automatic structures. 

In the second main
result, we use 
tree-automaticity as a tool. We show that collapsible
pushdown graphs of level $2$ are tree-automatic. Thus, they 
inherit the decidability of the  first-order model checking problem
from the general theory of tree-automatic structures.

\section{Logics and Interpretations}

In this section we briefly recall the definitions of the logics we are
concerned with. These are classical first-order logic and its
extensions by monadic second-order quantifiers, certain generalised
quantifiers, reachability predicates, or least fixpoint operators.
Furthermore, we present basic modal logic and the modal
$\mu$-calculus (denoted by $L\mu$), which is the extension of modal
logic by least fixpoint operators. 
The last part of this section also fixes our notation concerning
logical interpretations. 

\subsection{First-Order Logic, Locality and Ehrenfeucht-\Fraisse Games}

\paragraph{Vocabularies and Structures}

For reasons of convenience, we only introduce \emph{relational vocabularies}
and \emph{relational structures} because we are only concerned
with such structures.  
A \emph{vocabulary} (or \emph{signature})
$\sigma=( (R_i)_{i\in I})$ consists of
relation symbols $R_i$. Each relation symbol $R_i$ has a fixed arity
$\arity(R_i)\in\N$.

A \emph{$\sigma$-structure $\mathfrak{A}$} is a tuple
$(A,(R_i^{\mathfrak{A}})_{i\in I})$ where $A$ is a set called the \emph{universe}
of $\mathfrak{A}$, and
$R_i^{\mathfrak{A}}\subseteq A^{\arity(R_i)}$ is a relation of arity
$\arity(R_i)$ for each $i\in I$. 
We denote structures with the letters
$\mathfrak{A}$, $\mathfrak{B}$, $\mathfrak{C}$, and so on. We
silently assume that the universe of $\mathfrak{A}$ is a set $A$,
the universe of $\mathfrak{B}$ is a set $B$, etc. 

We introduce the following notation concerning elements of the
universe of a structure. 
For  some structure $\mathfrak{A}$, we use the notation
$a\in\mathfrak{A}$ for stating that $a$ is some element of the
universe $A$ of $\mathfrak{A}$. Furthermore, we use a sloppy notation
for tuples of elements. We write $\bar a:= a_1, a_2, \dots, a_n \in
\mathfrak{A}$ for $\bar a:=(a_1, a_2, \dots, a_n)\in A^n$. 

\paragraph{First-Order Logic}
Let $\sigma$ be a vocabulary. We denote by
$\FO{}(\sigma)$ \emph{first-order logic} over the
vocabulary $\sigma$. Formulas of $\FO{}(\sigma)$ are composed by
iterated use of the following rules:
\begin{enumerate} 
\item for $x, y$ some variable symbols, $x=y$ is a formula in
  $\FO{}(\sigma)$,
\item for $R_i\in \sigma$ a relation of arity $r:=\arity(R_i)$ and
  variable symbols $x_1, x_2, \dots x_r$, $R_i x_1 x_2
  \dots x_r$ is a formula in $\FO{}(\sigma)$,
\item for $\varphi, \psi\in\FO{}(\sigma)$, $\varphi \land \psi$,
  $\varphi \lor \psi$, and $\neg \varphi$ are formulas in
  $\FO{}(\sigma)$, 
\item for $\varphi\in \FO{}(\sigma)$ and $x$ a variable symbol, 
  $\exists x \varphi$ and $\forall x \varphi$ are formulas in
  $\FO{}(\sigma)$. 
\end{enumerate}
Let $\varphi\in \FO{}(\sigma)$ be a formula. We write $\Var(\varphi)$
for the set of variable symbols occurring in $\varphi$. The semantics
of first-order formulas is defined as follows. Let 
$\mathfrak{A}$ be a $\sigma$-structure with universe $A$ and 
$I:\Var(\varphi) \to A$ some function (called the variable assignment
or interpretation)
we write $\mathfrak{A},I \models \varphi$, and say $\mathfrak{A}, I$ is
a model of $\varphi$ (or $\mathfrak{A}, I$ satisfies $\varphi$), if
one of the following holds. 
\begin{enumerate}
\item $\varphi$ is of the form $x=y$ where $x,y$ are variable symbols
  and $I(x)= I(y)$.
\item $\varphi$ is of the form $R_i x_1 x_2 \dots x_r$ and
  $(I(x_1), I(x_2), \dots I(x_r))\in R_i^{\mathfrak{A}}$.
\item $\varphi$ is of the form $\psi \lor \chi$ for $\psi,\chi\in
  \FO{}(\sigma)$ and  
  $\mathfrak{A}, I{\restriction}_\Var(\psi) \models \psi$ or 
  $\mathfrak{A}, I{\restriction}_\Var(\chi) \models \chi$.
\item $\varphi$ is of the form $\psi \land \chi$ for $\psi, \chi\in
  \FO{}(\sigma)$ and 
  $\mathfrak{A}, I{\restriction}_\Var(\psi) \models \psi$ and
  $\mathfrak{A}, I{\restriction}_\Var(\chi) \models \chi$.
\item $\varphi$ is of the form $\neg \psi$ for $\psi\in \FO{}(\sigma)$
  and $\mathfrak{A}, I \not\models \psi$.
\item $\varphi$ is of the form $\exists x \psi$ and there is some
  $a\in A$ such that
  $\mathfrak{A},I_{x\mapsto a} \models \psi$ where 
  \mbox{$I_{x\mapsto a}:\Var(\psi) \to A$} with 
  $I_{x\mapsto a}(y):= \begin{cases}
    I(y) &\text{for } y\neq x,\\
    a    &\text{for } y=x.
  \end{cases}$
\item $\varphi$ is of the form $\forall x \psi$ and 
  $\mathfrak{A},I_{x\mapsto a} \models \psi$ for all 
  $a\in A$. 
\end{enumerate}
We denote by $\Freevar(\varphi)\subseteq\Var(\varphi)$ the set of
variables occurring free in $\varphi$. A variable $x$ does not occur
free if it only occurs under the scope of quantifiers $\exists x$ or
$\forall x$. 
If \mbox{$\Freevar(\varphi)\subseteq \{x_1, x_2, \dots, x_n\}$} we 
use the notation 
\begin{align*}
  \mathfrak{A}, a_1, a_2, \dots, a_n \models \varphi(x_1, x_2, \dots,
  x_n)  
\end{align*}
 for 
$\mathfrak{A}, I \models \varphi$ if $I$ maps $x_i$ to
$a_i$. Furthermore, if $x_1, x_2, \dots, x_n$ are clear from the
context, we also use the notation
\mbox{$\mathfrak{A} \models \varphi(a_1, a_2, \dots, a_n)$}. 

In the following, we  write $\FO{}$ instead of $\FO{}(\sigma)$
whenever $\sigma$ is clear from the context or if a  statement does
not depend on the concrete $\sigma$.  
We may assign to each formula in $\FO{}$ its quantifier
rank. This is the maximal nesting depth of existential and universal
quantifications in this formula. 
We write $\FO{\rho}$
for the restriction of $\FO{}$ to formulas of quantifier rank up to
$\rho$.  

Let $\mathfrak{A}$ and $\mathfrak{B}$ be structures. For 
$n$ parameters $\bar a\in \mathfrak{A}$ 
and $n$ parameters $\bar b\in \mathfrak{B}$, we write 
$\mathfrak{A},\bar a \equiv_\rho \mathfrak{B},\bar b$ for the fact
that  $\mathfrak{A}\models \varphi(\bar a)$ if and only if
$\mathfrak{B}\models\varphi(\bar b)$ 
for all $\varphi\in\FO{\rho}$ with free variables among 
$x_1, x_2, \dots, x_n$.

We conclude the section on first-order logic by recalling two
important concepts for the analysis of first-order theories. Firstly,
we present the concepts of Gaifman locality and
Gaifman graphs. Afterwards, we present Ehrenfeucht-\Fraisse games
which are a classical tool for the analysis of $\equiv_\rho$. In
this thesis, we develop a nonstandard application of these games for
the design of model checking algorithms. 

\paragraph{Gaifman-Locality}

First-order logic has a local nature, i.e., first-order formulas can
only express properties about local parts of structures. For example,
reachability along a path of some relation $E$ is not first-order
expressible. Gaifman introduced the notions of Gaifman graphs and
local neighbourhoods in order to give a precise notion of the 
local nature of first-order logic. Let us start by recalling 
these notions.

\begin{definition}
  Let $\sigma=(R_1, R_2, \dots, R_n)$ be a finite relational signature and
  let $\mathfrak{A}$ be a 
  $\sigma$-structure. The \emph{Gaifman graph $\mathcal{G}(\mathfrak{A})$} is
  the graph $(A,E)$ where $A$ is the universe of $\mathfrak{A}$ and
  $E\subseteq A^2$ is the relation defined as follows. 
  $E$ connects two distinct elements
  of $A$ if 
  they appear together in a tuple of some relation of
  $\mathfrak{A}$, i.e., $(a_1,a_2)\in E$ for $a_1\neq a_2$ if and only
  if there is an $1\leq 
  i \leq n$ and tuples $\bar x, \bar y, \bar z \in A$ such that
  \begin{align*}
    &\mathfrak{A} \models R_i \bar x a_1 \bar y a_2 \bar z\text{ or}\\
    &\mathfrak{A} \models R_i \bar x a_2 \bar y a_1 \bar z.    
  \end{align*}
  For $a_1, a_2 \in A$ we say $a_1$ and $a_2$ have distance $n$ in
  $\mathfrak{A}$, written $\mathrm{dist}(a_1,a_2)=n$, if their distance
  in $\mathcal{G}(\mathfrak{A})$ is $n$. Analogously, we use the
  terminology $\mathrm{dist}(a_1,a_2)\leq n$ with the obvious
  meaning. Note that $\mathrm{dist}(x_1,x_2)\leq n$ is first-order
  definable for each fixed $n\in \N$.\footnote{Note that the
    restriction to finite vocabulary is essential for this statement.}

  We define the \emph{$n$-local neighbourhood} of some tuple $\bar a = a_1,
  a_2, \dots, a_n \in A$
  inductively by 
  \begin{align*}
    \mathcal{N}_0(\bar a)&:=\{a_1, a_2, \dots, a_n\} \text{ and}\\
    \mathcal{N}_{n+1}(\bar a)&:= \{a\in A: \exists
    a'\in\mathcal{N}_n(\bar a) \text{ such that }
    \mathrm{dist}(a,a')\leq 1\}.
  \end{align*}
\end{definition}
When we say that  ``first-order logic  is Gaifman-local'', we refer to
the fact that for each quantifier rank $\rho$, there is a
natural number $n$ such that for each formula $\varphi\in\FO{\rho}$
the question whether a structure $\mathfrak{A}$ is a model of $\varphi$ only
depends on the $n$-local neighbourhoods of the elements in the
structure $\mathfrak{A}$. More precisely, any first-order formula is a
boolean combination of \emph{local formulas} and \emph{local
  sentences} which we introduce next. 

\begin{definition}
  Let $\sigma$ be a finite vocabulary and let
  $\varphi\in\FO{}$ be some formula with free variables $\bar
  x$. We write $\varphi^n$ for the \emph{relativisation of $\varphi$ to the
  $n$-local neighbourhood of $\bar x$}. $\varphi^n$ is
  obtained from $\varphi$ by replacing each quantifier 
  $\exists y(\psi)$ by 
  \begin{align*}
    \exists y (\mathrm{dist}(y,\bar x)\leq n \wedge \psi)    
  \end{align*}
  and each quantifier $\forall y(\psi)$ by 
  \begin{align*}
    \forall y (\mathrm{dist}(y,\bar
    x)\leq n \rightarrow \psi).
  \end{align*}
  This means that for each variable assignment $I$ that maps $\bar x
  \mapsto \bar a\in A$, $\mathfrak{A}, I \models \varphi^n(\bar
  x)$ if and only if  
  $\mathfrak{A}{\restriction}_{\mathcal{N}_n(\bar a)}, I\models \varphi$.

  We call $\varphi(x)$ an \emph{$n$-local formula} if 
  $\varphi(x) \equiv \varphi^n(x)$ and we call it \emph{local} if it
  is $n$-local for  some $n\in\N$.

  We call a sentence local if it is of the form 
  \begin{align*}
    \exists x_1 \dots \exists x_n \bigwedge_{1\leq i<j\leq n}
    \mathrm{dist}(x_i,x_j) > 2l \wedge \bigwedge_{1\leq i \leq n}
    \psi^l(x_i)    
  \end{align*}
  for some formula $\psi$ and some $l\in\N$. 
  Such a sentence asserts that there are $n$ elements far apart from
  one another each satisfying the $l$-local formula $\psi^l$.
\end{definition}
Using this notation we can state Gaifman's
Lemma. 
\begin{lemma}
  [\cite{Gaifman82}]
  Each first-order formula is equivalent to a boolean combination of
  local sentences and local formulas. 
\end{lemma}
This lemma has an interesting consequence. For each quantifier rank
$\rho\in\N$ there is a natural number $n\in\N$ such that the following
hold: if $\bar a\in A$ and $a,a'\in A$ are such that
$\mathrm{dist}(a,\bar a) > 2 \rho +1$ and $\mathrm{dist}(a',\bar
a)>2\rho+1$ and 
there is an isomorphism $\mathcal{N}_{\rho}(a) \simeq \mathcal{N}_{\rho}(a')$
mapping $a$ to $a'$, then 
$\mathfrak{A}, \bar a, a \equiv_\rho \mathfrak{A}, \bar a, a'$. 

In Section \ref{SecGaifman}, we develop a lemma of a similar style. 
But in contrast to Gaifman's Lemma,
this new lemma is tailored towards an application on
graphs of small diameter. Due to the small diameter,
all elements $a$ and $\bar a$
satisfy
$\mathrm{dist}(a,\bar a)\leq 2n+1$ whence we cannot use 
Gaifman's Lemma itself. Nevertheless, in that section we need a lemma
that provides $\equiv_\rho$-equivalence for certain tuples $\bar a, a$
and $\bar a, a'$ in certain graphs of small diameter. We obtain this
lemma using 
Ehrenfeucht-\Fraisse games which we introduce in the following.

\paragraph{Ehrenfeucht-\Fraisse Games and First-Order Model Checking}
\label{Sec:EFGame}

The equivalence $\equiv_\rho$ of first-order logic up to quantifier
rank $\rho$ has a nice characterisation via
Ehrenfeucht-\Fraisse  
games. Based on the work of \Fraisse \cite{Fraisse54},
Ehrenfeucht \cite{Ehren60} 
introduced these games which have become one of the most important
tools for proving inexpressibility of properties in first-order
logic. This tool is especially important in the context of finite
model theory where other 
methods, e.g. compactness, fail. The game is played by two 
players, who are 
called Spoiler and Duplicator. They play on two $\sigma$-structures
$\mathfrak{A_1}$ and $\mathfrak{A_2}$. The players alternatingly choose 
elements in the two structures. At the end of the game, Duplicator has
won if there is a partial isomorphism between the elements chosen in
each of the structures. 
Thus, 
Spoiler's goal is to choose elements in such a way 
that no choice of Duplicator yields a partial isomorphism between the
elements chosen so far. 
The precise definitions are as follows.

\begin{definition}
  Let $\mathfrak{A}_1$ and $\mathfrak{A}_2$ be $\sigma$-structures. 
  For 
  \begin{align*}
    &\bar a^1=a^1_1, a^1_2, \dots, a^1_m\in A_1^m\text{ and}\\
    &\bar a^2=a^2_1, a^2_2, \dots, a^2_m\in A_2^m
  \end{align*}
  we write $\bar a^1
  \mapsto \bar a^2$ for the map that maps $a^1_i$ to $a^2_i$ for all
  $1\leq i \leq m$. 

  In the  $n$-round Ehrenfeucht-\Fraisse game on
  $\mathfrak{A}_1, a^1_1, a^1_2, \dots, a^1_m$ and $\mathfrak{A}_2,
  a^2_1, a^2_2, \dots, a^2_m$ for $a^j_i\in A_j$ there are two
  players, Spoiler 
  and Duplicator, which play according to the following rules. 
  The game is played for $n$ rounds. The $i$-th round consists of the
  following steps 
  \begin{enumerate}
  \item Spoiler chooses one of the structures, i.e., he chooses $j\in\{1,2\}$.
  \item Then he chooses one of the elements of his structure, i.e., he
    chooses some $a^j_{m+i}\in A_j$.
  \item Now,  Duplicator chooses an element in the other structure,
    i.e., for $k:=3-j$, she chooses some $a^k_{m+i} \in A_k$. 
  \end{enumerate}
  Having executed $n$ rounds, Spoiler and Duplicator defined
  tuples 
  \begin{align*}
  &\bar a^1:= a^1_1, a^1_2, \dots, a^1_{m+n} \in A_1^{m+n}\text{ and}\\
  &\bar a^2:= a^2_1, a^2_2, \dots, a^2_{m+n} \in A_2^{m+n}.    
  \end{align*}
  Duplicator wins 
  the play if $f:\bar a^1 \mapsto \bar a^2$ is a partial isomorphism,
  i.e., if $f$ satisfies the following conditions.
  \begin{enumerate}
  \item $a^1_i = a^1_j$ if and only if  $a^2_i = a^2_j$ for all $1\leq
    i \leq j \leq m+n$ and
  \item for each $R_i\in \sigma$ of arity $r$ the following hold: if 
    $i_1, i_2, \dots i_r$ are numbers between $1$ and $m+n$, 
    $\mathfrak{A_1}, \bar a^1 \models R_i x_{i_1} x_{i_2} \dots x_{i_r}$
    if and only if
    $\mathfrak{A_2}, \bar a^2 \models R_i x_{i_1} x_{i_2} \dots x_{i_r}$.
  \end{enumerate}
\end{definition}

\begin{definition}
  Let $\mathfrak{A}_1$, $\mathfrak{A}_2$ be structures and 
  $\bar a^1 \in \mathfrak{A}^n$, $\bar a^2\in \mathfrak{A}_2$.
  We write 
  $\mathfrak{A}_1, \bar a^1 \simeq_\rho \mathfrak{A}_2, \bar  a^2$ if
  Duplicator has a winning strategy in the $\rho$-round 
  Ehrenfeucht-\Fraisse game on 
  $\mathfrak{A}_1, \bar a^1$ and $\mathfrak{A}_2, \bar a^2$.
\end{definition}

Our interest in Ehrenfeucht-\Fraisse games stems from the following
relationship of $\simeq_\rho$ and $\equiv_\rho$ (recall that
$\equiv_\rho$ is equivalence with respect to \FO{} formulas up to
quantifier rank $\rho$).

\begin{lemma}[\cite{Fraisse54},\cite{Ehren60}] 
  For  all $\sigma$-structures $\mathfrak{A}_1$, $\mathfrak{A}_2$, and 
  for all tuples $\bar a^1\in A_1^n$, and $\bar a^2\in A_2^n$, 
  \begin{align*}
      \mathfrak{A}_1, \bar a^1 \simeq_\rho \mathfrak{A}_2, \bar a^2
      \text{ iff } 
      \mathfrak{A}_1, \bar a^1 \equiv_\rho \mathfrak{A}_2, \bar a^2, 
  \end{align*}
  i.e., Duplicator has a winning strategy in the $\rho$ round
  Ehrenfeucht-\Fraisse game on 
  $\mathfrak{A}_1, \bar a^1$ and $\mathfrak{A}_2, \bar a^2$ if and only if 
  $\mathfrak{A}_1, \bar a^1$ and $\mathfrak{A}_2, \bar a^2$ are
  indistinguishable by  first-order formulas of 
  quantifier rank $\rho$. 
\end{lemma}
\begin{remark} 
  We want to give some brief comments on the proof. 
  
  If $\mathfrak{A}_1, \bar a^1 \not\equiv_\rho \mathfrak{A}_2, \bar a^2$, then
  there is a formula 
  $\varphi$ in negation normal form (i.e., negation only occurs in
  negated atomic formulas) such 
  that $\mathfrak{A}_1, \bar a^1 \models\varphi$ but
  $\mathfrak{A}_2, \bar a^2 \not\models\varphi$. By induction on the
  structure of 
  $\varphi$ one can prove that there is a winning strategy for
  Spoiler. Basically, for every subformula starting with an
  existential quantification, Spoiler chooses a witness for this
  quantification in 
  $\mathfrak{A}$ and for each universal quantification, he chooses an
  element in $\mathfrak{B}$ witnessing the negation of the
  subformula. 
  Due to the fact that the second structure does not satisfy
  $\varphi$, Duplicator must eventually respond  with an element
  not satisfying the existential claim made by Spoiler. 
  By clever choice of further elements, Spoiler can then point out
  this difference and Duplicator will loose the game. 

  On the other hand, if the two structures cannot be distinguished by
  quantifier rank 
  $\rho$ formulas, then Duplicator just has to preserve the
  equivalence of the quantifier rank $m$ types of 
  the elements chosen in both structures, where $m$ is the number of
  rounds left to play. 
  Note that the resulting
  partial map is a 
  partial isomorphism if and only if it preserves all quantifier-free
  formulas. Thus, Duplicator wins the game using the strategy
  indicated above. 
\end{remark}

Ehrenfeucht-\Fraisse games are usually used to show
that first-order logic cannot express certain  properties.  
We stress that our main application of these games is
nonstandard. Nevertheless, we first
present an example of this classical application. In Section
\ref{Chapter HONPT} we  use the result of this example.   

\begin{example} \label{Ex:colouredWordStructures}
  We present a  proof that there are only finitely many types of
  coloured finite successor structures that are distinguishable by
  $\FO{\rho}$. 
  This example will also illustrate how the concept of Gaifman
  locality can be fruitfully applied to the analysis of
  Ehrenfeucht-\Fraisse games.\footnote{Our example is in fact
    an application of Hanf's Lemma (cf. \cite{Hanf65}).} In general, the
  analysis of 
  Ehrenfeucht-\Fraisse games is difficult because one has to consider
  too many possible choices for Spoiler. But if the structure of the
  local neighbourhoods is simple, this can be used to analyse
  Duplicator's strategies in the game.
  
  A finite successor structure is up to isomorphism a
  structure of the form 
  \begin{align*}
    \mathfrak{A}:=(\{ 1, 2, \dots n\},
    \mathrm{succ}, P_1, P_2, \dots, P_m)    
  \end{align*}
  for some $n,m\in\N$ where
  $\mathrm{succ}=\{(k,k+1):1\leq k <n\}$ is the successor relation on
  the natural numbers up to $n$ and $P_1, \dots, P_m$ are unary
  predicates (which we call colours). 
  We are going to show that for fixed $m\in\N$ there are at most
  $(\rho+2^{\rho})^{(2^m+1)^{2^{\rho+1}+1}}$ successor structures with $m$
  colours that are pairwise not $\simeq_\rho$-equivalent. 

  In order to prove this claim, we consider a successor structure
  $\mathfrak{A}$ with $m$ colours and $n$ elements. 
  We will make use of the $2^l$-local neighbourhood
  $\mathcal{N}_{2^l}(a)$ of the elements $a\in\mathfrak{A}$.
  
  Note that $\mathcal{N}_{2^l}(a)$ is a successor structure with
  exactly $2^{l+1}+1$ 
  many elements unless $a\leq 2^l$ or $a\geq n-2^l$. 
  Since there are at most $2^m$ many possibilities to colour a node with
  $m$ colours, there are at most $(2^m+1)^{2^{l+1}+1}$ many distinct
  $2^l$-local neighbourhoods up to isomorphism. The base $2^m+1$ is due to
  the fact that elements may be undefined (if $a\leq 2^l$ or $a\geq
  n-2^l$) or coloured in one of the $2^m$ possibilities. 
  
  We claim that the number of occurrences of each $2^\rho$-local
  neighbourhood type counted up 
  to threshold $2^\rho+\rho$ determines the $\simeq_\rho$-type of 
  a successor structure. 
  
  In order to prove this, we use Ehrenfeucht-\Fraisse games. Before we
  explain Duplicator's strategy, note the
  following facts.
  \begin{enumerate}
  \item Counting the occurrences of each $2^l$-local neighbourhood
    type up to some threshold $t\in\N$ determines the occurrences 
    of \mbox{$2^{l-1}$-neighbourhood} types up to threshold $t$. 
    Furthermore, the $2^l$-local neighbourhood of an element
    $a\in\mathfrak{A}$ determines the $2^{l-1}$-local neighbourhood of
    the elements  $a-2^{l-1}, a-2^{l-1}+1, \dots, a+2^{l-1}-1,
    a+2^{l-1}$.    
  \item For any $k < l$, the union of the $2^{l-k-1}$-local
    neighbourhoods of $k$ elements 
    contains at most 
    $k(2^{l-k}+1) = k + 2^{\ln(k) + l -k} \leq  k+2^{l} < l+2^{l}$ many
    elements. 
  \item The $2^l$-local neighbourhood types of the $2^l$ first and the
    $2^l$ 
    last elements of a successor structure $\mathfrak{A}$ occur
    exactly once in $\mathfrak{A}$ because they are determined by the
    number of elements that exist to their left, respectively, right. 
  \end{enumerate}

  Let $\mathfrak{A}$ and $\mathfrak{B}$ be structures that have, 
  up to threshold $2^\rho+\rho$, the same number of occurrences of
  each $2^\rho$-local neighbourhood type.
  
  Duplicator has the following strategy in the $\rho$ round
  Ehrenfeucht-\Fraisse game. Without loss of generality,
  Spoiler chooses at first some element $a_1\in\mathfrak{A}$. 
  Duplicator may respond with any element $b_1\in\mathfrak{B}$ such
  that $\mathcal{N}_{2^{\rho-1}}(a)$  and
  $\mathcal{N}_{2^{\rho-1}}(b)$ are isomorphic. 
  
  For the following $\rho-1$ rounds, we distinguish between local and
  global moves of Spoiler. Assume that in the $i$-th round the game is
  in position $(a_1, a_2, \dots, a_{i-1}) \mapsto (b_1, b_2,
  \dots, b_{i-1})$ such that the following holds.
  \begin{enumerate}
  \item For each $1\leq j \leq i-1$,  $\mathcal{N}_{2^{\rho-(i-1)}}(a_j)$ and
    $\mathcal{N}_{2^{\rho-(i-1)}}(b_2)$ are isomorphic.
  \item Up to threshold $1+2^{\rho-(i-1)}$, the distance of $a_j$ from
    $a_k$ agrees with the distance of $b_j$ from $b_k$, i.e., 
    $a_j$ is the $n$-th successor of $a_k$ for some $n \leq 2^{\rho-(i-1)}$
    if and only if $b_j$ is the $n$-th successor of $b_k$. 
  \end{enumerate}
  Due to symmetry we may assume  that Spoiler chooses some
  $a_i\in\mathfrak{A}$. We call
  this move local, if there is some $j<i$ such that 
  the distance between $a_i$ and $a_j$ is at most $2^{\rho-i}$. 
  In this case, Duplicator chooses the element $b_i$ that has the same
  distance to $b_j$ as $a_i$ to $a_j$. 
  Since the $2^{\rho-(i-1)}$-local neighbourhood of $a_j$ and $b_j$ coincide,
  the $2^{\rho-i}$-local neighbourhood of $a_i$ and $b_i$
  agree. Furthermore, note that the distances of $a_i$ from each $a_k$
  and the distances of $b_i$ from the corresponding $b_k$ agree up to
  threshold $2^{\rho-i}$. 

  If Spoiler chooses some $a_i\in\mathfrak{A}$ such that the distance
  between $a_i$ and $a_j$ for all $j<i$ is more than $2^{\rho-i}$, we
  call the move global. 
  In this case,  Duplicator chooses an element $b_i$ such that 
  $\mathcal{N}_{2^{\rho-i}}(a_i)
  \simeq \mathcal{N}_{2^{\rho-i}}(b_i)$ and such that the distance
  from $b_i$ to any 
  $b_j$ is more than $2^{\rho-i}$ for all $j<i$. Such an element $b_i$
  exists due to the following facts.
  \begin{enumerate}
  \item The $2^{\rho-i}$-local neighbourhoods of $a_1, a_2, \dots, a_{i-1}$
    contain less than $\rho+2^\rho$ many elements.
  \item For each $j<i$, $\mathcal{N}_{2^{\rho-(i-1)}}(a_j)\simeq
    \mathcal{N}_{2^{\rho-(i-1)}}(b_j)$. Thus, the 
    $2^{\rho-i}$-local neighbourhoods of the elements of distance at most
    $2^{\rho-i}$ from one of the $a_j$ are isomorphic to the
    corresponding elements that are close to $b_j$. Hence, for each
    isomorphism-type of a $2^{\rho-i}$-local neighbourhood the number of
    elements that realise this type and that are close to one of the $a_j$
    coincide with the number of elements that realise this type and that are
    close to one of the $b_j$. Let $k$ be the number of elements $a$
    close  to one of the $a_j$ such that $\mathcal{N}_{2^{\rho-i}}(a)\simeq
    \mathcal{N}_{2^{\rho-i}}(a_i)$. 
  \item Since $a_i$ is far away from all $a_j$, there are at least
    $k+1 \leq \rho+2^\rho$ many elements of neighbourhood type
    $\mathcal{N}_{2^{\rho-i}}(a)$ in $\mathfrak{A}$. Due to our assumptions on
    $\mathfrak{A}$ and $\mathfrak{B}$ and on the neighbourhoods of the
    elements chosen so far, there are at least $k+1$ 
    elements of the  neighbourhood type $\mathcal{N}_{2^{\rho-i}}(a)$ in 
    $\mathfrak{B}$ of which exactly $k$ have distance at most
    $2^{\rho-i}$ of one of the $b_j$. Thus, there is an element
    $b_i\in\mathfrak{B}$ such that $\mathcal{N}_{2^{\rho-i}}(b_i)\simeq
    \mathcal{N}_{2^{\rho-i}}(a_i)$ that is far away from all the $b_j$
    for $j<i$.  
  \end{enumerate}
  It is straightforward to see that, after $\rho$ rounds, we end up
  with a partial map 
  \begin{align*}
    f:(a_1, a_2, \dots, a_\rho)\mapsto (b_1, b_2, \dots,
    b_\rho) \text{ such that}    
  \end{align*}
  \begin{enumerate}
  \item for each $1\leq j \leq \rho$,  
    $\mathcal{N}_{1}(a_j)\simeq \mathcal{N}_{1}(b_j)$ whence
    the colours of $a_j$ and
    the colours of $b_j$ are equal, and
  \item up to threshold $2$, the distance of $a_j$ from
    $a_k$ agrees with the distance of $b_j$ from $b_k$, i.e., 
    $f$ and $f^{-1}$ preserve the successor relation.
  \end{enumerate}
  Thus, $f$ is a partial isomorphism and Duplicator wins the game. 

  Note that counting $2^\rho$-local neighbourhoods up to threshold 
  $\rho+2^\rho$ assigns to each successor structure with $m$ colours a
  function ${(2^m+1)^{2^{\rho+1}+1}} \to (\rho+2^\rho)$. We have seen that
  if these
  functions agree for two structures $\mathfrak{A}$ and
  $\mathfrak{B}$, then 
  Duplicator wins the $\rho$ round Ehrenfeucht-\Fraisse game on these two
  structures whence $\mathfrak{A}\equiv_{\rho}\mathfrak{B}$. 
  Thus,
  there are at most $(\rho+2^\rho)^{(2^m+1)^{2^{\rho+1}+1}}$ many
  $m$-coloured successor structures that can be distinguished by
  quantifier rank $\rho$ first-order formulas.
\end{example}

We will use the result of the previous example in Section
\ref{subsec:EquivalenceonWords}. 
But beside this classical application of
Ehrenfeucht-\Fraisse games, a nonstandard application of Ehrenfeucht-\Fraisse
games plays a much more important
role in this thesis. This application gives rise to \FO{} model checking
algorithms on certain classes of structures. Ferrante and
Rackoff\cite{Fer79} were 
the first to mention the general approach 
of using Ehrenfeucht-\Fraisse analysis for the 
 the decidability of \FO{} theories.

We  
consider the game played on two 
copies of the same structure, i.e., the game on $\mathfrak{A}, \bar a^1$ and
$\mathfrak{A}, \bar a^2$ with identical choice of the initial
parameter 
$\bar a^1 = \bar a^2 \in \mathfrak{A}$. At a first glance, this
looks quite uninteresting because Duplicator has of course a winning
strategy in this setting: he can copy each move of Spoiler. 
But we want to look for winning strategies
with certain constraints. In our application the constraint will be that
Duplicator is only allowed to choose elements 
that are represented by short runs of certain automata, but the idea
can be formulated more generally. 
\begin{definition}
  Let $\mathcal{C}$ be a class of structures. Assume that
  $S^{\mathfrak{A}}(m) \subseteq \mathfrak{A}^m$ is a subset of the
  $m$-tuples of the structure $\mathfrak{A}$ for each
  $\mathfrak{A}\in\mathcal{C}$ and each $m\in\N$. Set
  $S:=(S^{\mathfrak{A}}(m))_{m\in\N, \mathfrak{A}\in\mathcal{C}}$.
  We
  call $S$ a \emph{constraint for Duplicator's strategy} and we say 
  Duplicator has an  
  \emph{$S$-preserving} winning strategy if she has a strategy for
  each game 
  played on two copies of $\mathfrak{A}$ for some $\mathfrak{A}\in
  \mathcal{C}$  with the
  following property. 
  Let $\bar a^1 \mapsto \bar a^2$ be a position 
  reached after $m$ rounds where Duplicator used her strategy.
  If $\bar a^2\in S(m)$, then Duplicator's
  strategy chooses 
  an element $a^2_{m+1}$ such that 
  $\bar a^2,a^2_{m+1}\in S^{\mathfrak{A}}(m+1)$  for each  challenge
  of Spoiler in the 
  first copy of $\mathfrak{A}$.
\end{definition}
\begin{remark}
  We write $S(m)$ for $S^{\mathfrak{A}}(m)$ if
  $\mathfrak{A}$ is clear from the context.
\end{remark}
Recall the following fact:
if Duplicator uses a winning strategy in the $n$ round game,
 her choice in the $(m+1)$-st round is an
element $a^2_{m+1}$ such that  
$\mathfrak{A},\bar a^1, a^1_{m+1} \equiv_{n-m-1} \mathfrak{A}, \bar a^2,
a^2_{m+1}$. 

This implies that if Duplicator has an $S$-preserving winning
strategy, then
for every formula
$\varphi(x_1, x_2, \dots, x_{m+1})\in \FO{n-m-1}$ and for all $\bar
a\in A^m$ with $\bar a\in S(m)$ the following holds:
\begin{align*}
  &\text{there is an element }a\in A\text{ such that }\bar a, a\in
  S(m+1)\text{ and } \mathfrak{A}, \bar a, a \models \varphi \\ 
  \text{iff }&\text{there is an element }a\in A\text{ such that }\mathfrak{A},
  \bar a, a \models  \varphi \\
  \text{iff }
  &\mathfrak{A}, \bar a \models \exists x_{m+1} \varphi.
\end{align*}
Replacing existential quantification with universal quantification we
obtain directly that this statement is equivalent to
\begin{align*}
  &\text{for all }a\in A\text{ such that } \bar a, a\in S(m+1)\text{
    we have } \mathfrak{A}, \bar a, a \models \varphi \\
  \text{iff }& \mathfrak{A}, \bar a \models \forall x \varphi(\bar y, x). 
\end{align*}

\begin{algorithm2e}[]
  \label{AlgoGeneralModelCheck}
  \caption{The general \FO{}-model checking as pseudo-code}
  \SetKwFunction{ModelCheck}{ModelCheck}
  \SetKw{accept}{accept}
  \SetKw{reject}{reject}
  \Titleofalgo{ \ModelCheck{$\mathfrak{A}, \bar a, \varphi(\bar x)$}}
  \KwIn{a structure $\mathfrak{A}$
    ,a formula $\varphi\in\FO{\rho}$, an
    assignment $\bar x \mapsto \bar a$}
  \uIf{ $\varphi$ is an atom or negated atom} {
    \lIf{ $\mathfrak{A}, \bar a \models \varphi(\bar x)$} 
    {\accept }\lElse{\reject\;}}
  \uIf{$\varphi = \varphi_1 \vee \varphi_2$}{
    \lIf{\ModelCheck{$\mathfrak{A}, \bar a, \varphi_1$} $=$ accept}{\accept}
    \uElse{ \lIf {\ModelCheck{$\mathfrak{A}, \bar a, \varphi_2$}$=$
        accept} {\accept}
      \lElse{\reject\;}
    }
  }
  \uIf{$\varphi = \varphi_1 \wedge \varphi_2$}{
    \lIf{
      \ModelCheck{$\mathfrak{A}, \bar a, \varphi_1$}$=$
      \ModelCheck{$\mathfrak{A}, \bar a, \varphi_2$}$=$ accept}{\accept}
    \lElse{\reject\;}}
  \uIf{$\varphi=\exists x \varphi_1(\bar x,x)$}{
    check whether there is an $a\in\mathfrak{A}$ such that
    \ModelCheck{$\mathfrak{A}, \bar a a, \varphi_1$}$=$ accept\;}
  \uIf{$\varphi=\forall x_i \varphi_1$}{
    check whether 
    \ModelCheck{$\mathfrak{A},\bar a a, \varphi_1$}$=$ accept
    holds for all $a\in\mathfrak{A}$\;}
\end{algorithm2e}

Now, we want to make use of this observation in a general approach to
first-order model checking. The pseudo-algorithm in Algorithm
\ref{AlgoGeneralModelCheck} is a correct description of 
first-order model checking as it just proceeds by syntactic induction
on the first-order formula in order to determine whether the given structure
is a model of the given formula. But of course, in general this is no
algorithm.
As soon as $\mathfrak{A}$ is infinite and a quantification occurs in
$\varphi$, this pseudo-algorithm would not terminate because it would
have to check infinitely many variable assignments. Nevertheless, it
is correct in the sense that if we consider a class of structures
where we could check these infinitely many variable assignments in
finite time, then it would
correctly determine the answer to the model checking problem. 
Using $S$-preserving strategies, we want to turn the pseudo-algorithm
into a proper algorithm for certain classes of structures. 
For this purpose, we first introduce the following notation. 
\begin{definition}
  Given a  class $\mathcal{C}$ of finitely represented structures, 
  we call a constraint $S$ for Duplicator's strategy 
  \emph{finitary} on $\mathcal{C}$, if
  for each $\mathfrak{A}\in\mathcal{C}$ we can compute a function
  $f_{\mathfrak{A}}$ such that for all $n\in \N$ 
  \begin{itemize}
  \item $S^{\mathfrak{A}}(n)$ is finite,
  \item for each $\bar a\in S^{\mathfrak{A}}(n)$, we can represent $\bar
    a$ in space $f_{\mathfrak{A}}(n)$, and
  \item $\bar a \in S^\mathfrak{A}(n)$ is effectively decidable. 
  \end{itemize}
\end{definition}
Using such a finitary constraint, we can rewrite the model checking
algorithm from above into Algorithm 
\ref{AlgoSPReservingModelCheck}. 
The condition of a finitary constraint is exactly what is needed to
guarantee termination of this algorithm. 
Furthermore, our observation on $S$-preserving constraints 
implies that this algorithm is correct for all structures
from a class $\mathcal{C}$ where Duplicator has an $S$-preserving
winning strategy for every $\mathfrak{A}\in\mathcal{C}$. We apply this
idea in Sections 
\ref{Chapter_FO-NPT} and \ref{Chapter HONPT}.  
There, we represent elements of certain structures 
$\mathfrak{A}$ by runs of some automaton. The sets
$S^\mathfrak{A}(n)$ consist of runs that have length bounded by
some function $f_\mathfrak{A}$ that is computable from the automaton
representing $\mathfrak{A}$. 

\begin{algorithm2e}[]
  \label{AlgoSPReservingModelCheck}
  \caption{\FO{}-model checking on $S$-preserving structures}
  \SetKwFunction{ModelCheck}{SModelCheck}
  \SetKw{accept}{accept}
  \SetKw{reject}{reject}
  \Titleofalgo{ \ModelCheck{$\mathfrak{A}, \bar a, \varphi(\bar x)$}}
  \KwIn{a structure $\mathfrak{A}$
    , a formula $\varphi\in\FO{\rho}$, an
    assignment $\bar x \mapsto \bar a$ for tuples $\bar x, \bar a$ of
    arity $m$ such that $\bar a\in S(m)$}
  \uIf{ $\varphi$ is an atom or negated atom} {
    \lIf{ $\mathfrak{A}, \bar a \models \varphi(\bar x)$} 
    {\accept }\lElse{\reject\;}}
  \uIf{$\varphi = \varphi_1 \vee \varphi_2$}{
    \lIf{\ModelCheck{$\mathfrak{A}, \bar a, \varphi_1$} $=$ accept}{\accept}
    \uElse{ \lIf {\ModelCheck{$\mathfrak{A}, \bar a, \varphi_2$}$=$
        accept} {\accept}
      \lElse{\reject\;}
    }
  }
  \uIf{$\varphi = \varphi_1 \wedge \varphi_2$}{
    \lIf{
      \ModelCheck{$\mathfrak{A}, \bar a, \varphi_1$}$=$
      \ModelCheck{$\mathfrak{A}, \bar a, \varphi_2$}$=$ accept}{\accept}
    \lElse{\reject\;}}
  \uIf{$\varphi=\exists x \varphi_1(\bar x,x)$}{
    check whether there is an $a\in\mathfrak{A}$ such that
    $\bar a, a\in S(m+1)$ and
    \ModelCheck{$\mathfrak{A}, \bar a a, \varphi_1$}$=$ accept\;}
  \uIf{$\varphi=\forall x_i \varphi_1$}{
    check whether 
    \ModelCheck{$\mathfrak{A},\bar a a, \varphi_1$}$=$ accept
    holds for all $a\in\mathfrak{A}$ such that $\bar a, a\in S(m+1)$\;}
\end{algorithm2e}

\subsection{Extensions of First-Order Logic}
In many cases  the expressive power of \FO{} is too weak. For example,
due to the local nature of first-order logic, simple reachability
questions cannot be formalised in \FO{}. In order 
to overcome this weakness, there have been proposed a lot of different
extensions tailored for different applications. 
In the following, we present those extensions that we use
later.

\paragraph{Monadic Second-Order Logic}
Perhaps the most classical extension of first-order logic is monadic
second-order logic (abbreviated \MSO). The formulas of this logic
are defined using the same rules as for \FO{}  but
additionally adding quantification over subsets. For this, we fix a
set $X_1, X_2, \dots$ of set variable symbols. We extend the
formation rules of first-order logic by the
following two rules.
\begin{itemize}
\item $X_i x$ is an \MSO formula for any variable symbol $x$ and any 
  set variable symbol $X_i$.
\item If $\varphi$ is an \MSO formula then $\exists X_i \varphi$ and 
  $\forall X_i \varphi$ are also  \MSO formulas.  
\end{itemize}
For the semantics, we extend the variable assignment $I$ to the set of set
variable symbols occurring in a formula $\varphi$. Now, $I$ maps each
symbol $X_i$ to a subset of the structure $\mathfrak{A}$. 
We then set 
\begin{itemize}
\item $\mathfrak{A}, I \models X_i x$ if $I(x)\in I(X_i)$,
\item $\mathfrak{A}, I \models  \exists X \varphi$ if there is some
  $M\subseteq A$ such that
  $\mathfrak{A},I_{X\mapsto M} \models \varphi$ where 
  $I_{X\mapsto M}$ is identical to $I$ but maps $X$ to $M$, and 
\item $\mathfrak{A}, I \models  \forall X \varphi$ if 
  $\mathfrak{A},I_{X\mapsto M} \models \varphi$
  for all $M\subseteq A$.
\end{itemize}

\MSO is the most expressive logic that we are going to
consider. But this expressive power comes at a prize. 
The \MSO model checking on collapsible
pushdown graphs and nested pushdown trees is undecidable. Thus, 
we look for weaker extensions of first-order logic that are still
decidable in our setting. 

\paragraph{Monadic Least Fixpoint Logic}

Another approach for extending first-order logic is the use of
fixpoint-operators. Here, we present the monadic least fixpoint
logic (MLFP). Consider an $\MSO$ formula $\varphi$ without
quantification over sets, with a free
variable $x$ and a free set variable $X$ that only
occurs positively, i.e., that only occurs under an even scope of
negations. For each structure $\mathfrak{A}$, each variable assignment
$I$, and $M\subseteq A$, we write $I^M$ for 
$I_{X\mapsto M}$, the variable assignment that is identical to $I$ but maps
the set variable $X$ to $M\subseteq A$.
Now, $\varphi$ defines a monotone operator
\begin{align*}
  &f^\varphi : \Pot{A} \to \Pot{A}\\
  &f^\varphi(M) := \{ a\in A: \mathfrak{A},I^M_{x\mapsto a}\models \varphi\}.
\end{align*}
Due to the theorem of Knaster and Tarski \cite{Knaster28}, $f^\varphi$
has a unique
least fixpoint $M\subseteq A$, i.e., 
there is a minimal
set $M^{\varphi}\subseteq A$ such that $f^{\varphi}(M^{\varphi})=M^{\varphi}$. 
MLFP is the extension of
\FO{} by the rule that $[\lfp_{x,X} \varphi](y)$ is an MLFP formula where
$\varphi$ is a formula as described above and where  $y$ is a free
variable. The semantics of 
$\psi=[\lfp_{x,X} \varphi](y)$ is defined by $\mathfrak{A}, I \models
\psi$ iff $I(y) \in M^{\varphi}$. 
 
It is clear that the expressive power of MLFP is between the
expressive power of \FO{} and that of \MSO. Each MLFP formula can be
translated into an equivalent \MSO formula because the least fixpoint of
$\varphi(x,X)$ is defined by the formula 
\begin{displaymath}
\psi(Z):= \left(\forall Y \forall y\ (Yy \leftrightarrow  \varphi(y,Y))\right)
\rightarrow Z\subseteq Y
\end{displaymath}
were $X\subseteq Y$ is an abbreviation for $\forall x (Xx \rightarrow
Yx)$. Thus, $[\lfp_{x,X} \varphi](y)$ can be translated into 
$\exists Z (\psi(Z) \wedge Zy)$. 
The expressive power of MLFP is strictly greater than that of \FO{} because
fixpoints can be used to formalise reachability queries. For example,
the fixpoint induced by the formula 
$\varphi(x,X):= Px \vee \exists y (Exy \wedge Xy)$ contains all
elements for which an $E$-path to an element in $P$ exists. Due to the
local nature of first-order logic, this is not expressible with an
\FO{} formula. 

The least fixpoint operator is a very strong extension of \FO{}
in the sense that
MLFP is much more expressive  than
\FO{}. As in the case of \MSO, 
MLFP is too powerful on those structures we are interested in. 
We will show that the MLFP-theory of
a certain collapsible pushdown graph of level $2$ is undecidable. 

Thus, in order to find logics with a decidable
model checking problems on collapsible pushdown graphs, 
we look at logics with strictly weaker expressive power than that of
MLFP. 

\paragraph{FO + Reachability Predicates}

During the last decade another extension of first-order logic
has been studied and successfully applied for  model checking. 
If one looks at verification problems, the most important properties
that one wants to verify often involve reachability of certain
states. Thus, for classes of graphs where $\MSO$ and MLFP are
undecidable, one may study the weakest extension of first-order
logic that allows to express reachability questions. We call this
logic $\FO{}(\Reach)$ and we introduce it formally in the following 
definition. 

\begin{definition} \label{DefReachPredicate}
  Let $\sigma=(E_1, E_2, \dots, E_n)$ be a relational signature and
  $E_i$ a binary relation symbol for each $1\leq i \leq n$. 
  Let $\FO{}(\Reach)$ denote the smallest set generated by the
  formation rules of first-order logic plus the rule that
  $\Reach x y$ is an $\FO{}(\Reach)$ formula for each pair of
  variables $x,y$. 
  
  $\FO{}(\Reach)$ inherits its semantics mainly from $\FO{}$. 
  If we consider $\mathfrak{A}$ as a graph with edge relation 
  $E:=\bigcup_{1\leq j \leq n} E_j^\mathfrak{A}$ where each edge is
  labelled with a nonempty subset of $\{1, 2, \dots, n\}$ then
  $\Reach$ is interpreted as the transitive closure of the edge
  relation $E$.  
  
  Similar to \FO{} extended by reachability we now introduce \FO{}
  extended by \emph{ regular reachability}. 
  Let $\sigma, \mathfrak{A}$ and $E$ be defined as before. 
  For simplicity, we assume that each edge in $\mathfrak{A}$ is
  labelled by exactly one label from $\{1, 2, \dots, n\}$, i.e., for
  all $a,a'\in A$, $(a,a')\in E_i$ and $(a,a')\in E_j$ implies $i=j$. 

  We write \FO{}(Reg) for the extension of $\FO{}$ by atomic formulas
  $\Reach_L xy$ for each regular language $L\subseteq \{1, 2, \dots,
  n\}^*$ and all variable symbols $x,y$. 

  For $a,b\in\mathfrak{A}$,  $\Reach_L ab$ holds if there is a
  path $a=a_1, a_2, a_3, \dots, a_k=b$ such that 
  $(a_i, a_{i+1})\in E$ for all $1\leq i <k$ and the word formed by
  the labels of the edges along this path form a word of $L$. 
\end{definition}

\paragraph{FO + Generalised Quantifiers}

Lastly, we present the  extensions of first-order logic 
by generalised quantifiers. The idea of generalised quantifiers was
first introduced by 
Mostowski \cite{Mostowski57}
and then further developed to full generality by
Lindstr\"om \cite{Lindstroem66}. We 
briefly recall the general notion. Afterwards,  we present the 
generalised quantifiers that occur in this thesis:
the
infinite existential quantifier $\exists^\infty$; the modulo counting
quantifiers $\exists^{(k,n)}$; and the Ramsey-  or Magidor-Malitz quantifier
$\RamQ{n}$, first introduced by Magidor and Malitz \cite{MagidorM77}. 

\begin{definition}[\cite{Lindstroem66}]
  Let $\sigma$ be some vocabulary. A collection of $\sigma$-structures
  $Q$ which is closed under isomorphisms is called a \emph{generalised
    quantifier}.
\end{definition}
Let $Q$ be some generalised quantifier. 
$\FO{}(Q)$ denotes the extension 
of first-order logic by this quantifier. 
The formulas of $\FO{}(Q)$
are defined using the formation rules of $\FO{}$ and the following rule.
If $\varphi_i(x^i_1, x^i_2, \dots x^i_{a_i}, \bar y)$ is a formula in
$\FO{}(Q)$ for each $1\leq i\leq n$, then 
\begin{align*}
  Q{x^1_1, \dots x^1_{a_1}, x^2_1, \dots, x^2_{a_2}, \dots,x^n_1,
    \dots,  x^n_{a_n}}
  (\varphi_1)(\varphi_2)\dots(\varphi_n)  
\end{align*}
 is a 
formula of $\FO{}(Q)$. 

The semantics of this formula is defined as
follows. Let $\mathfrak{A}$ be some structure and $I$ some
variable assignment. We set 
\begin{align*}
  &R_i^I:=\{\bar a \in A^{a_i}: \mathfrak{A},
  I_{\bar x^i \mapsto \bar a} \models \varphi_i\}\text{ where}\\
  &\bar x^i = x^i_1, x^i_2, \dots, x^i_{a_i}.    
\end{align*}
Now, we define the semantics of the quantifier by
\begin{align*}
  \mathfrak{A}, I \models Q{x^1_1, \dots x^1_{a_1}, x^2_1, \dots,
    x^2_{a_2}, \dots, x^n_1, \dots, x^n_{a_n}}
  (\varphi_1)(\varphi_2)\dots(\varphi_n)\text{ if } 
  (A, R_1^I, R_2^I, \dots, R_n^I)\in Q.  
\end{align*}
 
\begin{example}
  \begin{enumerate}
  \item[]
  \item Let $\exists$ consist of all structures $\mathfrak{A}=(A,P)$
    for $P$ a unary predicate $\emptyset\neq P \subseteq A$. The
    generalised quantifier defined by $\exists$ is the usual
    existential quantifier.
  \item Analogously, let $\forall$ consist of all structures
    $\mathfrak{A}=(A,P)$ for $P$ the unary predicate $P=A$. The
    generalised quantifier defined by $\forall$ is the usual universal
    quantifier.
  \item The infinite existential quantifier $\exists^\infty$ is
    defined by the collection of structures $\mathfrak{A}=(A,P)$ where
    $P\subseteq A$ is infinite. Some structure
    $\mathfrak{B}$ with some variable assignment $I$ satisfies
    $\mathfrak{B}, I \models \exists^\infty x \varphi$ if there are
    infinitely many pairwise distinct elements $b_1, b_2, b_3 \dots
    \in B$ such that $\mathfrak{B}, I_{x\to b_i} \models \varphi$ for
    all $i\in\N$.
  \item The modulo counting quantifier $\exists^{(k,m)}$ is defined by
    the collection of structures \mbox{$\mathfrak{A}=(A,P)$} where
    $P\subseteq A$ and $\lvert P \rvert = k \mod m$. 
    Some structure
    $\mathfrak{B}$ with some variable assignment $I$ satisfies
    $\mathfrak{B}, I \models \exists^{(k,m)} x \varphi$ if there are
    $k \mod m$ many pairwise distinct $b_i \in B$ such that
    $\mathfrak{B}, I_{x\to b_i} \models \varphi$.
  \item Finally, we introduce the Ramsey quantifier of arity $i$. Let
    us say that a set $S\subseteq A^n$ contains an infinite
    box if there is an infinite subset $A'\subseteq A$ such
    that $S$ contains all $n$-tuples of pairwise distinct elements
    from $A'$. Let $\RamQ{n}$ contain all structures
    $\mathfrak{A}=(A, P)$ where $P\subseteq A^n$ contains an infinite
    box. This means that $\RamQ{1} = \exists^\infty$ and 
    $\RamQ{2}xy(Exy)$ is the formula that states that there is an
    infinite clique with respect to the binary relation $E$. 
  \end{enumerate}
\end{example}
For a more detailed introduction to generalised quantifiers we
refer the reader to the survey of V\"a\"an\"anen
 \cite{Vaananen1997}. 

\begin{definition}
  We denote by $\FO{}(\exists^{\mathrm{mod}})$ the extension of $\FO{}$
  by modulo counting quantifiers. By $\FO{}((\RamQ{n})_{n\in\N})$ we denote
  the extension of $\FO{}$ by Ramsey quantifiers. Analogously, 
  $\FO{}(\exists^{\mathrm{mod}}, (\RamQ{n})_{n\in\N})$ denotes the
  extension of $\FO{}$ by both types of quantifiers. 
\end{definition}
\begin{remark}
  Note that $\exists^\infty$ is expressible in
  $\FO{}(\exists^{\mathrm{mod}})$: $\exists^\infty x(\varphi(x))$ is
  equivalent to $\neg(\exists^{0,2} x(\varphi(x)) \lor \exists^{1,2}
  x(\varphi(x)))$. Thus, we will use the quantifier $\exists^\infty$ as an
  abbreviation in $\FO{}(\exists^{\mathrm{mod}})$. 
\end{remark}

We come back to the generalised quantifiers $\exists^{\mathrm{mod}}$
and $\RamQ{n}$ in
Section \ref{SectionRamseyQuantifier}. We will show that first-order
logic extended by these quantifiers is decidable on tree-automatic
structures.

\subsection{Basic Modal Logic and  $L\mu$}

Beside the classical logics like \FO{}, \MSO, and their extensions
there is a another class of logics of great importance in the field of
model checking: Basic modal logic and its extensions. 

Almost
a century ago, C. I. Lewis \cite{Lewis18} introduced  a modal operator
for the first time.  
Since then, modal operators and modal logics have been studied
intensively and found applications in very different fields like 
philosophy, mathematics, linguistics, computer science, and economic
game theory. For an introduction to modal logics we refer the
reader to the introductory chapters of \cite{HandbookML06}.  

This thesis is mainly concerned with model checking for classical
logics. Nevertheless, we use  some results concerning model checking
for basic modal logic and modal $\mu$-calculus. Thus, 
we will briefly
recall the basic definitions and introduce our notation. 

We fix a signature $\sigma=(E_1, E_2, \dots, E_n, P_1, P_2, \dots,
P_m)$ of binary relations $E_i$ and unary relations $P_j$ called
propositions. 

\begin{definition}
  \emph{Modal logic} over the signature $\sigma$ consists of the formulas
  generated by iterated use of the following rules.
  \begin{enumerate}
  \item $\True$ and $\False$ are modal formulas. 
  \item $p_j$ is a modal formula for $1\leq j \leq m$.
  \item For $\varphi, \psi$ modal formulas, their conjunction,
    disjunction and negation are modal formulas, i.e., 
    $\varphi \land \psi, \varphi \lor \psi, \neg\varphi$ are modal
    formulas.
  \item If $\varphi$ is a modal formula, then $\langle E_i \rangle
    \varphi$ and $[E_i] \varphi$ are modal formulas. 
  \end{enumerate}
  In the modal terminology, one calls $\sigma$-structures \emph{Kripke
    structures}. A Kripke structure $\mathfrak{A}$ together with a
  distinguished element $a\in \mathfrak{A}$ is called a pointed Kripke
  structure. The semantics of modal formulas is inductively defined
  according to the following rules.
  \begin{enumerate}
  \item $\mathfrak{A}, a \models \True$ and     
    $\mathfrak{A}, a \not\models \False$ for all pointed Kripke structures 
    $\mathfrak{A},a$. 
  \item For $1\leq j \leq m$, $\mathfrak{A}, a \models p_j$ if $a\in P_j$.
  \item For formulas of the form $ \varphi \land \psi $, 
    $\varphi \lor \psi $, and $\neg \varphi$ we use the standard
    interpretation of the logical connectives.
  \item For $\varphi = \langle E_i \rangle \psi$, we set
    $\mathfrak{A}, a \models \varphi$ if there is some $a'\in A$ such
    that $(a,a')\in E_i$ and $\mathfrak{A}, a' \models \psi$. 
    For $\varphi = [E_i] \psi$, we set
    $\mathfrak{A}, a \models \varphi$ if $\mathfrak{A}, a' \models
    \psi$ for all $a'\in A$ such that $(a,a')\in E_i$. 
  \end{enumerate}
\end{definition}

The expressive power of modal logic is strictly contained in that of
first-order logic. This can be seen immediately when applying the
so-called standard translation. The basic idea
is that $\langle E_i \rangle \varphi$ is 
translated into a formula $\exists x ( E_i yx \land \hat \varphi(x))$
where $\hat\varphi(x)$ is the standard translation of $\varphi$.
$[E_i] \varphi$ is translated using the duality of $\langle E_i
\rangle$ and $[E_i]$, i.e., replacing $[E_i]\varphi$ by $\neg\langle
E_I \rangle \neg \varphi$. 
By clever reuse of variable names, it suffices to use $2$ variables in
this translation. The popularity of modal logic in model checking
stems from its algorithmic tractability. Each satisfiable modal
formula has a model which is a finite tree. Since trees are
algorithmically well-behaved, one can develop very efficient
algorithms for model checking of modal formulas. But this comes at the
cost that the 
expressive power of modal logic is quite low. Thus, there have been
many proposals how to extend the expressive power of modal logic while
keeping the good algorithmic behaviour. 
One of the most powerful extensions of modal logic is the modal
$\mu$-calculus. This is the extension of modal logic by
fixpoint operators analogously to the extension MLFP of first-order
logic. 
\begin{definition}
  In order to define the modal $\mu$-calculus, we fix set variables
  $X, Y, Z, \dots$ 
  The modal $\mu$-calculus (denoted as $L\mu$) over the signature $\sigma$ is 
  the set of formulas generated by the following rules.
  \begin{enumerate}
  \item We may use all the rules that are used to generate the
    formulas of  modal logic.
  \item Additionally, $X$ is a formula for each set variable $X$.
  \item Finally, if $X$ occurs only positively in an $L\mu$-formula
    $\varphi$, i.e., under the scope of an even number of negations,
    then $\mu X.\varphi(X)$ is a formula of $L\mu$. 
  \end{enumerate}
  Fix a $\sigma$-structure $\mathfrak{A}$, a variable assignment 
  $I: V \to A$ and a point $a\in A$. We say $\mathfrak{A}, I, a\models X$
  for $X\in V$ if $a\in I(X)$. 
  For $\varphi= \mu X.\psi(X)$, we say $\mathfrak{A}, I, a \models
  \varphi$ if $a\in M^\psi$ for $M^\psi\subseteq A$ the least fixpoint of the
  operator that maps any subset $B\subseteq A$ to
  $\{a\in A: \mathfrak{A}, I_{X\mapsto B}, a \models \psi(X)\}$. 
  The rules for all other formulas are inherited from the semantics of
  modal logic in the obvious way. 
\end{definition}
 $L\mu$ can be embedded into MLFP, i.e., for each $L\mu$-formula there is an
 equivalent MLFP formula. One extends the standard
 translation of modal logic to $\FO{}$ by the obvious translation of
 fixpoints in $L\mu$ to fixpoints in MLFP. 
 $L\mu$ is a very powerful modal logic. Its
 expressive power encompasses many modal logics
 like linear time logic (LTL), computation tree logic (CTL) or
 $\mathrm{CTL}^*$. 

\subsection{Logical Interpretations}

Logical interpretations are a formal framework to identify a structure
that ``lives'' in another structure. This concept is used widely in
mathematics. 
For instance,  if one investigates the multiplicative group of a
field, this is in fact
the interpretation of a group within a field. This is one of the
easiest examples of an interpretation. 
The following example from linear algebra illustrates a slightly more
involved application of the concept of interpretations.
\begin{example}
  Let $(V, +, \cdot)$ be some $n$-dimensional vectorspace. It is
  commonly  known that the endomorphisms of $V$ with concatenation
  $\circ$ and pointwise addition form a ring $\mathrm{End}(V)$. This
  ring is isomorphic 
  to the ring of $n\times n$-dimensional matrices with addition and
  multiplication. 

  In this representation, $\mathrm{End}(V)$ is interpretable in $V$:
  The domain of this interpretation are all $n^2$-tuples from
  $V$ where the $k$-th element of this tuple is considered as the
  entry in the $\left\lceil \frac{k}{n} \right\rceil$-th row and the
  $(k \mod 
  n)$-th column. Addition and composition of the endomorphisms can
  then be reduced to computations on these $n^2$-tuples. Addition of two
  morphisms corresponds to pointwise addition of the $n^2$-tuples and
  composition can be reduced using the known formulas for matrix
  multiplication. 

  In logical terms, this is an $n^2$-dimensional first-order
  interpretation of the ring $\mathrm{End}(V)$ in the vectorspace
  $V$. 
\end{example}

We call an interpretation \emph{logical} if it is defined using
formulas from some 
logic. The idea of using logical interpretations goes back to
Tarski who used this concept to obtain undecidability results. 
Since then, the use of interpretations for decidability or
undecidability proofs for the theories of certain structures has been
a fruitful approach. For a detailed survey on logical interpretations
we recommend the article of Blumensath et al.\ \cite{BlumensathCL07}.
We briefly introduce our notation concerning interpretations and
the important results that we are going to use. 

Given some structure $\mathfrak{A}$, we can use formulas of some
logic L to define 
a new structure $\mathfrak{B}$ from $\mathfrak{A}$.
The idea is to obtain the domain of $\mathfrak{B}$ as an L definable
subset of $A^n$. Then we define relations in this new structure via
formulas in the signature of the old structure. If we obtain some
structure $\mathfrak{B}$ in this way from another structure
$\mathfrak{A}$, we say that $\mathfrak{B}$ is interpretable in
$\mathfrak{A}$. If $\mathfrak{B}$ is interpretable in $\mathfrak{A}$
this can be used to reduce the model checking problem 
on input $\mathfrak{B}$ to the model checking problem on
$\mathfrak{A}$. In this thesis we will use $\FO{}$-interpretations and
one-dimensional $\MSO$-interpretations. Let us start with introducing
$\FO{}$-interpretations formally. 
\begin{definition}
  Let $\sigma:=(E_1, E_2, \dots, E_n)$ and 
  $\tau:=(F_1, F_2, \dots, F_m)$ be relational signatures. For  $n\in\N$,
  an ($n$-dimensional-$\sigma$-$\tau$) $\FO{}$-interpretation is given by a
  tuple of $\FO{}(\sigma)$-formulas  
  $I:=(\varphi, \psi_{F_1}, \psi_{F_2}, \dots, \psi_{F_n})$ where 
  $\varphi$ has $n$ free variables and each $\psi_{F_i}$ has $r_i\cdot
  n$ free variables where $r_i$ is the arity of $F_i$. 
  
  The interpretation $I$ induces two maps: one from
  $\sigma$-structures to $\tau$-structures and another from
  $\tau$-formulas to $\sigma$-formulas. 
  
  Let $\mathrm{Str}_I$ be the map that maps a $\sigma$-structure
  $\mathfrak{A}:=(A, E_1^\mathfrak{A}, \dots, E_n^\mathfrak{A})$ to
  the $\tau$-structure $\mathfrak{B}:=(B, F_1^\mathfrak{B}, \dots,
  F_m^\mathfrak{B})$ where
  \begin{align*}
    &B:=\left\{\bar a\in A^n: \mathfrak{A} \models \varphi(\bar a)\right\}
    \text{ and }\\
    &F_i^\mathfrak{B}:=\left\{(\bar a_1, \bar a_2,\dots, \bar
      a_{r_i})\in A^{n \cdot r_i}: \bar
      a_1, \bar a_2, \dots, \bar a_{r_i}\in B \text{ and }
      \mathfrak{A}\models \psi_{F_i}(\bar a_1, 
      \bar a_2, \dots, \bar a_{r_i})\right\}.  
  \end{align*}
  Let $\mathrm{Frm}_I$ be the map that maps an $\FO{}(\tau)$ formula
  $\alpha$ to the formula $\mathrm{Frm}_I(\alpha)$ which is obtained
  by the following rules:
  \begin{itemize}
  \item If $\alpha= F_i x^1 x^2 \dots x^k$ for  variable
    symbols $x^i$, then set 
    \begin{align*}
      \mathrm{Frm}_I(\alpha):= \psi_{F_i}(x^1_1, x^1_2,
      \dots, x^1_n, x^2_1, \dots, x^2_n, \dots, x^k_1, \dots x^k_n)      
    \end{align*}
    where $n$ is the  dimension of $I$.
  \item Boolean connectives are preserved, i.e., if $\alpha=\alpha_1
    \lor \alpha_2$ then
    \begin{align*}
      \mathrm{Frm}_I(\alpha)=
      \mathrm{Frm}_I(\alpha_1)\lor
      \mathrm{Frm}_I(\alpha_2)   
    \end{align*}
    and analogously for $\neg$ and $\land$.
  \item If $\alpha= \exists x \alpha_1$, then 
    $\mathrm{Frm}_I(\alpha) := \exists x_1 \exists x_2 \dots \exists x_n
    (\varphi(x_1, x_2, \dots, x_n) \land \mathrm{Frm}_I(\alpha_1))$. \\
    If $\alpha= \forall x \alpha_1$, then 
    $\mathrm{Frm}_I(\alpha) := \forall x_1 \forall x_2 \dots \forall x_n
    (\varphi(x_1, x_2, \dots, x_n) \rightarrow \mathrm{Frm}_I(\alpha_1))$. 
  \end{itemize}
\end{definition}
The well-known connection between $\mathrm{Str}_I$ and $\mathrm{Frm}_I$ is given
in the following lemma. 
\begin{lemma}
  Let $I$ be an $n$-dimensional-$\sigma$-$\tau$ $\FO{}$-interpretation,
  $\mathfrak{A}$ some $\sigma$-structure and $\varphi$ some
  $\FO{}(\tau)$ sentence. Then 
  \begin{displaymath}
    \mathrm{Str}_I(\mathfrak{A}) \models \varphi \text{ iff }
    \mathfrak{A} \models \mathrm{Frm}_I(\varphi).
  \end{displaymath}
\end{lemma}
The proof is by induction on the structure of $\varphi$. 

For $\FO{}$ model checking purposes, interpretations can be used as
follows. Fix an interpretation $I$ and two classes $\mathcal{C}_1$ and
$\mathcal{C}_2$ of structures. 
Assume that there is a 
computable function $\mathrm{Str}_I^{-1}$ that maps each
$\mathfrak{A}\in \mathcal{C}_1$ to a 
structure $\mathfrak{B}\in\mathcal{C}_2$ such that
$\mathrm{Str}_I(\mathfrak{B})=\mathfrak{A}$. 
Then we can reduce the model checking problem for $\mathcal{C}_1$
to  the model checking problem for  $\mathcal{C}_2$.
For $\mathfrak{A}\in\mathcal{C}_1$, we decide whether
$\mathfrak{A}\models \varphi$ as follows:  
Firstly, we compute 
$\mathfrak{B}:=\mathrm{Str}_I^{-1}(\mathfrak{A})$. 
Secondly, we solve the model checking problem
$\mathfrak{B}\models \mathrm{Frm}_I(\varphi)$. 

Similar to $\FO{}$-interpretations we can define
$\MSO{}$-interpretations: simply replace the $\FO{}$ formulas in $I$
by $\MSO$ formulas. Again, these can be used to reduce the
$\MSO$ model checking on one class of structures to another class, but
only if the interpretation is one-dimensional. If we use
an $n$-dimensional $\MSO$-interpretation for $n>1$, the resulting
transformation $\mathrm{Frm}_I$ translates $\MSO$ formulas into
second-order formulas as quantification over unary relations
is turned into quantification over $n$-ary relations.
As long as we stick to one-dimensional $\MSO$-interpretations, the
transformation $\mathrm{Frm}_I$ turns an $\MSO(\tau)$ formula into an
$\MSO(\sigma)$ formula and analogously to the previous lemma one
obtains the following statement.

\begin{lemma} \label{MSO-Interpretation-Properties}
  Let $I$ be a $1$-dimensional-$\sigma$-$\tau$ $\MSO$-interpretation,
  $\mathfrak{S}$ some $\sigma$-structure and $\varphi$ some
  $\MSO(\tau)$ sentence. Then 
  \begin{displaymath}
    \mathrm{Str}_I(\mathfrak{A}) \models \varphi \text{ iff }
    \mathfrak{A} \models \mathrm{Frm}_I(\varphi).
  \end{displaymath}
\end{lemma}

\section{Grids and  Trees}
\label{Sec:StrucandInt}
\subsection{A Grid-Like Structure}

Grid-like structures often play a crucial role in undecidability results
for model checking problems. 
In this section, we introduce a certain grid-like structure, namely, 
the bidirectional half-grid. 
It is a version of the upper half of the
$\N\times\N$ grid with an edge-relation for each direction, i.e.,
there are relations for the left, right, upward, and downward
successor. 

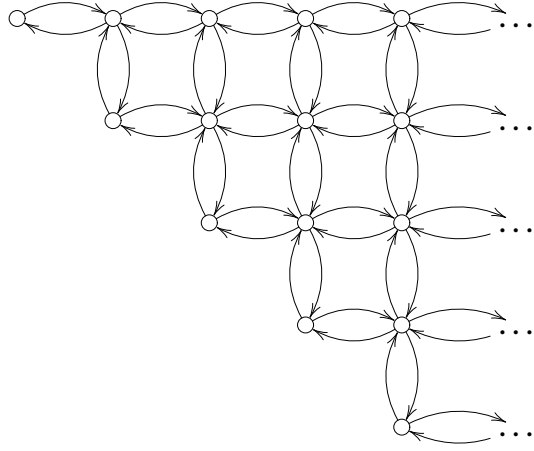
\begin{figure}
  \centering
  $
  \begin{xy}
    \xymatrix@R=30pt@C=30pt{
      *+[o][F]{ } \ar@/^/[r]^{} &        
      *+[o][F]{ } \ar@/^/[l]_{} \ar@/^/[r]^{} \ar@/^/[d]^{}&       
      *+[o][F]{ } \ar@/^/[l]_{} \ar@/^/[r]^{} \ar@/^/[d]^{}&         
      *+[o][F]{ } \ar@/^/[l]_{} \ar@/^/[r]^{} \ar@/^/[d]^{}&         
      *+[o][F]{ } \ar@/^/[l]_{} \ar@/^/[r]^{} \ar@/^/[d]^{}&  
      {\dots} \ar@/^/[l]_{} \\
      &         
      *+[o][F]{ } \ar@/^/[r]^{}  \ar@/^/[u]_{}&        
      *+[o][F]{ } \ar@/^/[l]_{} \ar@/^/[r]^{} 
      \ar@/^/[d]^{} \ar@/^/[u]_{}&         
      *+[o][F]{ } \ar@/^/[l]_{} \ar@/^/[r]^{} 
      \ar@/^/[d]^{} \ar@/^/[u]_{}&       
      *+[o][F]{ } \ar@/^/[l]_{} \ar@/^/[r]^{} 
      \ar@/^/[d]^{} \ar@/^/[u]_{}&
      {\dots} \ar@/^/[l]_{}\\
      & &         
      *+[o][F]{ } \ar@/^/[r]^{} \ar@/^/[u]_{}&       
      *+[o][F]{ } \ar@/^/[l]_{} \ar@/^/[r]^{} 
      \ar@/^/[d]^{} \ar@/^/[u]_{}&
      *+[o][F]{ } \ar@/^/[l]_{} \ar@/^/[r]^{} 
      \ar@/^/[d]^{} \ar@/^/[u]_{}& 
      {\dots} \ar@/^/[l]_{} \\ 
      &  & &       
      *+[o][F]{ } \ar@/^/[r]^{} \ar@/^/[u]_{}&       
      *+[o][F]{ } \ar@/^/[l]_{} \ar@/^/[r]^{}
      \ar@/^/[d]^{} \ar@/^/[u]_{}& 
      {\dots}\ar@/^/[l]_{} \\
      &  & & &       
      *+[o][F]{ }\ar@/^/[r]^{}  \ar@/^/[u]_{}& 
      {\dots}  \ar@/^/[l]_{} \\
    }    \end{xy}
  $ 
  \caption{The 
    bidirectional half-grid.}
  \label{HALFGRID}
\end{figure}

\begin{definition}
  The half-grid is the structure 
  $\BiHalfgrid:=(H,\rightarrow, \leftarrow, \downarrow,\uparrow)$
  where
  \begin{align*}
    &H:=\{(i,j)\in\N\times\N: i\leq j\},\\
    &\rightarrow:=\left\{\big((i,j),(k,l)\big)\in H^2: i=k, j=l-1\right\},\\
    &\leftarrow:=\left\{\big((i,j),(k,l)\big)\in H^2: i=k, j=l+1\right\},\\
    &\downarrow:=\left\{\big((i,j),(k,l)\big)\in H^2: i=k-1,
      j=l\right\},\text{ and}\\
    &\uparrow:=\left\{\big((i,j),(k,l)\big)\in H^2: i=k+1, j=l\right\},
  \end{align*}
  See Figure \ref{HALFGRID} for a
  pictures of $\BiHalfgrid$.   
\end{definition}

Many \MSO model checking results can be reduced to
the question of tree-likeness or grid-likeness of the underlying
graphs. On the one hand, if a class of structures consists only of
structures that are  
similar to trees, e.g, structures with small tree-width, then the
\MSO model checking is effectively decidable. On the other hand, if a
class contains a grid-like structure then the \MSO model checking is
undecidable. We do not want to go into the details what grid-likeness
means exactly. But for our purposes, the crucial observation is that
the upper half of a 
grid is, of course, grid-like whence $\BiHalfgrid$
has undecidable $\MSO$-theory. In fact, we can even 
show undecidability of the $L\mu$-theory of this structure. 

\begin{lemma}\label{Lmu-Half-grid-undecidability}
  $L\mu$ model checking is undecidable on the bidirectional half-grid
  $\BiHalfgrid$.  
\end{lemma}
\begin{remark}
  Note that we consider $L\mu$ on the naked half-grid, i.e., without any
  additional propositions. Although this result is not very surprising
  for people familiar with $L\mu$, we have not found any
  proof of this lemma in the literature. The interested reader may find
  a detailed proof of this result in Appendix
  \ref{Appendix_Lmuundecidabilty} where we reduce the halting problem
  for Turing machines to $L\mu$ model checking on $\BiHalfgrid$. 
\end{remark}

Since $L\mu$ may be seen as a fragment of MLFP and of $\MSO$, the following
corollaries follow immediately.

\begin{corollary}
  MLFP and \MSO are undecidable on  $\BiHalfgrid$. 
\end{corollary}

These results play a crucial role in Section \ref{Lmu-FO-CPG}, where
we investigate the $L\mu$-theory of \FO{}-interpretations of
collapsible pushdown graphs.

\subsection{Words and Trees}
\label{Sec:wordsAndTrees}
If $\Sigma$ is a finite set (called alphabet) then $\Sigma^*$ denotes
the set of finite words over the alphabet $\Sigma$. 
For words $w_1,w_2\in \Sigma^*$, 
we write $w_1 \leq w_2$ if $w_1$ is a prefix of $w_2$. We write
$w_1 < w_2$ for $w_1 \leq w_2$ and $w_1 \neq w_2$. We denote by
$w_1\circ w_2$ (or simply  $w_1 w_2$) the concatenation of $w_1$ and
$w_2$. Furthermore, we write $w_1 \sqcap w_2$ for the
greatest common prefix of $w_1$ and $w_2$.
If $\lvert w \rvert =n$, we set $w_{-i}$ for $0\leq i \leq n$ to be
the prefix of $w$ of length $n-i$.  

We now turn to trees. In this thesis we only consider
binary trees. Most of the time we are concerned with finite trees, but
in Section \ref{SectionRamseyQuantifier} we have to treat infinite
trees as well. We 
use the word ``tree'' only for finite trees unless we explicitly
say otherwise. 

We call a set $D\subseteq \{0,1\}^*$ a \emph{tree domain}, if $D$ is prefix
closed, i.e., for each $d\in D$ and $d'\in\Sigma^*$ we have $d'\in D$
if $d'\leq d$. 

A \emph{$\Sigma$-labelled tree} is a mapping $T:D\rightarrow \Sigma$
for $D$ some tree domain. $T$ is called \emph{finite}, if $D$ is
finite; otherwise $T$ is an \emph{infinite} tree. 

For $d\in D$ we denote the \emph{subtree rooted at $d$} by
$\inducedTreeof{d}{T}$. This is the tree defined by
$\inducedTreeof{d}{T}(e):= T(de)$. 
For $T_1, T_2$ trees, we write $T_1 \initialTreeq T_2$ if $T_1$ is an
initial segment of $T_2$, i.e., if $\domain(T_1) \subseteq
\domain(T_2)$ and $T_2{\restriction}_{\domain(T_1)} = T_1$. 

We denote the \emph{depth} of the tree $T$ by 
\mbox{$\depth{T}\coloneqq\max\left\{\lvert t\rvert: t\in\domain(T)\right\}$}.

\label{HochPlusEinfuehrung}
For $T$ some tree with domain $D$,  let $D^+$ denote the set of minimal
elements of the complement of $D$, i.e., 
\begin{align*}
D^+=\{e\in\{0,1\}^*\setminus D:\text{ all proper ancestors of }e\text{
  are contained in }D\}.   
\end{align*}
In particular, note that $\emptyset^+=\{\varepsilon\}$.
Under the same assumptions, we write $D^\oplus$ for $D\cup D^+$. 
Note that $D^\oplus$ is the extension of the tree domain $D$ by one
layer.  

Sometimes it is useful to define trees inductively by describing the
subtrees rooted at $0$ and $1$. For this purpose we fix the following
notation. 
Let $\hat T_0$ and $\hat T_1$ be $\Sigma$-labelled trees and
$\sigma\in\Sigma$. Then we write 
$T\coloneqq \treeLR{\sigma}{\hat T_0}{ \hat T_1}$
for the $\Sigma$-labelled tree $T$ with the following three
properties
\begin{align*}
  &1.\ T(\varepsilon) = \sigma,& &2.\ \inducedTreeof{0}{T} = \hat
  T_0 \text{, and }& &3.\ \inducedTreeof{1}{T} = \hat T_1.&
\end{align*}
We call $\inducedTreeof{0}{T}$ the \emph{left subtree of $T$} and
$\inducedTreeof{1}{T}$ the \emph{right subtree of $T$}. 

We denote by $\Trees{\Sigma}$ the set of all finite
$\Sigma$-labelled trees and by $\infTrees{\Sigma}$ the set of all
infinite $\Sigma$-labelled trees. 
We set 
$\allTrees{\Sigma}:=\Trees{\Sigma} \cup \infTrees{\Sigma}$ to be the
set of all finite or infinite 
$\Sigma$-labelled trees. 
We call the elements of $\Trees{\Sigma}$ \emph{trees}
without referring to finiteness. This convention is useful because we use
infinite trees only in Section \ref{SectionRamseyQuantifier}. In that
section we always clarify whether we talk about finite or infinite
trees. The elements of
$\infTrees{\Sigma}$ are called \emph{infinite trees}

For $T\in \allTrees{\Sigma}$  a finite or infinite tree, we write
$T^\Box$ for its 
lifting to the domain $\{0,1\}^*$ by padding with a special symbol
$\Box$, i.e.,  
\begin{align*}
T^\Box:\{0,1\}\rightarrow\Sigma\cup\{\Box\}, 
T^\Box(d):=
\begin{cases}
  T(d)&\text{if }d\in\domain(T),\\
  \Box & \text{otherwise.}
\end{cases}
\end{align*}

Note that we consider a $\Sigma$-word $w$ as a $\Sigma$-tree $t_w$
with domain $\{0^{i}: 0\leq i \leq \lvert w \rvert -1\}$ where
$t_w(0^{i-1})$ is labelled by the $i$-th letter of $w$.


\section{Generalised Pushdown Graphs}
\label{Chapter_InfiniteStructures}
In this section we introduce the objects of our study, namely, the
class of collapsible pushdown graphs and the class of nested pushdown
trees. Both classes generalise the class of pushdown graphs. It will
turn out that nested pushdown trees may be seen as a subclass of the
collapsible pushdown graphs which has a nicer algorithmic behaviour
than the class of all collapsible pushdown graphs. 
We start by recalling the well-known basics on pushdown systems.
Then we present nested pushdown trees (NPT) and  
in the last part we present collapsible pushdown
graphs (CPG). 

\subsection{Pushdown Graphs}

A pushdown system  is  a finite automaton extended
by a stack. These systems were first developed in formal language
theory. Used as word- or tree acceptors, pushdown systems recognise
exactly the context-free languages. 

We are interested in the model checking properties of graphs generated
by  generalisations of pushdown system. 
The graph of a pushdown system is the graph of all reachable
configurations where the edge-relation is induced by the transition
relation of the pushdown system. 

We briefly recall the  definitions and present some classical
results on pushdown systems. 

\begin{definition}
  A \emph{pushdown system}
  is a tuple \mbox{$\mathcal{S}=(Q,\Sigma,\Gamma, q_I, \Delta)$}
  satisfying the following conditions. 
  $Q$ is 
  finite and it is called the set of states. It contains the initial
  state $q_I\in Q$.  
  $\Sigma$ is finite and is called the set of stack symbols.
  There is a special symbol $\bot\in\Sigma$ which is called the  bottom-of-stack symbol.
  $\Gamma$ is finite and it is called the  input alphabet. 
  \begin{align*}
    \Delta\subseteq Q\times \Sigma \times \Gamma \times Q
    \times \Op    
  \end{align*}
  is the 
  transition relation 
  where 
  \begin{align*}
    \Op\coloneqq \{\Pop{1},\Id\}\cup\left\{\Push{\sigma}:
      \sigma\in\Sigma\setminus\{\bot\}\right\}.    
  \end{align*}
  The elements of $\Op$ are called stack operations. Each stack
  operation induces a function $\Sigma^+ \to \Sigma^*$ as follows. 
  \begin{itemize}
  \item Let $w, w'\in \Sigma^+$ be words and  $\sigma\in\Sigma$ a
    letter such that $w=w'\sigma$. Then $\Pop{1}(w):=w'$.
  \item $\Id$ is the identity on $\Sigma^+$.
  \item Let $w\in\Sigma^+$. For each
    $\sigma\in\Sigma\setminus\{\bot\}$, we set
    $\Push{\sigma}(w):=w\sigma$. 
  \end{itemize}

  A \emph{configuration} of $\mathcal{S}$ is a tuple $(q,s)\in Q\times
  \Sigma^+$. 
  Let $\delta=(q,\sigma,\gamma,q',\op)\in\Delta$. We call $\delta$ a 
  \emph{$\gamma$-labelled transition}.
  $\delta$
  connects the configuration 
  $(q,s)$ with the configuration $(q',s')$ if 
  $s=\op(s')$. We set $(q,s)\trans{\gamma} (q',s')$ if there is a
  $\gamma$-labelled transition $\delta\in\Delta$ that connects $(q,s)$
  with $(q',s')$. 
  
  We call $\trans{}:=\bigcup_{\gamma\in\Gamma} \trans{\gamma}$ the
  transition relation of $\mathcal{S}$

  The \emph{configuration graph of} $\mathcal{S}$ (also called the
  \emph{graph generated by} $\mathcal{S}$) consists of all
  configurations that are reachable from the initial configuration
  $(q_0,\bot)$ via a path along $\trans{}$. 
\end{definition}
\begin{remark}
  We call a graph $\mathfrak{A}$ a \emph{pushdown graph} 
  if it is the graph  generated by some pushdown system
  $\mathcal{S}$.
  
  Without loss of generality, we assume that there is no
  transition of the form $(q,\bot,\gamma,q',\Pop{1})\in\Delta$. This means
  that we never remove the bottom-of-stack symbol from the
  stack. Thus, we never have to deal with an empty stack. 
\end{remark}

\begin{definition}   \label{Def:PSRunDefinition}
  Let $\mathcal{S}$ be a pushdown
  system. Let $C$ be the set of configurations of $\mathcal{S}$ and
  $\trans{}$ its transition relation 
  A \emph{run} $\rho$ of $\mathcal{S}$ is a sequence of configurations
  that are connected by transitions, i.e., a sequence
  $c_0 \trans{\gamma_1} c_1\trans{\gamma_2} c_2 \trans{\gamma_3} \dots
  \trans{\gamma_n}c_n$. 
  
  We call $\rho(i):=c_i$ the 
  \emph{configuration of $\rho$ at position $i$}. 
  We call $\rho$ a run from $\rho(0)$ to $\rho(n)$ and say that the
  \emph{length} of $\rho$ is \mbox{$\length(\rho):=n$}.

  We write $\Runs{\mathcal{S}}$ for the \emph{set of all runs of}
  $\mathcal{S}$. 

  For runs $\rho_1, \rho_2$ of a pushdown system we write 
  $\rho_1 \preceq \rho_2$ for the fact 
  that $\rho_1$ is an initial segment of $\rho_2$.
  We write $\rho_1 \prec \rho_2$ if $\rho_1$ is a proper initial
  segment, i.e., $\rho_1\preceq \rho_2$ and
  $\length(\rho_1)<\length(\rho_2)$.

  For runs 
  $\rho=c_0 \trans{\gamma_1} c_1\trans{\gamma_2} c_2 \trans{\gamma_3} \dots
  \trans{\gamma_n}c_n$  and 
  $\rho'=  c'_0 \trans{\gamma'_1} c'_1\trans{\gamma'_2} c'_2
  \trans{\gamma'_3} \dots 
  \trans{\gamma'_m}c'_m$ where $c_n=c'_0$
  we 
  define 
  \begin{align*}
  \pi:=\rho\circ\rho':=c_0 \trans{\gamma_1} c_1\trans{\gamma_2} c_2
  \trans{\gamma_3} \dots \trans{\gamma_n}c_n 
  \trans{\gamma'_1} c'_1\trans{\gamma'_2} c'_2 \trans{\gamma'_3} \dots
  \trans{\gamma'_m}c'_m    
  \end{align*}
  and we call $\rho\circ\rho'$ 
  the \emph{composition of $\rho$ and $\rho'$}. We also say that 
  \emph{$\pi$ decomposes as $\pi = \rho \circ \rho'$}.
\end{definition}

\begin{remark}
  Note that a run does not necessarily start in the initial
  configuration. This convention is
  useful for the analysis of decompositions of runs because every
  restriction $\rho{\restriction}_{[i,j]}$ of a run $\rho$ with $0\leq
  i \leq j \leq \length(\rho)$ is again a run. 

  In the following, we will often identify a run $\rho$ of length $n$
  with a function 
  from $\{0, 1, 2, \dots, n\}$ to $C$ that maps $i$ to $\rho(i)$. 
  This is a sloppy notation because there may be two different
  transitions $(q, \sigma, \gamma, q', \op)$ and
  $(q,\sigma,\gamma',q',\op)$ with $\gamma\neq \gamma'$ that give
  rise to two runs $(q,w)\trans{\gamma} (q', \op(w))$ and $(q,w)
  \trans{\gamma} (q', \op(w))$. In this case we would identify both
  runs with the same function $f$ where $f(0)=(q,w)$ and
  $f(1)=(q',\op(w))$. For simplicity, we will always assume that the
  configurations of a run already determine the whole run. 
\end{remark}

Perhaps the most important theorem concerning pushdown systems and
formal languages is the
so-called uvxyz-theorem or pumping lemma of Bar-Hillel et
al.\ \cite{Bar-HillelPS61}. 
It is a classical tool for proving
that a language is not context-free.  
The uvxyz-theorem states the following. Given a context-free language
$L$ there is a natural number $n\in\N$ such that for all words from
$w\in L$ of length at least $n$, there is a decomposition $w=uvxyz$
such that $uv^ixy^iz\in L$ for all $i\in\N$. 
There are elegant proofs of this theorem using context-free grammars. 

In Chapter \ref{Chapter_FO-NPT}, we are interested in the runs of
pushdown systems. Especially, we need to find a short run that is
similar to a given long run. Thus, we are interested in a version of
the uvxyz-theorem where we look at the run corresponding to a word $w$.
We want to find a decomposition such that we can remove certain
parts from the run and obtain a valid run (corresponding to some word
$uxz$ where $w=uvxyz$).  

In this form, the proof of the lemma is slightly more complicated than
in the version of context free languages. Thus, we start by giving an
auxiliary lemma. It says that the run of a pushdown system does not
depend on a prefix of the stack that is never read. 
A generalised version of this lemma for higher-order
pushdown systems can be found in \cite{Blumensath2008}.

\begin{definition}
  Let $w\in\Sigma^*$. Let $\rho$ be a run of a pushdown system. We
  set $(q_i,w_i):=\rho(i)$ for all
  $i\in\domain(\rho)$. If $w\leq w_i$ for all
  $i\in\domain(\rho)$,  we write  $w\prefixeq \rho$ and say that
  $\rho$ is prefixed by $w$.
\end{definition}

\begin{lemma} \label{LemmaBlumensath}
  Let $\rho$ be a run of some pushdown system $\mathcal{S}$ and let 
  $w\in \Sigma^*$ be some word such that $w\prefixeq \rho$. For each
  $i\in\domain(\rho)$, let $v_i$ denote 
  the suffix of $\rho(i)$ such that $\rho(i)=(q_i,w v_i)$ for some
  state $q_i\in Q$. 
  
  If $w'\in\Sigma^*$ ends with the same letter
  as $w$ then the function
  \begin{align*}
    \rho[w/w']:\domain(\rho) &\rightarrow Q\times\Sigma^* \\
    \rho[w/w'](i)&:= (q_i,w' v_i)
  \end{align*}
  is a run of $\mathcal{S}$.
\end{lemma}

The proof of this lemma is straightforward: just observe that any
stack operation commutes with the prefix replacement. The claim follows
by induction on $\domain(\rho)$. 
We are now prepared to state the uvxyz-theorem in a version for
pushdown systems.

\begin{lemma}[\cite{Bar-HillelPS61}] \label{PumpingLemmaTemplate}
  Let $\mathcal{S}$ be some pushdown system. There is a constant $n\in \N$ such
  that for every run $\rho$ of length greater than $n$ that starts in the
  initial configuration at least one of the following holds.
  \begin{enumerate}
  \item There is a decomposition $\rho = \rho_1\circ \rho_2 \circ \rho_3$,
    words $w_1< w_2$, and a state $q\in Q$  such that
    $\rho_2(0)=(q,w_1)$, $\rho_3(0)=(q,w_2)$,
    $\length(\rho_2)\geq 1$, and $\rho':=\rho_1\circ \rho_3[w_2/w_1]$
    is a run of 
    $\mathcal{S}$. 
  \item There is a decomposition $\rho = \rho_1\circ \rho_2 \circ
    \rho_3 \circ \rho_4 \circ \rho_5$, words $w_1 < w_2$ with equal
    topmost letter and states
    $q,q'\in Q$ such that 
    $\rho_2(0)=(q,w_1)$, $\rho_3(0)=(q,w_2)$, $\rho_4(0)=(q',w_2)$,
    $\rho_5(0)=(q',w_1)$, 
    $\length(\rho_2)+\length(\rho_4)\geq 1$, and 
    $\rho':=\rho_1\circ \rho_3[w_2/w_1]\circ \rho_5$ is a run of $\mathcal{S}$.  
  \item There is a decomposition $\rho = \rho_1\circ \rho_2 \circ \rho_3$ with
    $\length(\rho_2)\geq 1$ such that $\rho':=\rho_1\circ \rho_3$ is a
    run of $\mathcal{S}$. 
  \end{enumerate}
\end{lemma}
\begin{proof}
  We assume that $n\in\N$ is some large natural number (what large
  means can be obtained from the proof). Let $\rho$ be some 
  run such that $m:=\length(\rho) > n$. 
  In order to prove this claim, we look for configurations in the run
  that share the 
  same state and share the same topmost element on their stack. There
  are the following cases.
  \begin{enumerate}
  \item The run ends with a large stack: assume that $\rho$ ends in a
    stack $w$ with $\lvert w \rvert > \lvert\Sigma\times Q\rvert$.
    For each $i\leq \lvert w\rvert$, let $w_i$ be
    the prefix of $w$ of length $i$. 
    Let $n_i\leq\length(\rho)$ be maximal such that the stack at
    $\rho(n_i)$ is $w_i$. Set $(q_i, w_i):=\rho(n_i)$. 
    By pigeon-hole principle there are $j<k<\length(\rho)$ such that
    $q_j=q_k$ and $\TOP{1}(w_j)=\TOP{1}(w_k)$. 
    Since $n_k$ is maximal, the run
    \begin{align*}
      \rho':= \rho{\restriction}_{[0,n_j]} \circ
      \rho{\restriction}_{[n_k,\length(r)]}[w_k/w_j]      
    \end{align*}
    is well defined and
    satisfies the lemma. 
  \item 
    The run passes a large stack but ends in a small one:
    assume that $\rho$ ends in some word of length at most
    $\lvert\Sigma\times Q\rvert$. Furthermore, assume that $\rho$
    passes a word of length greater than 
    $\lvert\Sigma\times Q\rvert + 
    \lvert Q\times Q  \times \Sigma \rvert$.  

    Let $i_{\max}\in\domain(\rho)$ be a position such that the word at
    $\rho(i_{\max})$ has maximal length in $\rho$. 

    By assumption, it follows that for
    $\rho(i_{\max})=:(q_{\max}, w_{\max})$, 
    \begin{align*}
      \lvert w_{\max}\rvert > \lvert\Sigma\times Q\rvert + 
      \lvert Q\times Q  \times \Sigma \rvert.  
    \end{align*}
    For each $i\leq\lvert w_{\max}\rvert$, let $w_i$ be the prefix of
    $w_{\max}$ of length $i$.  
    For each 
    \begin{align*}
      \lvert \Sigma\times Q\rvert \leq i \leq 
      \lvert\Sigma\times Q\rvert +
      \lvert Q\times Q  \times \Sigma \rvert,      
    \end{align*}
    let $n_i\leq i_{\max}$ be maximal such that $\rho(n_i)=(q_i,w_i)$
    for some $q_i\in Q$. 
    Analogously, let $m_i\geq i_{\max}$ be minimal such that 
    $\rho(m_i+1)=(\hat q_i, \Pop{1}(w_i))$ for some $\hat q_i\in Q$. 

    Note that
    $w_i$ is the stack at $\rho(m_i)$ and $w_i\prefixeq
    \rho{\restriction}_{[n_i,m_i]}$ due to the   
    definition of $n_i$ and $m_i$. 
    
    By the pigeon-hole principle, there are $\lvert \Sigma\times
    Q\rvert
    \leq j < k \leq \lvert\Sigma\times Q\rvert + 
    \lvert Q\times Q  \times \Sigma \rvert$ such that $q_j=q_k, 
    \hat q_j=\hat q_k$ and $\TOP{1}(w_j)=\TOP{1}(w_k)$. 

    Then the run $\rho':= \rho{\restriction}_{[0,n_j]} \circ
    \rho{\restriction}_{[n_k,m_k]}[w_k/w_j] \circ
    \rho{\restriction}_{[m_j,\length(\rho)]}$ satisfies the lemma. 
  \item The run never visits a large stack, i.e., a stack of size
    greater than 
    $\lvert Q\times \Sigma\rvert + \lvert Q \times Q \times \Sigma
    \rvert$. Since there are only finitely many stacks of size smaller
    than this bound, in a long run of this form there is a 
    configuration which is visited twice and the subrun in between may
    be omitted.\qedhere
  \end{enumerate}
\end{proof}
Pushdown graphs form a class of finitely represented infinite graphs
with good model checking
properties. Almost fifty years ago, Buchi \cite{Buchi64} showed that
the reachability 
problem on pushdown graphs is decidable. This result was notably
extended by Muller and Schupp in the 80's as follows.

\begin{theorem}[\cite{MullerS85}]
  The \MSO-theory of every pushdown graph is decidable.
\end{theorem}

This result was important for the development of
software verification because of the following fact. 
A pushdown graph  naturally arises as the
abstraction of some programme using (first-order recursive)
functions. Given a programme, one can design a pushdown system that 
simulates the behaviour of this programme. Every run
of the pushdown system corresponds to a possible execution of the
programme. Here, the state of the pushdown system stores the
programme counter. This means that the state of the pushdown system
stores the 
line number that is executed by the programme in this step. 
If a function call occurs, the pushdown system does the
following. It writes the programme counter onto the stack, and the new
state is the first line of the function which is called. When
this function eventually terminates, the programme counter is
restored by reading the stack. While the programme counter is restored,
the topmost element of the stack is deleted.

Using this reduction, many problems occurring in software 
verification can be reduced to model checking on pushdown graphs. 
But this approach has a severe limitation: in the language of the
pushdown 
graph, \MSO cannot be used to define a function return corresponding
to a function call. 
Defining a return that corresponds to a certain call is equivalent to
defining a subrun of the pushdown system that starts at this function
call and forms a well-bracketed word (where we interpret
push operations as opening brackets and pop operations as closing
brackets). But it is well known that \MSO cannot define the language
of well-bracketed words (the so-called Dyck-languages).

Thus, if one wants to verify properties of a programme that involves a
comparison of the situation just before a function call with the
situation exactly after the return of the function, one cannot reduce this
problem to a model checking problem on pushdown graphs. 

In the next section we present nested pushdown trees. 
These generalise trees generated by  pushdown systems in such a way
that pairs of corresponding calls and returns become definable even in
first-order logic. Therefore, nested pushdown trees are suitable
abstractions for programmes if one wants to verify properties involving
the pairs of corresponding function calls and returns.

\subsection{Nested Pushdown Trees}
\label{ChapterNPT}

Alur et 
al.\ \cite{Alur06languagesof} 
proposed the study of the model checking problem on nested pushdown trees. 
A nested pushdown tree is the  tree generated by a pushdown system
where the pairs of corresponding push and pop operations are
marked by a new relation $\hookrightarrow$. 
This new relation is called \emph{jump-relation}. We stress that due
to this new relation, a nested pushdown tree is no tree. 

\begin{definition}
  Let $\mathcal{S}=(Q,\Sigma, \Gamma, \Delta, q_0)$ be a pushdown
  system. Then the \emph{nested pushdown tree} generated by
  $\mathcal{S}$ is 
  \begin{align*}
    \NPT(\mathcal{S}):=(R,(\trans{\gamma})_{\gamma\in\Gamma},
    \hookrightarrow)    
  \end{align*}
  where 
  $(R,(\trans{\gamma})_{\gamma\in\Gamma})$ 
  is the unfolding of the configuration graph of $\mathcal{S}$. 
  $R$ is the 
  set of all runs of $\mathcal{S}$ starting at the
  configuration $(q_0, \bot)$. For two runs
  $\rho_1, \rho_2\in R$, we have $\rho_1\trans{\gamma} \rho_2$ if
  $\rho_2$ extends $\rho_1$ by 
  exactly one $\gamma$-labelled transition. The binary relation
  $\hookrightarrow$ is called 
  \emph{jump-relation} and is defined as follows:
  let $\rho_1,\rho_2\in R$ with $\length(\rho_i)=n_i$ and 
\mbox{$\rho_1(n_1)=(q,w)\in
  Q\times\Sigma^*$}. Then 
  $\rho_1\hookrightarrow \rho_2$ if $\rho_1$ is an initial segment of $\rho_2$,
  $\rho_2(n_2)=(q',w)$ for some $q'\in Q$ and $w$ is a proper prefix of
  all stacks between $\rho_1(n_1)$ and $\rho_2(n_2)$, i.e., 
  $w<\rho_2(i)$   for all $n_1 <i< n_2$. 
\end{definition}

Alur et al.\ proved the following results concerning the model checking
properties of the class of nested pushdown trees. 
\begin{theorem}[\cite{Alur06languagesof}]
  \label{NPTLmuModelCheckingComplexity}  
  The $L\mu$ model checking problem for nested pushdown trees is in
  EXPTIME. 
\end{theorem}

\begin{lemma}[\cite{Alur06languagesof}]
  \label{MSO-NPT-MC-undecidable}
  The \MSO model checking problem for nested pushdown trees in
  undecidable. 
\end{lemma}
\begin{proof}
  Let $\mathcal{S}:=(\{0,1\}, \{a,\bot\}, \{A,P\}, (0, \bot), \Delta)$ with
  \begin{align*}
  \Delta=\left\{ 
    (0, \bot, A, 0, \Push{a}),
    (0, a, A, 0, \Push{a}),
    (0, a, P, 1, \Pop{1}),
    (1, a, P, 1, \Pop{1})
  \right\}.      
  \end{align*}
  Figure \ref{fig:nptExample} shows the nested pushdown tree generated
  by $\mathcal{S}$. 
  \begin{figure}[h]
    \centering 
    $ 
    \begin{xy}
      \xymatrix@R=20pt@C=20pt{
        0 \bot \ar[r]^-{A} \ar@{^{(}->}[dr] \ar@{^{(}->}@/_18pt/[drdr]
        \ar@{^{(}->}@/_40pt/[drdrdr] &  
        0 \bot a  \ar[r]^-{A} \ar[d]^{P}
        \ar@{^{(}->}[dr] \ar@{^{(}->}@/_18pt/[drdr]& 
        0 \bot a a  \ar[r]^-{A} \ar[d]^-{P} \ar@{^{(}->}[dr]& 
        0 \bot a a a \ar[d]^{P}& \dots \\
        & 
        1 \bot  &
        1 \bot a \ar[d]^P&  
        1 \bot a a \ar[d]^P & \dots \\
        & 
        &
        1 \bot  &  
        1 \bot a  \ar[d]^P&  
        \dots \\
        & 
        &
        & 
        1 \bot  &  
        \dots \\
      }
    \end{xy}
    $
    \caption{Example of a nested pushdown tree.}
    \label{fig:nptExample}
  \end{figure}
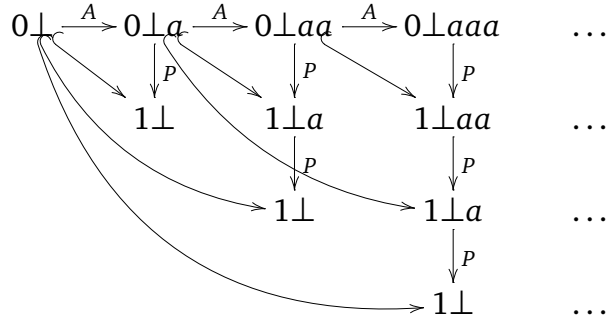
  We now show that the bidirectional halfgrid $\BiHalfgrid$ is
  \MSO-interpretable in 
  this graph. Application of Lemmas \ref{Lmu-Half-grid-undecidability}
  and \ref{MSO-Interpretation-Properties}
  then directly yields the claim.  
  
  As an abbreviation, we use the binary relation $\mathrm{REACH}_{P^*}$
  which holds for configurations $(c_1,c_2)$ if and only if $c_2$ is
  reachable from $c_1$ by a $P$-labelled path. This predicate is clearly
  \MSO-definable and in the structure $\NPT(\mathcal{S})$ it describes
  reachability along the columns.  
  Now, we define the next-column relation by 
  \begin{align*}
    \varphi_{nc}(x,y)\coloneqq\exists z_1, z_2 (\mathrm{REACH}_{P^*}
    z_1 x \wedge \mathrm{REACH}_{P^*} z_2 y \wedge z_1 \trans{A} z_2).    
  \end{align*}
  Similarly, we can define a next-diagonal relation by
  \begin{align*}
    \varphi_{nd}(x,y) \coloneqq\exists
    z_1, z_2 ( z_1 \hookrightarrow x \wedge z_2 \hookrightarrow y
    \wedge z_1 \trans{A} z_2).    
  \end{align*}
   We conclude that $\varphi_{nc}\land\varphi_{nd}(x,y)$ holds if and
   only if $y$ 
  is the right neighbour of $x$ in the  half-grid. 
  Thus, $\downarrow$ coincides with $\trans{\Pop{1}}$ and 
  $x \rightarrow y $ is defined by $\varphi_{nc} (x,y)\wedge
  \varphi_{nd}(x,y)$. 
  Switching the roles of $x$ and $y$, we can also define $\leftarrow$
  and $\uparrow$.

  This completes the interpretation of $\BiHalfgrid$ in
  $\NPT(\mathcal{S})$. Lemmas \ref{Lmu-Half-grid-undecidability}
  and \ref{MSO-Interpretation-Properties} then yield the claim. 
\end{proof}
\begin{remark}
  Even though $\MSO$ model checking for nested pushdown trees is
  undecidable, nested pushdown trees form an interesting class for
  software verification. Many interesting properties of programmes are
  expressible in $L\mu$. Moreover, the jump-relation allows to use
  $L\mu$ in order to express properties concerning corresponding push
  and pop operations. Such properties are not expressible when using 
  \MSO on pushdown graphs. 
\end{remark}

We have seen that \MSO is undecidable on nested pushdown trees while
$L\mu$ model checking is decidable. 
This difference concerning decidability of \MSO and $L\mu$
model checking turns nested 
pushdown trees into an interesting  class of structures from a model
theoretic point of view. Natural classes of graphs tend to have either
decidable MSO and $L\mu$ model checking or undecidable MSO and
$L\mu$ model checking.  
Beside the class of nested pushdown trees we only know of one other
natural 
class that does not follow this general rule: the class of collapsible
pushdown graphs. In 
Section \ref{Chapter_FO-NPT} we will show that nested pushdown trees
and collapsible pushdown graphs are
closely related via a simple \FO{}-interpretation. 
This relationship between nested pushdown trees and collapsible pushdown
graphs (of level $2$) will motivate our definition of the hierarchy of
higher-order nested pushdown trees in Section \ref{Chapter HONPT} in
analogy to the hierarchy of collapsible pushdown graphs. But before we
come to this generalisation of the concept of a nested pushdown tree,
let us introduce collapsible pushdown graphs.

\subsection{Collapsible Pushdown Graphs}
\label{STACS:SecCPG}
Before we introduce Collapsible pushdown graphs (CPG)
in
detail, we fix some notation. Then, we informally explain collapsible
pushdown systems. Afterwards, we formally introduce these
systems and the graphs generated by them. We conclude this section with some basic results on runs of collapsible pushdown systems. In Chapter
\ref{CPG-Tree-Automatic} we will then investigate \FO{} model
checking  on collapsible pushdown graphs. 

For some alphabet $\Sigma$, we inductively define $\Sigma^{*n}$ and
$\Sigma^{+n}$ for all $n\in\N\setminus\{0\}$ as follows. 
We set $\Sigma^{*1}:=\Sigma^*$, i.e., $\Sigma^{*1}$ is the set of all
finite words over
alphabet $\Sigma$. Then we set
\mbox{$\Sigma^{*(n+1)}\coloneqq(\Sigma^{*n})^*$}. 
Analogously, we write $\Sigma^{+1}:=\Sigma^+$ for the set of all
nonempty finite words over alphabet $\Sigma$ and we 
set \mbox{$\Sigma^{+(n+1)}\coloneqq(\Sigma^{+n})^+$}. Each element of
$\Sigma^{*n}$ is called an $n$-word. Stacks of a level $n$ collapsible
pushdown system are certain nonempty $n$-words over a special
alphabet.

Let us fix a  word $s\in \Sigma^{*(n+1)}$ of level $n+1$.
$s$  consists of an ordered list 
\mbox{$w_1, w_2, \dots, w_m$} of $n$-words, i.e., 
\mbox{$w_1, w_2, \dots, w_m \in\Sigma^{*n}$}. If 
we want to state this list of $n$-words explicitly, we
separate them by colons writing $s=w_1: w_2 : \dots : w_m$. 
By $\lvert  s \rvert$ we denote the number of $n$-words $s$
consists of, i.e., $\lvert  s\rvert = m$. We say $\lvert s \rvert$ is
the \emph{width} of $s$. We also use the notion of the height of an
$(n+1)$-word. The \emph{height} of $s$ is $\height(s):=\max\{\lvert
w_i \rvert: 1\leq i \leq m\}$ which is the width of the widest
$n$-word occurring in $s$.

Let $s'$ be another word of level $n+1$ such that
$s'=w_1':w_2':\dots w_l' \in \Sigma^{*(n+1)}$. We write 
$s:s'$ for the concatenation  $w_1: w_2 : \dots w_m:w_1':w_2':\dots:w_l'$.

If $s\in\Sigma^{*n}$, we denote by $[s]$ the $n+1$ word
that only consists of a list of one $n$ word which is 
$s$. We regularly omit the brackets if no confusion arises. 

Let $\Sigma$ be some finite alphabet. A level $n$ stack $s$ is an $n$-word
where each letter carries a link to some substack. Each link has a
certain level $1\leq i \leq n$. A level $i$ link points to some
$(i-1)$-word of the topmost level $i$ stack of $s$. 
Now, we first define the initial level $n$ stack; afterwards we
describe some stack operations that are used to generate all level $n$
stacks from the initial one. 

\begin{definition}
  Let $\Sigma$ be some finite alphabet with a distinguished
  bottom-of-stack symbol $\bot\in\Sigma$. 
  The \emph{initial stack} of level $l$ over $\Sigma$ is inductively
  defined as follows. 
  The initial level $1$ stack is $\bot_1\coloneqq \bot$. 
  For the higher levels, we set $\bot_n\coloneqq[\bot_{n-1}]$ to be
  the initial stack of level $n$. 
\end{definition}

We informally describe the operations that 
can be applied to a level $n$ stack. 
There are the following stack operations:
\begin{itemize}
\item The  push operation of level $1$, denoted by $\Push{\sigma,k}$ for
  $\sigma\in\Sigma$ and $1\leq k \leq n$, writes the symbol $\sigma$
  onto the topmost level $1$ stack and attaches a link of level
  $k$. This link points to the next to last entry
  of the topmost level $k$ stack.
\item For $2\leq i \leq n$, the push operation of level $i$ is denoted
  by $\Clone{i}$. It
  duplicates the topmost entry of the topmost 
  level $i$ stack.
  The links are preserved by $\Clone{i}$ in the following sense. 
  Let $s$ be some stack. Let  $a$ be a letter in the topmost level $i$
  stack of  $s$. Assume that $a$ has a link of
  level $j$. Let $a'$ be the copy of $a$ in
  $\Clone{i}(s)$. Then the link of $a'$ points to the unique level
  $j-1$ stack in the topmost level $j$ stack of $\Clone{i}(s)$ that is
  a clone of the $j-1$ stack to which the link of $a$ points. This
  means that for 
  $j\geq i$, $a$ and $a'$ carry links to the same stack. For $j<i$,
  the link of $a'$ points to the clone of the stack to which the
  link of $a$ points. 
\item The level $i$ pop operation $\Pop{i}$ for $1\leq i\leq n$ removes
  the topmost entry of the topmost level $i$ stack. Note
  that the
  $\Pop{1}$ operation corresponds to the ordinary pop in a pushdown
  system that just removes the topmost symbol from the stack.
\item The last operation is $\Collapse$. 
  The result of $\Collapse$ is determined by the link attached to the
  topmost letter of the stack. If we apply collapse to a stack $s$ where
  the link level of the topmost letter is $i$, then $\Collapse$
  replaces the topmost level $i$ stack of $s$ by the level $i$ stack
  to which the link points. 
  Note that the application of a collapse is equivalent to the
  application of a sequence of  
  $\Pop{i}$ operations where the link of the topmost letter controls
  how long this sequence is.
\end{itemize}

In the following, we formally introduce collapsible pushdown stacks
and the stack operations. 
We represent such a stack of letters with links as $n$-words over the
alphabet  
\mbox{$(\Sigma\cup (\Sigma\times\{ 2, \dots,
  n\}\times\N))^{+n}$}. We consider elements from $\Sigma$ as elements
with a link of level $1$ and elements $(\sigma,l,k)$ as letters with a
link of level $l$. 
In the latter case, the third component specifies the width of the
substack to which the link  points. 
For letters with link of level $1$, the position of this letter within
the stack already determines the stack to which the link points. Thus,
we need not explicitly specify the link in this case.

\begin{remark}
  Other equivalent definitions, for instance in \cite{Hague2008}, use
  a different 
  way of storing the links: they also store symbols $(\sigma,i,n)$ on
  the stack, but here $n$ denotes the number of $\Pop{i}$ transitions
  that are equivalent to performing the collapse operation at a stack
  with topmost element $(\sigma,i,n)$. The disadvantage of that
  approach is that the $\Clone{i}$ operation  cannot copy
  stacks. Instead, it can only copy the symbols stored
  in the topmost 
  stack and has to alter the links in the new copy. A clone of level
  $i$ must replace all links $(\sigma,i,n)$ by $(\sigma,i,n+1)$ in
  order to preserve the links stored in the stack. 
\end{remark}

Before we give a formal definition of the stack operations, we
introduce some auxiliary functions. 
\begin{definition}
  For $s=w_1:w_2:\dots: w_n\in(\Sigma\cup(\Sigma\times\{ 2, \dots
  l\}\times\N))^{+l}$, we define the 
  following  auxiliary functions:
  \begin{itemize}
  \item  For $1\leq k\leq l$, the \emph{topmost level $k-1$ word of
      $s$} is $\TOP{k}(s)\coloneqq 
    \begin{cases}
      w_n & \text{if } k=l,\\
      \TOP{k}(w_n) &\text{otherwise.}
    \end{cases}$
  \item For $\TOP{1}(s)=(\sigma,i,j)  \in
    \Sigma\times\{2,3,\dots,l\}\times \N$,
    we define the \emph{topmost symbol}
    \mbox{$\Sym(s)\coloneqq\sigma$}, the  \emph{collapse level of the
      topmost element}
    $\Lvl(s)\coloneqq i$, and the \emph{collapse link of the topmost element}
    $\Lnk(s)\coloneqq j$. 
    
    For $\TOP{1}(s)=\sigma  \in \Sigma$,
    we define the \emph{topmost symbol}
    \mbox{$\Sym(s)\coloneqq\sigma$}, the  \emph{collapse level of the
      topmost element}
    $\Lvl(s)\coloneqq 1$, and the \emph{collapse link of the topmost element}
    $\Lnk(s)\coloneqq \lvert \TOP{2}(s)\rvert-1$. 
  \item For $m\in\N$, we define
    $\mathrm{p}_{\sigma,k,m}(s)\coloneqq
    \begin{cases}
      s (\sigma,k,m) &\text{if } l=1,\\
      w_1: w_2 : \dots : w_{n-1}: \mathrm{p}_{\sigma,k,m}(w_n)& \text{otherwise.}
    \end{cases}$
  \end{itemize}
\end{definition}
These auxiliary functions are useful for the formalisation of the stack
operations. 
\begin{definition}
  For $s=w_1:w_2:\dots: w_n\in(\Sigma\cup(\Sigma\times\{ 2, 3,\dots
  l\}\times\N))^{+l}$,
  for $\sigma\in\Sigma\setminus\{\bot\}$,
  for $1\leq k \leq l$ and for $2\leq j \leq l$, we define
  the stack operations 
  \begin{align*}
    \Clone{j}(s)\coloneqq
    &\begin{cases}
      w_1: w_2 : \dots : w_{n-1}: w_n : w_n & \text{if } j = l\geq 2,\\
      w_1: w_2 : \dots : w_{n-1}: \Clone{j}(w_n) & \text{otherwise.}
    \end{cases}\\
    \Push{\sigma,k}(s)\coloneqq
    &\begin{cases}
      s\sigma  &\text{if } k=l=1,\\
      \mathrm{p}_{\sigma,k,n-1}(s) &\text{if } k=l\geq 2,\\
      w_1: w_2 : \dots : w_{n-1}: \Push{\sigma,k}(w_n)& \text{otherwise.}
    \end{cases} \\
    \Pop{k}(s)\coloneqq
    &\begin{cases}
      w_1: w_2 : \dots : w_{n-1}: \Pop{k}(w_n) & \text{if } k<l,\\
      w_1: w_2 : \dots : w_{n-1} & \text{if } k=l, n>1,\\
      \text{undefined} & \text{otherwise, i.e.,} k=l, n=1. 
    \end{cases}\\
    \Collapse{}(s)\coloneqq 
    &\begin{cases}
      w_1: w_2 : \dots : w_m & \text{if } \Lvl(s)=l, \Lnk(s)=m>0,\\
      w_1: w_2 : \dots : w_{n-1}: \Collapse{}(w_n) & \text{if }
      \Lvl(s)<l,\\
      \text{undefined}&\text{if } \Lnk(s)=0.
    \end{cases}    
  \end{align*}
  The \emph{set of level $l$ operations} is 
  \begin{align*}
    \Op_l\coloneqq\{(\Push{\sigma,k})_{\sigma\in\Sigma,k\leq l},
    (\Clone{k})_{2\leq k\leq l}, (\Pop{k})_{1\leq k\leq l},
    \Collapse{}\}.
  \end{align*}
  The \emph{set of level $l$ stacks}, $\Stacks_l(\Sigma)$, is the smallest set
  that contains $\bot_l$ and is closed under application of operations
  from $\Op_l$. 
\end{definition}
\begin{remark}
  It is sometimes convenient to assume that
  the identity
  \begin{align*}
    \Id:\Stacks_l(\Sigma)\to\Stacks_l(\Sigma)    
  \end{align*}
  is also a
  stack operation. Whenever this assumption is useful, we assume $\Id$
  to be a stack operation.  
\end{remark}
We illustrate the definition of the stack operations with the
following example. 
\begin{example}
  We start with the level $3$ stack 
  $s_0:=\left[ \bot  \right] : \left[  \bot   : \bot \right]$. 
  We have 
  \begin{eqnarray*}
    \Push{a,2}(s_0) &=& \left[\bot \right]  : 
      \left[ \bot  :  \bot (a,2,1) \right] =:s_1\\
    \Push{b,3}(s_1) &=&  \left[\bot  \right] :
     \left[ \bot  :  \bot (a,2,1) (b,3,1)
      \right]=:s_2 \\
    \Clone{3}(s_2) &=&  \left[\bot
      \right]  :  \left[ \bot :  \bot
        (a,2,1) (b,3,1) \right]:  \left[ \bot :
       \bot  (a,2,1) (b,3,1) \right] =: s_3\\
    \Clone{2}(s_3) &=& s_2 :  \left[ \bot  :
       \bot  (a,2,1) (b,3,1) :  \bot
      (a,2,1) (b,3,1) \right]=: s_4 \\
    \Collapse(s_4) &=& \left[ \bot \right]\\
    \Pop{1}(s_4) &= & s_2: \left[ \bot  :
       \bot  (a,2,1) (b,3,1) : \bot
      (a,2,1) \right] =:s_5\\
    \Collapse(s_5) &=&  s_2:[\bot] = \left[\bot
      \right]  :  \left[ \bot : \bot
        (a,2,1) (b,3,1) \right]:  \left[ \bot \right]. 
  \end{eqnarray*}
\end{example}

Note that $\Collapse$ and $\Pop{k}$ operations are only allowed if
the resulting stack is nonempty. This avoids the special treatment of
empty stacks. Furthermore, any $\Collapse$ that works on level $1$
is equivalent to one $\Pop{1}$ operation: 
level $1$ links always point to the preceding letter
because there is no $\Clone{1}$ operation. Furthermore,  every
$\Collapse$ that works on a level $i\geq 2$ is equivalent to a
sequence of $\Pop{i}$ operations. 

Let us now define the substack relation on collapsible pushdown
stacks. It is the natural generalisation of the prefix order on words.  
\begin{definition}
  Let $s,s'\in\Stacks_l(\Sigma)$. We say that $s'$ is a substack of $s$ if
  there are  $n_i\in\N$ for $1\leq i \leq l$ such that
  $s' = \Pop{1}^{n_1}( \Pop{2}^{n_2}(\dots(\Pop{l}^{n_l}(s))))$.
  We write $s' \leq s$ if $s'$ is a substack of $s$. 
\end{definition}

Now, it is time to formally define collapsible pushdown systems. These
are defined completely analogously to pushdown systems but using a
level $l$ stack and all the level $l$ stack operations. 

\begin{definition}
  A \emph{collapsible pushdown system} of level $l$ ($l$-\CPS) is
  a tuple 
  \begin{align*}
    \mathcal{S} = (Q,\Sigma, \Gamma, \Delta, q_0)    
  \end{align*}
  where $Q$ is a
  finite set of states, $\Sigma$  a
  finite stack  alphabet with a distinguished bottom-of-stack symbol
  $\bot\in\Sigma$, $\Gamma$  a finite 
  input alphabet,  $q_0\in Q$ the initial state, and 
  \begin{align*}
    \Delta\subseteq
    Q\times \Sigma \times\Gamma \times Q \times \Op_l    
  \end{align*}
  the transition relation.

  A level $l$ \emph{configuration} is a pair $(q,s)$ where $q\in Q$ and
  $s\in\Stacks_l(\Sigma)$. 
  For \mbox{$q_1,q_2\in Q$} and $s,t\in \Stacks_l(\Sigma)$ we 
  define a $\gamma$-labelled
  transition $(q_1,s) \trans{\gamma} (q_2, t)$ if there is
  a $(q_1, \sigma, \gamma, q_2, op)\in\Delta$ such that $\op(s)=t$ and
  \mbox{$\Sym(s) = \sigma$}. 

  We call $\trans{}:=\bigcup_{\gamma\in\Gamma} \trans{\gamma}$ the
  transition relation of $\mathcal{S}$. 
  We set $C(\mathcal{S})$ to be the set of all
  configurations that are reachable from $(q_0,\bot_l)$ via
  $\trans{}$
  and call $C(\mathcal{S})$ the set of \emph{reachable}
  or \emph{valid} configurations. The 
  \emph{collapsible pushdown graph (\CPG) generated by $\mathcal{S}$} is 
  \begin{align*}
    \CPG(\mathcal{S})\coloneqq\big(C(\mathcal{S}), (C(\mathcal{S})^2\cap
    \trans{\gamma})_{\gamma\in\Gamma}\big)
  \end{align*}
\end{definition}

\begin{remark} 
  \begin{itemize}
  \item[]
  \item Note that the transitions of a
  collapsible pushdown system only depend on the state and the topmost
  symbol, but not on the topmost collapse level and collapse link. The latter
  are only used to handle the result of a collapse operation. 
  \item
    In the following,  we always assume that the label of each
    transition carries information about the stack operation and the
    state that is reached, i.e., we assume that there is a map 
    $f:\Gamma\to Q\times \Op$ such that for each transition
    $(q,\sigma,\gamma,q',\op)\in \Delta$ we have
    $f(\gamma)=(q',\op)$. It is obvious that each collapsible
    pushdown system can be transformed into  one that satisfies this
    assumption: use $\Gamma\times Q\times\Op$ as new input alphabet;
    then $\trans{\gamma}=\bigcup_{q\in Q, \op\in\Op}
    \trans{(\gamma,q,\op)}$. In this sense, we will write
    $\trans{q,\op}:= \bigcup_{f(\gamma)=(q,\op)} \trans{\gamma}$ and
    also
    $\trans{q}:=\bigcup_{\op\in\Op} \trans{q,\op}$ and 
    $\trans{\op}:=\bigcup_{q\in Q} \trans{q,\op}$.
  \item
    An \emph{higher-order pushdown system} is a collapsible pushdown
    system that does not use the collapse operation. 

    To be more precise, we call a collapsible pushdown system with
    transition relation $\Delta$ 
    an \emph{higher-order  pushdown system} if
    \begin{align*}
      \Delta\subseteq
      Q\times \Sigma \times\Gamma \times Q \times
      \left(\Op_l\setminus\left(\{\Collapse\}\cup\{\Push{\sigma,i}: i\geq 2\}\right)\right),      
    \end{align*}
    i.e., if it does not use the collapse operation and the links of
    level $i$ for all $i>1$.
  \end{itemize}
\end{remark}

\begin{example}
  The following example of a collapsible pushdown graph $\mathfrak{G}$
  of level $2$  is taken from 
  \cite{Hague2008}. Let $Q\coloneqq\{0,1,2\}, \Sigma\coloneqq
  \{\bot,a\}$, 
  $\Gamma:=\{\mathrm{Cl}, A, A', P, \mathrm{Co}\}$.
  $\Delta$ is given by $(0,-,\mathrm{Cl},1,\Clone{2})$,
  $(1,-,A,0,\Push{a,2})$, $(1,-,A',2,\Push{a,2})$,
  $(2,a,P,2,\Pop{1})$, and $(2,a,\mathrm{Co},0,\Collapse)$, where $-$
  denotes any 
  letter from $\Sigma$. \\
  \begin{figure}[h]
    \centering 
    $
    \begin{xy}
      \xymatrix@R=12pt@C=9pt{ 
        0 \bot \ar[r]^-{\mathrm{Cl}}& 
          1 \bot:\bot \ar[r]^-{A} \ar[d]^{ A'}& 
          0 \bot:\bot a \ar[r]^-{\mathrm{Cl}} &
          1 \bot:\bot a:\bot a  \ar[r]^-{A} \ar[d]^-{A'}& 
          0\bot:\bot a: \bot aa \ar[r]^-{\mathrm{Cl}}  
          & 1 \bot:\bot a : \bot aa:\bot aa \ar[d]^{A'}& \dots \\
          & 
          2 \bot:\bot a  \ar[d]^P \ar[ul]^{\mathrm{Co}}&  
          & 
          2 \bot:\bot a:\bot aa \ar[d]^P \ar[ul]^{\mathrm{Co}}&  
          & 
          2 \bot:\bot a : \bot aa:\bot aaa\ar[d]^P
            \ar[ul]^{\mathrm{Co}}
          & \dots \\
        & 
          2 \bot:\bot  &  
          &
          2 \bot:\bot a:\bot a  \ar[d]^P \ar[uulll]^{\mathrm{Co}}&  
          & 
          2 \bot:\bot a : \bot aa:\bot aa
          \ar[d]^P\ar[uulll]^(.415){\mathrm{Co}}
          & \dots \\
        &                &  
          & 
          2 \bot:\bot a:\bot  &  
          & 
          2 \bot:\bot a : \bot aa:\bot a \ar[d]^P\ar[uuulllll]^(.3){\mathrm{Co}}
          & \dots \\
        &              
          &  
          &
          &
          &
          2 \bot:\bot a : \bot aa:\bot  & \dots 
      }
    \end{xy}
    $
    \caption{Example of the  $2$-\CPG $\mathfrak{G}$ (the level
      $2$ links of the letters $a$  are omitted due to space
      restrictions).} 
    \label{fig:CPGExample}
  \end{figure}
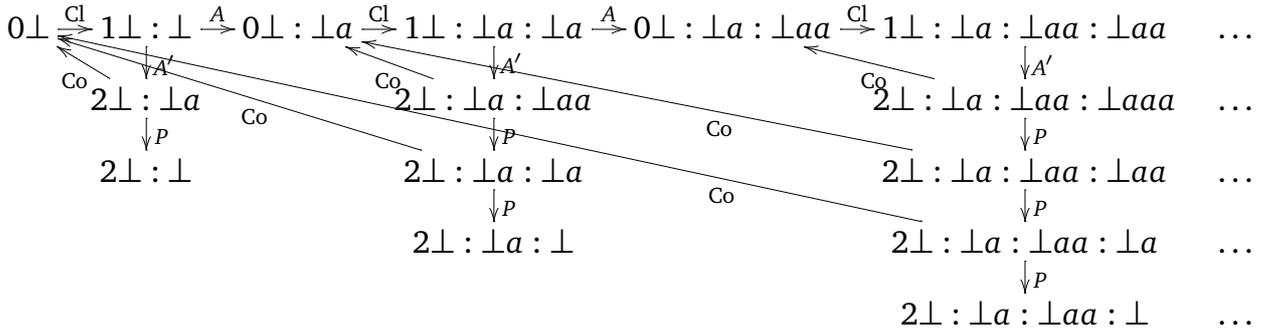
\end{example}

The next two theorems summarise the known results concerning model
checking on collapsible pushdown graphs.

\begin{theorem}[\cite{Hague2008}]
 There is a collapsible pushdown graph of level $2$ with undecidable
  \MSO model checking.
\end{theorem}
\begin{proof}
  The graph from figure \ref{fig:CPGExample} is an example. 
  Note that the graph from figure \ref{fig:nptExample} is clearly
  $\FO{}$-interpretable in this graph; one merely has to reverse the
  collapse-edges in order to obtain the jump-edges and to omit every
  second node in the topmost line. Hence, the corresponding
  \MSO undecidability result from theorem \ref{MSO-NPT-MC-undecidable}
  applies also to this collapsible pushdown graph. 
\end{proof}

\begin{theorem}[\cite{Hague2008}] \label{CPGLmuModelCheckingComplexity}
  $L\mu$ model checking on level $n$ collapsible
  pushdown graphs is \mbox{$n$-EXPTIME} complete. 
\end{theorem}
We briefly sketch the proof idea for this theorem. 
The proof uses parity-games on collapsible pushdown
graphs. It is commonly known that $L\mu$ model checking and the
calculation of winning regions in a parity-game are equivalent. 
In order to solve parity-games on a collapsible pushdown graph of
level $l+1$, Hague et al.\ reduce this problem to another parity game
on a collapsible pushdown graph of level $l$. Their proof consists of
two steps.
\begin{enumerate}
\item First, they prove that for each collapsible pushdown graph
  there is another one of the same level that is \emph{rank
    aware}. A level $l$ collapsible pushdown graph is rank aware if it
  ``knows'' at each configuration $(q,s)$ with $\Lvl(s)=l$ 
  the minimal rank (or priority) that was
  visited since the last occurrence of the stack $\Collapse(s)$. 
  One can show that for each parity game on a collapsible pushdown
  graph, one can construct a parity game on a rank aware collapsible
  pushdown graph such that every winning strategy in this new  game can be
  transformed into a winning strategy for the original game. 
\item In the second step, Hague et al.\ reduce the problem of solving a
  parity game on a rank aware collapsible pushdown graph of level $l+1$
  to the problem of solving a game on a graph of level $l$. The basic
  idea is to simulate only the topmost level $l$ stack of the level
  $l+1$ graph and to handle the attempt to use a $\Clone{l+1}$ at a
  certain stack $s$ in the
  following way: if one player would perform a $\Clone{l+1}$ operation
  in the original game, Verifier has to make a certain claim
  about her winning strategy in the original game. According to her 
  winning strategy, for each priority $i$, there is a set of
  states $Q_i$  such that whenever the game
  returns to $s$ and the 
  minimal priority between the $\Clone{l+1}$ operation and this new
  occurrence of the stack $s$ is $i$, then the stack $s$ is
  visited again in one of the states from $Q_i$. Now, in the new game,
  Verifier moves to a state representing the set $(Q_i)_{i\in
    P}$ where $P$ is the finite set of priorities. Falsifier now has two
  choices. Either he believes Verifier 
  or he does not
  believe that Verifier's claim is correct. 

  If he believes her, he
  chooses one of the
  $i\in P$ and a $q\in Q_i$. The new game continues in $(q,s)$
  after visiting an auxiliary state of priority $i$. 

  Otherwise, the new game continues
  with the stack $\TOP{l-1}(s)$ and Falsifier has to show that there is some
  position where he could use a $\Pop{l}$- or $\Collapse$ operation
  (in the original game) and return to some state $q$ and stack
  $s$ such that $q\notin Q_i$ for $i$ the least priority visited since
  Verifier had made this claim. 
  At this point, rank awareness comes into play. 
  At each position where a $\Pop{2}$- or $\Collapse$ operation may be
  performed, rank-awareness allows to determine the minimal priority
  $i$ since Verifier had made her claim.
  Thus, due to rank-awareness, we can check whether Falsifier managed to
  reach a position $(q,s)$ where $q\notin Q_i$. In this case,
  Falsifier wins the game. If $q\in Q_i$ then the game ends and
  Verifier wins. 
\end{enumerate}
Using this reduction $l-1$ times, one derives a parity game on a
level $1$ pushdown graph such that a strategy on this game can be used in order
to compute a strategy in the original parity
game. Walukiewicz \cite{Wal01} showed the 
solvability of parity games on level $1$ pushdown graphs. 
Now, the decidability of $L\mu$ model checking on collapsible pushdown
graphs follows by induction on the level of the graph. 

Since \MSO model checking is undecidable for collapsible pushdown
graphs, it is  interesting to investigate model checking for fragments of
\MSO. What is the largest fragment such that model checking on
collapsible pushdown graphs is decidable? 
In Chapter \ref{CPG-Tree-Automatic}, we make a first step towards an
answer to this question. We prove the decidability of
the first-order model checking on level $2$ collapsible pushdown
graphs extended by $L\mu$-definable predicates. 

Now, we come to the notion of a run of a collapsible pushdown
system. This definition is completely analogous to the definition of a
run of a pushdown system (cf. Definition \ref{Def:PSRunDefinition}).

\begin{definition}
  Let $\mathcal{S}$ be a collapsible pushdown system. 
  A \emph{run} $\rho$ of $\mathcal{S}$ is a sequence of configurations
  that are connected by transitions, i.e., a sequence
  \begin{align*}
    c_0 \trans{\gamma_1} c_1\trans{\gamma_2} c_2 \trans{\gamma_3} \cdots
    \trans{\gamma_n}c_n.    
  \end{align*}
\end{definition}
\begin{remark}
  As in the case of pushdown systems, we identify a run $\rho$
  of length $n$ with the function that maps a number $i$ to the 
  configuration occurring just after the $i$-th transition in $\rho$
  for each $0\leq i \leq n$, i.e., $\rho(i)$ denotes the configuration
  after the $i$-th transition of $\rho$ and especially  $\rho(0)$ is
  the first configuration of $\rho$.   
\end{remark}

The final part of this section consists of some basic results
concerning runs of collapsible pushdown systems of level $2$. 
We focus on level $2$ because all of our main results only treat
pushdown systems of level $2$. 

First, we come to the question whether certain runs can create links
to certain stacks. 
Consider some configuration $(q,s)$ of a level $2$ collapsible pushdown
system. If $\lvert s \rvert=n$ then 
a $\Push{\sigma,2}$ transition applied to $(q,s)$ 
creates a letter with a link to the substack
of width $n-1$. 
Thus, links to the substack of width $n-1$ in some word above the
$n$-th one are always created by a $\Clone{2}$ operation. A
direct consequence of this fact 
is the following lemma.

\begin{lemma} \label{Lem:Howtwolinksevolve}
 Let $s$ be some level $2$ stack with $\TOP{1}(s)=(\sigma,2,k)$. 
 Let $\rho$ be a run of a pushdown system of level $2$
 that starts with
 stack $s$, that passes 
 $\Pop{1}(s)$, and that ends in $s$. If $k<\lvert s \rvert -1$ then
 $\rho$  passes $\Pop{2}(s)$.  
\end{lemma}
The proof is left to the reader. 
Later we often use the contraposition of this statement. We use the
fact that a certain 
run to $s$ does not visit $\Pop{2}(s)$ and conclude that it cannot visit
$\Pop{1}(s)$. 

The next result deals with the 
\emph{decision problem for configurations}: given a  
collapsible pushdown system 
$\mathcal{S}$, and a configuration $(q,s)$, is
$(q,s)\in\CPG(\mathcal{S})$?
We can solve this problem using the decidability of $L\mu$ model
checking on collapsible pushdown systems.

In the following we reduce the decision problem for
configurations  for a level $2$ collapsible pushdown system
$\mathcal{S}$ to the 
$L\mu$ model checking on a variation of $\mathcal{S}$. 
The proof is based
on the idea that a stack is uniquely determined by its top element and
the information which substacks can be reached via $\Collapse$ and 
$\Pop{i}$. 

We can compute a variant $\mathcal{S}^{(q,s)}$ of a given
$\CPS$ $\mathcal{S}$ 
such that $\mathcal{S}^{(q,s)}$ satisfies a certain
$L\mu$ formula if and only if $(q,s)$ is a configuration of the graph
generated by $\mathcal{S}$. The new 
pushdown system is the extension of $\mathcal{S}$ by a testing device
for the configuration $(q,s)$. 
Let us describe this testing device. 

Assume that we want to define a  testing device for the configuration
$(q,s)$. 
Furthermore, assume that for each configuration $(q',s')$ where $s'$ is a
proper substack of $s$, there already is a testing device for configuration
$(q',s')$. 
The testing device for $(q,s)$ works as follows. 

Whenever the system is in some configuration $(q,\hat s)$, it switches to
$(q_s, \hat s)$ where $q_s$ is a new ``testing state''. In $q_s$, the
system checks 
whether $\TOP{1}(\hat s) = \TOP{1}(s)$. If this is the case, then 
the following happens. 
Let $\hat s'$ be the stack obtained from $\hat s$ by removing the
topmost element 
and let $s'$ be the stack obtained from $s$ by removing the
topmost element. 
Now, we start the testing device for
the substack $s'$ on the stack $\hat s'$.  
If this testing device returns that $\hat s'$ is $s'$, then $\hat s=s$
and the new testing device was started in $(q,s)$. 

For each configuration $(q,s)$, there is an $L\mu$ formula such that this
formula is satisfied at some configuration of $\mathcal{S}^{(q,s)}$ if
and only if this configuration is $(q,s)$. 

Before we go into the details of this proof, we recall the terminology
concerning $L\mu$ on collapsible pushdown graphs. 
The binary relations on such a graph are labelled by symbols from the input
alphabet $\Gamma$ and we  use expressions as $\langle \gamma
\rangle \varphi$ for the formula saying ``there is a $\gamma$-labelled
edge leading to a node where $\varphi$ holds''. 
As an abbreviation we use $\Diamond \varphi$ for the formula
saying ``there is an arbitrary labelled edge leading to a node where
$\varphi$ holds'', i.e., as an abbreviation for 
$\bigvee_{\gamma\in\Gamma} \langle \gamma \rangle \varphi$. 

\begin{lemma} \label{ReachableConfigsinMuCalcul}
  Given some \CPS 
  \mbox{$\mathcal{S}=(Q,\Sigma, \Gamma, \Delta, q_0)$} of level $2$,
  some $q\in Q$ and some  
  stack $s$,
  it is decidable  whether $(q,s)$ is a reachable
  configuration 
  of $\mathcal{S}$, i.e., whether $(q,s)$ is a vertex of $\CPG(\mathcal{S})$.
\end{lemma}

\begin{proof}
  For $q\in Q$ and $s$ a stack,  we define a system $\mathcal{S}^{(q,s)}$ and
  a formula $\psi_{(q,s)}\in L\mu$,  such that 
  \begin{align*}
    \mathcal{S}^{(q,s)}, (q_0,\bot_2)\models \psi_{(q,s)} \text{\quad iff\quad }
    (q,s)\in\CPG(\mathcal{S}).
  \end{align*}
  We set 
  \begin{align*}
    &\mathcal{S}^{(q,s)}\coloneqq(Q', \Sigma, \Gamma', \Delta^{(q,s)},
    q_0) \text{ with}\\ 
    &Q'\coloneqq Q \cup \{q_{t}: t\leq s\} \cup \{q_{\emptyset}\},\text{ and}\\
    &\Gamma':=\Gamma\cup\{(q_{t}, \op): t\leq s, \op\in\Op\} \cup (\{q_{\emptyset}\}\times\Op),
  \end{align*} 
  where  $q_t$ is a new state for every substack $t$ of the stack s we
  are looking for and $q_\emptyset$ is used for 
  checking that certain operations can or cannot be performed on a
  configuration.  

  In the following we define $\Delta^{(q,s)}\supseteq \Delta$ by
  induction on the size of $s$ such that
  \mbox{$\Delta^{(q_t,t)}\subseteq\Delta^{(q,s)}$} for all proper substacks
  $t<s$. 
  \begin{enumerate}
  \item For $s=\bot_2$ we set 
    \begin{align*}
      \Delta^{(q,s)}\coloneqq\Delta\cup
      \big\{(q,\bot,(q_{\emptyset},\Clone{2}),q_{\emptyset},\Clone{2}), 
      (q,\bot, (q_\emptyset, \Pop{2}),q_\emptyset, \Pop{2})\big\}.
    \end{align*}
    Additionally, we set
    $\varphi_{(q_0,s)}\coloneqq\langle q_{\emptyset}, \Clone{2}\rangle
    \True  \wedge
    [q_\emptyset,\Pop{2}] \False$. Note that the first part of this
    formula is satisfied in $\mathcal{S}^{(q,s)}$ at some
    configuration $c$ if the state is $q$ and the topmost symbol is
    $\bot$. At such a configuration $c$, the second part can only be
    satisfied if no $\Pop{2}$ is possible, i.e., if the width of the
    stack is $1$. 
  \item Assume that $\lvert s\rvert>1$ and $\Sym(s)=\bot$ for some
    stack $s$. Then we set $t=\Pop{2}(s)$ and 
    \begin{align*}
      \Delta^{(q,s)}\coloneqq\Delta^{(q_t,t)}\cup
      \big\{(q,\bot,(q_t, \Pop{2}), q_{t},\Pop{2})\big\}       
    \end{align*}
    and $\varphi_{(q,s)}\coloneqq\langle q_{t},\Pop{2}\rangle \varphi_{(q_t,t)}$. 
  \item The next case is $\Sym(s)\neq \bot$ and $\Lnk(s)=0$.
    Then we set $t\coloneqq\Pop{1}(s)$ and
    \begin{align*}
      \Delta^{(q,s)}\coloneqq&\Delta^{(q_t,t)}\\
      &\cup \big \{(q,\Sym(s),(q_t,\Pop{1}),q_t,\Pop{1}), 
      (q,\Sym(s), (q_\emptyset,\Collapse), q_\emptyset, \Collapse)\big\} 
    \end{align*}
    and 
    \begin{align*}
      \varphi_{(q,s)}\coloneqq\langle q_t, \Pop{1} \rangle
      \varphi_{(q_t,t)} \wedge 
      [q_\emptyset,\Collapse]\False.
    \end{align*}
  \item In all other cases we set $t\coloneqq\Pop{1}(s)$ and
    $u\coloneqq\Collapse(s)$. 
    We set
    \begin{align*}
      \Delta^{(q,s)}\coloneqq&\Delta^{(q_t,t)} \cup
      \big \{(q,\Sym(s),(q_t, \Pop{1}), q_t,\Pop{1}) \big\} \\
      &\cup\big\{ 
      (q,\Sym(s), (q_U, \Collapse), q_u, \Collapse)\big\}
      \end{align*}
      and
      \begin{align*}
        \varphi_{(q,s)}\coloneqq&\langle q_t,\Pop{1}\rangle
        \varphi_{(q_t,t)} \wedge 
        \langle q_u,\Collapse \rangle
        \varphi_{(q_u,u)}.
    \end{align*} 
  \end{enumerate}
  We show by induction that for all $q\in Q$ and stacks $s$ 
  \begin{align*}
    \CPG(\mathcal{S}^{(q,s)}), c \models \varphi_{(q,s)} \text{ iff }
    c=(q,s).
  \end{align*}
  The initial stack $s=\bot_2=[\bot]$ is characterised by the facts that
  the top symbol of the stack is the bottom-of-stack symbol and that
  $\Pop{2}$ is undefined. The first conjunct of 
  \begin{align*}
    \varphi_{(q_0,s)}= \langle q_{\emptyset}, \Clone{2}\rangle
    \True  \wedge [q_\emptyset,\Pop{2}] \False    
  \end{align*}
  is only satisfied if
  the top symbol is $\bot$ and the second conjunct is satisfied
  if and only if $\Pop{2}$ is undefined. Thus, $\varphi_{(q,\bot_2)}$
  and $\mathcal{S}^{(q,\bot_2)}$ satisfy our claim. 
  
  For the induction step, note that $\Collapse$ is defined if and only
  if the collapse link of the topmost symbol is not $0$. If
  $\Collapse$ is defined 
  for some stack $s$ and $u\coloneqq\Collapse(s), t=\Pop{1}(s)$ then
  $\Delta^{(q_u,u)}\subseteq \Delta^{(q_t,t)}\subseteq \Delta^{(q_s,s)}$ because
  $u\leq t$. With
  these observations the induction step is straightforward by case
  distinction on the topmost symbol of $s$ and on the fact whether
  $\Collapse(s)$ is defined. 
  Let $(q,s)$ be some configuration and let $l$ be the 
  minimal level such
  that $\Pop{l}(s)$ is defined. 
  On the graph generated by  $\mathcal{S}^{(q,s)}$,
  the formula $\varphi_{(q,s)}$ asserts
  that this $\Pop{l}$ operation is defined and, by induction
  hypothesis, results in the stack $\Pop{l}(s)$. The analogous argument
  applies to the result of a collapse operation if the operation is
  defined on $s$. If it is undefined, i.e., $\Lnk(s)=0$ then the
  formula $\varphi_{(q,s)}$ asserts that the collapse is undefined. 

  Now, we set $\psi_{(q,s)}\coloneqq\mu Z.(\Diamond Z \vee
  \varphi_{(q,s)})$ which is just the formula asserting reachability
  of some point where $\varphi_{(q,s)}$ holds. Thus, 
  \begin{align*}
    &(q,s)\in\CPG(\mathcal{S})&
    &\text{iff}&
    &\CPG(\mathcal{S}^{(q,s)}), (q_0, \bot_2)\models \psi_{(q,s)} \enspace.&
  \end{align*}
  The latter problem is decidable due to Theorem
  \ref{CPGLmuModelCheckingComplexity}.  
\end{proof}
\begin{remark}
  This lemma extends to systems of higher level. 
  But in the case of higher levels, the proof needs some further
  preparation. The underlying problem that one faces on higher levels
  is the following. Consider the level $3$ stacks $s_2:=[[\bot
  (\sigma,2,0)]]$ and $s_3:=[[\bot (\sigma, 3, 0)]]$. For any sequence of
  stack operations, the result of the application of this sequence to
  $s_2$ is defined if and
  only if its application to $s_3$ is defined. Furthermore, the
  resulting stacks are identical except for the replacement of 
  level $2$ links of value $0$  by level $3$ links of value $0$. 

  Thus, our approach cannot distinguish between $s_2$ and $s_3$. 

  In order to make our approach work, we have to transform
  $\mathcal{S}$ into a new pushdown system $\mathcal{S'}$ over a new
  alphabet $\Sigma'$ which is \emph{level aware}. 
  Level awareness is defined as follows. There is a mapping
  $f:\Sigma' \to \{1, 2, 3, \dots, l\}$ such that for each stack
  generated by $\mathcal{S'}$, $\Sym(s)=\sigma$ implies that $\Lvl(s)
  = f(\sigma)$. This system can be obtained by replacing $\Sigma$ by
  $\Sigma\times\{1,2,3,\dots, l\}$ and by using $\Push{(\sigma,k),k}$
  instead of $\Push{\sigma,k}$. 

  Then we can apply the generalisation of the approach of level $2$ to
  this new system and solve the decision problem for configurations.
\end{remark}

We now turn to a quantitative version of the  decision
problem for configurations. We want to compute how many runs to a
given configuration 
exist up to a given threshold $k\in\N$.  In the next lemma we show
that this question can be reduced to the  decision
problem for configurations.

\begin{lemma}\label{LemmaCountNumberOfRuns}
  There is an algorithm solving the following problem. Given a level
  $2$ collapsible pushdown system
  $\mathcal{S}=(Q,\Sigma,\Gamma,q_0,\Delta)$, 
  a state $q\in Q$, a stack $s\in\Stacks_2(\Sigma)$ and a threshold
  $k\in\N$, how many runs from the initial configuration to $(q,s)$
  exist up to threshold $k$?
\end{lemma}
\begin{proof}
  First of all, by Lemma \ref{ReachableConfigsinMuCalcul}, it is
  decidable whether $(q,s)$ is a node of $\Graph(\mathcal{S})$. If
  this is the case then there is at least one run of the desired
  form. Otherwise there are $0$ runs of this form. 
  
  Assume that $(q,s)\in\Graph(\mathcal{S})$. 
  Since the runs of $\mathcal{S}$ are recursively enumerable, we can
  compute the length-lexicographically smallest run
  $\rho_1$ to $(q,s)$.\footnote{We assume that the
    transition relation of  $\mathcal{S}$ is a totally ordered set.}  

  In the following we show how to decide whether there is a second run
  of the desired form. For this purpose, let $l:=\length(\rho_1)$ and
  let $\delta_i$  be the transition between $\rho_1(i)$ and 
  $\rho_1(i+1)$ for all $0\leq i < l$. Furthermore, 
  let $q_i$ be the state at $\rho_1(i)$.
  Now, we construct a new pushdown system $\tilde{\mathcal{S}}:=(\tilde
  Q, \Sigma ,\Gamma,\tilde q_0, \tilde\Delta)$ where
  \begin{align*}
    \tilde Q:= Q\cup\{\tilde q_0, \tilde q_1, \dots, \tilde q_l\}  
  \end{align*}
   for new states $\tilde q_0 ,\dots, \tilde q_l$ and
   \begin{align*}
     \tilde\Delta:=&\Delta \cup  
     \left\{(\tilde q_i, \sigma, \gamma, \tilde q_{i+1}, \op): \delta_i = (q,
     \sigma, \gamma, q', \op)\right\} \\
     &\cup 
     \left\{ (\tilde q_i, \sigma, \gamma, q', \op): (q_i, \sigma, \gamma,
     q', \op)\in\Delta\setminus\{\delta_i\},
     \right\}.  
   \end{align*}
   This system copies the behaviour of every initial segment of
   $\rho_1$ and stays within the new states. As soon as it simulates
   one of the transitions of $\Delta$ that do not extend the run
   to another initial segment of $\rho_1$, it changes to the correct
   original state in $Q$. From this point on, the system behaves
   exactly like $\mathcal{S}$. Note that the run corresponding to
   $\rho_1$ in 
   $\tilde{\mathcal{S}}$ ends in configuration $(\tilde q_l, s)$ (for $s$
   the final stack of $\rho_1$). Hence, the corresponding run is no
   witness for the reachability of $(q,s)$ in the new system. Thus, if
   $\CPG(\tilde{\mathcal{S}})$ contains $(q,s)$, then there are two
   different runs in $\mathcal{S}$ from the initial configuration to
   $(q,s)$. 

   Repeating this construction up to $k$ times, we compute the
   runs to $(q,s)$ up to threshold $k$. 
\end{proof}
\begin{remark}
  We have no elementary bound on the complexity of this
  algorithm. This is due to the fact that we cannot derive a
  polynomial bound on the length of the run $\rho_1$. Hence, the size
  of the 
  pushdown system under consideration may increase too much in
  each iteration. Since we use the $L\mu$ model checking algorithm on
  each of the pushdown systems we construct, the resulting algorithm
  is doubly exponential in the size of the largest pushdown system
  that we construct. 
\end{remark}

In the last part of this section, we recall a lemma of
Blumensath from \cite{Blumensath2008} concerning the substitution of
prefixes of stacks. The original
lemma was stated for higher-order pushdown systems (without collapse)
of arbitrary level. Here, we only recall the result for level $2$
pushdown systems and we present a
straightforward adaption to the  case of level $2$ collapsible
pushdown systems. We start by defining a prefix relation on
stacks. Note that this relation does not coincide with the substack
relation.

\begin{definition} \label{Def:Prefixrelation}
  For some level $2$ stack $t$ and some substack $s\leq t$ we say that
  $s$
  is a \emph{prefix} 
  of $t$ and write 
  $s\prefixeq t$, if there are $n\leq m
  \in\N$ such that $s=w_1:w_2:\dots : w_{n-1}: w_n$ and 
  $t=w_1:w_2 \dots :w_{n-1} : v_n : v_{n+1} : \dots : v_m$ such that $w_n\leq
  v_j$ for all $n\leq j \leq m$. 

  For  some run $\rho$, we write $s\prefixeq \rho$ if
  $s\prefixeq\rho(i)$ for all $i\in\domain(\rho)$. 
\end{definition}
\begin{remark}
  Note that $s\prefixeq t$ obtains if
  $s$ and $t$ agree on the
  first $\lvert s\rvert -1$ words and the last word of $s$ is a prefix of all
  other words of $t$. Especially, 
  $s$ has to be a substack of  $t$ ands
  $\lvert s \rvert \leq \lvert t \rvert$.
\end{remark}

Now, we introduce a function that replaces the prefix of some
stack by some other.
\begin{definition}
  Let $s,t,u$ be level $2$ stacks such that $s\prefixeq t$. Assume that 
  \begin{align*}
   &s=w_1:w_2:\dots : w_{n-1}: w_n,\\
   &t=w_1:w_2 \dots :w_{n-1} : v_n : v_{n+1} : \dots : v_m, \text{
     and}\\
   &u=x_1: x_2: \dots : x_p
  \end{align*}
  for numbers $n,m,p\in\N$ such that $n\leq m$. 
  For each $n\leq i \leq m$, let $\hat v_i$ be the unique word such that
  $v_i=w_n\circ \hat v_i$. We define
  \begin{align*}
    t[s/u]:= x_1: x_2: \dots: x_{p-1} : (x_p\circ \hat v_n) :
    (x_p \circ \hat v_{n+1}): \dots : (x_p\circ \hat v_m)
  \end{align*}
  and call $t[s/u]$ 
\emph{the stack obtained from $t$ by replacing the prefix $s$ by $u$}. 
\end{definition}
\begin{remark}
  Note that for $t$ some stack with level $2$ links, the resulting
  object $t[s/u]$ may be no stack. 
  Take for example the stacks 
  \begin{align*}
    &t=\bot (a,2,0)
    : \bot (a,2,0),\\ &s=\bot (a,2,0) : \bot\text{ and}\\ &u=\bot:\bot.     
  \end{align*}
  Then $t[s/u]=\bot: \bot (a,2,0)$.
  This list of words
  cannot be created from the initial stack using the stack operation
  because an element $(a,2,0)$ in the second word has to be a clone of
  some element in the first one. But $(a,2,0)$ does not occur in the
  first word. 
  
  If $t\in\Sigma^{+2}$, i.e., if $t$ does not contain links of level
  $2$,  then $t[s/u]$ is always a stack. Thus, the prefix 
  replacement for stacks of higher-order pushdown systems always
  results in a well-defined stack while prefix replacement for stacks
  of collapsible pushdown systems may result in objects that are not stacks. 
\end{remark}

In the following we study the compatibility of prefix replacement with
the stack operations. 

\begin{lemma} \label{Lem:PrefixesAndOperations}
  Let $s,t$ be stacks such that $s\prefixeq t$. Let $\op$ be some
  operation. If $s\notprefixeq \op(t)$, then one of the following
  holds:
  \begin{enumerate}
  \item $\op(t) = \Pop{2}^k(s)$ for some $k\in\N$ or
  \item $\op(t) = \Pop{1}(s)$, $\TOP{2}(t)=\TOP{2}(s)$ and
    $\TOP{2}(\op(t))=\Pop{1}(\TOP{2}(s))$. 
  \end{enumerate}
\end{lemma}
\begin{proof}
  If $\op$ is $\Clone{2}$ or $\Push{\sigma,i}$ for some
  $\sigma\in\Sigma$ and $i\in\{1,2\}$, then $s\prefixeq t$ implies
  $s\prefixeq \op(t)$. 

  If $\op = \Pop{2}$ and $s\notprefixeq\op(t)$ then $\lvert t \rvert =
  \lvert s \rvert$ and $\op(t) = \Pop{2}(s)$. 

  If $\op = \Pop{1}$, $s\prefixeq t$ and $s\notprefixeq\op(t)$ implies
  that $\TOP{2}(s)\leq \TOP{2}(t)$ but $\TOP{2}(s)\nleq \TOP{2}(\op(t))$. 
  One immediately concludes that $\TOP{2}(t)=\TOP{2}(s)$ and
  $\TOP{2}(\op(t))=\Pop{1}(\TOP{2}(s))$.

  If $\op = \Collapse$, we have to distinguish two cases. If $\Lvl(t)
  = 1$, then we apply the same argument as in the case of $\op=\Pop{1}$. 
  Otherwise, $\op(t) = \Pop{2}^m(t)$ for some $m\in\N$ and one
  reasons analogously to the case of $\op=\Pop{2}$. 
\end{proof}

\begin{lemma} \label{Lem:PrefixPop1Pop2wieder}
  Let  $s,t$ be stacks such that $s\prefixeq t$ and  $\lvert t \rvert
  > \lvert s \rvert$. 
  Then it holds that $s\prefixeq\Pop{2}(t)$. 
\end{lemma}
\begin{proof}
  Just note that $t=\Pop{2}(s): t'$ for $t'$ a stack where each word
  is prefixed by $\TOP{2}(s)$. Furthermore, $\lvert  t \rvert > \lvert
  s \rvert$ implies that $\lvert t' \rvert \geq 2$. 
  Hence, $\Pop{2}(t) = \Pop{2}(s): t'$ for $t'$ a stack of width at
  least $1$ where each word is $\TOP{2}(s)$ prefixed by
  $\TOP{2}(s)$. Thus, $s\prefixeq 
  t$ holds. 
\end{proof}

Blumensath showed the following important compatibility of prefix
replacement and stack operations in the case of level $2$  pushdown
systems (without collapse!). 


\begin{lemma}[\cite{Blumensath2008}] \label{Lem:BlumensathHOLevel2}
  Let $\rho$ be a run of some pushdown system
  $\mathcal{S}$ of level $2$ and let $s,u\in \Sigma^{+2}$ be stacks such that
  the following conditions are satisfied:
  \begin{enumerate}
  \item $s\prefixeq \rho$,
  \item $\TOP{2}(s)<\TOP{2}(\rho(i))$ for all $i<\length(\rho)$ 
    or  $\Sym(u)=\Sym(s)$. 
  \end{enumerate}
   Under these conditions,  the function 
  $\rho[s/u]$ defined by $\rho[s/u](i):=\rho(i)[s/u]$ is a run of
  $\mathcal{S}$.  
\end{lemma}
\begin{proof}[Proof (sketch).]
  One proves this
  lemma by induction on $\domain(\rho)$. The transitions performed in
  $\rho$ can be carried over one by one to the
  transitions of $\rho[s/u]$. 
\end{proof}

Now, we present an adaption of this idea to collapsible pushdown
systems. 

\begin{lemma} \label{Lem:BlumensathLevel2}
  Let $\rho$ be a run of some collapsible pushdown system
  $\mathcal{S}$ of level $2$ and let $s$ and $u$ be stacks such that
  the following conditions are satisfied:
  \begin{enumerate}
  \item $s\prefixeq \rho$,
  \item $\TOP{2}(s)<\TOP{2}(\rho(i))$ for all $i<\length(\rho)$ 
    or $\TOP{1}(u)=\TOP{1}(s)$, 
  \item $\lvert s \rvert = \vert u \rvert$, and
  \item for $\rho(0)=(q,t)$, $t[s/u]$ is a stack. 
  \end{enumerate}
  Under these conditions  the function 
  $\rho[s/u]$ defined by $\rho[s/u](i):=\rho(i)[s/u]$ is a run of
  $\mathcal{S}$.  
\end{lemma}
\begin{proof}
  The proof is again by induction on $\domain(\rho)$. 
  For all operations, except for $\Collapse$, the proof of this lemma is
  analogous to the proof of the previous lemma. For each such operation $\op$
  occurring at position $i$ in
  $\rho$ one shows that $\rho(i+1)[s/u] = \op(\rho(i)[s/u])$. 
  
  For the collapse operation, assume that there is a position $i$ such that 
  \begin{align*}
    \rho(i+1) = \Collapse(\rho(i))  
  \end{align*}
  and such that $\rho(i)[s/u]$
  is defined. Due to condition 2, the topmost symbol and the collapse
  level of $\rho(i)$ and   $\rho(i)[s/u]$ agree. 
  Thus, if the collapse level is $1$, then the collapse acts on both
  configurations like a $\Pop{1}$. In this case, the compatibility of
  this $\Collapse$ 
  with the prefix replacement follows from the proof of the 
  case of $\Pop{1}$. 
  Otherwise, the collapse level of the topmost element of both stacks is
  $2$.
  In this case the collapse links of the topmost elements also agree by
  definition. Furthermore, due to $\lvert s \rvert = \lvert u \rvert$
  the width of $\rho(i)$ and $\rho(i)[s/u]$ agrees. Hence, there is
  some $k\in\N$ such that the collapse applied to both configurations
  results in $\Pop{2}^k(\rho(i))$ and $\Pop{2}^k(\rho(i))[s/u]$,
  respectively. Thus,  the reduction to the iterated use of the
  case of $\Pop{2}$ proves the claim. 
\end{proof}


\section{Technical Results on the Structure of Collapsible Pushdown
  Graphs} 
\label{ChapterMilestons}
\label{ChapterLoops}
\label{ChapterReturns}

In this section, we develop the technical background for
our main results that are presented in  
Sections \ref{CPG-Tree-Automatic} and \ref{Chapter HONPT}. 

As in the end of the previous section, this section is only concerned
with collapsible pushdown systems of level $2$. Hence, if we write
collapsible pushdown system, we always mean one of level $2$. 

The overall goal of this section is the following:
finite automata can be used  to determine how many\footnote{For the
  rest of this section, the question ``how many?'' is meant up to a
  certain threshold $k\in\N$, i.e., ``how many runs to $(q,s)$ exist''
stands for ``given a threshold $k\in\N$, how many runs to $(q,s)$
exist up to threshold $k$?''.}
runs from the initial configuration to some configuration $(q,s)$
exist. 
In order to prove this result, we introduce three notions:
\emph{returns}, \emph{loops}, and 
\emph{generalised milestones}\footnote{The term ``generalised'' refers
  to the fact that this notion is a generalisation of the notion
  ``milestone'' which we introduced in \cite{Kartzow10}.}. We motivate
these notions from the last to the first.

Let $s$ and $s'$ be stacks. 
We call $s'$ a generalised milestone of $s$ if every run
from the 
initial configuration to a configuration with stack $s$ has to pass 
$s'$ at some intermediate step. Thus, it
follows directly from this definition that the reachability of a
certain stack from the initial configuration decomposes into the
analysis of the reachability of milestones from other milestones
of this stack. We will see that every run to $s$ passes all the
milestones of 
$s'$ in a certain order. Thus, the question ``how many runs to $s$
exist?'' can be reduced to the question ``how many runs from one
milestone of $s$ to the next exist?''. 

A closer analysis of this decomposition shows that the run from one
milestone to the next is always a \emph{loop} followed by exactly one
transition. A \emph{loop} is a run from some configuration $(q,s)$ to
some configuration $(q',s)$ not passing a substack of
$\Pop{2}(s)$. This means that a run starts and ends with the same
stack $s$ 
and it does not ``look into'' the content of $\Pop{2}(s)$. 

Using this result, the question ``how many runs to $(q,s)$ exist?''
can be reduced to the question ``how many loops of each generalised
milestone of $s$ exist?''. 

In order to show that a finite automaton can answer the last question,
we introduce the notion of a \emph{return}. A run $\rho$ is called
return if it is a run from some
stack $s$ to the stack $\Pop{2}(s)$ that satisfies the following
conditions:
\begin{enumerate}
\item before the last position, no substack of $\Pop{2}(s)$ is passed, and
\item the collapse links of level $2$ stored in $\TOP{2}(s)$ are not
  used by $\rho$. 
\end{enumerate}
It turns out that returns naturally appear as subruns of loops. 
In the following we first introduce generalised milestones and develop
their theory. Then we define loops and returns  and show their
connection to generalised milestones in Section
\ref{sec:LoopsAndReturns}.
In Section \ref{sec:CompReturns} we develop the theory of counting
returns. Finally, we develop the analogous theory of loops in Section
\ref{sec:CompLoops}. 

\subsection{Milestones and Loops}
\label{subsec:Milestones}
Recall that $w\sqcap v$ denotes the greatest common prefix of the
words $w$ and $v$ (cf. Section \ref{Sec:wordsAndTrees}). We start with
a formal definition of generalised milestones. Afterwards, we show that
this definition fits the informal description given before. 

\begin{definition}
  Let $s=w_1:w_2:\dots :w_k$ be a stack. We call a stack $m$ a
  \emph{generalised milestone of} $s$ if $m$ is  of the form 
  \begin{align*}
    &m=w_1:w_2:\dots: w_i:v_{i+1} \text{ where }
    0\leq i<k,\\
    &w_i \sqcap w_{i+1}\leq v_{i+1}\text{ and }\\
    &v_{i+1}\leq w_i\text{ or }v_{i+1}\leq w_{i+1}.
  \end{align*}
  We denote by $\genMilestones(s)$ the set of all generalised
  milestones of $s$.  
  
  For a generalised milestone $m$ of $s$, we call $m$ a
  \emph{milestone of }$s$ if  $m$ is a substack of $s$. 
  We write $\Milestones(s)$ for the set of all milestones of
  $s$.
\end{definition}
\begin{remark}
  In the following we are mainly concerned with generalised
  milestones. Only in Section \ref{CPG-Tree-Automatic} the concept of
  milestones appears as a useful concept on its own.
\end{remark}

A simple observation is that  we
can derive a bound on the number of generalised milestones from the
height and the width of a stack. 

\begin{lemma} \label{LemmaNumberMilestones}
  For each stack $s$ there are less than $2 \cdot \height(s) \cdot
  \lvert s\rvert$ many generalised milestones.
\end{lemma}

In our informal description of generalised milestones, we said that
the generalised milestones of $s$ are those  stacks that every run to
$s$ has to pass.
In order to show this, we use a result of Carayol \cite{Carayol05}.
He showed the following. For each
higher-order pushdown stack $s$
there is a 
unique minimal sequence of stack operations that creates $s$ from the
initial stack. 
On level two, this sequence creates the stack
word by word, i.e., it starts with a sequence of push operations
writing the first word onto the stack, then there is a clone
operation, after this there is a sequence of $\Pop{1}$ transitions
followed by a sequence of push transitions that create the second word
of the stack, then there follows a clone and so on. Furthermore, the
topmost word reached after the $n$-th of the $\Pop{1}$ sequences is exactly
the greatest common prefix of the $n$-th and the $(n-1)$-st word of
the stack. 
This result directly carries over to collapsible
pushdown stacks due to the following fact: 
on level two the result of a collapse operation is either the same as
applying a $\Pop{1}$ or a sequence of $\Pop{2}$ operations. 
Carayol's result shows that for any sequence containing a $\Pop{2}$
operation, there is a shorter one where this $\Pop{2}$ is
eliminated. Hence, we can eliminate in the same way any $\Collapse$
of link level $2$. Finally, any other $\Collapse$ can be treated like a
$\Pop{1}$ operation. 
We describe Carayol's result more formally in the following lemma.
\begin{lemma}[\cite{Carayol05}] \label{Lemma:Carayol05Milestones}
  \begin{itemize}
  \item[]
  \item For each collapsible pushdown stack $s$ of level $2$ there is
    a minimal sequence of operations $\op_1, \op_2, \dots, \op_n\in\{
    \Push{\sigma,i}, \Pop{1}, \Clone{2}\}$ such that 
    $s=\op_n(\op_{n_1}(\dots(\op_1(\bot_2))))$.
  \item   For $\op_1, \op_2,\dots, \op_n$ the minimal sequence
    generating a stack $s$, the stack
    \begin{align*}
    \op_j (\op_{j-1}(\dots\op_0(\bot_2)))      
    \end{align*}
    is a
    generalised milestone of $s$ for each $0\leq j\leq n$.
    
    Furthermore, for each generalised milestone $m$ of $s$ there is a
    $0\leq j \leq n$ such that 
    $m= \op_j (\op_{j-1}(\dots\op_0(\bot_2)))$.
  \item Every run $\rho$ to some stack $s$ passes all generalised
    milestones of $s$. 
  \end{itemize}
\end{lemma}
\begin{remark}
  From the minimality of the sequence 
  $\op_1, \op_2, \dots, \op_n$ generating $s$ it follows that
  there is a bijection between the initial subsequences 
  $\op_1, \op_2, \dots, \op_j$ and the milestones of $s$. 
  From now on, we call
  $\op_j (\op_{j-1}(\dots\op_0(\bot_2)))$ the $j$-th
  milestone of $s$. 

  Note that if $i\leq j$ then the $i$-th milestone  $m_i$ of $s$ is a
  milestone of the $j$-th milestone $m_j$ of $s$.  
  If we restrict this order to the set of milestones
  $\Milestones(s)$, then it
  coincides with the substack relation. 
\end{remark}

We want to conclude the analysis of generalised milestones with a
lemma that characterises
runs connecting generalised milestones in terms of loops. 
Thus, we first give a precise definition of loops. Then we prove this
characterisation. 
A loop
is a run that starts and ends in the same stack and which satisfies
certain restrictions concerning the substacks that are passed. 

\begin{definition} \label{DefLoop}
  A \emph{loop}
  from $(q,s)$ to $(q',s)$ is a run $\lambda$ that does not pass a
  substack of  $\Pop{2}(s)$ and that may pass $\Pop{1}^k(s)$ only if
  the $k$  topmost elements of $\TOP{2}(s)$ are letters with links of
  level $1$. This means that for all $i\in\domain(\lambda)$, if
  $\lambda(i)= (q_i,\Pop{1}^k(s))$ then $\Lvl(\Pop{1}^{k'}(s))=1$ for
  all $0\leq k' < k$. 

  If $\lambda$ is a loop from $(q,s)$ to $(q',s)$ such that
  $\lambda(1)=\Pop{1}(s)$ and
  $\lambda(\length(\lambda)-1)=\Pop{1}(s)$, then we call 
  $\lambda$ a \emph{low loop}. 

  If $\lambda$ is a loop from $(q,s)$ to  $(q',s)$ that never passes
  $\Pop{1}(s)$, then we call $\lambda$ a \emph{high loop}. 
\end{definition}
\begin{remark}
  If $\lambda$ is a loop from $(q,s)$ to $(q',s)$ such that the stack
  at $i$ and at $j$ is $s$ for  $i\leq j\in\domain(\lambda)$, then
  $\lambda{\restriction}_{[i,j]}$ is a loop. 
\end{remark}

We now characterise runs connecting milestones in terms of loops. 

\begin{lemma} \label{LemmaDecompositioninLoops}
  Let $\rho$ be a run from the initial configuration to the stack
  $s=w_1:w_2:\dots :w_k$. Furthermore, let $n$ be the number of
  generalised milestones of $s$. 
  For all $i\leq n$, let $m_i$ be the $i$-th generalised milestone of
  $s$.
  Furthermore, let  $n_i$ denote the maximal position such that the
  stack of 
  $\rho(n_i)$ is $m_i$. We write $q_i$ for the state of $\rho(n_i)$,
  i.e., $\rho(n_i)=(q_i,m_i)$. 
  For all $i<n$, there is some state $q'_{i+1}$ such that there is a
  transition  from $\rho(n_i)$ to $(q_{i+1}', m_{i+1})=\rho(n_i+1)$
  and $\rho{\restriction}_{[n_i+1, n_{i+1}]}$ is a loop of
  $m_{i+1}$. Furthermore, $\rho{\restriction}_{[0,n_1]}$ is a loop of
  $\bot_2$. 
\end{lemma}
\begin{proof}
  Fix some $i\in\N$. We prove the claim for $m_i$ and $m_{i+1}$. 
  We distinguish the following cases.
  \begin{itemize}
  \item Assume that $m_{i+1} = \Clone{2}(m_i)$. In this case 
    \mbox{$m_i = w_1 : w_2 : \dots : w_{\lvert m_i \rvert}$}. Thus, at the last
    position $j\in\domain(\rho)$ where $\lvert \rho(j) \rvert =\lvert
    m_i \rvert$, the stack at 
    $\rho(j)$ is $m_i$ (because $\rho$ never changes the first $\lvert
    m_i\rvert$ many words after passing $\rho(j)$). Hence, $j=n_i$ by
    definition. Since 
    $\lvert s \rvert > \lvert m_i\rvert$, it follows directly that the
    operation at 
    $n_i$ is a $\Clone{2}$ leading to $m_{i+1}$. Note that $\rho$ never
    passes  a stack of 
    width $\lvert m_i \rvert$ again. Thus, it follows from Lemma 
    \ref{Lem:Howtwolinksevolve} that
    $\rho{\restriction}_{[n_i+1,n_{i+1}]}$ satisfies the restriction
    that it never 
    visits $\Pop{1}^k(m_{i+1})$ if
    $\Lvl(\Pop{1}^{k-1}(m_{i+1}))=2$. Thus, we conclude that this
    restriction is a loop. 
  \item Assume that $m_{i+1} = \Pop{1}(m_i)$. In this case, 
    $m_i=w_1: w_2: \dots: w_{\lvert m_i\rvert -1}: w$ for some $w$ such that
    $w_{\lvert m_i\rvert-1} \sqcap w_{\lvert m_i\rvert} < w \leq
    w_{\lvert m_i \rvert -1}$. 
    Thus, $w\not\leq w_{\lvert m_i\rvert}$ and creating $w_{\lvert m_i\rvert}$
    as the $\lvert m_i\rvert$-th word 
    on the stack requires passing $w_1: w_2: \dots :w_{\lvert m_i \rvert-1}:
    w_{\lvert m_i\rvert-1} \sqcap  w_{\lvert m_i\rvert}$. This is only possible
    via applying $\Pop{1}$ or $\Collapse$ of level $1$ to $m_i$. Since
    we assumed $n_i$ 
    to be maximal, the operation at $n_i$ must be $\Pop{1}$ or
    $\Collapse$ of level $1$ and leads to $m_{i+1}$. 

    We still have to show that $\rho{\restriction}_{[n_i+1,n_{i+1}]}$
    is a loop. 
    By definition of $n_{i+1}$, $\rho{\restriction}_{[n_i+1,n_{i+1}]}$
    starts and ends in $m_{i+1}$. 
    By maximality of $n_i$, $\rho{\restriction}_{[n_i+1,n_{i+1}]}$
    does not visit the stack \mbox{$\Pop{2}(m_i) = \Pop{2}(m_{i+1})$.}
    Furthermore, 
    note that $\TOP{1}(m_{i+1})$ is a cloned element. Hence, 
    Lemma 
    \ref{Lem:Howtwolinksevolve} implies that
    $\rho{\restriction}_{[n_i+1,n_{i+1}]}$ may only visit
    $\Pop{1}^k(m_{i+1})$ in case that
    $\Lvl(\Pop{1}^{k-1}(m_{i+1}))=1$.
    Thus, $\rho{\restriction}_{[n_i+1,n_{i+1}]}$ is a loop.  
  \item The last case is $m_{i+1} = \Push{\sigma,l}(m_i)$ for
    $(\sigma,l)\in\Sigma\times\{1,2\}$. In this
    case, 
    \begin{align*}
      m_i = w_1: w_2: \dots: w_{\lvert m_i\rvert-1}: w      
    \end{align*}
    for some $w$
    such that 
    $w_{\lvert m_i\rvert-1} \sqcap w_{\lvert m_i\rvert} \leq w <
    w_{\lvert m_i\rvert}$. 
    Creating $w_{\lvert m_i\rvert}$ on the stack requires pushing the
    missing symbols onto the stack as they cannot be obtained via
    clone operation  from the previous word. Since $n_i$ is maximal,
    the operation at $n_i$ is some $\Push{\sigma,l}$ leading to
    $m_{i+1}$.
    $\rho{\restriction}_{[n_i+1,n_{i+1}]}$ is a high loop due
    to the maximality of $n_i$ (this part of $\rho$ never visits
    $m_i=\Pop{1}(m_{i+1})$ or any other proper substack of $m_{i+1}$).\qedhere
  \end{itemize}
\end{proof}

We conclude this section by rephrasing  this result in terms of
milestones. We will use it in this form in Chapter
\ref{CPG-Tree-Automatic}.

\begin{corollary}\label{ClonePopMilestoneRun}
  Let $\rho$ be a run from the initial configuration to the
  configuration $(q,s)$ where $s$ decomposes as
  \begin{align*}
    s=w_1:w_2:\dots :w_k.
  \end{align*}
  Let $n_i$ denote the maximal position such that $\rho(n_i)=(q,m_i)$ for
  some $q\in Q$ and $m_i$ the $i$-th  milestone of $s$. 
  We define $q_i\in Q$ such that $\rho(n_i)=(q_i,m_i)$. 
  Then one of the following applies.
  \begin{enumerate}
  \item There is a $\Push{\sigma,j}$ transition
    from $\rho(n_i)=(q_i,m_i)$ to $(q_{i+1}', m_{i+1}):=\rho(n_i+1)$
    and $\rho{\restriction}_{[n_i+1, n_{i+1}]}$ is a loop of
    $m_{i+1}$, or
  \item there is a $\Clone{2}$ transition followed by a sequence 
    $\lambda_0\circ \pi_1 \circ \lambda_1 \dots \circ \pi_n\circ
    \lambda_n$ where the $\lambda_i$ 
    are loops and the $\pi_i$ are runs that perform exactly one
    $\Pop{1}$ operation or collapse of level $1$ each. 
  \end{enumerate}
  Furthermore, we have
  \begin{enumerate}
  \item[3.] $\rho(n_1)=(q_1, [\bot])$, i.e.,
    $\rho{\restriction}_{[0,n_1]}$ is a loop of $[\bot]$. 
    If $m$ is the number of milestones of $s$, then $\rho(n_m)=(q,s)$
    is the final configuration of $\rho$. 
  \end{enumerate}
\end{corollary}

As another direct corollary of the lemma, we obtain that the linear
order of the 
milestones induced by the substack relation coincides with the order
in which the milestones appear for the last time in a given run.

\begin{corollary} \label{Cor:OrderEmbedding}
  For an arbitrary run $\rho$ from the initial configuration to some
  stack $s$, the
  function
  \begin{align*}
    f:\Milestones(s) &\rightarrow \domain(\rho)\\  
    s'&\mapsto
    \max\{i\in\domain(\rho): \rho(i)=(q,s')\text{ for some }q\in Q\}    
  \end{align*}
  is an order embedding.
\end{corollary}

We have seen that generalised milestones induce a uniform
decomposition of all runs to a given stack. Furthermore, the parts of
the run that connect generalised milestones always consist of a loop
plus one further transition. 
In order to understand the existence of runs to certain
configurations,  we investigate the theory of loops in the following. 

\subsection{Loops and Returns}
\label{sec:LoopsAndReturns}
Recall that we have already defined loops in Definition \ref{DefLoop}.
Next, we define \emph{returns} which are runs from a stack $s:w$ to $s$
without visiting substacks of $s$. 
Our interest in returns stems from the fact that they
appear as subruns of high loops whence they play an important role in
finding loops for a given stack.

\begin{definition}  \label{DefReturn}
  Let $t=s:w$ be some stack with topmost word $w$. A \emph{return from
    $t$ to $s$} is 
  a run $\rho$ from $t$ to $s$ such that $\rho$ never visits a
  substack of $s$ except for the last stack of $\rho$ and
  such that one of the following holds:
  \begin{enumerate}
  \item the last operation in $\rho$ is $\Pop{2}$,
  \item the last operation in $\rho$ is a $\Collapse$ and 
    $w < \TOP{2}(\rho(\length(\rho)-1))$, i.e., $\rho$ pushes at first some new
    letters onto $t$ and then performs a collapse of one of these new
    letters, or 
  \item  there is some
    $i\in\domain(\rho)$ such that
    $\rho{\restriction}_{[i,\length(\rho)]}$ is a return from
    $\Pop{1}(t)$ to $s$. 
  \end{enumerate}
\end{definition}
\begin{remark}
  A return from
    $t$ to $\Pop{2}(t)$ is 
  a run $\rho$ from $t$ to $\Pop{2}(t)$ such that $\rho$ never visits a
  substack of $\Pop{2}(t)$ except for the last stack of $\rho$ and that does
  not use the level $2$ links stored in $\TOP{2}(t)$. 
\end{remark}

We first give an example for this definition, afterwards we discuss
its motivation. 

\begin{example} \label{Example:ReturnExample}
  Consider a collapsible pushdown system $\mathcal{S}$ over the
  alphabet $\{\bot, \top, a, b\}$ with the
  transitions 
  $(q_0, a , \gamma_0,q_1,\Clone{2}), 
  (q_1,  a, \gamma_1, q_1, \Collapse)$ and 
  $(q_1, b, \gamma_2, q_1,  \Push{a,2})$. 
  Consider the stack 
  \begin{align*}
    s:=\bot(b,2,0)\Box:
    \bot(b,2,1)a.     
  \end{align*}
  The transitions induce a unique run $\rho$ from $(q_1, s)$ to 
  $(q_1,\Pop{2}(s))$ of length $3$. $\rho$ is depicted on
  the left side of
  Figure \ref{fig:FirstReturnExample}.
  \begin{figure}
    \centering
    \begin{xy}
      \xymatrix@R=0pt@C=5pt{ 
        &\Box & a\\
        &(b,2,0) & (b,2,1) \\
        q_1, &\bot \ar[dddd]^{\gamma_1}&\bot\\
        \\
        \\
        \\
        &\Box &\\
        &(b,2,0) & (b,2,1) \\
        q_1, &\bot \ar[dddd]^{\gamma_2}&\bot\\
        \\
        \\
        \\
        &\Box & (a,2,1)\\
        &(b,2,0) & (b,2,1) \\
        q_1, &\bot \ar[dddd]^{\gamma_1}&\bot\\
        \\
        \\
        \\
        &\Box &\\
        &(b,2,0) & \\
        q_1, &\bot &
        }\hskip 10cm
        \xymatrix@R=0pt@C=5pt{ 
          &      & a       & a\\
        &(b,2,0) & (b,2,0) & (b,2,0) \\
        q_1, &\bot \ar[dddd]^{\gamma_1}&\bot&\bot \\
        \\
        \\
        \\
        &         & a\\
        &(b,2,0) & (b,2,0) & (b,2,0)\\
        q_1, &\bot \ar[dddd]^{\gamma_2}&\bot&\bot\\
        \\
        \\
        \\
        &        &a       &(a,2,2)\\
        &(b,2,0) & (b,2,0)&(b,2,0) \\
        q_1, &\bot \ar[dddd]^{\gamma_1}&\bot&\bot\\
        \\
        \\
        \\
        &        & a \\
        &(b,2,0) & (b,2,0)\\
        q_1, &\bot &\bot
        }
      \end{xy}
    \caption{The run $\rho$ from $(q_1,s)$ to
      $(q_1,\Pop{2}(s))$ on the left side and
      the run $\rho'=\lambda{\restriction}_{[1,4]}$ from $(q_1,s')$ to
      $(q_1, \Pop{2}(s'))$ on the 
      right side.}
    \label{fig:FirstReturnExample}
  \end{figure}
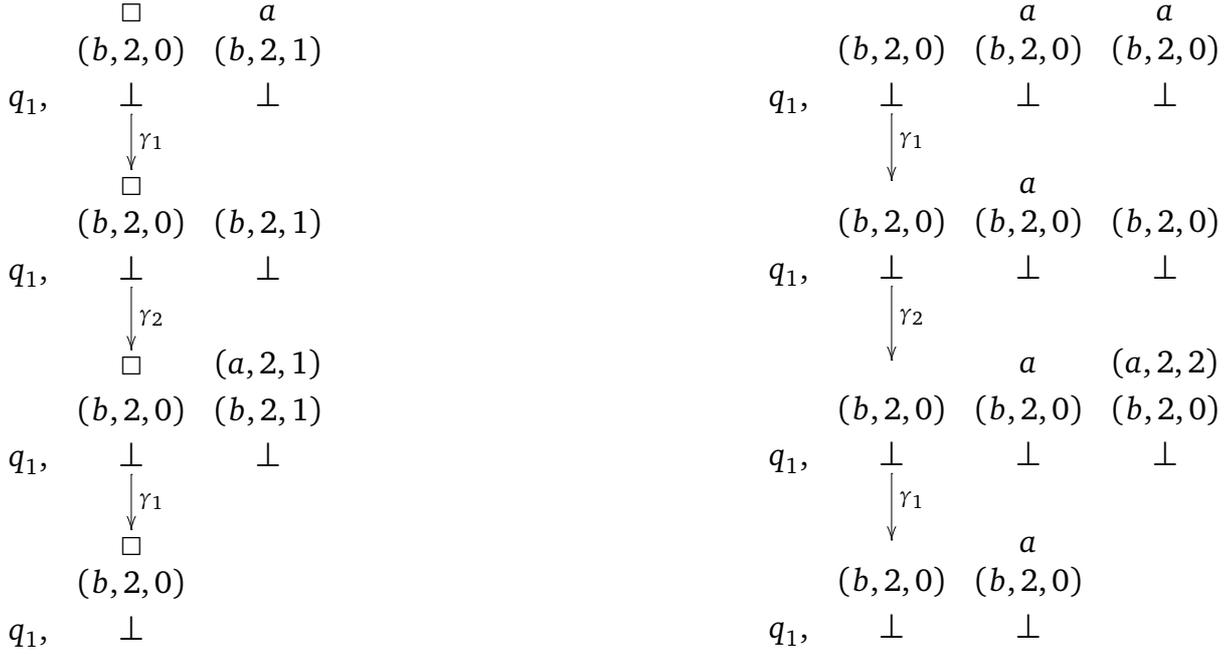
  $\rho{\restriction}_{[1,3]}$ is a return from 
  \begin{align*}
    &(q_1, \bot(b,2,0)\Box: \bot(b,2,1))  \text{ to}\\
    &(q_1, [\bot(b,2,0)\Box])
  \end{align*}
  because it satisfies the
  second item of the definition of a return. 
  Hence, $\rho$ is a return because it satisfies the third item
  of the definition.   
  
  We want to consider a second example that shows how returns occur
  as subruns of loops. 
  The transitions induce a loop $\lambda$ from 
  \begin{align*}
    &(q_0, \bot (b,2,0): \bot (b,2,0) a)\text{ to}\\  
    &(q_1, \bot    (b,2,0): \bot (b,2,0) a).     
  \end{align*} 
  The run passes 
  $(q_1, \bot (b,2,0): \bot (b,2,0) a: \bot (b,2,0) a)$ and
  continues from there as depicted on the right side of   
  Figure \ref{fig:FirstReturnExample} (the figure shows $\lambda$
  without its first configuration because this final part of
  $\lambda$ 
  plays a role in the next remark). 
  Note that $\lambda{\restriction}_{[2,4]}$ is a return starting from
  a stack with topmost word $\Pop{1}(\TOP{2}(\lambda(0)))$. 
  Later, when we analyse loops in detail we will see that this is a typical
  occurrence of a return. Any loop of a stack $s$ decomposes into
  parts prefixed by $s$ and parts that are returns of stacks with
  topmost word $\Pop{1}(\TOP{2}(s))$. 
\end{example}

\begin{remark}
  A return is a run from
  some stack $s$ to $\Pop{2}(s)$ that depends on the symbols and
  link levels of $\TOP{2}(s)$, but not on any other content of $s$
  in the following sense.  
  A return from $s$ to $\Pop{2}(s)$ consists of a sequence of
  transitions. For any stack $s'$ with $\lvert s' \rvert \geq
  2$ such that the topmost words of $s$ and $s'$ coincide on their
  symbols and link levels, this sequence can be applied to $s'$. The
  resulting run induced by this sequence is then a return from $s'$ to
  $\Pop{2}(s')$. 
 
  We explain this idea with some examples. 
  Let $\mathcal{S}$ be the pushdown system
  and $\rho$ the return from $s$ to $\Pop{2}(s)$ as in Example \ref{Example:ReturnExample}. 
  Consider the stack 
  \begin{align*}
    s':=\bot(b,2,0): \bot (b,2,0)a:
    \bot (b,2,0) a.    
  \end{align*}
  Note that the symbols and link levels of $\TOP{2}(s)$ and
  $\TOP{2}(s')$ agree  while their links  differ. 
  There is a return $\rho'$ from $(q_1,s')$ to 
  $(q_1, \Pop{2}(s'))$ which is obtained by
  starting in $(q_1,s')$ and  copying the transitions of $\rho$ one by
  one. The resulting return $\rho'$ is depicted on the right side of 
  Figure \ref{fig:FirstReturnExample}.
  
  This is not by accident, but by intention: whenever two stacks $s$
  and $s'$ coincide on the symbols and link levels of their topmost
  words, we can copy a return from $s$ to $\Pop{2}(s)$ transition by
  transition and obtain a return from $s'$ to $\Pop{2}(s')$. 
  This is due to two facts. 

  Firstly, a return from $s$ to $\Pop{2}(s)$
  never looks into $\Pop{2}(s)$ before its last configuration. Thus,
  the words below the topmost word have no influence on this run. 
  Secondly, the restriction of the use of collapse links ensures that a
  return  only
  uses  collapse links of level $2$ if these were created during the
  run $\rho$. If such a link points to $\Pop{2}(s)$, it is created by
  a push operation at some position $i$ in $\rho$ on a stack of width
  $\lvert s \rvert$. But then $\rho'$, the one to one copy of the
  transitions of $\rho$ with starting stack $s'$, uses a push transition
  at position $i$ on a stack of width $\lvert s' \rvert$. Thus, the link
  created in this step points to $\Pop{2}(s')$. 
  Hence, if $\rho$ uses the created collapse link and collapses the stack to
  $\Pop{2}(s)$, then $\rho'$ uses the copy of this link and collapses to
  $\Pop{2}(s')$. 
  
  We defined returns in such a way that they are  runs from some stack
  $s$ to $\Pop{2}(s)$ 
  that are independent of the links and the words below the topmost
  one. 
  The next example shows that the restricted use of the collapse
  operation  in the definition of returns is crucial for this property. 
  We present a  run from some stack $\hat s$ to $\Pop{2}(\hat s)$ that
  does not 
  look into the substacks of $\Pop{2}(\hat s)$ before the final position
  but that lacks the independence of the level $2$ links of the
  topmost word. 
  
  Consider the stacks
  \begin{align*}
    \hat s:=\bot (b,2,0)(b,2,0):\bot(a,2,1)a \text{ and}\\
    \hat s':=\bot(b,2,0):\bot(a,2,1):\bot(a,2,1)a. 
  \end{align*}
  We still consider the transitions given in Example 
  \ref{Example:ReturnExample}. Using these transitions, there are 
  runs $\hat\rho$ from $(q_1,  \hat s)$ to $(q_1, 
  \Pop{2}(\hat s))$ and $\hat\rho'$ from $(q_1, \hat s')$ to $(q_1,
  [\bot(b,2,0)])$ as depicted in Figure 
  \ref{fig:SecondReturnExample}.
  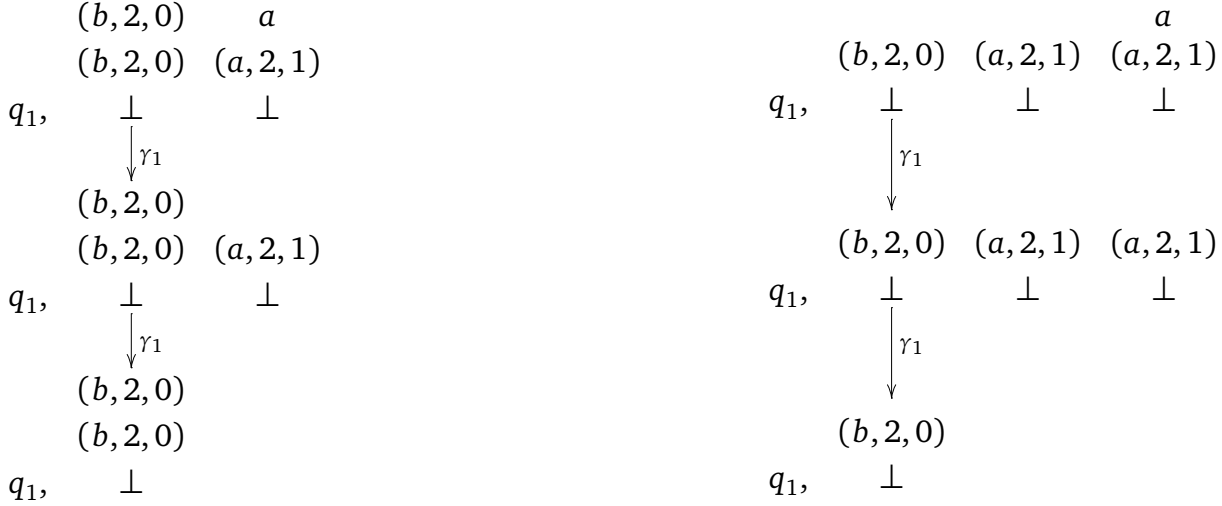
\begin{figure}
    \centering
    \begin{xy}
      \xymatrix@R=0pt@C=5pt{ 
        &(b,2,0) & a\\
        &(b,2,0) & (a,2,1) \\
        q_1, &\bot \ar[dddd]^{\gamma_1}&\bot\\
        \\
        \\
        \\
        &(b,2,0) &\\
        &(b,2,0) & (a,2,1) \\
        q_1, &\bot \ar[dddd]^{\gamma_1}&\bot\\
        \\
        \\
        \\
        &(b,2,0)\\
        &(b,2,0) & \\
        q_1, &\bot &
        }\hskip 10cm
        \xymatrix@R=0pt@C=5pt{ 
        & & & a\\
        &(b,2,0) & (a,2,1) & (a,2,1) \\
        q_1, &\bot \ar[dddddd]^{\gamma_1}&\bot&\bot \\
        \\
        \\
        \\
        \\
        \\
        & \\
        &(b,2,0) & (a,2,1) & (a,2,1)\\
        q_1, &\bot \ar[dddddd]^{\gamma_1}&\bot&\bot\\
        \\
        \\
        \\
        \\
        \\
        &  &\\
        &(b,2,0)  \\
        q_1, &\bot
        }
      \end{xy}
    \caption{The run $\hat \rho$ from $(q_1,\hat s)$ to
      $(q_1,\Pop{2}(\hat s))$ on the left side and
      the run $\hat\rho'$ from $(q_1,\hat s')$ to $(q_1, \Pop{2}(\Pop{2}(\hat
      s')))$ on the right side.}
    \label{fig:SecondReturnExample}
  \end{figure}

  Note that $\hat\rho$ is no return because it uses the level $2$
  collapse link stored in $\TOP{2}(\hat s)$. 
  Furthermore, $\TOP{2}(\hat s)= \TOP{2}(\hat s')$ and
  $\hat\rho'$ copies $\hat\rho$ transition by transition. 
  Nevertheless, $\hat\rho'$ does not end with the stack $\Pop{2}(\hat
  s')$ but with $\Pop{2}^2(\hat s')$. 
  
  Thus, if we drop the restriction on the use of collapse links, then
  we obtain runs from some stack $s$ to $\Pop{2}(s)$ that cannot be
  transferred into runs from stacks $s'$ to $\Pop{2}(s')$ even though
  $\TOP{2}(s)=\TOP{2}(s')$.  
\end{remark}

\subsection{Computing Returns} 
\label{sec:CompReturns}
As already mentioned in the previous section, the theory of returns is
important for the theory of loops. Thus, we first study
the theory of returns on its own. Later we apply this theory to the
theory of loops.  
Our main goal in this part is to provide  a finite
automaton that calculates on input $\TOP{2}(s)$ the number of returns
from $(q,s)$ to $(q',\Pop{2}(s))$ up to a given threshold $k\in\N$. 
We start by introducing appropriate  notation for this purpose.

\begin{definition} \label{Def:ReturnFunc}
  Let $\mathcal{S}$ be a collapsible pushdown system of level $2$. 
  We set
  \begin{align*}
    &\ReturnFunc{k}_{\mathcal{S}}(s):Q\times Q \rightarrow
    \{0,1,\dots,k\} \\ 
    &(q,q')\mapsto
    \begin{cases}
      i &\text{if there are exactly }i\leq k\text{ different returns of
        $\mathcal{S}$ 
        from }(q,s)\text{ to } (q',\Pop{2}(s))\\
      k & \text{otherwise.}
    \end{cases}
  \end{align*}
\end{definition}
\begin{remark}
  This function maps $(q,q')$ to the number $i$ of returns from $(q,s)$ to
  $(q',s)$ if $i\leq k$ and it maps $(q,q')$ to $k$ otherwise. In
  this sense $k$ stands for the class of at least $k$ returns. 
  Thus, the answer to the question ``how many returns from $(q,s)$ to
  $(q',\Pop{2}(s))$ exist up to threshold $k$?'' is exactly the value
  of $\ReturnFunc{k}_{\mathcal{S}}(s)(q,q')$. 
  If $\mathcal{S}$ is clear from the context, we will omit it and
  write $\ReturnFunc{k}$ instead of 
  $\ReturnFunc{k}_{\mathcal{S}}$. 
\end{remark}

As already indicated in the examples, it turns out that we can copy
returns between stacks which agree on their topmost words.
Lemma \ref{ReturnDependTopword} proves this fact. A corollary of
this lemma is that the number of returns from $(q,s)$ to
$(q',\Pop{2}(s))$ only depend on the topmost word of $s$. 
Hence, the following definition is well-defined. 
\begin{definition} \label{Def:returnsofWord}
  For $w$ an arbitrary word, let $\ReturnFunc{k}(w)$ be 
  $\ReturnFunc{k}(s)$ for an arbitrary stack $s$ with $\TOP{2}(s)=w$ and
  $\lvert s \rvert \geq 2$. 
\end{definition}

The next part of this section aims at a better understanding of the
dependence of $\ReturnFunc{k}(s)$ from $\TOP{2}(s)$. 
Let $w$ be the topmost word of the stack $s$. 
It will turn out that $\ReturnFunc{k}(w)$ only depends on
$\ReturnFunc{k}(\Pop{1}(w))$, on $\Sym(w)$ and on $\Lvl(w)$. This
means that the topmost element of $w$ and the number of returns of stacks
with topmost word $\Pop{1}(w)$ already determine the number of returns
of $s$. This implies that $\ReturnFunc{k}(s)$ can be computed as follows.
First, we compute the number of returns of stacks with topmost word
$\bot$. Then we compute $\ReturnFunc{k}(w_i)$ where $w_i$ is the
prefix of $w$ of length $i$ from $i=2$, $i=3$, \dots until we have
computed $\ReturnFunc{k}(w_i)$ for $w_i=w$ or equivalently, for
$i=\lvert w \rvert$.
Before we prove this claim in detail, let us give an
example. 

\begin{example} \label{Ex:DecompositionSubReturns}
  Consider the pushdown system $\mathcal{S}$ given by the transitions
  \begin{align*}
    &(q_0, a, \gamma_0, q_2, \Collapse),& 
    &(q_1, a, \gamma_1, q_2, \Collapse),& 
    &(q_2, b, \gamma_2, q_1, \Push{a,2}),\\
    &(q_0, a, \gamma_3, q_3, \Push{c,2}),& 
    &(q_3, c, \gamma_4, q_2, \Clone{2}),& 
    &(q_2, c, \gamma_5, q_2, \Pop{1}), \\
    &(q_2, b, \gamma_6, q_2, \Pop{1}),& 
    &(q_2, a, \gamma_7, q_2, \Pop{1}),& \text{ and }
    &(q_2, \bot, \gamma_8, q_2, \Pop{2}).
  \end{align*}
  Consider the stacks 
  \begin{align*}
  &s=\bot(b,2,0)(b,2,0):  \bot(b,2,1)a&\text{ and }& &s':=\Pop{1}(s).     
  \end{align*}
  There are exactly two returns of $\mathcal{S}$ from $(q_2,s')$ to
  $(q_2,\Pop{2}(s'))$. These are depicted in Figure
  \ref{fig:ThirdReturnExample}. We call them $\rho_1$ and $\rho_2$. 
  \begin{figure}
    \centering
    \begin{xy}
      \xymatrix@R=0pt@C=5pt{ 
        &\rho_1:\\
        &(b,2,0) & \\
        &(b,2,0) & (b,2,1) \\
        q_2, &\bot \ar[dddd]^{\gamma_2}&\bot\\
        \\
        \\
        \\
        &(b,2,0) & (a,2,1)\\
        &(b,2,0) & (b,2,1) \\
        q_1, &\bot \ar[dddd]^{\gamma_1}&\bot\\
        \\
        \\
        \\
        &(b,2,0)\\
        &(b,2,0) & \\
        q_2, &\bot &
        }\hskip 10cm
        \xymatrix@R=0pt@C=5pt{ 
        &\rho_2:\\
        &(b,2,0) & \\
        &(b,2,0) &(b,2,1)  \\
        q_2, &\bot \ar[dddd]^{\gamma_6}&\bot \\
        \\
        \\
        \\
        &(b,2,0) & \\
        &(b,2,0) & \\
        q_2, &\bot \ar[dddd]^{\gamma_8}&\bot\\
        \\
        \\
        \\
        &(b,2,0)  &\\
        &(b,2,0)  \\
        q_2, &\bot
        }
      \end{xy}
    \caption{The two returns from $(q_2,s')$ to $(q_2,\Pop{2}(s'))$.}
    \label{fig:ThirdReturnExample}
  \end{figure}
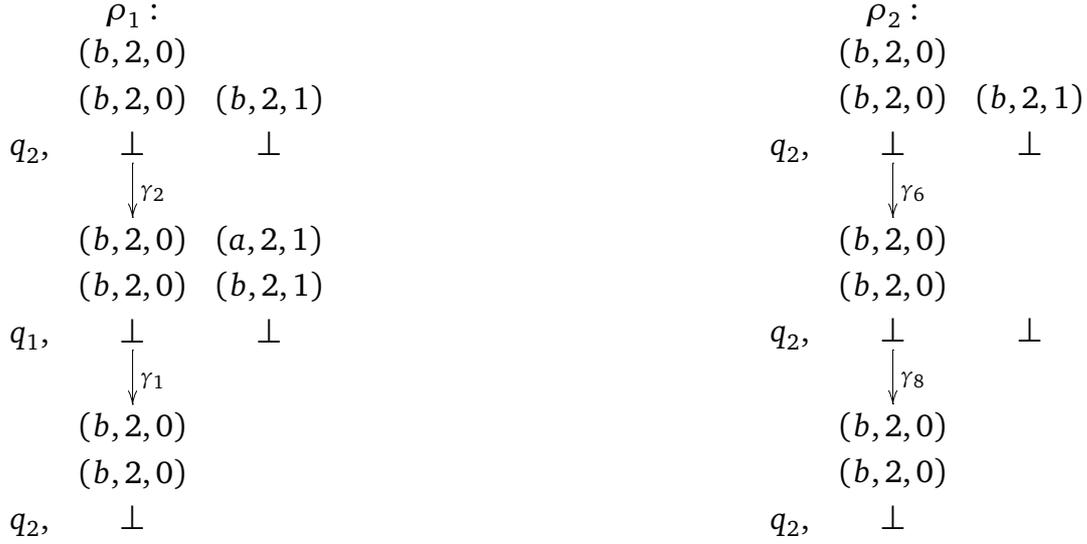

  We explain how  returns from $(q_0, s)$ to
  $(q_2,\Pop{2}(s))$ depend on those from $(q_2, s')$ to
  $(q_2,\Pop{2}(s'))$. 
  First of all note that there is a return $\pi_1$ from $(q_0,s)$ to
  $(q_2, \Pop{2}(s))$ as depicted on the left side of Figure
  \ref{fig:FourthReturnExample}.
  This return $\pi_1$ decomposes as $\pi_1= \pi_1{\restriction}_{[0,1]} \circ
  \rho_1$. If we replace $\rho_1$ by the other return $\rho_2$, then
  we obtain again a return which we call 
  $\pi_2:=\pi_1{\restriction}_{[0,1]}\circ \rho_2$. This run is
  depicted on the right side of Figure \ref{fig:FourthReturnExample}.
  In the following, we consider $\pi_1{\restriction}_{[0,1]}$ 
  as a representative for the returns $\pi_1$ and $\pi_2$ because both
  returns can be obtained from $\pi_1{\restriction}_{[0,1]}$ by
  attaching  a return with topmost word $\Pop{1}(\TOP{2}(s))$. 
  Furthermore, the existence of $\pi_1{\restriction}_{[0,1]}$ only
  depends on $\Sym(s)$ and $\Lvl(s)$: on any stack with topmost symbol
  $a$ of link level $1$, we can perform the sequence of transitions
  $\pi_1{\restriction}_{[0,1]}$ consists of. 
  \begin{figure}
    \centering
    \begin{xy}
      \xymatrix@R=0pt@C=5pt{ 
        &\pi_1:\\
        &(b,2,0) & a \\
        &(b,2,0) & (b,2,1) \\
        q_0, &\bot \ar[dddd]^{\gamma_0}&\bot\\
        \\
        \\
        \\
        &(b,2,0) & \\
        &(b,2,0) & (b,2,1) \\
        q_2, &\bot \ar[dddd]^{\gamma_2}&\bot\\
        \\
        \\
        \\
        &(b,2,0) & (a,2,1)\\
        &(b,2,0) & (b,2,1) \\
        q_1, &\bot \ar[dddd]^{\gamma_1}&\bot\\
        \\
        \\
        \\
        &(b,2,0)\\
        &(b,2,0) & \\
        q_2, &\bot &
        }\hskip 10cm
        \xymatrix@R=0pt@C=5pt{ 
        &\pi_2:\\
        &(b,2,0) & a \\
        &(b,2,0) &(b,2,1)  \\
        q_0, &\bot \ar[dddd]^{\gamma_0}&\bot \\
        \\
        \\
        \\
        &(b,2,0) & \\
        &(b,2,0) &(b,2,1)  \\
        q_2, &\bot \ar[dddd]^{\gamma_6}&\bot \\
        \\
        \\
        \\
        &(b,2,0) & \\
        &(b,2,0) & \\
        q_2, &\bot \ar[dddd]^{\gamma_8}&\bot\\
        \\
        \\
        \\
        &(b,2,0)  &\\
        &(b,2,0)  \\
        q_2, &\bot
        }
      \end{xy}
    \caption{Two returns from $(q_0,s)$ to $(q_2,\Pop{2}(s))$.}
    \label{fig:FourthReturnExample}
  \end{figure}
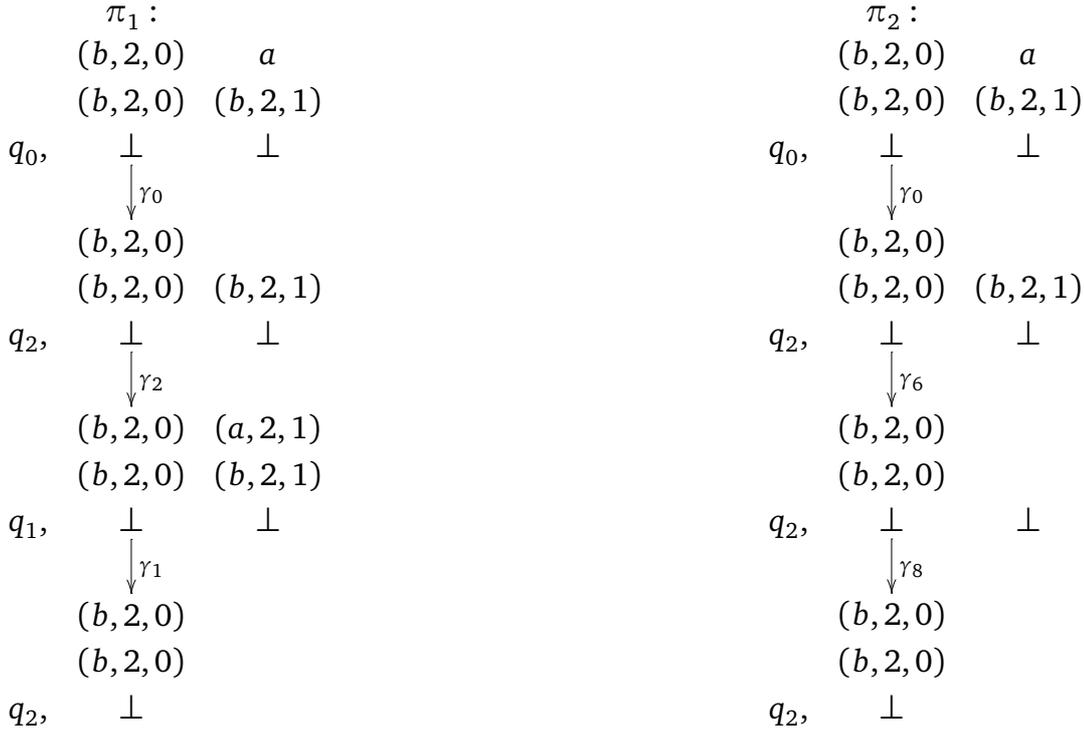
  Let us now turn to the other returns from $(q_0, s)$ to $(q_2,
  \Pop{2}(s))$.  Figure \ref{fig:FifthReturnExample} depicts another
  return $\pi_3$. 
  \begin{figure}
    \centering
    \begin{xy}
      \xymatrix@R=0pt@C=2pt{ 
        &(b,2,0) & a \\
        &(b,2,0) & (b,2,1) \\
        q_0, &\bot \ar[dddd]^{\gamma_3}&\bot\\
        \\
        \\
        \\
        &        & (c,2,1)\\
        &(b,2,0) & a\\
        &(b,2,0) & (b,2,1) \\
        q_3, &\bot \ar[dddd]^{\gamma_4}&\bot\\
        \\
        \\
        \\
        & &        (c,2,1)& (c,2,1) \\
        &(b,2,0) & a &  a \\
        &(b,2,0) & (b,2,1)& (b,2,1) \\
        q_2, &\bot \ar[dddd]^{\gamma_5}&\bot&\bot\\
        \\
        \\
        \\
        & &        (c,2,1)&  \\
        &(b,2,0) & a& a \\
        &(b,2,0) & (b,2,1)& (b,2,1) \\
        q_2, &\bot \ar[dddd]^{\gamma_7}&\bot&\bot\\
        \\
        \\
        \\
        & &        (c,2,1)&  \\
        &(b,2,0) & a&  \\
        &(b,2,0) & (b,2,1)& (b,2,1) \\
        q_2, &\bot \ar[dddd]^{\gamma_6}&\bot&\bot\\
        \\
        \\
        \\
        & &        (c,2,1)&  \\
        &(b,2,0) & a&  \\
        &(b,2,0) & (b,2,1)&  \\
        q_2, &\bot \ar[dddd]^{\gamma_8}&\bot&\bot\\
        \\
        \\
        \\
        & &        (c,2,1)&  \\
        &(b,2,0) &   a &  \\
        &(b,2,0) & (b,2,1)&  \\
        q_2, &\bot \ar[dddd]^{\gamma_5}&\bot\\
        \\
        \\
        \\
        &(b,2,0) & a&  \\
        &(b,2,0) & (b,2,1)&  \\
        q_2, &\bot &\bot\\
      }\hskip 10cm
      \xymatrix@R=0pt@C=2pt{ 
        &\ar[dddd]^{\gamma_7}  \\
        \\
        \\
        \\
        &(b,2,0) &&  \\
        &(b,2,0) & (b,2,1)&  \\
        q_2, &\bot \ar[dddd]^{\gamma_6}&\bot\\
        \\
        \\
        \\
         &(b,2,0) &&  \\
        &(b,2,0) & &  \\
        q_2, &\bot \ar[dddd]^{\gamma_8}&\bot\\
        \\
        \\
        \\
         &(b,2,0) &&  \\
        &(b,2,0) & &  \\
        q_2, &\bot \\
        }
      \end{xy}
    \caption{The  return $\pi_3$ from $(q_0,s)$ to $(q_2,\Pop{2}(s))$.}
    \label{fig:FifthReturnExample}
  \end{figure}
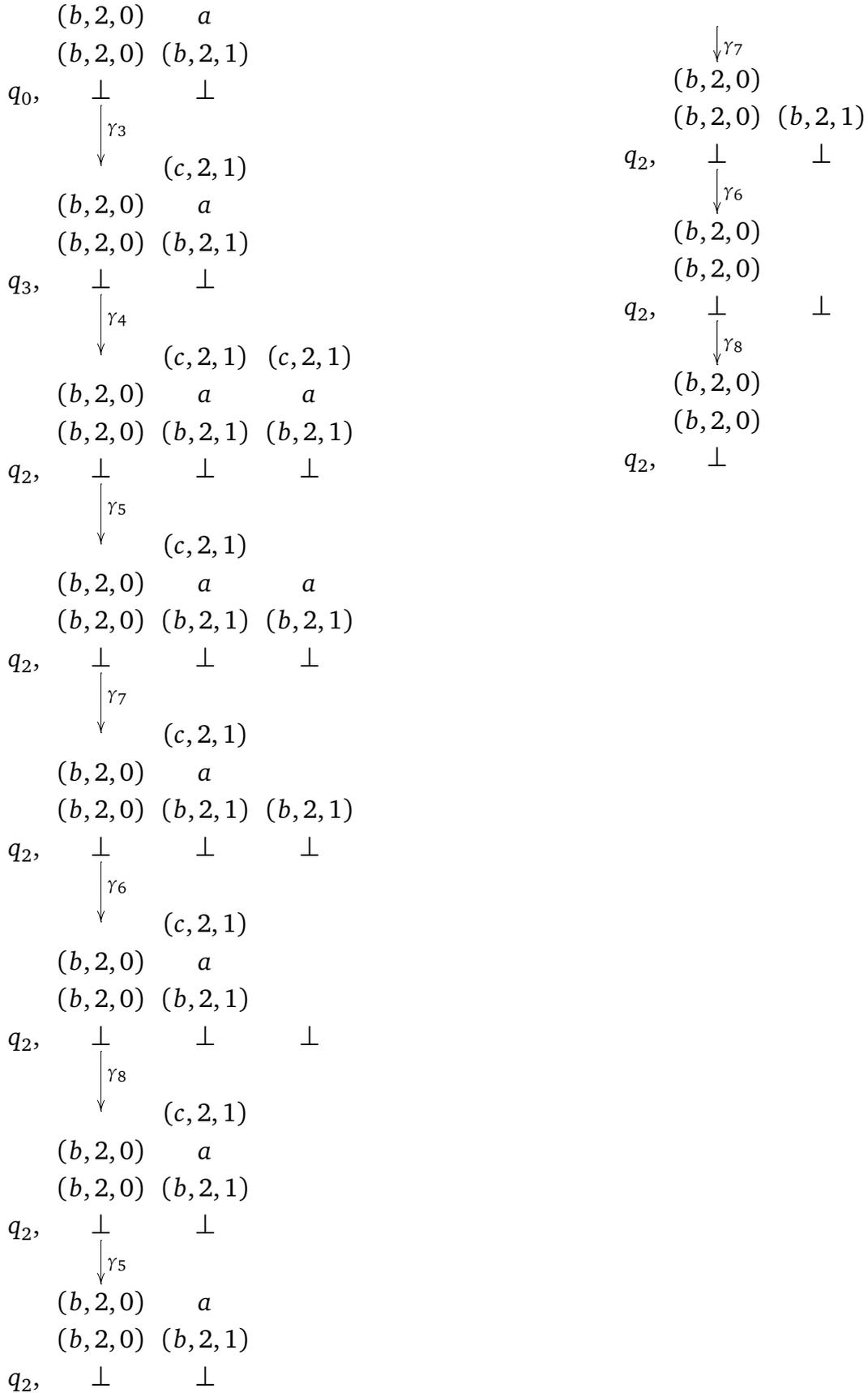
  Note that $\pi_3{\restriction}_{[4,6]}$ and
  $\pi_3{\restriction}_{[8,10]}$ are returns starting at stacks with
  topmost words $\Pop{1}(\TOP{2}(s))$. $\pi_3{\restriction}_{[8,10]}$
  is $\rho_2$ and $\pi_3{\restriction}_{[4,6]}$ copies the transitions
  of $\rho_2$ one by one. We can replace each of these parts
  of $\pi_3$ by the return $\rho_1$ (or by a one by one copy of its
  transitions) and obtain another return. We can also replace both
  parts by copies of the return 
  $\rho_1$ and obtain a fourth return. Thus, we obtain $4$ different
  returns from $(q_0, s)$ to 
  $(q_2, \Pop{2}(s))$ from the pair 
  $(\pi_3{\restriction}_{[0,4]}, \pi_3{\restriction}_{[6,8]})$ by
  plugging in different returns of 
  topmost word $\Pop{1}(\TOP{2}(s))$ after each element of this pair. 
  
  It is again an important observation that 
  $(\pi_3{\restriction}_{[0,4]}, \pi_3{\restriction}_{[6,8]})$ only 
  depends on $\Sym(s)$ and $\Lvl(s)$ in the following sense.
  Given any other stack $t$ with topmost symbol $a$ of link level $1$, 
  there is a run $\hat \pi_3{\restriction}_{[0,4]}$ that copies the
  transitions of $\pi_3{\restriction}_{[0,4]}$ and that ends in a
  stack $t'$ with $\TOP{2}(t')=\Pop{1}(\TOP{2}(t))$. 
  Similarly, we can copy the transitions of
  $\pi_3{\restriction}_{[6,8]}$ to a run starting at $\Pop{2}(t')$ and
  which ends again in a stack with topmost word
  $\Pop{1}(\TOP{2}(t))$. 
  
  It is easy to see that there are no other returns from $(q_0,s)$ to
  $(q_2, \Pop{2}(s))$ than the ones we discussed above. 

  Thus,  the tuple 
  $\pi_1{\restriction}_{[0,1]}, (\pi_3{\restriction}_{[0,4]},
  \pi_3{\restriction}_{[6,8]})$ represents all returns from the 
  configuration $(q_0,s)$
  to $(q_2, \Pop{2}(s))$ in
  the following sense.
  \begin{enumerate}
  \item   $\pi_1{\restriction}_{[0,1]}$
    can be turned into a return from $(q_0,s)$ to
    $(q_2, \Pop{2}(s))$ by appending a return of a stack with topmost
    word $\Pop{1}(\TOP{2}(s))$.
  \item $(\pi_3{\restriction}_{[0,4]}, \pi_3{\restriction}_{[6,8]})$
    can be turned into a return from $(q_0, s)$ to $(q_2, \Pop{2}(s))$
    by plugging in  one return of a stack with topmost word
    $\Pop{1}(\TOP{2}(s))$ between the two runs and by appending such a
    return to the end of $\pi_3{\restriction}_{[6,8]}$.
  \item  All returns from $(q_0, s)$ to $(q_2, \Pop{2}(s))$ are
    induced by this tuple in the sense of items 1 and 2. 
  \end{enumerate}
  Since there are $2$ returns from state $(q_2,t)$ to
  $(q_2, \Pop{2}(t))$ for any stack $t$ of width at least $2$ and
  topmost word $\TOP{2}(t)=\Pop{1}(\TOP{2}(s))$, we conclude that 
  there are $2 + 2\cdot 2=6$ different returns. 

  This form of computing the number of returns works for all stacks. 
  Take for example the stack 
  \begin{align*}
    \hat s:=\bot: \bot(b,2,1)(b,2,1)a.     
  \end{align*}
  This stack has the same topmost symbol and link level as $s$. 
  Thus, we can copy transition by transition the runs
  $\pi_1{\restriction}_{[0,1]}$  and 
  $\pi_3{\restriction}_{[0,4]}$  to runs $\hat \pi_1$ and $\hat\pi_2$
  starting from $\hat s$.  
  Furthermore, note that the stack of $\pi_3(6)$ is
  $\Pop{2}(\pi_3(4))$. For $\hat t$ the stack obtain via a $\Pop{2}$
  from the last stack of $\hat\pi_2$ we can copy 
  $\pi_3{\restriction}_{[6,8]}$ 
  transition by transition to a run $\hat\pi_3$ starting at $\hat t$.
  The resulting runs are depicted in Figure
  \ref{fig:SixthReturnExample}.  
  \begin{figure}
    \centering
    \begin{xy}
      \xymatrix@R=0pt@C=2pt{ 
        &\hat\pi_1:\\
        \\
        & & a \\
        & & (b,2,1) \\
        & & (b,2,1) \\
        q_0, &\bot \ar[dddd]^{\gamma_0}&\bot\\
        \\
        \\
        \\
        & & (b,2,1)\\
        & & (b,2,1) \\
        q_3, &\bot &\bot\\
        \\
        \\
        \\
        \\
        \\
        \\
        & \hat\pi_3:\\
        \\
        & & (c,2,1)& \\
        & & a& \\
        & & (b,2,1)& \\
        & & (b,2,1)& \\
        q_2, &\bot\ar[dddd]^{\gamma_5} &\bot
        \\
        \\
        \\
        \\
        & & & \\
        & & a& \\
        & & (b,2,1)& \\
        & & (b,2,1)& \\
        q_2, &\bot\ar[dddd]^{\gamma_7} &\bot
        \\
        \\
        \\
        \\        
        & & & \\
        & & & \\
        & & (b,2,1)& \\
        & & (b,2,1)& \\
        q_2, &\bot &\bot
      }\hskip 10cm
      \xymatrix@R=0pt@C=2pt{ 
        &\hat\pi_2:\\
        \\
        & & a \\
        & & (b,2,1) \\
        & & (b,2,1) \\
        q_0, &\bot \ar[dddd]^{\gamma_3}&\bot
        \\
        \\
        \\
        \\
        & & (c,2,1) \\
        & & a \\
        & & (b,2,1) \\
        & & (b,2,1) \\
        q_3, &\bot \ar[dddd]^{\gamma_4}&\bot
        \\
        \\
        \\
        \\
        & & (c,2,1)& (c,2,1) \\
        & & a      & a \\
        & & (b,2,1)& (b,2,1) \\
        & & (b,2,1)& (b,2,1) \\
        q_2, &\bot \ar[dddd]^{\gamma_5}&\bot&\bot
        \\
        \\
        \\
        \\
        & & (c,2,1)& \\
        & & a      &  a  \\
        & & (b,2,1)& (b,2,1) \\
        & & (b,2,1)& (b,2,1) \\
        q_2, &\bot \ar[dddd]^{\gamma_7}&\bot&\bot
        \\
        \\
        \\
        \\
        & & (c,2,1)& \\
        & &    a   & \\
        & & (b,2,1)& (b,2,1) \\
        & & (b,2,1)& (b,2,1) \\
        q_2, &\bot &\bot&\bot
        }
      \end{xy}
    \caption{$\hat\pi_1$ corresponding to
      $\pi_1{\restriction}_{[0,1]}$, 
      $\hat\pi_2$ corresponding to $\pi_3{\restriction}_{[0,4]}$, and
      $\hat\pi_3$ corresponding to $\pi_3{\restriction}_{[6,8]}$. }
    \label{fig:SixthReturnExample}
  \end{figure}

  Again, we can  turn $\hat \pi_1$ and the pair $(\hat\pi_2,
  \hat\pi_3)$ into returns from $(q_0, \hat s)$ 
  to $(q_2, \Pop{2}(s))$. For this purpose,  we have to plug in
  returns with topmost word 
  $\bot(b,2,1)(b,2,1)$ from state $q_2$ to state $q_2$ after
  $\hat\pi_1, \hat\pi_2$, and 
  $\hat\pi_3$.  

  There are exactly three returns
  from \mbox{$(q_2,  \bot:\bot(b,2,1)(b,2,1))$} to 
  $(q_2, [\bot])$. The first performs $\trans{\gamma_2}$ and then
  $\trans{\gamma_1}$, the second performs $\trans{\gamma_6},
  \trans{\gamma_2}$ and $\trans{\gamma_1}$, and the last one
  performs $\trans{\gamma_6}$, $\trans{\gamma_6}$ and 
  $\trans{\gamma_8}$. 

  Since we have to append such a return to $\hat\pi_1$  in order to
  obtain a return of $\hat s$, $\hat\pi_1$ induces $3$ different
  returns from $(q_0, \hat s)$ to $(q_2, \Pop{2}(\hat s))$. 
  Moreover, using $2$ of these returns, we can turn the pair
  $(\hat\pi_2, \hat \pi_3)$ into a return from $(q_0, \hat s)$ to
  $(q_2, \Pop{2}(\hat s))$. Hence, there are $3\cdot 3=9$
  possibilities to turn this pair into a return. 
  We conclude that there are $3+9=12$ returns from $(q_0, \hat s)$ to
  $(q_2, \Pop{2}(\hat s))$. We leave it as an exercise to figure out
  that there are exactly $12$ returns from ($q_0, \hat s)$
  to $(q_2, \Pop{2}(\hat s))$. 
\end{example}

The previous example pointed to a connection between returns of a
stack with topmost word $w$ and the returns of stacks with topmost
word $\Pop{1}(w)$. 
The main result of this section is that this connection can be used to
define a finite automaton that calculates on input $\TOP{2}(s)$ the
function $\ReturnFunc{k}(s)$ for a given $k\in\N$. 
Furthermore, this dependence can be used to calculate a bound on the
length of returns in dependence of the length of the topmost word of a
stack. 
We first state these two results, afterwards we provide the technical
background for the proofs. 

\begin{proposition}  \label{Prop:AutomatonForReturns}
  There is an algorithm that, given a collapsible pushdown system
  $\mathcal{S}$ of level $2$, computes a
  deterministic finite automaton  
  $\mathcal{A}_{\mathrm{ret}}$ with the following property. 
  $\mathcal{A}_{\mathrm{ret}}$ computes $\ReturnFunc{k}(s:w)$ on input
  $\pi(w)$ where $\pi(w)$ denotes the projection of $w$ to its symbols
  and link levels.  
\end{proposition}

\begin{proposition} \label{Prop:FuncBoundReturnLengthLemma}
  There is an algorithm that, on input some $2$-\CPG
  $\mathcal{S}$ and a natural number $k$,  computes  a
  function
  $\FuncBoundReturnLength{\mathcal{S}}{k}:\N \to\N$ with the following
  properties. 
  \begin{enumerate}
  \item 
    For each stack $s$, for states $q_1,q_2$ and for
    $i:=\ReturnFunc{k}(s)(q_1,q_2)$, the
    length-lexicographically shortest returns 
    $\rho_1,\dots, \rho_i$ from
    $(q_1,s)$ to $(q_2,\Pop{2}(s))$ satisfy 
    \begin{align*}
      \length(\rho_j)\leq \FuncBoundReturnLength{\mathcal{S}}{k}(\lvert
      \TOP{2}(s)\rvert)\text{ for all }1\leq j \leq i.      
    \end{align*}
  \item If there is a return $\rho$ from $(q_1,s)$ to $(q_2, \Pop{2}(s))$ with
    $\length(\rho) > \FuncBoundReturnLength{\mathcal{S}}{k}(\lvert
    \TOP{2}(s)\rvert)$, then there are $k$ returns from $(q_1,s)$ to
    $(q_2, \Pop{2}(s))$ of length at most 
    $\FuncBoundReturnLength{\mathcal{S}}{k}(\lvert \TOP{2}(s)\rvert)$.  
  \end{enumerate}
\end{proposition}

For any stack $s$ with topmost word $\bot$, we will calculate the
number of returns of $s$ using Lemma \ref{LemmaCountNumberOfRuns}. 
 
We can inductively calculate the returns of some stack $s$ as
follows. Assume that we already know how to calculate returns of
stacks with topmost word of size $\lvert \TOP{2}(s)\rvert -1$. 
Any return of $s$ splits into those parts that only
depend on its topmost symbol and link level and those parts that are
returns from stacks with smaller topmost word 
(cf. Example \ref{Ex:DecompositionSubReturns}). By induction
hypothesis we already counted the latter parts. Hence,  we have to
focus on the other parts. Here again, we can reduce the counting of
these runs to an application of Lemma \ref{LemmaCountNumberOfRuns}.
This reduction to Lemma \ref{LemmaCountNumberOfRuns} is uniformly in
the length of the topmost word from $s$. Due to this uniformity, we
can then compute a finite automaton that calculates the number of
returns. 

The reader who is not interested in the technical details of the
proofs of the propositions may safely skip this part and continue
reading Section \ref{sec:CompLoops}. 

We start the analysis of returns with a general observation. 
By definition, there are returns $\rho$ 
where  there is some $i\in\domain(\rho)$ such that
$\rho{\restriction}_{[i,\length(\rho)]}$ is a return from
$\Pop{1}(\rho(0))$. 
Our first lemma shows that a run $\rho$ which visits
$\Pop{1}(\rho(0))$ is a return if and only if a suffix of $\rho$ is a 
return from $\Pop{1}(\rho(0))$ to $\Pop{2}(\rho(0))$. 

\begin{lemma} \label{Lemma_ReturnDecomposition}
  For $\rho$ a return from $s$ to $\Pop{2}(s)$ and for
  $i\in\domain(\rho)$ minimal such that
  $\rho(i)=\Pop{1}(s)$, the restriction
  $\rho{\restriction}_{[i,\length(\rho)]}$ is a return from $\Pop{1}(s)$
  to $\Pop{2}(s)$. 
\end{lemma}
\begin{proof}
  If $\rho$ ends with a $\Pop{2}$ transition, there is nothing to
  show. Now assume that $\rho$ ends with a $\Collapse$ operation. 
  If  $\TOP{2}(s)\leq \TOP{2}(\rho(\length(\rho)-1))$ then
  \begin{align*}
    \TOP{2}(\Pop{1}(s)) < \TOP{2}(s) \leq
    \TOP{2}(\rho(\length(\rho)-1))    
  \end{align*}
  immediately yields the claim. 
  The last possible case is that there is some $j\in\domain(\rho)$
  such that $\rho{\restriction}_{[j,\length(\rho)]}$ is a return of
  $\Pop{1}(s)$. Since $i\leq j$, this immediately implies the claim as
  $\rho{\restriction}_{[i,j]}$ is a run from $\Pop{1}(s)$ to $\Pop{1}(s)$
  that never visits $\Pop{2}(s)$. But the class of returns is  closed
  under prefixing by such runs.  
\end{proof}

The previous observation gives rise to a classification of returns
into \emph{low} and \emph{high} ones. 
\begin{definition}
  Let $\rho$ be some return. We call $\rho$ a \emph{low return}, if there
  is some 
  \mbox{$i\in\domain(\rho)$} such that $\rho(i)=\Pop{1}(\rho(0))$. Otherwise
  we call $\rho$ a \emph{high return}. 
\end{definition}
\begin{remark}
  Due to \ref{Lemma_ReturnDecomposition}, a low return decomposes as a
  run to $\Pop{1}(\rho(0))$ followed by a return of
  $\Pop{1}(\rho(0))$. High returns never pass
  $\Pop{1}(\rho(0))$. Hence, low returns pass ``lower'' stacks than
  high returns. 
\end{remark}

In fact, the analysis of high returns and low returns is very
similar. But there are small differences which provoke a lot of case
distinctions when dealing with both types at the same time. 
In order to avoid these case distinctions, 
some of our lemmas will concentrate on high returns and we will only
remark the differences to the case of low returns.

Next, we show that the notion $\ReturnFunc{k}(w)$
(cf. Definition \ref{Def:returnsofWord}) is  well-defined for every
word $w$. For this purpose let us first introduce auxiliary notation. 

\begin{definition} \label{def:downarrowzero}
  The word $w{\downarrow_0}$ is obtained from
  $w\in(\Sigma\cup(\Sigma\times\{2\}\times\N))^*$ 
  by replacing every occurrence of $(\sigma,2,j)$ in $w$ by $(\sigma,2,0)$
  for all $\sigma\in \Sigma$ and all $j\in \N$.
\end{definition}
\begin{remark}\label{WortLevelProjektion}
  Later, it is important that 
  $w{\downarrow_0}$ is a word over the finite alphabet
  $\Sigma\cup(\Sigma\times\{2\}\times\{0\})$. 
\end{remark}

\begin{definition} \label{Def:ReturnEquivalent}
  Let $s, s'$ be stacks such that $\TOP{2}(s){\downarrow_0} =
  \TOP{2}(s'){\downarrow_0}$. Let $\rho$ be a return from $(q_1,s)$ to
  $(q_2,\Pop{2}(s))$ and $\rho'$ be a return from $(q_1, s')$ to
  $(q_2, \Pop{2}(s'))$. 
  We say $\rho$ and $\rho'$ are \emph{equivalent} returns if they consist of
  the same sequence of transitions.
\end{definition}
For an example, note that the returns $\rho$ and $\rho'$ in Figure
\ref{fig:FirstReturnExample} are equivalent.
The crucial observation is that different stacks whose topmost words
agree on their symbols and link levels have the same returns modulo
this equivalence of returns. 

\begin{lemma}\label{ReturnDependTopword}
  Let $s$ and $s'$ be stacks of width at least $2$ such that
  $\TOP{2}(s){\downarrow_0} = \TOP{2}(s'){\downarrow_0}$. If $\rho$ is a
  return from $(q_1, s)$ to $(q_2, \Pop{2}(s))$ then there is an
  equivalent return $\rho'$ from $(q_1, s')$ to $(q_2, \Pop{2}(s'))$. 
\end{lemma}
\begin{proof} 
  We assume that there
  is a symbol $\Box\in\Sigma$ not occurring in any of the transitions
  of the pushdown system.  Let $s$ be some stack and $w:=\TOP{2}(s)$. 
  
  Now, we define $s':=w{\downarrow_0} \Box:
  w{\downarrow_0}$.
  This definition is tailored towards the fact that $s'$ is minimal
  with the following two properties.
  \begin{enumerate}
  \item The assumption that $\Box$ does not appear
    within the transitions implies that an arbitrary run from $s'$ to
    $\Pop{2}(s')$ is a return.\footnote{This fact is not important for
    the proof of this lemma, but this fact gets important in the
    next lemmas.}
  \item $\TOP{2}(s){\downarrow_0} = \TOP{2}(s'){\downarrow_0}$.
  \end{enumerate}
  In order to prove the lemma, it suffices to show that for every return
  starting in $s$ there is an equivalent one starting in $s'$ and vice
  versa. 

  The proof of this lemma is as follows. 
  Let $t:=\Pop{2}(s)$, $t':=\Pop{2}(s')$ and let 
  \mbox{$k:=\lvert t \rvert - \lvert t' \rvert = \lvert t \rvert -1$}. 
  Furthermore, let
  $\rho$ be some 
  return from $(q,s)$ to 
  $(\hat q, t)$. 
  We define $\rho'$ to be the largest run which starts
  in $(q,s')$ and copies an initial part of $\rho$ transition by
  transition. 
  
  We prove that $\domain(\rho')=\domain(\rho)$ by showing a
  stronger claim. For some word $v$ let $v^{-k}$ denote the word that
  is obtained from $v$ by replacing every link $l$ of level $2$ by the
  link $l-k$. 
  \begin{claim}
    The domains of $\rho'$ and $\rho$ agree. Furthermore, 
    for each $i\in\domain(\rho)$ the following holds.
    \begin{enumerate}
    \item  The states of $\rho(i)$ and
      $\rho'(i)$ agree
    \item  The stack of $\rho(i)$ decomposes as
      $t: w_1v_1:w_2v_2:\dots:w_nv_n$ and the stack of $\rho'(i)$ decomposes as
      $t': w_1{\downarrow_0} v_1^{-k}: w_2{\downarrow_0} v_2^{-k}: \dots:
      w_n{\downarrow_0} v_n^{-k}$
      where the $w_i$ are chosen in such a way that 
      $w_i{\downarrow_0}$ is a maximal  prefixes of
      $\TOP{2}(s){\downarrow_0}$.
    \item If the operation at $i$ is a collapse of link level
      $2$, then $v_n$ is nonempty. 
    \end{enumerate}
  \end{claim}
  This claim can be proved by induction on the domain of $\rho$. This is 
  tedious but straightforward.
  The construction of $\rho'$ can be seen as the application of a
  prefix replacement as follows: $\rho'=\rho[t:\bot/t'\bot]$ where we
  manually repair the links of level $2$. 
  
  The lemma then follows as a direct corollary of the claim: just note
  that the last operation of $\rho$ is a $\Collapse$ of level
  $2$ or a 
  $\Pop{2}$ and yields a stack of width $\lvert s \rvert -1$. In both
  cases it follows directly from the statements of 
  the claim that the same transition is applicable to $\rho'(i)$
  and results in a stack of width $\lvert s\rvert -1 - \lvert s \rvert
  +2=1$. Since $\rho'$ never changed the first word of the stack, this
  stack is $t'=\Pop{2}(s')$.  
\end{proof}

The previous lemma shows that $\ReturnFunc{k}(w)$ is well defined
(cf. Definition \ref{Def:returnsofWord}) and 
$\ReturnFunc{k}(s) = \ReturnFunc{k}(\TOP{2}(s){\downarrow_0})$ 
for all stacks $s$ of width at least $2$.  
Thus, if we want to compute $\ReturnFunc{k}(s)$, we can concentrate of
the returns of a fixed stack with topmost word
$w:=\TOP{2}(s){\downarrow_0}$.
We will do this by choosing the stack
$s':=w{\downarrow_0} \Box:w{\downarrow_0}$ to be the representative of any
stack with $\TOP{2}(s)=w$. 
$s'$ is the smallest
stack with topmost word $w{\downarrow_0}$ such that any run from $s'$ to
$\Pop{2}(s')$ is a return.

The following lemma contains the observation that every return from
some stack $s$ to $\Pop{2}(s)$ decomposes into parts that are prefixed
by $s$  and parts that are returns of stacks with topmost word
$\Pop{1}(\TOP{2}(s))$. This lemma shows that the decomposition of
the returns in Example \ref{Ex:DecompositionSubReturns} can be
generalised to decompositions of all returns. 

\begin{lemma} \label{Lemma_Return_Decomposition_For_Simulation}
  Let $\rho$ be some high return of some stack $s$ with topmost word
  $w:=\TOP{2}(s)$. Then there is a well-defined sequence
  \begin{align*}
    0:=j_0<i_1<j_1<i_2<j_2<\dots<i_n<j_n < i_{n+1}:=\length(\rho)-1    
  \end{align*}
  with the following properties.
  \begin{enumerate}
  \item
    \label{Lemma_Return_Decomposition_For_Simulation_PrefixCondition}
    For $1\leq k \leq n+1$,  
    $s\prefixeq \rho{\restriction}_{[j_{k-1}, i_k]}$.
  \item For all $1\leq k \leq n$, $\TOP{2}(\rho(i_k))=w$ and the
    operation at $i_k$ in $\rho$ is a $\Pop{1}$ or a collapse of level
    $1$.
  \item Either $w$ is a proper prefix of $\TOP{2}(\rho(i_{n+1}))$ and
    the  operation at $i_{n+1}$ is a collapse of level $2$ or
    $w$ is a prefix of $\TOP{2}(\rho(i_{n+1}))$ and the operation at
    $i_{n+1}$ is a $\Pop{2}$. 
  \item 
    \label{Lemma_Return_Decomposition_For_Simulation_ReturnCondition}
    For each $1\leq k \leq n$, there is a stack $s_k$ with
    $\TOP{2}(s_k)=\Pop{1}(w)$ such that
    $\rho{\restriction}_{[i_{k}+1,j_k]}$ is a 
    return from $s_k$ to $\Pop{2}(s_k)$. 
  \end{enumerate}
\end{lemma}
\begin{remark}
  If $\rho$ is a low return, a completely analogous lemma holds. We
  just have to omit $i_{n+1}$, i.e., the sequence ends with
  $j_n=\length(\rho)$. Then statements 1, 2, and 4 hold for this sequence
  $0:=j_0<i_1<j_1<i_2<j_2<\dots<i_n<j_n=\length(\rho)$. 
\end{remark}
\begin{proof}
  Set $j_0:=0$.
  Let
  $i_1\in\domain(\rho)$ be the minimal position such that 
  $s\prefixeq \rho(i_1)$ but $s\notprefixeq\rho(i_1+1)$. If
  $i=\length(\rho)-1$ we set $n:=0$ and we are done: due to the
  definition of 
  a high return, the last operation is either a collapse of level $2$ and 
  $w<\TOP{2}(\rho(i_1))$ or it is a $\Pop{2}$ and $w\leq
  \TOP{2}(\rho(i_1))$. 
  
  So let us assume that $i_1<\length(\rho)-1$.
  By definition, $s\prefixeq\rho{\restriction}_{[j_0,i_1]}$.
  
  Since $\rho$ is a return
  of $s$ and $i_1+1<\length(\rho))$, 
  $\lvert \rho(i_1+1) \rvert \geq \lvert s \rvert$.
  Lemma \ref{Lem:PrefixesAndOperations} then implies that
  $\TOP{2}(\rho(i_1))=w$ and $\TOP{2}(\rho(i_1+1))=\Pop{1}(w)$. 
  
  Since $\lvert \rho(i_1+1) \rvert \geq \lvert s \rvert > \lvert
  \rho(\length(\rho)) \rvert$, there is some minimal $j_1$ such that
  $i_1<j_1$ and $\lvert
  \rho(j_1) \rvert < \lvert \rho(i_1) \rvert$. We want to prove the
  following claim. 
  \begin{claim}
    For $s'$ the stack at $\rho(i_1+1)$, $\rho{\restriction}_{[i_1+1,j_1]}$
    is a return from  $s'$ to $\Pop{2}(s')$. 
  \end{claim}
    
  First observe that by definition of $j_1$ for all $i_1+1\leq k <j_1,
  \lvert\rho(k)\rvert \geq \lvert s' \rvert$. 
  The operation at $j_1-1$ has to be a $\Pop{2}$ or
  a $\Collapse$ (of link level $2$) because it decreases the width of
  the stack.

  If it is
  $\Pop{2}$, then we conclude that $\rho(j_1)=\Pop{2}(s')$ and the claim
  is satisfied.
  
  Now, we consider the case that the operation before $j_1$ is a collapse
  of level $2$. 
  Since $\rho$ is a high return, $\rho(i_1+1)\neq\Pop{1}(s)$. 
  Thus, $\lvert s' \rvert > \lvert s \rvert$. 
  Since $\TOP{2}(s') =\Pop{1}(w)$, all elements
  in $\TOP{2}(s')$ are clones of elements in the topmost word $w$ of
  $s$. Thus, their level $2$ links 
  point to stacks $t$ with $\lvert t \rvert < \lvert s \rvert <
  \lvert s' \rvert$. 
  Heading for a contradiction, let us first assume that
  $\TOP{2}(\rho(j_1-1))$ is a prefix of $\TOP{2}(s')$,
  i.e.,  $\TOP{2}(\rho(j_1-1))\leq \TOP{2}(s') < w$. 
  In this case,  
  $\lvert \rho(j_1) \rvert < \lvert s \rvert$ whence
  $j_1=\length(\rho)$ and the 
  operation at $j_1-1$ is the last collapse operation in $\rho$. But
  $\Pop{1}(s)$ does not occur within $\rho$ 
  because $\rho$ is a high return.
  We conclude that
  the last operation of $\rho$ is $\Collapse$, but 
  neither $w < \TOP{2}(\rho(j_1-1))$ nor a final segment of $\rho$ is a
  return of $\Pop{1}(s)$. This contradicts the definition of a
  return. 
    
  Thus, we conclude that the topmost element
  of $\rho(j_1-1)$ was pushed onto the stack between $i_1+1$ and $j_1-1$. Since
  $\lvert \rho(k)\rvert \geq \lvert s' \rvert$ for  all $\rho(k)$ with
  $i_1+1\leq k \leq j_1-1$, the  link of this element is at least $\lvert s'
  \rvert -1$. But by definition of $j$ this link also points below
  $s'$, whence the link is $\lvert s' \rvert -1$. But then
  $\rho{\restriction}_{[i,j]}$ satisfies all requirements of a
  return of $s'$ and we are done. 

  This completes the claim. 

  Thus, we have obtained that $i_1$ and $j_1$ are candidates for
  the initial elements of the sequence required by the lemma. 
  Note that the proof yields even more
  information. We have seen that $j_1<\length(\rho)$. 
  Thus, $\lvert \rho(i_1) \rvert > \lvert \Pop{2}(\rho(i_1)) \rvert =
  \lvert \rho(j_1) \rvert \geq \lvert s \rvert$.  
  Lemma \ref{Lem:PrefixPop1Pop2wieder} implies that
  $s\prefixeq\rho(j_1)$ because $s\prefixeq\rho(i_1)$. 
  
  Hence, we can use the same arguments
  (restricted to $\rho{\restriction}_{[j_1,\length(\rho)]}$)
  to show that for  the minimal $i_2>j_1$ such that
  $s\notprefixeq \rho(i_2+1)$, we have
  $\TOP{2}(\rho(i_2+1))=\Pop{1}(w)$.  By induction
  one concludes that the whole run $\rho$ decomposes into parts
  prefixed by $s$ and returns of stacks with topmost word $\Pop{1}(w)$
  as desired.
\end{proof}
\begin{remark}
  For low returns the proof is analogous. The only difference is the
  following. When defining inductively $0=j_0< i_1< j_1 < \dots <
  i_k$ at some point, we will obtain that $\rho(i_k)=\Pop{1}(s)$. 
  In this case, we set $j_k:=\length(\rho)$. 
  Lemma \ref{Lemma_ReturnDecomposition} shows that
  $\rho{\restriction}_{[i_k,j_k]}$ is a return. 
  Thus, this definition satisfies the claim of the lemma for the case
  of low returns. 
\end{remark}

\begin{lemma}
  In Lemma \ref{Lemma_Return_Decomposition_For_Simulation}, the
  sequence $0=j_0<i_1<j_1<i_2<j_2<\dots<i_{n+1}$ is uniquely 
  defined by conditions 
  \ref{Lemma_Return_Decomposition_For_Simulation_PrefixCondition}
  and
  \ref{Lemma_Return_Decomposition_For_Simulation_ReturnCondition}:
  assume that there is another sequence 
  \begin{align*}
    0=l_0<k_1< l_1 < k_2 < l_2 <
    \dots < k_{n+1}=\length(\rho)    
  \end{align*}
  satisfying these conditions. 
  Then $l_0=j_0, i_1=k_1, j_1=l_1, \dots, i_{n+1} = k_{n+1}$.
\end{lemma}
\begin{proof}
  If $i_1<k_1$ then
  $\rho{\restriction}_{[l_0,k_1]}$ contains $\rho(i_1)$ but
  $s\notprefixeq\rho(i_1)$ which is a contradiction. If $k_1< i_1$,
  we derive the contradiction $s\notprefixeq\rho(k_1)$ analogously.
  Now, $l_1=j_1$ follows from the fact that a return of
  $\rho(k_1)$ has to visit a stack $s'$ with $\lvert s' \rvert <
  \lvert \rho(k_1) \rvert$ at its last position but it is not allowed
  to do so before. But $j_1$ is the minimal position where such a
  stack is reached whence $l_1 = j_1$. The claim follows by induction.
\end{proof}

Before we continue our analysis of returns, it is useful to fix an
enumeration of all runs of a pushdown system. 
\begin{assumption}
  Let $\mathcal{S}$ be some pushdown system and $\Delta$ its
  transition relation. From now on, we assume that $\Delta$ is a linearly
  ordered set. Thus, all runs of $\mathcal{S}$ that start in a fixed
  configuration  are well-ordered via the length-lexicographic
  ordering of the transitions that they use. 
\end{assumption}

\begin{figure}
  \centering
  \begin{xy}
    \xymatrix@R=0pt@C=2pt{ 
             & \Box\\
             &a    & a \\
             &\top & \top \\
        q_0, &\bot \ar[dddd]^{\gamma_0}&\bot\\
        \\
        \\
        \\
        &\Box \\
        &a  & \\
        &\top & \top \\
        q_2, &\bot \ar[dddd]^{\mathrm{Rt}_1}&\bot\\
        \\
        \\
        \\
        & \Box\\
        & a\\
        & \top & \\
        q_2, &\bot &
        }\hskip 10cm
      \xymatrix@R=0pt@C=2pt{ 
        & \Box\\
        &a   & a \\
        &\top & \top \\
        q_0, &\bot \ar[dddd]^{\gamma_3}&\bot\\
        \\
        \\
        \\
        & \Box   & (c,2,1)\\
        &a       & a\\
        &\top    & \top \\
        q_3, &\bot \ar[dddd]^{\gamma_4}&\bot\\
        \\
        \\
        \\
        & \Box  & (c,2,1)& (c,2,1) \\
            & a & a      &  a \\
        & \top  & \top   & \top \\
        q_2, &\bot \ar[dddd]^{\gamma_5}&\bot&\bot\\
        \\
        \\
        \\
        &  \Box&   (c,2,1)&  \\
        & a & a& a \\
        & \top & \top& \top \\
        q_2, &\bot \ar[dddd]^{\gamma_7}&\bot&\bot\\
        \\
        \\
        \\
        &  \Box&   (c,2,1)&  \\
        & a & a&  \\
        & \top & \top& \top \\
        q_2, &\bot \ar[dddd]^{\mathrm{Rt}_2}&\bot&\bot\\
        \\
        \\
        \\
        & \Box &(c,2,1)&  \\
        & a &         a&  \\
        &\top &    \top&  \\
        q_2, &\bot \ar[dddd]^{\gamma_5}&\bot\\
        \\
        \\
        \\
        & \Box \\
        &a & a&  \\
        &\top & \top&  \\
        q_2, &\bot \ar[dddd]^{\gamma_7}&\bot\\
        \\
        \\
        \\
        & \Box \\
        &a  &&  \\
        & \top & \top &  \\
        q_2, &\bot \ar[dddd]^{\mathrm{Rt}_2}&\bot\\
        \\
        \\
        \\
        & \Box \\
         &a &&  \\
        &\top & &  \\
        q_2, &\bot 
        }
  \end{xy}
  \caption{Simulation of $\pi_1$ on the left and $\pi_3$ on the right.}
  \label{fig:SimulationsOfReturns}
\end{figure}
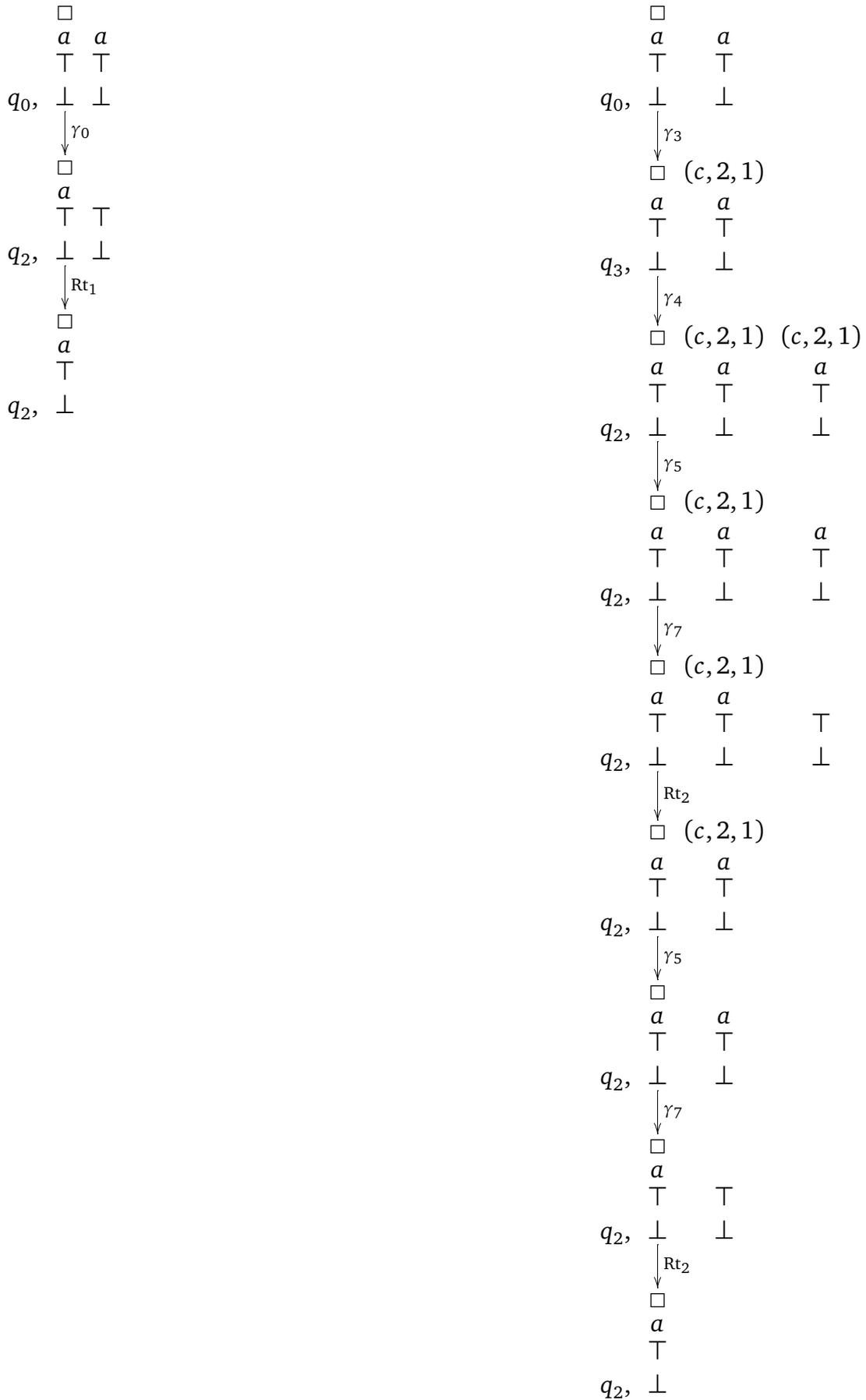

The rest of this section is concerned with the question ``How can we
determine $\ReturnFunc{k}(s)$ for some collapsible pushdown system
$\mathcal{S}$ using a finite automaton?''. The technical tools that we
use in order to answer this question are the notions of a \emph{return
simulator} and a \emph{simulation of a return}. 
We start with an informal description. Afterwards, we precisely
define these notions. 
A return simulator is a copy
of the pushdown system $\mathcal{S}$ enriched by transitions that
simulate each return of $\Pop{1}(s)$ in one transition. The simulation
of a return from $s$ to $\Pop{2}(s)$ is a return of this return
simulator from the special stack 
\begin{align*}
  s':=\bot\top \TOP{1}(s)\Box:
  \bot\top\TOP{1}(s)
\end{align*}
to $\Pop{2}(s')$.  
$\top$ is a new symbol
representing $\TOP{2}(\Pop{1}(s))$ and $\Box$ is a symbol not
occurring in the transitions of the return simulator. $\Box$ is used to
stop the computation once we reached $\Pop{2}(s')$. This guarantees
that any run from $s'$ to $\Pop{2}(s')$ is a return. 
Figure \ref{fig:SimulationsOfReturns} shows the simulations of the
return $\pi_1$ and $\pi_3$ from figures \ref{fig:FourthReturnExample}
and \ref{fig:FifthReturnExample}. 

Before we introduce simulations and simulators formally, we want to
explain the connection between a run and its 
simulation. For this purpose, we fix some notation. 
Let $\rho$ be some return and $\rho'$ its simulation (which is also a
return). 
According to Lemma \ref{Lemma_Return_Decomposition_For_Simulation},
there is a sequence 
\begin{align*}
  0=j_0 < i_1 < j_1 < \dots j_n\leq i_{n+1} = \length(\rho)-1  
\end{align*}
such that, for all $1\leq k \leq n$, 
$s\prefixeq \rho{\restriction}_{[j_{k-1},  i_k]}$ and 
$\rho{\restriction}_{[i_k+1,j_k]}$ is a return from some stack with
topmost word $\TOP{2}(\Pop{1}(s'))$. 
Analogously,  there is a sequence 
\begin{align*}
  0=j'_0  < i'_1 < j'_1 < \dots j'_n\leq i'_{n+1} =
  \length(\rho')-1  
\end{align*}
such that, for $1\leq k \leq n$, 
$s'\prefixeq \rho'{\restriction}_{[j'_{k-1},  i'_k]}$ and 
$\rho'{\restriction}_{[i'_k+1,j'_k]}$ is a return from some stack with
topmost word $\bot\top = \TOP{2}(\Pop{1}(s'))$. 
The run $\rho$ and its simulation $\rho'$ are connected as follows:

$\rho'{\restriction}_{[j'_{k-1},  i'_k]}=\rho{\restriction}_{[j_{k-1},
  i_k]}[s/s']$, i.e., $\rho'{\restriction}_{[j'_{k-1},  i'_k]}$
copies  $\rho{\restriction}_{[j_{k-1},  i_k]}$ transition by
transition but starts in a different stack. 
Furthermore, 
$\rho'{\restriction}_{[i'_k+1,j'_k]}$ is a return of length $1$, i.e.,
it is a run that only consists of one $\Pop{2}$ operation. 

Thus, the simulation induces a decomposition of a run into those parts
prefixed by its initial stack $s$ and those parts that form returns which
are equivalent to a return from $\Pop{1}(s)$ to $\Pop{2}(s)$. Using
this decomposition we 
prove the inductive computability of $\ReturnFunc{k}(s)$ from
$\ReturnFunc{k}(\Pop{1}(s))$. We first define the notion of a return
simulator. Afterwards, we introduce the notion of a simulation of a
return. 

\begin{definition} \label{Def:ReturnSimulator}
  Let $k\in\N$ be a threshold,
  $\mathcal{S}=(Q,\Sigma,\Gamma,q_0,\Delta)$ a collapsible 
  pushdown system, and $w$ some word. 
  Let $s:=w{\downarrow_0}\Box:w{\downarrow_0}$. 
  The \emph{return simulator} with respect
  to $(k,\mathcal{S}, s)$, denoted by
  $\ReturnSimulator{\mathcal{S}}{s}{k}$,  
  is the tuple
  $(Q,\Sigma,\Gamma\cup\{\mathrm{Rt}_i: i\leq k\},q_0,\hat\Delta)$ 
  where the $\mathrm{Rt}_i\notin\Gamma$ are new edge labels and
  \begin{align*}
    \hat\Delta:=\Delta\cup \{ (q_1,\top,\mathrm{Rt}_i, q_2,\Pop{2}):
    i\leq \ReturnFunc{k}(\Pop{1}(w))(q_1, q_2)\}.
  \end{align*}
  We also use the notation $\ReturnSimulator{\mathcal{S}}{\rho}{k}$ for 
  $\ReturnSimulator{\mathcal{S}}{s}{k}$ if $\rho$ is a return starting at
  stack $s$.  
\end{definition}
\begin{remark} \label{RemarkReturnSimulator}
  Before we continue, let us make some remarks concerning this
  definition.
  \begin{itemize}
  \item The return simulator copies the behaviour of $\mathcal{S}$ as
    long as the topmost symbol of a stack is not $\top$.
  \item We consider $\top$ as an abbreviation for the word
    $\Pop{1}(\TOP{2}(s))$. A run starting in the stack 
    \begin{align*}
      s':=\bot\top \TOP{1}(s)\Box : \bot\top\TOP{1}(s)      
    \end{align*}
    reaches a stack with topmost symbol $\top$ if and only if the
    equivalent run that starts in $s$ reaches a stack with topmost
    word $\Pop{1}(\TOP{2}(s))$. Recall that a return of $s$ always
    continues with a return if it reaches a stack with topmost word
    $\Pop{1}(\TOP{2}(s))$.
  \item \label{Rem:ReturnsOfSimulator}
    By definition
    $\ReturnFunc{k}_{\ReturnSimulator{\mathcal{S}}{s}{k}}(\bot\top)$
    agrees with
    $\ReturnFunc{k}_{\mathcal{S}}(\Pop{1}(s))$. 
    On topmost symbol $\top$ the applicable transitions of the return
    simulator are only
    $\Pop{2}$ transitions. 
    Hence,  a return from 
    $(q_1,\Pop{1}(s'))$ to $(q_2, \Pop{2}(s'))$ consists by definition
    of only one $\Pop{2}$ transition. If
    $\ReturnFunc{k}_{\mathcal{S}}(\Pop{1}(s))(q_1,q_2)=i$, then there are $i$ such
    transitions which are labelled by $\mathrm{Rt}_1, \dots,
    \mathrm{Rt}_i$. Each of these induces exactly one return whence
    $\ReturnFunc{k}_{\ReturnSimulator{\mathcal{S}}{s}{k}}(\bot\top)(q_1,q_2)=i$.
  \end{itemize}
  The last two observations will lead to the result that  the number
  of returns of 
  $\mathcal{S}$ from  $s$ and the returns of the simulator from $s'$
  agree up to threshold $k$. 
\end{remark}

\begin{definition}
  Let $\mathcal{S}$, $s$, $s'$ and $k$ as in Definition \ref{Def:ReturnSimulator}. 
  We call any run of $\ReturnSimulator{\mathcal{S}}{s}{k}$ from
  $(q_1,s')$ to $(q_2,\Pop{2}(s'))$ a \emph{simulation of a return}
  from $(q_1,s)$ to $(q_2, \Pop{2}(s))$. 
\end{definition}

\begin{lemma}
  Let $\mathcal{S}, k, s$ and $s'$ be as in 
  Definition \ref{Def:ReturnSimulator}.  
  If $\rho$ is a simulation of a return
  from 
  $(q,s')$ to $(q',\Pop{2}(s'))$, then
  $\rho$ is in fact a return of the return simulator.
\end{lemma}
\begin{proof}
  Let $i$ be minimal in $\domain(\rho)$ such that the stack at
  $\rho(i)$ is a substack of $\Pop{2}(s')$. Due to
  $\lvert s' \rvert = 2$, $\rho(i) = \Pop{2}(s')$. Furthermore
  $\Sym(\Pop{2}(s'))= \Box$. Since
  $\ReturnSimulator{\mathcal{S}}{s}{k}$ does not contain any
  transition of the form $(q_1, \Box, \gamma, q_2, \op)$, $\rho(i)$
  cannot be extended. Thus, $i=\length(\rho)$. 

  Furthermore, $\TOP{2}(s')=\bot \top \TOP{1}(s)$. 
  If $\Lvl(s)=2$, then $\Lnk(s)=0$ by definition of $s$. 
  Thus, $\TOP{2}(s')$ does not contain any defined level $2$
  link. One easily concludes that $\rho$ is a return. 
\end{proof}

In the following, we justify the term 
\emph{simulation of a return}. To each simulation of the return
simulator $\ReturnSimulator{\mathcal{S}}{s}{k}$ with initial state $q$
and final  state $q'$ we associate a return from $(q,s)$ to 
$(q',\Pop{2}(s))$. 

\begin{definition}\label{DefSimToReturn}
  Let $\mathcal{S}$ be a collapsible pushdown system, $k\in\N$ a
  threshold, and $w$ some word. 
  Let 
  \begin{align*}
    &s:=w{\downarrow_0}\Box:w{\downarrow_0} \text{ and}\\
    &s':=\bot\top\TOP{1}(s)\Box :\bot\top\TOP{1}(s)
  \end{align*}
  We define a function $\SimToReturnSym{s}$ that maps every
  run $\rho'$ of $\ReturnSimulator{\mathcal{S}}{s}{k}$ from 
  $(q,s')$ to $( q',\Pop{2}(s'))$ to a return
  $\SimToReturnSym{s}(\rho')$ of $\mathcal{S}$ from $(q,s)$ to
  $(q',\Pop{2}(s))$. 

  In order to explain $\SimToReturnSym{s}$, we fix a run $\rho'$ of
  $\ReturnSimulator{\mathcal{S}}{s}{k}$ from  $(q_1,s')$ to
  $(q_2,\Pop{2}(s'))$.  
  Due to the previous lemma, $\rho'$ is a return. 
  We assume that it is a high return. 
  Let 
  \begin{align*}
    0=j_0<i_1<j_1<\dots< j_n< i_{n+1}    
  \end{align*}
  be the sequence
  corresponding to $\rho'$ according to Lemma
  \ref{Lemma_Return_Decomposition_For_Simulation}. Without loss of
  generality, we assume that $n>0$.
  
  We set $\pi'_k:=\rho'{\restriction}_{[j_{k-1},i_k]}$ for all 
  $1\leq k \leq n+1$ and
  $\rho'_k:=\rho'{\restriction}_{[i_k+1,j_k]}$ for all $1\leq k \leq n$.
  Now, we write $c'_k=(q'_k,s'_k)$ for the last configuration of $\pi'_k$ and
  $\hat c'_{k}=(\hat q'_k, \hat s'_k)$ for the configuration following
  $c'_k$ in $\rho'$.  
  Lemma \ref{Lemma_Return_Decomposition_For_Simulation} implies the
  following. 
  \begin{enumerate}
  \item $s'\prefixeq \pi'_k$ for all $1\leq k \leq
    n$.
  \item For all $1\leq k \leq n$, 
    $\TOP{2}(s'_k) = \bot \top \TOP{1}(s)$
    and  $\hat s'_k=\Pop{1}(s'_k)$. Thus, 
    $\TOP{2}(\hat s'_{k})=\bot\top$. 
  \item $\TOP{2}(s')\leq\TOP{2}(c'_{n+1})$ and $c'_{n+1}$ is connected
    to $\hat c'_{n+1}$ via a $\Pop{2}$ or $\Collapse$ of level $2$. In the
    latter case, $\TOP{2}(s') < \TOP{2}(c'_{n+1})$.
  \item For $1\leq i \leq n$, $\rho'_k$ is a return from $\hat c'_k$ to
    $\Pop{2}(\hat c'_k)$.  
    Thus, $\rho'_k$ is a return of a stack with topmost word 
    $\bot\top$. 
  \end{enumerate}
  We now define iteratively runs $\pi_i, \xi_i$ and $\rho_i$ whose
  composition then forms $\SimToReturnSym{s}(\rho')$. 
  
  Due to condition 1, $\pi_1:=\pi_1[s'/s]$ is well-defined. 
  $\pi_1$ ends with stack $s_1:=s'_1[s'/s]$. 
  Due to condition 2, $\TOP{2}(s'_1)=\TOP{2}(s')$ whence
  $\TOP{2}(s_1)=w$. This implies that $\TOP{1}(s'_1) =\TOP{1}(w) =
  \TOP{1}(s_1)$.   
  Furthermore, the transition connecting $\rho'(i_1)=(q_1',s_1')$ with
  $\rho'(i_1+1)=(\hat q_1', \hat s_1')$ performs  a $\Pop{1}$ or a
  collapse of level $1$.  
  Let $\xi_1$ be the run of length $1$ that applies this transition
  to the last configuration of $\pi_1$, i.e., $\xi_1$ is a run 
  $(q_1', s_1)\trans{} (\hat q_1', \Pop{1}(s_1))$. 
  
  Due to the observation in Remark \ref{Rem:ReturnsOfSimulator}, 
  the form of
  $\rho'_1$ is $(\hat q'_1,\hat s'_1)\trans{\mathrm{Rt}_{n}}
  (\tilde q'_1, \Pop{2}(\hat s'_1))$ for some
  $n \leq\ReturnFunc{k}(\Pop{1}(s))(\hat q'_1,\tilde q'_1)$. 
  Thus, we can define
  $\hat\rho_1$ to be the $n$-th return from $(\hat q'_1, \Pop{1}(s))$ 
  to $(\tilde q'_1, \Pop{2}(s))$ in length-lexicographic order.

  Recall that $\TOP{2}(\Pop{1}(s_1))=\Pop{1}(s)$ is the topmost word of
  the last configuration of $\xi_1$. Hence, there is a return $\rho_1$
  that is equivalent to $\hat\rho_1$ and starts in the last configuration
  of $\xi_1$. 

  $\rho_1$ ends in configuration $(\tilde q'_1,\Pop{2}(s_1))$ where
  $\tilde q'_1$ is by definition the state of the initial
  configuration of $\pi_2'$. Furthermore the stack of $\pi_2'$ is
  $\Pop{2}(s_1')$. Since $s'\prefixeq \pi_2'$, $s_1'[s'/s]$ is a 
  well-defined stack. Due to $s_1=s_1'[s'/s]$, we conclude that
  $\Pop{2}(s_1)=\Pop{2}(s_1')[s'/s]$. 

  Thus, we can repeat this construction for $2, 3, 4, \dots, n$ and
  obtain runs $\pi_k, \xi_k, \rho_k$ such that 
  $\pi_1\circ\xi_1\circ\rho_1\circ\pi_2\circ\xi_2\circ\rho_2\circ\dots\circ
  \pi_{n}\circ\xi_{n}\circ\rho_{n}$ is a well-defined run
  from $(q,s)$ to $\pi_{n+1}'(0)[s'/s]$. 
  
  We set $\pi_{n+1}:=\pi_{n+1}[s'/s]$. Due to condition 3, $w\leq
  s_{n+1}:=s'_{n+1}[s'/s]$ which is the last stack of $\pi_{n+1}$.  
  As in the cases $i\leq n$, it follows that
  $\TOP{1}(s_{n+1})=\TOP{1}(s'_{n+1})$. 
  Thus, the last transition $\delta$ of $\rho'$ is also applicable to
  the last 
  configuration of $\pi_{n+1}$. 
  By definition, $\delta$ connects the last configuration
  of $\pi_{n+1}'$ with $(q', \Pop{2}(s'))$. 
  Since $\lvert s_{n+1}' \rvert = \lvert s_{n+1} \rvert$ and
  $\Lvl(s_{n+1}) = \Lvl(s'_{n+1})$, the application of this transition
  to the last configuration of $\pi_{n+1}$ results in $(q', \tilde
  s)$ where $\tilde s$ is a stack  of width $1$ such that 
  $\tilde s = \Pop{2}^m(s_{n+1})$ for some $m\in\N$. But this is by
  definition  $(q', w{\downarrow_0}\Box)$. Let $\xi_{n+1}$ be the
  run that applies $\delta$ to the last configuration of $\pi_{n+1}$. 
  We define 
  \begin{align*}
    \SimToReturnSym{s}(\rho'):=\pi_1 \circ \xi_1 \circ \rho_1
    \circ \pi_2 \circ \xi_2 \circ \rho_2 \circ \dots \circ \pi_n
    \circ \xi_n \circ \rho_n \circ \pi_{n+1} \circ \xi_{n+1}.     
  \end{align*}
  We say $\SimToReturnSym{s}(\rho')$ is the return 
  \emph{simulated by  $\rho'$}.  
\end{definition}
\begin{remark}
  For low returns $\rho'$,
  $\SimToReturnSym{s}(\rho')$ is defined completely analogous. We define
  $\pi_i, \xi_i$, and $\rho_i$ for all $1\leq i \leq n$ as
  before. Then 
  \begin{align*}
    \SimToReturnSym{s}(\rho'):=
    \pi_1 \circ \xi_1 \circ \rho_1
    \circ \pi_2 \circ \xi_2 \circ \rho_2 \circ \dots \circ \pi_n
    \circ \xi_n \circ \rho_n.     
  \end{align*}
\end{remark}

\begin{lemma}
  Let $s= w{\downarrow_0}\Box:w{\downarrow_0}$ and $\rho'$ a
  simulation of a return as in the previous definition. Then 
  $\SimToReturnSym{s}(\rho')$ is a return from $s$ to $\Pop{2}(s)$. 
\end{lemma}
\begin{proof}
  By definition, $\SimToReturnSym{s}(\rho')$ is a run from $s$ to
  $\Pop{2}(s)$ that does not pass any substack of $\Pop{2}(s)$ before
  its final configuration. 
  If its last operation is $\Pop{2}$ we are done. 
  
  Otherwise, the last operation is a $\Collapse$ of level $2$.
  Then we distinguish the following
  cases:

  First consider the case that $\rho'$ is a high return. 
  Recall
  that by definition of $\pi_{n+1}$, we have that $\TOP{2}(s)$ is a
  proper prefix of the topmost word of the last stack of
  $\pi_{n+1}$. But then the use of the last collapse in
  $\SimToReturnSym{s}(\rho')$ satisfies the restrictions from the 
  definition of a return. 
  
  Now, consider the case that $\rho'$ is a low return. By definition, 
  $\SimToReturnSym{s}(\rho')$ ends with $\rho_n$. But $\rho_n$ was
  defined to be a return from $\Pop{1}(s)$ to $\Pop{2}(s)$. 
  Thus, $\SimToReturnSym{s}(\rho')$ is a return due to 
  Lemma \ref{Lemma_ReturnDecomposition}. 
\end{proof}

\begin{lemma} \label{LemmaSimToReturn}
  Let $s=w{\downarrow_0}\Box:w{\downarrow_0}$  as in 
  Definition \ref{DefSimToReturn}. Then  
  $\SimToReturnSym{s}$ is injective. 
\end{lemma}
\begin{proof}[Proof (Sketch).]
  
  The proof is by contradiction. Assume that there are two runs
  $\rho_1'$ and $\rho_2'$ such that $\rho_1'\neq\rho_2'$. 
  We write $\rho_i:=\SimToReturnSym{s}(\rho_i)$ for $i\in\{1,2\}$. 
  
  Then there is a minimal $i'\in\domain(\rho_1')$ such that the
  transition $\delta_1'$ at position $i'$ in $\rho_1'$ is not the
  transition $\delta_2'$ at position $i'$ in $\rho_2'$. 
  Set
  $\pi':=\rho_1'{\restriction}_{[0,i']}=\rho_2'{\restriction}_{[0,i']}$. 
  Now, $\pi'$ induces a
  common initial segment $\pi$ of 
  $\rho_1$ and 
  $\rho_2$ of length $i$. 
  By this we mean that $\rho_1{\restriction}_{[0,i]} =
  \rho_2{\restriction}_{[0,i]}$ and that $\delta_1'$ and $\delta_2'$ 
  determine the transition of $\rho_1$ and $\rho_2$ at position $i$. 
  
  We distinguish two cases.
  \begin{enumerate}
  \item Assume that
    $\TOP{1}(\rho_1(i))=\TOP{1}(\rho_2(i))\neq\top$.  
    By definition of $\rho_1$ and $\rho_2$, this implies that the
    transition at $i$ in $\rho_j$ is $\delta_j'$ for
    $j\in\{1,2\}$. Since $\delta_1'\neq\delta_2'$, this implies that 
    $\rho_1$ and $\rho_2$ differ in the transition applied at $i$
    whence $\rho_1\neq\rho_2$. 
  \item Otherwise, assume that
    $\TOP{1}(\rho_1(i))=\TOP{1}(\rho_2(i))=\top$. 
    Then we directly conclude that \mbox{$\delta_1'=(q, \top,
      \mathrm{Rt}_{j_1}, q'_1, 
      \Pop{2})$} and  
    $\delta_2'=(q,\top, \mathrm{Rt}_{j_2}, q'_2, \Pop{2})$ where
    either  $q'_1\neq q'_2$ or $j_1\neq j_2$. 
    
    If $q'_1\neq q'_2$, then there are  $i_1>i$ and $i_2>i$ 
    such that $\rho_1{\restriction}_{[i,i_1]}$ is a return from state
    $q$ to state $q'_1$ while 
    $\rho_2{\restriction}_{[i,i_i]}$ is a return from state $q$ to
    state $q'_2$. Thus, $\rho_1 \neq \rho_2$. 

    Otherwise, $j_1\neq j_2$ and $q'_1 = q'_2$. By definition, 
    $\rho_1$ continues with the $j_1$-th return from $\rho_1(i)$ to
    $(q'_1, \Pop{2}(\rho_1(i)))$ and $\rho_2$ continues with the
    $j_2$-th return from $\rho_2(i)=\rho_1(i)$ to $(q'_1,
    \Pop{2}(\rho_1(i)))$. Since $j_1\neq j_2$, these returns differ
    whence the runs $\rho_1$ and $\rho_2$ differ. \qedhere
  \end{enumerate}
\end{proof}

The last fact  that we prove about $\SimToReturnSym{s}$ is a
characterisation of its image. 
Consider a run $\rho$ of $\ReturnSimulator{\mathcal{S}}{s}{k}$ in the
domain of $\SimToReturnSym{s}$. 
By definition of $\SimToReturnSym{s}$, $\SimToReturnSym{s}(\rho)$ is a
return from $s$ to $\Pop{2}(s)$
that satisfies the following restriction: 
let $\rho'$ be a subrun of $\SimToReturnSym{s}$
that is a return from some stack $s'$ with
$\TOP{2}(s') = \TOP{2}(\Pop{1}(s))$. Then $\rho'$ is one of the
$k$ smallest returns of $s'$ (with respect to length-lexicographic
order). 

The following lemma shows that this condition already defines the
image of $\SimToReturnSym{s}$. 
We only state the lemma for high returns, but for low returns the
analogous statement holds. 

\begin{lemma} \label{Lem:SimToRetUsesSmall}
  Let $\mathcal{S}$ be some collapsible pushdown system, $s$ some
  stack of the form \mbox{$s=w{\downarrow_0}\Box:w{\downarrow_0}$} and $k$
  some threshold.  
  Furthermore, let $\rho$ be a high return from $(q,s)$ to
  $(q',\Pop{2}(s))$ 
  of $\mathcal{S}$ and
  let 
  \begin{align*}
    0=j_0< i_1 < j_1 < \dots < i_n < j_n < i_{n+1}    
  \end{align*}
  be the
  sequence corresponding to $\rho$ according to Lemma 
  \ref{Lemma_Return_Decomposition_For_Simulation}. 
  Let $\rho_m:=\rho{\restriction}_{[i_m,j_m]}$ for each $1\leq m \leq
  n$. $\rho_m$ is a return from some $(q_m, s_m)$ to $(q_m',
  \Pop{2}(s_m))$. 
  If for all $1\leq m \leq n$, 
  $\rho_m$ is one of the $m$ length-lexicographically smallest returns
  from $(q_m, s_m)$ to $(q_m', \Pop{2}(s_m))$,
  then there is a run
  $\rho'$ of $\ReturnSimulator{\mathcal{S}}{s}{k}$ such that 
  $\SimToReturnSym{s}(\rho')=\rho$. 
\end{lemma}
\begin{proof}
  Let $\rho$ be a high return from $s$ to $\Pop{2}(s)$ satisfying the
  properties required in the 
  lemma. For each $1\leq m\leq n+1$ let $\delta_m$ be the transition
  connecting $\rho(i_m)$ and $\rho(i_m+1)$. 
  Furthermore, for $1\leq m\leq n$ let $l_m$ be the number such that
  $\rho{\restriction}_{[i_m+1,j_m]}$ is the $l_m$-th return from
  $\rho(i_m+1)$ to $\rho(j_m)$ in length-lexicographic order. 
  Set $s':=\bot\top \TOP{1}(s) \Box: \bot\top \TOP{1}(s)$.
  For $1\leq m \leq n+1$, set 
  $\pi'_m:=\rho{\restriction}_{[j_{m-1}, i_m]}[s/s']$. 
  
  We define $\rho'$ to be the run 
  \begin{align*}
    \rho':=\pi'_1 \circ \xi'_1  \circ \rho'_1 \circ \pi'_2 \circ \dots
    \circ \rho'_n \circ \pi'_{n+1} \circ \xi'_{n+1}    
  \end{align*}
  where $\xi'_m$ applies
  $\delta_m$ to the last configuration of $\pi'_m$ and $\rho'_m$
  applies an $\mathrm{Rt}_{l_m}$-labelled transition to the last
  configuration of $\xi'_{m}$. 

  It is now easy to check that $\rho'$ is a well defined run of
  $\ReturnSimulator{\mathcal{S}}{s}{k}$ from 
  $s'$ to $\Pop{2}(s)$ and that
  $\rho= \SimToReturnSym{s}(\rho')$. 
\end{proof}

A corollary of the previous lemma is that there are at least as many
returns of a pushdown system $\mathcal{S}$ from $(q,s)$ to
$(q',\Pop{2}(s))$  as there are runs of
$\ReturnSimulator{\mathcal{S}}{s}{k}$ from $(q,s')$ to
$(q',\Pop{2}(s'))$ for 
$s'=\bot\top \TOP{1}(s)\Box:\bot\top\TOP{1}(s)$. 

In fact, we want to prove that these two numbers agree up to
threshold $k$.  
We obtain this result as a corollary of the following lemma.
Again we only formulate the lemma for high returns, but the corresponding
statement for low returns is proved analogously. 

\begin{lemma} \label{Lem:LargeSubreturnImpliesManyReturns}
  Let $\mathcal{S}$ be a collapsible pushdown system and
  $s=w{\downarrow_0}\Box:w{\downarrow_0}$ for some word $w$. 
  Let $\rho$ be a return of $\mathcal{S}$ from $(q,s)$ to $(q',\Pop{2}(s))$. 
  Let 
  \begin{align*}
    0=j_0<i_1< j_1 < \dots < i_n < j_n < i_{n+1}    
  \end{align*}
  be the sequence
  corresponding to $\rho$ according to Lemma 
  \ref{Lemma_Return_Decomposition_For_Simulation}. 
  If there is a $1\leq k \leq n$ such that
  $\rho_k:=\rho{\restriction}_{[i_k+1,j_k]}$  is not one of the
  minimal $k$ returns from $\rho(i_k+1)$ to $\rho(j_k)$, then there
  are more than $k$ returns of $\mathcal{S}$ from $(q,s)$ to 
  $(q', \Pop{2}(s))$. 
\end{lemma}
\begin{proof}
  Let $1\leq k \leq n$ be a number such that
  $\rho_k:=\rho{\restriction}_{[i_k+1,j_k]}$ satisfies the
  requirements of the lemma. 
  Then there are $k$ returns from $\rho(i_k+1)$ to $\rho(i_j)$ that
  are length-lexicographically smaller than $\rho_k$. 
  Now, let $\hat\rho_1, \hat\rho_2, \dots, \hat\rho_k$ be an enumeration of
  these runs. For $1\leq i \leq k$, the run
  $\pi_i:=\rho{\restriction}_{[0,i_k+1]} \circ \hat \rho_i \circ
  \rho{\restriction}_{[j_k, i_{n+1}+1]}$ is a return from $(q,s)$ to
  $(q', \Pop{2}(s))$. The $\pi_i$ are pairwise distinct and distinct
  from $\rho$. Thus, there are at least $k+1$ returns from $(q,s)$ to
  $(q', \Pop{2}(s))$. 
\end{proof}

As a direct corollary of the previous two lemmas, we obtain that the
runs of the return simulator 
and the returns from a configuration $(q,s)$ to $(q',\Pop{2}(s))$
agree up to threshold $k$. 

\begin{corollary} \label{Cor:InductiveComputability}
  Let $\mathcal{S}$ be a collapsible pushdown system, $k\in\N$ some
  threshold and $w$ some word. 
  For $s:=w{\downarrow_0}\Box:w{\downarrow_0}$ and 
  $s':=\bot\top \TOP{1}(s) \Box: \bot\top \TOP{1}(s)$ and
  for all $q,q'\in Q$ let $M_{q,q'}$ be the set of runs of
  $\ReturnSimulator{\mathcal{S}}{s}{k}$ from $(q,s')$ to $(q',
  \Pop{2}(s'))$. 
  For all $q,q'\in Q$,
  \begin{align*}
    \ReturnFunc{k}_{\mathcal{S}}(w)(q,q') = 
    \min\{k, \lvert M_{q,q'} \rvert\} =
    \ReturnFunc{k}_{\ReturnSimulator{\mathcal{S}}{s}{k}}(\TOP{2}(s'))(q,q').
  \end{align*}
\end{corollary}
The last corollary shows that we can count simulations of a return
simulator in order to calculate
$\ReturnFunc{k}_{\mathcal{S}}(s)$. Now, we use this 
result in order to obtain a proof of Proposition \ref{Prop:AutomatonForReturns}.
Recall that this proposition asserts that there is a finite automaton
$\mathcal{A}_{\mathrm{ret}}$ that calculates $\ReturnFunc{k}(s)$ for
each stack $s$ on input $\TOP{2}(s){\downarrow_0}$. 

\begin{proof}[Proof of Proposition \ref{Prop:AutomatonForReturns}]
  In the following, we define the finite automaton 
  \begin{align*}
    \mathcal{A}_{\mathrm{ret}}:=(Q_{\mathrm{ret}}, 
    \Sigma \cup (\Sigma \times\{2\}\times\{0\}), a_0,
    \Delta_{\mathrm{ret}}). 
  \end{align*}
  Let $Q_{\mathrm{ret}}:=\{a_0\}\cup \{0, 1, \dots, k\}^{Q\times Q}$
  where $Q$ is the set of states of
  $\mathcal{S}$ and $a_0$ is an extra initial state distinct from 
  all other states. Thus, beside the initial state  all functions from
  $Q\times Q$ to $\{0, 1,\dots,k\}$ are states of
  $\mathcal{A}_{\mathrm{ret}}$.

  We define $\Delta_{\mathrm{ret}}$ in such a way that the run of
  $\mathcal{A}_{\mathrm{ret}}$ on some word $w=\TOP{2}(s){\downarrow_0}$
  ends in a state $a=\ReturnFunc{k}(s)$.
  We compute the transitions of
  $\mathcal{A}_{\mathrm{ret}}$ iteratively. 

  We start with the transitions from the
  initial state. Recall that all words occurring in some stack start
  with the 
  letter $\bot$. Thus, the only transition at $a_0$ should be of the
  form $(a_0, \bot, a)$ where
  $a$ must satisfy  $a=\ReturnFunc{k}(\bot)=\ReturnFunc{k}( [\bot
  \Box: \bot])$. 
  Due to Lemma \ref{LemmaCountNumberOfRuns}, the value of $a$ is
  computable. 

  Now, we repeat the following construction.
  Assume that for all reachable states
  \mbox{$a\in Q_{\mathrm{ret}}\setminus\{a_0\}$} every path from $a_0$ to $a$ is
  labelled by some word $w$ such that $a=\ReturnFunc{k}(w)$. 

  For each $\sigma\in\Sigma$, we want to compute the
  value $a'=\ReturnFunc{k}(w\sigma)$.
  Let \mbox{$s:=w{\downarrow_0} \sigma \Box :w{\downarrow_0} \sigma$.} 
  Recall that $\ReturnSimulator{\mathcal{S}}{s}{k}$ is computable from
  $\mathcal{S}, a=\ReturnFunc{k}(w)$, and $k$. 
  Due to Lemma \ref{LemmaCountNumberOfRuns}, we can count the number
  $i_{q_1,q_2}$ 
  of runs of $\ReturnSimulator{\mathcal{S}}{s}{k}$ from 
  $(q_1,\bot\top\sigma\Box:\bot\top\sigma)$ to $(q_2,[\bot\top\sigma\Box])$
  for each pair $q_1,q_2\in Q$ up to threshold $k$. 
  Finally, Corollary \ref{Cor:InductiveComputability} shows that 
  \mbox{$i_{q_1,q_2}=\ReturnFunc{k}(w\sigma)(q_1,q_2)$}. 

  Thus, $a':=\ReturnFunc{k}(w\sigma)$ is computable and we add the
  transition $(a,\sigma,a')$ to $\Delta_{\mathrm{ret}}$. 
  By induction hypothesis and by Corollary
  \ref{Cor:InductiveComputability}, all nonempty paths to some state
  $\hat a$ are now labelled by words $v$ such that
  $\hat a=\ReturnFunc{k}(v)$ (this can be proved by induction on the
  length of the path). 

  For words of the form $w(\sigma,2,0)$ the transitions
  $(a,(\sigma,2,0),a')$ are defined completely analogous. 

  After finitely many iterations of this process, we cannot add any new
  transitions to $\mathcal{A}_{\mathrm{ret}}$. 
  Then the construction of $\mathcal{A}_{\mathrm{ret}}$ is finished. 

  We claim that the resulting automaton $\mathcal{A}_{\mathrm{ret}}$
  calculates $\ReturnFunc{k}(s)$ on input
  $\TOP{2}(s){\downarrow_0}$ for every  stack $s$. 

  The claim is proved by contradiction. Assume that there is some
  stack $s$ such 
  that there is no run of $\mathcal{A}_{\mathrm{lp}}$ on
  $w:=\TOP{2}(s){\downarrow_0}$. By minimality there is a run on
  $\Pop{1}(w)$. Now, we could add 
  a transition from the final state of this run which is labelled by
  $\TOP{1}(w)$. This contradicts the assumption that we added all
  possible transitions. 

  Thus, there is a run of $\mathcal{A}_{\mathrm{lp}}$ on
  $w:=\TOP{2}(s){\downarrow_0}$ for all stacks $s$. 
  By construction, the run on $w$ calculates $\ReturnFunc{k}(w)$. 
\end{proof}

We conclude this section by proving Proposition
\ref{Prop:FuncBoundReturnLengthLemma}.
Recall that this proposition asserts the existence of a function
$\FuncBoundReturnLength{\mathcal{S}}{k}$ that bounds the length of the
shortest returns of every stack. 
We first define
$\FuncBoundReturnLength{\mathcal{S}}{k}$, then we prove the
properties asserted in the lemma.

  Let $\mathcal{S}$ be some collapsible pushdown system and let 
  $\mathcal{A}_{\mathrm{ret}}$ be the corresponding finite automaton
  that calculates the returns of $\mathcal{S}$ up to threshold $k$. 
  Recall that $\Delta_{\mathrm{ret}}$  denotes the
  transition relation of $\mathcal{A}_{\mathrm{ret}}$.

  Recall that for each $\delta=(a,\tau,b)\in \Delta_{\mathrm{ret}}$,
  it holds that 
  $\tau\in \Sigma\cup(\Sigma\times\{2\}\times\{0\})$ and 
  $a,b:Q\times Q \to \{0,1,\dots, k\}$ 
  are functions such that
  there is some word $w_\delta\in (\Sigma\cup
  (\Sigma\times\{2\}\times\{0\}))^+$ with $\ReturnFunc{k}(w_\delta)=a$
  and $\ReturnFunc{k}(w_\delta\tau)=b$. 
  In the following, we fix a $w_\delta$ for each
  $\delta\in \Delta_{\mathrm{ret}}$. 
  For each $w_\delta$, we define the stack 
  $s'_\delta:=\bot\top\tau\Box:\bot\top\tau$. 
  Due to Corollary \ref{Cor:InductiveComputability}, there are (up to
  threshold $k$) $b(q_1, q_2)$ many simulations of returns from $(q_1, s'_\delta)$
  to $(q_2, 
  \Pop{2}(s'_\delta))$. Using Lemma \ref{LemmaCountNumberOfRuns} we
  can compute the $b(q_1,q_2)$ many lexicographically smallest such
  simulations. We call these $\rho_1^{\delta,q_1,q_2},
  \rho_2^{\delta,q_1,q_2}, \dots,
  \rho_{b(q_1,q_2)}^{\delta,q_1,q_2}$. 

  Let 
  \begin{align*}
    &l_{q_1,q_2}^\delta:=\max\left\{ \length(\rho_1^{\delta,q_1,q_2}), 
    \length(\rho_2^{\delta,q_1,q_2}), \dots, 
    \length(\rho_{b(q_1,q_2)}^{\delta,q_1,q_2})\right\} \text{ and}\\
    &\#\top_{q_1,q_2}^\delta:=
    \max\left\{\lvert\{j\in\domain(\rho_i^{\delta,q_1,q_2}):
    \Sym(\rho_i^{\delta,q_1,q_2}(j))=\top \}\rvert : 1\leq i \leq
    b(q_1,q_2) \right\}      
  \end{align*}
  be the maximal length of any of these simulations and  the 
  maximal number of occurrences of $\top$ as topmost symbol in any of these
  returns, respectively. 
  Now, set 
  \begin{align*}
    &l:=\max\{ l_{q_1,q_2}^\delta: q_1,q_2\in Q,
    \delta\in\Delta_{\mathrm{ret}}\}\text{ and}\\
    &\#\top:=\max\{ \#\top_{q_1,q_2}^\delta:q_1,q_2\in Q,
    \delta\in\Delta_{\mathrm{ret}}\}. 
  \end{align*}
\begin{definition} 
  We define 
  \begin{align*}
    &\FuncBoundReturnLength{\mathcal{S}}{k}:\N \to \N \text{ by}\\
    &\FuncBoundReturnLength{\mathcal{S}}{k}(0)=0\text{ and}\\
    &\FuncBoundReturnLength{\mathcal{S}}{k}(n+1):=
    l + \#\top \cdot \FuncBoundReturnLength{\mathcal{S}}{k}(n).    
  \end{align*}
\end{definition}
\begin{remark} \label{Rem:BoundingFunctionsMotivation}
  The following idea underlies this definition. Assume that there is
  some word $w$ of length $n-1$ such that the length of the shortest
  $\ReturnFunc{k}(w)(q_1,q_2)$ returns from $(q_1, w\Box:w)$ to
  $(q_2,w\Box)$  is bound by
  $\FuncBoundReturnLength{\mathcal{S}}{k}(n-1)$. 
  Furthermore, let $s$ be a stack such that $w=\TOP{2}(\Pop{1}(s))$. 
  
  Due to the definition of $l$, the lexicographically smallest
  $\ReturnFunc{k}(s)(q,q')$ many returns from $(q,s)$ to
  $(q',\Pop{2}(s))$ have simulations of length at most $l$. 
  
  Now, $\SimToReturnSym{s}$ translates these simulations into
  $\ReturnFunc{k}(s)(q,q')$ many returns by copying all transitions
  one by one except for transitions on topmost symbol $\top$. The
  latter are  
  replaced by
  lexicographically small returns equivalent to those from 
  $(q_1, w\Box:w)$ to $(q_2,w\Box)$. Since this replacement happens at
  at most $\#\top$ many positions, we obtain $\ReturnFunc{k}(s)(q,q')$
  many returns from $(q,s)$ to $(q',\Pop{2}(s))$ of length at most
  $l + \#\top \cdot \FuncBoundReturnLength{\mathcal{S}}{k}(n-1) = 
  \FuncBoundReturnLength{\mathcal{S}}{k}(n)$. 
\end{remark}

Next, we prove Proposition \ref{Prop:FuncBoundReturnLengthLemma}. 
\begin{proof}
  Recall that we have to show the following two properties of
  $\FuncBoundReturnLength{\mathcal{S}}{k}$.
  \begin{enumerate}
  \item 
    For each stack $s$, for all states $q_1,q_2$ and for
    $i:=\ReturnFunc{k}(s)(q_1,q_2)$, the
    length-lexicographically shortest returns 
    $\rho_1,\dots, \rho_i$ from
    $(q_1,s)$ to $(q_2,\Pop{2}(s))$ satisfy 
    $\length(\rho_j)\leq \FuncBoundReturnLength{\mathcal{S}}{k}(\lvert
    \TOP{2}(s)\rvert)$ for all $1\leq j \leq i$.
  \item If there is a return $\rho$ from $(q_1,s)$ to $(q_2, \Pop{2}(s))$ with
    $\length(\rho) > \FuncBoundReturnLength{\mathcal{S}}{k}(\lvert
    \TOP{2}(s)\rvert)$, then there are $k$ returns from $(q_1,s)$ to
    $(q_2, \Pop{2}(s))$ of length at most 
    $\FuncBoundReturnLength{\mathcal{S}}{k}(\lvert \TOP{2}(s)\rvert)$.  
  \end{enumerate}
  Note that the previous remark already contains a proof of the first
  part. The second part is proved by induction on $k$. 

  Let  $\rho$ be a return from $(q_1,s)$ to
  $(q_2,\Pop{2}(s))$ with  
  $\length(\rho) > \FuncBoundReturnLength{\mathcal{S}}{k}(\lvert
  \TOP{2}(s)\rvert)$. Then we conclude that $\ReturnFunc{k}(s)\geq
  1$. Due to the first statement, the lexicographically shortest return
  $\rho_1$ from $(q_1,s)$ to $(q_2, \Pop{2}(s))$ satisfies 
  $\length(\rho_1) \leq \FuncBoundReturnLength{\mathcal{S}}{k}(\lvert
  \TOP{2}(s)\rvert)$. Thus, $\rho_1\neq \rho$ and we conclude that 
  $\ReturnFunc{k}(s)\geq 2$ (if $k\geq 2$). 

  We can iterate this argument $k$ times and obtain 
  $\rho_1, \rho_2 \dots, \rho_k$ many short returns from 
  $(q_1, s)$ to  $(q_2, \Pop{2}(s))$ as desired. 
\end{proof}

\subsection{Computing Loops}
\label{sec:CompLoops}
This section investigates the computability of loops. In fact, it
lifts the results on returns to analogous results on loops. 
Again, we start by defining the functions we are interested in. 

\begin{definition} \label{Def:LoopFunc}
  Let $\mathcal{S}$ be some collapsible pushdown system of level $2$,
  $k\in\N$ some threshold and $s$ some stack.  
  We define 
  \begin{align*}
    &\LoopFunc{k}_{\mathcal{S}}(s):Q\times Q \rightarrow \{0,1,\dots,k\} \\
    &(q,q')\mapsto
    \begin{cases}
      i &\text{if there are exactly }i\leq k\text{ different loops of
        $\mathcal{S}$ 
        from }(q,s)\text{ to } (q',s)\\
      k & \text{otherwise.}
    \end{cases}
  \end{align*}
  This function maps $(q,q')$ to the number $i$ of loops from $(q,s)$ to
  $(q',s)$ if $i\leq k$ and it maps $(q,q')$ to $k$ otherwise. In
  this sense $k$ stands for the class of at least $k$ loops. 
  
  Analogously, we define 
  $\LowLoopFunc{k}_{\mathcal{S}}(s):Q\times Q \rightarrow
  \{0,1,\dots,k\}$ to be the 
  function that maps $(q,q')$ to the number $i$ of low loops from $(q,s)$ to
  $(q',s)$ if $i\leq k$ and that maps $(q,q')$ to $k$ otherwise. 

  Finally, we define 
  $\HighLoopFunc{k}_{\mathcal{S}}(s):Q\times Q \rightarrow
  \{0,1,\dots,k\}$ to be the 
  function that maps $(q,q')$ to the number $i$ of high loops from $(q,s)$ to
  $(q',s)$ if $i\leq k$ and that maps $(q,q')$ to $k$ otherwise. 
  
  If $\mathcal{S}$ is clear from the context, we omit it and write
  $\LoopFunc{k}$ for $\LoopFunc{k}_{\mathcal{S}}$, etc. 
\end{definition}

Analogously to the theory of returns, we want to show that
$\LoopFunc{k}$, $\HighLoopFunc{k}$, and $\LowLoopFunc{k}$ can be
calculated by a finite automaton. Furthermore, we also want to prove
bounds on the length of short loops analogously to Proposition
\ref{Prop:FuncBoundReturnLengthLemma}. 
We start by stating these two propositions.

\begin{proposition} \label{Prop:AutomatonForLoops}
  There is an algorithm that, given a collapsible pushdown system
  $\mathcal{S}$ of level $2$, computes a
  deterministic finite automaton  
  $\mathcal{A}_{\mathrm{loop}}$ that computes $\LoopFunc{k}(s:w)$ on
  input $w{\downarrow_0}$.

  In the same sense, there are  automata that compute
  $\HighLoopFunc{k}$ and $\LowLoopFunc{k}$. 
\end{proposition}

\begin{proposition} \label{Prop:FuncBoundLoopLengthLemma}  
  There is an algorithm that computes on input some $2$-\CPG
  $\mathcal{S}$ and a natural number $k$ a
  function
  $\FuncBoundLoopLength{\mathcal{S}}{k}:\N \to\N$ such that the
  following holds. 
  \begin{enumerate}
  \item 
    For every stack $s$, for $q_1,q_2\in Q$ and for
    $i:=\LoopFunc{k}(s)(q_1,q_2)$, the
    length-lexicographically shortest loops
    $\lambda_1,\dots, \lambda_i$ from
    $(q_1,s)$ to $(q_2,s)$ satisfy 
    \begin{align*}
       \length(\lambda_j)\leq
      \FuncBoundLoopLength{\mathcal{S}}{k}(\lvert 
      \TOP{2}(s)\rvert)      
    \end{align*}
    for all $1\leq j \leq i$.
  \item If there is a loop $\lambda$ from $(q_1,s)$ to $(q_2, s)$ with
    $\length(\lambda) > \FuncBoundLoopLength{\mathcal{S}}{k}(\lvert
    \TOP{2}(s)\rvert)$, then there are $k$ loops from $(q_1,s)$ to
    $(q_2, \Pop{2}(s))$ of length at most 
    $\FuncBoundLoopLength{\mathcal{S}}{k}(\lvert \TOP{2}(s)\rvert)$.  
  \end{enumerate}
  
  Analogously, there are functions 
  $\FuncBoundHighLoopLength{\mathcal{S}}{k}:\N \to\N$ and
  $\FuncBoundLowLoopLength{\mathcal{S}}{k}:\N \to\N$ that satisfy the
  same assertions but for the set of high loops or low loops,
  respectively. 
\end{proposition}

Before we prove these propositions, we present two corollaries of the
previous Proposition that play a crucial role in Section \ref{Chapter
  HONPT}. 

\begin{corollary} \label{Cor:GlobalBoundRun}
  Let $\mathcal{S}$ be some level $2$ collapsible pushdown
  system. Furthermore, let  $(q,s)$ be some configuration and 
  $\rho_1, \dots, \rho_n$ be pairwise distinct runs from the initial
  configuration to $(q,s)$. 
  There is a run $\hat\rho_1$ from the initial configuration to
  $(q,s)$ such that the following holds.
  \begin{enumerate}
  \item $\hat\rho_1\neq\rho_i$ for $2\leq i \leq n$ and
  \item $\length(\hat\rho_1)\leq 2 \cdot \width(s) \cdot \height(s) 
    ( 1+ \FuncBoundLoopLength{\mathcal{S}}{n}(\height(s)))$.
  \end{enumerate}
\end{corollary}
\begin{proof}
  If $\length(\rho_1)\leq
  2 \cdot \lvert s \rvert \cdot \height(s) 
  ( 1+ \FuncBoundLoopLength{\mathcal{S}}{n}(\height(s)))$, set
  $\hat\rho_1:=\rho_1$ and we are done. 
  Assume that this is not the case. 
  Due to Lemmas 
  \ref{LemmaNumberMilestones} and
  \ref{Lemma:Carayol05Milestones}, $\rho_1$ decomposes as
  \begin{align*}
    \rho_1 = \lambda_0 \circ \op_1 \circ \lambda_1
    \circ \dots \circ \lambda_{m-1} \circ \op_m \circ \lambda_{m}    
  \end{align*}
  where every $\lambda_i$ is a loop and every $\op_i$ is a run of
  length $1$ such that $m\leq 2\cdot \lvert s \rvert \cdot
  \height(s)$. Proposition \ref{Prop:FuncBoundLoopLengthLemma} implies
  the following: 
  If $\length(\lambda_i) >
  \FuncBoundLoopLength{\mathcal{S}}{n}(\height(s)))$, then there are
  $n$ loops from $\lambda(0)$ to $\lambda(\length(\lambda))$ of length 
  at most $\FuncBoundLoopLength{\mathcal{S}}{n}(\height(s)))$. 
  At least one of these can be plugged into the position of
  $\lambda_i$ such that the resulting run does not coincide with any
  of the $\rho_2, \rho_3, \dots, \rho_n$. In other words, there is
  some loop $\lambda_i'$ of length at most
  $\FuncBoundLoopLength{\mathcal{S}}{n}(\height(s)))$ such that
  \begin{align*}
    \hat\rho_1:=
    \lambda_0 \circ \op_1 \circ \lambda_1
    \circ \dots \circ \op_i \circ \lambda_i'\circ\op_{i+1}\circ
    \lambda_{i+1} \circ\dots \circ  \lambda_{m-1} \circ \op_m \circ
    \lambda_{m}     
  \end{align*}
  is a run to $(q,s)$ distinct from $\rho_2, \rho_3 \dots, \rho_n$ and
  shorter than $\rho_1$.
  Iterated replacement of large loops results in a run 
  $\rho_1'$ with the desired properties. 
\end{proof}

Now, we state a second corollary that is quite similar to the previous
one but deals with runs of a different form. 

\begin{corollary}\label{Cor:GlobalBoundRun2}
  Let $\hat\rho_1, \hat\rho_2, \dots \hat\rho_n$ be runs from the
  initial configuration to some configuration $(q,s)$. 
  Furthermore, let $w$ be some word and 
  $\rho_1, \rho_2, \dots \rho_n$ be runs from $(q,s)$ to $(q',s:w)$
  that do not visit proper substacks of $s$. 
  If $\hat\rho_1 \circ \rho_1, \hat\rho_2\circ\rho_2, \dots,
  \hat\rho_n\circ\rho_n$ are pairwise distinct, then there is a run
  $\rho_1'$ from $(q,s)$ to $(q',s:w)$ that satisfies the following.
  \begin{enumerate}
  \item $\rho_1'$does not visit a proper
    substack of $s$,
  \item   
    $\length(\rho_1') \leq 2\cdot \height(s:w)  \cdot
    ( 1+ \FuncBoundLoopLength{\mathcal{S}}{n}(\height(s:w)))$, and
  \item $\hat\rho_1\circ\rho_1'$ is distinct from each
    $\hat\rho_i\circ\rho_i$ for $2\leq i \leq n$. 
  \end{enumerate}
\end{corollary}
\begin{proof}
  It is straightforward to see that $\rho_1$ decomposes as 
  \begin{align*}
    \rho_1= \lambda_0 \circ \op_1 \circ \lambda_1
    \circ \dots \circ \lambda_{m-1} \circ \op_m \circ \lambda_{m}    
  \end{align*}
  where every $\lambda_i$ is a loop and every $\op_i$ is a run of
  length $1$ such that $m\leq 2\cdot \height(s:w)$. 
  We then proceed completely analogous to the previous corollary.
\end{proof}

We now come to the proofs of the main propositions on loops. 
The proofs of these two propositions are analogous to the proofs for
the return case. The reader who is not interested in these rather
technical proofs should skip
the rest of this section  and continue reading Section 
\ref{Chapter_AutomaticStructures}.

We now prepare the proofs of the two main propositions on
loops. Analogously to the return case, the first important observation
is that $\LoopFunc{k}(s)$, $\HighLoopFunc{k}(s)$, and
$\LowLoopFunc{k}(s)$ only depend on the symbols and link levels of the
topmost word of the stack $s$. In order to show this, we first define
the notion of equivalent loops analogously to the notion of equivalent
returns in Definition \ref{Def:ReturnEquivalent}. 

\begin{definition}
  Let $s, s'$ be stacks such that $\TOP{2}(s){\downarrow_0} =
  \TOP{2}(s'){\downarrow_0}$. Let $\lambda$ be a loop from $(q_1,s)$ to
  $(q_2,s)$ and $\lambda'$ be a loop from $(q_1, s')$ to
  $(q_2, s')$. 
  We say $\lambda$ and $\lambda'$ are \emph{equivalent} loops if they
  consist of the same sequence of transitions.
\end{definition}

The crucial observation is that different stacks whose topmost words
agree on their symbols and link levels have the same loops modulo
this equivalence relation. 

\begin{lemma}\label{LoopDependTopword}
  Let $s$ and $s'$ be stacks such that
  $\TOP{2}(s){\downarrow_0} = \TOP{2}(s'){\downarrow_0}$. If $\lambda$ is a
  loop from $(q_1, s)$ to $(q_2, s)$ then there is an
  equivalent loop $\lambda'$ from $(q_1, s')$ to $(q_2, s')$. 
\end{lemma}
\begin{remark}
  The proof of this lemma is analogous to the proof of Lemma
  \ref{ReturnDependTopword}.  
  Furthermore, it is straightforward to see that a low loop can only
  be equivalent to a low loop (analogously, a high loops can only
  be equivalent to a high loop). 
\end{remark}
Since the lemma shows that loops of a given
stack only depend on its topmost word, it is a meaningful concept to
speak about the loops of some word. 

\begin{definition} \label{Def:loopsofWord}
  For $w$ some word, let $\LoopFunc{k}(w)$ be 
  $\LoopFunc{k}(s)$ for some stack $s$ with $\TOP{2}(s)=w$. 
  Analogously, we define the notions $\HighLoopFunc{k}(w)$ and 
  $\LowLoopFunc{k}(w)$. 
\end{definition}

The next step towards the proof of our main propositions is a
characterisation of $\LoopFunc{k}(w)$ in terms of
$\LoopFunc{k}(\Pop{1}(w))$ and $\ReturnFunc{k}(\Pop{1}(w))$
analogously to the result of Lemma
\ref{Lemma_Return_Decomposition_For_Simulation} for returns. 
We do this in the following three lemmas. First, we present a
unique decomposition of loops into high and low loops. Afterwards, we
characterise low loops and high loops. 

\begin{lemma} \label{Lem:LoopDecomposition}
  Let $\lambda$ be a loop from $(q,s)$ to $(q',s)$. 
  $\lambda$ is either a high loop or it has a unique decomposition as
  $\lambda = \lambda_0 \circ \lambda_1 \circ \lambda_2$ where
  $\lambda_0$ and $\lambda_2$ are
  high loops and $\lambda_1$ is a low loop.
\end{lemma}
\begin{proof}
  Assume that $\lambda$ is no
  high loop. Since it is a loop, it
  visits $\Pop{1}(s)$ at some position. 
  Let $i\in\domain(\lambda)$ be the minimal position just
  before the first occurrence of
  $\Pop{1}(s)$  and $j\in\domain(\lambda)$ be the position
  directly after the last occurrence of $\Pop{1}(s)$. It is
  straightforward to see that 
  $\lambda{\restriction}_{[i+1,j-1]}$ is
  by definition a loop of $\Pop{1}(s)$ and the initial and final part
  of $\lambda$ are loops of $s$. We conclude by noting that
  $\lambda{\restriction}_{[i,j]}$ is then a low loop of $s$. 
\end{proof}
\begin{remark}\label{Rem:LoopsDeterminedByLowAndHighOnes}
  An important consequence of this lemma is the fact that
  $\LoopFunc{k}(s)$ is determined by $\HighLoopFunc{k}(s)$ and
  $\LowLoopFunc{k}(s)$. 
  $\LoopFunc{k}(s)(q,q')$ counts  the high loops from $(q,s)$ to
  $(q',s)$ and those loops that consists of a high loop from $(q,s)$
  to $(\hat q, s)$ followed by a low loop from $(\hat q, s)$ to $(\hat
  q', s)$ followed by a loop from $(\hat q', s)$ to $(q', s)$. 
  Thus, writing
  $H(s)$ for $\HighLoopFunc{k}(s)$ and $L(s)$ for
  $\LowLoopFunc{k}(s)$, we obtain that 
  \begin{align*}
    \LoopFunc{k}(s)(q,q') = 
    \min\left\{k, H(s)(q,q') + 
    \sum_{\hat q, \hat q'\in Q}  H(s)(q,\hat q) \cdot
    L(s)(\hat q, \hat q') \cdot 
    H(s)(\hat q', q') \right\}.
  \end{align*}
\end{remark}

In the following, we first explain how low loops depend on the loops
of smaller stacks, afterwards we explain how high loops depend on
returns of smaller stacks. 

\begin{lemma} \label{Lem:LowLoopDecomposition}
  Let $\lambda$ be a low loop starting and ending in stack $s$. Then
  $\lambda{\restriction}_{[1,\length(\lambda)-1]}$ is a loop starting
  and ending in $\Pop{1}(s)$. The operation at $0$ is a $\Pop{1}$ or a
  $\Collapse$ of level $1$. The operation at $\length(\lambda)-1$ is
  a $\Push{\sigma}$ where $\TOP{1}(s)=\sigma\in\Sigma$. 
\end{lemma}
\begin{proof}
  Let $\TOP{1}(s)=\sigma$. 
  The lemma follows directly from the observation that a $\Push{\sigma}$
  transition followed by a loop of $s$ followed by a $\Pop{1}$ or
  $\Collapse$ forms a loop of $\Pop{1}(s)$. 
\end{proof}

The following lemma provides the analysis of high loops. Every high
loop decomposes into parts that are prefixed by its initial stack
$s$ and parts that are returns of stacks with topmost word
$\Pop{1}(\TOP{2}(s))$. Note the similarity of this characterisation
and its proof with the characterisation of returns in Lemma
\ref{Lemma_Return_Decomposition_For_Simulation}. 

\begin{lemma}
  \label{Lemma_Loop_Decomposition_For_Simulation}
  Let $\lambda$ be some high loop of some stack $s$ with topmost word
  $w=\TOP{2}(s)$. Then 
  there is a sequence $0=:j_0<i_1<j_1< i_2 < j_2 <\dots < i_n < j_
  n \leq i_{n+1}:=\length(\lambda)$ such that
  \begin{enumerate}
  \item for $1\leq k\leq n+1$,   $s\prefixeq \lambda{\restriction}_{[j_{k-1},i_k]}$ and
  \item for each $1\leq k \leq n$, there is a stack $s_k$ with
    $\TOP{2}(s_k) = \Pop{1}(w)$ such that
    $\lambda{\restriction}_{[i_k+1,j_k]}$ is a return of $s_k$. 
  \end{enumerate}
\end{lemma}

\begin{proof}
  This is completely analogous to the proof of Lemma
  \ref{Lemma_Return_Decomposition_For_Simulation}. Assume that
  $\lambda$ is a high loop and let $i_1$ be the minimal position in
  $\domain(\lambda)$ such that $s\notprefixeq \lambda(i_1+1)$. For
  exactly the same reasons as in the return case,
  $\TOP{2}(\lambda(i_1+1))=\Pop{1}(w)$. Since $\lambda$ is a high loop,
  $\lambda(i_1+1)\neq \Pop{1}(s)$ whence 
  $\lvert \lambda(i_1+1) \rvert >  \lvert s \rvert$. Since $\lambda$ ends
  in stack $s$, there is a minimal 
  $j_1>i_1+1$ with $\lvert\lambda(j_1)\rvert <
  \lvert\lambda(i_1+1)\rvert$. Now, completely analogous to the return
  case one concludes that $\lambda{\restriction}_{[i_1+1,j_1]}$ is a
  return: just note that all level $2$ links in
  $\TOP{2}(\lambda(i_1+1))$ are clones of $\TOP{2}(s)$ whence they point
  to stacks of width smaller than $s$. By definition of a loop,
  $\lambda$ cannot use any of these links. Thus, it is clear
  that $\lambda(j_1)=\Pop{2}(\lambda(i_1+1))$. Furthermore, it is easy to
  see that $s\prefixeq \lambda(j_1)$. Thus, an inductive definition of
  the $i_k$ and $j_k$ provides a proof of the lemma. 
\end{proof}
\begin{remark}
  Completely analogous to the decomposition of returns, one proves
  that the sequence $j_0<i_1<\dots< j_n \leq i_{n+1}$ is unique. 
\end{remark}

Having obtained this decomposition of high loops we show that the
computation of $\HighLoopFunc{k}$ can be done analogously to the
computation of $\ReturnFunc{k}$: we use certain runs of the return
simulator $\ReturnSimulator{\mathcal{S}}{s}{k}$ as simulations of high
loops with initial and final stack $s$. 

\begin{definition}
  Let $\mathcal{S}$ be a collapsible pushdown system of level $2$,
  $k\in\N$ some threshold, and 
  $s$  some stack of the form $s=w{\downarrow_0}\Box:w{\downarrow_0}$. 
  We set $s':=\bot\top\TOP{1}(s)\Box:\bot\top\TOP{1}(s)$. 
  We call any run of $\ReturnSimulator{\mathcal{S}}{s}{k}$ from
  $(q_1,s')$ to $(q_2,s')$ a \emph{simulation} of a high loop
  from $(q_1,s)$ to $(q_2, s)$. 
\end{definition}

This terminology is justified for the same reasons as in the case of
returns. Analogously, to the function $\SimToReturnSym{s}$, we next
define a function $\SimToLoopSym{s}$ that translates simulations of
loops into loops with initial and final stack $s$. 

\begin{definition}
  Let $\mathcal{S}, k$, $s$ and $s'$ be as in the previous definition. Let
  $\lambda'$ be a simulation of a high loop from $(q_1,s')$ to
  $(q_2,s')$ where $q_1,q_2\in Q$. 
  Due to the definition of the return simulator
  $\ReturnSimulator{\mathcal{S}}{s}{k}$, the run $\lambda'$ cannot
  pass any substack of $\Pop{1}(s')$ (there are no transitions that
  allow to return to $s'$ once the run reaches $\Pop{1}(s')$ or
  $\Pop{2}(s')$). Thus, $\lambda'$ is a high loop and there is a
  sequence $j_0<i_1<j_1<\dots<j_n\leq i_{n+1}$ according to Lemma
  \ref{Lemma_Loop_Decomposition_For_Simulation}. 
  
  We set $\pi'_k:=\lambda'{\restriction}_{[j_{k-1},i_k]}$ for all 
  $1\leq k \leq n+1$ and
  $\rho'_k:=\lambda'{\restriction}_{[i_k+1,j_k]}$ for all $1\leq k \leq n$.
  
  Completely analogous to what we did in Definition
  \ref{DefSimToReturn} we can define runs 
  $\pi_k, \xi_k$, and $\rho_k$ and set
  \begin{align*}
    \SimToLoopSym{s}(\rho'):= \pi_1 \circ \xi_1 \circ \rho_1
    \circ \pi_2 \circ \xi_2 \circ \rho_2 \circ \dots \circ \pi_n
    \circ \xi_n \circ \rho_n \circ \pi_{n+1}.    
  \end{align*}
  We say $\SimToLoopSym{s}(\lambda')$ is the high loop
  \emph{simulated by  $\lambda'$}.  
\end{definition}

We omit the details of the following claims because they are
completely analogous to the return case. 
One can show that $\SimToLoopSym{s}$ is an injective function
(cf. Lemma \ref{LemmaSimToReturn}). 
The image of $\SimToLoopSym{s}$ contains exactly all those high loops
with initial and final stack $s$ that use length-lexicographic small
returns of stacks with topmost word $\TOP{2}(\Pop{1}(s))$ in the
following sense. Let $\lambda$ be in the image of $\SimToLoopSym{s}$. 
If $\lambda$ contains a subrun that is a return from $(q_1,\hat s)$ to
$(q_2, \Pop{2}(\hat s))$ with $\TOP{2}(\hat s) = \Pop{1}(\TOP{2}(s))$
then this subrun is equivalent to one of the $k$
length-lexicographically smallest returns from $(q_1, \Pop{1}(s))$ to
$(q_2, \Pop{2}(s))$. The proof of this claim is analogous to the proof of  
Lemma \ref{Lem:SimToRetUsesSmall}. 
Finally, if there is a high loop from $(q,s)$ to $(q',s)$ that is not
the image of $\SimToLoopSym{s}$, then there are $k$ high loops from
$(q,s)$ to $(q',s)$ in the image of $\SimToLoopSym{s}$
(cf. Lemma \ref{Lem:LargeSubreturnImpliesManyReturns}). 

Analogous to Lemma \ref{Cor:InductiveComputability}, these results
imply that the number of simulations of high loops is up to threshold
$k$ the number of high loops. 

\begin{corollary} \label{Cor:InductiveComputabilityHighLoops}
  Let $\mathcal{S}$ be a collapsible pushdown system, $k\in\N$ some
  threshold and $w$ some word. 
  For $s:=w{\downarrow_0}\Box:w{\downarrow_0}$, 
  $s':=\bot\top \TOP{1}(s) \Box: \bot\top \TOP{1}(s)$ and
  for $q,q'\in Q$,  let $M_{q,q'}$ be the set of runs of
  $\ReturnSimulator{\mathcal{S}}{s}{k}$ from $(q,s')$ to 
  $(q',s')$. 
  For all $q,q'\in Q$,
  \begin{align*}
    \HighLoopFunc{k}_{\mathcal{S}}(w)(q,q') = \min\{k, \lvert M_{q,q'}
    \rvert\}=
    \HighLoopFunc{k}_{\ReturnSimulator{\mathcal{S}}{s}{k}}(\TOP{2}(s'))(q,q').
  \end{align*}
\end{corollary}

We will soon see that this corollary can be used to define a finite
automaton that computes $\HighLoopFunc{k}$. But before we come to this
result, we briefly examine how we can compute the number of low
loops of a given stack. 

In Lemma \ref{Lem:LowLoopDecomposition} we proved that the
number of low loops of a stack $s$ depends on the number of loops of
$\Pop{1}(s)$. The following lemma shows how we can use this dependence
in order to compute $\LowLoopFunc{k}(s)$ from
$\LoopFunc{k}(\Pop{1}(s))$. 

\begin{lemma}\label{inductiveComputationLowLoops}
  Let $\mathcal{S}$ be some collapsible pushdown system of level $2$,
  $k\in\N$ some threshold. 
  There is a function that computes $\LowLoopFunc{k}(w)$ on input
  $\Sym(w)$, $\Lvl(w)$, $\Sym(\Pop{1}(w))$, $\LoopFunc{k}(\Pop{1}(w))$
  and the transition relation $\Delta$ of $\mathcal{S}$. 
\end{lemma}
\begin{proof}
  A low loop from $(q,s)$ to $(q',s)$ with $w=\TOP{2}(s)$ can only
  exists if $\Lvl(w)=1$. Hence, we only have to consider the case
  $\Lvl(w)=1$.
  Due to Lemma \ref{Lem:LowLoopDecomposition}, a low loop from $(q,s)$
  to 
  $(q',s)$ starts with a 
  $\Pop{1}$ or collapse  (of level $1$) and ends with
  $\Push{\Sym(w)}$. Between these two transitions, the low loop performs a
  loop of $\Pop{1}(s)$. 
  We set 
  \begin{align*}
    &M_{q_1,q_2}:=\{(q_1, \Sym(w), \gamma, q_2, \op)\in\Delta:
    \op=\Pop{1}\text{ or }\op=\Collapse\}\text{ and}\\
    &N_{q_1,q_2}:=\{(q_1, \Sym(\Pop{1}(w)), \gamma, q_2,
    \Push{\Sym(w)})\in\Delta \}.   
  \end{align*}
  Then,  $\LowLoopFunc{k}(s)(q,q')= \min\left\{ k,\sum_{\hat q,\hat q'\in Q} 
  \lvert M_{q,\hat q} \rvert \cdot \LoopFunc{k}(s)(\hat q, \hat q')
  \cdot \lvert N_{\hat q', q'} \rvert\right\}$.  
\end{proof}

By now, we are prepared to prove Proposition
\ref{Prop:AutomatonForLoops}. Recall 
that we have to provide automata that calculate $\LoopFunc{k}$,
$\HighLoopFunc{k}$ and $\LowLoopFunc{k}$. In fact, we provide one
automaton that calculates these functions and $\ReturnFunc{k}$ at the
same time. 

\begin{proof}[Proof of Proposition \ref{Prop:AutomatonForLoops}.]
  Let $\mathcal{S}=(Q,\Sigma,\Gamma,q_0,\Delta)$ be a collapsible
  pushdown system of level $2$.  

  We want to define a finite automaton
  $\mathcal{A}_{\mathrm{lp}}:=(Q_{\mathrm{lp}},
  \Sigma\cup(\Sigma\times\{2\}\times\{0\}), a_0,
  \Delta_{\mathrm{lp}})$ that computes 
  $\ReturnFunc{k}(s),\HighLoopFunc{k}(s),\LowLoopFunc{k}(s)$, and
  $\LoopFunc{k}(s)$ on input $w:=\TOP{2}(s){\downarrow_0}$. 
  
  Recall the following facts.
  \begin{enumerate}
  \item $\ReturnFunc{k}(s) = \ReturnFunc{k}(w)$, $\HighLoopFunc{k}(s)
    = \HighLoopFunc{k}(w)$, etc.
  \item Due to Proposition \ref{Prop:AutomatonForReturns},
    $\ReturnFunc{k}(w)$ is 
    computable by some automaton.
  \item Due to the proof of corollary
    \ref{Cor:InductiveComputabilityHighLoops}, we
    can compute a function $f$ such that for all words $w$ and all
    $\tau\in \Sigma\cup(\Sigma\times\{2\}\times\N)$, we have
    \begin{align*}
    \HighLoopFunc{k}(w\tau)= f(\ReturnFunc{k}(w), \Sym(\tau),
    \Lvl(\tau)). 
    \end{align*}
  \item  Due to Lemma \ref{inductiveComputationLowLoops}, we can
    compute 
    a function $g$ such that
    \begin{align*}
      \LowLoopFunc{k}(w\tau)= g(  \Sym(\tau), \Lvl(\tau), \Sym(w),
      \LoopFunc{k}(w)). 
    \end{align*}
  \item Due to Remark \ref{Rem:LoopsDeterminedByLowAndHighOnes}
    there is a function $h$ such that
    \begin{align*}
    \LoopFunc{k}(w\tau)=h(\HighLoopFunc{k}(w\tau),
    \LowLoopFunc{k}(w\tau)).      
    \end{align*}
  \end{enumerate}
  Analogously to the proof of Proposition \ref{Prop:AutomatonForReturns}, 
  we can use these observations in order to define an automaton that
  computes 
  $\ReturnFunc{k}(s),\HighLoopFunc{k}(s),\LowLoopFunc{k}(s)$, and
  $\LoopFunc{k}(s)$ on input $w:=\TOP{2}(s){\downarrow_0}$. 
\end{proof}

Finally, we have to prove Proposition
\ref{Prop:FuncBoundLoopLengthLemma}. Analogously to the case of returns, we
first define functions 
$\FuncBoundLoopLength{\mathcal{S}}{k}$, 
$\FuncBoundHighLoopLength{\mathcal{S}}{k}$, and
$\FuncBoundLowLoopLength{\mathcal{S}}{k}$.
Then we show that these functions satisfy the conditions of Proposition
\ref{Prop:FuncBoundLoopLengthLemma}.
We first prepare the definition of
$\FuncBoundHighLoopLength{\mathcal{S}}{k}$. Afterwards, we define the
functions mentioned above.

Let $\mathcal{S}$ be some collapsible pushdown system and  
$\mathcal{A}_{\mathrm{lp}}$ the corresponding finite automaton
that calculates the returns, high loops, low loops, and loops of
$\mathcal{S}$.  
Recall that we write $\Delta_{\mathrm{lp}}$  for
the transition relation of $\mathcal{A}_{\mathrm{lp}}$.

Recall that the transitions are labelled by elements of
$\Sigma\cup(\Sigma\times\{2\}\times\{0\})$. 
The run of the automaton on some word $w$ over this alphabet leads
to a state $a$ such that $a$ encodes $\Sym(w), \ReturnFunc{k}(w),
\HighLoopFunc{k}(w)$, etc. 

For each transition $\delta(a,\tau,b)$ we fix a word $w_\delta$ such
that the  run on $w_\delta$ ends in state $a$. 
  
For each $w_\delta$, we define the stack 
$s'_\delta:=\bot\top\tau\Box:\bot\top\tau$. 
Due to Corollary \ref{Cor:InductiveComputability}, there are (up to
threshold $k$) $\HighLoopFunc{k}_{\mathcal{S}}(s'_\delta)(q_1, q_2)$
many 
simulations of high loops
from $(q_1, s'_\delta)$  to $(q_2, s'_\delta)$. 
Using Lemma \ref{LemmaCountNumberOfRuns} we
can compute the $\HighLoopFunc{k}_{\mathcal{S}}(s'_\delta)(q_1, q_2)$ many
lexicographically smallest such 
simulations. We call these $\lambda_1^{\delta,q_1,q_2},
\lambda_2^{\delta,q_1,q_2}, \dots,
\lambda_{\HighLoopFunc{k}(s'_\delta)(q_1, q_2)}^{\delta,q_1,q_2}$. 
Let 
\begin{align*}
  l_{q_1,q_2}^\delta&:=\max\left\{ \length(\lambda_1^{\delta,q_1,q_2}), 
    \length(\lambda_2^{\delta,q_1,q_2}), \dots, 
    \length(\lambda_{\HighLoopFunc{k}(s'_\delta)(q_1, q_2)}^{\delta,q_1,q_2})\right\} \text{and}\\
  \#\top_{q_1,q_2}^\delta&:=
  \max\left\{\left\lvert\left\{j\in\domain(\lambda_i^{\delta,q_1,q_2}):
        \Sym(\lambda_i^{\delta,q_1,q_2}(j))=\top \right\}\right\rvert:  \leq
    \HighLoopFunc{k}(s'_\delta)(q_1, q_2) \right\}      
\end{align*}
be the maximal length of any of these simulations and the 
number of occurrences of $\top$ as topmost symbol in any of these
simulations, respectively. 
Now, set 
\begin{align*}
  &l:=\max\{ l_{q_1,q_2}^\delta: q_1,q_2\in Q,
  \delta\in\Delta_{\mathrm{lp}}\}\text{ and}\\
  &\#\top:=\max\{ \#\top_{q_1,q_2}^\delta:q_1,q_2\in Q,
  \delta\in\Delta_{\mathrm{lp}}\}. 
\end{align*}
\begin{definition}
  We define
  \begin{align*}
    &\FuncBoundHighLoopLength{\mathcal{S}}{k}:\N \to \N\text{ via}\\
    &\FuncBoundHighLoopLength{\mathcal{S}}{k}(0)= 0 \text{ and}\\
    &\FuncBoundHighLoopLength{\mathcal{S}}{k}(n+1)=
    l + \#\top \cdot \FuncBoundReturnLength{\mathcal{S}}{k}(n).     
  \end{align*}
  Furthermore, we define 
  $\FuncBoundLowLoopLength{\mathcal{S}}{k}:\N \to \N$ and
  $\FuncBoundLoopLength{\mathcal{S}}{k}:\N \to \N$ simultaneously
   via 
   \begin{align*}
     &\FuncBoundLowLoopLength{\mathcal{S}}{k}(0):=0,\\ 
     &\FuncBoundLoopLength{\mathcal{S}}{k}(0):=0,\\
     &\FuncBoundLowLoopLength{\mathcal{S}}{k}(n+1):=
     2+\FuncBoundLoopLength{\mathcal{S}}{k}(n)\text{ and}\\
     &\FuncBoundLoopLength{\mathcal{S}}{k}(n+1):=
     \FuncBoundLowLoopLength{\mathcal{S}}{k}(n+1)  
     + 2\cdot 
     \FuncBoundHighLoopLength{\mathcal{S}}{k}(n+1).
   \end{align*}
\end{definition}
\begin{remark}
  The following idea underlies the definition of 
  $\FuncBoundHighLoopLength{\mathcal{S}}{k}$. Let $w$ be 
  some word $w$ of length $n-1$. 
  Assume that 
  we have already proved that
  the lengths of the shortest
  $\ReturnFunc{k}(w)(q_1,q_2)$ many returns from $(q_1, w\Box:w)$ to
  $(q_2,[w\Box])$  are  bound by
  $\FuncBoundReturnLength{\mathcal{S}}{k}(n-1)$. 

  Now, let $s$ be a stack such that $w=\TOP{2}(\Pop{1}(s))$. 
  Due to the definition of $l$, the lexicographically smallest
  $\HighLoopFunc{k}(s)(q,q')$ many  high loops from $(q,s)$ to
  $(q',s)$ have simulations of length at most $l$. 
  
  $\SimToLoopSym{s}$ translates these simulations into
  $\HighLoopFunc{k}(s)(q,q')$ many high loops by copying all transitions
  one by one but by replacing transitions on topmost symbol $\top$ by
  lexicographically small returns equivalent to those from 
  $(q_1, w\Box:w)$ to $(q_2,w\Box)$. Since this replacement happens at
  at most $\#\top$ many positions, we obtain $\HighLoopFunc{k}(s)(q,q')$
  many returns from $(q,s)$ to $(q',\Pop{2}(s))$ of length at most
  $l + \#\top \cdot \FuncBoundReturnLength{\mathcal{S}}{k}(n-1) = 
  \FuncBoundHighLoopLength{\mathcal{S}}{k}(n)$. 

  The other two functions are motivated as follows. A low loop of
  a word $w\sigma$ consists of its initial and final transition plus a
  loop of $w$. Hence, a short low loop consists of a short loop of $w$
  plus $2$ transitions. 

  Due to Lemma \ref{Lem:LoopDecomposition}, a loop of $w\sigma$ is
  either 
  a high loop or consists of a high loop 
  followed by a low loop followed by a high loop. Thus, short loops
  consists of at most three short loops, one a low the two others
  high ones.   
\end{remark}

In analogy to Remark \ref{Rem:BoundingFunctionsMotivation}, 
the previous remark already contains the first half of the proof of 
Proposition \ref{Prop:FuncBoundLoopLengthLemma}. The second half is proved
completely analogous to the return case.


\section{Automatic Structures}
\label{Chapter_AutomaticStructures}

For over 50 years finite automata have been playing a crucial role in
theoretical computer science and have found various applications in
very different fields. In this chapter we recall the basic notions and
techniques concerning finite tree-automata and tree-automatic
structures. In general finite automata come in different flavours. On
one hand automata can be used as acceptors for strings or for
trees and on the other hand one can consider the variants for inputs
of finite or infinite length. 
We mainly focus on finite tree-automata for finite binary trees because we will
use these automata as one of the crucial tools in 
Section \ref{CPG-Tree-Automatic}. 
Nevertheless, in Section \ref{SectionRamseyQuantifier} we will 
also use finite $\omega$-tree-automata on infinite trees as tools for
our proof. 
But only basic facts concerning $\omega$-tree-automata are
actually needed to understand that section.  

For automata on strings, most of the facts we present here are
folklore. Their analogues for tree-automata are mostly 
straightforward generalisations.

This section is organised as follows: we first recall the notions of a
finite tree-automaton and a finite $\omega$-tree-automaton, then we
introduce tree-automatic structures as a form of internal
representation for infinite structures. Finally,  
we recall the known decidability results for model checking on
tree-automatic structures.  

\subsection{Finite Automata}
In this section, we present
the basic theory of tree-automata and tree-automatic
structures. For a more detailed introduction, we refer the reader to 
\cite{tata2007}.
We start by fixing our notation concerning tree-automata.

\begin{definition}
  A \emph{finite tree-automaton} is a tuple
  $\mathcal{A}=(Q,\Sigma,q_I,F,\Delta)$ 
  where $Q$ is a finite nonempty set of states, $\Sigma$ is a finite
  alphabet, $q_I\in Q$ is the initial  state, $F\subseteq Q$ is the
  set of final states,
  and $\Delta\subseteq Q\times Q\times \Sigma \times Q$ is the
  transition relation.
\end{definition}
\begin{remark}
  In the following, we simply write \emph{automaton} for ``finite
  tree-automaton''. 
\end{remark}

We next define the concept of a run of an  automaton on a tree.
Before we state the
definition, recall that for any tree $t$, $t^+$ denotes the minimal
elements of $\{0,1\}^*\setminus\domain(t)$ and $t^\oplus = \domain(t)
\cup t^+$ (cf. Section \ref{HochPlusEinfuehrung}). 
\begin{definition}
  A \emph{run} of $\mathcal{A}$ on a binary $\Sigma$-labelled tree $t$
  is a map  
  $\rho:\domain(t)^\oplus \rightarrow Q$ such that  
  \begin{itemize}
  \item 
    $\rho(d)=q_I $ for all $d\in\domain(t^+)$, and
  \item $\big(\rho(d0), \rho(d1), t(d), \rho(d)\big)\in \Delta$ for
    all $d\in\domain(t)$. 
  \end{itemize}
  $\rho$ is called \emph{accepting} if
  $\rho(\varepsilon)\in F$. 
  We say $t$ is accepted by $\mathcal{A}$ if there is an accepting run
  of $\mathcal{A}$ on $t$. 
  With each automaton $\mathcal{A}$, we associate the \emph{language} 
  \begin{align*}
    L(\mathcal{A}):=\{t: t\text{ is accepted by }\mathcal{A}\}    
  \end{align*}
  \emph{accepted} (or
  recognised)  by $\mathcal{A}$. 
  The class of languages accepted by automata is called the
  class of \emph{regular} languages. 
\end{definition}
\begin{remark}
  Recall that we can consider any string as a tree where each node has
  at most one successor. 
  Using this idea a
  finite string-automaton is just the corresponding special case of a
  finite automaton. 
\end{remark}

One of the reasons for the success of the concept of automata
in computer science is the robustness of this model with
respect to determinisation. We call an automaton $\mathcal{A}$ 
\emph{bottom-up deterministic}, if $\Delta$ is the graph of a function
$Q\times Q\times \Sigma \rightarrow Q$. 

\begin{lemma}[see \cite{tata2007}]
  For each automaton $\mathcal{A}$ there is a bottom-up
  deterministic automaton $\mathcal{A'}$ that accepts exactly the same
  trees as $\mathcal{A}$. 
\end{lemma}

This correspondence between deterministic and nondeterministic
automata is one of the reasons why the class of regular languages has
very strong closure properties. 
\begin{lemma}[see \cite{tata2007}]
  \label{AbschlusseigenschaftenTreeAutomata}
  The regular languages are closed under conjunction, disjunction,
  complementation, and projection. 
\end{lemma}
While the proof of closure under complementation is straightforward
using deterministic automata (just use $Q\setminus F$ as set of
accepting states), the closure under projection is easily shown using
nondeterministic automata. 

In the next section we will see how these closure properties can be
used to turn automata into a useful tool for first-order model
checking purposes using the concept of an automatic structure. 
Beforehand, we recall some more facts about
automata. 
First, we present the generalisation of the pumping lemma for
regular string languages to the tree case. 
Instead of the length of
a string, one uses the  depth of a tree. One obtains
the completely analogous result which states that if a
regular language of trees contains a tree of large depth, then the
language contains infinitely many trees and some of these have smaller
depth than the tree considered initially. 

\begin{lemma}[see \cite{tata2007}]
  Let $\mathcal{A} = (Q, \Sigma, q_I, F, \Delta)$ be an automaton
  recognising the language $L$. For each tree $t\in L$ with
  $\depth{t}> \lvert Q \rvert$, there are nodes
  $d,d'\in \domain(t)$ with $d\leq d'$ such that the following holds. 
  If we replace in $t$ the subtree rooted at $d$ by the subtree rooted
  at $d'$, the tree $t_0$ resulting from this replacement satisfies
  $t_0\in L$. Furthermore, let $t_1, t_2, t_3, \dots $ be the
  infinite sequence of trees where $t_1=t$ and $t_{i+1}$ arises from
  $t_i$ by 
  replacing the subtree rooted at $d'$ in $t_i$ by the subtree rooted
  at $d$ in $t_i$, then $t_i\in L$ for all $i\in\N$. 
\end{lemma}
\begin{proof}
  Take an accepting run $\rho_1$ of $\mathcal{A}$ on $t$. Since 
  $\depth{t} > \lvert Q \rvert$, there are $d,d'\in\domain(t)$
  such that $\rho_1(d)=\rho_1(d')$ and $d\leq d'$. 
  Now, we have to show that the trees $t_0, t_1, t_2, \dots$ are
  accepted by $\mathcal{A}$, i.e., we have to define accepting runs
  for these trees. 
  For $t_0$ consider the run 
  \begin{align*}
    \rho_0:\domain(t_0)^\oplus\rightarrow Q\text{ where }
    \rho_0(e):=
    \begin{cases}
      \rho_1(e) &\text{if }d\not\leq e, \\
      \rho_1(d'f) &\text{if }e=df.
    \end{cases}    
  \end{align*}
  It is easy to see that $\rho_0$ is a run of $\mathcal{A}$ on $t_0$.
  It is
  accepting because $\rho$ was accepting and we did not alter the
  label of the root. 
  For $i\geq 1$ we use the same trick the other way round, setting
  \begin{align*}
    \rho_{i+1}(e):=
    \begin{cases}
      \rho_i(e)  &\text{if }d'\not\leq e,\\
      \rho_i(df) &\text{if }e=d'f.
    \end{cases}    
  \end{align*}
   Again one easily sees that this defines an accepting run of
  $\mathcal{A}$ on $t_{i+1}$. 
\end{proof}

As a direct corollary of the pumping lemma we obtain that finiteness
of regular languages is decidable because finiteness of such a language is
equivalent to not containing a tree of depth between $\lvert Q \rvert$
and $2\lvert Q \rvert$. The latter can be checked by exhaustive search.
\begin{corollary}\label{CorAutomata_FinitenessDecidable}
  Given an automaton $\mathcal{A}$, it is decidable whether
  $L(\mathcal{A})$ is 
  finite. If this is the case, we can compute $\lvert L(\mathcal{A}) \rvert$.
\end{corollary}

We conclude this introduction to automata on trees by recalling
a well known characterisation of
regular classes of trees in terms of \MSO-definability. 

\begin{lemma}[\cite{ThatcherW68}, \cite{Doner70}]
  \label{RegularEqualsMSOonTrees}
  For a set $T$ of finite $\Sigma$-labelled trees,
  there is an  
  automaton recognising $T$ if
  and only if $T$ is \MSO definable.
\end{lemma}

Beside the successful applications of automata on finite
strings or trees
in many areas of computer science, the lifting of the underlying ideas
to the case of infinite inputs had a mayor impact on the importance of
automata theory for computer science. Rabin\cite{Rab69} played a
prominent role in the development of this theory. 

In order to give a meaningful definition of an automaton processing an
infinite tree, we have to reverse the direction in which the automaton
works. 
Up to now, we have considered bottom-up automata, i.e., automata which
start to label a tree at the leaves and then process the tree up to
the root. Of course, one can also imagine an automaton that starts
labelling the root and then labels top-down all the nodes from the
root to the leaves. For the determinisation result we presented, it is
very important to think of a bottom-up automaton. Deterministic
top-down automata are strictly weaker than nondeterministic ones:
there is a regular language which is not 
the language recognised by any top-down deterministic automaton.
Top-down automata become important as soon as we look at infinite
trees. Since there are infinite trees without leaves, the bottom-up
approach is not meaningful anymore. 
But the top-down approach generalises from  finite to  infinite
trees. Considered as a device working  top-down, an automaton is a
device labelling the root with the 
initial state and then
forks into two copies of this automaton -- one for each successor of
the root. Each of these copies now repeats the same procedure on the
corresponding subtree but
starting from a different state according to the transition
relation. With this view, it is straightforward to
generalise the notion of a run from finite trees to infinite
trees. We only have to come up with a new concept of an accepting
run. We now introduce \emph{finite $\omega$-tree automata}. 
Recall that we write $t^\bot$ for the lifting of a (possibly
infinite) tree to the domain $\{0,1\}^*$ by padding with $\bot$. 

\begin{definition}
  A finite \emph{$\omega$-tree automaton} is a tuple 
  $\mathcal{A} =(Q, \Sigma, Q_I, \Delta, \Omega)$, where $Q$ and $\Sigma$,
   and $\Delta$ are as in the case of 
  a finite tree-automaton, $Q_I$ is a set (called the set of initial
  states), and $\Omega$ is 
  a  function $\Omega:Q\rightarrow \N$ (called priority function). 

  A function $\rho:\{0,1\}^*\rightarrow Q$ is a \emph{run} of $\mathcal{A}$
  on an infinite tree $t^\bot$ if $\rho(\varepsilon)\in Q_I$ and $\rho$
  respects $\Delta$. 
  We say $\rho$ is a run on an arbitrary finite or infinite tree $t$
  if it is a run on $t^\bot$. 

  Given some run $\rho$ of $\mathcal{A}$ on $t$, we
  call $\rho$ accepting if
  $\liminf_{n\to\infty}\Omega(\rho(b_1 b_2 b_3 \dots b_n))$ is even for all
  infinite branches $b_1b_2b_3\dots \in\{0,1\}^*$. 
\end{definition}
\begin{remark}
  The acceptance condition that we present here is called
  parity condition. 
  In the literature,  several other acceptance
  conditions for automata on infinite trees are studied, e.g., Buchi-,
  Muller-, 
  Rabin- or Street-conditions. The parity condition turned out to be the
  strongest of all these in the sense that all other 
  conditions mentioned can be reformulated in terms of parity conditions, while
  the parity condition is weak enough in order to transfer most of the
  important results from the theory of finite trees to the infinite
  tree case.  

  In the following we use the term $\omega$-automaton for ``finite
  $\omega$-tree-automaton''. 
\end{remark}

Even though the determinisation result for finite automata does not
carry over to $\omega$-automata, the languages accepted by
$\omega$-automata have the same good closure properties as in the
finite case. Rabin was the first who gave a construction for the
complementation of a nondeterministic $\omega$-automaton. 
In analogy to the finite case, we call the class of languages of
(finite and infinite) trees accepted  by $\omega$-automata 
\emph{$\omega$-regular languages}.

\begin{lemma}[\cite{Rab69}]
  \label{omegaAutomatischAbschlusseigenschaften} 
  The $\omega$-regular languages are closed under conjunction, disjunction,
  complementation, and projection. 
\end{lemma}
The proof of this lemma is through effective constructions of the
corresponding $\omega$-automata. 
Furthermore, the tight correspondence between automata and $\MSO$
carries over 
from the finite to the $\omega$-case.

\begin{theorem}[\cite{Rab69}]
  \label{Rab69AutomaticEqualsMSO}
  A subset $S\subseteq\allTrees{\Sigma}$ is $\omega$-regular if and
  only if it is $\MSO$-definable. 
\end{theorem}

We conclude this brief introduction of $\omega$-automata by recalling
the connection between regular and $\omega$-regular sets of trees. We
show that each regular set has an $\omega$-regular representation via
padding with some label $\bot$. Recall that for a $\Sigma$-labelled
tree $t$, we write 
$t^\bot$ for the full binary tree which coincides with $t$ on
$\domain(t)$ and is labelled by $\bot$ at all other positions. In the
following lemmas, we
assume that $\bot\notin\Sigma$. 

\begin{lemma} \label{Lemma:FinToOmegaAutomaton}
  Given an automaton $\mathcal{A}=(Q,\Sigma,q_I,F,\Delta)$, one
  can construct an 
  $\omega$-automaton $\mathcal{A}^\infty$ such that for all finite
  $\Sigma$-labelled trees $t$,
  \begin{align*}
    \mathcal{A} \text{ accepts } t \text{ iff } \mathcal{A}^\infty
    \text{ accepts } t^\bot.
  \end{align*}
\end{lemma}
\begin{proof}
  The construction of 
  $\mathcal{A}^\infty:=(Q^\infty, \Sigma\cup\{\bot\}, Q^\infty_I, \Delta^\infty,
  \Omega)$ is as follows. We add a 
  new  state  $q_{\mathrm{acc}}$ to the set of
  states by setting  
  $Q^\infty:=Q\cup \{q_{\mathrm{acc}}\}$. 
  Set $Q^\infty_I:=F$ (since we change from the bottom-up view to the
  top-down view, the final states of the automaton become the
  initial state of the $\omega$-automaton). 
  $\Delta^\infty$ is a copy of $\Delta$ enriched by the following 
  transitions:
  $\{ (q_{\mathrm{acc}}, q_{\mathrm{acc}}, \bot,  q_I), 
  (q_{\mathrm{acc}}, q_{\mathrm{acc}}, \bot, q_{\mathrm{acc}}) \}$. 
  The priority function $\Omega: Q^\infty_I \rightarrow \{1,2\}$  is
  defined by
  \begin{align*}
    \Omega(q) =
    \begin{cases}
      1 &\text{if } q\in Q,\\
      2 &\text{if } q=q_{\mathrm{acc}}.
    \end{cases}    
  \end{align*}
  Using these definitions a tree $t$ is accepted if and only if there is
  a finite  initial part $D\subseteq\{0,1\}^*$ such that $\mathcal{A}$
  accepts $t{\restriction}_D$, i.e., it labels $t{\restriction}_{D^+}$
  only with the initial state $q_I$ and all descendants of $D^+$ 
  are nodes labelled by $\bot$. 
  Thus, $\mathcal{A}$ accepts a tree $t'$ if and only if it is of the
  form $t' = t^\bot$ for some finite tree $t$ such that $\mathcal{A}$
  accepts $t$.
\end{proof}

\begin{lemma} \label{Lemma:OmegaToFinAutomaton}
  Given an $\omega$-automaton 
  $\mathcal{A}=(Q, \Sigma, Q_I, \Delta, \Omega)$, one
  can construct an
  automaton $\mathcal{A}^\mathrm{fin}$ such that for all finite
  $\Sigma$-labelled trees $t$,
  \begin{align*}
    \mathcal{A} \text{ accepts } t^\bot \text{ iff } \mathcal{A}^\mathrm{fin}
    \text{ accepts } t.
  \end{align*}
\end{lemma}
\begin{proof}
  We construct
  $\mathcal{A}^{\mathrm{fin}}:=(Q^{\mathrm{fin}},\Sigma,q_I^{\mathrm{fin}},F^{\mathrm{fin}},
  \Delta^{\mathrm{fin}})$ as follows:
  \begin{itemize}
  \item $Q^{\mathrm{fin}}:=Q\cup\{q_{\mathrm{init}}\}$ for a new state
    $q_{\mathrm{init}}$ not contained in $Q$,
  \item $F^{\mathrm{fin}}:=Q_I$,
  \item $q_I^{\mathrm{fin}}:=q_{\mathrm{init}}$, and
  \item $\Delta^{\mathrm{fin}}$ is
    constructed as follows. 
    For each $q\in Q$ we consider the runs of the automaton
    $\mathcal{A}$ with initial state $q$, i.e., the automaton
    $\mathcal{A}_q:=(Q, \Sigma, \{q\}, \Delta, \Omega)$, on the
    $\{\bot\}$-labelled full binary tree $\emptyset^\bot$. We call $q$
    good if there is an accepting run of $\mathcal{A}_q$ on
    $\emptyset^\bot$. Now, for each transition $(q_1, q_2, \sigma,
    q_3)$ we add a new
    transition $(q_{\mathrm{init}}, q_2, \sigma, q_3)$ to
    $\Delta^{\mathrm{fin}}$ if $q_1$ is
    good. Analogously, we add a transition
    $(q_1, q_{\mathrm{init}}, \sigma, q_3)$ to     $\Delta^{\mathrm{fin}}$ if $q_2$ is
    good. Finally, we add
    $(q_{\mathrm{init}},q_{\mathrm{init}},\sigma, q_3)$ to
    $\Delta^{\mathrm{fin}}$ if both $q_1$  and $q_2$ are good. 
    Furthermore, $\Delta^{\mathrm{fin}}$ contains a copy of each
    transition in $\Delta$, i.e., $\Delta\subseteq
    \Delta^{\mathrm{fin}}$. 
  \end{itemize}    
  Now, $\mathcal{A}^\mathrm{fin}$ copies the behaviour of
  $\mathcal{A}$ but at any position where one of the successor nodes
  is labelled by a good state, it can nondeterministically guess
  that the tree it processes is not defined on this successor. 
  If this guess is right, then $\mathcal{A}$ processes at this
  successor a tree which is completely labelled by $\bot$. Since the
  state at this successor is good, the partial run up to this
  position can be extended in such a way that each path starting at
  this successor is accepting. 

  Now, if $\mathcal{A}^\mathrm{fin}$ labels some node by its initial
  state, there is no transition that is applicable at this node. Thus,
  any run of $\mathcal{A}^\mathrm{fin}$ on a tree $t$ labels only
  those positions  by $q_{\mathrm{init}}$ that are in
  $\domain(t)^+$. 

  By definition, a run of $\mathcal{A}^\mathrm{fin}$ on some tree
  $t$ labels all elements of $\domain(t)^+$ by $q_{\mathrm{init}}$
  if and only if there is a run of $\mathcal{A}$ on $t^\bot$ that
  labels all elements of $\domain(t)^+$ by good states. But this is
  equivalent to the fact that $\mathcal{A}$ accepts $t^\bot$ by the
  definition of good states. 
  
  We conclude that  $\mathcal{A}^\mathrm{fin}$ satisfies the claim
  of this lemma. 
\end{proof}

\subsection{Automatic Structures}

As already mentioned the algorithmic tractability of problems on an
infinite structure depends on a good finite representation. 
In this section we recall how automata can be used for this
purpose. The general underlying idea is the following.

Given some structure $\mathfrak{A}$, 
one defines  a tuple of machines from some 
fixed model of computation such that these machines can be used to
evaluate atomic formulas on $\mathfrak{A}$. 
 
A presentation of some structure $\mathfrak{A}=(A, E_1, E_2, \dots,
E_n)$ using a certain model of computation $\mathcal{M}$ is a tuple of
machines $M, M_1, M_2, \dots, M_n$ from $\mathcal{M}$  and some map
$f$ such that the following 
holds.
\begin{itemize}
\item $M$ accepts a set $L$ of strings or trees.
\item  $f$ is a bijective
  map from $L$ to $A$.
\item  $M_i$ accepts a tuple of elements from $L$
  if and only if the image of this tuple under $f$ is in $E_i$. 
\end{itemize}

The first model of computation that was considered for this approach
is that of Turing machines (cf. Appendix \ref{Appendix_Lmuundecidabilty}). 
If one uses Turing machines for representing a structure in this way,
one obtains the so-called class of 
recursive structures (cf. \cite{Harel94}).
But for algorithmic issues, Turing machines turned out to
be far too strong, resulting in the undecidability of model checking on
recursive structures for most logics. 
In general, it is only possible to evaluate 
quantifier-free formula on recursive structures.  

Automata can be used
much more fruitfully as underlying model of computation. This is
due to their good computational behaviour. The resulting
structures are called automatic structures. 
Hodgson \cite{Hodgson82,Hodgson83} first proposed this idea. But it
took more than 10 years until 
the systematic investigation of the general notion of automatic
structures started. 
Khoussainov and Nerode \cite{KhoussainovN94} reintroduced the
notion of string-automatic structures. They obtained the
first 
important results. For instance, they proved that
these structures have decidable $\FO{}$ model checking due to the good
closure properties of regular languages. 
Another boost to the study of automatic structures came from 
the work of Blumensath \cite{Blumensath1999} who
developed the theory further and lifted the idea from the finite
string case to the cases of finite or infinite strings or trees.
Since then, the field of automatic structures has been an active area
of research and many new results have been collected over the years by 
Blumensath, Gr\"adel, Khoussainov, Kuske, Lohrey,  Rubin, et
al.(e.g., \cite{KhoussainovN94, 
  Blumensath00automaticstructures, 
  DBLP:journals/jalc/KhoussainovR01,
  DBLP:conf/lpar/Kuske03,
  BlumensathGra04,
  KhoussainovRS04,
  DBLP:journals/lmcs/KhoussainovNRS07,
  Rubin2008,
  DBLP:journals/apal/KhoussainovM09,
  Kuske2009,DBLP:conf/csl/KuskeL09,
  DBLP:conf/lics/KuskeLL10}).
In the following, we recall the definitions and important results 
with a focus on
tree-automatic structures (which we simply call automatic-structures
in the following). 
String-automatic structures are obtained
by restriction of the accepted languages to languages of strings.
We start by introducing the convolution of trees. This is a tool
for representing an $n$-tuple of $\Sigma$-trees as a single tree over the
alphabet $(\Sigma\cup\{\Box\})^n$ where $\Box$ is a padding symbol
$\Box\notin\Sigma$. 

\begin{definition}
  The \emph{convolution} of two
  $\Sigma$-labelled trees $t$ and $s$ is given by a function
  \begin{align*}
    t\otimes s : \domain(t)\cup\domain(s) \rightarrow
    (\Sigma \cup \{\Box\} )^2
  \end{align*}
  where $\Box$ is some new padding symbol, and
  \begin{align*}
    (t\otimes s)(d) \coloneqq
    \begin{cases}
      (t(d),s(d)) & \text{ if } d\in \domain(t)\cap
      \domain(s), \\ 
      (t(d), \Box) & \text{ if }d\in \domain(t)\setminus
      \domain(s), \\ 
      (\Box, s(d)) & \text{ if }d\in \domain(s) \setminus
      \domain(t).
    \end{cases}
  \end{align*}
  We also use the notation $\bigotimes(t_1, t_2, \dots, t_n)$ for 
  $t_1 \otimes t_2 \otimes \dots \otimes t_n$.
\end{definition}

Using convolutions of trees we can use a single automaton for
defining $n$-ary relations on a set of trees.
Thus, we can then use automata to represent a set and a tuple of
$n$-ary relations on this set. If we can represent the domain of some
structure and all its relations by automata, we call the structure
automatic.  
 
\begin{definition}
  We say a relation $R\subseteq \Trees{\Sigma}^n$ is automatic
  if there 
  is an automaton $\mathcal{A}$ such that 
  $L(\mathcal{A})=\{ \bigotimes(t_1 ,t_2, \dots,
  t_n)\in\Trees{\Sigma}^n: (t_1, t_2, \dots, t_n)\in R\}$. 

  A structure $\mathfrak{B}=(B,E_1, E_2, \dots, E_n)$ with relations
  $E_i$ is \emph{automatic} if there are automata $\mathcal{A}_B,
  \mathcal{A}_{E_1}, \mathcal{A}_{E_2}, 
  \dots, \mathcal{A}_{E_n}$ such that for the language
  $L(\mathcal{A}_B)$ accepted by
  $\mathcal{A}_B$  the following holds:
  \begin{enumerate}
  \item There is a bijection $f: L(\mathcal{A}_B)\rightarrow B$.
  \item For $c_1, c_2, \dots, c_n \in L(\mathcal{A}_B)$, the
    automaton $\mathcal{A}_{E_i}$
    accepts $\bigotimes(c_1, c_2, \dots, c_n)$ if and only if 
    \mbox{$(f(c_1), f(c_2), \dots, f(c_n))\in E_i$.}  
  \end{enumerate}
  In other words, $f$ is a bijection between $L(\mathcal{A}_B)$ and $B$ and
  the relations $E_i$ are automatic via the automata $\mathcal{A}_{E_i}$. 
  We call $f$ a tree presentation of $\mathfrak{B}$.
\end{definition}
Automatic structures form a nice class because automata theoretic
techniques may be used to decide first-order formulas on these
structures:

\begin{theorem}[\cite{Blumensath1999}, \cite{Rubin2008}]
  \label{Thm:FOTreeAutomaticDecidable} 
  If $\mathfrak{B}$ is automatic, then its 
  $\FO{}(\exists^{\mathrm{mod}})$-theory is decidable.
\end{theorem} 
\begin{proof}  
  Given some $\FO{}$ formula $\varphi(x_1, \dots, x_n)$, we can
  construct effectively 
  an automaton $\mathcal{A}_\varphi$ such that $\mathcal{A}_\varphi$
  accepts $t_1 \otimes \cdots \otimes t_n$ if and only if
  $\mathfrak{B},f(t_1), 
  \dots, f(t_n) \models \varphi$ 
  for $f$ the bijection from the previous definition. 
  For atomic formulas, this is clear from the definition of an
  automatic structure because the automata for the relations are
  already given in the definition. 
  Conjunction and negation transform into the classical automata
  constructions of product and complementation. Finally, existential
  quantification corresponds to the closure of regular languages under
  projection. 

  The decidability of the modulo counting quantifier was first proved
  for the string-automatic case in \cite{KhoussainovRS04}.
  Our presentation follows the ideas of Rubin \cite{Rubin2008}.
  He provided a proof for the string-automatic case that 
  allows a straightforward adaption to the case of trees. 

  For simplicity in the presentation, we assume that 
  $\mathfrak{B}$ is an automatic structure whose presentation is
  the identity $\Id$. 
  Let
  \begin{align*}
    \varphi(x,y_1,y_2,\dots, y_n)\in\FO{}(\exists^{\mathrm{mod}})    
  \end{align*}
  be some formula which is represented on $\mathfrak{B}$ by the automaton
  \mbox{$\mathcal{A} = (Q, \Sigma, q_I, F, \delta)$}, i.e.,
  \begin{align*}
    \mathfrak{B}, t,
    t_1, t_2, \dots, t_n \models \varphi(x,y_1,y_2,\dots, y_n)      
  \end{align*}
  if and
  only if $\mathcal{A}$ accepts $\bigotimes(t,t_1,t_2,\dots, t_n)$. 
  
  Given a tuple $t_1\otimes\dots\otimes t_n$ representing the
  assignment of the free
  variables in a formula $\exists^{(k,m)} x (\varphi(x,y_1, \dots,
  y_n))$ for which we want to evaluate the formula, 
  we have to construct an automaton that counts modulo $j$ the number
  of  trees $t$ such that $\mathcal{A}$ accepts 
  \mbox{$t\otimes \bar t:= t\otimes t_1\otimes\dots\otimes t_n$.}  
  Without loss of generality, we assume that $\mathcal{A}$ is a bottom-up
  deterministic automaton.  
  In this case the number of trees $t$ such that $t\otimes \bar t$ is
  accepted by $\mathcal{A}$ coincides with the number of
  accepting runs on trees of the form $t\otimes \bar t$. 
  We now construct an automaton $ \mathcal{\hat A}$ that does this counting. 
  The states $\hat Q$ of $\mathcal{\hat A}$ are functions
  $Q\rightarrow \{0, 1, \dots, m-1, \infty\}$. $\mathcal{\hat A}$
  will label a node $d$ of 
  $\bar t$ with a function $f$ such that there are (modulo $m$) $f(q)$
  different trees $t$ such that the unique run of $\mathcal{A}$ on 
  $t\otimes \inducedTreeof{d}{\bar t}$ labels the root with state
  $q$. By this we mean that $f(q)=\infty$ iff there are infinitely
  many such trees $t$ and otherwise $f(q)$ determines the number of
  such trees modulo $m$. 
  If we know how to label the successors of some node according to
  this rule, then some automaton can update this information. The
  details of  the construction are as follows. 

  We set $\hat Q:=\{0, 1, \dots, m-1, \infty\}^Q$, 
  the initial state is $\hat q_I:=f$ where -- modulo $m$ -- 
  $f(q)$ is 
  \begin{align*}
    \lvert \{t: \rho(\varepsilon)=q \text{ for } \rho \text{ the run
      of } \mathcal{A} \text{ on }  t\otimes \emptyset^n \}\rvert.
  \end{align*}
  Note that $q_I$ is computable due to Corollary  
  \ref{CorAutomata_FinitenessDecidable}. 
  The set of final states $\hat F$ consists of those function
  $f: Q\rightarrow \{0, 1, \dots, m-1, \infty\}$ such that
  $\sum\limits_{q\in F} f(q) = k \mod m$.
  The transition relation $\hat \Delta$ consists of all tuples
  $(f_0, f_1, \sigma, f)$ where $f_0, f_1, f$ are functions 
  $Q\rightarrow \{0, 1, \dots, m-1, \infty\}$ such that
  \begin{align*}
    f(q) = \left(\sum\limits_{(q_0,q_1,\sigma,q)\in\Delta} f_0(q_0) \cdot
    f_1(q_1) \right)\mod k
  \end{align*}
  holds for all $q\in Q$. Note that for fixed $f_0, f_1$, and $\sigma$
  the function $f$ is uniquely determined. Thus,  the resulting
  automaton is a deterministic bottom-up automaton. 

  By an easy induction, one sees that for all $d\in\bar t$ the run of 
  $\mathcal{\hat A}$ on $\inducedTreeof{d}{\bar t}$ labels the root with some
  function $f:Q\rightarrow\{0, 1, \dots, m-1, \infty\}$ such that
  there are  $f(q)$ many different trees $t$ (modulo $m$) such that the
  run $\rho_t$ of $\mathcal{A}$ on $t\otimes \inducedTreeof{d}{\bar t}$
  satisfies $\rho_t(\varepsilon)=q$. 

  From this fact, we directly obtain the desired result, namely, that
  $\mathcal{\hat A}$ accepts $\bar t$ if and only if there are $k$
  modulo $m$ many 
  trees $t$ such that $\mathcal{A}$ accepts $t\otimes \bar t$. 
\end{proof}

\begin{remark}
  The complexity of the \FO{} model checking algorithm for automatic
  structures is nonelementary. 
  This is due to the following facts.
  \begin{itemize}
  \item  Applying a projection to some
    deterministic automaton yields a nondeterministic one.
  \item  Complementation of an automaton can only be done
    efficiently if the automaton is deterministic.
  \item Determinisation of a nondeterministic automaton yields an
    exponential blow-up. 
  \end{itemize}
  Thus, the size of the automaton obtained by
  the construction in the proof is an  exponential tower in the number
  of alternations of existential quantification and negation in the
  formula (or equivalently the number of alternations of existential
  and universal quantifications), i.e., if there are $n$ alternations
  between existential and universal quantification in $\varphi$, then the
  corresponding automaton $\mathcal{A}_\varphi$ may have 
  $\exptower{n}(c)$
  many states (for $c$ some constant)\footnote{
    We denote by $\exptower{i}$ the following function.
    $\exptower{0}(m)\coloneqq c$ and $\exptower{n+1}(c)\coloneqq 
    2^{\exptower{n}(m)}$, i.e., $\exptower{n}(m)$ is an exponential tower
    of height $n$ with topmost exponent $m$.}. 
  
  On the other hand, 
  the algorithm cannot be improved essentially:
  the extension $(\mathbb{N}, +, |_p)$\footnote{$x |_p y$ if $x$ is a
    power of $p$ dividing $y$.} of 
  Presburger Arithmetic is a  string-automatic 
  structure \cite{Blumensath00automaticstructures} 
  which has nonelementary model checking complexity \cite{Gr90d}.
  This implies that there is a nonelementary lower bound for the
  $\FO{}$ model checking complexity on automatic structures. 
\end{remark}

The theory of automatic structures can be naturally extended to the
theory of structures that are represented by automata for infinite trees. 
The structures obtained in this way are called
$\omega$-automatic structures. We conclude this section by precisely
defining $\omega$-automatic structures. 

 \begin{definition}
   We say a relation $R\subseteq (\allTrees{\Sigma})^n$ is $\omega$-automatic
  if there 
  is an $\omega$-automaton $\mathcal{A}$ such that 
  \begin{align*}
    L(\mathcal{A})=\left\{ \bigotimes(t_1 ,t_2, \dots,
    t_n)\in(\allTrees{\Sigma})^n: (t_1, t_2, \dots, t_n)\in R\right\}. 
  \end{align*}
  A structure $\mathfrak{A}=(A,E_1, E_2, \dots, E_n)$ with relations
  $E_i$ is \emph{$\omega$-automatic} if there are $\omega$-automata
  $\mathcal{A}_A,   \mathcal{A}_{E_1}, \mathcal{A}_{E_2}, 
  \dots, \mathcal{A}_{E_n}$ such that for $L(\mathcal{A}_A)$, the
  language accepted by 
  $\mathcal{A}_A$, the following holds:
  \begin{enumerate}
  \item There is a bijection $f: L(\mathcal{A}_A)\rightarrow A$.
  \item For $t_1, t_2, \dots, t_n \in L(\mathcal{A}_A)$, the
    automaton $\mathcal{A}_{E_i}$
    accepts $\bigotimes(t_1, t_2, \dots, t_n)$ if and only if 
    \mbox{$(f(t_1), f(t_2), \dots, f(t_n))\in E_i$.}  
  \end{enumerate}
  We call $f$ an $\omega$-presentation of $\mathfrak{A}$.
 \end{definition}

The similarity of Lemma \ref{AbschlusseigenschaftenTreeAutomata} and Lemma
\ref{omegaAutomatischAbschlusseigenschaften} yields the straightforward
extension of Theorem \ref{Thm:FOTreeAutomaticDecidable} to the
$\omega$-automatic case. 

\begin{theorem}[\cite{Blumensath1999,Blumensath00automaticstructures,
    BarKR09}]  
  \label{DecidabilityOmegaStructures}
  The $\FO{}(\exists^\infty)$-theory of every $\omega$-automatic
  structure is decidable.
\end{theorem}


\chapter{Main Results}
\label{Chapter_MainResults}

\setcounter{section}{-1}
\refstepcounter{section}

In this chapter, we present four results concerning model checking on
certain graph structures. The first and the last involve
automaticity\footnote{We stress that the term ``automaton'' stands for
  ``finite tree-automaton'' and ``automatic'' stands for
  ``tree-automatic''.} 
while the other two are based on modularity arguments for 
Ehrenfeucht-\Fraisse games.

Our first result concerns \FO{} model checking 
on collapsible pushdown graphs of level $2$. 
The expansion of every level $2$ collapsible pushdown graph by regular
reachability and $L\mu$-definable predicates is automatic.
From the general decidability result for \FO{} on automatic
structures, we obtain the following theorem. 
\begin{theorem} \label{Thm:CPGRegular}
 Let 
 $\mathcal{S} = (Q,\Sigma, \Gamma, \Delta, q_0)$ 
 be a collapsible pushdown system of level $2$. Let  
 \begin{align*}
 \mathfrak{G}:=(\CPG(\mathcal{S}), \Reach_{L_1}, \Reach_{L_2}, \dots,
 \Reach_{L_n}, P_1, \dots, P_m)   
 \end{align*}
 be an expansion of the collapsible
 pushdown graph 
 $\CPG(\mathcal{S})$ where $L_1, L_2, \dots, L_n$ are arbitrary
 regular languages over $\Gamma$ and $P_1, \dots, P_m$ are arbitrary
 $L\mu$-definable  predicates. Then the 
 $\FO{}$-theory of $\mathfrak{G}$ is decidable. 
\end{theorem}
  Using our fourth main theorem, 
  even the
  \mbox{\FO{}($\exists^{\infty},
  \exists^{\mathrm{mod}},(\RamQ{n})_{n\in\N}$)-theory} of $\mathfrak{G}$
  is decidable.   
A preliminary version of this theorem was published in
\cite{Kartzow10} and 
we present the proof of this theorem in Section
\ref{CPG-Tree-Automatic}.  

Next, we turn to modularity arguments for Ehrenfeucht-\Fraisse games 
on  nested pushdown trees. The analysis of restricted strategies in
these games lead to model checking algorithms on the class of nested
pushdown trees (cf. Section \ref{Sec:EFGame}). We obtain the following
two results. 

\begin{theorem}
  $\FO{}(\Reach)$ model checking on nested pushdown trees is
  decidable. 
  Furthermore, there is an $\FO{}$ model checking algorithm
  on nested pushdown trees with the following complexities:
  Its structure complexity is in
  $\mathrm{EXPSPACE}$, while its expression complexity and its
  combined complexity are in $2$-$\mathrm{EXPSPACE}$. 
\end{theorem}

\begin{theorem}
  $\FO{}$ model checking on level $2$ nested pushdown trees is
  decidable. 
\end{theorem}
  The concept of a level $2$ nested pushdown tree is a combination of
  the concepts of 
  higher-order pushdown systems and nested pushdown trees. One takes
  a level $2$ pushdown system (without collapse) and
  enriches the unfolding of its graph by jump-edges connecting
  corresponding clone and pop operations (of level $2$). 
  We formally introduce the hierarchy of higher-order nested
  pushdown trees in Section \ref{Chapter HONPT}. 

The proof of the first theorem, which was published in 
\cite{Kartzow09}, is contained in Section
\ref{Chapter_FO-NPT}. The second theorem is proved in Section 
\ref{Chapter HONPT}.

Finally, motivated by the automaticity of collapsible pushdown
graphs of level $2$. We study the model checking problem on the class
of automatic structure. 
We extend the automata-based approach for $\FO{}(\exists^{\infty},
 \exists^{\mathrm{mod}})$ model
checking on automatic structures to
$\FO{}(\exists^{\infty},
 \exists^{\mathrm{mod}},(\RamQ{n})_{n\in\N})$ model checking. In Section
\ref{SectionRamseyQuantifier} we prove the following theorem which
was developed by Dietrich Kuske and the author.

\begin{theorem} \label{Thm:RamseyQuantifierDecidable}
  The $\FO{}(\exists^\infty, \exists^{\mathrm{mod}},
  (\RamQ{n})_{n\in\N})$-theory of automatic structures 
  is decidable. 
\end{theorem}


\section{Level 2 Collapsible Pushdown Graphs are Tree-Automatic}
\label{CPG-Tree-Automatic}
In this section, we focus on collapsible pushdown graphs of level
$2$. Thus, whenever we talk about collapsible pushdown systems
or graphs, we mean those of level $2$.

The main result of this section is the following theorem. 

\begin{theorem}\label{ThmCPGTreeAutomatic}
  Given a collapsible pushdown graph 
  $\CPG(\mathcal{S})=(C(\mathcal{S}),(\trans{\gamma})_{\gamma\in\Gamma})$,
  regular languages $L_1, L_2, \dots, L_n \subseteq \Gamma^*$, and
  $L\mu$-definable predicates $P_1, P_2, \dots, P_m \subseteq C(\mathcal{S})$,
  its expansion
  $(\CPG(\mathcal{S}), \Reach_{L_1}, \Reach_{L_2},
  \dots,\Reach_{L_n}, P_1, P_2, \dots, P_m)$ is
  automatic\footnote{Recall that ``automatic'' is an abbreviation for
    ``tree-automatic''.}.
\end{theorem}

A direct consequence of this result is the automaticity of the
second level of the Caucal hierarchy. 

\begin{corollary}
  The  second level of the Caucal hierarchy is automatic. 
\end{corollary}
\begin{proof}
  The second level of the Caucal hierarchy is obtained by
  $\varepsilon$-contraction\footnote{
    An $\varepsilon$-contraction of a higher-order pushdown graph $G=(V,
    E_1, E_2, \dots, E_n)$ is a graph $(V, E'_1, E'_2, \dots, E'_m)$ for
    $m\leq n$ where $E'_i:=\Reach_{L_i}$ for $L_i=L( (E_{m+1} + E_{m+2} +
    \dots + E_n)^* E_i)$.}
  from the class of higher-order pushdown graphs of level
  $2$ (cf. \cite{cawo03}). 
\end{proof}

Using Theorem \ref{Thm:RamseyQuantifierDecidable} we obtain the
decidability of the first-order theory of collapsible pushdown graphs
of level $2$:

\begin{corollary}
  Let $\mathcal{S}$ be a collapsible pushdown system of level $2$. 
  Let $\mathfrak{G}$ be the expansion of $\CPG(\mathcal{S})$ 
  by $L\mu$-definable predicates and by regular reachability predicates. 
  Under these conditions, the \FO{}(Reg,
  $\exists^{\mathrm{mod}},(\RamQ{n})_{n\in\N}$)-theory 
  of $\mathfrak{G}$ is decidable.
\end{corollary}
\begin{remark}
  Note that this corollary is just a reformulation of Theorem
  \ref{Thm:CPGRegular}. 
\end{remark}

The main part of this section consists of  
a proof of theorem \ref{ThmCPGTreeAutomatic}. 
Furthermore, we discuss the limitations of our
approach and the limitations of first-order model checking on
collapsible pushdown graphs in general. 

The section is organised as follows.
In Section \ref{STACS:SecEncoding}, we present a function $\Encode$
which translates configurations of collapsible pushdown systems into
trees. This function $\Encode$ yields an automatic representation
for every collapsible pushdown graph. 
We show in Section \ref{STACS:SecCertificates} that the reachable
configurations of a collapsible pushdown system are turned into a
regular set of trees by $\Encode$.
The proof of this statement takes the results on loops from section
\ref{ChapterLoops} as a main ingredient.
Recall that the loops of a given stack can be calculated by a
string-automaton reading the topmost word of the stack. This result carries
over to an automaton reading the encoding of a given stack. 
Since runs from the initial configuration to some configuration $c$
mainly consist of loops, this kind of
regularity of loops can be used to show the regularity of the
set of reachable configurations. 
In Section \ref{Reg_Stack_OP}, we prove that the stack operations
are regular via $\Encode$. Hence, for each
transition relation $\trans{\gamma}$ there is an automaton recognising
those encodings of pairs of configurations that are related by
$\trans{\gamma}$.  
Then we show  that regular reachability predicates
over $\Gamma^*$ are regular sets via $\Encode$.  
This is done in Section \ref{Reachability_Tree-Automatic} 
as follows. First, we prove that the image of the  ``ordinary''
reachability predicate $\Reach$
is a regular relation via $\Encode$. 
Then we show that collapsible
pushdown graphs are closed under products with string-automata.
Finally, we reduce
the predicate $\Reach_L$ to the predicate $\Reach$ on the product of
the collapsible pushdown system and the string-automaton corresponding to $L$.

Afterwards, we  relate our result to other known results. In
Section \ref{Lmu-FO-CPG}, we first 
investigate combinations of the  known $L\mu$ model checking algorithm
with our $\FO{}$ model checking algorithm. Then, in Section 
\ref{sec:LowerBounds}, we provide a lower bound for \FO{} model
checking on level $2$ collapsible pushdown graphs. 
Recall that the first-order model checking on automatic
structures has nonelementary complexity.  
We show that the complexity of the first-order model checking on
collapsible pushdown graphs  is also nonelementary. Thus, our model
checking algorithm cannot be improved essentially.
In the final part we discuss first-order
model checking on higher-order collapsible pushdown graphs. 
Recently, Broadbent \cite{Broadbent2010Mail} showed that first-order
model checking 
is undecidable on level $3$ collapsible pushdown graphs. 
Thus, there is no hope to extend our technique to higher levels of the
collapsible pushdown hierarchy.

\subsection{Encoding of Level 2 Stacks in Trees}
\label{STACS:SecEncoding}
 
In this section we present an encoding of level $2$ stacks in trees. 
The idea is to divide a stack into blocks and to
encode different blocks in different subtrees.
The crucial observation is that every stack is a list of words that
share the same first letter. A block is a maximal list of words occurring
in the stack which share the same two first letters. If we remove the first
letter of every word of such a block, the resulting  $2$-word
decomposes again as a list of blocks. Thus, we can inductively
carry on to decompose parts of a stack into blocks and encode every
block in a different subtree. The roots of these subtrees are
labelled with the first letter of the block. 
This results in a tree where every initial left-closed
path in the tree represents one word of the stack. 
A path of a tree is left-closed if its last element has no left
successor.

As we already mentioned, the encoding works by dividing stacks
into blocks.  
The following notation is useful for the formal definition of blocks. 
Let $w\in\Sigma^*$ be some word and \mbox{$s=w_1:w_2:\dots:w_n
\in\Sigma^{*2}$} some stack. We
write $s'\coloneqq w\mathrel\backslash s$ 
for $s'=ww_1 : ww_2 : \dots: ww_n$. 
Note that $[w]$ is a
prefix of $s'$, i.e., in the notation from Definition
\ref{Def:Prefixrelation},  $[w] \prefixeq w\mathrel\backslash s$. We
say that $s'$ is \emph{$s$ prefixed by $w$}. 

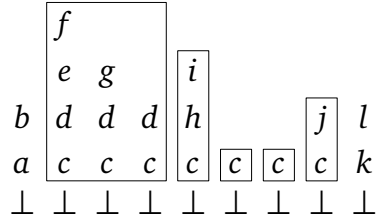
\begin{figure}
  \centering
  $
  \begin{xy}
    \xymatrix@=0.2mm{
      & f \\
      & e & g & & i\\
      b& d & d & d & h & & &j & l\\
      a & c & c & c & c & c & c&c & k\\
      \bot & \bot & \bot & \bot & \bot & \bot & \bot & \bot&\bot
      \save "1,2"."4,4"*[F-]\frm{}
      \save "2,5"."4,5"*[F-]\frm{}
      \save "4,6"."4,6"*[F-]\frm{}
      \save "4,7"."4,7"*[F-]\frm{}
      \save "3,8"."4,8"*[F-]\frm{}
    }
  \end{xy}
  $
  \caption{A stack with blocks forming a $c$-blockline.}
  \label{STACSfig:Blocks}
\end{figure}

\begin{definition}
  Let $\sigma\in\Sigma$ and
  $b\in\Sigma^{*2}$. We call $b$  a \emph{$\sigma$-block} if
  $b=[\sigma]$ or $b=\sigma\tau\mathrel{\backslash} s'$ for some
  $\tau\in\Sigma$ and some $s'\in\Sigma^{*2}$.
  If $b_1, b_2,\dots, b_n$ are $\sigma$-blocks, then we call
  $b_1:b_2:\dots :b_n$ a \emph{$\sigma$-blockline}. 
  See Figure \ref{STACSfig:Blocks} for an example of a blockline with
  its blocks. 
\end{definition}

Note that every 
stack in $\Stacks_{2}(\Sigma)$ forms a $\bot$-blockline. 
Furthermore, every blockline $l$ decomposes uniquely as
\mbox{$l=b_1: b_2: \dots: b_n$} of maximal blocks $b_i$ in $l$. We
will call these 
maximal blocks the blocks of $l$. 

Another crucial observation is that a $\sigma$-block
$b\in\Sigma^{*2}\setminus \Sigma$
decomposes as $b=\sigma\mathrel{\backslash}l$ for some blockline $l$
and we call $l$  the  blockline \emph{induced} by $b$. For a block of
the 
form $[b]$ with $b\in\Sigma$, we define the
blockline induced by $[b]$ to be $\emptyset$.

Recall that the symbols of a collapsible pushdown stack (of level $2$)
come from the set $\Sigma \cup (\Sigma \times\{2\}\times\N)$ where
$\Sigma$ is the stack alphabet. 

We are now going to define our encoding of stacks in trees. 
For $\tau\in\Sigma \cup (\Sigma \times\{2\}\times\N)$, we encode a
$\tau$-blockline $l$ in a tree as follows. The root of the tree is
labelled by $(\Sym(\tau), \Lvl(\tau))$. 
The
blockline induced by the first block of 
$l$ is encoded in the left subtree  and  the
rest of $l$ is encoded in the right subtree. This means that we only
encode explicitly the symbol and the 
collapse level of each element of the stack, but not the
collapse link. We will later see how to decode the collapse links from the
encoding of a stack. 
When we
encode a part of a blockline in the right subtree, 
we do not repeat the label $(\Sym(\tau), \Lvl(\tau))$, but 
replace it by the empty word $\varepsilon$. 

\begin{definition}
  Let $\tau\in \Sigma\cup(\Sigma\times\{2\}\times\N)$. Furthermore,
  let 
  \begin{align*}
    s =
    w_1:w_2:\dots:w_n\in(\Sigma\cup(\Sigma\times\{2\}\times\N))^{+2}    
  \end{align*}
  be some  
  $\tau$-blockline. Let $w_i'$ be words
  such that 
  \mbox{$s= \tau \mathrel\backslash [w_1': w_2' :\dots
    :w_n']$} and set \mbox{$s'\coloneqq w_1': w_2' : \dots : w_n'$}. As an
  abbreviation we write  
  $_is_k\coloneqq w_i:w_{i+1}:\dots:w_k$.
  Furthermore, let  
  $w_1:w_2:\dots: w_j$ be a maximal block of $s$. Note that $j>1$ implies 
  that there is some $\tau'\in
  \Sigma\cup(\Sigma\times\{2\}\times\N)$ and there are words $w''_{j'}$ for each
  $j'\leq j$ such that
  $w_{j'}= \tau\tau' w_{j'}''$. 

  Now, for arbitrary
  $\sigma\in(\Sigma\times\{1,2\})\cup\{\varepsilon\}$, we define
  recursively the
  $(\Sigma\times\{1,2\})\cup\{\varepsilon\}$-labelled
  tree $\Encode(s,\sigma)$ via
  \begin{align*}
    \Encode(s,\sigma)\coloneqq
    \begin{cases}
      \sigma & \text{if } \lvert w_1\rvert=1, n=1\\
      \treeR{\sigma}{\Encode( _2s_n,\varepsilon)}
      &\text{if } \lvert w_1 \rvert =1, n>1\\
      \treeL{\sigma}{\Encode( _1s_n',(\Sym(\tau'), \Lvl(\tau')))}
      &\text{if }  \lvert w_1\rvert >1, j=n \\
      \treeLR{\sigma}{\Encode(
        _1s_j',(\Sym(\tau'),\Lvl(\tau')))}{\Encode( 
        _{j+1}s_n,\varepsilon)}      
      &\text{otherwise}
    \end{cases}
  \end{align*}
For every $s\in\Stacks_2(\Sigma)$, $\Encode(s)\coloneqq
\Encode(s,(\bot,1))$ is called the \emph{encoding of the stack} $s$.
\end{definition}
Figure \ref{STACSfig:Encoding} shows a configuration and its encoding. 
\begin{figure}[t]
  \centering
  $
  \begin{xy}
    \xymatrix@R=0pt@C=0pt{
      &       & (c,2,1) &         & e & \\
      &(b,2,0)& (b,2,0) &       c & (d,2,3) & \\
      &(a,2,0)& (a,2,0) & (a,2,2) & (a,2,2) & (a,2,2) \\
      & \bot  & \bot    &\bot     & \bot    & \bot
      }
  \end{xy}$
  \hskip 1cm
  $\begin{xy}
    \xymatrix@R=9pt@C=3pt{
           & c,2 &       & e,1      &  \\
      b,2 \ar[r]& \varepsilon\ar[u] & c,1 & d,2\ar[u] & \\
      a,2 \ar[u] &  & a,2\ar[r]\ar[u] & \varepsilon\ar[r]\ar[u]
      &\varepsilon \\
       \bot,1 \ar[rr]\ar[u] &  &\varepsilon\ar[u] &  & &      
      }
  \end{xy}
  $
  \caption{A stack $s$ and its Encoding $\Encode(s)$: right
    arrows lead to $1$-successors (right successors), upward arrows
    lead to $0$-successors (left successors).}
  \label{STACSfig:Encoding}
\end{figure}

\begin{remark}
  Fix some stack $s$. 
  For $\sigma\in\Sigma$ and $k\in\N$, every $(\sigma,2,k)$-block of
  $s$ is
  encoded in a subtree whose root $d$ is labelled $(\sigma,2)$. 
  We can restore
  $k$ from the position of $d$ in the tree $\Encode(s)$ as follows.
  \begin{align*}
    k = \lvert \{d'\in \domain\Encode(s) \cap \{0,1\}^*1:
    d'\leq_{\mathrm{lex}} d\}\rvert,     
  \end{align*}
  where $\leq_{\mathrm{lex}}$ is the lexicographic order. This
  is due to the fact that every right-successor corresponds to the
  separation of some block from
  some other. 

  This correspondence can be seen as a bijection.
  \label{milestonesInEncoding}
  Let \mbox{$s=w_1:w_2:\dots:w_n$} be some stack. We define the set 
  \mbox{$R\coloneqq \domain(\Encode(s)) \cap(\{\varepsilon\}\cup
    \{0,1\}^*1)$}. Then there is a bijection
  \mbox{$f: \{ 1,2, 3,   \dots, n\} \rightarrow R$} 
  such that $i$ is mapped to the $i$-th
  element of $R$ in lexicographic order. Each $1\leq i \leq n$
  represents the $i$-th word of $s$. 
  $f$ maps the
  first word of $s$ to the root of $\Encode(s)$ and every other word in
  $s$ to the element of $\Encode(s)$ that separates this word from its left
  neighbour in $s$. 

  If we interpret $\varepsilon$ as empty word, the word from the root
  to $f(i)$ in $\Encode(s)$ is the greatest common prefix of $w_{i-1}$
  and $w_i$. 
  More precisely,
  the word read along this path is the projection onto the letters and
  collapse levels of $w_{i-1}\sqcap w_i$.

  Furthermore, set $f'(i)\coloneqq d0^m\in\Encode(s)$ for $d:=f(i)$
  such that 
  $m$ is maximal with this property, i.e., $f'(i)$ is the leftmost
  descendent of $f(i)$. Then the path from $f(i)$ to $f'(i)$ is the
  suffix $w_i'$ such that $w_i=(w_{i-1}\sqcap w_i) \circ w_i'$ (here
  we set $w_0\coloneqq \varepsilon$). More precisely, the word read along this
  path is the projection onto the symbols and collapse levels of $w_i'$.
\end{remark}

Having defined the encoding of a stack, we want to encode whole
configurations, i.e., a stack together with a state. To this end, we
just add the state as new root of the tree and attach
the encoding of the stack as left subtree, i.e.,
for some configuration $(q,s)$ we set
\begin{align*}
  \Encode(q,s)\coloneqq \treeL{q}{\Encode(s)}.
\end{align*}

The image of this encoding function contains only trees of a very
specific type. We call this class $\EncTrees$. 
In the next definition we state the characterising properties of
$\EncTrees$. This class is  \MSO-definable whence
automata-recognisable (cf. Lemma \ref{RegularEqualsMSOonTrees}). 
\begin{definition} \label{STACS:DefEncodingTrees}
  Let $\EncTrees$ be the class of trees $T$ that satisfy
  the following conditions.
  \begin{enumerate}
  \item  The root of $T$ is labelled by some element of $Q$
    ($T(\varepsilon)\in Q$).
  \item Every element of the form $\{0,1\}^*0$ is labelled by some
    $(\sigma,l)\in\Sigma\times\{1,2\}$, especially $T(0)=(\bot,1)$. 
  \item Every element of the form $\{0,1\}^*1$ is labelled by
    $\varepsilon$. 
  \item $1\notin\domain(T)$,   $0\in\domain(T)$.
  \item \label{STACS:fifthofDefEnc} For all $t\in T$ we have that
    $T(t0)=(\sigma,1)$ implies 
    $T(t10)\neq(\sigma,1)$. 
  \end{enumerate}
\end{definition}
\begin{remark}
  Note that all trees in the image of $\Encode$ satisfy condition
  \ref{STACS:fifthofDefEnc} due to the following. 
  \mbox{$T(t0)=T(t10)=(\sigma,1)$} would imply that the 
  subtree rooted at $t$ encodes a blockline $l$ such that the first
  block $b_1$ of $l$ induces a $\sigma$-blockline and the second
  block $b_2$
  induces also a $\sigma$-blockline.  This
  contradicts the maximality of the blocks used in the encoding because
  all words of $b_1$ and $b_2$ have $\sigma$ as
  second letter whence $b_1:b_2$ forms a larger block.
  Note that for letters with links of level $2$ the analogous
  restriction does not hold. In Figure \ref{STACSfig:Encoding} one
  sees the encoding of a stack $s$ where $\Encode(s)(0) =
  \Encode(s)(10) = (a,2)$. Here, the label $(a,2)$ represents two
  different letters. $\Encode(s)(0)$ encodes the element $(a,2,0)$,
  while 
  $\Encode(s)(10)$ encodes the element $(a,2,2)$,
  i.e., the first element encodes a letter $a$ with undefined link and
  the second  encodes the  
  letter $a$ with a link to the substack of width $2$. 
\end{remark}

Having defined the encoding function $\Encode$, we next show that it
induces a bijection between the configurations of \CPG and $\EncTrees$. 
The rest of this section is a formal proof of the following
lemma. 

\begin{lemma} \label{STACS:Bijective}
  $\Encode:Q\times \Stacks_{2}(\Sigma) \rightarrow \EncTrees$ is a 
  bijection. We denote its inverse by $\Decode$.   
\end{lemma}

The formal proof of this lemma is rather technical. 
The reader who is not interested in the technical details of this
proof may continue with reading Section
\ref{STACS:SecCertificates} directly.

We start our proof of the lemma by explicitly constructing 
the inverse of $\Encode$. This inverse is called \mbox{$\Decode$}. Since
$\Encode$ removes the collapse links of 
the elements in a stack, we have to restore these now. For restoring
the collapse links, we use the following auxiliary 
function. For $g\in N$ and $\tau\in
\{\varepsilon\}\cup(\Sigma\times\{1,2\})$, we set 
\begin{align*}
  f_g(\tau)\coloneqq
  \begin{cases}
    \sigma &\text{if } \tau=(\sigma,1),\\
    (\sigma,2,g)  &\text{if } \tau=(\sigma,2),\\
    \varepsilon &\text{if } \tau=\varepsilon.
  \end{cases}
\end{align*}
Later, $g$ will be the width of the stack decoded so far.
\begin{definition}
  Let $\Gamma:=(\Sigma\times\{1,2\})\cup\{\varepsilon\}$.
  We define the following function 
 \mbox{ $\Decode: \Trees{\Gamma}\times \N \rightarrow
  (\Sigma\cup(\Sigma\times\{2\}\times\N))^{*2}$} by recursion. Let
  \begin{align*}
      \Decode(T,g) = 
      \begin{cases}
        f_g(T(\varepsilon)) &
        \text{if }\domain(T)=\{\varepsilon\},\\
        f_g(T(\varepsilon)) \mathrel{\backslash}
        \Decode(\inducedTreeof{0}{T},g) &\text{if } 1\notin\domain(T),\\
        f_g(T(\varepsilon)) \mathrel{\backslash}
        (\varepsilon:\Decode(\inducedTreeof{1}{T},g+1)) &
        \text{if } 0\notin\domain(T),\\
        f_g(T(\varepsilon)) \mathrel{\backslash}
        (\Decode(\inducedTreeof{0}{T},g):
        \Decode(\inducedTreeof{1}{T},g+G(\inducedTreeof{0}{T}))) &
        \text{otherwise,}
      \end{cases}
    \end{align*} where
    $G(\inducedTreeof{0}{T})\coloneqq
    \lvert\Decode(\inducedTreeof{0}{T},0)\rvert$ is 
    the width of the 
    stack encoded in $\inducedTreeof{0}{T}$. 
    For a tree $T\in\EncTrees$, the decoding of $T$ is
    \begin{align*}
      \Decode(T)\coloneqq(T(\varepsilon),
      \Decode(\inducedTreeof{0}{T},0))\in Q\times
      (\Sigma\cup(\Sigma\times\{2\}\times\N))^{+2}.       
    \end{align*}
\end{definition}
\begin{remark}
  Obviously, for each $T\in\EncTrees$, $\Decode(T)\in
  Q\times(\Sigma\cup (\Sigma\times\{2\}\times\N))^{+2}$. 
  In fact, the image of $\Decode$ only consists of configurations,
  i.e., $\Decode(T)=(q,s)$ such that $s$ is a level $2$ stack. 
  The verification of this claim relies on two important
  observations. 

  Firstly, $T(0)=(\bot,1)$ due to condition 2 of Definition
  \ref{STACS:DefEncodingTrees}. Thus, all words in $s$ start with
  letter $\bot$. 

  Now, 
  $s$ is a stack if and only if the link structure of $s$ can be
  created using the push, clone and $\Pop{1}$ operations. 
  The proof of this claim can be done by a tedious but straightforward
  induction. We only sketch the most important observations for this
  fact. 

  Every letter $a$ of the form
  $(\sigma,2,l)$ occurring in $s$ is either a clone or can be created
  by the $\Push{\sigma,2}$ operation.   
  We call $a$ a clone if $a$
  occurs in $s$ in some word $waw'$ such that the word to the left of
  this word has $wa$ as prefix. Note that cloned elements are those
  that can 
  be created by use of the $\Clone{2}$ and $\Pop{1}$ operations from a
  certain substack of $s$. 
  
  If $a$ is not a clone in this sense, then $\Decode$ creates the
  letter $a$ because there is some $(\sigma,2)$-labelled node in $T$
  corresponding to $a$. Now, the important observation is that
  $\Decode$ defines
  $a=f_g((\sigma,2))$ where $g+1$ is the width of the stack decoded from
  the lexicographically smaller nodes. Hence, the letter $a$ occurs in
  the $(g+1)$-st word of $s$ and points to the $g$-th word. 
  Such a letter $a$ can clearly be created by a $\Push{\sigma,2}$
  operation. Thus, all $2$-words in the image of $\Decode$ can be
  generated by stack operations from the initial stack. A
  reformulation of this observation is that the image of $\Decode$
  only contains stacks.
\end{remark}

Now, we prove that 
$\Decode$ is injective on $\EncTrees$.  Afterwards, we show
that 
$\Decode\circ\Encode$ is 
the identity on the set of all configurations. This implies that
$\Decode$ is a surjective map from $\EncTrees$ to 
$Q\times\Stacks_2(\Sigma)$. Putting both facts together, we obtain the
bijectivity of $\Encode$. 

\begin{lemma} \label{DecInjective}
  $\Decode$ 
  is injective on $\EncTrees$.
\end{lemma}
\begin{proof}
  Assume that there are trees $T', U'\in\EncTrees$ with
  $\Decode(T')=\Decode(U')=(q,s)$. Then by definition 
  $T'(\varepsilon) = U'(\varepsilon) = q$. Thus, we only have to compare
  the subtrees
  rooted at $0$, i.e., $T\coloneqq\inducedTreeof{0}{T'}$ and 
  $U\coloneqq\inducedTreeof{0}{U'}$. From our assumption it follows that
  \mbox{$\Decode(T,0)=\Decode(U,0)$.} 

  Note that the roots of $T$ and of $U$ are both labelled by 
  $(\bot, 1)$. 

  Now, the lemma follows from the 
  following claim. 
  \begin{claim}
    Let $T$ and $U$ be trees such that there are $T', U'\in\EncTrees$
    and $d\in\domain(T')\setminus \{\varepsilon\}$, 
    $e\in\domain(U')\setminus \{\varepsilon\}$ such that
    $T=\inducedTreeof{d}{T'}$ and 
    $U=\inducedTreeof{e}{U'}$. 
    If $\Decode(T,  m) = \Decode(U, m)$ and either
    $T(\varepsilon) = U(\varepsilon) = \varepsilon$ or
    $T(\varepsilon)\in\Sigma\times\{1,2\}$ and
    $U(\varepsilon)\in\Sigma\times\{1,2\}$, then 
    $U=T$. 
  \end{claim}

  The proof is by induction on the depth of the trees $U$ and
  $T$. 
  If \mbox{$\depth{U}=\depth{T}=0$}, $\Decode(U,m)$ and $\Decode(T,m)$
  are uniquely determined by
  the label of their roots. 
  A straightforward consequence of the definition of $\Decode$ is
  that $U(\varepsilon)=T(\varepsilon)$ whence $U=T$.

  Now, assume that the claim is true for all trees of depth at most
  $k$ for some fixed $k\in\N$. Let $U$ and $T$ be trees of
  depth at most $k+1$. 

  We proceed by a case distinction on whether the left or right
  subtree of $T$ and $U$ are defined. 
  In fact, 
  $\Decode(T, m) = \Decode(U, m)$ implies that
  \begin{enumerate}
  \item $\inducedTreeof{0}{T}\neq \emptyset$ if and only if 
    $\inducedTreeof{0}{U}\neq \emptyset$ and
  \item $\inducedTreeof{1}{T}\neq \emptyset$ if and only if
    $\inducedTreeof{1}{U}\neq \emptyset$. 
  \end{enumerate}
  We  first prove that $\Decode(T,  m)=\Decode(U,m)$ implies
  $U=T$ in the cases satisfying these conditions. 
  Afterwards, we show that all possible combinations that do not
  satisfy this condition imply $\Decode(T,  m) \neq
  \Decode(U,m)$. 
  
  \begin{enumerate}
  \item \label{Case1111}
    Assume that $\inducedTreeof{0}{U} = \inducedTreeof{1}{U} =
    \inducedTreeof{0}{T} = \inducedTreeof{1}{T} = \emptyset$. Then
    $\depth{T}= \depth{U} = 0$. For trees of depth $0$ we have already
    shown that  $\Decode(U, 0) = \Decode(T,0)$ implies $U=T$.
  \item \label{Case1010} 
    Assume that
    $\inducedTreeof{0}{U} = \emptyset$, 
    $\inducedTreeof{1}{U} \neq \emptyset$, 
    $\inducedTreeof{0}{T} = \emptyset$ and
    $\inducedTreeof{1}{T} \neq \emptyset$.
    In this case 
    \begin{align*}
      &\Decode(U, m) = f_m(U(\varepsilon)) \mathrel{\backslash}
      (\varepsilon: \Decode(\inducedTreeof{1}{U}, m+1))\text{
        and}\\
      &\Decode(T, m) = f_m(T(\varepsilon)) \mathrel{\backslash}
      (\varepsilon: \Decode(\inducedTreeof{1}{T}, m+1)). 
    \end{align*}
    Since $U(\varepsilon)=\varepsilon$ if and only if 
    $T(\varepsilon)=\varepsilon$, we can directly conclude that
    $U(\varepsilon) = T(\varepsilon)$. 
    But then $\Decode(T, m)  = \Decode(U, m)$ implies that 
    $\Decode(\inducedTreeof{1}{T}, m+1) =
    \Decode(\inducedTreeof{1}{U}, m+1)$.  
    Since $\depth{\inducedTreeof{1}{T}}\leq k$ and 
    $\depth{\inducedTreeof{1}{U}}\leq k$, the induction hypothesis
    implies that $\inducedTreeof{1}{T}=\inducedTreeof{1}{U}$. 
    We
    conclude that $T=U$.
  \item \label{Case0101}
    Assume that
    $\inducedTreeof{0}{U} \neq \emptyset$, 
    $\inducedTreeof{1}{U} = \emptyset$, 
    $\inducedTreeof{0}{T} \neq \emptyset$, and
    $\inducedTreeof{1}{T} =\emptyset$.
    In this case, 
    \begin{align*}
      &\Decode(U, m) = f_m(U(\varepsilon)) \mathrel{\backslash}
      \Decode(\inducedTreeof{0}{U}, m)\text{ and}\\
      &\Decode(T, m) = f_m(T(\varepsilon)) \mathrel{\backslash}
      \Decode(\inducedTreeof{0}{T}, m).      
    \end{align*}
    Since $U(\varepsilon)=\varepsilon$ if and only if
    $T(\varepsilon)=\varepsilon$, we conclude that
    $U(\varepsilon)=T(\varepsilon)$ and  
    $\Decode(\inducedTreeof{0}{U}, m) = 
    \Decode(\inducedTreeof{0}{T}, m)$. 
    Since the depths of $\inducedTreeof{0}{U}$ and of
    $\inducedTreeof{0}{T}$ are at most $k$, the induction hypothesis
    implies $\inducedTreeof{0}{U}=\inducedTreeof{0}{T}$ whence $U=T$.
  \item \label{Case0000} 
    Assume that
    $\inducedTreeof{0}{U}\neq \emptyset$, 
    $\inducedTreeof{1}{U}\neq \emptyset$, 
    $\inducedTreeof{0}{T}\neq \emptyset$,  and
    $\inducedTreeof{1}{T}\neq \emptyset$. 
    Then we have 
    \begin{align*}
      & 
      \Decode(U,m)= f_m(U(\varepsilon)) \mathrel{\backslash} 
      \left(\Decode(\inducedTreeof{0}{U},m):
        \Decode(\inducedTreeof{1}{U},m+m')\right) \text{ and }\\
      & 
      \Decode(T,m)= f^n_m(T(\varepsilon)) \mathrel{\backslash} 
      \left( \Decode(\inducedTreeof{0}{T},m):
        \Decode(\inducedTreeof{1}{T},m+m'')\right)
    \end{align*}
    for some natural numbers  $m',m''>0$. 
    
    Since $U(\varepsilon)=\varepsilon$ if and only if
    $T(\varepsilon)=\varepsilon$ this implies that the roots of $U$
    and $T$ coincide. 
    Hence,
    \begin{align*}
    \Decode(\inducedTreeof{0}{U},m):
    \Decode(\inducedTreeof{1}{U},m+m') 
    =
    \Decode(\inducedTreeof{0}{T},m):
    \Decode(\inducedTreeof{1}{T},m+m'')      
    \end{align*}
    If $\Decode(\inducedTreeof{0}{U}, m) = 
    \Decode(\inducedTreeof{0}{T}, m)$, then the induction
    hypothesis yields $\inducedTreeof{0}{U}=\inducedTreeof{0}{T}$. 
    Furthermore, this implies 
    $\Decode(\inducedTreeof{1}{U},m+m') =
    \Decode(\inducedTreeof{1}{T},m+m'')$ and $m'=m''$ whence
    by induction hypothesis 
    $\inducedTreeof{1}{U}=\inducedTreeof{1}{T}$.  
    In this case we conclude immediately that $T=U$. 
    
    The other case is that  the width of
    $\Decode(\inducedTreeof{0}{U},m)$ and the width of 
    $\Decode(\inducedTreeof{0}{T},m)$ do not coincide. 
    
    We prove that this case contradicts the assumption that
    $\Decode(U,m)=\Decode(T,m)$. 
    
    Let us assume that
    $\Decode(\inducedTreeof{0}{U},m) = 
    \Pop{2}^z\left(\Decode(\inducedTreeof{0}{T},m)\right)$ for
    some $z\in\N\setminus\{0\}$.  Note that this implies that the
    first word of $\Decode(\inducedTreeof{1}{U},m+m')$ is a word
    in $\Decode(\inducedTreeof{0}{T},m)$.
    
    Since $U(0)$ is a left successor in some tree from $\EncTrees$, it
    is labelled by some \mbox{$(\sigma,l)\in \Sigma\times\{1,2\}$.} 
    We make a case distinction on $l$.
    \begin{enumerate}
    \item Assume that $U(0) = (\sigma,2)$ for some
      $\sigma\in\Sigma$. Then 
      all words in $\Decode(\inducedTreeof{0}{T},m)$ start with
      the letter $(\sigma, 2, m)$. Thus, the first word of
      $\Decode(\inducedTreeof{1}{U},m+m')$ must also start with
      $(\sigma, 2, m)$. But all collapse links of level $2$ in
      $\Decode(\inducedTreeof{1}{U}, m+m')$ are at least
      $m+m'>m$. This is a contradiction.
    \item Otherwise,  $U(1) = (\sigma, 1)$ for some $\sigma\in\Sigma$. Thus, 
      all words in $\Decode(\inducedTreeof{0}{T},m)$ start with
      the letter $\sigma$. Thus, the first word of 
      $\Decode(\inducedTreeof{0}{U},m)$ and the first word of 
      $\Decode(\inducedTreeof{1}{U},m+m')$ have to start with
      $\sigma$. But this requires that
      \mbox{$U(0)=U(10)=(\sigma, 1)$}. This contradicts the assumption that
      $U$ is a proper subtree of a tree from $\EncTrees$
      (cf. condition 5 of Definition \ref{STACS:DefEncodingTrees}). 
    \end{enumerate}
    Both cases result in contradictions. Thus,  it is not the fact
    that there is  some $z\in\N\setminus\{0\}$ such that
    \begin{align*}
      \Decode(\inducedTreeof{0}{U},m) = 
      \Pop{2}^z\left(\Decode(\inducedTreeof{0}{T},m)\right)      
    \end{align*}
    By symmetry, we obtain that there is no $z\in\N\setminus\{0\}$
    such that
    \begin{align*}
      \Decode(\inducedTreeof{0}{T},m) = 
      \Pop{2}^z\left(\Decode(\inducedTreeof{0}{U},m)\right).
    \end{align*}
    Thus, we conclude that
    $\Decode(\inducedTreeof{0}{T},m) = 
    \Decode(\inducedTreeof{0}{U},m)$ whence $U=T$ as shown above. 
  \end{enumerate}
    If $\Decode(T,m)=\Decode(U,m)$, one of the previous cases applies:
  the following case distinction shows that 
  all other cases for 
  the defined or undefined subtrees of $T$ and $U$ imply
  $\Decode(T,m)\neq\Decode(U,m)$. 
  \begin{enumerate}
  \item \label{Case1110}
    Assume that $\inducedTreeof{0}{U}=\inducedTreeof{1}{U}=
    \inducedTreeof{0}{T}=\emptyset$ and 
    $\inducedTreeof{1}{T} \neq\emptyset$. In this case, 
    $\Decode(U,m)$ is $[\varepsilon]$ or
    $[\tau]$ for some
    $\tau\in \Sigma\cup(\Sigma\times\{2\}\times\N$). 
    Furthermore, 
    \begin{align*}
      \Decode(T, m) = f_m(T(\varepsilon)) \mathrel{\backslash} (\varepsilon:
      \Decode(\inducedTreeof{1}{T},m+1)).       
    \end{align*}
    It follows  that 
    $\lvert \Decode(T, m) \rvert \geq 2 > 
    \lvert \Decode(U, m) \rvert =  1$ whence
    \mbox{$\Decode(T,m)\neq\Decode(U,m)$.}
  \item \label{Case1101}
    Assume that 
    $\inducedTreeof{0}{U} = \inducedTreeof{1}{U} = \emptyset$, 
    $\inducedTreeof{0}{T} \neq \emptyset$, and 
    $\inducedTreeof{1}{T} = \emptyset$. 
    In this case,  $\Decode(U,m)$ is again $[\varepsilon]$ or
    $[\tau]$ for some
    $\tau\in\Sigma\cup(\Sigma \times\{2\}\times\N)$. 
    Since we assumed that $U(\varepsilon)= T(\varepsilon)$, 
    \begin{align*}
      \Decode(T, m) = f_m(T(\varepsilon)) \mathrel{\backslash}
      f_m(T(0)) \mathrel{\backslash} s      
    \end{align*}
    for some
    $2$-word $s$. Since $T$ is a subtree of a tree in $\EncTrees$,
    $T(0)\in\Sigma\times\{1,2\}$. Thus,
    $f_m(T(0))
    \in\Sigma\cup(\Sigma\times\{1,2\}\times\N)$.
    We conclude that the length of the first word of $\Decode(T, m)$
    is greater than the length of the first word of 
    $\Decode(U, m)$. Thus,
    \mbox{$\Decode(T, m)  \neq \Decode(U, m)$.}
  \item \label{Case1100} Assume that 
    $\inducedTreeof{0}{U} = \inducedTreeof{1}{U} = \emptyset$, 
    $\inducedTreeof{0}{T} \neq \emptyset$, and 
    $\inducedTreeof{1}{T} \neq \emptyset$. 
    Completely analogous to case
    \ref{Case1110}, we conclude that 
    $\lvert \Decode(T,  m) \rvert \geq 2 > 
    \lvert \Decode(U, m) \rvert =  1$ whence
    \mbox{$\Decode(T,m) \neq \Decode(U,m)$.} 
  \item \label{Case1011}
    Assume that 
    $\inducedTreeof{0}{U} = \emptyset$, 
    $\inducedTreeof{1}{U} \neq \emptyset$, and
    $\inducedTreeof{0}{T} =
    \inducedTreeof{1}{T} = \emptyset$. Exchanging the roles of $U$ and
    $T$, this is exactly the same as case \ref{Case1110}.
  \item \label{Case1001}
    Assume that
    $\inducedTreeof{0}{U} = \emptyset$, 
    $\inducedTreeof{1}{U} \neq \emptyset$, 
    $\inducedTreeof{0}{T} \neq \emptyset$, and
    $\inducedTreeof{1}{T} = \emptyset$.
    Analogously to case \ref{Case1101}, we derive that the length of
    the first word of $\Decode(T, m)$ is greater than the length of
    the first word of $\Decode(U, m)$. Thus,
    $\Decode(T, m) \neq \Decode(U, m)$.
  \item \label{Case1000}
    Assume that
    $\inducedTreeof{0}{U} = \emptyset$, 
    $\inducedTreeof{1}{U} \neq \emptyset$, 
    $\inducedTreeof{0}{T} \neq \emptyset$, and
    $\inducedTreeof{1}{T} \neq \emptyset$.
    Analogously to case \ref{Case1101}, we derive that the length of
    the first word of $\Decode(T, m)$ is greater than the length of
    the first word of $\Decode(U, m)$. Thus,
    $\Decode(T, m) \neq \Decode(U, m)$.
  \item \label{Case0111}
    Assume that
    $\inducedTreeof{0}{U} \neq \emptyset$, and
    $\inducedTreeof{1}{U} = 
    \inducedTreeof{0}{T} =
    \inducedTreeof{1}{T} =\emptyset$.
    Exchanging the roles of $U$ and
    $T$, this is exactly the case \ref{Case1101}.
  \item \label{Case0110}
    Assume that
    $\inducedTreeof{0}{U} \neq \emptyset$, 
    $\inducedTreeof{1}{U} = 
    \inducedTreeof{0}{T} =\emptyset$, and
    $\inducedTreeof{1}{T} \neq\emptyset$.
    Exchanging the roles of $U$ and $T$, this is exactly the case
    \ref{Case1001}.
  \item \label{Case0100}
    Assume that
    $\inducedTreeof{0}{U} \neq \emptyset$, 
    $\inducedTreeof{1}{U} = \emptyset$, 
    $\inducedTreeof{0}{T} \neq \emptyset$, and
    $\inducedTreeof{1}{T} \neq \emptyset$.
    In this case, 
    \begin{align*}
      &\Decode(U,m) = f_m(U(\varepsilon))\mathrel{\backslash}
      \Decode(\inducedTreeof{0}{U},m)\\
      \text{and }     
      &\Decode(T,m) = f_m(T)\mathrel{\backslash}
      \left(\Decode(\inducedTreeof{0}{T},m) :     
        \Decode(\inducedTreeof{1}{T},m+m')\right)      
    \end{align*}
    for some $m'\in\N\setminus\{0\}$. 
    Since $U(\varepsilon) = \varepsilon$ if and only if 
    $T(\varepsilon)= \varepsilon$, we conclude that  $U(\varepsilon) =
    T(\varepsilon)$.  
    Now, 
    \begin{align*}
      \Decode(\inducedTreeof{0}{U},  m) =
      \tau\mathrel{\backslash} u'
    \end{align*}
    for
    $\tau=f_m(U(0))
    \in\Sigma\cup(\Sigma\times\{2\}\times\{m\})$ and $u'$ some level
    $2$-word. We distinguish the following cases. 
    
    First assume that $\tau=(\sigma,2,m)$. 
    For all letters in $T'\coloneqq
    \Decode(\inducedTreeof{1}{T},m+m')$ of 
    collapse level $2$, the collapse link is greater or equal to $m+m'$. Hence,
    $T'$ does not contain a symbol $(\sigma,2,m)$ whence
    $\Decode(U,m) \neq \Decode(T,m)$. 
    
    Otherwise, $\tau\in\Sigma$. But then 
    $\Decode(U,m) = \Decode(T,m)$ would imply that
    \begin{align*}
      &\Decode(\inducedTreeof{0}{T},m) =
      \tau \mathrel\backslash T'\\
      \text{and }
      &\Decode(\inducedTreeof{10}{T},m+m') =
      \tau \mathrel\backslash T''      
    \end{align*}
    for certain nonempty level
    $2$-words 
    $T'$ and $T''$. But then $T(0)=T(10)=(\tau,1)$ 
    which contradicts the fact that $T$ is a subtree of some tree from
    $\EncTrees$. 
    
    Thus, we conclude that 
    $\Decode(T, m) \neq \Decode(U,m)$.
  \item  \label{0011} 
    Assume that 
    $\inducedTreeof{0}{U} \neq \emptyset$, 
    $\inducedTreeof{1}{U} \neq \emptyset$, and
    $\inducedTreeof{0}{T} =
    \inducedTreeof{1}{T} = \emptyset$. 
    Exchanging the roles of $U$ and $T$, this is the same as case
    \ref{Case1100}.
  \item \label{Case0010}
    Assume that
    $\inducedTreeof{0}{U} \neq \emptyset$, 
    $\inducedTreeof{1}{U} \neq \emptyset$,
    $\inducedTreeof{0}{T} = \emptyset$, and
    $\inducedTreeof{1}{T} \neq \emptyset$. 
    Exchanging the roles of $U$ and $T$, this is the same as case
    \ref{Case1000}.
  \item \label{Case0001}
    Assume that
    $\inducedTreeof{0}{U} \neq \emptyset$, 
    $\inducedTreeof{1}{U} \neq \emptyset$,
    $\inducedTreeof{0}{T} \neq \emptyset$, and
    $\inducedTreeof{1}{T} = \emptyset$. 
    Exchanging the roles of $U$ and $T$, this is the same as case 
    \ref{Case0100}.
  \end{enumerate}
  Hence, we have seen that 
  $\Decode(T, m)= \Decode(U, m)$ implies that each of the
  subtrees of $T$ is defined if and only if the corresponding subtree
  of $U$ is defined. 
  Under this condition, we concluded that $U=T$.
  Thus, the claim holds
  and the lemma follows as indicated above. 
\end{proof}

Next, we prove that $\Decode$ is a surjective map from
$\EncTrees$ to $Q\times\Stacks_2(\Sigma)$. This is done by induction
on the blocklines used to encode a stack.  
In this proof we use the notion of \emph{left-maximal} blocks and
\emph{good} blocklines. Let 
\begin{align*}
  s:\left(w\mathrel\backslash (w':b)\right) : s'  
\end{align*}
be a stack  
where $s$ and $s'$ are  $2$-words, $w$, and $w'$ are words, and 
$b$ is a $\tau$-block. We call $b$  \emph{left maximal} in this
stack if either $b=[\tau]$ or 
$b=\tau\tau'\mathrel{\backslash} b'$ such that $w'$ does not
start with $\tau\tau'$ for some
$\tau'\in\Sigma\cup(\Sigma\times\{2\}\times\N)$.  
We call a blockline in some stack \emph{good}, if its first block is left
maximal. Furthermore, we call the blockline starting with the block
$b$ \emph{left maximal} if $w'$ does not start with $\tau$.
Recall that the encoding of stacks works on left maximal blocks
and good blocklines.

\begin{lemma}
  $\Decode\circ\Encode$ is the identity on $Q\times\Stacks_2(\Sigma)$,
  i.e., $\Decode(\Encode(c)) = c$, for all $c\in Q\times\Stacks_2(\Sigma)$.
\end{lemma}
\begin{corollary} \label{DecSurjective}
  $\Decode:\EncTrees\rightarrow Q\times\Stacks_2(\Sigma)$ is surjective.
\end{corollary}
\begin{proof}[Proof of Lemma]
  Let $c=(q,s)$ be a configuration. 
  Since $\Decode$ and $\Encode$ encode and decode the state of $c$ in the
  root of $\Encode(c)$, it suffices  to show that
  \begin{align*}
    \Decode(\Encode(s, (\bot,1)),0) = s    
  \end{align*}
  for all stacks $s\in \Stacks_2(\Sigma)$. 
  We proceed by induction on
  blocklines of the stack $s$. For this purpose we reformulate the lemma
  in the following claim.
  \begin{claim}
    Let $s'$ be some stack which decomposes as
    $s'=s'':(w\mathrel\backslash b):s'''$ such that
    \mbox{$b\in(\Sigma\cup(\Sigma\times\{2\}\times\N))^{+2}$} is a good  
    $\tau$-blockline for some
    $\tau\in\Sigma\cup(\Sigma\times\{2\}\times\N)$.   
    Then
    \begin{enumerate}
    \item $\Decode(\Encode(b,\varepsilon),
      \lvert s''\rvert )= b'$ for the unique $2$-word $b'$ such that
      \mbox{$b=\tau\mathrel{\backslash}b'$} and 
    \item if $b$ is left maximal, then
      $\Decode(\Encode(b,(\sigma,l)),\lvert s''\rvert)=
      b$ where $\sigma=\Sym(\tau)$ and \mbox{$l=\Lvl(\tau)$.} 
    \end{enumerate}
  \end{claim}
  Note that the conditions in the second part require that 
  either $\tau\in\Sigma$ or $\tau=(\sigma,2,\lvert s'' \rvert)$ for
  some $\sigma\in\Sigma$. 

  The lemma follows from the second part of the claim because every stack
  is a left maximal $\bot$-blockline. 

  We prove both claims by parallel induction on the size of $b$. 
  As abbreviation we set
  \mbox{$g\coloneqq\lvert s'' \rvert$}. We
  write $\overset{(1)}{=}$($\overset{(2)}{=}$, respectively)  when
  some equality is due to the induction hypothesis of the first claim
  (the second claim , respectively). 
  The arguments for the first claim are as follows.
  \begin{itemize}
  \item If $b=[\tau]$ for
    $\tau\in\Sigma\cup(\Sigma\times\{2\}\times\N)$,  the claim  is trivially
    true because  $\Decode(\Encode(b,\varepsilon), g)=
    \Decode(\varepsilon,  g) = \varepsilon$. 
  \item  If there are $b_1,b_1'\in(\Sigma\cup(\Sigma\times\{2\}\times\N))^{*2}$
    such that 
    \begin{align*}
      b= [ \tau ] : b_1 = [\tau]
      :\left(\tau\mathrel{\backslash} b_1'\right)\text{ then}
    \end{align*}
    \begin{align*}
      &\Decode(\Encode(b,\varepsilon), g) = \Decode(
      \treeR{\varepsilon}{ \Encode(b_1,\varepsilon)},g)\\
      =   
      &f_g(\varepsilon) \mathrel{\backslash} \left(\varepsilon:
      \Decode(\Encode(b_1,\varepsilon),g+1)\right) \\
      \overset{(1)}{=}
      &\varepsilon \mathrel{\backslash}(\varepsilon : b_1') =
      \varepsilon:b_1'=b'.
    \end{align*}
  \item Assume that there is some
    $\tau'\in\Sigma\cup(\Sigma\times\{2\}\times\N)$ and some 
    $b_1\in (\Sigma\cup (\Sigma\times\{2\}\times\N))^{*2}$ such that
    \begin{align*}
      b=\tau\tau' \mathrel{\backslash} b_1.      
    \end{align*}
    The assumption that $b$ is good implies that  the blockline
    $\tau'\mathrel\backslash b_1$ is left maximal 
    whence
    \begin{align*}
      &\Decode(\Encode(b,\varepsilon),g) = 
      \Decode(    
      \treeL{\varepsilon}{
      \Encode( \tau'\mathrel{\backslash}
      b_1, (\Sym(\tau'), \Lvl(\tau')))}, g)\\
      = &f_g(\varepsilon) \mathrel{\backslash} 
      \Decode(\Encode( \tau'\mathrel{\backslash} b_1,
      (\Sym(\tau'), \Lvl(\tau')), g))\\
      \overset{(2)}{=}  &\tau' \mathrel{\backslash} b_1 =
      b'.
    \end{align*}
  \item The last case is that
    \begin{align*}
      b = \tau \mathrel{\backslash}
      \left((\tau' \mathrel{\backslash} b_1):b_2\right)
    \end{align*}
    for  $b_2$ a blockline of $s$ not starting with
    $\tau'$. By this we mean that
    $b_2\neq \tau' w':b_2'$ for any word $w'$ and any $2$-word
    $b_2'$. 
    Since $b$ is good,  
    \mbox{$\tau' \mathrel{\backslash} b_1$} is a left maximal
    blockline. Furthermore, $\tau \mathrel{\backslash} b_2$ is a
    good blockline.  
    Thus, 
    \begin{align*}
    &\Decode(\Encode(b,\varepsilon), g)\\
    =
    &\Decode\left(
    \treeLR{\varepsilon}{
      \Encode\left( \tau'\mathrel{\backslash}
      b_1, (\Sym(\tau'), \Lvl(\tau'))\right)}{
      \Encode(\tau\mathrel{\backslash}b_2,\varepsilon)},g\right)\\
    = 
    &f_g(\varepsilon)\mathrel{\backslash} \\
    &\left(\Decode(\Encode(\tau'\mathrel{\backslash} b_1,
    ( \Sym(\tau'), \Lvl(\tau'))),g) : 
    \Decode(\Encode(\tau\mathrel{\backslash} b_2,
    \varepsilon),g + f)\right), 
    \end{align*}
    where
    \begin{align*}
      f = \lvert
    \Decode(\Encode(\tau'\mathrel{\backslash} b_1, 
    (\Sym(\tau'), \Lvl(\tau'))),g)  \rvert  \overset{(2)}{=}
    \lvert b_1\rvert.  
    \end{align*}
    From this, we obtain that
    \begin{align*}
      &\Decode(\Encode(b,\varepsilon),g) \\
      \overset{(1)}{=}
      &\varepsilon \mathrel{\backslash}\big(
      (\tau'\mathrel{\backslash} b_1) :
      \Decode(\Encode(\tau \mathrel{\backslash} b_2,
      \varepsilon), g+f)\big) \\
      \overset{(2)}{=}
      &(\tau'\mathrel{\backslash} b_1) : b_2 = b'.
    \end{align*}
  \end{itemize} 
  For the proof
  of the second claim, note that the calculations are basically the
  same, but $f_g(\varepsilon)$ is replaced by $f_g(\sigma,l)$.
  Thus, if $l=1$ then $f_g(\sigma,l) = \sigma = \tau$. 
  For the case  $l=2$, recall that $g=\lvert s'' \rvert$ whence
  $f_g(\sigma,l) = (\sigma,2,\lvert s'' \rvert)$. 
  Note that $\Lnk(\tau)=\lvert s''\rvert$ due to the left maximality
  of $b$. 
  
  Thus, one proves the second case using
  the same calculations, but replacing
  $\varepsilon$ by $\tau$.
\end{proof}

The previous lemmas provide a proof of
Lemma \ref{STACS:Bijective}: we have shown that $\Decode$ is
bijective and it is the inverse of $\Encode$.

\subsection{Recognising Reachable Configurations}
\label{STACS:SecCertificates}

In this section, we show that $\Encode$ maps the reachable
configurations of a given collapsible pushdown system to a regular set. 

Fix a configuration $c=(q,s)$.
Recall that every run from the initial configuration to some stack $s$
has to pass each of the generalised  milestones $\genMilestones(s)$ of
$s$ (cf. Section \ref{subsec:Milestones}). Especially, the set of
milestones 
$\Milestones(s)\subseteq\genMilestones(s)$ has a close connection to
our encoding:
with every $d\in \Encode(c)$, we can associate a subtree of
$\Encode(c)$ which encodes a milestone. 
Via this correspondence, 
the substack relation  
on the milestones corresponds exactly to the lexicographic order of
the elements of $\Encode(c)$. 

We show the regularity of the set of encodings of reachable
configurations as follows. Given the tree $\Encode(c)$, we annotate
each node $d\in\Encode(c)$ with a state $q_d$. This annotation
represents the claim that there is a run from the initial
configuration to $c$ that passes the milestone associated with $d$ in
state $q_d$. Then we show that an automaton can check the
correctness of such an annotation. 
Since this annotation can be generated nondeterministically by an
automaton, it follows that the set of encodings of reachable
configurations is regular. 

The correspondence between nodes of $\Encode(s)$ and milestones of $s$
is established via the notion of \emph{the left stack induced by}
$d\in\domain(\Encode(s))$.  
This left stack is the decoding of the subtree of
$\Encode(s)$ which contains all nodes that are lexicographically
smaller than $d$. We show that these left stacks always
form milestones and that each milestone can be represented by such an
element.

\begin{definition}\label{STACS:DefLeftTree}
  Let $T\in\EncTrees$ be a tree and $d\in
  T\setminus\{\varepsilon\}$. Then the \emph{left and
  downward closed tree of $d$} is
  $\LeftTree{d,T} \coloneqq T{\restriction}_D$ where 
  \mbox{$D\coloneqq\{d'\in T: d'\leq_{\mathrm{lex}} d\}\setminus\{\varepsilon\}$}. 
  Then we denote by 
  $\LeftStack(d,T)\coloneqq\Decode(\LeftTree{d,T}, 0)$ the 
  \emph{ left stack induced by $d$}. 
  If $T$ is clear from the context, we omit it.
\end{definition}
\begin{remark}
  \label{TOP2Determined}
  We exclude the case $d=\varepsilon$ from the definition because the
  root encodes the state of the configuration and not a part of the
  stack. In the following, we are often interested in the stack
  encoded in a tree, whence we will consider all nodes except for the
  root of the encoding tree. 

  Recall that
  \mbox{$w\coloneqq\TOP{2}(\LeftStack(d,s)){\downarrow_0}$} is
  $\TOP{2}(\LeftStack(d,s))$ where all level $2$ links are set to $0$
  (cf. Definition \ref{def:downarrowzero}). 
  Due to the definition of the encoding, for every
  $d\in\domain(\Encode(s))$,
  $w$ is
  determined by the path 
  from the root to $d$: 
  interpreting $\varepsilon$ as empty word, the word along this
  path contains the pairs of      
  stack symbols and collapse levels of the letters of 
  $\TOP{2}(\LeftStack(d,s))$. Since all level $2$ links in $w$ are
  $0$, $w$ is determined by this path. 
  Thus, Proposition \ref{Prop:AutomatonForLoops} implies that there is
  an automaton that calculates at each position $d\in\Encode(q,s)$ the
  number  of possible loops of $\LeftStack(d, \Encode(q,s))$ with given
  initial and final state. 
\end{remark}
\begin{remark}
  $\LeftStack(d,\Encode(q,s))$ is a substack of $s$ for all
  $d\in\domain(\Encode(q,s))$.
  This observation follows from Remark \ref{milestonesInEncoding}
  combined with the fact that the left stack is induced by a  
  lexicographically downward closed subset. 
\end{remark}

\begin{lemma}
  Let $q\in Q$ and $s\in \Stacks_2(\Sigma)$. For each $d\in
  \Encode(q,s)\setminus\{\varepsilon\}$ we have
  $\LeftStack(d,\Encode(q,s))\in\Milestones(s)$.
  Furthermore, for each
  $s'\in\Milestones(s)$ there is some
  \mbox{$d\in\Encode(q,s)\setminus\{\varepsilon\}$} such that
  \mbox{$s'=\LeftStack(d,\Encode(q,s))$.}
\end{lemma}

\begin{proof}
  For the first claim, let
  $d\in\domain(\Encode(q,s))\setminus\{\varepsilon\}$.  We already
  know that \mbox{$s_d:=\LeftStack(d, \Encode(q,s))$} is a substack of $s$. 

  Recall that the path from the root to $s_d$ encodes
  $\TOP{2}(s_d)$. Furthermore, by definition of $\Encode$, $d$
  corresponds to some maximal block $b$ occurring in $s$ in the
  following sense:
  there are $2$-words $s_1, s_2$ and a word $w$ such that
  $s=s_1: (w\mathrel{\backslash} b) : s_2$ and such that the subtree
  rooted at $d$ encodes $b$. Moreover, $d$ encodes the first letter of
  $b$, i.e., if $b$ is a $\tau$-block, then the path from the
  root to $d$ encodes $w\tau$. 

  Note that by maximality of $b$, the greatest common prefix of the
  last word of $s_1$ and the first word of $w\mathrel{\backslash} b$
  is a prefix of $w\tau$. 

  Since the elements that are lexicographically smaller than $d$
  encode the 
  blocks to the left of $b$, one sees that $s_d=s_1:w\tau$. 
  Setting $k:=\lvert s_d \rvert$, we conclude that $s_d$ is a substack
  of $s$ such that the greatest common prefix of the $(k-1)$-st and
  the $k$-th word of $s$ is a prefix of $\TOP{2}(s_d)$. 

  Recall that this is exactly a characterisation of a milestone of
  $s$. Thus, $s_d$ is a milestone of $s$ and  we completed the proof of
  the first claim.

  Now, we turn to the second claim. 
  The fact that every milestone $s'\in\Milestones(s)$ is indeed represented by
  some node of $\Encode(q,s)$ can be seen by induction on the block
  structure of $s'$. 
  Assume that $s'\in\Milestones(s)$ and that $s'$  decomposes as 
  $s'=b_0: b_1 : \dots: b_{m-1}: b_m'$ into maximal blocks. We claim
  that $s$ then 
  decomposes as $s=  b_0: b_1: \dots:  b_{m-1}:b_m:\dots:b_n$ into
  maximal  blocks. 
  In order to verify this claim, we have to prove
  that $b_{m-1}$ cannot be the initial segment of a larger block
  $b_{m-1}:b_m$ in $s$. Note that if $b_{m}'$ only contains one letter, then 
  by definition of a milestone the last word of $b_{m-1}$ and the
  first word occurring in $s$ after $b_{m-1}$, which is the first word of
  $b_m$,  can only have a common prefix of length at most $1$. Hence,
  their composition does not form a block. 
  Otherwise, the first word of $b_m'$ contains two letters which do
  not coincide with the first two letters of the words in $b_{m-1}$. 
  Since this word is by definition a prefix of the first word in
  $b_m$, we can conclude again that $b_{m-1}:b_m$ does not form a block.
  

  Note that all words in the blocks $b_i$ for $1\leq i \leq n$ and in
  the block $b'_m$ share the same first letter which is encoded at the
  position $0$ in $\Encode(q,s)$ and in $\Encode(q,s')$. 
  By the definition of $\Encode(q,s)$ the blockline induced by 
  $b_i$ is encoded in the subtree rooted at $01^i0$ in $\Encode(q,s)$. 
  For $i<m$ the same holds in $\Encode(q,s')$. We set $d:=01^m$. 
  Note that $\Encode(q,s')$ and $\Encode(q,s)$ coincide on all
  elements that are lexicographically smaller than $d$ (because these
  elements encode the blocks $b_1:b_2:\dots b_{m-1}$.  

  Now, we distinguish the following cases.
  \begin{enumerate}
  \item Assume that $b_m'= [\tau]$ for
    $\tau\in\Sigma\cup(\Sigma\times\{2\}\times\N)$. Then the block $b_m'$
    consists of only one letter. In this case
    $d$ is the lexicographically largest element of $\Encode(q,s')$
    whence
    $s'=\LeftStack(d,\Encode(q,s')) = \LeftStack(d,\Encode(q,s))$. 
  \item Otherwise, there is a
    $\tau\in\Sigma\cup(\Sigma\times\{2\}\times\N)$ such that
      \begin{align*}
        b_m=  \tau \mathrel{\backslash}(c_0: c_1:\dots:
        c_{m'}:\dots: c_{n'}) \text{ and}\\
        b_m'= \tau\mathrel{\backslash}(c_0: c_1: \dots:
        c_{m'-1}:c'_{m'})        
      \end{align*}
     for some $m' \leq n'$ such that $c_0:c_1:\dots: c_{n'}$ are the
     maximal blocks of the blockline 
     induced by $b_m$  
     and \mbox{$c_0:c_1:\dots c_{m'-1}: c'_{m'}$} are the maximal
     blocks of the 
     blockline induced by $b'_{m}$.
     Now, \mbox{$c_1:c_2:\dots: c_{m'-1}$} are encoded in the subtrees
     rooted at $d01^i0$ for $0\leq i \leq m'-1$ in $\Encode(q,s)$ as
     well as in $\Encode(q,s')$.  $c_{m'+1}: c_{m'+2}:\dots :c_{n'}$ is
     encoded in the subtree rooted at $d01^{m'+1}$ in $\Encode(q,s)$
     and these elements are all lexicographically larger than
     $d01^{m'}0$. Hence, we can set $d':=d01^{m'}$ and repeat this
     case distinction on $d', c'_{m'}$ and $c_{m'}$ instead of $d,
     b'_m$ and $b_m$. 
  \end{enumerate}
  Since $s'$ is finite, by repeated application of the case distinction,
  we will eventually end up in the first case where we find a
  $d\in\Encode(q,s)$ such that $s'=\LeftStack(d,\Encode(q,s))$. 
\end{proof}

The next lemma states the tight connection between milestones of a stack
(with substack relation) and elements in the encoding of this stack 
(with lexicographic order).

\begin{lemma} \label{LemmaOrderIso}
  The map  
  \begin{align*}
    g: \domain(\Encode(q,s))\setminus\{\varepsilon\}
    &\rightarrow  \Milestones(s)\\  
    d &\mapsto \LeftStack(d,s)
  \end{align*}
  is an order isomorphism between
  $\left(\domain(\Encode(q,s))\setminus\{\varepsilon\},
    \leq_{\mathrm{lex}}\right)$ and 
  $\big(\Milestones(s),\leq\big)$.     
\end{lemma}  
\begin{proof}
  If the successor of $d$ in lexicographic order is $d0$, then the
  left stack of the latter
  extends the former by just one letter. Otherwise, the
  left and downward closed tree of the successor of $d$
  contains more elements ending in 
  $1$, whence it encodes a stack of larger width. Since each left and
  downward closed tree
  induces a milestone, it follows that $g$ is an order isomorphism. 
\end{proof}

Recall that by Lemma \ref{Cor:OrderEmbedding}, each run to a
configuration $(q,s)$ visits the milestones of $s$ in the order given
by the substack relation. With the previous lemma, this translates
into the fact that the left stacks induced by the elements of
$\Encode(q,s)$ are visited by the run in lexicographical order of the
elements of $\Encode(q,s)$. 

This gives rise to the following algorithm for identifying reachable
configurations of a collapsible pushdown system $\mathcal{S}$:
we label each node $d$ of the encoding with a state $q_d$. 
Let $s_d$ be the left stack induced by each $d$. Fix a $d$ and let $d'$
be the lexicographical successor of $d$. 
Then we
check whether 
there is a run from $(q_d,s_d)$ to $(q_{d'}, s_{d'})$. 

In the next section we show that
this check depends only on the local structure of the encoding 
of a configuration. Hence, an automaton can do this check.

\paragraph{Detection of Reachable Configurations}
We have already seen that every run to a valid configuration $(q,s)$
passes all the milestones of $s$. Now, we use the last states in 
which a run $\rho$ to $(q,s)$ visits the milestones as a certificate for
the reachability of $(q,s)$.

\begin{definition}
  Let $(q,s)$ be some configuration and
  $\rho$ a run from the initial configuration to $(q,s)$. 
  \emph{The certificate for the
    reachability of $(q,s)$} induced by $\rho$ is the map
  \mbox{$C_\rho:\domain(\Encode(q,s))\setminus\{\varepsilon\}
    \rightarrow Q$}  such
  that $d\mapsto \hat q$ if and only if $\rho(i)=(\hat q,\LeftStack(d))$
  and $i$ is 
  the maximal position in $\rho$ where $\LeftStack(d, \Encode(q,s))$
  is visited. 
\end{definition}
\begin{remark}
  In the following, we identify a function
  $f:\domain(\Encode(q,s))\setminus\{\varepsilon\} \rightarrow Q$ with
  the $Q\cup\{\Box\}$-labelled tree 
  \mbox{$\hat f:\domain(\Encode(q,s)) \rightarrow Q \cup\{\Box\}$} where
  \mbox{$\hat f(d) =
  \begin{cases}
    \Box &\text{if }d=\varepsilon,\\
    f(d) &\text{otherwise.}
  \end{cases}$}
\end{remark}

In the following, we analyse the existence of certificates for
reachability. The existence of certain loops plays an important role
in this analysis.   
Thus, we first fix some notation concerning the existence of returns
and loops. 
Recall that we defined the functions $\ReturnFunc{k}$, $\LoopFunc{k}$, 
etc. (cf. Definitions \ref{Def:ReturnFunc} and \ref{Def:LoopFunc})
that count up to threshold 
$k$ the number of returns and loops starting in a given
configuration. Recall that there is a loop from $(q,s)$ to $(q',s)$ if
and only if $\LoopFunc{k}(s)(q,q')\geq 1$ for all $k\geq 1$. 
\begin{definition}
  We set 
  \begin{align*}
    \exLoops(s):=\{(q,q')\in Q\times Q:\LoopFunc{1}(s)(q,q')=1\}.    
  \end{align*}
  $\exLoops(s)$ contains those pairs of
  states $q,q'$ such that there exists at least one loop from $(q,s)$
  to $(q',s)$. 
  Completely analogously, we set 
  \begin{align*}
    &\exHighLoops(s):=\{(q,q')\in Q\times Q:
    \HighLoopFunc{1}(s)(q,q')=1\},\\ 
    &\exLowLoops(s):=\{(q,q')\in Q\times Q:
    \LowLoopFunc{1}(s)(q,q')=1\}\text{ and}\\ 
    &\exReturns(s):=\{(q,q')\in Q\times Q: \ReturnFunc{1}(s)(q,q')=1\}.    
  \end{align*}
  These sets contain the pairs of initial and final states of
  low loops, high loops and returns starting with stack $s$. 
\end{definition}
\begin{remark}
  Due to Remark \ref{TOP2Determined} and
  due to  Proposition
  \ref{Prop:AutomatonForLoops}, 
  the function  that assigns
  \begin{align*}
    d\mapsto \exLoops(\LeftStack(d,\Encode(q,s)))   
  \end{align*}
  is calculated by some automaton for all
  configurations $(q,s)$. 
  Analogous, the function that  assigns
  \begin{align*}
    d\mapsto \exHighLoops(\LeftStack(d,\Encode(q,s)))   
  \end{align*}
  is also calculated by some automaton.
\end{remark}
Using this notation, we can prove the first important lemma concerning
certificates for reachability. 

\begin{lemma}\label{STACS:CertificateCheck}
  For every \CPG $G$, there is an automaton $\mathcal{A}$
  that checks for each map
  \begin{align*}
   f:\domain(\Encode(q,s))\setminus\{\varepsilon\} &\rightarrow Q 
  \end{align*}
  whether $f$ is a certificate for the reachability of $(q,s)$. This
  means 
  that $\mathcal{A}$ accepts  $\Encode(q,s)\otimes f$ if $f=C_\rho$
  for some run $\rho$ from the initial configuration to $(q,s)$. 
\end{lemma}
\begin{proof}
  As before, we identify $f$ with a $Q\cup\{\Box\}$-labelled tree
  encoding $f$.  
  Due to the previous remark, it is sufficient to prove that there is
  an automaton 
  which accepts 
  \begin{align*}
    &\Encode(q,s)\otimes f \otimes T_{\mathrm{Lp}} \otimes
    T_{\mathrm{HLp}}\\
    &\text{if and only if }  f=C_\rho\text{ for some run }\rho,
  \end{align*}
  where $T_{\mathrm{Lp}}$ is  a tree encoding the value of 
  $\exLoops(\LeftStack(d,\Encode(q,s)))$ at each node
  \mbox{$d\in\domain(\Encode(q,s))$} and 
  $T_{\mathrm{HLp}}$ is a tree encoding the value of
  $\exHighLoops(\LeftStack(d,\Encode(q,s))$. 
  We write 
  $T:=\Encode(q,s)$ 
  as an abbreviation. 
  We start with an informal description what we have to check at  some
  node $d\in\domain(T)$.  
  According to Corollary \ref{ClonePopMilestoneRun}, it is sufficient
  to check the following facts.
  \begin{enumerate}
  \item Assume that $d,d0\in \domain(t)$. 
    We know that $\LeftStack(d,T)=\Pop{1}(\LeftStack(d0,T))$. 
    By definition, we know that $f$
    can only be a certificate for reachability if there is a run
    $\rho'$ from  
    $\left(f(d),\LeftStack(d,T)\right)$ to
    $\left(f(d0),\LeftStack(d0,T)\right)$ that starts with some 
    push operation followed by a high loop of $\LeftStack(d0,T)$.
    This requirement can be checked by an automaton  
    when it reads the labels $t(d)$ and $t(d0)$ as follows.

    We assume that the automaton has stored the information about the
    topmost symbol $\sigma$ of $\LeftStack(d,T)$.
    When it reads $t(d)$ it guesses nondeterministically a pair
    $(q', (\sigma',i))$ for $q'\in Q$, $\sigma'\in\Sigma$ and
    $i\in\{1,2\}$ such that there is a $\Push{\sigma',i}$ transition from 
    state $f(d)$ and topmost symbol $\sigma$ going to state $q'$. 
    Reading the label $t(d0)$ it  checks whether $\Encode(q,s)(d0)
    =(\sigma',i)$ and whether
    $\left(q',f(d0)\right)\in
    \exLoops\left(\LeftStack(d0,T)\right)$. If this is the     
    case then the automaton  guessed the right push transition and
    there is a run from $(f(d), \LeftStack(d,T))$ 
    to 
    $(f(d0), \LeftStack(d0, T))$.
  \item  Consider the case where $d\in \domain(T)$ but
    $d0\notin\domain(T)$ and where $d$ has a successor $d'$ in
    lexicographic order. This implies that the direct successor of
    $\LeftStack(d,T)$ in $\Milestones(s)$ is of the form 
    \begin{align*}
      s' \coloneqq \Pop{1}^m(\Clone{2}(\LeftStack(d,T))).      
    \end{align*}
    In this case there is a maximal prefix $d_0\leq d$ and some
    $d_1\in\{0,1\}^*$ such that 
    $d=d_0 0 d_1$ and $d' = d_0 1\in\domain(T)$.
    Due to Lemma \ref{LemmaOrderIso}, 
    we know that $s'=\LeftStack(d',T)$. 
    
    From our observations about milestones we know that we have to
    verify that there 
    is some run $\rho' \coloneqq \rho_0 \circ \lambda_0 \circ \rho_1
    \circ \lambda_1 \circ \rho_2 \circ \lambda_2 \dots \rho_{m} \circ
    \lambda_m$ where $\lambda_i$ is a loop for all $0\leq i \leq m$  and 
    $\rho_0$ is a run that 
    performs  one clone operation and for $j>0$ the run $\rho_j$ 
    performs  either one $\Pop{1}$ or one $\Collapse$ of level $1$
    such that $\rho'$ starts in 
    $\left(f(d),\LeftStack(d,T)\right)$ and ends in  $\left(f(d'),
      \LeftStack(d',T)\right)$.  
    
    An automaton can verify this because the 
    path from $d_0$ to
    $d$ encodes the topmost stack symbols and
    collapse levels of 
    $\Pop{1}^{m'}(\Clone{2}(\LeftStack(d,T)))$ for $m'\leq m$.
    Since the
    existence of loops only depends on the topmost word, an automaton can
    check the existence of $\rho'$ while processing the path from $d$ to
    $d_0$.
  \item  Finally, we have to consider the lexicographically minimal and
    maximal element in the encoding of the stack. 
    Let $d$ be the rightmost leaf of $T$.
    Recall that \mbox{$\LeftStack(d, T)=s$}. $f$ can
    only be a certificate for reachability for $(q,s)$ if it labels $d$
    with the last state in which $s$ is visited. But if $\rho$ is a
    run to $(q,s)$ then this last state must be $q$. Thus, the
    condition for the rightmost leaf $d$ is that $f(d) = q$. 

    Recall that $\LeftStack(0, T)=[\bot]$. 
    Due to Corollary \ref{ClonePopMilestoneRun}, the run starts with a
    loop from the initial configuration to some configuration 
    $(\hat q, [\bot])$. Hence, we have to check whether $f(0)=\hat q$. 
  \end{enumerate}
  The lemma claims that there is an automaton $\mathcal{A}$  checking
  these conditions. 
  Instead of a concrete
  construction of $\mathcal{A}$, we
  present an \MSO formula $\chi$ that checks at each node
  $d\in\domain(\Encode(q,s))$ the corresponding condition.  Due to the
  correspondence between  
  \MSO definability and automata recognisability, 
  the automaton $\mathcal{A}$ can be constructed from this formula
  using standard constructions. 
  \begin{enumerate}
  \item For the first condition consider the formula 
    \begin{align*}
      &\chi_1:=\forall x \forall y (\neg \mathrm{Root}(x) \land y=x0)
      \rightarrow\\ 
      &\left(\bigvee_{(q_1,\sigma,\gamma, q_2, \Push{\tau,i})\in \Delta}
        \mathrm{Sym}(x)=\sigma \land \mathrm{Top}(y)=(\tau,i)
        \land
        f(x)=q_1 \land  (q_2,f(y))\in \mathrm{HLp}(y)
      \right ),
    \end{align*}
    where 
    \begin{itemize}
    \item $\mathrm{Root}(x)$ is the formula stating that $x$ is the
      root of the tree, i.e., $x$ has no predecessor,
    \item $\mathrm{Sym}(x)=\sigma$ is  an \MSO formula stating
      that the maximal $1$-ancestor $z$ of $x$ satisfies
      $\Encode(q,s)(z)=(\sigma,i)$ for some $i\in\{1,2\}$, 
    \item $\mathrm{Top}(y)=(\tau,i)$ is a formula stating that 
      $\Encode(q,s)(y)=(\tau, i)$, and
    \item $(q_2,f(y))\in \mathrm{HLp}(y)$  asserts that
      $T_{\mathrm{HLp}}(q_2,f(y)) = 1$, i.e., it asserts that
      $(q_2, f(y))\in \exHighLoops(\LeftStack(y,T))$. 
    \end{itemize}
    This formula asserts exactly the
    conditions of the first case at all nodes $x$ that have a left
    successor. Note that we exclude the root of the tree because it
    encodes the state of the configuration and not a part of the
    stack. 
  \item For the second case, 
    let 
    $\varphi(x,y,X)$ be an
    \MSO formula that is valid if $x$ does not have a left successor,
    if $y$ is the successor of  $x$ with respect to lexicographic
    ordering and if $X$
    contains the path connecting the predecessor of $y$ with $x$. 
    
    Assume that there is a triple $(x,y,X)$ that satisfies $\varphi$ on
    $\Encode(q,s)$. Then  there
    are  a  node $z\in\domain(\Encode(q,s))$,  a number $k\in\N$ and
    numbers $n_1, n_2, \dots, n_k\in\N$   
    such that $y=z1$ and $x=z01^{n_1}01^{n_2} \dots 01^{n_k}$. 
    Then
    $X=\{a: z\leq a \leq x\}$. 
    For each node $a\in X$, there is some $0\leq l \leq k$ and a
    number $n_l'\leq n_l$ such that
    $a=y01^{n_1}01^{n_2} \dots 01^{n_{l-1}}01^{n'_l}$. 
    Since the path to $a$ encodes the topmost word of the left stack
    induced by $a$, setting $k_a:=k-l$ we obtain that
    \begin{align*}
      \TOP{2}(\LeftStack(a, \Encode(q,s))) =
      \TOP{2}(\Pop{1}^{k_a}(\Clone{2}(\LeftStack(x,\Encode(q,s))))).       
    \end{align*}
    Furthermore, 
    \begin{align*}
      \LeftStack(y,\Encode(q,s)) =
      \Pop{1}^k(\Clone{2}(\LeftStack(x,\Encode(q,s)))).       
    \end{align*}
    We will use the following abbreviations: 
    \begin{align*}
      &s_a:=\LeftStack(a, \Encode(q,s))\text{ and}\\
      &\hat s_a:=\Pop{1}^{k_a}(\Clone{2}(\LeftStack(x,\Encode(q,s)))).      
    \end{align*}
    By definition, $\exLoops(s_a) =
    \exLoops(\hat s_a)$. 
    We use $a$ as  the representative for 
    $\hat  s_a$.  
  
    We next define a formula $\chi_2$. $\chi_2$ asserts the existence of a
    function 
    \mbox{$g:X\to Q$} that labels each node $a\in X$ with a state $q_a$ such
    that there is a $\Pop{1}$ or $\Collapse$ of level $1$ followed by
    a loop which connects  
    $(q_a,\hat s_a)$ with 
    $(q_b,\hat s_b)$ for
    $b\leq a$ some node such that $k_b = k_a+1$. 
    Furthermore, the formula asserts that there is a run from 
    $\left(f(x), s_x\right)$ to \mbox{$\left(g(x), \hat s_x\right) =
      \left(g(x), \Clone{2}(s_x)\right)$} and 
    it asserts that $g(z)=f(y)$.  
    Note that such a labelling $g$ is exactly a witness for a run 
    $\rho'=\rho_0 \circ \lambda_0 \circ \rho_1
    \circ \lambda_1 \circ \rho_2 \circ \lambda_2 \dots \rho_{m} \circ
    \lambda_m$ as described above. 

    Let $\chi_2$ be the formula
    \begin{align*} 
      &\forall x,y \forall X 
      \Big( \varphi(x,y,X) \rightarrow\\
      &\exists g: X\to Q\ \left(
        \bigvee_{(q_1,\sigma,\gamma,q_2,\Clone{2})\in\Delta} 
        \left(\mathrm{Sym}(x)=\sigma \land f(x)=q_1\land
          (q_2,g(x))\in\mathrm{Lp}(x)\right)
      \right.\\
      &\land \psi(g,X) \land \exists z (z1=y\land f(y)=g(z))\big)\Big)
    \end{align*}
    where  
    \begin{align*}
      &\psi(g,X):=  \forall v,z\in X
      \left((z=v1 \rightarrow g(z)=g(v)) 
      \land 
      \left(z=v0 \rightarrow (\psi_{p} \lor
          \psi_{c})\right)\right), \\
      &\psi_p(v,z):=\bigvee\limits_{(q_1,\sigma,\gamma,q_2, \Pop{1})\in\Delta} \left(
      \mathrm{Sym}(z)=\sigma \land g(z)=q_1 \land
      (q_2,g(v))\in \mathrm{Lp}_s(v)\right) \text{ and}\\
      &\psi_p(v,z):=\bigvee\limits_{(q_1,\sigma,\gamma,q_2, \Collapse)\in\Delta} \left(
      \mathrm{Top}(z)=(\sigma,1) \land g(z)=q_1 \land
      (q_2,g(v))\in \mathrm{Lp}_s(v)\right).
    \end{align*}
    Note that the function $g$ has finite range whence it may be
    encoded in a finite number of set-variables. Thus, $\chi_2$ can be
    formalised in $\MSO$.
  \item Let $\chi_3$ be the formula asserting that
    \begin{enumerate}
    \item the rightmost leaf
      $d$  of $\Encode(q,s)$ satisfies $f(d)=q$, and that
    \item  $(q_0,f(0))\in \exLoops([\bot])$, i.e., if
      $T_{\mathrm{Lp}}(0)(q_0, f(0)) = 1$. 
    \end{enumerate}

  \end{enumerate}
  Now, $\Encode(q,s) \otimes f\otimes T_{\mathrm{Lp}} \otimes
  T_{\mathrm{HLp}} \models \chi:= \chi_1\land \chi_2 \land \chi_3$ if
  and only if  
  $f=C_\rho$ for some run $\rho$ from the initial configuration to
  $(q,s)$. 
\end{proof}

Since regular tree-languages are closed under projection, 
there is an automaton that nondeterministically guesses the
existence of a 
certificate for reachability for each encoding of a reachable configuration. 

\begin{corollary}
 For every  collapsible pushdown system $\mathcal{S}$ of level $2$,
 there is an automaton $\mathcal{A}$ that accepts a 
 tree $T$ if and only if $T=\Encode(q,s)$ for a reachable 
 configuration $(q,s)$ of $\mathcal{S}$. 
\end{corollary}
\begin{proof}
  Note that $T=\Encode(q,s)$ for an arbitrary
  configuration if and only if $T\in \EncTrees$ which is a regular set.
  Furthermore, the set of encodings of reachable configurations forms
  a regular 
  subset of $\EncTrees$ due to the previous lemma and due to the closure of
  regular languages under projection. 
\end{proof}

\subsection{Regularity of the Stack Operations}
\label{Reg_Stack_OP}

In the previous section, we have seen that the function $\Encode$
translates the reachable configurations of a collapsible pushdown
graph $\mathcal{S}$ (of level $2$) into a regular
tree language.
In order to prove that
$\CPG(\mathcal{S})$ is automatic, we have to define automata
recognising the transition relations $\trans{\gamma}$ for
every $\gamma\in\Gamma$. 
In fact, we will prove that for each transition
$(q,\sigma,\gamma,q',\op)\in Q\times\Sigma\times\Gamma\times
Q\times\Op$ the set 
\begin{align*}
  \left\{(\Encode(q,s),\Encode(q',s')): \Sym(s)=\sigma\text{ and
    }\op(s)=s'\right\}  
\end{align*}
is regular.    
In preparation of this proof, we
analyse the relationship between 
the encodings of the stack $s$ and the stack $s':=\Pop{2}(s)$. 

\begin{figure}
  \centering
    \psset{unit=.9mm} 
    \begin{pspicture}(-15,0)(130,60)
      \psline(-15,15)(0,0)
      \psline(15,15)(0,0)
      \psline(15,15)(15,20) \rput(15,23){$(\sigma,l)$}
      \psline(15,25)(15,30)
      \psline(15,30)(10,35)
      \psline(15,30)(20,35)
      \psline(21,23)(25,23) \rput(27,23){$\varepsilon$}
      \psline(27,25)(27,30)
      \psline(27,30)(22,35)
      \psline(27,30)(32,35)
      \psline(30,23)(34,23) \rput(37,23){$\dots$}
      \psline(40,23)(45,23)\rput(47,23){$\varepsilon$}
      \psline(47,25)(47,30)
      \psline(47,30)(52,35)
      \psline(47,30)(42,35)
      \psline(50,23)(55,23)\rput(57,23){$\varepsilon$}
      \psline(57,25)(57,29)\rput(57,32){$(\tau,k)$}
      \psline(57,35)(57,39)\rput(57,44){$\vdots$}
      \psline(57,46)(57,50)\rput(57,53){$(\sigma',l')$}

      \psline(65,15)(80,0)
      \psline(95,15)(80,0)
      \psline(95,15)(95,20) \rput(95,23){$(\sigma,l)$}
      \psline(95,25)(95,30)
      \psline(95,30)(90,35)
      \psline(95,30)(100,35)
      \psline(101,23)(105,23) \rput(107,23){$\varepsilon$}
      \psline(107,25)(107,30)
      \psline(107,30)(102,35)
      \psline(107,30)(112,35)
      \psline(110,23)(114,23) \rput(117,23){$\dots$}
      \psline(120,23)(125,23)\rput(127,23){$\varepsilon$}
      \psline(127,25)(127,30)
      \psline(127,30)(132,35)
      \psline(127,30)(122,35)
    \end{pspicture}   
  \caption{$\Pop{2}$ operation in the tree-encoding.}
  \label{fig:Pop2}
\end{figure}
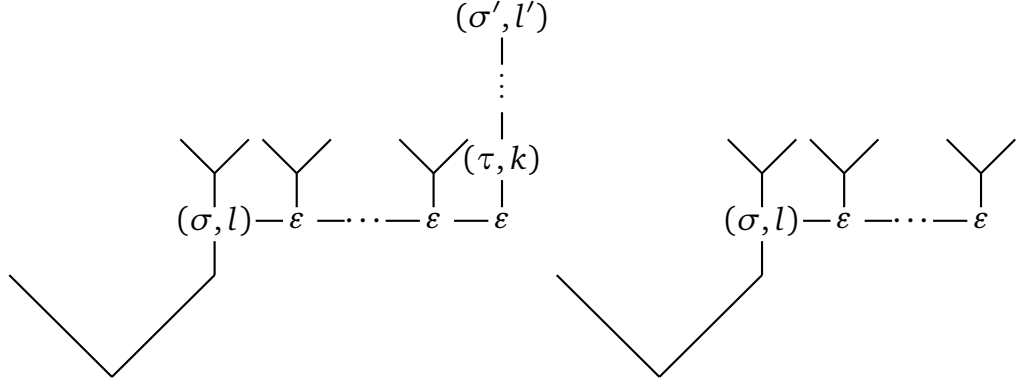

\begin{lemma} \label{LemmaStructurePop}
  Let $c=(q,s)$ and $c'=(q',s')$ be configurations of a pushdown
  system $\mathcal{S}$ such that
  \mbox{$s'=\Pop{2}(s)$}. 
  There is a unique element $t\in\Encode(c')$  
  such that  $t\in\Encode(c)\setminus\Encode(c')$  
  For
   \mbox{$D\coloneqq\{d\in\domain(\Encode(c)): t1\not\leq
    d\}$}, we have
  \begin{align*}
    &\domain(\Encode(c))\setminus
    \domain(\Encode(c'))\subseteq t10^*  
    \\
    \text{and }
    &\Encode(c')=\Encode(c){\restriction}_D
    \text{  (see  Figure \ref{fig:Pop2}).}
  \end{align*}
\end{lemma}

\begin{proof}
  The  proof is by induction on the structure of $\Encode(s, (\bot,1))$
  and $\Encode(s', (\bot, 1))$. 
  In fact, we prove the following stronger claim. 
  \begin{claim}
    Let $\tau\in(\Sigma\times\{1,2\} )\cup\{\varepsilon\}$.
    Let $t$ be the maximal
    element such that 
    \begin{itemize}
    \item $t$ is in the rightmost path of 
      $\domain(\Encode(s, \tau))$,
    \item    $t\in\domain(\Encode(s',\tau))$ and
    \item     $t1\in\domain(\Encode(s,
      \tau))\setminus\domain(\Encode(s',\tau))$.
    \end{itemize}
    Set
   $D\coloneqq\{d\in\domain(\Encode(s,\tau)): t1\not\leq
    d\}$. It holds that
  \begin{align*}
    &\domain(\Encode(s,\tau))\setminus
    \domain(\Encode(s',\tau))\subseteq t10^*  
    \\
    \text{and }
    &\Encode(s',\tau)=\Encode(s,\tau){\restriction}_D
    \text{  (see  Figure \ref{fig:Pop2}).}
  \end{align*}    
  \end{claim}
  Recall that for some stack consisting of just one word $w_1$,
  its encoding $\Encode(w_1,\varepsilon)$ is a path with $0$-edges
  only, i.e., $\domain(\Encode(w_1,\varepsilon))\subseteq \{0\}^*$.
  
  Let $s \coloneqq w_1: w_2: \dots: w_n : w_{n+1}$ and correspondingly 
  \mbox{$s'\coloneqq w_1: w_2: \dots: w_n$}. 
  In the case that $\lvert w_1 \rvert\geq 1$, 
  let $\tau_1, \tau_2\in \Sigma\cup(\Sigma\times\{2\}\times\N)$ 
  and $w_1'$ some word such that $w_1= \tau_1\tau_2 w_1'$.
  We prove the lemma by induction on the size of $s$. 
  We distinguish the following cases.
  \begin{enumerate}
  \item For all $i\leq n$ there are words $w_i'$ such that 
    $w_i= \tau_1 \tau_2 w_i'$, but 
    $\tau_1\tau_2 \not\leq w_{n+1}$.
    Then the root in $\Encode(s',\tau)$ has only a left successor
    and $\Encode(s,\tau)$ extends $\Encode(s',\tau)$ by a right
    subtree of the root which is $\Encode(w_{n+1},\varepsilon)$. 
    Due to our initial remark on the structure of the encoding of a
    single word, the claim follows immediately.
  \item For all $i\leq n+1$, there are words $w_i'$ such that 
    $w_i=\tau_1\tau_2 w_i'$.
    In this case $\Encode(s,\tau)$ and $\Encode(s',\tau)$
    coincide on their roots, these roots do not have right successors
    and the 
    subtrees induced by the left successor are 
    \begin{align*}
      &\Encode\left( \tau_2 \mathrel{\backslash} 
      (w_1': \dots : w_n': w_{n+1}'), \left(\Sym(\tau_2),
        \Lvl(\tau_2)\right)\right)\\  
      \text{and } &\Encode\left( \tau_2 \mathrel{\backslash} (w_1': \dots :
      w_n'), \left(\Sym(\tau_2), \Lvl(\tau_2)\right)\right).
    \end{align*} 
    Now,  
    we apply again the same case distinction to the subtrees encoding
    these parts of the stacks.  
  \item There is some $j<n$ such that 
    $\tau_1\tau_2 \leq w_i$ for all $i\leq j$ and
    $\tau_1\tau_2 \not\leq w_i$ for all $i>j$.
    In this case the claim of the lemma reduces to the claim that the
    lemma holds for $t\coloneqq w_{j+1}: w_{j+2} : \dots : w_n: w_{n+1}$ and
    $t'\coloneqq w_{j+1}: w_{j+2} : \dots : w_{n}$. Since the left
    subtrees of the encodings of $s$ and $s'$ agree and their right
    subtrees encode $t$ and $t'$, respectively, we can apply again this
    case distinction to $t$ and $t'$. 
  \item The last case is that $\lvert w_1 \rvert = 1$. If $n>1$, 
    the claim reduces 
    to the claim that the lemma holds for $t\coloneqq w_2 : w_3 :
    \dots : w_n:w_{n+1}$ and $t'\coloneqq w_2: w_3 : \dots : w_{n}$
    because $w_1$ is encoded in the root of $\Encode(s,\tau)$ and
    $\Encode(s',\tau)$ and the right subtree of the trees encode $t$
    and $t'$, respectively.

    If $n=1$, this leads to the fact that
    $\Encode(s',\sigma)$ is only a tree of one element and
    $\Encode(s,\sigma)$ extends this root by a right subtree, namely
    $\Encode(w_{n+1}, \varepsilon)$. In this case the lemma holds due
    to our initial remark.
  \end{enumerate}
  In each iteration of the case distinction, the stacks get smaller. 
  Thus, we eventually reach the first case or the last
  case with condition $n=1$. This observation completes the proof of
  the lemma.  
\end{proof}

\begin{figure}
  \centering
    \psset{unit=.9mm} 
    \begin{pspicture}(-15,0)(95,40)
      \psline(-15,15)(0,0)
      \psline(0,0)(10,10)
      \rput(10,15){$\vdots$}
      \psline(10,10)(12,10)
      \rput(15,10){$\varepsilon$}
      \psline(15,13)(15,15)
      \rput(15,20){$\vdots$}
      \psline(15,22)(15,27)
      \rput(15,30){$(\sigma,2)$}
      \psline(15,32)(15,37)
      \psline(15,37)(10,42)
      \psline(15,37)(20,42)
      \psline(21,30)(25,30) \rput(27,30){$\varepsilon$}
      \psline(27,32)(27,37)
      \psline(27,37)(22,42)
      \psline(27,37)(32,42)
      \psline(30,30)(34,30) \rput(37,30){$\dots$}
      \psline(40,30)(45,30)\rput(47,30){$\varepsilon$}
      \psline(47,32)(47,37)
      \psline(47,37)(52,42)
      \psline(47,37)(42,42)
      \psline(50,30)(55,30)\rput(57,30){$\varepsilon$}

      \psline(65,15)(80,0)
      \psline(90,10)(80,0)
      \rput(90,15){$\vdots$}
    \end{pspicture}   
  \caption{$\Collapse$ operation of level $2$ (if the collapse is of
  level $1$ then it is identical to the $\Pop{1}$ operation).}
  \label{fig:Col}
\end{figure}
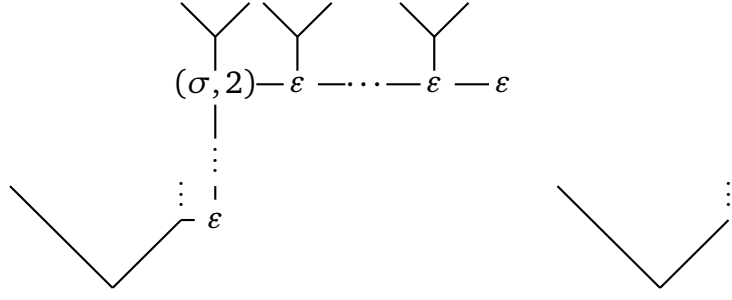

Analogously to the case of $\Pop{2}$, one proves a similar result 
for the $\Collapse$ operation:

\begin{lemma}\label{LemmaStructureCol}
  Let $s,s'$ be stacks of a pushdown system $\mathcal{S}$ such that
  $\Lvl(s)=2$ and 
  $s'\coloneqq\Collapse(s)$. Let $t'$ be the
  maximal element in the rightmost path of $\Encode(s, (\bot,1))$ which
  is labelled by some $(\sigma,2)$ for $\sigma\in\Sigma$. 
  Furthermore, let $t$ be the maximal ancestor of $t'$ such that
  $t1\leq t'$. 
  For 
  $D\coloneqq\left\{d\in\domain(\Encode(s, (\bot,1))): t1\not\leq
    d\right\}$, it 
  holds that 
  \begin{align*}
  \Encode(s',(\bot,1))=\Encode(s, (\bot,1)){\restriction}_D\text{ (see Figure
    \ref{fig:Col}).}
  \end{align*}
\end{lemma}
\begin{proof}
    Note that the rightmost leaf of $\Encode(s,(\bot,1))$ is of the
    form $t'1^n$ for some $n\in \N$. Hence, the topmost element of $s$
    is a clone of the element encoded at $t'$. Thus, 
    \begin{align*}
      \Collapse(s) =
      \Pop{2}(\LeftStack(t', \Encode(s,(\bot,1)))).      
    \end{align*}
    Using the previous
    lemma, the claim follows immediately.
\end{proof}

With these auxiliary lemmas we can now prove that $\Encode$ turns the
relations of collapsible pushdown graphs into automatic relations.  
\begin{lemma} \label{OperationsRegular}
  Let $\mathcal{S}=(Q, \Sigma,\Gamma, \Delta, q_0)$ be a collapsible
  pushdown system. 
  For each $\delta\in\Delta$, there is an automaton
  $\mathcal{A}_\delta$
  such that  for all configurations $c_1$ and  $c_2$
  \begin{align*}
    \mathcal{A}_{\delta} \text{ accepts } \Encode(c_1)\otimes \Encode(c_2)
    \text{\quad iff\quad } c_1\trans{\gamma} c_2.
  \end{align*}

\end{lemma}
\begin{figure}
  \centering
    \psset{unit=.9mm} 
    \begin{pspicture}(-15,0)(90,40)
      \psline(-15,15)(0,0)
      \psline(15,15)(0,0)
      \psline(15,15)(20,15) \rput(22,15){$\varepsilon$}

      \psline(35,15)(50,0)
      \psline(65,15)(50,0)
      \psline(65,15)(70,15) \rput(72,15){$\varepsilon$}
      \psline(74,15)(79,15) \rput(81,15){$\varepsilon$}
    \end{pspicture}\\
    \begin{pspicture}(90,0)(185,40)
      \psline(95,15)(110,0)
      \psline(125,15)(110,0)
      \rput(125,17){$(\sigma,l)$}

      \psline(135,15)(150,0)
      \psline(165,15)(150,0)
      \rput(165,17){$(\sigma,l)$}
      \psline(172,17)(177,17)
      \rput(179,17){$\varepsilon$}
    \end{pspicture}   
  \caption{The two versions of $\Clone{2}$ operations.}
  \label{fig:Clone}
\end{figure}
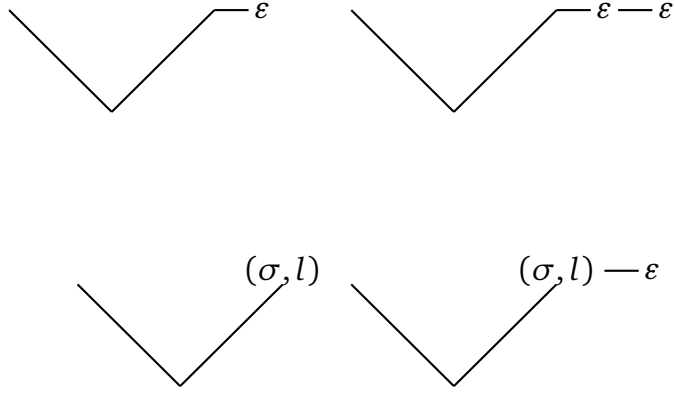

\begin{figure}
  \centering
    \psset{unit=.9mm} 
    \begin{pspicture}(-15,0)(135,40)
      \psline(-15,15)(0,0)
      \psline(15,15)(0,0)
      \psline(15,15)(15,20) \rput(15,23){$(\sigma,1)$}
      \psline(15,25)(15,30)
      \psline(15,30)(10,35)
      \psline(15,30)(20,35)
      \psline(21,23)(25,23) \rput(27,23){$\varepsilon$}
      \psline(27,25)(27,30)
      \psline(27,30)(22,35)
      \psline(27,30)(32,35)
      \psline(30,23)(34,23) \rput(37,23){$\dots$}
      \psline(40,23)(45,23)\rput(47,23){$\varepsilon$}
      \psline(47,25)(47,30)
      \psline(47,30)(52,35)
      \psline(47,30)(42,35)
      \psline(15,15)(50,15)\rput(52,15){$\varepsilon$}

      \psline(65,15)(80,0)
      \psline(95,15)(80,0)
      \psline(95,15)(95,20) \rput(95,23){$(\sigma,1)$}
      \psline(95,25)(95,30)
      \psline(95,30)(90,35)
      \psline(95,30)(100,35)
      \psline(101,23)(105,23) \rput(107,23){$\varepsilon$}
      \psline(107,25)(107,30)
      \psline(107,30)(102,35)
      \psline(107,30)(112,35)
      \psline(110,23)(114,23) \rput(117,23){$\dots$}
      \psline(120,23)(125,23)\rput(127,23){$\varepsilon$}
      \psline(127,25)(127,30)
      \psline(127,30)(132,35)
      \psline(127,30)(122,35)
      \psline(130,23)(135,23)\rput(137,23){$\varepsilon$}

    \end{pspicture}   
  \caption{$\Push{\sigma,1}$ operation with
    $\TOP{2}(\Push{\sigma,1}(s_1)) \leq \TOP{2}(\Pop{2}(s_1))$.}
  \label{fig:Push1}
\end{figure}
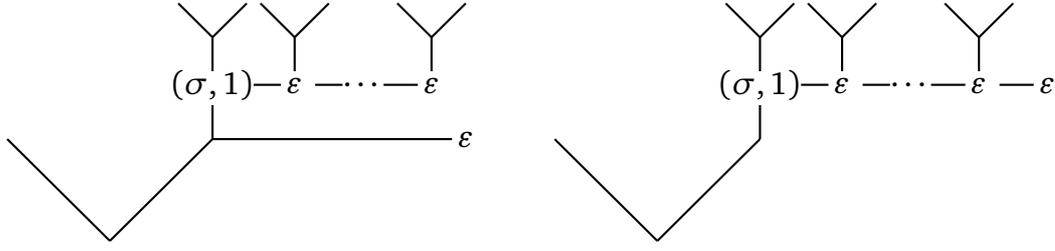

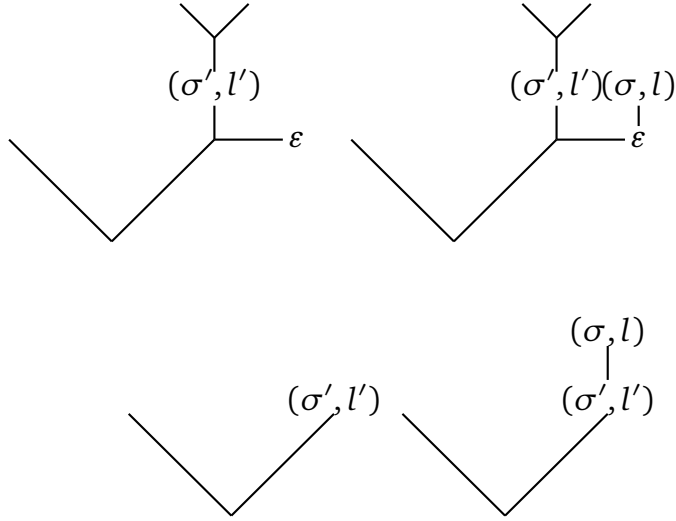
\begin{figure}
  \centering
    \psset{unit=.9mm} 
    \begin{pspicture}(-15,0)(90,40)
      \psline(-15,15)(0,0)
      \psline(15,15)(0,0)
      \psline(15,15)(15,20) \rput(15,23){$(\sigma',l')$}
      \psline(15,25)(15,30)
      \psline(15,30)(10,35)
      \psline(15,30)(20,35)
      \psline(15,15)(25,15) \rput(27,15){$\varepsilon$}

      \psline(35,15)(50,0)
      \psline(65,15)(50,0)
      \psline(65,15)(65,20) \rput(65,23){$(\sigma',l')$}
      \psline(65,25)(65,30)
      \psline(65,30)(60,35)
      \psline(65,30)(70,35)
      \psline(65,15)(75,15) \rput(77,15){$\varepsilon$}
      \psline(77,17)(77,20) \rput(77,23){$(\sigma,l)$}
    \end{pspicture}   
    \begin{pspicture}(90,0)(170,40)
      \psline(95,15)(110,0)
      \psline(125,15)(110,0)
      \rput(125,17){$(\sigma',l')$}

      \psline(135,15)(150,0)
      \psline(165,15)(150,0)
      \rput(165,17){$(\sigma',l')$}
      \psline(165,20)(165,25)
      \rput(165,27){$(\sigma,l)$}
    \end{pspicture}   
  \caption{The two versions of $\Push{\sigma,l}$ operation otherwise.}
  \label{fig:Push1otherwise}
\end{figure}

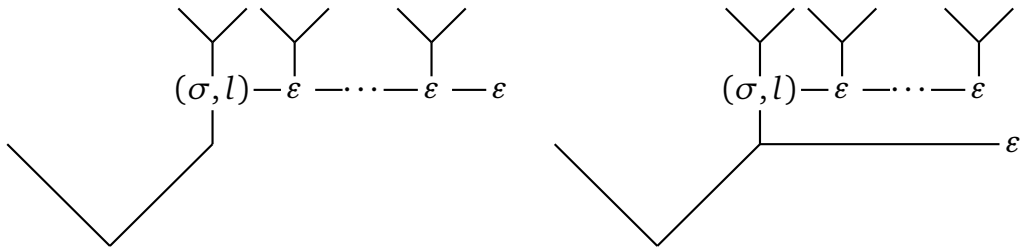
\begin{figure}
  \centering
    \psset{unit=.9mm} 
    \begin{pspicture}(-15,0)(135,40)
      \psline(-15,15)(0,0)
      \psline(15,15)(0,0)
      \psline(15,15)(15,20) \rput(15,23){$(\sigma,l)$}
      \psline(15,25)(15,30)
      \psline(15,30)(10,35)
      \psline(15,30)(20,35)
      \psline(21,23)(25,23) \rput(27,23){$\varepsilon$}
      \psline(27,25)(27,30)
      \psline(27,30)(22,35)
      \psline(27,30)(32,35)
      \psline(30,23)(34,23) \rput(37,23){$\dots$}
      \psline(40,23)(45,23)\rput(47,23){$\varepsilon$}
      \psline(47,25)(47,30)
      \psline(47,30)(52,35)
      \psline(47,30)(42,35)
      \psline(50,23)(55,23)\rput(57,23){$\varepsilon$}
      \psline(65,15)(80,0)
      \psline(95,15)(80,0)
      \psline(95,15)(95,20) \rput(95,23){$(\sigma,l)$}
      \psline(95,25)(95,30)
      \psline(95,30)(90,35)
      \psline(95,30)(100,35)
      \psline(101,23)(105,23) \rput(107,23){$\varepsilon$}
      \psline(107,25)(107,30)
      \psline(107,30)(102,35)
      \psline(107,30)(112,35)
      \psline(110,23)(114,23) \rput(117,23){$\dots$}
      \psline(120,23)(125,23)\rput(127,23){$\varepsilon$}
      \psline(127,25)(127,30)
      \psline(127,30)(132,35)
      \psline(127,30)(122,35)
      \psline(95,15)(130,15)\rput(132,15){$\varepsilon$}

    \end{pspicture}   
  \caption{$\Pop{1}$ operation on a cloned element.}
  \label{fig:Pop11}
\end{figure}

\begin{figure}
  \centering
    \psset{unit=.9mm} 
    \begin{pspicture}(-15,0)(85,40)
      \psline(-15,15)(0,0)
      \psline(15,15)(0,0)
      \psline(15,15)(15,20) \rput(15,23){$(\sigma,l)$}
      \psline(65,15)(80,0)
      \psline(95,15)(80,0)

    \end{pspicture}   
  \caption{$\Pop{1}$ operation otherwise.}
  \label{fig:Pop12}
\end{figure}
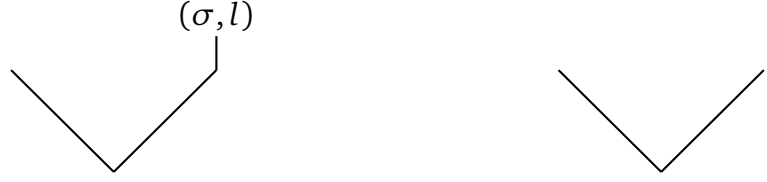

\begin{proof}
  Consider a transition $\delta:=(q, \sigma,\gamma, q', \op)$. 
  We show that there is an automaton that accepts
  $\Encode(c_1)\otimes\Encode(c_2)$ if and only if $\delta$ induces a
  transition from $c_1$ to $c_2$. 
  Thus, we have to define an automaton that accepts 
  $\Encode(c_1)\otimes\Encode(c_2)$ if and only if the following
  conditions are satisfied.
  \begin{enumerate}
  \item $c_1=(q,s_1)$ for some stack $s_1$,
  \item $c_2=(q',s_2)$ for some stack $s_2$,
  \item $\Sym(c_1)=\sigma$, and
  \item $\op(s_1) = s_2$. 
  \end{enumerate}
  The states of $c_1$ and $c_2$ may be checked directly at the
  root of \mbox{$\Encode(c_1) \otimes \Encode(c_2)$}. 
  $\Sym(c_1)$ is encoded in the last node of the rightmost path in
  $\Encode(c_1)$ that is not labelled $\varepsilon$. Hence, 
  the remaining problem is to construct an automaton for each stack
  operation $\op$ which recognises  $\Encode(s_1,
  (\bot,1))\otimes\Encode(s_2, (\bot, 1))$
  if and only if $s_2=\op(s_1)$.   

  We proceed by a case distinction on the stack operation. 
  \begin{itemize}
  \item   If $s_2 =\Push{\sigma,2}(s_1)$  or
    $s_2=\Clone{2}(s_1)$, then   
    $\Encode(s_1, (\bot, 1))$ and $\Encode(s_2, (\bot, 1))$ differ
    only in one node, which is 
    the rightmost leaf of $\Encode(s_2, (\bot, 1))$ (cf. Figures \ref{fig:Clone}
    and \ref{fig:Push1otherwise}). 
    This can easily be checked by an automaton.
  \item   If $s_2=\Push{\sigma,1}(s_1)$, we have to distinguish two
    cases.  In most
    cases, this operation behaves analogous to $\Push{\sigma,2}$ and
    $\Encode(s_2, (\bot, 1))$ is the extension of $\Encode(s_1, (\bot, 1))$ by  a left
    successor of the rightmost leaf of $\Encode(s_1, (\bot, 1))$. This new node
    is labelled $(\sigma,1)$. 

    But there is one case that is different, namely, when
    $s_1$ decomposes as 
    \begin{align*}
      s_1=s_1':
      \left(w \mathrel{\backslash}\left(\sigma w_1:\sigma w_2:
        \dots: \sigma w_n: \varepsilon\right)\right).      
    \end{align*}
    This case is depicted in
    Figure \ref{fig:Push1}. In this case, 
    \begin{align*}
      \TOP{1}(w) \sigma w_1: \TOP{1}(w) \sigma w_2:
      \dots: \TOP{1}(w)\sigma w_n      
    \end{align*}
    forms a block $b$ of the stack $s_1$. 
    $ \TOP{1}(w)\mathrel{\backslash}\varepsilon$ forms
    another block which is encoded in the rightmost leaf of
    $\Encode(s_1, (\bot, 1))$. Now, 
    \begin{align*}
      s_2=s_1':(  w 
      \mathrel{\backslash} \sigma w_1: \sigma w_2:  
      \dots: \sigma w_n: \sigma)
    \end{align*}
    i.e., in $\Encode(s_2, (\bot, 1))$ the
    whole block $\TOP{1}(w)
    \mathrel{\backslash} \sigma w_1: \sigma w_2:  
    \dots:  \sigma w_n: \sigma $ is encoded in a single
    subtree. This subtree extends the subtree encoding the block $b$
    by exactly one $\varepsilon$-labelled node as depicted in Figure
    \ref{fig:Push1}. 
    Thus, $s_2=\Push{\sigma,1}(s_1)$ if the following conditions are
    satisfied:
    \begin{enumerate}
    \item there is a node $d1\in\Encode(s_1, (\bot,
      1))\otimes\Encode(s_2, (\bot, 1))$ such 
      that $d1$ is the rightmost leaf of $\Encode(s_1, (\bot, 1))$,
    \item $d1\notin\Encode(s_2, (\bot, 1))$,
    \item $\Encode(s_2, (\bot, 1))$ extends $\Encode(s_1, (\bot, 1))$ by one
      node of the form $d01^m$, 
    \item  $d0$ is labelled by $(\sigma,1)$ in
      $\Encode(s_1, (\bot, 1))$ and $\Encode(s_2, (\bot, 1))$, and
    \item the two trees  
      coincide on all nodes but $d1$ and $d01^m$.
    \end{enumerate}
    These conditions are clearly \MSO-definable whence there is an
    automaton recognising these pairs of trees.
    
    Note that  the case distinction is also $\MSO$-definable. 
    For $d1$ the rightmost leaf of $\Encode(s_1, (\bot,1))$, the second case
    applies if and only if
    $d0$  has label $(\sigma,1)$ in $\Encode(s_1, (\bot,1))$.
    Again, the correspondence between \MSO and automata yields
    an automaton that accepts 
    $\Encode(s_1, (\bot,1))\otimes\Encode(s_2,(\bot,1))$ if and only if
    \mbox{$s_2=\Push{\sigma,1}(s_1)$}.  
  \item  Consider $s_2= \Pop{1}(s_1)$.
    Since $\Pop{1}$  is a kind of inverse of $\Push{\sigma,i}$, we 
    make a similar case distinction as in that case. 
    
    The different possibilities are depicted in the Figures
    \ref{fig:Pop11} and \ref{fig:Pop12}.
    Note the similarity of Figure \ref{fig:Pop11} and of Figure
    \ref{fig:Push1}, as well as the similarity of Figure
    \ref{fig:Pop12} and Figure \ref{fig:Push1otherwise}. 
    
    Both cases can be distinguished by an automaton. $\Encode(s_1,
    (\bot,1))$ and $\Encode(s_2, (\bot, 1))$ are as in Figure
    \ref{fig:Pop11} if and only if the rightmost leaf of $\Encode(s_1,
    (\bot,1))$ is a right successor. 

    Analogously to the push case, we conclude that there is an
    automaton that recognises $\Encode(s_1,
    (\bot,1))\otimes\Encode(s_2, (\bot,1))$ if and only if
    $s_2=\Pop{1}(s_1)$. 
  \item   For the case of $\Pop{2}$, recall Lemma \ref{LemmaStructurePop}
    and Figure \ref{fig:Pop2}. 
    An automaton recognising the $\Pop{2}$ operation only has to guess
    the set $D$ from Lemma \ref{LemmaStructurePop} and check whether the
    second tree is the restriction of
    the first tree to $D$. Note that the last element of $D$ along the
    rightmost path may be guessed nondeterministically and then the
    automaton may check that its guess was right. 
  \item   For the case of $\Collapse$, we have a case distinction due to
    the collapse level of the stack $s_1$.  Either $\Lvl(s_1)=1$ or
    $\Lvl(s_2)=2$. If it is $1$, the collapse operation on $s_1$ is
    equivalent to a $\Pop{1}$ operation.     
    Otherwise,  the collapse level of $s_1$ is $2$. 
    This case can be treated as in the 
    case of a $\Pop{2}$, but using Lemma \ref{LemmaStructureCol}
    instead of Lemma \ref{LemmaStructurePop}.

    Since the case distinction only
    depends on the collapse level stored in the label of the maximal node
    in the rightmost path of 
    $\Encode(s_1, (\bot,1))$
    which is not labelled $\varepsilon$, an automaton may 
    nondeterministically guess which case applies and verify its guess
    during the run on $\Encode(s_1, (\bot,1))\otimes\Encode(s_2,
    (\bot,1))$.\qedhere
  \end{itemize}
\end{proof}

We have seen that for each collapsible pushdown system the class of
encodings of valid 
configurations of this system is a set of regular
trees. Furthermore, 
all operations of a collapsible pushdown system are
automata-recognisable in this encoding. Putting these facts together
we obtain 
the following theorem.

\begin{theorem} \label{STACS:CPGTreeAutomatic}
 Given a collapsible pushdown system $\mathcal{S}$ of level $2$,
 one can effectively compute an automatic presentation of 
 the collapsible pushdown graph generated by
 $\mathcal{S}$. 
\end{theorem}

A direct corollary of this theorem is
the decidability of the first-order model checking on 
collapsible pushdown graphs (cf. Theorem
\ref{Thm:RamseyQuantifierDecidable}). 

\begin{corollary}
  The \mbox{$\FO{}(\exists^\infty,
    \exists^{\mathrm{mod}},(\RamQ{n})_{n\in\N})$-theory} of 
  every level $2$ collapsible pushdown graph is decidable.
\end{corollary}

\subsection{Tree-Automaticity of  Regular Reachability Predicates} 
\label{Reachability_Tree-Automatic}

In this section we show that 
regular reachability predicates are also automatic via 
$\Encode$. 
In the first part, 
we expand a collapsible pushdown graph $\CPG(\mathcal{S})$
 by the binary relation
$\Reach$ (cf. Definition \ref{DefReachPredicate}) and prove that this
predicate is 
automatic in our encoding. 
In the second part, we use the closure of collapsible pushdown systems
under products with automata in order to provide the
automaticity of all regular reachability predicates.

In order to show the regularity of the reachability predicate, 
we start with an observation about the general form of a run between
two configurations. 
Let $c_1=(q_1,s_1)$ and
\mbox{$c_2=(q_2,s_2)$}. For every run $\rho$ from $c_1$ to $c_2$ 
there are configurations $(q_3,s_3)$, $(q_4,s_4)$, $(q_5,s_5)$,
positions $i_3\leq i_4 \leq i_5\in\domain(\rho)$, and numbers
$m_2,m_3,m_5\in\N$ 
such that the following holds:
\begin{enumerate}
\item $\rho(i_3)=(q_3,s_3)$ and $s_3=\Pop{2}^{m_2}(s_1)$, 
\item $\rho(i_4)=(q_4,s_4)$, $s_4=\Pop{1}^{m_3}(s_3)$ and $s_4$ is a
  common substack of $s_1$ and $s_2$, 
\item $\rho(i_5)=(q_5,s_5)$,  
  $s_4=\Pop{1}^{n_5}(s_5)$ and $s_5=\Pop{2}^{m_5}(s_2)$, and
\item $\rho$ does not visit any proper substack of $s_4$. 
\end{enumerate}

For any run $\rho$, $s_4$ is found as follows: it is the minimal
substack of $s_1$ that is visited 
by $\rho$. $i_4$ is then an arbitrary position in $\rho$ that visits $s_4$. 
The existence of $i_5$ follows directly from the fact that
$s_5$ is a milestone of $s_2$ and the fact that $\rho$ visits $s_4$,
which is a substack of $s_5$. The existence of $i_3$ is clear from
the fact that the run $\rho$ has to reach a stack of width 
$\lvert s_4 \rvert$ at first, before it can change the $\lvert s_4
\rvert$-th word of the stack, i.e., before it can reach $s_3$. 
 
We use this decomposition for proving the regularity of $\Reach$ as
follows. 

\begin{definition} \label{ABCD-Definition}
  Given a collapsible pushdown system $\mathcal{S}$, we define the
  following four relations on the configurations of $\mathcal{S}$:
  \begin{enumerate}
  \item   Let $A\subseteq \CPG(\mathcal{S})\times\CPG(\mathcal{S})$ be the
    relation containing those pairs of configurations $(c_1,c_2)$ with
    $c_1=(q_1, s_1)$ and $c_2=(q_2, s_2)$ such that
    \begin{enumerate}
    \item  $s_2=\Pop{2}^m(s_1)$,
    \item there is a run  $\rho$ from 
      $c_1$  to $c_2$ and
    \item  $\rho$  does not
    visit a proper substack of $s_2$. 
    \end{enumerate}
  \item Let $B\subseteq \CPG(\mathcal{S})\times\CPG(\mathcal{S})$ be
    the relation containing those pairs of configurations $(c_1, c_2)$
    with $c_1=(q_1, s_1$) and $c_2=(q_2, s_2)$ such that
    \begin{enumerate}
    \item  $s_2=\Pop{1}^m(s_1)$,
    \item  there is a run  $\rho$ from  $c_1$ to  $c_2$  and
    \item  $\rho$ does not
      visit a proper substack of $s_2$. 
    \end{enumerate}
  \item Let $C\subseteq \CPG(\mathcal{S})\times\CPG(\mathcal{S})$ be
    the relation containing those pairs of configurations $(c_1, c_2)$
    with $c_1=(q_1, s_1$) and $c_2=(q_2, s_2)$ such that
    \begin{enumerate}
    \item $s_1=\Pop{1}^m(s_2)$,
    \item there is a run  $\rho$ from $c_1$ to $c_2$ and
    \item $\rho$  does not visit  a proper substack of $s_1$. 
    \end{enumerate}
  \item   Let $D\subseteq \CPG(\mathcal{S})\times\CPG(\mathcal{S})$ be the
    relation containing those pairs of configurations $(c_1,c_2)$ with
    $c_1=(q_1, s_1)$ and $c_2=(q_2, s_2)$ such that
    \begin{enumerate}
    \item 
      $s_1=\Pop{2}^m(s_2)$ ,
    \item there is a run  $\rho$ from $c_1$ to $c_2$ and
    \item $\rho$ does not visit a substack of $s_1$ after its initial
      configuration.  
    \end{enumerate}
  \end{enumerate}
\end{definition}
\begin{remark}
  Since we allow runs of length $0$, the relations $A$, $B$, $C$ and
  $D$ are reflexive, i.e., for all configurations $c$, $(c,c)\in A$,
  $(c,c)\in B$, $(c,c)\in C$ and $(c,c)\in D$. 
\end{remark}
The relation $\Reach$ can be expressed via $A,
B, C$ and $D$ in the sense that for arbitrary configurations $c_1,
c_2$, $(c_1, c_2)\in\Reach$ holds if and only if there are configurations
$x, y, z$ such that $(c_1,x)\in A$, $(x,y)\in B$, $(y,z)\in C$ and
$(z,c_2)\in D$. Since projections of regular sets are regular,
$\Reach$ is an automatic relation via the encoding $\Encode$ if
the relations $A, B, C$ and $D$ are automatic via $\Encode$. 
Proving the regularity of these relations is our next goal. We first
prove the  
regularity of $A$. This proof requires an analysis of runs from some
stack $s$ to some stack $\Pop{2}^n(s)$ for every $n\in\N$. We obtain a
characterisation of 
these runs that can be checked by an automaton.

\paragraph{Regularity of the Relation $\mathbf{A}$}

At a first glance one
might think that  a run from some stack $s$ to a stack
$\Pop{2}^n(s)$ 
only consists of a sequence of returns. But this is only true if
we do not use the collapse operation. 
A collapsible pushdown system may start by writing a lot of
information with $\Clone{2}$ 
and $\Push{\sigma,l}$ operations onto the stack, then use a couple of
$\Pop{1}$ operations to come to an element with a small collapse link
and finally use the collapse to jump to a very small substack of $s$
without using any other substack of $s$ in between. Such a run does
not contain any returns at all. 

In order to cope with such runs, 
we introduce the notion of a \emph{level-$1$-loop}. A level-$1$-loop
is a 
kind of loop of the topmost word which increases the number of
words on the level $2$ stack. We prove that the pairs of
initial and final states of these new loops
are computable in a similar way
as for ordinary loops. 
Furthermore, we show that 
every run from $s$ to $\Pop{2}^n(s)$ decomposes mainly into parts that are
basically returns, loops or $1$-loops. These parts are connected by
application of either a $\Pop{1}$  or a $\Collapse$ operation.
First, we introduce $1$-loops. Then we show the decomposition
result we mentioned above. Finally, we use this decomposition for
showing the regularity of 
the relation $A$. For this purpose, we introduce
\emph{certificates for substack reachability}. 
We consider a certificate as the  abstract representation of 
the  decomposition of some (potentially existing) run. The certificate
consists of the final state of each part of the
decomposition of this run. Using these certificates, we
reduce the problem whether a run exists to the problem whether the subruns
that form the parts of the decomposition exist. 
This is a much simpler problem because each of these subruns can only
have a very special form. 
Finally, we show that an automaton can check the existence of these
subruns while processing the certificate and the trees encoding the
initial and final configuration of the run.

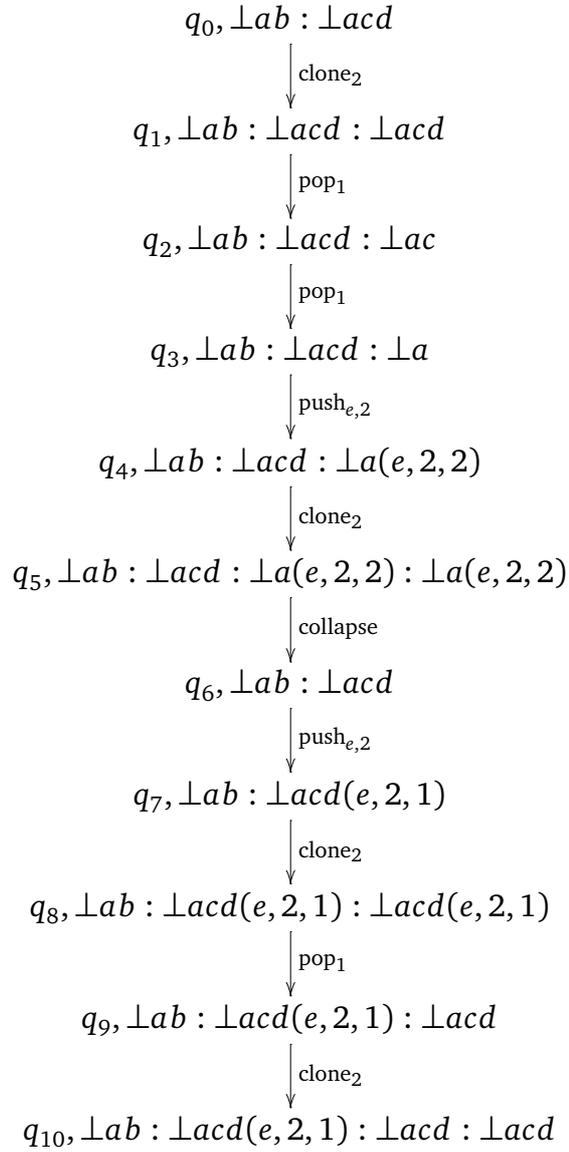
\begin{figure}
  \hskip 5cm
  \begin{xy}
    \xymatrix{ q_0, \bot a b : \bot a c d \ar[d]^{\Clone{2}}\\
      q_1, \bot a b : \bot acd : \bot acd \ar[d]^{\Pop{1}}\\
      q_2, \bot ab : \bot acd: \bot ac \ar[d]^{\Pop{1}}\\
      q_3, \bot ab: \bot acd: \bot a \ar[d]^{\Push{e,2}}\\
      q_4, \bot ab: \bot acd : \bot a(e,2,2) \ar[d]^{\Clone{2}}\\
      q_5, \bot ab: \bot acd: \bot a(e,2,2): \bot a(e,2,2)
      \ar[d]^{\Collapse}\\ 
      q_6, \bot ab: \bot acd \ar[d]^{\Push{e,2}}\\
      q_7, \bot ab: \bot acd(e,2,1) \ar[d]^{\Clone{2}}\\
      q_8, \bot ab: \bot acd(e,2,1) : \bot acd(e,2,1) \ar[d]^{\Pop{1}}\\
      q_9, \bot ab: \bot acd(e,2,1) : \bot acd \ar[d]^{\Clone{2}}\\
      q_{10}, \bot ab: \bot acd(e,2,1) : \bot acd: \bot acd
    }
  \end{xy}
  \caption{Example of a $1$-loop of $s:=\bot ab: \bot acd$. The part
    between $q_2$ and $q_6$ forms a return of a stack with topmost
    word $\TOP{2}(\Pop{1}(s))$. Note that the run up to $q_9$ 
    also forms a $1$-loop.}
  \label{fig:OneLoop}
\end{figure}

\begin{definition}
  Let $s$ be some stack and $w$ some word. 
  A run $\lambda$ of length $n$ is called a \emph{level-$1$-loop} (or
  $1$-loop) of $s:w$ 
  if the following conditions are satisfied.
  \begin{enumerate}
  \item $\lambda(0)=(q_0, s:w)$ for some $q_0\in Q$,
  \item $\lambda(n)=(q_n, s:s':w)$ for some nonempty stack $s'$ and
    some state $q_n\in Q$, \label{OneLoopCondTwo}
  \item for every $i\in\domain(\lambda)$, $\lvert \lambda(i)\rvert >
    \lvert s \rvert$,  \label{OneLoopCondThree} and
  \item for every $i\in\domain(\lambda)$ such that $w\leq
    \TOP{2}(\lambda(i-1))$ and 
    $\TOP{2}(\lambda(i))=\Pop{1}(w)$, there is some $j>i$ such that
    $\lambda{\restriction}_{[i,j]}$ is a return. 
  \end{enumerate}
\end{definition}
\begin{remark}
   Under condition \ref{OneLoopCondTwo}, condition
   \ref{OneLoopCondThree} is equivalent to the condition that
   $\lambda$  never passes the stack  $s$.
    An example of a $1$-loop can be found in Figure \ref{fig:OneLoop}.   
    Note that the last two conditions imply that a $1$-loop does never
    visit a proper substack of $s:w$. 
\end{remark}

\begin{definition}
  For a fixed collapsible pushdown system  $\mathcal{S}$ and some
  stack $s$ we denote by 
  $\exOneLoops_{\mathcal{S}}(s)$ the set 
  \begin{align*}
  \left\{(q_1,q_2)\in Q\times Q:\text{ there is an
  }s'\in\Stacks({\Sigma}) \text{ and a $1$-loop of $\mathcal{S}$ from
  }(q_1,s)\text{ to 
  }(q_2,s')\right\}.       
  \end{align*}
  If $\mathcal{S}$ is clear
  from the context, we omit it. 
\end{definition}

We use this rather technical definition of a $1$-loop due to two 
important properties. Firstly, we obtain a similar computational
behaviour of 
$1$-loops  as for loops and returns: 
$\exOneLoops(s)$
only depends on the returns of $\Pop{1}(s)$,
 $\Lvl(s)$ and $\Sym(s)$. 
Secondly, this notion is strong enough to capture all parts of a run
from a stack $s$ to $\Pop{2}^n(s)$ that are not captured by the
notions of loops and returns. This idea is made precise in Lemma
\ref{FormLemma}.

\begin{lemma}\label{DetermineOneLoops}
  There is an algorithm that determines for every stack $s$ the set
  $\exOneLoops(s)$ from the input $\Sym(s)$, $\Lvl(s)$ and
  $\exReturns(\TOP{2}(\Pop{1}(s)))$.
\end{lemma}
\begin{proof}[Proof (sketch).]
  First of all note the similarity of the claim to the corresponding
  lemmas concerning returns, loops, low loops and high loops. 
  The main ingredients of this proof are variants of Lemma
  \ref{ReachableConfigsinMuCalcul} and Lemma 
  \ref{Cor:InductiveComputability}.
  \begin{itemize}
  \item Analogously to Lemma \ref{ReachableConfigsinMuCalcul}, it is
    decidable whether there is some reachable 
    configuration of the form $c=(q,s')$ with $\lvert s' \rvert \geq 3$
    and $\lvert \TOP{2}(s') \rvert = 3$. Note that by definition of the
    return-simulator $\ReturnSimulator{\mathcal{S}}{s}{k}$, $\lvert
    \TOP{2}(s') \rvert = 3$ is equivalent to 
    \mbox{$\TOP{2}(s')= \bot\top
    (\Sym(s), \Lvl(s), \kappa(\Lvl(s)))$} for all configurations of
    some return simulator for all stacks $s$ and $s'$. 
    The decidability follows by reduction to $L\mu$ model checking. We
    equip the pushdown system with a testing device. This testing
    devices first tries to perform two $\Pop{1}$ operations. If it
    then reaches the bottom of stack, it tries to perform two
    $\Pop{2}$-operations. If this is possible, then the stack is of
    the desired form.
  \item Analogously to the return case 
    \ref{Cor:InductiveComputability}, one proves that 
    $(q_1, q_2)\in\exOneLoops(s)$ if and only if the graph of the
    return simulator 
    $\ReturnSimulator{\mathcal{S}}{s}{k}$ contains some run from 
    $(q_1, \bot\top \TOP{1}(s)\Box: \bot\top\TOP{1}(s))$ to 
    $(q_2, s')$ where $\lvert s' \rvert \geq 3$ and
    $\TOP{2}(s')=\bot\top \TOP{1}(s)$.
    Analogously to the return case, 
    each such run corresponds to a $1$-loop starting in $(q_1,s)$ and
    ending in state $q_2$. 
    Again, we copy the transitions of such a simulation one to one to
    a run starting in $(q_1,s)$. Whenever we come to a transition on
    topmost symbol $\top$, we replace the following
    $\Pop{2}$-transition by a return of some stack with topmost word
    $\Pop{1}(\TOP{2}(s))$. 
  \end{itemize}
  Putting these two facts together, we obtain that
  $\exOneLoops(s)$ can be computed from
  $\ReturnSimulator{\mathcal{S}}{s}{k}$, $\Sym(s)$ 
  and $\Lvl(s)$. But the definition of the return simulator only
  depends on $\exReturns(\Pop{1}(s))$. This concludes the proof. 
\end{proof}

Of course, we can turn the previous proof into a definition of an
automaton calculating the $1$-loops of all milestones of a stack. This
is completely analogous to Propositions
\ref{Prop:AutomatonForReturns} and \ref{Prop:AutomatonForLoops}.

\begin{corollary}\label{OneLoopsAutomaton}
  For each collapsible pushdown system $\mathcal{S}$ of level $2$,
  we can compute an automaton $\mathcal{A}$ that calculates for
  each configuration $c$ at
  each $d\in\Encode(c)$ the set
  $\exOneLoops_{\mathcal{S}}(\LeftStack(d,\Encode(c)))$. 
\end{corollary}

Now, we analyse the form of any run from some stack $s$ to
some substack $s'=\Pop{2}^n(s)$. 

\begin{lemma} \label{FormLemma}
  Let  $s$ and $s'$ be stacks such that
  $s'=\Pop{2}^m(s)$ for some $m\in\N$. 
  Let $\rho$ be a run from $s$ to $s'$
  such that $\rho$ does not visit a proper substack of $s'$. 
  Then $\rho$ decomposes as $\rho_1\circ \rho_2 \circ \dots \circ
  \rho_n \circ \lambda$ where
  $\lambda$ is a high loop of $s'$ and each $\rho_i$ is of one of
  the following forms.
  \begin{enumerate}[F1.]
  \item \label{FormReturn}$\rho_i$ is a return, 
  \item \label{1LoopCol} $\rho_i$ is a
    $1$-loop followed 
    by a $\Collapse$ of collapse level $2$, 
  \item \label{LoopCol} 
    $\rho_i$ is a loop followed by a 
    $\Collapse$ of collapse level $2$,
  \item \label{LoopPop} $\rho_i$ is a loop followed by
    a $\Pop{1}$ (or a 
    $\Collapse$ operation of collapse level $1$) and 
    there is a $j>i$ such that 
    $\rho_j$ is of the form F\ref{1LoopCol} or 
    F\ref{LoopCol} and there
    is no $i<k<j$ such that $\rho_k$ is
    of the form F\ref{FormReturn}, 
  \item \label{1LoopPop} $\rho_i$ is a 
    $1$-loop followed by a 
    $\Pop{1}$ operation (or a 
    $\Collapse$ of collapse level $1$) and 
    there is a $j>i$ such that 
    $\rho_j$ is of the form F\ref{1LoopCol} or 
    F\ref{LoopCol} and there
    is no $i<k<j$ such that $\rho_k$ is
    of the form F\ref{FormReturn}. 
  \end{enumerate}
\end{lemma}

\begin{proof}
  Let $s' = \Pop{2}^n(s)$ for $n\geq 0$ and $\rho$ a run from $s$ to
  $s'$ not passing any proper substack of $s'$. 
  
  First of all, note that the case $n=0$ is trivial. If $n=0$, $\rho$
  is by definition a high loop of $s$. 
  
  For the case $n>0$, 
  we proceed by induction on the length of $\rho$. We write
  $(q_i,s_i)$ for the configuration $\rho(i)$.  
  Firstly, consider the case where there is some $m\in\domain(\rho)$
  such that  $\rho_1:=\rho{\restriction}_{[0,m]}$ is a return.
  Then $\rho_1$ is of the form
  F\ref{FormReturn}. By induction hypothesis,
  $\rho{\restriction}{[m,\length(\rho)]}$ decomposes as desired. 

  Otherwise, 
  assume that there is no $m\in\domain(\rho)$ such that
  $\rho{\restriction}_{[0,m]}$ is a return. 
  
  Nevertheless, there is a minimal $m\in\domain(\rho)$ such that for
  all $i<m$, it holds that
  $\lvert s_i \rvert \geq \lvert s \rvert$ and $\lvert s_m\rvert <
  \lvert s \rvert$.  
  The last operation of $\hat\rho:=\rho{\restriction}_{[0,m]}$ is a
  $\Collapse$ such that \mbox{$\TOP{2}(s_{m-1}) \leq \TOP{2}(s)$} (otherwise
  $\hat\rho$ would be a return). 
  
  Writing $w:=\TOP{2}(s_{m-1})$, we distinguish two cases.
  \begin{enumerate}
  \item 
  First consider the case that $w=\TOP{2}(s)$. 
  Note that this  implies $\Lvl(s)=2$ because the  last operation of
  $\hat\rho$ is a collapse of level $2$.
  
  Furthermore, we claim that  $\hat\rho$ does not
  visit $\Pop{1}(s)$. Heading for a contradiction, assume that 
  $\hat\rho(i) = \Pop{1}(s)$ for some $i\in\domain(\hat\rho)$.
  Since $\hat\rho$ does not visit $\Pop{2}(s)$ between $i$ and $m-1$,
  $\TOP{2}(\hat\rho(m-1))=w$ is only possible if $\Lnk(w) = \lvert s
  \rvert -1$ (smaller links cannot be restored by $\hat\rho$). 
  But then $\hat\rho{\restriction}_{[i,m]}$ is a return of
  $\Pop{1}(s)$ whence $\hat\rho$ is a return of $s$. This contradicts
  the assumption that $\hat\rho$ is no return. 

  Hence, $\hat\rho$ does not pass $\Pop{1}(s)$ and we distinguish the
  following cases
  \begin{itemize}
  \item Assume that the stack of $\hat\rho(m-1)$ is $s$. Then
    $\hat\rho$ is a high loop followed by a collapse: 
    the stack at $\hat\rho(0)$ and $\hat\rho(m-1)$ is $s$ and the run does
    not visit $\Pop{2}(s)$ or $\Pop{1}(s)$ in between whence its
    restriction to $[0,m-1]$ is a
    high loop. Thus, 
    $\rho_1:=\hat\rho$ is of the form F\ref{LoopCol} and the claim
    follows by induction hypothesis.
  \item Assume that the stack of $\hat\rho(m-1)$ is
    $s' = s:t :w$ for some nonempty
    $2$-word $t$. We claim that $\hat\rho$  
    is a $1$-loop plus a
    $\Collapse$ operation: 
    We have already seen that $\hat\rho$ does not visit any proper
    substack of $s$. 
    Thus, it suffices to show that $\hat\rho$ reaches a stack with
    topmost word $\Pop{1}(w)$ only at positions where a return
    starts. 

    Let $i$ be some position such that $w\leq \hat\rho(i-1)$ and
    $\TOP{2}(\hat\rho(i))=\Pop{1}(w)$. 
    Recall that $\TOP{2}(s_{m-1})=w$, $\Lnk(w)=2$ and 
    $\Lvl(w)\leq\lvert s \rvert -1$. Since 
    $\lvert \hat\rho(i) \rvert > \lvert s \rvert$, 
    we cannot restore $\TOP{1}(w)$ by a push operation. Thus, there is
    some minimal position $j>i$ such that 
    $\lvert \hat\rho(j) \rvert < \lvert \hat\rho(i) \rvert$. Since the
    level $2$ links of $w$ point 
    below $\Pop{2}(s)$ and no proper substack of $s$ is reached by
    $\hat\rho{\restriction}_{[0,m-1]}$, the links stored in $w$ are
    not used in  
    $\hat\rho{\restriction}_{[i,j]}$. It follows immediately that
    $\hat\rho{\restriction}_{[i,j]}$ is a return. 

    Thus, $\rho_1:=\hat\rho$ is of the form F\ref{1LoopCol}.
  \end{itemize}
\item 
  For the other case, assume that $w<\TOP{2}(s)$. 
  Then there is a
  minimal \mbox{$i\in\domain(\hat\rho)$} such that
  $\TOP{2}(\hat\rho(i))=\Pop{1}(w)$ and 
  there is no $j>i$ such that $\hat\rho{\restriction}_{[i,j]}$
  is a return. 

  We claim the following: if $\hat\rho(i) = \Pop{1}(s)$, then $\hat\rho_1$
  is of the form F\ref{LoopPop}, otherwise $\hat\rho$ is of the form F\ref{1LoopPop}.
  Due to the definition, 
  $\rho_1:=\hat\rho{\restriction}_{[0,i]}$ is a loop or $1$-loop
  followed by a $\Pop{1}$ or a collapse. 
  Hence, it suffices to check the side conditions on the segments
  following in the decomposition of $\rho$. 
  For this purpose set $\rho':=\rho{\restriction}_{[i,\length(\rho)]}$. 
  By induction hypothesis $\rho'$ decomposes as 
  $\rho'=\rho_2\circ\rho_3\circ\dots\circ\rho_n\circ\lambda$ where
  the $\rho_i$ and $\lambda$ satisfy the claim of the lemma. 

  Now, by definition of $i$, $\rho'$ does not start with a
  return. Thus, $\rho_2$ is of one of the forms
  F\ref{1LoopCol}--F\ref{1LoopPop}. But all these 
  forms require that there is some $j\geq 2$ such that $\rho_j$ is of
  form F\ref{1LoopCol} or F\ref{LoopCol} and for all $2\leq k < j$,
  $\rho_k$ is not of the form 
  F\ref{FormReturn}. 
  
  From this condition, it follows directly that $\rho=\rho_1\circ\rho'
  = \rho_1\circ\rho_2\circ\rho_3\circ\dots\circ\rho_n\circ\lambda$
  and $\rho_1$ is of the form F\ref{LoopPop} or F\ref{1LoopPop}.\qedhere
  \end{enumerate}
\end{proof}

The following example illustrates the lemma. 
\begin{example} \label{Example:ReachDecomposition}
  Consider the stack $s$ in Figure \ref{fig:Stack_s} and the 
  following transitions:
  \begin{enumerate}
  \item $(q_1,a,q_2,\Clone{2})$,
  \item $(q_2,a,q_3,\Collapse)$,
  \item $(q_3,a,q_2,\Clone{2})$,
  \item $(q_3,b,q_2,\Push{b,2})$,
  \item $(q_1,c,q_4,\Push{b,1})$,
  \item $(q_2,b,q_1,\Collapse)$,
  \item $(q_2,c,q_1,\Collapse)$,
  \item $(q_4,b,q_2,\Pop{1})$.
  \end{enumerate}
  Since these transitions form a deterministic relation, there is a
  unique run 
  starting in $(q_1,s)$. This run $\rho$ is generated by
  using the transitions in the following order: $(1)$, $(2)$, $(3)$,
  $(2)$, $(4)$, $(6)$, $(5)$, $(8)$, $(7)$, $(5)$, $(8)$, $(7)$. 
  The run $\rho$ ends in the configuration 
  \mbox{$(q_1,\bot_2)=(q_1,\Pop{2}^4(s))$}. According to the
  decomposition of 
  Lemma 
  \ref{FormLemma}, 
  $\rho{\restriction}_{[0,2]}$ is of the form F\ref{1LoopPop}, 
  $\rho{\restriction}_{[2,4]}$ is of the form F\ref{1LoopCol}, 
  $\rho{\restriction}_{[4,6]}$ is of the form F\ref{FormReturn}, 
  $\rho{\restriction}_{[6,9]}$ is of the form F\ref{LoopPop}, and
  $\rho{\restriction}_{[9,12]}$ is of the form F\ref{LoopCol}.
\end{example}
\begin{figure}
  \centering
  $
   \begin{matrix}
             &     &           &          &  a & a \\
             &     &  c   &   b &  (a,2,3) & (a,2,3) \\
             &     & (c,2,1)   &  (c,2,1) &  (c,2,1) & (c,2,1) \\
  s\coloneqq &\bot & \bot     &\bot      &\bot     &\bot 
  \end{matrix}
  $ 
  \caption{The stack $s$ of example \ref{Example:ReachDecomposition}.
  }
  \label{fig:Stack_s}
\end{figure}

The previous lemma tells us that any run to a substack decomposes into
subruns of three forms: 
\begin{enumerate}
\item returns,
\item  subruns that decrease the length of the topmost word by one, or
\item subruns that end in a collapse of level $2$ applied to some stack
  with the same topmost word as their initial stack. 
\end{enumerate}
If subruns of the second case occur, then they are followed by a
subrun of the third form before any return occurs. 
Since runs of the third form end in a $\Collapse$ of level $2$, it
does not matter whether runs of the second or third form have
increased the width of the stack in between:
eventually we perform a collapse operation on a prefix of
the initial topmost word. This collapse then deletes all the new words
that were created in between.

The decomposition of a run according to Lemma \ref{FormLemma} is the
starting point for deciding whether there is a run from some
configuration $(q,s)$ to some
$(q',\Pop{2}^n(s))$. The basic idea is that we guess the form and the
final state of each segment the run consists of. 
We then attach this guess to the encoding of the two configurations.
We will call such a guess \emph{ certificate for substack
  reachability}. 
Finally, we prove that there is an automaton that can check
whether a certificate for substack reachability actually encodes some
run from $(q,s)$ to $(q', \Pop{2}^n(s))$. 

This approach is quite similar to the proof that the reachable
configurations of a given collapsible pushdown system form a regular
set. 
Let us first recall the basic idea of that proof.
We used each node of $d\in \Encode(q,s)$ as representative for the
milestone $\LeftStack(d, \Encode(q,s))$ and a certificate for
reachability labelled every node with the state in which some run
visited the 
corresponding milestone. 

Now, we do a similar thing. Given a run $\rho$ from $(q_1, s)$
to $(q_2, \Pop{2}^m(s))$, let 
\mbox{$\rho=\rho_1 \circ \rho_2 \circ \dots 
\circ \rho_n \circ\lambda$} be its decomposition according to Lemma
\ref{FormLemma}. We want to find a representative for the initial
configuration of 
each of the $\rho_j$ and label this representative with a description
of
$\rho_j$. In fact, we label the representative with the
final state of 
$\rho_j$ and the type of $\rho_j$ according to the classification from
Lemma \ref{FormLemma}. 

Let us first explain the system of representation. Let
$d\in\Encode(q,s)$ be some node. We write
$s_d:=\LeftStack(d,\Encode(q,s))$ for the milestone induced by $d$. 
Now, we will use $d$ as a representative for any stack $\hat s_d$ that has
the following two properties:
\begin{enumerate}
\item $\Pop{2}(s_d) = \Pop{2}^k(\hat s_d)$ for some $k\in \N$ and
\item $\TOP{2}(s_d) = \TOP{2}(\hat s_d)$. 
\end{enumerate}
This implies that $d$ may represent $s_d$ or
some stack $\Pop{2}(s_d): s': \TOP{2}(s_d)$ for $s'$ an arbitrary
$\TOP{2}(s_d)$ prefixed stack. 

Let us explain why this form of representation is sufficient for our
purpose. Recall 
that the existence of $1$-loops and loops only depends on the topmost
word of a 
stack. Thus, we only need to know the topmost word of some stack in
order to verify the existence of $1$-loops or loops for certain pairs of
initial and final states. Furthermore, if we know the topmost word of
some stack, we can easily derive the topmost word of the stack reached
via  $\Pop{1}$ or $\Collapse$ of level $1$. 
Moreover, if $d$ is a representative for some stack $\hat s_d$, then a
collapse of level $2$ from $\hat s_d$ and from $s_d$ result in the same
stack: if $d$ represents $\hat s_d$ then $\TOP{2}(s_d)$ and
$\TOP{2}(\hat s_d)$
coincide. Thus, level $2$ collapse links in the topmost word of $\hat s_d$
point to some substack of $\Pop{2}(s_d)$ (because the links of $s_d$
have this property by definition of a stack). Since 
$\Pop{2}(s_d) = \Pop{2}^k(\hat s_d)$,
the collapse link of
$s_d$ and of $\hat s_d$ point to the same substack of $\Pop{2}(s_d)$. 

Hence, the representatives that we use are sufficiently similar to the
represented stacks in the following sense.
The existence of subruns of the forms
F\ref{1LoopCol}--F\ref{1LoopPop}  can 
be decided by considering the representatives. 
Note that subruns of the form F\ref{FormReturn}, which are  returns,
occur as an 
initial part of the run or
after the application of some
$\Collapse$ of level $2$. At such positions, the
corresponding node $d$ represents a stack $\hat s_d$ such that
$\hat s_d=s_d$. 
Thus, 
$\Pop{2}(\hat s_d)=\Pop{2}(s_d)$ is determined by $d$. 
Hence, we can find  a node $d'$ such that \mbox{$s_{d'} = \hat s_{d'} =
\Pop{2}(\hat s_d)$}.  

Having explained the system of representation, let us introduce 
\emph{certificates for substack reachability}. Before we come to the
formal definition, we explain the underlying idea. 

Given a tree $\Encode(q_1,s_1)\otimes\Encode(q_2,s_2)$ such that
$s_2=\Pop{2}^n(s_1)$, we want to label this tree with information
witnessing the existence of a run from $(q_1,s_1)$ to $(q_2, s_2)$. 
Assume that there is such a run $\rho$. Let $\rho=\rho_1 \circ \rho_2
\circ \dots \rho_n \circ \lambda$ be its decomposition into parts
according to 
Lemma \ref{FormLemma}. Recall that the rightmost leaf $d$ of
$\Encode(q_1,s_1)$ represents the stack $s_1$. Since $\rho_1$ starts with stack
$s_1$, this $d$ is the position in $\Encode(q_1,s_1)\otimes\Encode(q_2, s_2)$
which we want to label with information concerning $\rho_1$. We will
label this node with the final state of $\rho_1$ and the type of
this run according to the classification from Lemma \ref{FormLemma}. If
$\rho_1$ is of the form F\ref{FormReturn} (i.e., $\rho_1$ is a return), we
label it by \ref{FormReturn}, if it is of the form F\ref{1LoopCol}, we
label it by \ref{1LoopCol}, etc.  

Now, assume that there is some node $d$ that represents some stack
$s'$ such that $s'$ is the initial stack of $\rho_j$ for some $1\leq j
\leq n$. The type of $\rho_j$ defines a representative for the initial
stack of $\rho_{j+1}$, which is the final stack of
$\rho_j$, as follows.
\begin{enumerate}
\item  If $\rho_j$ is a return, the initial 
  stack of $\rho_{j+1}$ is $\Pop{2}(s_d)$. There is a node $d'$ such
  that $s_{d'}=\Pop{2}(s_d)$. This node is the representative of
  $\rho_{j+1}$.
\item If $\rho_j$ ends in a
  collapse of level $2$ (from a stack with topmost word $\TOP{2}(s_d)$)
  the initial stack of $\rho_{j+1}$ is $\Collapse(s_d)$. There is a node
  $d'$ such that 
  $s_{d'}=\Collapse(s_d)$.
\item Finally, if 
  $\rho_j$ ends in a $\Pop{1}$ or a $\Collapse$ of level $1$, we need to
  find a representative  $d'$ such that 
  $\TOP{2}(s_{d'}) = \TOP{2}(\Pop{1}(s_d))$. We take the
  lexicographically maximal node $d'$ such that $\LeftStack(d',
  \Encode(q_1,s_1))$ is a milestone of $s_d$ with topmost word
  $\TOP{2}(\Pop{1}(s_d))$. 
\end{enumerate}

We call the representative $d'$ of the initial stack of $\rho_{j+1}$
the \emph{successor of $d$}. Keep in mind that this successor depends on
the label of $d'$. Furthermore, note that the successor is \MSO-definable
on $\Encode(q_1, s_1)$ if the label of $d$ is known: the $\Pop{2}$ or
$\Collapse$ successor of $s_d$ is clearly definable due to the
regularity of the operations $\Pop{2}$ and $\Collapse$. For the third
case, note that the successor of $d$ is the unique ancestor of $d$
such that $d=d'01^m$ for some $m\in\N$. 

Since we have found a representative for $\rho_{j+1}$, we label
it again by the final state of $\rho_{j+1}$ and by the type of
$\rho_{j+1}$. We continue this process  until we have defined a
representative for each segment of the run $\rho$. 
We will soon see
that an automaton can check whether an arbitrary labelling of the
nodes of $\Encode(q_1, s_1) \otimes \Encode(q_2, s_2)$ is indeed a
labelling corresponding to an existing run from $(q_1, s_1)$ to
$(q_2,s_2)$ in this sense.

Let us now formally introduce \emph{certificates for substack
reachability}. We will call such a certificate valid if it witnesses
the existence of a run from the larger configuration to the smaller
one.

\begin{definition}
  Let $c_1=(q_1,s_1), c_2=(q_2,s_2)$ be configurations such that 
  $s_2=\Pop{2}^n(s_1)$ for some $n\in \N$. We call a function
  \begin{align*}
    \ReachguessSym:\domain(\Encode(c_1))\setminus\domain(\Encode(c_2))
    &\rightarrow \{1, 2, 3, 4, 5\} \times Q
  \end{align*} 
  a \emph{certificate for substack reachability} for $c_1$ and $c_2$. 
\end{definition}
\begin{remark} \label{CertificatesofSubstackReachabilityMSODefinable}
  Due to the finite range of a certificate for substack reachability,
  we can express quantification over
  certificates for substack reachability 
  for $c_1$ and $c_2$   on the structure
  $\Encode(c_1)\otimes\Encode(c_2)$ in \MSO. 
  
  Even though these certificates are defined
  on domain $\domain(\Encode(c_1))\setminus\domain(\Encode(c_2))$, we will only
  use some of the information, 
  namely those labels assigned to nodes that represent one of the stacks we
  pass on some run from $c_1$ to $c_2$. The first component represents 
  a guess on the kind of segment starting at the corresponding
  stack. The numbers correspond to the enumeration in Lemma
  \ref{FormLemma}. The second component asserts  the final state
  of the corresponding segment. 
\end{remark}

As already mentioned, for each  encoding of
two configurations and each  certificate for substack reachability on
this encoding, there is a successor function. This successor function
chooses, according to the label of one representative, the
representative for the next stack. 

\begin{definition}
  Let $c_1, c_2$ be configurations such that $c_2=\Pop{2}^n(c_1)$ for some $n\in\N$. 
  Furthermore, let $f$ be a partial function from
  $\domain(\Encode(c_1))$ to $\{1,2,3,4,5\}\times Q$. 
  For \mbox{$d\in\domain(f)$}, 
  the \emph{successor of $d$ with respect to
    $f$} 
  is defined by case distinction on the first component of
  $f(d)$, denoted by $\pi_1(f(d))$, as follows. 
  \begin{enumerate}
  \item $\pi_1(f(d))=1$:
    If $d \in \{0\}^*$, there is no successor of $d$ with respect to
    $f$. 

    Otherwise, let $d'$ be the ancestor of $d$ such that
    $d=d'10^m$ for 
    some number \mbox{$m\in\N$}. Let 
    \mbox{$d''\in\{\varepsilon\}\cup\{0\{0,1\}^*\}$}
    be the lexicographically maximal word such that
    \mbox{$d'd''\in\domain(\Encode(c_1))$} ($d'd''$ is the maximal element in the
    subtree rooted at $d'0$ if $d'0$ is in the tree, otherwise we have
    $d'd''=d'$). We say $d'd''$ is the successor of $d$ with respect to
    $f$.  
  \item $\pi_1(f(d))\in\{2,3\}$:
    If $d\in\{0\}^*\{1\}^*$ the successor of $d$ with respect to
    $f$ is undefined. 
    
    Otherwise, let $d'$ be the element
    such that \mbox{$d=d'10^m1^n$} where $m>0$ and \mbox{$n\in\N$}.
    Let 
    $d''\in\{\varepsilon\}\cup\{0\{0,1\}^*\}$ 
    be the lexicographically maximal word such that
    \mbox{$d'd''\in\domain(\Encode(c_1))$.} We say $d'd''$ is the successor
    of $d$ with respect to 
    $f$.  
  \item $\pi_1(f(d))\in\{4,5\}$:
    If $d\in \{\varepsilon\}\cup \{0\}\{1\}^*$, then the successor of
    $d$ with respect to $\ReachguessSym$ is undefined. 

    Otherwise, 
    let $d'$ be the unique element such that $d=d'01^n$ for some
    $n\in\N$. Then $d'$ is the successor of $d$ with respect to
    $f$. 
  \end{enumerate}
\end{definition}
\begin{remark}
As already said in the informal description, the motivation of the
previous definition are the 
  following observations.
  \begin{enumerate}
  \item If $\pi_1(f(d))=1$, then
    the successor $\hat d$ of $d$ is chosen such that
    \begin{align*}
    \LeftStack(\hat d,\Encode(c_1))=\Pop{2}(\LeftStack(d,\Encode(c_1))).      
    \end{align*}
    If this is not possible, i.e., if 
    $\lvert \LeftStack(d, \Encode(c_1)) \rvert = 1$, the successor is
    undefined.    
  \item If $\pi_1(f(d))\in\{2,3\}$, then the successor
    $\hat d$ is chosen such that 
    \begin{align*}
    \LeftStack(\hat d,\Encode(c_1)) =
    \Collapse(\LeftStack(d),\Encode(c_1))      
    \end{align*}
    (assuming 
    that $\Lvl(\LeftStack(d,\Encode(c_1)))$ is $2$). If such an element
    does not exist, then the successor is undefined.  
  \item If $\pi_1(f(d))\in\{4,5\}$, then 
    the successor $\hat d$ of $d$ is chosen such that 
    $\LeftStack(\hat d,\Encode(c_1))$ is the maximal milestone of
    $\LeftStack(d,\Encode(c_1))$ satisfying 
    \begin{align*}
      \TOP{2}(\LeftStack(\hat d,\Encode(c_1))) =
      \TOP{2}(\Pop{1}(\LeftStack(d,\Encode(c_1)))).      
    \end{align*}
       If this is not possible, i.e., if $\TOP{2}(\LeftStack(d,
    \Encode(c_1)))=\bot$, then the successor is undefined. 
  \end{enumerate}
\end{remark}
\begin{example}\label{Example:CertificateSubstackReachability}
  Recall the run $\rho$ from example
  \ref{Example:ReachDecomposition}. 

  The decomposition of $\rho$ induces a
  certificate for substack reachability on
  \begin{align*}
    \Encode(q_1,s)\otimes \Encode(q_1,\bot_2).    
  \end{align*}
  This certificate is
  depicted in Figure \ref{fig:ValidReachability} (the bold labels are the
  values of the 
  certificate). 
  We only state the values of the certificate on the rightmost leaf of
  $\Encode(q_1,s)$ and the chain of successors with respect to this
  certificate. 
  These are the important values that 
  witness the existence of $\rho$.
\end{example}

\begin{figure}
  \centering
  $
  \begin{matrix}
             &  &           &          &  a & a \\
             &  &  c    &   b &  (a,2,3) & (a,2,3) \\
             &  & (c,2,1)   &  (c,2,1) &  (c,2,1) & (c,2,1) \\
  s\coloneqq &    \bot & \bot&\bot&\bot&\bot 
  \end{matrix}
  $ \\
  $
  \begin{xy}
    \xymatrix{
      &            &             & (a,1),\Box, \ar[r] &
      *+\txt{$\varepsilon$,$\Box$ \\
        $\mathbf{(\ref{1LoopPop},q_3)}$}\\
      &  *+\txt{$(c,1)$,$\Box$ \\
        $\mathbf{(\ref{LoopPop},q_1)}$}
      & *+\txt{$(b,1)$,$\Box$ \\
        $\mathbf{(\ref{FormReturn},q_1)}$}  
      & *+\txt{$(a,2)$,$\Box$ \\
        $\mathbf{(\ref{1LoopCol},q_3)}$} \ar[u] \\
      &  *+\txt{$(c,2)$,$\Box$ \\
        $\mathbf{(\ref{LoopCol},q_1)}$} \ar[u] \ar[r] &
      \varepsilon, \Box \ar[r]\ar[u] &
      \varepsilon, \Box \ar[u]\\
      (\bot, 1), (\bot, 1) \ar[r] & \varepsilon, \Box \ar[u] \\
      (q_1,q_1) \ar[u]
    }
  \end{xy}
  $
  \caption{Example of a valid certificate for substack reachability for
    $\Encode(q_1, s) \otimes \Encode(q_1, \bot_2)$.}
  \label{fig:ValidReachability}
\end{figure}

We will show that there is a close connection between runs
from some configuration $(q,s)$ to another configuration $(\hat q,
\hat s)$ where $\hat s=\Pop{2}^n(s)$ and certificates for substack
reachability. We prove that every such run induces a
certificate with certain properties. 
Since these properties are rather technical, we postpone the detailed
description of these properties for a short while.  
In the following, we first explain how to obtain such a certificate
from the run. Then we present the characterising properties of these
certificates. Finally, 
we show that each certificate with these properties actually
represents a run from $(q,s)$ to $(\hat q, \hat s)$. Hence, 
deciding the existence of such a run  reduces to deciding whether
there is such
a certificate. We then show that the latter problem is
$\MSO$-definable on the encoding of the configurations. 
From this, the regularity of the relation $A$ follows.

\begin{lemma} \label{Lemma:ValidCertificatefromRun}
  Let $c=(q,s)$ and $\hat c=(\hat q, \hat s)$ be configurations such
  that $\hat s=\Pop{2}^m(s)$ for some $m\in\N$. Let $\rho$ be a run
  from $c$ to $\hat c$ such that $\rho$ does not pass a proper
  substack of $\hat s$. 
  Assume that $\rho$ decomposes as 
  $\rho=\rho_1 \circ \rho_2 \circ \dots \circ \rho_n \circ\lambda$
  according to 
  Lemma \ref{FormLemma}. Then there 
  is a certificate for substack reachability $\ReachguessSym^\rho$ on
  $\Encode(c) \otimes \Encode(\hat c)$ such that the
  following conditions hold.
  
  There is a finite sequence $\bar t= t_1, t_2, \dots, t_n
  \in\domain(\Encode(c)) \setminus \domain(\Encode(\hat c))$ such that  
  \begin{enumerate}
  \item $t_1$ is the rightmost leaf of $\Encode(c) \otimes \Encode(\hat c)$,
  \item for each $1\leq i \leq n$, $\ReachguessSym^\rho(t_i) =
    (k_i,q_i)$ where $k_i$ is the form 
    of $\rho_i$ according to Lemma \ref{FormLemma} and $q_i$ is the
    final state of $\rho_i$. 
  \item $t_{i+1}$ is the successor of $t_i$ with respect to
    $\ReachguessSym^\rho$ for every $1\leq i<n$, and
  \item the successor of $t_n$ with respect to $\ReachguessSym^\rho$ is
    the rightmost leaf of $\Encode(\hat c)$.
  \end{enumerate}
\end{lemma}
\begin{proof}
  First of all, we define inductively the sequence $t_1, t_2, \dots,
  t_n$ and the values of $\ReachguessSym^\rho$ on these elements. 

  We write $k_i$ for the form of $\rho_i$ and $q_i$ for the final
  state of $\rho_i$. 
  Let $t_1$ be the rightmost leaf of $\Encode(q,s)$. 
  We define $\ReachguessSym^\rho(t_1):=(k_1,q_1)$. 
  Now assume that we have already defined $t_i$ for some $1\leq i<n$. We
  define $\ReachguessSym^\rho(t_i):=(k_i,q_i)$. If the successor of
  $t_i$ with respect to $\ReachguessSym^\rho$ exists, it is uniquely
  defined and we call it $t_{i+1}$. We proceed with this definition
  until $i=n$ or until there is some $l<n$ such that the successor
  of $t_l$ with respect to $\ReachguessSym^\rho$ is not defined. 


  In order to prove that $\ReachguessSym^\rho$ can be extended to a
  well-defined certificate for substack reachability satisfying
  conditions 1--4, we show a stronger 
  claim. 
  For $1\leq i \leq n+1$ such that $t_i$ is defined, set
  $s_{t_i}:=\LeftStack(t_i,\Encode(c))$.  
  For $1\leq i \leq n$ let $s_i$ be the stack of $\rho_i(0)$. Let
  $s_{n+1}$ be the stack of $\lambda(0)$.  
  \begin{claim}
    For all $1\leq i \leq n$ such that $t_i$ is defined, $t_{i+1}$ is
    also defined. 
    Furthermore, 
    if $t_i$ is defined for some $1\leq i \leq n+1$, then  
    $s_i$ and $s_{t_i}$ are similar in the following sense:
    \begin{enumerate}
    \item $\TOP{2}(s_i)=\TOP{2}(s_{t_i})$ and \label{FirstPartofClaim}
    \item there is an $n_i\geq 1$ such that $\Pop{2}(s_{t_i})=\Pop{2}^{n_i}(s_i)$.
    \item Moreover, if $i=1$ or $\rho_{i-1}$ is of the form F\ref{FormReturn},
      F\ref{1LoopCol}, or F\ref{LoopCol}. (with respect to Lemma
      \ref{FormLemma}), then $s_i=s_{t_i}$. This means
      that the substack represented by $t_i$ is the initial stack of
      $\rho_{i}$ whenever $\rho_{i-1}$ ended in a $\Pop{2}$ or $\Collapse$
      of level $2$.  \label{LastPartofClaim}
    \end{enumerate}
  \end{claim}
  Before we prove the claim, let us explain how the lemma follows from the claim. Note that $t_1$ is defined and
  $s_{t_1}= \LeftStack(t_1,\Encode(q,s)) = s = s_1 = \rho_1(0)$. 
  Due to the claim, $t_1, t_2, t_3, \dots, t_{n+1}$ are
  defined. According to Lemma \ref{FormLemma}, $\rho_n$ is of the form  
  F\ref{FormReturn},  F\ref{1LoopCol}, or F\ref{LoopCol}. Thus, the 
  claim implies that 
  $\LeftStack(t_{n+1}, \Encode(c)) = s_{t_{n+1}} = s_{n+1} =
  \lambda(0) = \hat s$.  
  Thus, $t_{n+1}$ is the rightmost leaf of $\Encode(\hat c)$. Since
  the successor with respect to $\ReachguessSym^\rho$ of  
  some node is always lexicographically smaller than  this node,
  $t_i>_{\mathrm{lex}}t_{n+1}$ for all $1\leq i \leq n$.  
  Hence, $\{t_1, t_2, \dots,
  t_n\}\subseteq\domain(\Encode(c))\setminus\domain(\Encode(\hat c))$.  
  Thus, we can extend the partial definition of $\ReachguessSym^\rho$
  to a map from  
  $\domain(\Encode(c))\setminus\domain(\Encode(\hat c))$ to
  $\{1,2,3,4,5\}\times Q$. Furthermore,  
  \mbox{$\bar t=t_1, t_2, \dots, t_n$} satisfies items 1--3 by
  definition of the $t_i$. $\bar t$ also satisfies  
  item 4 because  we proved that the successor $t_{n+1}$ of $t_n$
  is the rightmost leaf of $\Encode(\hat c)$.  

  Now, we prove  the claim. 
  Assume 
  that there is  some  $i \leq n$ such that $t_i$ is defined.
  Furthermore, assume that $s_i$ and $s_{t_i}$ are similar, i.e., 
  $s_i$ and $s_{t_i}$ satisfy conditions 1--3 of the claim.   
  We distinguish the following cases according to the
  form of $t_i$:
  \begin{enumerate}
  \item Consider the case $k_i=1$. In this case,  $\rho_i$ is a return. 
    If $i>1$, then Lemma \ref{FormLemma}
    implies that $\rho_{i-1}$ is of the form  F\ref{FormReturn},
    F\ref{1LoopCol}, or F\ref{LoopCol}. Thus,  
    item \ref{LastPartofClaim} of the claim implies that $s_i=s_{t_i}$. 

    Due to $k_i=1$, the successor
    of $t_i$ -- if defined -- is  a node $t_{i+1}$ such that
    \begin{align*}
      \LeftStack(t_{i+1}, \Encode(c)) = \Pop{2}(\LeftStack(t_i,
      \Encode(c)))=\Pop{2}(s_{t_i}).      
    \end{align*}
    $\rho_i$ is a return
    starting at $s_i = s_{t_i}$ whence $\rho_i$ ends in
    $s_{i+1}:=\Pop{2}(s_{t_i})$.  
    Thus, 
     we conclude that 
    $\lvert s_{t_i} \rvert \geq 2$ whence $t_{i+1}$ is defined.
    Furthermore, note that 
    \begin{align*}
      s_{t_{i+1}}=\LeftStack(t_{i+1},
      \Encode(c)) =  \Pop{2}(s_{t_i})  = s_{i+1}      
    \end{align*}
    whence
    $t_{i+1}$ satisfies item \ref{LastPartofClaim} of the claim. 
  \item Consider the case $k_i\in\{2,3\}$. This means that $\rho_i$ ends with a collapse
    of level $2$ from a stack with topmost word $\TOP{2}(s_i)$. Thus,
    $\Lvl(s_i)=2$. 
    Since $k_i\in\{2,3\}$, the
    successor of $t_i$  -- if defined -- is a node $t_{i+1}$ such that
    \begin{align*}
      \LeftStack(t_{i+1}, \Encode(c)) = \Collapse(\LeftStack(t_i,
      \Encode(c)) = \Collapse(s_{t_i}).      
    \end{align*}
    Due to the form of
    $\rho_i$,  
    $\Collapse(s_i)$ is defined. 
    Due to item \ref{FirstPartofClaim} of
    the claim,
    $\TOP{1}(s_{t_i})$ and $\TOP{1}(s_i)$ coincide. Thus,
    $\Collapse(s_{t_i})$ is also 
    defined. But then $t_{i+1}$ is defined whence the first part of
    the claim holds.
    Furthermore, item 2 of the claim implies that
    that 
    \begin{align*}
      s_{t_{i+1}} = \Collapse(s_{t_i}) = \Collapse(s_i) = s_{i+1}      
    \end{align*}
    whence $t_{i+1}$ satisfies the second part of the claim. 
  \item Consider the case $k_i\in\{4,5\}$. This means that $\rho_i$ is a loop or
    $1$-loop followed by a $\Pop{1}$ or $\Collapse$ of level $1$. 
    Since the topmost word of $s_i$ and the topmost word
    before the last operation of $\rho_i$ agree, 
    $\lvert \TOP{2}(s_i) \rvert >1$ holds. Due to item 1 of the claim,
    $\lvert\TOP{2}(s_{t_i})\rvert > 1$ follows. This implies
    that there is some node $d\in\Encode(c)$ such that $t_i=d01^l$
    for some $l\in\N$. Since $k_i\in\{4,5\}$, $d$ is the successor of
    $t_i$ with respect to $\ReachguessSym^\rho$, i.e., $t_{i+1}=d$
    whence the first part of the claim holds. 

    For the second part, note that there is some $m\geq 1$ such that
    $\Pop{2}(s_i) = \Pop{2}^m(s_{i+1})$ due to the definition of loops and $1$-loops. 

    Since the left stack induced by $t_{i+1}$ is a milestone of the one induced by $t_i$,
    there is an $n\geq 1$ such that 
    $\Pop{2}(s_{t_{i+1}}) = \Pop{2}^n(s_{t_{i}})$. Thus, by item 2 of the claim, 
    we obtain that $\Pop{2}(s_{t_{i+1}}) = \Pop{2}^{n+m+n_i-2}(s_{i+1})$. 
    Since $n_{i+1}:=n+m+n_i-2\geq 1$, item 2 of the claim holds for $t_{i+1}$. 
     
    Furthermore, since $\TOP{2}(s_{i+1}) = \TOP{2}(\Pop{1}(s_i))$  
    and $\TOP{2}(s_{t_{i+1}}) = \TOP{2}(\Pop{1}(s_{t_i}))$, item 1 of
    the claim carries  
    over from $t_i$ to $t_{i+1}$. 
    This completes the proof that $s_{i+1}$ and $s_{t_{i+1}}$ are
    similar in the sense of the claim.  \qedhere
  \end{enumerate}
\end{proof}
\begin{remark}
  In the following we say that a certificate for substack reachability
  $\ReachguessSym$  \emph{represents} $\rho$ if 
  it coincides with $\ReachguessSym^\rho$ on $\{t_1, t_2, \dots, t_n\}$.
\end{remark}

In the next lemma, we collect important properties of a certificate
which represents some run. Afterwards, we turn these properties into the
defining conditions of \emph{valid} certificates. This terminology is
justified because each valid certificate represents in fact some run.

\begin{lemma} \label{ValidCertificateConditions}
  Let $(q, s), (\hat q, \hat s)$ be configurations such that $\hat s =
  \Pop{2}^m(s)$ for some $m\in\N$. 
  Let $\ReachguessSym$ be a certificate representing a run $\rho$ from
  $(q,s)$ to $(\hat q, \hat s)$. 

  Then there is an $n\in\N$ and a finite sequence $t_1, t_2, \dots, t_n \in
  \domain(\Encode(q,s)) \setminus \domain(\Encode(\hat q, \hat s))$ with the
  following properties (setting $(k_i,q_i):=\ReachguessSym(t_i)$ and
  $q_0:=q$): 
  \begin{enumerate}[\normalfont{A}1.]
  \item\label{Axiom1} $t_1$ is the rightmost leaf of
      $\Encode(q,s)$,  
  \item\label{Axiom2} for all $1\leq i < n$, the successor of $t_i$
    with respect to 
    $\ReachguessSym$ is $t_{i+1}$,
  \item\label{Axiom3} the successor of $t_n$ with respect to
    $\ReachguessSym$ is the 
    rightmost leaf of $\Encode(\hat q, \hat s)$,
  \item\label{Axiom4} $k_n\in\{1,2,3\}$,
  \item\label{Axiom5} $(q_n, \hat q)\in \exHighLoops(\hat s)$, i.e.,
    there is a high 
    loop from $(q_n, \hat s)$ to $(\hat q, \hat s)$,
  \item\label{Axiom6} if $k_i\in\{4,5\}$ for some  $i<n$ then there is
    a $j>i$ such 
    that $k_j\in\{2,3\}$ and $k_l\neq 1$ for all $i<l<j$,
  \item\label{Axiom7} For each $1\leq i \leq n$, the stack induced by $t_i$ 
    satisfies in dependence of the value of $k_i$ a certain assertion
    as follows:
    \begin{enumerate}
    \item if $k_i=1$ then $(q_{i-1}, q_{i})\in
      \exReturns(\LeftStack(t_i, \Encode(q,s)))$, 
    \item if $k_i=2$ then $\Lvl(\LeftStack(t_i, \Encode(q,s)))=2$ and there is a
      $q'\in Q$ and a $\gamma\in\Gamma$ such that 
      \begin{align*}
        &(q_{i-1},q')\in\exOneLoops(\LeftStack(t_i, \Encode(q,s))) \text{ and}\\
        &(q',\Sym(\LeftStack(t_i, \Encode(q,s))),\gamma,
        q_{i},\Collapse)\in \Delta,    
      \end{align*}
    \item if $k_i=3$ then $\Lvl(\LeftStack(t_i, \Encode(q,s)))=2$ and
      there is some 
      $q'\in Q$ and some $\gamma\in\Gamma$ such that 
      \begin{align*}
        &(q_{i-1},q')\in\exLoops(\LeftStack(t_i, \Encode(q,s))) \text{ and}\\
        &(q',\Sym(\LeftStack(t_i,
        \Encode(q,s))),\gamma,q_{i},\Collapse)\in \Delta,  
      \end{align*}
    \item if $k_i=4$ then there is some $q'\in Q$ and some
      $\gamma\in\Gamma$ such that
      \begin{align*}
        &(q_{i-1},q')\in\exLoops(\LeftStack(t_i, \Encode(q,s))))
        \text{ and either}\\
        &(q',\Sym(\LeftStack(t_i,
        \Encode(q,s))),\gamma,q_{i},\Pop{1})\in 
        \Delta \text{ or }\\ 
        &\Lvl(\LeftStack(t_i, \Encode(q,s)))=1 \text{ and } 
        (q',\Sym(\LeftStack(t_i,
        \Encode(q,s))),\gamma,q_{i},\Collapse)\in \Delta.        
      \end{align*}
    \item if $k_i=5$ then there is a $q'\in Q$ such that
      \begin{align*}
        &(q_{i-1}, q') \in\exOneLoops(\LeftStack(t_i, \Encode(q,s))))
        \text{ and either}\\ 
        &(q',\Sym(\LeftStack(t_i,
        \Encode(q,s))),\gamma,q_{i},\Pop{1})\in 
        \Delta \text{ or}\\ 
        &\Lvl(\LeftStack(t_i, \Encode(q,s)))= 1 \text{ and } 
        (q',\Sym(\LeftStack(t_i,
        \Encode(q,s))),\gamma,q_{i},\Collapse)\in \Delta. 
      \end{align*}
    \end{enumerate}
  \end{enumerate}
\end{lemma}
\begin{proof}
  Let $\ReachguessSym$ represent a run $\rho$. 
  
  There is a unique sequence $t_1, t_2, \dots, t_n$ of maximal length
  that satisfies
  A\ref{Axiom1} and A\ref{Axiom2}. Furthermore, the
  previous 
  lemma showed that $t_n$ satisfies A\ref{Axiom3}.  
  
  From the previous lemma we also know that $\ReachguessSym(t_i)$
  encodes the form and the final state of $\rho_i$ where 
  $\rho=\rho_1 \circ \rho_2 \circ \dots \circ \rho_n\circ \lambda$ is
  the decomposition of $\rho$ according to Lemma \ref{FormLemma}. 
  Thus, $k_n$ is the form of $\rho_n$. Hence, Lemma \ref{FormLemma}
  implies that $k_n\in\{1,2,3\}$. 

  $q_n$ is the final state of $\rho_n$
  and due to Lemma \ref{FormLemma}, $\lambda$ is a high loop from 
  $(q_n,\hat s)$ to $(\hat q, \hat s)$. Thus, $\lambda$  witnesses that
  $(q_n, \hat q)\in\exHighLoops(\hat s)$. This is exactly the assertion
  of A\ref{Axiom5}. 

  A\ref{Axiom6} is also a direct consequence of Lemma
  \ref{FormLemma}:
  if there is some  $\rho_i$ of the form F\ref{LoopPop} or
  F\ref{1LoopPop}, then there is a $j>i$ such that 
  $\rho_j$ is of the form F\ref{1LoopCol} or F\ref{LoopCol}, and for all 
  $i < k < j$, $\rho_k$ is not of the
  form F\ref{FormReturn}. From the correspondence between the form  of 
  $\rho_l$ and the 
  value of $k_l$ for all $1\leq l \leq n$,   A\ref{Axiom6} follows directly. 

  A\ref{Axiom7} is a consequence of the claim in the previous
  proof. There we showed that  
  \begin{align} \label{Certifictetopwordagree}
    \TOP{2}(\LeftStack(t_i, \Encode(q,s)))=\TOP{2}(\rho_i(0)).
  \end{align}
  Thus,
  $\exReturns, \exLoops$ and  $\exOneLoops$ agree on the stacks
  $\LeftStack(t_i, \Encode(q,s))$ and $\rho_i(0)$. We conclude by case
  distinction on $k_i$ as follows. 
  \begin{enumerate}
  \item[$k_i=1$] This implies that $\rho_i$ is a return from $\rho_i(0)$
    to $(q_i, \Pop{2}(\rho_i(0)))$. By definition, the state of
    $\rho_i(0)$ is $q_{i-1}$. Thus, $\rho_i$ witnesses
    \begin{align*}
      (q_{i-1},q_i)\in \exReturns(\rho_i(0))=
      \exReturns(\LeftStack(t_i,\Encode(q,s))).      
    \end{align*}
  \item[$k_i=2$] This implies that $\rho_i$ is a $1$-loop followed by
    a $\Collapse$ of level $2$. Let $j$ be the position just before
    this $\Collapse$, i.e., $j:=\length(\rho_i)-1$. 
    Let $q'$ be the state of $\rho_i(j)$. 
    Now, $\rho_i{\restriction}_{[o,j]}$ witnesses the existence of a
    $1$-loop from state $q_{i-1}$ to state $q'$ on topmost word
    $\TOP{2}(\rho_i(0))$. Thus, 
    $(q_{i-1},q')\in\exOneLoops(\LeftStack(t_i, \Encode(q,s)))$. 
   
    Due to (\ref{Certifictetopwordagree}) and the definition of
    $1$-loops, we have
    \begin{align*}
      \TOP{1}(\LeftStack(t_i, \Encode(q,s))) =\TOP{1}(\rho_i(0)) =
      \TOP{1}(\rho_i(j)). 
    \end{align*}
    By definition of  $\rho_i$, $\Lvl(\rho_i(j))=2$. 
    We conclude that 
    \begin{align*}
      \Lvl(\LeftStack(t_i, \Encode(q,s)))=\Lvl(\rho_i(j))=2.      
    \end{align*}
    Since $\rho_i$ performs a collapse at
    $j$, there is some transition 
    \begin{align*}
      (q', \Sym(\rho(j)), \gamma, q_i, \Collapse)\in\Delta.      
    \end{align*}
    Due to
    \mbox{$\Sym(\rho(j)) = \Sym(\LeftStack(t_i, \Encode(q,s)))$}, this
    transition witnesses that 
    \begin{align*}
      (q',\Sym(\LeftStack(t_i,
      \Encode(q,s))),\gamma, q_i,\Collapse)\in\Delta.
    \end{align*}  
  \item[$k_i=3$] Replacing the role of $1$-loops by loops, we can copy
    the proof from the previous case word by word.
  \item[$k_i=4$] This implies that $\rho_i$ is a loop followed by a
    $\Pop{1}$ transition or a $\Collapse$ of level $1$. 
    We set $j:=\length(\rho_i)-1$ and $q'$ to be the state of
    $\rho_i(j)$. Completely analogous to the previous case, one derives
    that $(q_{i-1}, q')\in\exLoops(\LeftStack(t_i, \Encode(q,s)))$. 
    
    The transition at $\rho_i(j)$ is a
    $\Pop{1}$ or a $\Collapse$ of level $1$. Thus, this transition is
    either $(q', \Sym(\rho_i(j)), \gamma, q_i, \Pop{1})$ or
    $(q', \Sym(\rho_i(j)), \gamma, q_i, \Collapse)$ and
    \mbox{$\Lvl(\rho_i(j))=1$.} 
    Due to (\ref{Certifictetopwordagree}), 
    $\TOP{2}(\LeftStack(t_i, \Encode(q,s))) =  \TOP{2}(\rho_i(j))$.
    Thus, this transition is also applicable to $\LeftStack(t_i,
    \Encode(q,s))$. This completes the proof in the case $k_i=4$.  
  \item[$k_i=5$] This case is completely analogous to the
    previous one: we only have to replace loops by $1$-loops. \qedhere
  \end{enumerate}
\end{proof}

\begin{definition}
  Let  $c=(q, s)$ and $\hat c=(\hat q, \hat s)$ be configurations such that
  $\hat s = \Pop{2}^n(\bar s)$. 
  Let $\ReachguessSym :
  \domain(\Encode(c))\setminus\domain(\Encode(\hat c) \to
  \{1,2,3,4,5\}\times Q$ be a certificate for substack reachability on
  $c$ and $\hat c$. 
  Setting $q_0:=q$, 
  we call $\ReachguessSym$ \emph{valid} if 
  there is an $n\in\N$ and a finite sequence $t_1, t_2, \dots, t_n \in
  \domain(\Encode(c)) \setminus \domain(\Encode(\hat c))$ 
  which satisfies conditions  A\ref{Axiom1} -- A\ref{Axiom7} from Lemma
  \ref{ValidCertificateConditions}.  
\end{definition}

The next lemma shows the tight correspondence 
between valid certificates of substack
reachability  and runs.  
For $q,\hat q \in Q$ two states, $s$ some stack and $\hat s=\Pop{2}^n(s)$
for some $n\in\N$, 
there is a run from $(q,s)$ to $(\hat q, \hat s)$ if and only if there
is a valid certificate for substack reachability for $(q,s)$ and
$(\hat q, \hat s)$. 

\begin{lemma}
  Let  $c=(q,s)$ and $\hat c=(\hat q,\hat s)$ be configurations
  such
  that $\hat s=\Pop{2}^m(s)$ for some $m\in\N$. 
  There is a run from $c$
  to $\hat c$ which does not visit proper substacks of $\hat s$ if and
  only if there is a valid certificate for substack 
  reachability 
  \begin{align*}
    \ReachguessSym:\domain(\Encode(c))\setminus
    \domain(\Encode(\hat c)) \to \{1,2,3,4,5\}\times Q.     
  \end{align*}
\end{lemma}
\begin{proof}
  The implication from left to right follows from Lemma
  \ref{Lemma:ValidCertificatefromRun}. 
  
  For the proof from right to left assume that $\ReachguessSym$ is a valid
  certificate for substack reachability for $c$ and $\hat c$. 

  Then there is a sequence $t_1, t_2, \dots, t_n$ in
  $\domain(\ReachguessSym)$ that witnesses the conditions
  A\ref{Axiom1}--A\ref{Axiom7}. 
  
  We  now construct runs $\rho_0, \rho_1, \rho_2, \dots, \rho_n,
  \rho_{n+1}$ such 
  that $\rho_i$ is an initial segment of $\rho_{i+1}$ for each $i \leq
  n$.
  The run $\rho_{n+1}$ is  then  a run from $c$ to $\hat c$. 

  Before we start the construction, let us define some notation. 
  For all $1\leq i \leq n$, let $(q_i,
  k_i):=\ReachguessSym(t_i)$ and let 
  $s_{t_i}:=\LeftStack(t_i, \Encode(c))$.
  Furthermore, for reasons of convenience, we set $q_0:=q$ and
  we set $t_{n+1}$ to be the rightmost leaf of $\Encode(\hat c)$.   
  As soon as $\rho_{i-1}$ is defined for some $1\leq i \leq n$, we
  denote by $s_i$ the last 
  stack of $\rho_{i-1}$.   
  
  We  define $\rho_0$ to be a run of length $0$ with $\rho_0(0):=c$.

  During the construction of $\rho_i$ for $1\leq i \leq n$, we
  preserve the following conditions:
  \begin{enumerate}
  \item the last state of $\rho_{i-1}$ is $q_{i-1}$,
  \item $\TOP{2}(s_i)=\TOP{2}(s_{t_i})$, and 
  \item there is an $n_i\geq 1$ such that
    $\Pop{2}(s_{t_i})=\Pop{2}^{n_i}(s_i)$. 
  \item Moreover, if $i=1$ or $k_{i-1}\in\{1,2,3\}$, then
    $s_i=s_{t_i}$. 
  \end{enumerate}
  Note that $\rho_0(0)=(q,s_0)= (q_0,s_{t_0})$ by definition whence
  for $i=1$ these conditions are satisfied. 
  
  Now assume that $\rho_{i-1}$ is defined for some $1\leq
  i\leq n$ such that these conditions are satisfied. 
  By case distinction on the value of $k_i$ we define $\rho_i$ as
  follows.
  \begin{enumerate}
  \item[$k_i=1$] Since $\ReachguessSym$ satisfies A\ref{Axiom6},
    $i=1$ or $k_{i-1}\in\{1,2,3\}$. Thus, $s_i=s_{t_i}$ by induction
    hypothesis. Due to A\ref{Axiom7}, we have $(q_{i-1},
    q_i)\in\exReturns(s_{t_i})=\exReturns(s_i)$. Hence, there is a
    return $\hat\rho$ from $(q_{i-1},s_i)$ to $(q_i, \Pop{2}(s_i))$. 
    We set 
    \begin{align*}
      \rho_i:=\rho_{i-1}\circ \hat\rho.       
    \end{align*}
    Due to A\ref{Axiom2} and A\ref{Axiom3},  the successor of $t_i$
    with respect to  
    $\ReachguessSym$ is $t_{i+1}$. 
    Since $k_i=1$, 
    $s_{t_{i+1}}=\LeftStack(t_{i+1},\Encode(c)) =
    \Pop{2}(s_{t_i}) = \Pop{2}(s_i)=s_{i+1}$.
  \item[$k_i=2$] Due to A\ref{Axiom7}, $\Lvl(s_{t_i})=2$ and there
    is a $q'\in Q$ and a $\gamma\in\Gamma$ such that
    \begin{align*}
      &(q_{i-1}, q')\in\exOneLoops(s_{t_i})\\
      \text{and } &\delta:=(q', \Sym(s_{t_i}), 
      \gamma, q_i, \Collapse)\in\Delta.       
    \end{align*}
    Since the topmost words of $s_i$ and $s_{t_i}$ agree, there is some
    stack $s'$ such that there is a $1$-loop $\hat\lambda$ from
    $(q_{i-1}, s_i)$ to $(q', s')$. By definition of a $1$-loop,
    \begin{align*}
      \TOP{2}(s') = \TOP{2}(s_i) = \TOP{2}(s_{t_i}).      
    \end{align*}
    Thus,
    $\Lvl(s')=2$ and $\hat\lambda$ can be extended by $\delta$. 
    We write $\hat\lambda^+$ for $\hat\lambda$ extended by one
    application of $\delta$. 
    Since $\hat\lambda$ is a $1$-loop, it does not visit any 
    substacks of $\Pop{2}(s_i)$. Thus, 
    \mbox{$\Collapse(s_i)=\Collapse(s')$.} 
    Set $\rho_i:=\rho_{i-1}\circ \hat\lambda^+$. 
    By assumption, we conclude that 
    \begin{align*}
      s_{t_{i+1}} = \Collapse(s_{t_i}) = \Collapse(s_i) =
      \Collapse(s') = s_{i+1}.
    \end{align*}
    Thus, the last configuration of $\rho_i$ is $(q_i, s_{t_{i+1}})$.
  \item[$k_i=3$] We can copy the argument from the case $k_i=2$: just
    replace $1$-loops by loops. Then we obtain a run $\rho_i$ that
    ends in $(q_i, s_{t_{i+1}})$.
  \item[$k_i=4$] 
    Due to condition A\ref{Axiom7}, there
    is a $q'\in Q$ and a $\gamma\in\Gamma$ such that
    $(q_{i-1}, q')\in\exLoops(s_{t_i})$ and 
    $\delta_p:=(q', \Sym(s_{t_i}), \gamma, q_i, \Pop{1})\in\Delta$ 
    or 
    $\delta_c:=(q', \Sym(s_{t_i}), \gamma, q_i, \Collapse)\in\Delta$ and
    $\Lvl(s_{t_i})=1$.
 
    Since the topmost words of $s_i$ and $s_{t_i}$ agree, there is 
    a loop $\hat\lambda$ from
    $(q_{i-1}, s_i)$ to $(q', s_i)$. Due to
    $\TOP{2}(s_i) = \TOP{2}(s_{t_i})$, it follows that
    $\Lvl(s_i)=\Lvl(s_{t_i})$. We conclude that $\delta_c$ (or
    $\delta_p$, respectively) can be applied to the last 
    configuration of $\hat\lambda$. In both cases, the resulting configuration
    is $(q_i, \Pop{1}(s_i))$. Writing $\hat\lambda^+$ for $\hat\lambda$
    extended by $\delta_p$ or $\delta_c$, we set
    $\rho_i:=\rho_{i-1}\circ \hat\lambda^+$. 
    
    By definition of $t_{i+1}$, $s_{t_{i+1}}$ is the maximal milestone
    of $s_{t_i}$ such that 
    \begin{align*}
      \TOP{2}(s_{t_{i+1}})= \TOP{2}(\Pop{1}(s_{t_i})) = \TOP{2}(s_{i+1}).   
    \end{align*}
    We still have to show that there is some $n_{i+1}\geq 1$ such that
    $\Pop{2}(s_{t_{i+1}})= \Pop{2}^{n_{i+1}}(s_{i+1})$. 
    
    Since $s_{t_{i+1}}$ is a milestone of $s_{t_i}$, there is some
    $j\geq 1$ such that $\Pop{2}(s_{t_{i+1}}) =
    \Pop{2}^j(s_{t_{i+1}})$. 
    Since $\hat\lambda$ is a loop, we have $\Pop{2}(s_{i+1}) =
    \Pop{2}(s_{i})$. By assumption on $\rho_{i-1}$, we obtain that
    \begin{align*}
      \Pop{2}(s_{t_{i+1}}) = \Pop{2}^{j}(s_{t_{i}}) =
      \Pop{2}^{j+n_i-1}(s_i) = \Pop{2}^{j+n_i-1}(s_{i+1}).  
    \end{align*}
    Let
    $n_{i+1}:=j+n_i-1$. We conclude by noting that $n_{i+1}\geq 1$. 
  \item[$k_i=5$] This case is analogous to the previous
    one. We can replace the loop $\hat\lambda$ in the previous case by
    some $1$-loop from $(q_{i-1}, s_i)$ to some $(q', s')$ where
    $\Pop{2}(s_i)$ is a substack of $s'$ and
    $\TOP{2}(s')=\TOP{2}(s_i)$. The rest of the argument is then
    completely analogous. 
  \end{enumerate}
  Repeating this construction for all $i\leq n$, we define a run
  $\rho_n$ with last state $q_i$. Due to A\ref{Axiom4}, the last step
  in this 
  construction uses one of the first three cases. Thus, 
  \begin{align*}
    s_{n+1}=
    s_{t_{n+1}}=\LeftStack(t_{n+1}, \Encode(c)) = \hat s.    
  \end{align*}
  Note that the
  last equality is due to A\ref{Axiom3}. Thus, $\rho_n$ ends in 
  $(q_n, \hat s)$. 

  Due to A\ref{Axiom5}, there is a high loop $\lambda$ from  $(q_n,
  \hat s)$ to $(\hat q, \hat s)$. We set
  $\rho_{n+1}:=\rho_n\circ\lambda$. This completes the proof because
  $\rho_{n+1}$  is a  run from $(q, s)$ to $(\hat q, \hat s)$ that
  does not visit any proper substack of $\hat s$. 
\end{proof}

We have seen that there is a run from $c=(q,s)$ to $\hat c= (\hat q,
\hat s)$ 
for $\hat s = \Pop{2}^n(s)$ if and only if there is a valid
certificate for substack reachability on
$\Encode(c)\otimes \Encode(\hat c)$. 
The final step of the analysis of runs of this form is the following
lemma. We prove that the set of all pairs of configurations
$(c,\hat c)$ of the form mentioned 
above is an automatic relation via the encoding $\Encode$. 
We show that the set of encodings of such pairs is
\MSO-definable. The automaticity of the relation follows from the
correspondence of \MSO and automata on trees.

\begin{lemma} \label{MSODefinabilityofA}
  There is a formula in \MSO that defines the set
  \begin{align*}
    S:=\{(\Encode(c),\Encode(\hat c)): \exists
    \ReachguessSym:\domain(\Encode(c))\setminus\domain(\Encode(\hat
    c)) \to \{1,\dots,5\} \times Q, \ReachguessSym \text{ is valid}\}
  \end{align*}
\end{lemma}
\begin{proof}
  Certificates for substack reachability are only defined for
  configurations $c=(q,s)$ and $\hat c=(\hat q, \hat s)$ where $\hat s
  = \Pop{2}^m(s)$ for some $m\in\N$.  
  Note that this necessary condition is satisfied by a pair $(c,\hat
  c)$ if and only if there is some node
  $d\in\domain(\Encode(c))$ such that $d0\notin\domain(\Encode(c))$
  and $\LeftStack(d, \Encode(c)) = \hat s$. 
  These pairs of configurations are obviously $\MSO$ definable. 

  In this proof, we use the following claim.
  \begin{claim}
    There is an \MSO formula $\varphi$ such that for each
    certificate for substack reachability $\ReachguessSym$ on $c$ and
    $\hat c$ the following
    holds. 
    $\ReachguessSym$ is valid if and only if
    \begin{align*}
    \Encode(c)\otimes\Encode(\hat c)  \otimes \ReachguessSym \otimes
    T_{HL} \otimes T_{L} \otimes  T_{1L} \otimes T_{R} \models \varphi      
    \end{align*}
    where $T_{HL}, T_{L}, T_{1L}, T_{R}$ are trees 
    such that 
    \begin{align*}
      &T_{HL} \text{ encodes the mapping } d\mapsto
      \exHighLoops(\LeftStack(d, \Encode(c))), \\
      &T_{L} \text{ encodes the mapping } d\mapsto
      \exLoops(\LeftStack(d, \Encode(c))), \\
      &T_{1L} \text{ encodes the mapping } d\mapsto
      \exOneLoops(\LeftStack(d, \Encode(c))), \text{ and}\\
      &T_{R} \text{ encodes the mapping } d\mapsto
      \exReturns(\LeftStack(d, \Encode(c))).
    \end{align*}
  \end{claim}
  Before we prove this claim, we show that it implies the lemma.
  Due to Propositions
  \ref{Prop:AutomatonForReturns},
  \ref{Prop:AutomatonForLoops}, and
  \ref{OneLoopsAutomaton}, 
  the trees $T_{HL}, T_L, T_{1L}$ and $T_R$ are definable on
  $\Encode(c)$ using $\MSO$.  
  
  Furthermore, in Remark  
  \ref{CertificatesofSubstackReachabilityMSODefinable}  we saw that
  $\MSO$ can 
  express the existence of a certificate for substack reachability on
  $\Encode(c)\otimes \Encode(\hat c)$. Thus, given the formula
  $\varphi$ from the claim, we can construct a formula $\psi$
  asserting that ``there is a certificate for substack reachability
  $\ReachguessSym$ on 
  $\Encode(c)\otimes \Encode(\hat c)$ such that 
  \begin{align*}
    \Encode(c) \otimes \Encode(\hat c) \otimes \ReachguessSym \otimes
    T_{HL}\otimes T_L \otimes T_{1L} \otimes T_R \models
    \varphi\text{''}.     
  \end{align*}
  This \MSO formula defines the set
  $S$. 

  Let us now prove the claim. As an abbreviation, we write
  \begin{align*}
    \mathfrak{A}:= \Encode(c)\otimes\Encode(\hat c)  \otimes
    \ReachguessSym \otimes T_{HL} \otimes T_{L} \otimes  T_{1L} \otimes
    T_{R}.     
  \end{align*}
  We provide formulas that assert the conditions
  A\ref{Axiom1}--A\ref{Axiom7}.  
  The rightmost leaf of $\Encode(c)$ is of course \MSO-definable in
  $\mathfrak{A}$. Furthermore, the successor of a given node $d$ with
  respect to $\ReachguessSym$ is also \MSO-definable. Thus, 
  the uniquely defined maximal set $T:=\{t_1, t_2,
  \dots, t_n\}$  such that the sequence $t_1, t_2, \dots, t_n$
  satisfies condition  A\ref{Axiom1} and A\ref{Axiom2} is \MSO-definable. 

  Note that $i<j$ is equivalent to 
  $t_j \leq_{\mathrm{lex}} t_i$. Since the lexicographic order on
  $\domain(\Encode(c))$ is \MSO-definable, the successor of $t$ with
  respect to $\ReachguessSym$ is definable for each $t\in T$. 
  But this implies directly that the lexicographically minimal element
  in $T$ is $t_n$. Thus,  $t_n$ is definable and we can express ``the
  successor of $t_n$ with 
  respect to $\ReachguessSym$ is the rightmost leaf of $\Encode(\hat q,
  \hat s)$'' in an \MSO formula. This formula expresses A\ref{Axiom3}. 

  Furthermore, we conclude that
  ``$\ReachguessSym(t_n)=(k_n,q_n)$ such that $k_n\in\{1,2,3\}$'' is
  definable because $t_n$ is definable. Thus, we can express A\ref{Axiom4}.
  
  Next, we define a formula expressing A\ref{Axiom5}. 
  The label of the root of $\mathfrak{A}$ encodes the state
  $\hat q$.  ``$(q_n, \hat q) \in
  \exHighLoops(\hat s)$'' is expressible because $\hat q$ is encoded
  in the label of $t_n$ and $\exHighLoops(\hat
  s)$ is encoded in the label of the rightmost leaf of $\Encode(\hat
  q, \hat s)$. This  leaf is  definable in $\mathfrak{A}$. Thus, we
  conclude that A\ref{Axiom5} is expressible by some \MSO formula. 

  Condition A\ref{Axiom6} says that after some $t_i$ with $k_i\in\{4,5\}$
  there is a $t_j$ with $k_j\in\{2,3\}$ before the next occurrence of
  some $k_l=1$. Since the order of the $t_i$ is definable and since the
  $k_i$ are encoded in the labels of the $t_i$, this is clearly
  \MSO-definable. 
  
  Finally, note that condition A\ref{Axiom7} only depends on the
  values of 
  $\ReachguessSym$ on the $t_i$ and on the values of $\exReturns$,
  $\exLoops$ and $\exOneLoops$ for $\LeftStack(t_i, \Encode(c))$. But
  all these information are encoded in the labels of the $t_i$ whence
  condition A\ref{Axiom7} is \MSO-definable

  Thus, we conclude that the validity of a certificate for substack
  reachability is expressible in an \MSO formula on $\mathfrak{A}$. 
  This proves the claim. The lemma follows from the claim as indicated
  above. 
\end{proof}

The following corollary summarises the results obtained so far. 

\begin{corollary}
  Let $\mathcal{S}$ be some collapsible pushdown system. 
  The relation $A$ from Definition \ref{ABCD-Definition} is a
  regular relation via $\Encode$.  
\end{corollary}

Having shown the regularity of $A$, we prove the regularity of $B$,
$C$ and $D$ in the following. We then obtain the regularity of
$\Reach$ as a corollary. 

\paragraph{Regularity of the Relation $\mathbf{B}$}

$B$ contains pairs $(c,\hat c)$ where $c=(q,s)$,
$\hat c=(\hat q,\hat s)$, and $\hat s=\Pop{1}^m(s)$ such that there
is a run from $c$ 
to $\hat c$ not passing any proper substack of $\hat s$. 
By definition, such a run is a composition of high loops and
$\Pop{1}$ or $\Collapse$ of level $1$. 

Recall that we already dealt with a similar problem.
In the previous section we investigated  milestones
$m_1, m_2$ where $m_2=\Pop{1}^n(\Clone{2}(m_1))$ and proved that 
the existence of a run from $m_1$ to $m_2$ with a given initial and
final state is \MSO-definable (cf. Lemma \ref{STACS:CertificateCheck}). 
The following lemma adapts the same idea and proves the regularity of
$B$. 

\begin{figure}[t]
  \centering
  $
  \xymatrix@=1mm{
    &&& g\\
    &&& f\\
    & b & d &e \\
    &a & c & c \\
    s= & \bot & \bot & \bot
  }
  $
  \hskip 2cm
  $
  \xymatrix@=1mm{
    \\
    \\
    & &  \\
    & b & d & e\\
    &a & c & c \\
    \hat s= & \bot & \bot & \bot
  }
  $ 
  \vskip 1cm
  $
  \xymatrix@=7mm{
   && g,\Box:\mathbf{d}\\
   && f,\Box \ar[u]\\
    b,b & d,d &e,e:\hat{\mathbf{d}} \ar[u]\\
   a, a\ar[u] & c,c\ar[r]\ar[u] & \varepsilon, \varepsilon\ar[u] \\
    \bot, \bot \ar[u] \ar[r] & \varepsilon, \varepsilon\ar[u] 
  }
  $
\caption{Illustration of the first case of the proof of lemma
  \ref{Lemma:BReg}. For better orientation, we have marked the
  rightmost leaves of the two encodings by $\mathbf{d}$ and
  $\hat{\mathbf{d}}$.}
  \label{fig:BRegularIllustration1}
\end{figure}

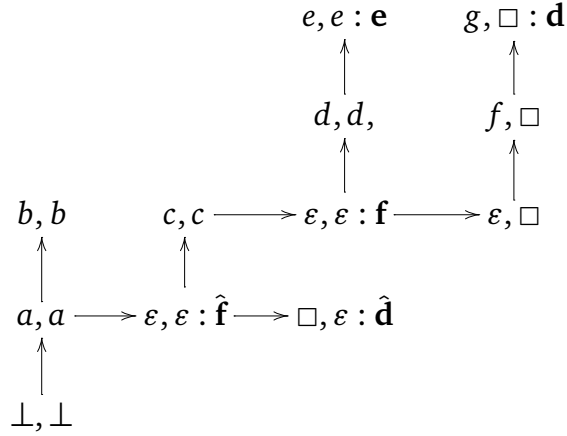
\begin{figure}[t]
  \centering
  $
  \xymatrix@=1mm{
    \\
         &      &      &   e  &  g \\
         &      &      &   d  &  f \\
         &   b  &   c  &   c  &  c \\
         &   a  &   a  &   a  &  a \\
    s= & \bot & \bot & \bot & \bot
  }
  $
  \hskip 2cm
  $
  \xymatrix@=1mm{
    \\
         &      &      &  e     \\
         &      &      &  d     \\
         &  b   &   c  &  c     \\
         &  a   &   a  &  a   &a\\
    \hat s= & \bot & \bot & \bot &\bot
  }
  $
  \vskip 1cm
  $
  \xymatrix@=7mm{
    & & e, e:\mathbf{e} & g, \Box:\mathbf{d}\\
    & & d, d, \ar[u]  & f, \Box \ar[u]\\
    b, b  & c, c \ar[r] & \varepsilon, \varepsilon:\mathbf{f} \ar[u]
    \ar[r] & \varepsilon, \Box \ar[u]\\
    a,a \ar[u] \ar[r]& \varepsilon,
    \varepsilon:\hat{\mathbf{f}} 
    \ar[r] \ar[u] & \Box, \varepsilon:\hat {\mathbf{d}} \\
     \bot, \bot \ar[u] 
  }
  $
\caption{Illustration of the second case of the proof of lemma
  \ref{Lemma:BReg}. For better orientation, 
  we have marked the nodes $\mathbf{d}$, $\hat{\mathbf{d}}$,
  $\mathbf{e}$, $\mathbf{f}$, and $\hat{\mathbf{f}}$.}
  \label{fig:BRegularIllustration2}
\end{figure}

\begin{lemma} \label{Lemma:BReg}
  $B$ is regular via $\Encode$.
\end{lemma}
\begin{proof}
  Let us first recall the structure of
  $\Encode(c)\otimes\Encode(\hat c)$ where $c=(q,s)$, 
  $\hat c=(\hat q, \hat s)$ and $\hat s=\Pop{1}^n(s)$ for
  some $n\in\N$. 
  There are two cases.
  \begin{enumerate}
  \item Let us first assume that $\hat s$ is a milestone of $s$. 
    Figure
    \ref{fig:BRegularIllustration1} shows an
    example of this 
    case. Let $\hat d$
    be the rightmost leaf in $\Encode(\hat c)$. $\domain(\Encode(c))$
    extends $\domain(\Encode(\hat c))$ by the nodes $\hat d0, \hat
    d00, \dots, 
    \hat d0^n$. The labels on the path from $\hat d0$ to the rightmost leaf
    $d:=\hat d0^n$ of $\Encode(c)$ encode the suffix $w$ such that
    $\TOP{2}(s) = \TOP{2}(\hat s)\circ w$. 

    Note that this condition on the domains of $\Encode(c)$ and
    $\Encode(\hat c)$ is $\MSO$-definable whence the pairs
    of configurations of this form are regular via $\Encode$.
  \item Now assume that $\hat s$ is not a milestone of $s$. 
    Figure
    \ref{fig:BRegularIllustration2} shows an
    example of this case.
    Let $d$ be the rightmost leaf of $\Encode(c)$ 
    and  analogously let $\hat d$ be the rightmost leaf of
    $\Encode(\hat c)$. Let $e$ be the second 
    rightmost element in 
    $\Encode(c)$ without left successor. In fact, $e$ is the
    second rightmost leaf of $\hat s$. 
    There are nodes $\hat f < f \leq e$ such that
    \begin{enumerate}
    \item $\hat d = \hat f1$,
    \item $d = f 1 0^m$ for some $m< n$ and
    \item $\LeftStack(e, \Encode(c)) = 
      \LeftStack(e, \Encode(\hat c))$, i.e., the encodings of $s$ and
      $\hat s$ agree on the elements that are lexicographically smaller
      than $e$.  
    \end{enumerate}
    Moreover, the path from $\hat f0$ to $d$ encodes the suffix $v$
    such that $\TOP{2}(s) = \TOP{2}(\hat s)\circ v$. 
    
    Note that these conditions on the domains of $\Encode(c)$ and
    $\Encode(\hat c)$ are $\MSO$-definable whence the pairs
    of configurations of this form are regular via $\Encode$.
  \end{enumerate}
  We conclude that the encodings of pairs of configurations such that
  the stack of the second one is obtained from the first one by a
  sequence of $\Pop{1}$ operations forms a regular set $S$. 

  We show
  that there is an \MSO formula $\psi$ that defines the 
  relation $B$ from Lemma \ref{ABCD-Definition}
  relatively to $S$. 

  We only present the proof for configurations of the second form. The proof
  for the first case is analogous by replacing $\hat f0$ by $\hat d0$.

  Recall that $\hat s=\Pop{1}^n(s)$. 
  There are nodes
  \begin{align*}
    \hat f=: g_0 < g_1 < g_2 < \dots < g_{n-1} < g_n \leq d    
  \end{align*}
  (uniquely
  determined) such that
  $g_i\in\{0,1\}^*0$. These are uniquely determined because
  there are exactly $n$ letters encoded on the path from $\hat f$ to
  $d$ and each left-successor on this path corresponds to one of the
  letters.  

  Let $w_i$ be the topmost word of  $\LeftStack(g_i, \Encode(c))$.
  These form a chain 
  \begin{align*}
    \TOP{2}(\hat c) = w_0 < w_1 < w_2 < w_3 < \dots <
    w_n=\TOP{2}(c)    
  \end{align*}
  where $w_{i+1}$ extends $w_i$ by exactly one
  letter.
  Thus, 
  \begin{align*}
    \exHighLoops(\Pop{1}^m(c)) = 
    \exHighLoops(\LeftStack(g_{n-m}, \Encode(c)))\text{ for all }m\leq n.    
  \end{align*}
  Since $\exHighLoops(\LeftStack(g_{n-m}, \Encode(c)))$ is definable at
  $g_{n-m}$ in $\Encode(c)$, we can \MSO-definably access the
  pairs of initial and final states of all $\Pop{1}^m(c)$. 
  
  Recall that we are looking for a run from $c$ to $\hat c$ that do not visit
  proper substacks of $s$. Such a run consists of a sequence of high
  loops combined with $\Pop{1}$ or $\Collapse$ of level $1$. 
  
  Since the set $\{g_i: 0\leq i \leq n\}$ is \MSO-definable and
  since their order is also \MSO-definable, there is a formula which
  is satisfied by $\Encode(c)\otimes \Encode(\hat c)$ if and only if
  there is a run from $c$ to $\hat c$: given a function $f$ that labels
  each $g_i$ with some state $q_i$, we can check whether there is a loop
  followed by a $\Pop{1}$ or $\Collapse$ of level $1$ 
  from $(q_{n-m}, \Pop{1}^m(s))$ to $(q_{n-m-1},
  \Pop{1}^{m+1}(s))$ such that the following holds:
  \begin{enumerate}
  \item   $q_n=q$ i.e., $q_n$ is the state of $c=(q,s)$ and
  \item  there is  loop from $(q_0, \hat s)$ to $\hat c=(\hat q,\hat s)$. 
  \end{enumerate}
  The function $f$ can be encoded by $\lvert Q \rvert$ many
  sets. Thus, there is a formula asserting that there is a function $f$
  which satisfies the conditions mentioned above. 
  
  For all configurations $(c, \hat c)\in S$, 
  $\Encode(c)\otimes \Encode(\hat c)$ satisfies this formula if and
  only if $(c, \hat c)\in B$, i.e., if there is a run from $c$ to $\hat c$
  that does not visit a proper substack of $\hat c$. 

  Since we have already seen that $S$ is also \MSO-definable, we conclude that
  the Relation $B$ is regular via $\Encode$. 
\end{proof}

\paragraph{Regularity of the Relation $\mathbf{C}$}

Recall that the relation $C$ contains a pair of configurations
$(c,\hat c)$ with $c=(q,s)$ and $\hat c=(\hat q, \hat s)$ 
if and only if the following holds:
$s=\Pop{1}^m(\hat s)$ for some $m\in\N$ and there is a run from
$c$ to $\hat c$ that does not pass a proper substack of $s$. 
We prove the regularity of $C$ analogously to the proof of Lemma
\ref{Lemma:BReg}. 

\begin{lemma}
  The relation $C$ from Definition \ref{ABCD-Definition} is a regular
  relation via $\Encode$. 
\end{lemma}
\begin{proof}
  We proceed completely analogous to Lemma \ref{Lemma:BReg}. 

  Let us first recall the structure of
  $\Encode(c)\otimes\Encode(\hat c)$ where $c=(q, s)$, 
  $\hat c=(\hat q, \hat s)$ and $\hat s$ can be generated from $s$
  by a sequence of 
  $\Push{\sigma,i}$ of length $m\in\N$.  
  There are two cases.
  \begin{enumerate}
  \item Let us first assume that $s$ is a milestone of $\hat s$. 
    Let $d$
    be the rightmost leaf in $\Encode(c)$. $\domain(\Encode(\hat c))$
    extends $\domain(\Encode(c))$ by the nodes $d0, d00, \dots,
    d0^m$. The labels on the path from $d0$ to
    the rightmost leaf 
    $\hat d:=d0^m$ of $\Encode(\hat c)$ encode the suffix $w$ such that
    $\TOP{2}(\hat s) = \TOP{2}(s)\circ w$. 

    Note that this condition on the domains of $\Encode(c)$ and
    $\Encode(\hat c)$ is $\MSO$-definable whence the pairs
    of configurations of this form are regular via $\Encode$.
  \item Now assume that $s$ is not a milestone of $\hat s$. 
    Let $d$ be the rightmost leaf of $\Encode(c)$ and $\hat d$ be the
    rightmost leaf of $\Encode(\hat c)$. Let $e$ be the second
    rightmost element in 
    $\Encode(\hat c)$ without left successor. In fact, $e$ is the
    second rightmost leaf of $\Encode(c)$. 
    There are nodes $f < \hat f \leq e$ such that
    \begin{enumerate}
    \item $d = f1$,
    \item $\hat d = \hat f 1 0^m$ for some $m< n$, and
    \item $\LeftStack(e, \Encode(c)) = 
      \LeftStack(e, \Encode(\hat c))$, i.e., the encodings of $s$ and
      $\hat s$ agree on the elements that are  lexicographically less
      or equal to $e$. 
    \end{enumerate}
    Moreover, the path from $f0$ to $\hat d$ encodes the suffix $w$
    such that $\TOP{2}(\hat s) = \TOP{2}(s)\circ w$. 
    Note that these conditions are similar to those in the 
    proof of lemma \ref{Lemma:BReg} with exchanged roles for $s$ and
    $\hat s$. 
    
    But in this proof there is one further condition on the encodings
    of $c$ and $\hat c$.
    The path
    from $f0$ to $\hat f$ may only encode letters with link
    level $1$. This stems from the following fact.

    Since $\TOP{2}(s)$ is a prefix of $\TOP{2}(\hat s)$ and $s$ is no
    milestone of $\hat s$, $\TOP{2}(s)$ is a proper prefix of the
    greatest common prefix of the two topmost words of $\hat s$, i.e.,
    \begin{align*}
      \TOP{2}(s) \leq \TOP{2}(\hat s) \sqcap \TOP{2}(\Pop{2}(\hat s)).  
    \end{align*}
    Furthermore, $\hat f$ is defined in such a way that
    $\LeftStack(\hat f1, \Encode(\hat c))$ is the minimal milestone of
    $\hat s$
    that has width
    $\lvert \hat s \rvert$. Thus, the elements encoded along the path
    from $f0$ to $\hat f$ are also contained in the second topmost word
    of $\hat s$. Thus, if any of these is of link level $2$, then it
    points strictly below $\Pop{2}(\hat s)$.  But such a link
    cannot be constructed from $s$ by application of
    push operations because a push applied to a stack of width $\lvert
    \hat s \rvert$ cannot generate an
    element that points below $\Pop{2}(\hat s)$. 

    Note that these conditions  are  $\MSO$-definable on
    $\Encode(c)\otimes\Encode(\hat c)$. 
  \end{enumerate}
  Thus, the pairs of
  configurations $(c, \hat c)$ where the stack of $\hat c$ can be
  generated from the stack of $c$ by a sequence of push operations 
  form a regular set via $\Encode$. Let $S$ denote the set of these
  pairs of configurations. 
  Next, we show
  that there is an \MSO formula $\psi$ defining the  relation $C$ from
  Definition \ref{ABCD-Definition}
  with respect to $S$.

  We only present the proof for configurations of the second form. The proof
  for the first case is analogous, just replace $f0$ by $d0$.

  Recall that $n=\lvert \TOP{2}(\hat s) \rvert - \lvert \TOP{2}(s)\rvert$.
  By definition of the encoding, there are nodes
  \begin{align*}
    f0 = g_0 < g_1 < g_2 < \dots < g_{n-1} < g_n \leq \hat d    
  \end{align*}
  (uniquely
  determined) such that
  $g_i\in\{0,1\}^*0$ for all $1\leq i \leq n$. 

  Let $w_i$ be the topmost word of  $\LeftStack(g_i, \Encode(\hat c))$.
  We obtain a chain 
  \begin{align*}
    \TOP{2}(s) =: w_0 < w_1 < w_2 < w_3 < \dots <
    w_n=\TOP{2}(\hat s)    
  \end{align*}
  where $w_{i+1}$ extends $w_i$ by exactly one
  letter.
  Thus, 
  \begin{align*}
    \exHighLoops(\Pop{1}^m(\hat s)) = 
    \exHighLoops(\LeftStack(g_{n-m}, \Encode(\hat c)))\text{ for all
    }0\leq m\leq n.     
  \end{align*}
  Since $\exHighLoops(\LeftStack(g_{n-m}, \Encode(\hat c)))$ is
  \MSO-definable at 
  $g_{n-m}$ in $\Encode(\hat c)$, we can \MSO-definably access the
  pairs of initial and final states of high loops of all
  $\Pop{1}^m(\hat s)$.  
  
  Since we are looking for runs from $c$ to $\hat c$ that do not visit
  a proper substack of the stack of $c$, these consist of a sequence
  of high loops and 
  push-operations.  
  
  Since the set of the $g_i$, $0\leq i \leq n$, is \MSO-definable and
  since their order is also \MSO-definable, there is a formula which
  is satisfied by $\Encode(c)\otimes \Encode(\hat c)$ if and only if
  there is such a run from $c$ to $\hat c$ not passing 
  a substack of $\Pop{2}(c)$: given a function $f$ that labels
  $g_i$ with a state $q_i$, a formula can assert that there is a high loop
  followed by a push operation   
  from $(q_{n-m}, \Pop{1}^m(\hat c))$ to $(q_{n-(m-1)},
  \Pop{1}^{m-1}(c))$ for each $0\leq m < n$ such that the following holds:
  \begin{enumerate}
  \item   $q_0=q$, i.e., $q_0$ is the state of $c$ and
  \item  there is a high loop from $(q_n, \hat s)$ to 
    $\hat c=(\hat q, \hat s)$. 
  \end{enumerate}
  The function $f$ can be encoded by $\lvert Q \rvert$ many
  sets. Thus, there is a formula $\psi$ which asserts that there is a
  function $f$ 
  which satisfies the conditions mentioned above. 
  
  For all configurations $c, \hat c$ such that the stack of $\hat c$
  can be created 
  from $c$ by a sequence of push transitions, 
  $\Encode(c)\otimes \Encode(\hat c)$ satisfies $\psi$ if and
  only if there is a run from $c$ to $\hat c$ not passing a proper
  substack of $c$. 

  Since we have already seen that the set $S$ of pairs $(c, \hat c)$
  such that $\hat c$ can be created from $c$ by a sequence of
  push transitions is also \MSO-definable, we conclude that
  the Relation $C$ is regular via $\Encode$.   
\end{proof}

\paragraph{Regularity of the Relation $\mathbf{D}$}

Recall that the relation $D$ contains a pair $(c,\hat c)$ of
configurations for
$c=(q,s)$ and \mbox{$\hat c=(\hat q,\hat s)$} if and only if the following
holds: $s=\Pop{2}^m(\hat s)$ and there 
is a run from $c$ to $\hat c$ not passing any substack of $s$ after
its initial configuration. 

Recall that $s=\Pop{2}^n(\hat s)$ implies that $s$ is a milestone of
$\hat s$. Hence, the existence of a 
run from $c$ to $\hat c$ can be checked in a similar manner as the
existence of a 
run from the initial configuration to $\hat c$. Any run of the latter
form passes $s$. If it passes $s$ in state $q$ this is a witness for
$(c,\hat c)\in D$.

\begin{lemma}
  $D$ is regular via $\Encode$. 
\end{lemma}
\begin{proof}
  Fix a collapsible pushdown system $\mathcal{S}$.
  The set
  \begin{align*}
    S := \{\Encode(c)\otimes\Encode(\hat c):
    c=(q,s), \hat c=(\hat 
    q,\hat s) \in\CPG(\mathcal{S})
    \text{ and } s=\Pop{2}^n(\hat s)\text{ for some }n\in\N\}
  \end{align*}
  is regular (cf. the proof of \ref{MSODefinabilityofA}). 

  Let $c, \hat c$ be configurations such that
  $(c,\hat c)\in S$. We write $c=(q,s)$ and
  $\hat c=(\hat q,\hat s)$. Furthermore, since
  $c\in\CPG(\mathcal{S})$, 
  there is a run $\rho_0$ from $(q_0, \bot)$ to $c$.
  
  There is a run $\rho$ from $c$ to $\hat c$ if and only if there is a run 
  $\hat \rho:=\rho_0\circ\rho$
  from $(q_0, \bot_2)$ to $\hat c$ passing $c$. 
  
  From the previous section,
  we know that there is a certificate for reachability $C_{\hat\rho}$
  induced by $\hat\rho$ which labels each  node $e\in\Encode(\hat c)$ by the
  last state in which $\hat\rho$ passes $\LeftStack(e, \Encode(\hat c))$. 
  Since $s$ is a milestone of $\hat s$, the rightmost leaf $d$ of
  $\Encode(c)$ is a 
  node in $\Encode(\hat c)$ such that $s=\LeftStack(d,\Encode(\hat c))$. 

  If $\rho$ does not visit any substack of $c$ after its initial
  configuration, then $C_{\hat\rho}(d)=q$. 
  On the other hand, if $C_{\hat\rho}(d)=q$ then there is a run
  $\hat\rho$ and some $i\in\domain(\hat\rho)$ such that
  $\hat\rho(i)=(q,s)=c$, $\hat\rho$ ends in $\hat c$ and after $i$ no
  substack of $s$ is visited. Thus,
  $\hat\rho{\restriction}_{[i+1,\length(\hat\rho)]}$ witnesses
  $(c,\hat c)\in D$.

  Since $d$ is $\MSO$-definable, there is a formula $\psi$ such that
  $\Encode(c)\otimes\Encode(\hat c)\models\psi$ for some $(c,\hat c)\in
  S$, if
  and only if there is a certificate $C_{\hat\rho}$ for $\hat c$ such that
  $C_{\hat\rho}(d) = q$. This means that $\psi$ defines $D$ relatively
  to $S$.
  Since $S$ is regular, we conclude that
  $D$ is also regular. 
 \end{proof}

\paragraph{Regularity of $\mathbf{Reach}$}

As already indicated, the regularity of $A$, $B$, $C$, and $D$
directly implies the regularity of $\Reach$. 
We obtain the following corollary. 
\begin{corollary}
  Let $\mathcal{S}$ be a collapsible pushdown system  of  level $2$. 
  The expansion of the graph of $\mathcal{S}$ by the reachability
  predicate is automatic, i.e., the graph $(\CPG(\mathcal{S}), \Reach)$ is 
  automatic. Thus, 
  the $\FO{}(\Reach)$-theory of $\CPG(\mathcal{S})$ is decidable.
\end{corollary}

\paragraph{Regularity of $\mathbf{Reach_L}$}

In the previous section, we proved that the reachability predicate
$\Reach$ on collapsible pushdown graphs is automatic via $\Encode$.
We improve this result and show that reachability by a path
that satisfies a regular expression is automatic. 

Recall that for $L\subseteq\Gamma^*$ some regular language, $\Reach_L$
is the binary relation that contains configurations $(c,\hat c)$ if and
only if there is a run $\rho$ from $c$ to $\hat c$ such that the labels of
the transitions used in $\rho$ form a word $w$ such that $w\in L$. 

Let $L$ be some regular language. We show that $\Reach_L$ is automatic
via $\Encode$ by
constructing a version of the 
product of $\mathcal{S}$ with the automaton $\mathcal{A}_L$
corresponding to $L$. We show that $\CPG(\mathcal{S})$ is first-order
interpretable in this 
product and we show that the predicate $\Reach_L$ on
$\CPG(\mathcal{S})$ can be expressed via  $\Reach$ on this product. 

Before we state the lemma, we introduce some abbreviations. For $x$ a
variable and $q$ a state of some collapsible pushdown system
$\mathcal{S}$, we write $x\in q$ for the $\FO{}$ 
formula stating that $x$ is a configuration with state $q$. This is
definable because we assume that the label of an incoming transition
encodes the state of the node. Furthermore, all configurations but the
initial one have at least one incoming edge. 
Since the set $Q$ of states is finite, we also write $x\in Q'$ where
$Q'\subseteq Q$ for the formula $\bigvee_{q\in Q'} x\in q$.

\begin{lemma} \label{ReachLRegular}
  Let $\mathcal{S}=(Q,\Sigma, \Gamma ,q_i,\Delta)$ be a
  $2$-\CPS. Furthermore, let
  $L_1, L_2, \dots, 
  L_n \subseteq\Gamma^*$ be regular languages. Then $(\CPG(\mathcal{S}),
  (\Reach_{L_i})_{1\leq i \leq n})$ is automatic via $\Encode$.
\end{lemma}
\begin{proof}
  Without loss of generality, we assume that $n=1$ and write $L$ for
  $L_1$. The general case 
  is proved by iterating the following construction. 
  Let $\mathcal{A}_L=(F,\Gamma,f_i,f_f, \Delta_L)$ be the finite
  string-automaton corresponding to $L$. 
  
  We define the product of $\mathcal{S}$ and $\mathcal{A}_L$ to be the
  collapsible pushdown system 
  \begin{align*}
    \mathcal{\bar S}=(\bar
    Q, \Sigma, \Gamma\cup\{\varepsilon_i, \varepsilon_f\}, q_i, \bar
    \Delta)    
    \text{  where}
  \end{align*}
  \begin{itemize}
  \item $\bar Q\coloneqq  Q \cup (Q\times F)$ and
  \item $\bar \Delta$ is the union 
    \begin{align*}
      &\Delta \\
      &\cup \{(q,\sigma, \varepsilon_i, (q,f_i),\Id): \sigma\in\Sigma,
      q\in Q\} \\ 
      &\cup \{(q, \sigma, \varepsilon_f, (q,f_f), \Id):
      \sigma\in\Sigma, 
      q\in Q\} \\ 
      &\cup \{((q,f),\sigma,\gamma, (q',f'), \op):
      (q,\sigma,\gamma,q',\op) \in \Delta 
      \text{ and } (f, \gamma, f')\in \Delta_L\}.      
    \end{align*}
  \end{itemize}

  Note that $\CPG(\mathcal{S})$ is \FO{} definable in 
  $\CPG(\mathcal{\bar S})$:
  both graphs 
  have the same initial configuration and $\mathcal{\bar S}$ extends
  $\mathcal{S}$ only
  by transitions that lead to configurations with states in $Q\times
  F$. Hence, the restriction of $\CPG(\mathcal{\bar S})$ to the set
  $\{x\in\CPG(\mathcal{\bar S}): x\in Q\}$ is isomorphic to
  $\CPG(\mathcal{S})$.

  On the other hand, by construction there is a path
  from $((q,f_i), s)$ to $((q',f_f), s')$ in $\CPG(\mathcal{\bar S})$
  if and 
  only if there is a path from $(q,s)$ to $(q',s')$ in
  $\CPG(\mathcal{\bar S})$ 
  whose path corresponds to an accepting word of $\mathcal{A}_L$. 
  Hence, $(x,y)\in \Reach_{L}$ on $\CPG(\mathcal{S})$ corresponds to
  \begin{align*}
   \exists x' \exists y' \big ( x \trans{\varepsilon_i} x' \wedge
    y \trans{\varepsilon_f} y' \wedge
  (x',y')\in \Reach\big)
  \end{align*}
  on $\CPG(\mathcal{ \bar S})$. The closure of automaticity under
  first-order interpretations yields the desired result. 
\end{proof}

\subsection{Combination of \FO{} and  $L\mu$ Model Checking}
\label{Lmu-FO-CPG}

We have obtained an $\FO{}$ model checking algorithm for
collapsible pushdown graphs of level two. Recall that Hague et
al.\ \cite{Hague2008} have shown that there is an  $L\mu$ model
checking algorithm for the class of all collapsible pushdown graphs. 
It is a natural question whether these two results can be combined. In
order to give an answer to this question, we investigate the
following three questions. 
\begin{enumerate}
\item Let $\mathcal{C}_{L\mu}$ be the class of graphs obtained by
  $L\mu$-interpretation from the class of level $2$ collapsible
  pushdown  graphs. 
  Is the $\FO{}$ model checking problem on $\mathcal{C}_{L\mu}$
  decidable?
\item  Let $\mathcal{C}_{\FO{}}$ be the class of graphs obtained by
  $\FO{}$-interpretation from the class of level $2$ collapsible
  pushdown graphs. Is the $L\mu$ model checking problem on
  $\mathcal{C}_{\FO{}}$ decidable?
\item 
  Is MLFP\footnote{
    Monadic least fixpoint logic (MLFP) is the smallest logic
    encompassing the expressive power of 
    $L\mu$ and $\FO{}$ that has sensible closure properties. 
  } model checking decidable on level $2$ collapsible pushdown graphs?
\end{enumerate}

Due to a recent result of Broadbent et al.\ \cite{BroadbentCOS10}, the
first question can be answered 
positively. They proved the following result: 
\begin{theorem}[\cite{BroadbentCOS10}] The global $L\mu$ model
  checking for collapsible pushdown graphs is decidable.\footnote{The
    global model checking problem asks the following: given a formula
    $\varphi\in L\mu$ and a graph $\mathfrak{G}$, what are the nodes of
    $\mathfrak{G}$ that satisfy $\varphi$, i.e., what is the set
    $\{g\in\mathfrak{G}: \mathfrak{G}, g \models \varphi\}$?}
\end{theorem}
For the proof of this theorem, Broadbent et al.\ introduced an encoding
of a collapsible pushdown stack as a word with back-edges. A word with
back-edges looks similar to a nested word, but the back-edges are not
well nested. Via this encoding, $L\mu$ definable sets of stacks are
turned into sets of words with back-edges that are recognised by
deterministic automata on such words with back-edges. These
automata work like finite automata on ordinary words, but they
propagate the state at a position in the word to the next position and
to those positions reachable via a back-edge. Broadbent et al.\ provide
a construction of an 
automaton on words with back-edges that corresponds to a given formula
$\varphi\in L\mu$. Using this construction one can then decide the global
model checking problem. 

For collapsible pushdown graphs of level two, their techniques imply that
$L\mu$-definable subsets  are automatic via $\Encode$ as follows. 
Let $S$ be a set of level $2$ stacks. If the set of words with back-edges
encoding $S$ is regular in the sense of Broadbent et al., then $S$ is
automatic via $\Encode$. 

Hence, from the results of Broadbent et al.\ the next corollary follows
immediately.

\begin{corollary}[\cite{BroadbentCOS10}]
  The $L\mu$-definable subsets in a $2$-\CPG are transformed into 
  regular sets of trees by the encoding function $\Encode$. 
\end{corollary}

From the decidability of model checking on automatic structures, 
we directly conclude that first-order logic on collapsible pushdown
graphs expanded by $L\mu$-definable predicates is decidable. 

\begin{corollary}
  The graph of a collapsible 
  pushdown system of level $2$ enriched by $L\mu$-definable predicates is
  automatic. Hence, 
  its \FO{}(Reg, $(\RamQ{n})_{n\in\N},
  (\exists^{k,m})_{k,m\in\N}$)-theory is decidable.  
\end{corollary}
Of course, this result is compatible with Lemma
\ref{ReachLRegular}. Thus, Theorem \ref{Thm:CPGRegular} follows
directly from these two results.

After the positive answer to our first question, we give
negative answers to the other two questions. 
We show the undecidability of the $L\mu$ model
checking  on graphs obtained by first-order interpretations from
collapsible pushdown graphs of level $2$. 

This negative answer to our second question implies also a negative
answer to the third: MLFP encompasses $\FO{}$ and $L\mu$ whence for
each $\FO{}$-interpretation $I$ and each $L\mu$-formula $\varphi$
there is a MLFP formula $\psi$ such that 
$\mathfrak{A}  \models \psi$ if and only if
$I_{\mathrm{Str}}(\mathfrak{A}) \models \varphi$ for all structures
$\mathfrak{A}$.  
Since we show the undecidability of the second problem, 
the first one is also undecidable.

Recall that Lemma \ref{Lmu-Half-grid-undecidability} shows the
undecidability of the $L\mu$ model checking on the bidirectional
half-grid (recall Figure \ref{HALFGRID}). 
Due to this result, the following lemma implies that $L\mu$ model
checking is
undecidable on $\mathcal{C}_{\FO{}}$.

\begin{lemma} \label{InterpretationBiGrid}
  The bidirectional half-grid $\BiHalfgrid$ is
  \FO{} interpretable in a certain \CPG of level $2$.
\end{lemma}
\begin{proof}
  Extending the idea for the MSO-undecidability result of Hague et
  al.\ \cite{Hague2008}, we consider the 
  following collapsible pushdown graph.

  Let $Q:=\{0,1,2\}, \Sigma:= \{\bot,a\}$, and
  $\Delta$ is given by 
  $(0,-,\mathrm{Cl}_1, 1,\Clone{2})$,
  $(1,-,\mathrm{A},0,\Push{a,2})$, 
  $(0,-,\mathrm{Cl}_2, 2,\Clone{2})$,
  $(2,a,\mathrm{P}_1,2,\Pop{1})$, 
  $(2,a,\mathrm{Co},0,\Collapse)$, and 
  $(2,a,\mathrm{P}_2,0,\Pop{2})$ where
  ``$-$'' denotes any
  letter from $\Sigma$. 
  We call this example
  graph $\mathfrak{G}$ (cf.  Figure \ref{fig:CPGExampleII}).\\
  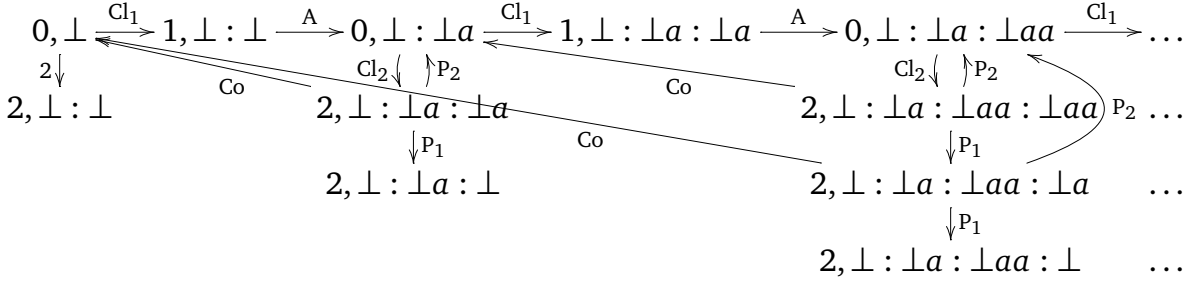
\begin{figure}[h]
    \centering
    $
    \begin{xy}
      \xymatrix@R=12pt@C=12pt{
        0, \bot \ar[r]^-{\mathrm{Cl}_1} \ar[d]_{ 2}& 
          1, \bot:\bot \ar[r]^-{\mathrm{A}} & 
          0, \bot:\bot a \ar[r]^-{\mathrm{Cl}_1} \ar@/_/[d]_{\mathrm{Cl}_2}&
          1, \bot:\bot a:\bot a  \ar[r]^-{\mathrm{A}} & 
          0, \bot:\bot a: \bot aa \ar[r]^-{\mathrm{Cl}_1}
          \ar@/_/[d]_{\mathrm{Cl}_2} & \dots \\
        %
        %
          2, \bot:\bot&  
          & 
          2, \bot:\bot a:\bot a
          \ar[d]^{\mathrm{P}_1}\ar@/_/[u]_{\mathrm{P}_2}
          \ar[ull]^{\mathrm{Co}}&   
          & 
          2, \bot:\bot a : \bot aa:\bot aa \ar[d]^{\mathrm{P}_1}
          \ar@/_/[u]_{\mathrm{P}_2} 
            \ar[ull]^{\mathrm{Co}}
          & \dots \\
        %
        %
          &
          & 
          2, \bot:\bot a : \bot  &  
          &
          2, \bot:\bot a:\bot aa: \bot a  \ar[d]^{\mathrm{P}_1}
          \ar@/_4.8pc/[uu]_{\mathrm{P}_2} 
          \ar[uullll]^(.4){\mathrm{Co}}&  
          \dots \\
        %
        %
          &&&&
          2, \bot:\bot a: \bot a a : \bot &
          \dots
      }
    \end{xy}
    $
    \caption{The collapsible pushdown graph $\mathfrak{G}$.}
    \label{fig:CPGExampleII}
  \end{figure}
  
  In order to interpret $\BiHalfgrid=(H,\rightarrow,\downarrow, \leftarrow,
  \uparrow)$ in $\mathfrak{G}$, we first  have to define the 
  domain of this interpretation. Let 
  \begin{align*}
    \varphi(x):= \exists y\ x\vdash^{\mathrm{P}_2} y.    
  \end{align*}
  This formula defines all elements that
  are not in the first row 
  of $\mathfrak{G}$ and
  which have a $\Pop{2}$ and a $\Collapse$ successor. 
  Set 
  \begin{align*}
    \varphi_{nd}(x,y):=\exists z \exists z' \left(x\vdash^{\mathrm{Co}} z \wedge
    y\vdash^{\mathrm{P}_1} z' \wedge z' \vdash^{\mathrm{Co}} z\right).    
  \end{align*}
  $\varphi_{nd}$
  defines the relation ``$y$ is on the diagonal to the right of the diagonal
  of $x$''.  Set 
  \begin{align*}
    \varphi_{nc}(x,y):=\exists z \exists z' \exists z'' \left(z
      \vdash^{\mathrm{Cl}_1} 
    z' \vdash^{\mathrm{A}} z'' \wedge x \vdash^{\mathrm{P}_2} z \wedge
    y \vdash^{\mathrm{P}_2} z''\right).    
  \end{align*}
  This formula defines the relation 
  ``$y$ is on the column to the right of the column of $x$''.

  Now, $y$ is the right neighbour of $x$ if and only if 
  \begin{align*}
    \varphi_{\rightarrow}:=\varphi_{nd}(x,y)\wedge \varphi_{nc}(x,y)  
  \end{align*}
  holds. 

  Hence, the $\FO{}$-interpretation  $I:=(\varphi, \varphi_{\rightarrow},
  \vdash^{\mathrm{P}_1}, \varphi_{\rightarrow}^{-1},
  (\vdash^{\mathrm{P}_1})^{-1})$ yields
  $\BiHalfgrid=\mathrm{Str}_I(\mathfrak{G})$.
\end{proof}

In the following, we summarise observations concerning the
optimality of our result. 
A first question is whether the
complexity of the $\FO{}$ model checking algorithm
may be
improved. As we mentioned in Chapter \ref{Chapter_InfiniteStructures},
using automatic representations for model checking purposes leads
to nonelementary complexity of the model checking algorithm. In the
first part of this section we show a matching lower bound: 
any $\FO{}$ model checking algorithm has nonelementary complexity. 
Then we briefly discuss a negative result concerning model checking on 
higher levels of the collapsible pushdown hierarchy:
Broadbent \cite{Broadbent2010Mail} has shown that
\FO{} is undecidable on the third level of the collapsible
pushdown hierarchy. 

\subsection{Lower Bound for FO Model Checking}
\label{sec:LowerBounds}

Recall that $\FO{}$ model checking on automatic structures has
nonelementary complexity. 
In the following theorem we show that 
there is no elementary algorithm for \FO{} model checking
on collapsible pushdown graphs. 
The proof is by reduction to the
nonemptyness problem for star-free regular expressions. 
As an auxiliary step, we prove that \FO{} model checking on the
full infinite binary tree is nonelementary.

\begin{lemma} \label{LemmaLowerBoundCPGModelChecking}
  The expression complexity of $\FO{}$ model checking on the full
  infinite binary tree $\mathfrak{T}:=(T, \prec, S_1, S_2)$ with
  prefix order $\prec$ and successor relations $S_1, S_2$ is nonelementary.  
\end{lemma}
\begin{proof}
  For each first-order sentence $\varphi$, there is a first-order sentence
  $\varphi'(x)$ 
  such that for all $t\in\mathfrak{T}$, 
  $\mathfrak{T}\models \varphi'(t)$
  if and only if $\mathfrak{T}{\restriction}_{\{ t': t'\prec t \}}
  \models \varphi$.
  Note that $\mathfrak{T}{\restriction}_{\{ t': t'\prec t \}}$ can be
  considered as a finite word structure over the alphabet $\{1,2\}$:
  We identify an incoming $S_1$ edge with the label $1$ and an 
  incoming $S_2$ edge with the label $2$.

  In this sense, the model checking problem for the formula $\exists x
  \varphi'(x)$ on $\mathfrak{T}$ is equivalent to the satisfiability
  problem for 
  $\varphi$ with respect to  the class of word-structures. Via the
  classical result of McNaughton and Papert \cite{McNaughtonP71}
  this problem is equivalent to the nonemptyness
  problem for languages defined by star-free regular expressions. 
  Since the latter problem has nonelementary
  complexity \cite{phd-stockmeyer}, the claim follows. 
\end{proof}

Now, we present a reduction of the \FO{} model checking on the
full infinite binary tree to the \FO{} model checking on
collapsible pushdown graphs. 

\begin{theorem}\label{thm:FOCPGnonelementary}
  The expression complexity of  any \FO{} model checking algorithm for
  level $2$ collapsible pushdown graphs is  nonelementary.
\end{theorem}
\begin{proof}
  For the proof of this theorem we modify the graph of 
  example \ref{fig:CPGExampleII}. Note that $(\omega, <)$ is first order
  definable in this graph: restrict the domain to all elements with
  state $0$. The order $<$ is then defined via $\varphi_<(x,y):=\exists z\
  z\trans{\Pop{2}} y \land z\trans{\Collapse} x$. 
  
  In order to obtain a binary tree from a collapsible pushdown graph
  we create an infinite tree-like graph where every branch is a
  copy of the graph from example \ref{fig:CPGExampleII}. The copies are
  ordered in such a way that the first-order interpretation from above
  yields the full binary tree when applied to this graph.
  
  To this end, we duplicate the letter $a$ and the label $A$. 
  We introduce a new letter
  $a'$ and a new label $A'$. Furthermore, 
  for each transition where $a$ occurs, we add the corresponding
  transition where $a$ is replaced by $a'$ as follows:
  we add the
  transitions
  $(1, -, A', 0, \Push{a',2})$,
  $(2, a', P_1, 2, \Pop{1})$, 
  $(2, a', P_2, 0, \Pop{2})$, and
  $(2, a', \mathrm{Co}, 0, \Collapse)$  
  where $A'$ is a new edge-label. 

  On the resulting graph restricted to the
  configurations with states $0$, the formula $\varphi_<(x,y)$
  from above defines the
  prefix order of the full infinite binary tree. Furthermore, the
  formulas 
  $\varphi_L(x,y):=\exists z\ x\trans{\mathrm{Cl}} z\trans{A} y$
  and  
  $\varphi_R(x,y):=\exists z\ x\trans{\mathrm{Cl}} z\trans{A'} y$
  define the left successor, respectively, the right successor
  relation.

  Lemma \ref{LemmaLowerBoundCPGModelChecking}
  implies the desired result.  
\end{proof}

\subsection{Model Checking  on Higher-Order Collapsible Pushdown Graphs} 
\label{FOReachUndecidableonFourCPG}

Recently, Broadbent \cite{Broadbent2010Mail} developed a reduction of
Post's
correspondence problem (PCP,cf. \cite{Post46}) to
\FO{} model checking on collapsible pushdown graphs of level $3$. 
Since the PCP is undecidable, it follows that the \FO{} model
checking problem on level $3$ collapsible pushdown graphs is
undecidable.

In fact, Broadbent's proof comes in two variants:
firstly, there is a fixed level $3$ collapsible
pushdown graph with undecidable $\FO{}$-theory. On this fixed graph,
there is a 
first-order formula
for each instance of the PCP with the following property. The graph
satisfies this 
formula if and 
only if the corresponding  instance of the PCP has a solution. 
Secondly, Broadbent provides a fixed formula $\varphi\in\FO{}$ such
that there is a class $\mathcal{C}$ of level $3$ collapsible pushdown
graphs such that the following holds. For each instance of the PCP
there is a graph $\mathfrak{G}\in\mathcal{C}$ such that
$\mathfrak{G}\models\varphi$ if and only if this instance of the PCP
has a solution. 

Thus, 
\FO{} model checking on level $3$ collapsible pushdown graphs is
undecidable even for either fixed structure or fixed formula.


\section{An FO Model Checking Algorithm on Nested Pushdown Trees}
\label{Chapter_FO-NPT}

This section analyses the \FO{} model checking problem on the class of
nested pushdown trees. In the first part we reduce the 
$\FO{}(\Reach)$ model checking problem to the $\FO{}(\mathrm{Reg})$ model
checking problem for level $2$ collapsible pushdown automata. We show
that there is a first-order interpretation $I$ such that for each
nested pushdown tree $\mathfrak{N}$ there is a
collapsible pushdown 
graph $\mathfrak{G}$ of level $2$ such that
$\mathrm{Str}_I(\mathfrak{G}) = \mathfrak{N}$. Furthermore, 
$I$ transfers the reachability predicate on 
nested pushdown trees into a certain regular
reachability predicate on the collapsible pushdown graph. 

In Sections \ref{SecGaifman}--\ref{SectionNPTModelChecking} we have a
closer look at the complexity of $\FO{}$ 
model checking on nested pushdown trees. 
We develop several versions of the pumping lemma for pushdown systems
which are
compatible with the jump edges in the following sense: application of
these lemmas to a run yields a short run with equivalent first-order
type. The bounds obtained by this lemma can be used as a constraint
for Duplicator's strategy in the Ehrenfeucht-\Fraisse game on two
identical copies of some nested pushdown tree. 
As indicated in Section \ref{Sec:EFGame}, this result can be turned
into a model checking algorithm. Using this approach, we show that the
complexity of \FO{} model  checking on nested pushdown trees is in
$2$-$\mathrm{EXPSPACE}$. 

\subsection{Interpretation of NPT in CPG}
\label{sectionNPTSimulation}

In this section, we show that any nested pushdown tree can be
first-order interpreted in some collapsible pushdown graph of level
$2$. For this purpose, we fix a pushdown system
$\mathcal{N}=( Q, \Sigma, \Gamma, \Delta, q_0)$. 
We show that there is a
collapsible pushdown system of level $2$ and a first-order
interpretation that yields the nested pushdown tree generated by
the pushdown system.

The basic idea is the following: every vertex of the nested pushdown
tree generated by $\mathcal{N}$ is a run,
i.e., a list of configurations that are passed by this run. Every
configuration is a 
level $1$-stack $s$ and a state $q$. We write the state $q$ on top of
the stack $s$ and obtain the stack $\Push{q}(s)$. Then we represent a run
$(q_1,s_1)\trans{}(q_2,s_2)\trans{} \dots\trans{} (q_n,s_n)$ by the stack
$\Push{q_1}(s_1):\Push{q_2}(s_2):\dots:\Push{q_n}(s_n)$. Using this
encoding, we can
simulate every transition of the pushdown system by at most four stack
operations of the collapsible pushdown system and the nesting edges can be simulated by
reverse collapse edges. 


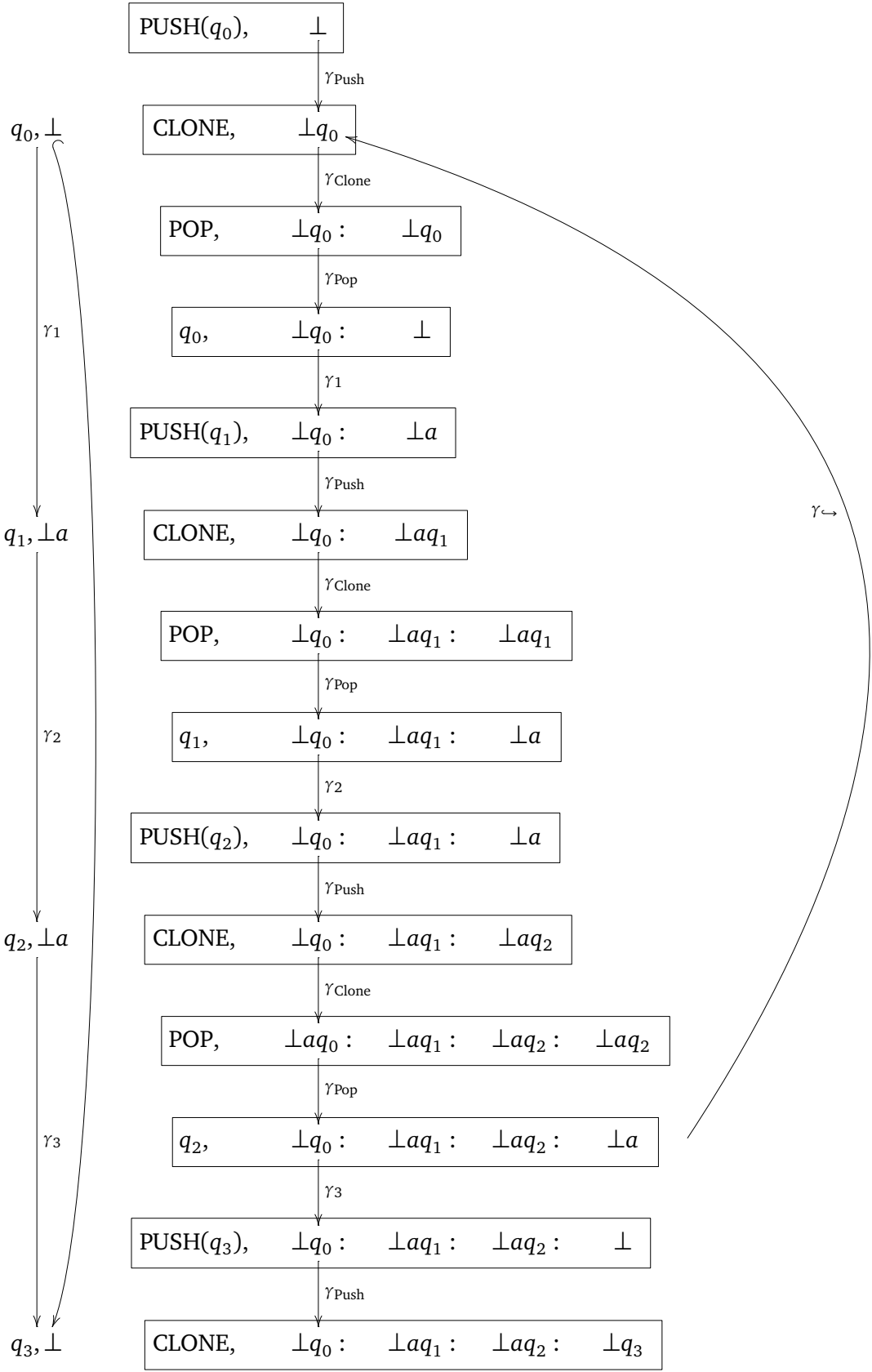
\begin{figure}  
  \begin{xy}{\small
    \xymatrix@C=10pt@R=30pt{
      & & \mathrm{PUSH}(q_0),& \bot 
      \ar[d]^{\gamma_{\mathrm{Push}}}\\
      q_0,\bot \ar[dddd]^{\gamma_1}
      \ar@{^{(}->}@(r,r)[dddddddddddd]   
      & & \mathrm{CLONE},& \bot  q_0
      \ar[d]^{\gamma_{\mathrm{Clone}}}\\ 
      & & \mathrm{POP},& \bot  q_0:
      \ar[d]^{\gamma_{\mathrm{Pop}}}
      & \bot q_0
       \\
      & & q_0,& \bot  q_0 : \ar[d]^{\gamma_1}
      & \bot 
      \\
      & & \mathrm{PUSH}(q_1),& \bot  q_0
      :\ar[d]^{\gamma_{\mathrm{Push}}} & \bot a \\
      q_1, \bot a \ar[dddd]^{\gamma_2}& &  
      \mathrm{CLONE}, & \bot q_0 : \ar[d]^{\gamma_{\mathrm{Clone}}} &
      \bot a q_1 \\
      & & \mathrm{POP},& \bot q_0 :
      \ar[d]^{\gamma_{\mathrm{Pop}}} &
      \bot a q_1 : & \bot a q_1 
       \\
      & & q_1,& \bot q_0 : \ar[d]^{\gamma_2} &
      \bot a q_1 : & 
      \bot a \\
      & & \mathrm{PUSH}(q_2),& \bot q_0 :
      \ar[d]^{\gamma_{\mathrm{Push}}} &
      \bot a q_1 : & 
      \bot a \\
      q_2, \bot a \ar[dddd]^{\gamma_3}& & 
      \mathrm{CLONE}, &\bot  q_0 : \ar[d]^{\gamma_{\mathrm{Clone}}} &
      \bot a q_1 : & \bot a q_2
      \\
      & &  \mathrm{POP},&
      \bot a q_0 : \ar[d]^{\gamma_{\mathrm{Pop}}} 
      &\bot a  q_1 : &
      \bot a  q_2 : & 
      \bot a q_2 
      \\
      & & q_2,
      & \bot q_0  :  
      \ar[d]^{\gamma_3}
      & 
      \bot a q_1 : & 
      \bot a q_2 : & \bot a
      &\ar@/_14pc/[uuuuuuuuuullll]^(.5){\gamma_{\hookrightarrow}}  
      \\
      & &\mathrm{PUSH}(q_3), &\bot q_0  :
      \ar[d]^{\gamma_{\mathrm{Push}}} &
      \bot a q_1 : &
      \bot a q_2 : & \bot  
      \\
      q_3, \bot & & \mathrm{CLONE},& \bot q_0 : &  \bot a 
      q_1 : & \bot a q_2 : & \bot q_3
      \save "1,3"."1,4"*+[F]\frm{}
      \save "2,3"."2,4"*+[F]\frm{}
      \save "3,3"."3,5"*+[F]\frm{}
      \save "4,3"."4,5"*+[F]\frm{}
      \save "5,3"."5,5"*+[F]\frm{}
      \save "6,3"."6,5"*+[F]\frm{}
      \save "7,3"."7,6"*+[F]\frm{}
      \save "8,3"."8,6"*+[F]\frm{}
      \save "9,3"."9,6"*+[F]\frm{}
      \save "10,3"."10,6"*+[F]\frm{}
      \save "11,3"."11,7"*+[F]\frm{}
      \save "12,3"."12,7"*+[F]\frm{}
      \save "13,3"."13,7"*+[F]\frm{}
      \save "14,3"."14,7"*+[F]\frm{}
    }}
  \end{xy}
  
  \caption{Simulation of a nested pushdown tree in a collapsible
    pushdown graph of level $2$.}
  \label{fig:SimNPTCPG}
\end{figure}

The following definition provides the  details of this simulation.
\begin{definition}
  Let $\mathcal{N}=(Q,\Sigma,\Gamma,\Delta,q_0)$ be a pushdown system
  generating $\NPT(\mathcal{N})$. We define a corresponding collapsible
  pushdown system 
  \begin{align*}
    C(\mathcal{N}):=(Q_C,\Sigma_C, \Gamma_C,\Delta_C, PUSH(q_0))    
  \end{align*}
  of level $2$
  as follows:
  \begin{itemize}
  \item $\Sigma_C:= Q \cup \Sigma $.
  \item $\Gamma_C:=\Gamma \cup \{\gamma_{\hookrightarrow},
    \gamma_{\mathrm{Clone}}, \gamma_{\mathrm{Pop}},
    \gamma_{\mathrm{Push}}\}$ for $\gamma_{\hookrightarrow},
    \gamma_{\mathrm{Clone}}, \gamma_{\mathrm{Pop}},
    \gamma_{\mathrm{Push}}$ new symbols   
    not contained in $\Gamma$. 
  \item $Q_C:=\{\mathrm{POP}, \mathrm{CLONE}\}\cup Q \cup
    \{\mathrm{PUSH}(q):q\in Q\}$, 
    where $\mathrm{POP}$ and $\mathrm{CLONE}$, and $\mathrm{PUSH}(q)$
    are new auxiliary states 
    used to perform exactly the stack operation indicated by the
    name. 
  \item $\Delta_C$ consists of the following transitions:
    \begin{itemize}
    \item For $q\in Q$ and $\sigma\in\Sigma$, let
      \begin{align*}
       &(\mathrm{PUSH}(q), \sigma, \gamma_{\mathrm{Push}},
       \mathrm{CLONE}, \Push{q,1}),\\  
       &(\mathrm{CLONE},q,\gamma_{\mathrm{Clone}},
       \mathrm{POP},\Clone{2}), \text{ and}\\ 
       &(\mathrm{POP},q,\gamma_{\mathrm{Pop}}, q,\Pop{1})
      \end{align*}
      be in $\Delta_C$. These transitions are auxiliary transitions
      that write the state of the run onto the topmost word and
      create a clone of the topmost word preparing the simulation
      of the next transition. 
    \item For $(q,\sigma,\gamma,p,\Id)\in\Delta$, set
      $(q,\sigma, \gamma, \mathrm{PUSH}(p), \Id)\in\Delta_C$. 
    \item For $(q, \sigma, \gamma, p, \Push{\tau})\in\Delta$ add
      $(q,\sigma, \gamma, \mathrm{PUSH}(p), \Push{\tau,2})\in
      \Delta_C$.
    \item For $(q,\sigma, \gamma, p, \Pop{1})\in \Delta$, set
      $(q,\sigma,\gamma,\mathrm{PUSH(p)},\Pop{1})\in\Delta_C$. 
      This transition simulates the $\Pop{1}$ transition. Moreover,
      whenever a $\Pop{1}$ occurs, we also have to simulate the
      jump-edge. For this purpose, we set
      $(q,\sigma,\gamma_{\hookrightarrow},\mathrm{CLONE},\Collapse)\in\Delta_C$.  
    \end{itemize}
  \end{itemize}
\end{definition}

Figure \ref{fig:SimNPTCPG} shows a path in a nested pushdown tree
generated by a pushdown system $\mathcal{N}$ and the corresponding path in
$C(\mathcal{N})$.
The
following lemma shows that the original nested pushdown 
tree is first-order definable in the graph generated by $C(\mathcal{N})$. 

\begin{lemma} \label{Lem:NPTInterpretation}
  If $\mathcal{N}$ is a pushdown system that generates a nested
  pushdown tree $\NPT(\mathcal{N})$, then $\NPT(\mathcal{N})$ is 
  $\FO{3}$-interpretable in  $\CPG(C(\mathcal{N}))$.
\end{lemma}
\begin{proof}
  First of all, note that $C(\mathcal{N})$ is deterministic whenever
  it is in one of the
  states  $\{\mathrm{POP}, \mathrm{CLONE}\}\cup \{\mathrm{PUSH}(q):q\in Q\}$. 
  
  For all $q\in Q, w\in\Sigma^*$ and $s$ a stack, we
  say that $(\mathrm{CLONE},s)\in\CPG(C(\mathcal{N}))$ represents a run to
  $(q,w)$ of  $\mathcal{N}$ if $\TOP{2}(s)= w q$ (in this equality we
  forget about the links stored in $s$, of course).

  The following holds for all configurations $(\mathrm{CLONE},s)$ that
  represent some run to 
  some configuration $(q,w)$.
  \begin{itemize}
  \item $(q,\sigma,p,id)\in \Delta$ iff there is a path from
    $(\mathrm{CLONE},s)$ to $(\mathrm{CLONE}, s:\bar w)$ for $\bar w$
    a word such that 
    $(\sigma, s:\bar w)$ represents a run to $(p,w)$. If such a path
    exists, it consists of the operations
    $\Clone{2};\Pop{1};\Id;\Push{p}$.
  \item $(q,\sigma,p,\Push{\tau})\in \Delta$ iff there is a path from
    $(\sigma,s)$ to $(\sigma, s:\bar w)$ for $\bar w$ a word such that
    $(\sigma, s:\bar w)$ represents a run to $(p,w\tau)$. If such a path
    exists, it consists of the operations
    $\Clone{2};\Pop{1};\Push{(\tau,2)};\Push{p}$. Furthermore
    note, that $\tau$ has a link to the stack $s$. 
  \item $(q,\sigma,p,\Pop{1})\in \Delta$ iff there is a path from
    $(\sigma,s)$ to $(\sigma, s:\bar w)$ for $\bar w$ a word such that
    $(\sigma, s:\bar w)$ represents a run to $(p,\Pop{1}(w))$. If such
    a path 
    exists, it consists of the operations
    $\Clone{2};\Pop{1};\Pop{1};\Push{p}$.
  \end{itemize}

  From these observations, an easy induction shows that there
  is a bijection 
  from the domain of $\NPT(\mathcal{N})$ to those configurations of
  $\CPG(C(\mathcal{N}))$ 
  which are in state $\mathrm{CLONE}$. Furthermore, the transition relation
  of $\NPT(\mathcal{N})$ is $\FO{3}$-definable on this subset of
  $\CPG(C(\mathcal{N}))$.  

  Finally, we have to show the $\FO{3}$-definability of the jump-edges
  of $\NPT(\mathcal{N})$ in $\CPG(C(\mathcal{N}))$. For this purpose,
  note that a 
  $\Push{\tau}$-transition in $\NPT(\mathcal{N})$ corresponds to a
  $\Push{\tau,2}$-transition in the collapsible pushdown graph. From the
  analysis of the 
  existence of $\Push{\tau}$-transitions in $\mathcal{N}$, we obtain directly
  that this $(\tau,2)$ has a pointer to the configuration
  representing the run to the configuration precisely before this
  $\Push{\tau}$-transition is simulated. 
  When we later simulate a $\Pop{1}$-transition of $\NPT(\mathcal{N})$ that
  corresponds to this $\Push{\tau}$-transition, then we remove one of
  the clones of the corresponding $(\tau,2)$ from the stack. 
  From this, one easily sees that if $(\mathrm{CLONE},s)$ represents a run to
  some configuration $(q,w)$ such that the last operation of this run
  was a $\Pop{1}$, then the prefix of the run up to the step before the
  corresponding $\Push{\tau,2}$-transition is encoded in the unique
  configuration $(\mathrm{CLONE}, s')$ such that there are configurations
  $c,d$ such that $c\trans{\gamma_{\hookrightarrow}} (\mathrm{CLONE},s')$ and
  $c\trans{\gamma} d \trans{\gamma_{\mathrm{Push}}} (\mathrm{CLONE}, s)$. 
  It is also easy to see that all configurations that satisfy 
  this condition correspond to positions that simulate corresponding
  $\Push{\tau}$ and $\Pop{1}$-transitions. Hence, the jump-edges are
  actually $\FO{2}$-definable in $\CPG(C(\mathcal{N}))$. 
\end{proof}

\begin{corollary}
  The $\FO{}$ model checking of nested pushdown trees is decidable.
\end{corollary}

A closer look at the pushdown system $C(\mathcal{N})$ even
gives a better result: $\FO{}(\Reach)$ model checking on nested
pushdown trees is decidable.  First of all observe that
reachability in a nested pushdown tree $\NPT(\mathcal{N})$ coincides
with reachability 
in $\NPT(\mathcal{N})$ without the use of jump-edges because
jump-edges only connect 
vertices $x$ and $y$ where $y$ is a run extending the run $x$. 
But there is a one-to-one correspondence between 
reachability along the transitions of the pushdown system $\mathcal{N}$
and reachability in $\CPG(C(\mathcal{N}))$ without use
of the collapse transitions. This is due to the fact that all
transitions in $C(\mathcal{N})$ that do not perform a $\Collapse$  are
used to simulate at least one 
of the transitions of $\mathcal{N}$. Hence, the predicate $\Reach$ on
$\NPT(\mathcal{N})$ 
reduces to $\Reach_{(\Gamma_C\setminus\{\gamma_{\hookrightarrow}\})^*}$
on $\CPG(C(\mathcal{N}))$. Thus, we obtain the following extension of
the previous corollary.

\begin{theorem}
  $\FO{}(\Reach)$ model checking on nested pushdown trees is decidable. 
\end{theorem}
\begin{remark}
  Moreover, 
  $\FO{}(\Reach_{L_1}, \Reach_{L_2}, \dots, \Reach_{L_n})$ is
  decidable on $\NPT(\mathcal{N})$ if the $L_i$ are regular languages
  over $\Gamma$ (i.e., not using $\hookrightarrow$). 
  This is due to the fact that each $\gamma\in\Gamma$ has
  a direct translation into a fixed sequence of labels in the
  simulating collapsible pushdown graph. 
\end{remark}

Having shown that nested pushdown trees are first-order interpretable
in collapsible pushdown graphs of level $2$, the question arises 
 whether the reverse statement also holds. Are
collapsible pushdown graphs interpretable in the class of nested
pushdown trees? 
The answer to this question is negative if we restrict our attention
to uniform first-order interpretations.

In Lemma
\ref{LemmaLowerBoundCPGModelChecking}, we proved that the first-order
model checking on collapsible pushdown graphs of level $2$ has
nonelementary complexity.  In the next section, we present an
elementary first-order model checking algorithm for nested pushdown
trees. Since first-order interpretations can be used to transfer the
first-order model checking problem, we obtain the following theorem. 

\begin{theorem} \label{CPGinNPTnotFOInterpretable}
  There is no first-order interpretation $I$ such that for each
  collapsible pushdown graph $\mathfrak{G}$ of level $2$, there is a nested
  pushdown tree $\NPT(\mathcal{N})$ such that
  \begin{align*}
    \mathfrak{G}=\mathrm{Str}_I(\NPT(\mathcal{N})).      
  \end{align*}
\end{theorem}
\begin{proof}
  Heading for a contradiction, assume that such an interpretation $I$
  exists. Fix some collapsible pushdown graph $\mathfrak{G}$ 
  such that its first-order model checking has nonelementary expression
  complexity.
  Set $\mathfrak{N}:=\NPT(\mathcal{N})$ 
  such that $\mathfrak{G}=\mathrm{Str}_I(\mathfrak{N})$.
  By definition of a first-order interpretation, for each
  sentence $\varphi\in\FO{}$ over the vocabulary of $\mathfrak{G}$,
  there is a 
  formula $\mathrm{Frm}_I(\varphi)$ such that
  $\mathfrak{G}\models\varphi$ if and only if 
  $\mathfrak{N}\models \mathrm{Frm}_I(\varphi)$. 
  As we will see in the following section, the question 
  ``$\mathfrak{N}\models \mathrm{Frm}_I(\varphi)$?'' has elementary
  expression complexity.  
  By definition of $I$, $\mathrm{Frm}_I(\varphi)$ has length linear in
  the length 
  of $\varphi$ which implies that the algorithm has also elementary
  complexity in the size of $\varphi$. 
  But then we obtain an elementary algorithm deciding
  $\mathfrak{G}\models \varphi$ by 
  just calculating $\mathrm{Frm}_I(\varphi)$ and solving $\mathfrak{N}\models
  \mathrm{Frm}_I(\varphi)$. This contradicts our assumption on $\mathfrak{G}$. 
\end{proof}
\begin{remark}
  More generally, we can weaken our assumption on the interpretation
  $I$. Assume that there is an elementary algorithm that, on input a
  collapsible pushdown graph $\mathfrak{G}$ of level $2$, computes  an
  interpretation $I$ and a pushdown system $\mathcal{N}$ such that
  \mbox{$\mathfrak{G}=\mathrm{Str}_I(\NPT(\mathcal{N}))$}. Let $f$ be an
  elementary bound on the running time of this algorithm in terms of
  the size of the pushdown system and the formula. 
  Then we obtain the following elementary model checking algorithm on
  the class of collapsible pushdown graphs of level $2$. 
  Given $\mathfrak{G}$ and a formula $\varphi$, we compute
  $\mathcal{N}$, $I$  and $\mathrm{Frm}_I(\varphi)$ such that
  $\mathfrak{G}\mathrm{Str}_I(\NPT(\mathcal{N}))$ in time $f(\lvert
  \mathfrak{G} \rvert, \lvert \varphi \rvert)$.
  Note that $\lvert \mathcal{N} \rvert$ and the size of 
  $\mathrm{Frm}_I(\varphi)$ are bound by 
  $f(\lvert  \mathfrak{G} \rvert, \lvert \varphi \rvert)$.
  Using the model checking algorithm on nested pushdown trees, we can
  decide whether $\NPT(\mathcal{N}) \models \mathrm{Frm}_I(\varphi)$
  in $\exp(\exp(\exp(f(\lvert  \mathfrak{G} \rvert, \lvert \varphi
  \rvert))))$. 

  This solves the model checking problem on
  collapsible pushdown graphs in running time three-fold exponential
  in the elementary function $f$. This contradicts the result that
  $\FO{}$ model checking on collapsible pushdown graphs has
  nonelementary complexity. 
\end{remark}

We have seen that first-order interpretations cannot be used to define
collapsible pushdown graphs in nested pushdown trees. The question
remains open 
whether there is another logical interpretation that allows
to interpret all collapsible pushdown graphs in the class of nested
pushdown trees. Before one could give a precise answer to this
question, we would have to specify what kind of interpretation we
would like to consider. 
Nevertheless, we conjecture that the answer to this question is
negative for all meaningful concepts of logical interpretation. 
We want to point out two facts that make it hard to imagine  an
interpretation of all collapsible pushdown graphs in nested pushdown
trees. 

We already mentioned the gap in the complexity of 
$L\mu$ model checking between the two classes. Recall Theorem
\ref{CPGLmuModelCheckingComplexity} which states that 
the $L\mu$ model checking problem
of level $2$ collapsible pushdown graphs is $2$-EXPTIME complete. On
the other hand, recall that Theorem
\ref{NPTLmuModelCheckingComplexity} states that the $L\mu$ model
checking problem for nested pushdown trees is in EXPTIME. This implies
that any such interpretation would have to imply an exponential blowup
in the size of the nested pushdown tree that is used to interpret some
graph or the interpretation cannot preserve $L\mu$ formulas. 

The second fact relies on comparison of the unfoldings of collapsible
pushdown graphs and  nested pushdown trees. Recall that the class of
collapsible pushdown graphs of level $2$ encompasses also all 
higher-order pushdown graphs and these graphs are contained in the
second level of the Caucal hierarchy. Furthermore, recall
that the third level of the Caucal hierarchy is generated by applying
graph unfoldings followed by $\MSO$-interpretations to all graphs in
the second level. Hence, applying unfoldings followed by
$\MSO$-interpretations to the collapsible pushdown graphs of level
$2$, we generate a class of graphs that contains the third level of
the Caucal hierarchy. 
If we apply the same transformation to nested pushdown trees, we end
up in the second level of the Caucal hierarchy due to the following
lemma. 

\begin{lemma}
  The unfolding $\mathfrak{U}$ of a nested pushdown tree
  $\mathfrak{N}$ is the $\varepsilon$-contraction of the
  unfolding of a pushdown graph. Thus, any
  $\MSO$-interpretation on $\mathfrak{U}$ yields a graph in the second
  level of the Caucal hierarchy.
\end{lemma}
\begin{proof}
  Recall that a nested pushdown tree $\mathfrak{N}$ is almost
  unfolded, in the sense 
  that it is a tree except for the jump-edges. Thus, the unfolding of
  $\mathfrak{N}$ is
  obtained by the following operation. 
  We remove each jump-edge $\rho\hookrightarrow \pi$ and we append a new copy
  of the subtree rooted at $\pi$ to $\rho$ via a $\hookrightarrow$-edge.
  Due to the definition of $\hookrightarrow$, the stacks in the last
  configuration of $\rho$ and $\pi$ 
  agree and the run from $\rho$ to $\pi$ does only ``see''
  the topmost element of this stack. Hence, generating the unfolding
  boils down to the generation of the right number of copies of the
  configuration $(q_2,s)$ for each run $\rho\in \mathfrak{N}$ ending in
  $(q_1,s)$ and to
  attaching the subtrees induced by this configuration via
  $\hookrightarrow$ to $\rho$. As we already mentioned, the number of
  outgoing jump-edges from $\rho$ to some position with state $q_2$ only
  depends on the topmost symbol of $\rho$ and the pair $(q_1,q_2)$. 
  Using new states and $\varepsilon$-contraction, we can easily design
  a pushdown system $\mathcal{S}$ that behaves as the one generating
  $\mathfrak{N}$, but 
  which furthermore generates the right number of copies of
  $\pi$ at each configuration (by writing and removing nondeterministically
  sufficiently many dummy symbols onto/from the stack). 
\end{proof}

We
next show that first-order model checking on nested pushdown
trees has elementary complexity. More precisely, we present an
algorithm that uses doubly exponential space in the size of the
pushdown system and the size of the formula. 
For this purpose, we first investigate variants of the pumping lemma
for pushdown systems that are compatible with nested pushdown trees in
the following sense. Application of the pumping lemma to some run yields
a shorter run such that both runs share the same first-order theory up
to a certain quantifier rank. In Section \ref{SectionNPTModelChecking}
we apply 
these lemmas in order to derive a dynamic small-witness property for
nested pushdown trees. This means that for any existential quantification
that is satisfied by some nested pushdown tree, there is a short run
witnessing this quantification. As explained in Section
\ref{Sec:EFGame}, this property gives rise to a model checking
algorithm. We prove that this algorithm is in $2$-$\mathrm{EXPSPACE}$.

\subsection{A Modularity Result for Games on  Graphs of Small Diameter}
\label{SecGaifman}

We prepare the pumping lemmas mentioned above by a general result on 
Ehrenfeucht-\Fraisse games on certain graphs. 
We  show that certain tuples
of a given graph have the same $\simeq_\rho$-type.
This argument forms the 
back-bone of the modification of the pumping lemma (Lemma
\ref{PumpingLemmaTemplate}) in order to obtain
$\simeq_\rho$-preserving pumping lemmas. 

Our lemma looks like a Gaifman-locality argument, but it can be used
in situations where ordinary locality arguments fail. 
It uses a locality 
argument on induced substructures whence it can be applied to certain graphs
that have a small diameter. The  
crucial property of these graphs is that there are some generic edges
that make the diameter small in the sense that a lot of vertices are
connected to the same vertex, but when these edges are removed the
diameter becomes large. Therefore, on the graph obtained by
removing these generic edges we can apply Gaifman-like arguments in
order to establish partial isomorphisms and
$\simeq_\rho$-equivalence. Since disjoint but isomorphic neighbourhoods in
such a
graph have generic edges to the same vertices (in the full graph),
moving a tuple from one neighbourhood to the other does not change the
$\simeq_\rho$-type of the tuple. 

We use the following notation. 

For some structure $\mathfrak{G}=(V,E_1, E_2, \dots, E_n)$ with
binary relations $E_1, E_2, \dots, E_n$ and sets $A,B\subseteq V$ we say
that $A$ and $B$ \emph{touch} if $A\cap B\neq \emptyset$ or if there are $a\in A$,
$b\in B$ such that $(a,b)\in E_i$ or $(b,a)\in E_i$ for some $i\leq n$. 
For a tuple $\bar a\in A$ we define inductively the \emph{$l$-neighbourhood} of
$\bar a$ with respect to $A$, denoted $A_n(\bar a)$, by  setting
$A_0(\bar a):=\{a_i\in\bar a\}$, and 
\begin{align*}
  &A_{l+1}(\bar a):=A_l(\bar a)\cup\{ b\in A: \text{there are
  }i\leq n \text{ and }
  c \in A_l(\bar a)\text{ s.t. } (b,c)\in E_i \text{ or }
  (c,b)\in E_i\}.
\end{align*}
In terms of Gaifman-neighbourhoods, $A_{l}(\bar a)$ is the $l$-local
neighbourhood of $\bar a$ with 
respect to $\mathfrak{G}{\restriction}_{A}$. 

We say that $A$ and $B$ are isomorphic over $C\subseteq V$
and write $A\simeq_C B$ if there is some isomorphism
$\varphi:G{\restriction}_A\simeq G{\restriction}_B$ such that for
all $a\in A$, all $c\in C$, and all $1\leq i \leq n$, 
\begin{align*}
  &(a,c)\in E_i\text{ iff } (\varphi(a),c)\in E_i &\text{ and }&
  &(c,a)\in E_i\text{ iff } (c,\varphi(a))\in E_i.
\end{align*}

\begin{lemma} \label{GenericGameArgument}
  Let $\mathfrak{G} = (V,E_1, E_2, \dots, E_n)$ be some structure, 
  $A,B\subseteq V$ not touching
  and let \mbox{$\varphi: A \simeq B$} be an isomorphism of the
  induced subgraphs. 
  Let $\bar a\in A$ and \mbox{$\bar c\in C:= V\setminus \big(A_{2^\rho}(\bar a)
  \cup B_{2^\rho}(\varphi(\bar a))\big)$.} Then
  \begin{align*}
    &\varphi{\restriction}{A_{2^\rho-1}(\bar a)}:
    A_{2^\rho-1}(\bar a)\simeq_C B_{2^\rho-1}(\varphi(\bar a)) 
    &\text{ implies }&
    &  \mathfrak{G}, \bar a, \varphi(\bar a), \bar c \simeq_\rho 
    \mathfrak{G}, \varphi(\bar a),\bar a, \bar c.    
  \end{align*}
\end{lemma}
\begin{proof}
  If $\rho=0$, the claim holds trivially: since $A$ and $B$ do
  not touch, there are no edges between the elements from $\bar a$ and
  $\varphi(\bar a)$; furthermore $\varphi$ preserves all edges between
  $\bar a$ and $\bar c$. 

  We prove the lemma by induction on $\rho$. 
  We consider the first round of the 
  Ehrenfeucht-\Fraisse-game on 
  $\mathfrak{G}, \bar a, \varphi(\bar a), \bar c$ and $\mathfrak{G}, \varphi(\bar
  a),\bar a, 
  \bar c$. By symmetry, we may
  assume that Spoiler extends 
  the left-hand side $\bar a, \varphi(\bar a), \bar c$, by some $d\in
  V$. We present a winning strategy for Duplicator. 
  The general idea is the following. 

  If Spoiler  has chosen an element
  in $A\cup B$ that is close to $\bar a$ or $\varphi(\bar a)$, then
  Duplicator responds with applying the isomorphism
  $\varphi$. Otherwise, Duplicator just responds choosing the same
  element as Spoiler. The details are as follows:\\
  \emph{Local case:} if $d\in A_{2^{\rho-1}}(\bar a)$ set $a':=d$ and if
  $d\in \varphi(A_{2^{\rho-1}}(\bar a))$ set
  $a':=\varphi^{-1}(d)$. We set $\bar a':=\bar 
  a, a'$.

  Since $A_{2^{\rho-1}}(\bar a')\subseteq 
  A_{2^{\rho}}(\bar a)$, we have 
  \begin{align*}
    \bar c\in
    C':=V\setminus \big( A_{2^{\rho-1}}(\bar a')\cup
    \varphi(A_{2^{\rho-1}}(\bar a'))\big).    
  \end{align*}
  By definition, there is some set    
  \begin{align*}
    &D\subseteq \big(A\setminus A_{2^{\rho-1}}(\bar a')) \cup 
    (B\setminus B_{2^{\rho-1}}(\varphi(\bar a'))\big)\\
    \text{such that } &C'=C\cup D.     
  \end{align*}
  We claim that there is no edge between any element in $D$ and any
  element in $A_{2^{\rho-1}-1}(\bar a')$. If some $d\in D$ satisfies
  $d\in A$, then by definition it has distance at least $2$ from any
  \mbox{$a\in A_{2^{\rho-1}-1}(\bar a')$}. If $d\in D$ satisfies $d\in B$
  then it has distance at least $2$ from $a\in A_{2^{\rho-1}-1}(\bar
  a')$ because $A$ and $B$  do not touch. 
  
  Analogously, one proves that there is no edge between elements in
  $D$ and elements in $\varphi(A_{2^{\rho-1}-1}(\bar a'))$. 

  Thus, 
  we conclude that $A_{2^{\rho-1}-1}(\bar a')
  \simeq_{C'} \varphi(A_{2^{\rho-1}-1}(\bar a'))$. By
  induction hypothesis, it follows that
  \begin{align*}
    \mathfrak{G},\bar a', \varphi(\bar a'), \bar c
    \simeq_{\rho-1}   
    \mathfrak{G},\varphi(\bar a'), \bar a', \bar c.
  \end{align*}
  \emph{Nonlocal case:} 
  otherwise, 
  \begin{align*}
    d\in C' := V\setminus \big( A_{2^{\rho-1}}(\bar a)\cup
    \varphi(A_{2^{\rho-1}}(\bar a))\big)    
  \end{align*}
  and we set 
  \mbox{$\bar c':=\bar c, d$.}\\
  Similarly to the local case, we conclude that
  $A_{2^{\rho-1}-1}(\bar a) \simeq_{C'}
  \varphi(A_{2^{\rho-1}-1}(\bar a))$ because $A$ and $B$ do
  not touch and the distance between elements in $A_{2^{\rho-1}-1}(\bar
  a)$ and elements in $C'\cap A$ is at least $2$. 
  Hence, by induction hypothesis 
  \begin{align*}
    &\mathfrak{G},\bar a, \varphi(\bar a), \bar c'
    \simeq_{\rho-1} \mathfrak{G},\varphi(\bar a), \bar a, \bar c'.& 
  \end{align*}  
  Thus, this strategy is winning for Duplicator in the $\rho$-round
  game. 
\end{proof}

\subsection{$\simeq_\alpha$-Pumping on \NPT}
\label{SectionNPTPumpingLemma}

Recall that $\simeq_\alpha$ coincides with $\equiv_\alpha$. Thus, it 
describes equivalence with respect to
$\FO{\alpha}$ formulas.
In this section we want to develop a version of the pumping Lemma
for pushdown systems (Lemma \ref{PumpingLemmaTemplate}) that preserves
$\simeq_\alpha$-types in the following sense. 
Given a tuple $\bar \rho$ of runs and another run $\rho$ such that
$\rho$ is very long compared to the runs of $\bar \rho$, then we want
to apply the pumping lemma in such a way that
the resulting run $\hat\rho$ is shorter than $\rho$ and such that
$\rho \bar \rho \simeq_\alpha \hat\rho \bar \rho$. 

In order to achieve this, we use the game argument developed in the
previous 
section and we make a clever choice in the pumping argument. Let us
first explain this choice:  we want to apply the pumping lemma to
$\rho$ and obtain 
a shorter run $\hat\rho$. We apply the lemma in such a way that $\rho$
and 
$\hat\rho$ share a long prefix and they share a long suffix
in the sense that the last $n$ transitions of $\rho$ and $\hat\rho$
agree for some large $n\in\N$. Later we specify what long exactly
means, but we first want to explain how this enables us to
use the 
general game argument in order to show that 
$\rho \bar \rho \simeq_\alpha \hat\rho \bar \rho$. 

The $2^\alpha$-neighbourhood of $\rho$ divides into two parts. 
The first part, denoted by
$A_\rho$, consists 
of runs $\rho'$ that are very similar to $\rho$ in the sense that
there is a large common prefix of $\rho$ and $\rho'$. 
The other part, denoted by $C_\rho$, consists of runs that are only reachable from $\rho$
via paths that pass a very small prefix of $\rho$. 
Now, the $2^\alpha$-neighbourhood of $\hat\rho$ is isomorphic to the one of
$\rho$ in the following sense. 

The elements in $A_\rho$ are reachable from $\rho$ via a path such that
every edge  of this path only changes a small final part of the runs
connected by 
this edge. Thus,  every intermediate step shares a large initial
prefix with 
$\rho$. Since $\hat\rho$ coincides with $\rho$ on the final
transitions, the path from $\rho$ to an element in $A_\rho$ can be
copied edge by edge. We obtain an element in the neighbourhood of
$\hat\rho$ that has a large common prefix with $\hat\rho$ because each edge
that we use only changes a small final part of the runs connected by
this edge. Since this argument applies to all runs in $A_\rho$, we
obtain an isomorphic copy $B_{\hat\rho}$ in the neighbourhood of
$\hat\rho$. 

Now, we consider an element $\pi\in C_\rho$. Any path from $\rho$ to
$\pi$ starts with an initial part that is contained in $A_\rho$ and
then at some point we use a $\hookrightarrow$-edge that connects an
element $\pi'\in A_\rho$ 
with a short prefix $\pi''$ of this element. Since all elements in
$A_\rho$ share a large common prefix, $\pi''$ is a prefix of $\rho$. 
Since $\rho$ and $\hat\rho$ agree on an initial part, $\pi''$ is also a
prefix of $\hat\rho$. Now, the crucial observation is that we can copy
the path from $\rho$ to $\pi'$ edge by edge to a path from $\hat\rho$
to some $\hat\pi'$ such that $\pi''$ and $\hat\pi'$ are connected by
an $\hookrightarrow$-edge. Since this argument applies to all elements
in $C_\rho$, one 
derives that $C_\rho$ is also part of the neighbourhood of
$\hat\rho$. 

Using the game argument from the previous section, the isomorphism
between $A_\rho$ and $B_{\hat\rho}$ can be used to show that 
$\rho \bar\rho \simeq_\alpha \hat\rho \bar\rho$. 

In fact, we divide this $\simeq_\alpha$-preserving pumping lemma into three
steps. The first translates a given run  into an equivalent run
that ends in a configuration with small stack. The second
step translates such a run with small final stack into an equivalent
run that only passes small stacks. The last step translates a run that
only uses small stacks into an equivalent short run. 

Later, we use the $\simeq_\alpha$-preserving pumping argument in order to
derive an elementary bound for
the complexity of \FO{} model checking on nested pushdown trees.

In the following, we first state the three pumping lemmas that we want
to prove in this section. Afterwards, we will present the proof of each
of these lemmas. 

Before we state the first pumping lemma, we want to
recall the necessary notation. 
Let $\rho$ be some run of a pushdown system ending in configuration
$c=(q,w)$ where $q\in Q$ and $w\in\Sigma^*$. 
Recall that, e.g., we write $\Pop{1}(c)$ for $\Pop{1}(w)$ and
similarly we write $\TOP{1}(\rho)$ for $\TOP{1}(c)=\TOP{1}(w)$. 
Since we only consider level $1$ pushdown systems, $\TOP{2}(\rho)$ is
the final stack of $\rho$. 
Recall that $\width(\rho)=\width(w)$ denotes the width of the stack,
i.e., $\width(\rho)=\lvert w \rvert$. 
Now, the first pumping lemma reduces the size of the
last configuration of a given run, while preserving its
$\simeq_\alpha$-type. 

\begin{lemma}\label{FirstPumpingLemma}
  Let $\mathfrak{N}:=\NPT(\mathcal{N})$ be a nested pushdown tree. 
  Let \mbox{$\bar \rho=\rho_1, \rho_2, \dots, \rho_m \in\mathfrak{N}$} be
  runs and $\rho\in\mathfrak{N}$ another run such that 
  \begin{align*}
    \width(\rho) > \width(\rho_i) +(2+2^{\alpha+1})\lvert
    Q\rvert\cdot\lvert \Sigma\rvert + 2^{\alpha}+1\quad 
    \text{ for
    all } i\leq m.
  \end{align*}
  There is a $\hat\rho\in\mathfrak{N}$ such that
  $\width(\hat\rho) < \width(\rho)$ and 
  $\mathfrak{N},\bar \rho, \rho \simeq_\alpha \mathfrak{N},\bar \rho,
  \hat\rho$. 
\end{lemma}

In the second pumping lemma, we want to bound the size of all the
stacks occurring in a run. For this purpose, we define the following
notation. 
\begin{definition}
  Let   $\max(\rho)$ denote the size of the largest stack occurring
  within $\rho$, i.e.,   
  \begin{align*}
    \max(\rho):=\max\{ \width(\rho(i)): i\in\domain(\rho)\}.    
  \end{align*}
\end{definition}

The second pumping lemma takes a run $\rho$ and transforms $\rho$ into
an equivalent run $\hat\rho$ such that $\max(\hat\rho)$ is bounded in
terms of $\width(\hat\rho)=\width(\rho)$. 

\begin{lemma} \label{SecondPumpingLemma}
  Let $\bar \rho = \rho_1, \rho_2, \dots, \rho_m\in \mathfrak{N}$ and 
  $\rho\in \mathfrak{N}$ such that 
  \begin{align*}
    &\max(\rho)>\max(\rho_i) + \lvert Q\rvert^2 \lvert
    \Sigma\rvert+1\text{ for all 
    }1 \leq i \leq m,\text{ and such that }\\
    &\max(\rho)  > \lvert \width(\rho) \rvert + \lvert Q\rvert^2 \lvert
    \Sigma\rvert+2^\alpha+1.    
  \end{align*}
  Then there is some 
  $\hat\rho\in \mathfrak{N}$ such that
  \begin{enumerate}
  \item $\hat\rho$ and $\rho$ agree on
    their final configuration,
  \item $\max(\hat\rho)<\max(\rho)$, and
  \item $\mathfrak{N}, \bar \rho, \rho \simeq_\alpha \mathfrak{N},\bar
    \rho, \hat\rho$. 
  \end{enumerate}
\end{lemma}

In the third pumping lemma, we want to translate a run $\rho$ into an
equivalent run $\hat\rho$ such that the length of $\hat\rho$ is
bounded in terms of $\max(\hat\rho)\leq\max(\rho)$. 
For this purpose we introduce a new measure $\OccNPT$ for the length
of a run.  
We  first define $\OccNPT$. Then we present the pumping lemma
that transforms a run $\rho$ into an equivalent run $\hat\rho$ such
that $\OccNPT(\hat\rho)$ is bounded in terms of
$\max(\hat\rho)\leq\max(\rho)$. 
Afterwards, we show that the length of a run $\hat\rho$ is polynomially
bounded in $\max(\hat\rho)$ and $\OccNPT(\hat\rho)$. 

\begin{definition}
  Let $\rho$ be a run of length $n$ of some pushdown system. 
  We denote the number of occurrences of a stack $w$ in $\rho$ by 
  $\lvert \rho \rvert_w:=\big\lvert \{i\in \N:\exists q\ \rho(i)=(q,w)
  \}\big\rvert$.  
  By $\mathrm{SR}(w,\rho):=\{\hat\rho:\exists i,j\
  \hat\rho=\rho{\restriction}_{[i,j]}, w\prefixeq \rho \}$ we denote the
  set of subruns of $\rho$ whose stacks are all prefixed by
  $w$. Then we define the  
  \emph{maximal number of connected occurrences of some stack} to be
  \begin{align*} 
    &\OccNPT(\rho,w):=\max\big\{\lvert \hat\rho\rvert_w:
    \hat\rho\in\mathrm{SR}(w,\rho)\big\} \text{ and}\\
    &\OccNPT(\rho):=\max\{\OccNPT(\rho,w):w\in\Sigma^*\}. 
  \end{align*}
\end{definition}

We first state the third pumping lemma. Then we show that it indeed
bounds the length of a run. 
The lemma is based on the fact that a long run $\rho$ that does not visit
large stacks has to visit some configuration a lot of times. We can
then safely delete a subrun $\rho{\restriction}_{[i,j]}$ that connects
this configuration with itself.  The crucial observation is that this
does not change the isomorphism type of the neighbourhood 
if $\rho{\restriction}_{[i,j]}$ is approximately the middle part of $\rho$. 

\begin{lemma}\label{ThirdPumpingLemma}
  Let $\bar \rho=\rho_1, \rho_2, \dots, \rho_n
  \in\mathfrak{N}:=\NPT(\mathcal{N})$ such that  
  there is a $B_\Xi\in\N$ satisfying \mbox{$\OccNPT(\rho_i)\leq B_\Xi$} for all
  $1\leq i \leq n$. For
  $\rho\in\mathfrak{N}$, there is some $\hat\rho\in\mathfrak{N}$ such that 
  \begin{enumerate}
  \item $\max(\hat\rho) \leq \max(\rho)$,
  \item    $\rho$ and $\hat\rho$ agree on
    their final configuration,
  \item $\OccNPT(\hat\rho)\leq B_\Xi+ (2^{\alpha+1}+2) \lvert Q\rvert
    +2^\alpha+1$, and
  \item $\mathfrak{N}, \bar \rho, \rho \simeq_\alpha \mathfrak{N},\bar
    \rho, \hat\rho.$ 
  \end{enumerate}
\end{lemma}
We derive a bound on the length of $\hat\rho$ from the bound on
$\OccNPT(\hat\rho)$ by using the following
lemma. 

\begin{lemma} \label{NPTDepthBound}
  Let $\mathcal{N}$ be a pushdown system and $\rho$ a run
  of $\mathcal{N}$ such that $\max(\rho) =h$ and $\OccNPT(\rho)=b$, then
  $\length(\rho) \leq \frac{b^{h+2}-b}{b-1}$. 
\end{lemma}
\begin{proof}
  Set $m_h:=b$. For every $w\in\Sigma^h$ and some 
  subrun $\pi\in\mathrm{SR}(\rho,w)$ we have $\length(\pi) \leq m_h$
  because the width of all stacks in $s$ is $h$, which implies that
  all elements in $s$ have stack $w$. 

  Now assume that every subrun $\pi'\in\mathrm{SR}(\rho,v)$ for some
  $v\in \Sigma^{n+1}$ has $\length(\pi')\leq m_{n+1}$. Let
  $w\in\Sigma^n$ be an arbitrary word and let 
  $\pi\in\mathrm{SR}(\pi',w)$. 
  Then there are 
  \begin{align*}
    0=e_1 < e_2 <\dots < e_f < e_{f+1}=\length(\pi)    
  \end{align*}
  such that for $0\leq i \leq f$, the stack at $e_i$ in $\pi$ is $w$
  and $\pi{\restriction}_{[e_i+1, e_{i+1}-1]}$ is
  \mbox{$w_i$-prefixed} for some $w_i\in\Sigma^{n+1}$. We have $f\leq b$ due
  to $\OccNPT(\pi)\leq \OccNPT(\rho)\leq b$. 
  By assumption we get 
  $\length(\pi)\leq (1+m_{n+1}) b$.  
  Note that $\rho\in\mathrm{SR}(\rho,\varepsilon)$ whence
  \begin{align*}
  \length(\rho) \leq m_0 = b + b m_1 = b + b^2 + b^2 m_2 =
  \dots =  m_h \sum\limits_{i=0}^{h} b^i =
  \frac{b^{h+2}-b}{b-1}.    
  \end{align*}  
\end{proof}

The rest of this section is concerned with the proofs of the pumping
lemmas. The reader who is not interested in these technical details
may skip the rest of this section and continue reading Section
\ref{SectionNPTModelChecking}.

We start with some auxiliary lemmas.
These are concerned with the structure of runs that are connected by a
path of a given length $n$. 

The first observation is that the final stack of runs $\rho$ and
$\hat \rho$ that are connected by an edge differ in at most one letter. 
Using this observation inductively, we obtain the following lemma. 
\begin{lemma} \label{Lem:NPTHeigthDifferenceAlongPath}
  Let $\rho$ and $\hat\rho$ be runs that are connected by a path of
  length $n$ in  some nested pushdown trees. Then
  $\lvert \width(\rho) - \width(\hat\rho)\rvert \leq n$. 
\end{lemma}

Next, we state another auxiliary lemma concerning prefixes of
connected runs. Recall that, for $w$ some word and $\rho$ some run,
$w\prefixeq \rho$ holds if $w$ is a prefix of all stacks occurring in
$\rho$. 

\begin{lemma}  
  Let $\rho$ and $\hat\rho$ be runs of a pushdown
  system  such that  the following holds. 
  Setting $n:=\length(\rho)$, there is a word 
  $w\in \Sigma^*$, a letter $\sigma\in\Sigma$, and numbers
  $i<j\in\domain(\rho)$  such that 
  $\TOP{2}(\rho(i))=w$, 
  $w\prefixeq \rho{\restriction}_{[i,n]}$ and 
  $w\sigma \prefixeq \rho{\restriction}_{[j,n]}$.
  
  For every \mbox{$*\in \{\hookrightarrow,
    \hookleftarrow, \trans{}, \invtrans{}\}$}, 
  if $\rho \mathrel{*} \hat\rho$  
  then $\hat\rho=\rho{\restriction}_{[0,i]} \circ \hat\rho'$
  for some $\hat\rho'$ with $w\prefixeq \hat\rho'$. 
\end{lemma}
\begin{proof}
  \begin{itemize}
  \item[]
  \item 
    If  $\hat\rho\trans{}\rho$, then
    it follows immediately from $i<j\leq n$ that 
    $w\prefixeq \hat\rho':=\hat\rho{\restriction}_{[i,n-1]}$.
  \item   
    If $\rho\trans{}\hat\rho$, then 
    $\hat\rho$ extends $\rho$ by one configuration. 
    Since each stack operation alters the height of the stack by at most
    one, $\TOP{2}(\rho(n))=w\sigma$ implies directly that
    $w\prefixeq \hat\rho'{\restriction}_{[i,\length(\hat\rho)]}$.
  \item 
    If $\rho\hookrightarrow\hat\rho$, a similar argument as in the
    previous case applies.  $\hat\rho$ extends $\rho$ only by
    configurations that are prefixed by $\TOP{2}(\rho)$. 
    Since the last stack of $\rho$ is prefixed by
    $w$, the claim follows immediately.
  \item 
    Finally, consider the case that $\hat\rho\hookrightarrow \rho$.
    By
    definition of $\hookrightarrow$, we have $w \sigma \leq \rho(i)$ for 
    all $i\in 
    \domain(\rho)\setminus\domain(\hat\rho)$. 
    Furthermore, 
    $\hat\rho$ is an initial segment of $\rho$.
    Thus, $\rho{\restriction}_{[0,i]}$ is an initial segment of
    $\hat\rho$. The claim follows because
    $\hat\rho{\restriction}_{[i,\length(\hat\rho)]}$ is an initial
    segment of $\rho{\restriction}_{[i,n]}$ whence it is $w$ prefixed. \qedhere
  \end{itemize}  
\end{proof}

Iterated use of the previous lemma yields the following corollary.

\begin{corollary} \label{NPTUmgebungPrefixed}
  Let $\rho$ and $\hat\rho$ be runs of a pushdown
  system  such that  the following holds. 
  Setting $n:=\length(\rho)$, there are words 
  $w, v\in \Sigma^*$ with $\lvert v \rvert \geq m$, and numbers
  $i<j\in\domain(\rho)$  such that 
  $\TOP{2}(\rho(i))=w$, 
  $w\prefixeq \rho{\restriction}_{[i,n]}$ and 
  $wv \prefixeq \rho{\restriction}_{[j,n]}$.
  
  If $\rho$ and $\hat\rho$ are connected by a path of length $m$, then 
  $\hat\rho=\rho{\restriction}_{[0,i]}\circ\hat\rho'$ such that
  $w\prefixeq\hat\rho'$. 
\end{corollary}
\begin{proof}
  The proof is by induction on $m$. The case $m=0$ is trivial and the
  case $m=1$ is exactly the previous lemma. 
  Assume that the claim holds for some $m\in\N$. 
  Let $\rho$ and $\hat\rho$ be connected by a path of length
  $m+1$, i.e., $\rho=\rho_1 * \rho_2 * \dots * \rho_m=\hat\rho$ where
  each $*$ can be replaced by an element of
  $\{\hookrightarrow,\hookleftarrow, \trans{}, \invtrans{}\}$. 

  For $u:=\Pop{1}(\rho)$, 
  let $k\in\domain(\rho)$ be maximal such that
  $\TOP{2}(\rho(k))=u$ for some $q\in Q$. 
  By definition  $u \prefixeq \rho{\restriction}_{[k,n]}$. Due to the
  previous lemma, $\rho{\restriction}_{[0,k]}$ is an initial segment
  of $\rho_2$ and $\rho_2=\rho{\restriction}_{[0,k]} \circ \rho_2'$
  with $u\prefixeq \rho_2'$. 

  Now, $\rho_2$ and $\hat\rho$ are connected by a path of length
  $m-1$. Furthermore, $\rho{\restriction}_{[0,i]}$ is a prefix of
  $\rho_2$ and
  $w\prefixeq\rho_2{\restriction}_{[i,\length(\rho_2)]}$. Moreover,
  there is some $j$ such that $\Pop{1}(wv) \prefixeq
  \rho_2{\restriction}_{[j, \length(\rho_2)]}$.  By induction
  hypothesis we conclude that $\rho{\restriction}_{[0,i]}$ is a prefix
  of $\hat\rho$ and $w\prefixeq
  \hat\rho':=\hat\rho{\restriction}_{[i,\length(\hat\rho)]}$.
\end{proof}

We now prove the first pumping lemma that translates a given run
$\rho$ into an equivalent one with small final stack.

\begin{proof}[Proof of Lemma \ref{FirstPumpingLemma}.]
  Let $v:=\TOP{2}(\rho)$. 
  Using the proof of Lemma \ref{PumpingLemmaTemplate}, we find
  \mbox{$w_1 <w_2 < v$} and numbers $n_1 < n_2 \leq \length(\rho)$ such that
  $\rho(n_1)=(q_1,w_1)$, $\rho(n_2)=(q_2,w_2)$ and such that
  \begin{align*}
    \hat\rho:=\rho{\restriction}_{[0,n_1]} \circ
    \rho{\restriction}_{[n_2,\length(r)]}[w_2/w_1]    
  \end{align*}
  is a valid run.  
  Because of the length of $v$, we can furthermore choose
  $w_1$ and $w_2$ such that the following holds:
  \begin{enumerate}
  \item $\lvert w_1 \rvert> \width(\rho_i)$ for each $i$, 
  \item $\lvert v \rvert > \lvert w_2 \rvert +2^\alpha$, and\label{fif} 
  \item $\lvert w_2\rvert -\lvert w_1\rvert > 1+2^{\alpha+1}$.
  \end{enumerate}
  We show that $\mathfrak{N},\bar \rho, \rho \simeq_\rho
    \mathfrak{N},\bar \rho, \hat\rho$. 

  Recall that we write $\mathfrak{N}_{2^\alpha}(\rho)$ for the
  $2^\alpha$-neighbourhood of $\rho$. 
  Note that
  \begin{align*}
    \width(\rho)  -  \width(\hat\rho)  = 
    \lvert w_2\rvert -\lvert w_1\rvert  > 1+2^{\alpha+1}.    
  \end{align*}
  Using Lemma \ref{Lem:NPTHeigthDifferenceAlongPath}, one concludes
  that $\mathfrak{N}_{2^\alpha}(\rho)$ and
  $\mathfrak{N}_{2^\alpha}(\hat\rho)$ do not touch.
  
  Furthermore, due to condition \ref{fif} and
  Lemma \ref{NPTUmgebungPrefixed} it
  follows that for all
  $\pi\in \mathfrak{N}_{2^\alpha}(\rho)$ we have 
  $\pi=\rho{\restriction}_{[0,n_2]} \circ \pi'$ for some run $\pi'$ with 
  $w_2 \prefixeq  \pi'$. 
  Analogously, for all 
  $\pi\in \mathfrak{N}_{2^\alpha}(\hat \rho)$ we have 
  $\pi=\rho{\restriction}_{[0,n_1]} \circ \pi'$ for some run $\pi'$ with 
  $w_1 \prefixeq  \pi'$. 
  Lemma \ref{LemmaBlumensath} and a straightforward induction on the
  neighbourhoods of $\rho$ and $\hat\rho$ show
  that the function
  \begin{align*}
    \varphi: \mathfrak{N}_{2^\alpha}(\rho) &\to
    \mathfrak{N}_{2^\alpha}(\hat\rho)\\ 
    \pi &\mapsto \rho{\restriction}_{[0,n_1]}\circ \pi'[w_2/w_1] \text{
      where}\\
    \pi'&:=\pi{\restriction}_{[n_2,\length(\rho)]}
  \end{align*}
  is a well-defined isomorphism between 
  $\mathfrak{N}_{2^\alpha}(\rho)$ and
  $\mathfrak{N}_{2^\alpha}(\hat\rho)$. 
  
  Finally, since $\width(\rho) > \width(\hat\rho) \geq
  \lvert w_1 \rvert >  \width(\rho_i) +
  2^\alpha$, again by Lemma \ref{NPTUmgebungPrefixed}, $\rho_i$ cannot be
  in the $2^\alpha$-neighbourhood of $\rho$ or $\hat\rho$. Hence, we
  apply 
  Lemma \ref{GenericGameArgument} and obtain that
  $\mathfrak{N},\bar \rho, \rho \simeq_\alpha \mathfrak{N},\bar \rho,
  \hat\rho$. 
\end{proof}

Next,  we  prove the  second $\simeq_\alpha$-type preserving
pumping lemma 
that preserves the last configuration of a run $\rho$, but reduces
$\max(\rho)$. Recall that $\max(\rho)$ denotes the size of the largest
stack occurring in $\rho$.

\begin{proof}[Proof of Lemma \ref{SecondPumpingLemma}.]
  Let $\rho_1, \rho_2, \dots, \rho_m$, and $\rho$ be runs such that 
  \begin{align*}
    &\max(\rho)>\max(\rho_i) + \lvert Q\rvert^2 \lvert
    \Sigma\rvert+1\text{ for all 
    }1 \leq i \leq m,\text{ and such that }\\
    &\max(\rho)  > \lvert \width(\rho) \rvert + \lvert Q\rvert^2 \lvert
    \Sigma\rvert+2^\alpha+1.    
  \end{align*}

  We construct $\hat\rho$ as follows.

  Let $i\in\domain(\rho)$ be such that $\rho(i)=(q,w)$ for some $q\in
  Q$ and $w\in\Sigma^*$ with $\lvert w \rvert = \max(\rho)$. This
  implies  
  $\lvert w \rvert > \lvert Q^2\rvert \lvert \Sigma\rvert + 2^\alpha +1 
  + \width(\rho)$. 
  
  Now, using the proof of Lemma  \ref{PumpingLemmaTemplate}  we find
  $w_1 < w_2 \leq w$ and numbers 
  \begin{align*}
    n_1 < n_2 < m_2 < m_1    
  \end{align*}
  such that
  \begin{enumerate}
  \item $\max(\rho_i) < \lvert w_1 \rvert $,
  \item $\lvert w_1\rvert > \width(\rho) + 2^\alpha+1$, and
  \item $\hat\rho:= \rho{\restriction}_{[0,n_1]} \circ
    \rho{\restriction}_{[n_2,m_2]}[w_2/w_1] \circ
    \rho{\restriction}_{[m_1,\length(\rho)]}$ is a valid run.
  \end{enumerate}
  Now, we set 
  \begin{align*}
    &m_1':=m_1-(n_2-n_1) - (m_1-m_2),\\
    &\rho_A:=\rho{\restriction}_{[0,m_1+1]} \text{ and}\\
    &\rho_B:=\hat\rho{\restriction}_{[0,m'_1+1]}.
  \end{align*}
  Note that $\hat\rho=\rho_B\circ
  \rho{\restriction}_{[m_1+1,\length(\rho)]}$.

  We use Lemma \ref{GenericGameArgument} to show that
  $\mathfrak{N},\bar \rho, \rho \simeq_\alpha \mathfrak{N}, \bar \rho,
  \hat\rho$. For this purpose we set
  \begin{align*}
    &A:=\left\{\pi\in \mathfrak{N}_{2^\alpha}(\rho): \pi=
    \rho_A \circ \pi', 
    \pi'\text{ some run}\right\}  \text{ and}\\
    &B:=\left\{\pi\in \mathfrak{N}_{2^\alpha}(\hat\rho): 
      \pi=\rho_B \circ \pi',
      \pi'\text{ some run}\right\}. 
  \end{align*}
  Observe that $\rho_A\notin A$ and $\rho_B\notin B$:
  this is due to Lemma \ref{Lem:NPTHeigthDifferenceAlongPath} 
  and the fact that 
  \begin{align*}
    \width(\rho_A) = \lvert w_1 \rvert -1
    > \width(\rho) + 2^\alpha.
  \end{align*}
  The proof for $\rho_B$ and $B$ is analogous. 
  Furthermore, for all $\pi\in A$ and all $\pi'\in B$ we have 
  \begin{align*}
    \pi(m_1'+1)
    = \rho_A(m_1'+1) \neq \rho_B(m_1'+1) = \pi'(m_1'+1).    
  \end{align*}
  This is due to the fact that  $\rho_B$ ends in stack $\Pop{1}(w_1)$
  (at position $m_1'+1$) and 
  \mbox{$w_1\prefixeq \rho_A(m'_1+1)$} because 
  $n_1 \leq m'_1+1 \leq m_1$. 

  We conclude that the greatest common prefix of some $a\in A$ and
  some $b\in B$ is a proper initial prefix of both runs. Hence, $a$
  and $b$ are not connected by an edge whence $A$ and $B$ do not
  touch. 

  Furthermore, note that $\rho_i \notin A\cup B$ because for all
  $\pi\in A\cup B$, we have 
  \begin{align*}
    \max(\pi) \geq \max(\rho_B) \geq \lvert w_1 \rvert > \max(\rho_i).    
  \end{align*}
  Recall that $A_{2^\alpha}(\rho)$ denotes the
  $2^\alpha$-neighbourhood of $\rho$ in the subgraph induced by $A$. 
  We claim that there is an isomorphism $\varphi$ of the induced
  subgraphs
  \begin{align*}
    &\varphi:A_{2^\alpha}(\rho) \simeq B_{2^\alpha}(\hat\rho)\\
    &\rho_A\circ\pi \mapsto \rho_B\circ\pi.    
  \end{align*}
  For the proof of this claim,  note that
  for any two runs $\pi', \pi''$ of length at least $1$, and for 
  $*\in\{\trans{}, \invtrans{}, \hookrightarrow,\hookleftarrow\}$  
  we have 
   \begin{align*}
    &\rho_A \circ \pi' \mathrel{*}
    \rho_A \circ \pi''\text{ iff }\\
    & \rho_B \circ
    \ \pi' \mathrel{*}
    \rho_B \circ \pi''.     
  \end{align*}
  From this observation it follows by induction on the distance from
  $\rho$ that 
  \begin{align*}
    \varphi(A_{2^\alpha}(\rho)) \subseteq
    B_{2^\alpha}(\hat\rho).    
  \end{align*}
  Analogously, by induction on the distance from $\hat\rho$ one shows
  that 
  \begin{align*}
    B_{2^\alpha}(\hat\rho) \subseteq \varphi(A_{2^\alpha}(\rho)).     
  \end{align*}
  One concludes immediately that $\varphi$ is an isomorphism. 

  In order to apply the game argument, we finally have to show that
  $\varphi$ and $\varphi^{-1}$ preserve edges between
  $\mathfrak{N}\setminus (A_{2^\alpha}(\rho) \cup
  B_{2^\alpha}(\hat\rho))$ and  $A_{2^\alpha-1}(\rho)$ or
  $B_{2^\alpha-1}(\hat\rho)$, respectively. 
  Assume that $a\in A_{2^\alpha-1}(\rho)$
  and $c\in \mathfrak{N} \setminus \big( A_{2^\alpha}(\rho)\cup
  B_{2^\alpha}(\hat\rho) \big)$. We claim that if $a$ and $c$ are
  connected by some edge, then we have $c\hookrightarrow a$. 

  Note that 
  $a\trans{} c$ or $a\hookrightarrow c$ implies that $a$ is a
  subrun of $c$ and therefor $c\in A_{2^\alpha}(\rho)$ by definition of $A$. 
  If $c\trans{} a$, then 
  $\width(c) \leq \width(\rho) \rvert +2^\alpha <
  \lvert w_1 \rvert-1$. Hence, $c\neq  \rho_A$. Since
  $\rho_A$ is a proper initial segment of $a$, this
  implies $c\in A_{2^\alpha}(\rho)$. 

  Thus, if $c\in\mathfrak{N}\setminus \left(A_{2^\alpha}(\rho) \cup
  B_{2^\alpha}(\hat\rho)\right) $ is connected
  to $a$ then $c\hookrightarrow a$ and $c$
  is a proper initial segment of $\rho_A$. Since the
  last stack of $a$ and $c$ agree and 
  $\width(a) < \lvert w_1 \rvert$, $c$ is an initial
  segment of $\rho{\restriction}_{[0,n_1]}$. 
  Furthermore, 
  if the stack at $a(i)$  is prefixed by some \mbox{$v<w_1$} for all 
  \mbox{$n_1  \leq i \leq \length(a)$}, then 
  the stack of $\varphi(a)(j)$ is prefixed by some
  $v<w_1$ for all \mbox{$n_1 \leq j \leq \length(\varphi(a))$.} 
  Moreover, $\rho{\restriction}_{[0,n_1]}$ 
  is an initial segment of $\varphi(a)$ whence
  $c\hookrightarrow \varphi(a)$. 
  
  An completely analogous analysis of $\varphi^{-1}$ shows that
  $\varphi^{-1}$ preserves edges between $B_{2^\alpha-1}(\hat\rho)$ and
  $\mathfrak{N}\setminus \left( A_{2^\alpha}(\rho) \cup
    B_{2^\alpha}(\hat\rho)\right)$.  
  
  Thus, we can apply Lemma \ref{GenericGameArgument} and obtain that
  \begin{align*}
    \mathfrak{N}, \bar\rho, \rho \simeq_\alpha \mathfrak{N}, \bar\rho, \hat\rho
  \end{align*}
  and $\length(\hat\rho)<\length(\rho)$.
  
  Now, either  $\max(\hat\rho) < \max(\rho)$ or we can apply the same
  construction again to $\hat\rho$. 
  Since $\length(\rho)$ is finite and the length decreases in every
  step, we eventually construct a run $\hat\rho$ with $\max(\hat\rho)
  < \max(\rho)$. 
\end{proof}

By now, we have shown how to preserve the $\simeq_\alpha$-type of a
run while bounding the size of all stacks that occur. 

Recall the statement of Lemma \ref{NPTDepthBound}:
if the size of the stacks that occur in a run $\rho$ is bounded, then
a bound on $\OccNPT(\rho)$ can be used to calculate a bound on the
length of $\rho$. $\OccNPT(\rho)$ is the maximal number of occurrences
of a word $w$ in a $w$ prefixed subrun of $\rho$. 

For the proof of the third pumping lemma, we need some insight into the
relationship
of $\OccNPT(\rho,w)$ and $\OccNPT(\pi,w)$ for runs $\rho$ and $\pi$
that are connected in $\NPT(\mathcal{N})$. 
Before we come to these insights, we introduce the following notation. 

\begin{definition}
  For $\hat\rho=\rho{\restriction}_{[i,j]}$ we call
  $\hat\rho$ a \emph{left maximal subrun} of $\rho$ if
  $\hat\rho\in\mathrm{SR}(w,\rho)$  
  and $w\not\leq \rho(i-1)$. Analogously, we call $\hat\rho$ a 
  \emph{right maximal subrun} of $\rho$ if $\hat\rho\in\mathrm{SR}(w,\rho)$ 
  and $w\not\leq \rho(j+1)$. 
  We call $\hat\rho$ \emph{maximal} if it is left and right maximal.
\end{definition}

\begin{lemma} \label{occpmeins}
  Let $\rho=\rho_1\circ \rho_2\circ \rho_3$ be a run such that
  $\rho_2\in\mathrm{SR}(w,\rho)$  
  is maximal for some $w\in\Sigma^*$.
  If $\rho\hookrightarrow \pi$ or $\rho\trans{} \pi$ for some run
  $\pi$, then $\pi$ decomposes as
  $\pi =\rho_1\circ \pi_2 \circ \pi_3$ for
  $\pi_2\in\mathrm{SR}(w,\pi)$  maximal. In this case, we have
  \begin{align*}
  \lvert \pi_2 \rvert_w - \lvert \rho_2 \rvert_w \in \{0,1\}.     
  \end{align*}
\end{lemma}

\begin{proof}
  For $\rho\trans{} \pi$, the proof is trivial because $\pi$ extends
  $\rho$ by 
  exactly one configuration. 

 It remains to consider the
  case $\rho\hookrightarrow \pi$. 
  Due to the maximality of $\rho_2$, we have $\length(\rho_3)=0$ or
  $\rho_3(1) < w$. 
  If $\rho_3(1) < w $, then $\pi=\rho_1\circ \rho_2\circ \rho_3 \circ
  \pi'$ for some run $\pi'$ which implies $\pi_2=\rho_2$.

  Otherwise, if 
  $\length(\rho_3)=0$, then $\rho=\rho_1\circ \rho_2$. Hence,
  $\pi=\rho_1\circ \rho_2\circ \pi'$ such that the last
  stacks of $\rho_2$ and $\pi'$ agree and 
  $w \leq \rho_2(\length(\rho_2)) = \pi'(\length(\pi')) < \pi'(i)$ for
  all $1 \leq i < \length(\pi')$. Thus,
  if $w$ is the stack of $\rho_2(\length(\rho_2))$ then 
  $\lvert \rho_2 \circ \pi' \rvert_w = \lvert \rho_2 \rvert_w +1$. 
  Furthermore, if $w<\rho_2(\length(\rho_2))$, then 
  $\lvert \rho_2 \circ \pi' \rvert_w = \lvert \rho_2 \rvert_w$. 
\end{proof}

This lemma has two corollaries that we are going to use in the
proof of the third pumping lemma. 

\begin{corollary} \label{occpmeins1}
  Let $\rho, \rho'$ be runs such that $\rho \trans \rho'$ or
  $\rho\hookrightarrow \rho'$. If $\rho$ decomposes as
  \mbox{$\rho=\rho_1\circ\rho_2$} where $\rho_2$ is a maximal, $w$-prefixed
  subrun, 
  then $\rho'$ decomposes as \mbox{$\rho'=\rho_1\circ\rho_2\circ\rho_3'$}
  such that $\rho_2\circ\rho_3'$ is maximal and $w$-prefixed such that
  \begin{align*}
  \lvert \rho_2\circ\rho_3' \rvert_w - \lvert \rho_2 \rvert_w \in \{0,1\}.     
  \end{align*}
\end{corollary}

\begin{corollary} \label{occpmeins2}
  Let $\rho=\rho_1\circ \rho_2\circ \rho_3$ be a run such that
  $\rho_2\in\mathrm{SR}(w,\rho)$  
  is maximal for some $w\in\Sigma^*$.
  Let $\pi$ be  a run that is connected to $\rho$ via a path of length
  $n$ that only
  visits runs $\pi'$ such that $\rho_1$ is a prefix of $\pi'$,  then
  $\pi$ decomposes as 
  $\pi =\rho_1\circ \pi_2 \circ \pi_3$ for
  $\pi_2\in\mathrm{SR}(w,\pi)$  maximal. In this case, we have
  \begin{align*}
  \lvert \pi_2 \rvert_w - \lvert \rho_2 \rvert_w \leq n.     
  \end{align*}
\end{corollary}
A straightforward induction proves this corollary.  

Using these results, we can prove
the third pumping lemma, which bounds $\OccNPT(\rho)$. The proof relies on
the fact that for some large run $\rho$, we find initial segments
$\rho_1$ and $\rho_2$ of $\rho$ 
ending in the same configuration $(q,w)$ such that
$\lvert \rho_1 \rvert_w$ is much smaller than $\lvert \rho_2\rvert_w$
for some word $w$ and some state $q$. 

\begin{proof}[Proof of Lemma \ref{ThirdPumpingLemma}]

  Assume $\OccNPT(\rho)$ is too big in the sense that there is a word
  $w\in\Sigma^*$
  such that $\OccNPT(\rho',w) > B_\Xi+ (2^{\alpha+1}+2) \lvert Q\rvert
  +2^\alpha+1$ for some $\rho'\in \mathrm{SR}(w,\rho)$, i.e., for some
  $w$ prefixed subrun $\rho'$ of $\rho$.
  
  Then there is a decomposition of $\rho$ as 
  $\rho=\rho_1\circ \rho_2\circ \rho_3\circ \rho_4 \circ \rho_5$ such that
  the following holds.
  \begin{enumerate}
  \item $\rho_2\circ \rho_3\circ \rho_4 \in \mathrm{SR}(w,\rho)$,
  \item $\rho_2(0)=\rho_3(0)=(q,w)$ for some $q\in Q$, 
  \item $\lvert \rho_2 \rvert_w \geq 2^{\alpha+1}+2$,
  \item $\lvert \rho_3 \rvert_w > B_\Xi$,
  \item $\lvert \rho_4 \rvert_w = 2^\alpha$, and
  \item $\rho_4$ is right maximal in $\mathrm{SR}(w,\rho)$, (this implies
    $\length(\rho_5)=0$ or $w< \rho_5(1)$). 
  \end{enumerate}
  We set $\hat\rho:=\rho_1\circ \rho_3\circ \rho_4\circ \rho_5$
  omitting $\rho_2$ in $\rho$ and claim that
  $\mathfrak{N}, \bar \rho, \rho \simeq_\alpha \mathfrak{N}, \bar
  \rho, \hat\rho$. The  
  proof uses again Lemma \ref{GenericGameArgument}. Let
  \begin{align*}
    &B:=\{\hat\pi\in\mathfrak{N}: \hat\pi=\rho_1\circ \rho_3\circ \hat\pi_1\circ
    \hat\pi_2,
    \hat\pi_1\in\mathrm{SR}(w,\hat\pi) \text{ right maximal and}
    \lvert \hat\pi_1\rvert_w \leq 2^{\alpha+1}\} \text{ and}\\
    &A:=\{\pi \in\mathfrak{N}: \pi=\rho_1\circ \rho_2\circ \rho_3\circ \pi_1\circ
    \pi_2,
    \pi_1\in\mathrm{SR}(w,\pi) \text{ right maximal and}
    \lvert \pi_1\rvert_w \leq 2^{\alpha+1}\}. \\
  \end{align*}
  First note that for all $1\leq i \leq n$, $\rho_i\notin A\cup B$
  because $\OccNPT(\rho_i)<B_\Xi < \lvert \rho_3\rvert_w \leq \OccNPT(\pi)$ for
  all $\pi\in A\cup B$. 

  Now, we show that $A$ and $B$ do not touch. 
  Let
  \begin{align*}
    a = \rho_1\circ \rho_2 \circ \rho_3  \circ \pi_1\circ\pi_2 \in A    
  \end{align*}
  such that $\pi_1$ is right maximal in $\mathrm{SR}(a,w)$
  and 
  \begin{align*}
    b = \rho_1 \circ \rho_3 \circ \hat\pi_1 \circ \hat\pi_2\in B    
  \end{align*}
  such that $\hat\pi_1$ is right maximal in $\mathrm{SR}(b,w)$.
  
  Heading for a contradiction, we assume that there is some edge
  connecting $a$ and $b$. There are the following cases.
  \begin{enumerate}
  \item Assume that $a\hookrightarrow b$ or $a\trans{} b$. 
    In both cases we have $b=a\circ\pi'$ for some run $\pi'$. 
    The assumption  implies that  $\rho_2\circ\rho_3$ is a prefix of 
    $\rho_3\circ\hat\pi_1\circ\hat\pi_2$. 
    Note that $\rho_2\circ\rho_3$ is $w$ prefixed, while 
    $\hat\pi_2(1)$ is not $w$ prefixed (if $\length(\hat\pi_2)\geq 1$). 
    Thus, we conclude that $\rho_2\circ\rho_3$ is a prefix of 
    $\rho_3\circ\hat\pi_1$. But this clearly contradicts 
    \begin{align*}
      \lvert \rho_2\circ\rho_3\rvert_w \geq B_\Xi+2^{\alpha+1}+2 >
      B_\Xi+2^{\alpha+1} \geq \lvert \rho_3\circ\hat\pi_1\rvert_w.      
    \end{align*}
  \item Assume that $b\hookrightarrow a$ or $b\trans{} a$. 
    Due to $\lvert \rho_2\circ\rho_3\rvert_w  >\lvert
    \rho_3\circ\hat\pi_1\rvert_w$, $\rho_3\circ\hat\pi_1$ is a proper
    prefix of $\rho_2\circ\rho_3$.       
    
    It follows that $\length(\hat\pi_2)=0$: otherwise,
    $\hat\pi_2(1) =\rho_3(j)$ for some $j\in\domain(\rho_3)$. 
    But this leads to the contradiction that $w\not\leq
    \hat\pi_2(1)=\rho_3(j)$ due to the right maximality of $\hat\pi_1$
    but $w\leq \rho_3(j)$ by definition of $\rho_3$. 
    
    Hence, Corollary \ref{occpmeins1}  shows that 
    \begin{align*}
      \lvert \rho_2 \circ \rho_3 \circ \pi_1 \circ \pi_2 \rvert_w \leq
      \lvert \rho_3 \circ \hat \pi_1\rvert +1.
    \end{align*}
    But this contradicts the fact that 
    \begin{align*}
      \lvert \hat\pi_1\rvert_w +1\leq
      2^{\alpha+1} +1 < 2^{\alpha+2}+2\leq   \lvert \rho_2 \rvert_w.
    \end{align*}
  \end{enumerate}
  Thus, $A$ and $B$ do not touch. Now,   
  the map 
  \begin{align*}
    &\varphi:A\rightarrow B \\
    &\rho_1\circ \rho_2 \circ \pi \mapsto \rho_1\circ \pi
  \end{align*}
  is clearly well-defined. Furthermore, it is an isomorphism. For 
  $*\in\{\trans, \invtrans{},\hookrightarrow, \hookleftarrow\}$ and
  for runs $\pi, \pi'$ with 
  $\rho_2(0)=\pi(0)=\pi'(0)$ we have
  \begin{align*}
    & (\rho_1\circ \rho_2\circ \pi) \mathrel{*} (\rho_1\circ \rho_2\circ \pi')\\
    \text{iff } &\pi \mathrel{*} \pi' \\
    \text{iff } &(\rho_1\circ \pi) \mathrel{*} (\rho_1\circ \pi').
  \end{align*}
  In order to apply Lemma \ref{GenericGameArgument}, we
  have to show that $\varphi$ and $\varphi^{-1}$ preserve edges between
  $\mathfrak{N} \setminus 
  \left( A_{2^\alpha}(\rho)\cup B_{2^\alpha}(\hat\rho) \right)$ and
  $A_{2^\rho-1}(\rho)$ or $B_{2^\rho-1}(\hat\rho)$, respectively.
    
  Note that for $k < 2^\alpha$,  Corollary
  \ref{occpmeins2} states that $a\in A_k(\rho)$ implies 
  \begin{align*}
    a=\rho_1\circ \rho_2\circ \rho_3 \circ \pi_1 \circ \pi_2    
  \end{align*}
  for some
  right maximal 
  $\pi_1\in\mathrm{SR}(w,a)$  
  such that
  \mbox{$\lvert \pi_1 \rvert_w\in[2^\alpha-k,2^\alpha+k]$.}

  One immediately concludes that $a\trans{}c, c\trans{}a$, or
  $a\hookrightarrow c$ implies that $c\in A_{2^\alpha}(\rho)$ because
  \begin{align*}
    c=\rho_1\circ\rho_2\circ\rho_3\circ\hat\pi_1\circ\hat\pi_2    
  \end{align*}
  for some
  right maximal subrun $\hat\pi_1\in\mathrm{SR}(w,c)$ with 
  \mbox{$\lvert \hat\pi_1 \rvert_w\in[2^\alpha-k-1,2^\alpha+k+1]$.} 
  Since this contradicts the assumption that $c\notin A_{2^\alpha}(\rho)$,
  we only have to 
  consider the case 
  $c\hookrightarrow a$. We analyse three possibilities.
  \begin{enumerate}
  \item If the last stack of $a$ is $w$ prefixed, then
    Corollary \ref{occpmeins1} implies that $c\in A_{2^\alpha}(\rho)$
    which contradicts the assumption on $c$.
  \item If the last stack of $a$ is not $w$ prefixed and 
    $c$ is not a proper prefix of $\rho_1$, then
    \begin{align*}
      c=\rho_1\circ\rho_2\circ\rho_3\circ\pi_1\circ\hat\pi_2      
    \end{align*}
    where
    $\hat\pi_2(1)=\pi_2(1)$. But then $c\in A_{2^\alpha}(\rho)$
    which again contradicts the assumption on $c$.
  \item Finally, we consider the case that $c$ is a proper prefix of
    $\rho_1$. Since the last stack of $c$ is then a proper prefix of
    $w$, one concludes immediately that 
    \begin{align*}
      c\hookrightarrow
      \varphi(a)=\rho_1\circ\rho_2\circ\pi_1\circ\pi_2.      
    \end{align*}
  \end{enumerate}
  Using the 
  analogous arguments with reversed roles for $A$ and $B$, one  shows that
  $\varphi^{-1}$ also preserves the edges from $B_{2^\alpha-1}(\hat\rho)$ to
  $\mathfrak{N}\setminus\left( A_{2^\alpha}(\rho) \cup
    B_{2^\alpha}(\hat\rho)\right)$.  
  
  Hence, Lemma \ref{GenericGameArgument} shows
  that 
  \begin{align*}
    \mathfrak{N}, \bar \rho, \rho\simeq_\alpha \mathfrak{N} \bar\rho,
    \hat\rho.      
  \end{align*}
  Iteration of this construction eventually leads to the construction of
  some $\hat\rho$ that satisfies the lemma.
\end{proof}

\subsection{First-Order Model Checking on NPT is in 
  2-EXPSPACE} 
\label{SectionNPTModelChecking}
Using the three pumping lemmas we can now establish a
dynamic small witness property for nested pushdown trees: 
let $\varphi(x_1, x_2, \dots, x_n)$ be an \FO{} formula that is
satisfied by some 
nested pushdown tree $\NPT(\mathcal{N})$ with parameters $\rho_1,
\rho_2, \dots, \rho_n\in\NPT(\mathcal{N})$. Then the outermost
existential quantification occurring in $\varphi$ is witnessed by a
small run $\rho$ such that the length of $\rho$ is bounded in terms of
the length of $\rho_1, \rho_2, \dots, \rho_n$. 
In order to state this fact in a precise manner, we first define the
appropriate notion of a small run.  

\begin{definition}
  Let $\mathcal{N}=(Q,\Sigma,\Gamma,\Delta,q_0)$ be a pushdown system. 
  For $j\leq k\in\N$ 
  we say that some $\rho\in
  \NPT(\mathcal{N})$ is \emph{$(j,k)$-small} if 
  \begin{align*}
    & \width(\rho) \leq 6\lvert \mathcal{N}\rvert^2 j 2^k,&
    &\max(\rho)\leq 8 
    \lvert \mathcal{N} \rvert^3 j 2^k,& \text{and } 
    &\OccNPT(\rho)\leq 6  \lvert \mathcal{N}  \rvert j 2^k. 
  \end{align*}
\end{definition}

Now, we can put all the pumping lemmas together in order to prove the
existence of a small $\simeq_\alpha$-equivalent tuple for every
tuple of elements. 

\begin{lemma} \label{LemmaikSmall}
  Let $\mathcal{N}=(Q,\Sigma,\Gamma,\Delta,q_0)$ be a pushdown system and
  \begin{align*}
    \bar \rho = \rho_1, \rho_2, \dots, \rho_{i-1}\in\NPT(\mathcal{N})    
  \end{align*}
  such that $\rho_j$ is $(j,\alpha)$-small for all $1\leq j \leq
  i-1$  and $1\leq i\leq \alpha \in\N$. 
  For each $\rho_{i}\in\NPT(\mathcal{N})$, there is an
  $(i,\alpha)$-small $\rho'_{i}\in\NPT(\mathcal{N})$ such that
  \begin{align*}
    \NPT(\mathcal{N}), \bar \rho, \rho_i
    \simeq_{\alpha-i} \NPT(\mathcal{N}), \bar \rho, \rho'_i.    
  \end{align*}
\end{lemma}
\begin{proof}
  Given $\rho_i$, the first pumping lemma (Lemma
  \ref{FirstPumpingLemma}) shows that there is some
  \mbox{$a\in\NPT(\mathcal{N})$} such that
  \begin{align*}
    &\mathfrak{N}, \bar \rho, \rho_i
    \simeq_{\alpha-i} 
    \mathfrak{N}, \bar \rho,a  \text{ and} \\ 
    &\width(a) \leq 6\lvert \mathcal{N}\rvert^2 i 2^\alpha 
    + \lvert Q\rvert\lvert \Sigma\rvert (2+2^{(\alpha-i)+1}) +
    2^{(\alpha-i)}+1 
    \leq 6\lvert \mathcal{N}\rvert^2 i 2^\alpha.
  \end{align*}
  Due to the second pumping lemma (Lemma \ref{SecondPumpingLemma}),
  there is some 
  $b\in\mathfrak{N}$ such that
  \begin{align*}
    &\mathfrak{N}, \bar \rho,  a
    \simeq_{\alpha-i} 
    \mathfrak{N}, \bar \rho, b,\\ 
    &b(\length(b)) = (q, w) = a(\length(a))\text{ for some } q\in Q,
    w\in\Sigma^* \text{, and }\\
    &\max(b) \leq 8\lvert \mathcal{N}\rvert^3 i 2^\alpha + \lvert Q\rvert^2
    \lvert \Sigma\rvert+1 
    \leq 8\lvert \mathcal{N}\rvert^3 i 2^\alpha.
  \end{align*}
  Finally, we  apply the third pumping lemma (Lemma
  \ref{ThirdPumpingLemma}) and find some
  $c\in\mathfrak{N}$ such that
  \begin{align*}
    &\mathfrak{N},  \bar \rho, b
    \simeq_{\alpha-i} 
    \mathfrak{N}, \bar \rho,
     c,\\
    &c(\length(c)) = (q, w) = b(\length(b))\text{ for some } q\in Q,
    w\in\Sigma^*,\\
    &\max(c) \leq \max(b) \text{, and}\\
    &\OccNPT(c) \leq
    6 \lvert \mathcal{N} \rvert i 2^\alpha + (2^{\alpha-i+1} + 2) \lvert
    Q\rvert+2^{\alpha-i}+1
    \leq 6 \lvert \mathcal{N} \rvert i 2^\alpha.             
  \end{align*}
\end{proof}

In the terminology of Section \ref{Sec:EFGame}, the previous lemma
shows that there is a finitary constraint $S$ for Duplicator's strategy
in the Ehrenfeucht-\Fraisse game.
We set 
\begin{align*}
  S^{\NPT(\mathcal{N})}(m):=\{\rho_1, \rho_2, \dots, \rho_m
  \in\NPT(\mathcal{N})^m: \rho_i \text{ is } (i, \alpha)\text{-small }
  \text{for all } i\leq m\}  
\end{align*}
and $S:=(S_i)_{i\leq \alpha}$. 
With this notation, the previous lemma shows that Duplicator has an
$S$-preserving winning strategy in the $\alpha$-round
Ehrenfeucht-\Fraisse-game on two copies of $\NPT(\mathcal{N})$. 
As explained in Section \ref{Sec:EFGame}, such a strategy has a direct
translation into a model checking algorithm. 

\begin{theorem} \label{MCEXPTIME}
  The Algorithm \ref{Algo:NPTModelCheck} (see next page)
  solves the  
  $\FO{}$ model checking problem on nested pushdown trees, i.e., given
  a pushdown system $\mathcal{N}$ and a sentence
  $\varphi\in\FO{\alpha}$, \mbox{NPTModelCheck} accepts  the input 
  $(\mathcal{N}, \alpha, \emptyset,
  \varphi)$, if and only if $\NPT(\mathcal{N})\models \varphi$.
  The structure complexity of this algorithm  is in
  $\mathrm{EXPSPACE}$, while its expression and
  combined complexity are in $2$-$\mathrm{EXPSPACE}$. 
\end{theorem}
\begin{algorithm2e}[]
  \label{Algo:NPTModelCheck}
  \caption{\FO{} model checking on nested pushdown trees}
  \SetKwFunction{ModelCheck}{NPTModelCheck}
  \SetKw{accept}{accept}
  \SetKw{reject}{reject}
  \Titleofalgo{ \ModelCheck{$\mathcal{N}, \alpha, \bar a, \varphi(\bar x)$}}
  \KwIn{a pushdown system $\mathcal{N}$ generating
    $\mathfrak{N}:=\NPT(\mathcal{N})$, $\alpha\in\N$,
    $\varphi\in\FO{\alpha}$, an 
    assignment $\bar x \mapsto \bar a$ for tuples $\bar x, \bar a$ of
    arity $m$ such that $\bar a$ is $(m,\alpha)$-small}
  \uIf{ $\varphi$ is an atom or negated atom} {
    \lIf{ $\mathfrak{N}, \bar a \models \varphi(\bar x)$} 
    {\accept }\lElse{\reject\;}}
  \uIf{$\varphi = \varphi_1 \vee \varphi_2$}{
    \lIf{\ModelCheck{$\mathfrak{N}, \alpha, \bar a, \varphi_1$} $=$
      accept}{\accept} 
    \uElse{ \lIf {\ModelCheck{$\mathfrak{N}, \alpha, \bar a, \varphi_2$}$=$
        accept} {\accept}
      \lElse{\reject\;}
    }
  }
  \uIf{$\varphi = \varphi_1 \wedge \varphi_2$}{
    \lIf{
      \ModelCheck{$\mathfrak{N}, \alpha, \bar a, \varphi_1$}$=$
      \ModelCheck{$\mathfrak{N}, \alpha, \bar a, \varphi_2$}$=$
      accept}{\accept} 
    \lElse{\reject\;}}
  \uIf{$\varphi=\exists x \varphi_1(\bar x,x)$}{
    check whether there is an $a\in\mathfrak{N}$ such that
    $a$ is $(m+1, \alpha)$-small  and
    \ModelCheck{$\mathfrak{N}, \alpha, \bar a a, \varphi_1$}$=$ accept\;}
  \uIf{$\varphi=\forall x_i \varphi_1$}{
    check whether 
    \ModelCheck{$\mathfrak{N},\alpha, \bar a a, \varphi_1$}$=$ accept
    holds for all $(m+1, \alpha)$-small
    $a\in\mathfrak{N}$\;} 
\end{algorithm2e}

\begin{proof}
  The correctness of the algorithm follows directly from the
  correctness of  Algorithm \ref{AlgoSPReservingModelCheck} and from
  Lemma \ref{LemmaikSmall}. 

  We analyse the space consumption of this
  algorithm. Due to Lemma \ref{NPTDepthBound} an \mbox{$(i,\alpha)$-small}
  run $\rho$ has bounded length. It can be stored as a list of 
  $\exp( O( i \lvert \mathcal{N} \rvert^4 \alpha \exp(\alpha)))$
  many transitions. Thus, we need
  $\exp( O( i \lvert \mathcal{N} \rvert^4 \alpha \exp(\alpha)))
  \log(\mathcal{N})$ space for storing 
  one run. Additionally, we need space for checking whether such a
  list of transitions forms a valid run and for checking
  the atomic type of the runs. We can do this by simulation of
  $\mathcal{N}$. The size of the stack is bounded by the size of
  the runs. Since we have to store up to $\alpha$ many runs at the same
  time and $i$ is bounded by $\alpha \leq \lvert \varphi \rvert$, the
  algorithm is in 
  \begin{align*}
   &\mathrm{DSPACE}
   \big( \lvert \varphi \rvert \log(\lvert\mathcal{N}\rvert) \exp(
    O( \lvert \mathcal{N} \rvert^4 \lvert \varphi \rvert^2
    \exp(\lvert\varphi\rvert)))\big)
    \subseteq \\
    &\mathrm{DSPACE}
      \big( \exp(O( \lvert \mathcal{N} \rvert^4
      \exp(2 \lvert \varphi \vert)))\big) \subseteq
      2\text{-}\mathrm{EXPSPACE}(\lvert \mathcal{N} \rvert + \lvert \varphi
      \rvert).      
  \end{align*}
  If the formula $\varphi$ is fixed, the space consumption of the
  algorithm is exponential in the size of $\mathcal{N}$. Thus, the  
  structure complexity of first-order model checking on nested
  pushdown trees is in EXPSPACE. 
\end{proof}
\begin{remark}
  Recall that we proved the existence of a nonelementary
  $\FO{}(\Reach)$ model checking algorithm for nested pushdown trees. 
  There is no hope in finding an elementary
  algorithm. A straightforward  adaption of the proof of Theorem
  \ref{thm:FOCPGnonelementary} shows this. As
  in the case of collapsible 
  pushdown graphs, one can define a nested pushdown tree that is the
  full binary tree 
  where each branch looks like the graph in Example
  \ref{fig:nptExample}. For similar arguments as in the proof of Theorem
  \ref{thm:FOCPGnonelementary}, 
  $\FO{}$ model checking on the full infinite binary tree can be
  reduced to $\FO{}(\Reach)$ model
  checking  on this nested pushdown tree. 
\end{remark}


\section{Higher-Order   Nested   Pushdown Trees}
\label{Chapter HONPT}

In this chapter, we propose the study of a new hierarchy of graphs. We
combine the idea underlying the definition of nested pushdown trees
with the idea of higher-order pushdown systems and obtain a notion of
a higher-order nested pushdown tree. 
We first give a formal definition of this hierarchy. 
Afterwards, we  compare this new hierarchy with the
hierarchies of higher-order pushdown graphs and collapsible pushdown
graphs. 

Recall that nested pushdown trees are \FO{}-interpretable in
collapsible pushdown graphs of level $2$. We show that this result
extends to the whole hierarchy. 
Every nested
pushdown tree of level $n$ is \FO{}-interpretable in 
some collapsible pushdown
graph of level $n+1$. 

In the final part of this chapter we then prove the decidability
of the first-order model checking on level $2$ nested pushdown
trees. The approach is an adaption of the idea underlying
the decidability proof of the level $1$ case: 
we prove that there is a strategy in the Ehrenfeucht-\Fraisse game
such that Duplicator always chooses small runs. 
But the techniques involved in the proof of the existence of such a
strategy are very different from those in the level $1$ case.

\subsection{Definition of Higher-Order Nested Pushdown Trees}
\label{sec:HONPT}

We want to define the notion of higher-order nested pushdown
trees. Recall 
that a nested pushdown tree is the unfolding of a pushdown graph
extended by a jump-relation $\hookrightarrow$ that connects
corresponding push- and pop operations. Extending this idea to higher
levels, one has to define what corresponding push- and pop operations
in a level $n$ pushdown system are. 
In order to obtain well-nested jump-edges, 
we concentrate on the push- and pop operations of the
highest level, i.e., for a level $n$ pushdown system we look at
corresponding $\Clone{n}$ and $\Pop{n}$ operations.

\begin{definition}
  Let $\mathcal{N}=(\Sigma,\Gamma,Q,q_0,\Delta)$ be a pushdown system of
  level $n$.\footnote{We stress that $\mathcal{N}$ is a pushdown system
    without links and without collapse-transitions.} 
  Then the  \emph{level $n$ nested pushdown tree}
  $\mathfrak{N}:=\NPT(\mathcal{N})$ is the 
  unfolding of the pushdown graph of $\mathcal{N}$ expanded by the
  relation $\hookrightarrow$ which connects 
  each $\Clone{n}$ operation with the corresponding $\Pop{n}$
  operation, i.e.,  
  for runs $\rho_1,\rho_2$ of $\mathcal{N}$ we have 
  $\rho_1\hookrightarrow \rho_2$
  if $\rho_2$ decomposes as
  $\rho_2=\rho_1\circ \rho$ for some run $\rho$ from $(q,s)$ to $(q',s)$
  of length $n$ such that
  \begin{align*}
    &\rho(0)\trans{\Clone{n}} \rho(1), \\  
    &\rho(n-1) \trans{\Pop{n}} \rho(n) \text{, and}\\
   &\rho(i)\neq(\hat q, s)\text{ for all } 1\leq i < n 
   \text{ and all }\hat q\in Q.
  \end{align*}
\end{definition}
\begin{remark}
  Another view on the jump edges is the following. 
  Some run $\rho_1$ is connected via $\hookrightarrow$ to
  some other run $\rho_2$ if $\rho_2$ decomposes as
  $\rho_2=\rho_1\circ\rho$ where $\rho$  consists of a
  $\Clone{n}$ operation 
  followed by a ``level $n$ return''. 
  It is straightforward to show that $\rho_1\hookrightarrow \rho_2$ if
  and only if $\rho_2=\rho_1\circ\rho$ for some run $\rho$ of length
  at least $2$ such that
  $\lvert \rho(0) \rvert = \lvert \rho(\length(\rho) \rvert$ and
  $\lvert \rho(i) \rvert > \lvert \rho(0)\rvert$ for all
  $0<i<\length(\rho)$. 
\end{remark}

 In the following, we write $n$-\HONPT for ``nested pushdown tree of
 level $n$''.

\subsection{Comparison with Known Pushdown Hierarchies} 
\label{sec:CompHierarchies}

The hierarchy of higher-order nested pushdown trees is a hierarchy
strictly extending the hierarchy of trees generated by higher-order
pushdown systems. Furthermore, it is first-order interpretable in the
collapsible pushdown hierarchy. In fact, this relationship of the
hierarchies 
is level by level.
In the following, we prove these claims. 

We start by
adapting the first-order interpretation of nested pushdown trees in
collapsible pushdown graphs of level $2$ to the
interpretation of nested pushdown trees of level $n$ in collapsible
pushdown graphs of level $n+1$.  
The approach is completely analogous. 
First of all, each configuration $(q,s)$ of a level $n$ pushdown
system $\mathcal{N}$ is identified with the level $n$ stack $\Push{q,1}(s)$. 
A run $\rho$ of $\mathcal{N}$ is a list of configurations $\rho(0),
\rho(1), \dots, \rho(\length(\rho))$. 
This run is identified with the level $n+1$ stack
$s_\rho:=\rho(0):\rho(1):\dots: \rho(\length(\rho))$. 

Each extension of $\rho$ by one transition
$\delta:=(q,\sigma,\gamma,q',\op)$ can be 
simulated by a level $n+1$ pushdown system by changing the stack to 
\begin{align*}
  s_{\rho'}:=\Push{q',1}(\op(\Pop{1}(\Clone{n+1}(s_\rho)))).  
\end{align*}
It is a
straightforward observation that 
$s_{\rho'}$ represents the run $\rho'$ which is $\rho$ extended by
$\delta$. Hence, the unfolding of a level $n$ pushdown system can be
simulated by some level $n+1$ collapsible pushdown system. 

In order to simulate the nested pushdown tree generated by
$\mathcal{N}$, we also have to simulate the jump-edges. 
A jump-edge connects a $\Clone{n}$ transition with the corresponding
$\Pop{n}$ transition. Thus, the collapsible pushdown system simulating
$\mathcal{N}$ has to keep track of the positions where a $\Clone{n}$
transition was performed. 

For this purpose we introduce a clone-marker $\#$. Before the
collapsible pushdown system performs a $\Clone{n}$ transition, it
applies a $\Push{\#,n+1}$ operation. This means that it writes the symbol
$\#$ onto the stack. This symbol carries a link to the stack
representing the run up to the configuration before the $\Clone{n}$
transition was applied.  

Later, when the system simulates a $\Pop{n}$ transition, it finds a
clone of this marker $\#$ on top of the stack reached by this $\Pop{n}$. 
The link of this clone still points to the position in the run where
the corresponding $\Clone{n}$ was performed. Thus, using the collapse
operation, we can connect any position simulating a $\Pop{n}$
transition with the position that simulated the corresponding
$\Clone{n}$. 

The following proposition provides the detailed construction of the
simulating collapsible pushdown system. 

\begin{proposition} \label{Prop:HONPTinCPG}
  Let $\mathcal{N}$ be a pushdown system of level $n\geq 2$. 
  We can effectively compute a collapsible pushdown system
  $\mathcal{S}$ of level $n+1$ and a first-order interpretation
  $I_{\mathcal{N}}$ 
  such that  $\NPT(\mathcal{N})$ is first-order interpretable in
  $\CPG(\mathcal{S})$ via $I_{\mathcal{N}}$. 

  Moreover, there is a uniform bound on the length of the formulas of 
  $I_{\mathcal{N}}$ for all higher-order pushdown systems  $\mathcal{N}$. 
\end{proposition}

\begin{figure}
  \begin{xy}{\tiny
    \xymatrix@C=30pt{
      q_0, [[\bot]] \ar@{_{(}->}@/_1.5pc/[rrr] \ar[r]^{\gamma_1} &
      q_1, [\bot]:[\bot] \ar[r]^{\gamma_2} &
      q_2, [\bot] : [\bot:\bot] \ar[r]^{\gamma_3} &
      q_3, [[\bot]] \ar[r]^{\gamma_4} &
      q_4, [[\bot a]] 
      }
    }
  \end{xy} \vskip 0.5cm
  \hskip -2.5cm
  \begin{xy}{\tiny 
      \xymatrix@C=0pt@R=18pt{ 
        \mathrm{PUSH}(q_0), &[[[\bot]]]
        \ar[d]^{\gamma_{\mathrm{Push}}} 
        &\\
        \mathrm{CLONE}, &[[[\bot q_0]]]
        \ar[d]^{\gamma_{\mathrm{Clone}}}
        &\\
        \mathrm{POP}, &[[[\bot q_0]]]
        \ar[d]^{\gamma_{\mathrm{Pop}}}
        &:&[[[\bot q_0]]]
        &\\
        q_0, &[[[\bot q_0]]] \ar[d]^{\gamma_1}&:&[[[\bot ]]] \\
        \mathrm{CPP}(q_1), &[[[\bot q_0]]] 
        \ar[d]^{\gamma_{\mathrm{CPP}}}
        &:&[[[\bot(\#,4,1)]]]&
        \\
        \mathrm{PP}(q_1), &[[[\bot q_0]]] 
        \ar[d]^{\gamma_{\mathrm{PP}}}
        &:&
        [[\bot(\#,4,1)]]:[[\bot(\#,4,1)]]
        \\
        \mathrm{PUSH}(q_1), &[[[\bot q_0]]] 
        \ar[d]^{\gamma_{\mathrm{Push}}}
        &:&
        [[\bot(\#,4,1)]]:[[\bot]] &
        \\
        \mathrm{CLONE},  &[[[\bot q_0]]] 
        \ar[d]^{\gamma_{\mathrm{Clone}}}
        &:&
        [[\bot(\#,4,1)]]:[[\bot q_1]]&
        \\ 
        \mathrm{POP},  &[[[\bot q_0]]] 
        \ar[d]^{\gamma_{\mathrm{Pop}}}
        &:&
        [[\bot(\#,4,1)]]:[[\bot q_1]]
        &:&
        [[\bot(\#,4,1)]]:[[\bot q_1]]&
        \\ 
        q_1,  &[[[\bot q_0]]] 
        \ar[d]^{\gamma_2}
        &:&
        [[\bot(\#,4,1)]]:[[\bot q_1]]
        &:&
        [[\bot(\#,4,1)]]:[[\bot]]&
        \\ 
        \mathrm{PUSH}(q_2),  &[[[\bot q_0]]] 
        \ar[d]^{\gamma_{\mathrm{Push}}}
        &:&
        [[\bot(\#,4,1)]]:[[\bot q_1]]
        &:&
        [[\bot(\#,4,1)]]:[[\bot]:[\bot]]
        \\ 
        \mathrm{CLONE}, &[[[\bot q_0]]] 
        \ar[d]^{\gamma_{\mathrm{Clone}}}
        &:&
        [[\bot(\#,4,1)]]:[[\bot q_1]]
        &:&
        [[\bot(\#,4,1)]]:[[\bot]:[\bot q_2]]&&
        \\ 
        \mathrm{POP}, &[[[\bot q_0]]] 
        \ar[d]^{\gamma_{\mathrm{Pop}}}
        &:&
        [[\bot(\#,4,1)]]:[[\bot q_1]]
        &:&
        [[\bot(\#,4,1)]]:[[\bot]:[\bot q_2]]
        &:&
        [[\bot(\#,4,1)]]:[[\bot]:[\bot q_2]]
        &\\ 
        q_2, &[[[\bot q_0]]] 
        \ar[d]^{\gamma_3}
        &:&
        [[\bot(\#,4,1)]]:[[\bot q_1]]
        &:&
        [[\bot(\#,4,1)]]:[[\bot]:[\bot q_2]]
        &:&
        [[\bot(\#,4,1)]]:[[\bot]:[\bot]]&
        \\ 
        \mathrm{PUSH}(q_3), &[[[\bot q_0]]] 
        \ar[d]^{\gamma_{\varepsilon}}
        &:&
        [[\bot(\#,4,1)]]:[[\bot q_1]]
        &:&
        [[\bot(\#,4,1)]]:[[\bot]:[\bot q_2]]
        &:&
        [[[\bot(\#,4,1)]]] 
        \ar `r/3pc [ur]
        `/20pc[uuuuuuuuuuuuulllll]
        [uuuuuuuuuuuuullllll]^{\gamma_{\hookrightarrow}}  
        \\ 
        \mathrm{PUSH}(q_3), &[[[\bot q_0]]] 
        \ar[d]^{\gamma_{\mathrm{Push}}}
        &:&
        [[\bot(\#,4,1)]]:[[\bot q_1]]
        &:&
        [[\bot(\#,4,1)]]:[[\bot]:[\bot q_2]]
        &:&
        [[[\bot]]]
        \\ 
        \mathrm{CLONE}, &[[[\bot q_0]]] 
        \ar[d]^{\gamma_{\mathrm{Clone}}}
        &:&
        [[\bot(\#,4,1)]]:[[\bot q_1]]
        &:&
        [[\bot(\#,4,1)]]:[[\bot]:[\bot q_2]]
        &:&
        [[[\bot q_3]]]
        \\  
        \mathrm{POP}, &[[[\bot q_0]]] 
        \ar[d]^{\gamma_{\mathrm{Pop}}}
        &:&
        [[\bot(\#,4,1)]]:[[\bot q_1]]
        &:&
        [[\bot(\#,4,1)]]:[[\bot]:[\bot q_2]]
        &:&
        [[\bot q_3]]
        &:
        [[\bot q_3]]
        \\ 
        q_3, &[[[\bot q_0]]] 
        \ar[d]^{\gamma_4}
        &:&
        [[\bot(\#,4,1)]]:[[\bot q_1]]
        &:&
        [[\bot(\#,4,1)]]:[[\bot]:[\bot q_2]]
        &:&
        [[[\bot q_3]]]
        &:
        [[[\bot]]]
        \\ 
        \mathrm{PUSH}(q_4), &[[[\bot q_0]]] 
        &:&
        [[\bot(\#,4,1)]]:[[\bot q_1]]
        &:&
        [[\bot(\#,4,1)]]:[[\bot]:[\bot q_2]]
        &:&
        [[[\bot q_3]]]
        &:
        [[[\bot a]]]
        \save "1,1"."1,2"*+[F]\frm{}
        \save "7,1"."7,4"*+[F]\frm{}
        \save "11,1"."11,6"*+[F]\frm{}
        \save "15,1"."15,8"*+[F]\frm{}
        \save "20,1"."20,9"*+[F]\frm{}
    }}
  \end{xy}
  \caption[Simulation of $3$-\HONPT in level $4$ collapsible pushdown
    graphs]{Simulation of $3$-\HONPT in level $4$ collapsible pushdown
    graphs; $\gamma_1$ is a $\Clone{3}$ transition, $\gamma_2$ is a
    $\Clone{2}$ transition, $\gamma_3$ is a $\Pop{3}$ transition and 
  $\gamma_4$ a $\Push{a}$ transition.} 
  \label{fig:HONPTinCPG}
\end{figure}

\begin{proof}[Proof (Sketch).]
  We prove this fact by a straightforward extension of the $n=1$
  case (cf. Lemma \ref{Lem:NPTInterpretation}). 
  Figure \ref{fig:HONPTinCPG} illustrates the simulation of 
  a
  $3$-\HONPT in a collapsible pushdown graph of level $4$. 
  
  Let $\mathcal{N}=(Q,\Sigma,\Gamma,\Delta,q_0)$ be
  a pushdown system 
  of level $n>1$ 
  generating $\mathfrak{N}:=\NPT(\mathcal{N})$. Then we define a collapsible
  pushdown system of level $n+1$
  $C(\mathcal{N}):=(Q_C,\Sigma_C, \Gamma_C,\Delta_C, I)$ 
  as follows.
  \begin{itemize}
  \item $\Sigma_C:= Q \cup \Sigma \cup\{\#\}$ for a new symbol
    $\#$ which is used to simulate the jump-edges. 
  \item $\Gamma_C:=\Gamma \cup \{\gamma_{\mathrm{Init}},
     \gamma_{\mathrm{Clone}},
    \gamma_{\mathrm{Pop}}, \gamma_{\mathrm{Push}},
    \gamma_{\mathrm{CPP}}, 
    \gamma_{\mathrm{PP}},
    \gamma_{\hookrightarrow},
    \gamma_{\varepsilon}\}$ for 
    new symbols
    not contained in $\Gamma$. 
  \item $Q_C:=Q \cup \Sigma \cup \{I, \mathrm{POP}, \mathrm{CLONE}\} \cup
    \{\mathrm{PUSH}(q):q\in Q\} \cup 
    \{\mathrm{CPP}(q): q\in Q\} \cup
    \{\mathrm{PP}(q): q\in Q\}$, 
    where $I$ is the new initial state, and the other states are 
    new auxiliary states for the simulation process. 
  \item $\Delta_C$ consists of the following transitions.
    \begin{enumerate}
    \item For the initialisation, we add the transition
      $(I,\bot,\gamma_{\mathrm{Init}},PUSH(q_0),\Id)\in
      \Delta_C$.
    \item 
      For $q\in Q$ and $\sigma\in\Sigma$, let
      \begin{align*}
        &(PUSH(q), \sigma, \gamma_{\mathrm{Push}}, \mathrm{CLONE}, \Push{q}),\\
        &(\mathrm{CLONE},q,\gamma_{\mathrm{Clone}},
        \mathrm{POP},\Clone{n+1}),\text{ and}\\
        &(\mathrm{POP},q,\gamma_{\mathrm{Pop}}, q,\Pop{1})        
      \end{align*}
      be in $\Delta_C$.\footnote{In the following, we write $\Push{q}$
        for $\Push{q,1}$.} These transitions are auxiliary transitions
      that write the state of the run onto the topmost level $n$
      stack and  create a clone of the topmost level $n$ stack
      preparing the simulation of the next transition. 
    \item For $\op\neq\Clone{n}$ and
      $(q,\sigma,\gamma,p,\op)\in\Delta$, set
      $(q,\sigma, \gamma, PUSH(p), \op)\in\Delta_C$.
    \item For 
      $(q,\sigma, \gamma, p, \Clone{n})\in\Delta$, set
      $(q,\sigma, \gamma , \mathrm{CPP}(p),
      \Push{\#,n+1})\in\Delta_c$.
    \item We handle the jump-edge marker $\#$ with the following
      transitions.  For all $q\in Q$, set 
      \begin{align*}
        &(\mathrm{CPP}(q), \#, \gamma_{\mathrm{CPP}}, \mathrm{PP}(q),
        \Clone{n}) \in \Delta_C,\\
        &(\mathrm{PP}(q), \#, \gamma_{\mathrm{PP}}, \mathrm{PUSH}(q),
        \Pop{1})\in \Delta_C,\\
        &(\mathrm{PUSH}(q), \#, \gamma_{\varepsilon}, \mathrm{PUSH}(q),
        \Pop{1})\in\Delta_C,\text{ and}\\
        &(\mathrm{PUSH}(q), \#, \gamma_{\hookrightarrow}, \mathrm{CLONE},
        \Collapse)\in\Delta_C.
      \end{align*}
      The first and the second transition are used to create the
      jump-edge 
      marker whenever a $\Clone{n}$ is simulated. The third transition
      is used to remove the marker after the simulation of a
      $\Pop{n}$. The last transition is used to simulate the
      jump-edge.  
    \end{enumerate}
  \end{itemize}
  We use those configurations with state $\mathrm{PUSH}(q)$ for all $q\in Q$ 
  that have no incoming $\gamma_{\varepsilon}$-edge 
  as representatives of the runs of $\mathcal{N}$.
  These configurations are defined by the formula
  \begin{align*}
    \varphi(x):=\exists y x\trans{\gamma_{\varepsilon}} y \lor
    \left( x\trans{\gamma_{\mathrm{Push}}} y \land \forall z \neg
    z\trans{\gamma_{\varepsilon}} x\right).       
  \end{align*}
  Now, we turn to the formulas that interpret the transitions
  $\trans{\gamma}$.  
  Let $\rho, \hat\rho\in\mathfrak{N}$ be connected by some transition
  $\delta=(q,\sigma,\gamma,p,\op)\in\Delta$. 
  We denote by
  $\rho'$  the representative of $\rho$ and by
  $\hat \rho'$ the representative of $\hat \rho$ in $C(\mathcal{N})$. 
  We distinguish the following cases.
  \begin{enumerate}
  \item Assume that the last transition of
    $\rho$ is not a $\Pop{n}$ transition and $\op\neq\Clone{n}$. 
    Then the transition $\rho \trans{\gamma} \hat\rho$ in
    $\mathfrak{N}$ corresponds to  
    a chain
    \begin{align*}
      \rho' \trans{\gamma_{\mathrm{Push}}} x_1
      \trans{\gamma_{\mathrm{Clone}}} x_2 \trans{\gamma_{\mathrm{Pop}}}
      x_3 \trans{\gamma} \hat \rho'      
    \end{align*}
    in $\CPG(C(\mathcal{N}))$. 
  \item 
    Assume that the last transition of
    $\rho$ is a $\Pop{n}$ transition and
    $\op\neq\Clone{n}$. 
    Then the transition $\rho \trans{\gamma} \hat\rho$ in
    $\mathfrak{N}$ corresponds to a chain
    \begin{align*}
      \rho' \trans{\gamma_\varepsilon} x_4 \trans{\gamma_{\mathrm{Push}}} x_1
      \trans{\gamma_{\mathrm{Clone}}} x_2 \trans{\gamma_{\mathrm{Pop}}}
      x_3 \trans{\gamma} \hat \rho'
    \end{align*}
    in $\CPG(C(\mathcal{N}))$.
  \item 
    Assume that the last transition of
    $\rho$ is not a $\Pop{n}$ transition and $\op=\Clone{n}$. 
    Then the transition $\rho \trans{\gamma} \hat\rho$ in
    $\mathfrak{N}$ corresponds to 
    a chain
    \begin{align*}
      \rho' \trans{\gamma_{\mathrm{Push}}} x_1
      \trans{\gamma_{\mathrm{Clone}}} x_2 \trans{\gamma_{\mathrm{Pop}}}
      x_3 \trans{\gamma}  x_5 \trans{\gamma_{\mathrm{CPP}}} x_6
      \trans{\gamma_{\mathrm{PP}}} \hat \rho'      
    \end{align*}
    in $\CPG(C(\mathcal{N}))$. 
  \item 
    Assume that the last transition of
    $\rho$ is  a $\Pop{n}$ transition and $\op=\Clone{n}$. 
    Then the transition $\rho \trans{\gamma} \hat\rho$ corresponds to 
    a chain
    \begin{align*}
      \rho' \trans{\gamma_\varepsilon} x_4 \trans{\gamma_{\mathrm{Push}}} x_1
      \trans{\gamma_{\mathrm{Clone}}} x_2 \trans{\gamma_{\mathrm{Pop}}}
      x_3 \trans{\gamma}  x_5 \trans{\gamma_{\mathrm{CPP}}} x_6
      \trans{\gamma_{\mathrm{PP}}} \hat \rho'      
    \end{align*}
    in $C(\mathcal{N})$. 
  \end{enumerate}
  Moreover, every chain that starts and ends in nodes defined by
  $\varphi$ and that is of one of the forms mentioned in the case
  distinction corresponds to a transition in $\mathfrak{N}$. 

  This claim is proved by induction on the length of the shortest path
  to some node satisfying $\varphi$. It is completely analogous to the
  corresponding proof in Lemma \ref{Lem:NPTInterpretation}. 
  
  Finally, we give an interpretation for the jump-edge relation
  $\hookrightarrow$.
  The jump-edges correspond to the edges defined by
  \begin{align*}
    \varphi_{\hookrightarrow}(x,y) :=\exists z
    (x\trans{\gamma_{\mathrm{Push}}} z \land
    y\trans{\gamma_{\hookrightarrow}} z). 
  \end{align*}
\end{proof}

The previous proposition shows that higher-order nested pushdown trees are
(modulo \FO{}-interpretations) contained in the collapsible
pushdown hierarchy. 
The hierarchy of nested pushdown trees is also an extension of the
pushdown tree hierarchy. This is shown in the following lemma. 
\begin{lemma}
  The unfoldings of graphs of level $n-1$ pushdown systems are
  contained in the $n$-th level of the nested pushdown tree
  hierarchy. 
\end{lemma}
\begin{proof}
  Consider any level $n-1$ pushdown system $\mathcal{S}$ as a
  level $n$ system that does not use $\Clone{n}$. Then
  $\mathcal{S}$ generates a level $n$ nested pushdown tree which
  coincides with the unfolding of the configuration graph of
  $\mathcal{S}$.  
\end{proof}
\begin{remark}
  Recall that the unfoldings of higher-order pushdown graphs form the
  pushdown tree hierarchy. 
  The previous lemma shows that the nested pushdown tree hierarchy is
  an extension of the pushdown tree hierarchy. 

  It is an interesting open question what the exact relationship
  between the  hierarchy of pushdown graphs and the hierarchy of
  nested trees is. 
  Since there are nested pushdown trees that have undecidable
  \MSO-theory (cf. Lemma \ref{MSO-NPT-MC-undecidable}), the hierarchy
  of nested pushdown trees is not contained in the hierarchy of
  pushdown graphs. 
  But it is an open question whether there is some logical
  interpretation that interprets every nested pushdown tree in some
  higher-order pushdown graph. 
  Lemma \ref{MSO-NPT-MC-undecidable} only implies that there is no
  $1$-dimensional \MSO interpretation that interprets nested pushdown
  trees in higher-order pushdown graphs. 
\end{remark}

The previous lemma and Proposition \ref{Prop:HONPTinCPG}
locate the 
hierarchy of nested pushdown trees between the hierarchy of pushdown
trees and the hierarchy of collapsible pushdown graphs. 
We propose the study of this new hierarchy in order to obtain new insights
into the relationship of the hierarchies of collapsible pushdown
graphs and higher-order pushdown graphs. 
In the following, we show that \FO{} model checking on $2$-\HONPT is
decidable. Via the interpretation of nested pushdown trees in
collapsible pushdown graphs, this can be seen as the first step
towards an characterisation of the largest subclass of the class of
collapsible pushdown graphs of level $3$ on which the \FO{} model
checking problem is decidable.

\subsection{Towards FO Model Checking on  Nested Pushdown Trees of Level 2}

In the following,  we develop an \FO{} model
checking algorithm on nested pushdown trees of level $2$. 

Before we continue, we want to stress that the rest of this chapter
deals exclusively with level $2$ pushdown systems and not with level
$2$ collapsible pushdown systems. Thus, stacks do not
carry any link structure and the systems never use collapse
operations. 
In this setting, loops and returns play an even more important role
than in the setting of collapsible pushdown systems. 
In runs of  pushdown systems of level $2$,  loops and returns occur almost
everywhere in the following sense:
\begin{enumerate}
\item every run $\rho$ from some stack $s$ to a substack of
  $\Pop{2}(s)$ has an initial part that is a return and
\item every run $\rho$ that starts and ends in stack $s$ and that
  never visits $\Pop{2}(s)$ is a loop. 
\end{enumerate}
We leave it as an easy exercise to check the correctness of these claims. 
In the following, we will use these facts without any further explanation.

We want to provide an $\FO{}$ model checking algorithm for
the  class of nested pushdown trees of level $2$.
We do this by 
adapting our approach for first-order model checking on nested
pushdown trees of level $1$. 
Fix some pushdown system  $\mathcal{N}$ of level $2$.
We show that
every formula of the form $\exists x 
\varphi$ such that \mbox{$\NPT(\mathcal{N}), \bar \rho \models\exists
x\varphi$}  has a short 
witness \mbox{$\rho\in \NPT(\mathcal{N})$} for the first existential
quantification. Here, the size of an  
element is given by 
the length of the run of $\mathcal{N}$ representing this element. 
We consider a run to be short, if
its size is bounded in terms of the length of the 
runs in the tuple $\bar\rho$ of parameters. 

As in the level $1$ case, we prove this dynamic small-witness property
via Ehrenfeucht-\Fraisse games.
The rough picture of the proof is as follows.

We analyse the $\alpha$-round Ehrenfeucht-\Fraisse game on two
copies of \mbox{$\mathfrak{N}:=\NPT(\mathcal{N})$}. We show that
Duplicator has a 
strategy that answers 
every move of Spoiler by choosing a small
element. An element is small if there is a bound on the size
of the element in terms of the size of the elements chosen so far in
the same copy of $\mathfrak{N}$. 
Using such a strategy, we obtain a model checking algorithm on nested
pushdown trees of level $2$ as explained in Section \ref{Sec:EFGame}. 

On this level of detail, the decidability proof on level $2$
is exactly the same as on level $1$. 
But the proof that Duplicator can always choose small runs is
completely different. 

The main technical tool for this proof is the concept of
\emph{relevant ancestors}. For each element of
$\mathfrak{N}$, the relevant 
$l$-ancestors are a finite set of initial subruns of this element. 
Intuitively, the relevant $l$-ancestors of a run $\rho$ are finitely
many ancestors of $\rho$ 
that give a description of
the $l$-local neighbourhood of $\rho$. Surprisingly, this finite
description is sufficiently complete for the purpose of preserving
partial isomorphisms during the Ehrenfeucht-\Fraisse game. 
We prove that there is a winning strategy for Duplicator with the
following property. 
Duplicator always  chooses small runs whose relevant $l$-ancestors are
isomorphic to the relevant $l$-ancestors of the element chosen by
Spoiler. 

In order to find such a strategy for Duplicator,
we analyse the structure of relevant ancestors. 
We show that a relevant
ancestor $\rho_1$ is connected to the next one, say $\rho_2$,
by either a single transition or 
by a run $\rho$ of a certain kind. 
This run $\rho$ satisfies the following conditions:
$\rho_2$ decomposes as 
$\rho_2=\rho_1\circ\rho$, 
the initial
stack of $\rho$ is $s:w$ where $s$ is some stack and $w$ is some
word. The final stack of $\rho$ is $s:w:v$ for some word $v$ and
$\rho$ does never pass a proper substack of $s:w$.

Due to this result, a typical set of relevant ancestors is of the form
\begin{align*}
  \rho_1 \prec \rho_2 \prec \rho_3 \prec \dots \prec \rho_m=\rho,
\end{align*}
where
$\rho_{n+1}$ extends $\rho_n$ by either one transition or by a run
that extends the last stack of $\rho_n$ by a new word $v$. 
If we want to construct a run $\rho'$ with isomorphic relevant
ancestor set, we have to provide runs 
\begin{align*}
  \rho'_1 \prec \rho'_2 \prec \rho'_3 \prec \dots \prec \rho'_m=\rho'  
\end{align*}
where $\rho_{n+1}'$ extends $\rho_n'$ in exactly the same manner as
$\rho_{n+1}$ extends $\rho_n$. 

We first concentrate on one step of this
construction.  Assume that $\rho_1$ ends in some configuration $(q,s:w)$ and
$\rho_2$ extends $\rho_1$ by a run creating the stack $s:w:v$.
How can we find another stack $s'$ and words $w',v'$ such that
there is a run $\rho_1'$ to $(q,s':w')$ and a run $\rho_2'$ that
extends $\rho_1'$ by a run from $(q,s':w')$ to the stack $s':w':v'$?

We introduce a family of equivalence relations on words
that preserves the existence of such runs. 
If we find some $w'$ that is equivalent to $w$
with respect to the $i$-th equivalence relation, then for any run from
$s:w$ to $s:w:v$ we can find a run from $s':w'$ to $s':w':v'$ for $v$
and $v'$ equivalent with respect to the $(i-1)$-st equivalence
relation. 

Let us explain the ingredients of these equivalence relations. 
Let $\rho_1$ be a run to some stack $s:w$ and let $\rho_2$ be a run
that extends $\rho_1$ and ends in a stack $s:w:v$. 
Recall that the theory of generalised milestones shows that the final
segment of $\rho_2$ is of the form
\begin{align*}
  \lambda_n\circ\op_n\circ\lambda_{n-1}\circ\op_{n-1}\circ\dots\circ
  \op_1\circ\lambda_0  
\end{align*}
where the $\lambda_i$ are loops and $\op_{n},
\op_{n-1}, \dots, \op_1$ is the minimal sequence generating $s:w:v$
from $s:w$. 
Thus, we are especially interested in the loops of each prefix
$\Pop{1}^k(w)$ of $w$ and 
each prefix $\Pop{1}^k(w')$ of $w'$. 
For this purpose we consider the word models of $w$ and $w'$ enriched by 
information on runs between certain prefixes of $w$ or $w'$.
Especially, each prefix is annotated with the number of possible loops
of each prefix.  
$w$ and $w'$ are equivalent with respect to the first equivalence
relation if the $\FO{k}$-types
of their enriched word structures coincide. 
The second, third, etc. equivalence relation is then defined as follows.
We enrich every element of the word model of some word $w$ by the
equivalence class of the corresponding prefix with respect to the
$(i-1)$-st equivalence relation. 
The $i$-th equivalence relation then compares the $\FO{k}$-types of
these enriched word models. This means that 
two words $w$ and $w'$ are equivalent with respect to the $i$-th
equivalence relation if the $\FO{k}$-types of their word models
enriched with the $(i-1)$-st equivalence class of each prefix coincide.

This iteration of equivalence of prefixes leads to the following
result. Let $w$ and $w'$ be equivalent with respect to the $i$-th
relation. Then we can transfer runs creating $i$ words in the
following sense: if
$\rho$ is a run creating $w:v_1:v_2:\dots:v_i$ from $w$, then there is
a run $\rho'$ creating $w':v'_1:v'_2:\dots:v'_i$ from $w'$ such that
$v_k$ and $v'_k$ are equivalent with respect to the $(i-k)$-th
relation. 
This property then allows to construct isomorphic relevant ancestors
for a given set of relevant ancestors of some run $\rho$. We only have
to start with a 
stack $s':w'$ such that $w'$ is $i$-equivalent to the topmost word of
the minimal element of the relevant ancestors of $\rho$ for some large
$i\in\N$. 

This observation reduces the problem of
constructing runs with isomorphic relevant ancestors to the problem of
finding runs whose last configurations have equivalent topmost
words (with respect to the $i$-th equivalence relation for some
sufficiently large $i$) such that one of these runs is always short. 

We solve this problem by application of several pumping constructions
that respect 
the equivalence class of the topmost word of the final configuration
of a run but which decrease the length of the run. 

Putting all these results together, we obtain that Duplicator has an
$S$-preserving strategy on every nested pushdown tree of level $2$
where $S$ is a finitary constraint bounding the length of the runs
that Duplicator may choose. Then we use the general model checking
algorithm from Section \ref{Sec:EFGame} in order to solve the $\FO{}$ 
model checking problem on nested pushdown trees of level $2$.

The outline of the next sections is as follows. 
In Section \ref{SectionRelevantAncestors} we define the important notion
of relevant ancestors and develop some theory concerning these sets. 
We then 
define a family of equivalence relations on words and stacks in
Section \ref{subsec:EquivalenceonWords}. 
In Section \ref{SectionEquivalenceonTuples} we put these things
together: the equivalence on stacks gives us a transfer property of  
relevant ancestors to isomorphic copies. Our analysis of loops
(cf. Section \ref{ChapterLoops})
yields the possibility to bound the length of the runs involved in the
isomorphic copy. Thus, preserving isomorphisms between relevant
ancestors while choosing small runs is a valid strategy for Duplicator
in the Ehrenfeucht-\Fraisse game. This gives us a small-witness
property which we use to  show the decidability of \FO{} model
checking on $2$-\HONPT in Section \ref{SectionFODecidability}.


\subsection{Relevant Ancestors}
\label{SectionRelevantAncestors}

This section aims at identifying those ancestors of a run $\rho$ in a
$2$-\HONPT $\mathfrak{N}$ that are relevant with respect to its
$\FO{k}$-type. We show that  
only finitely many ancestors of a certain kind fix the
$\FO{k}$-type of the
$l$-local neighbourhood of $\rho$.
We call these finitely many ancestors the
\emph{relevant $l$-ancestors} of $\rho$.

Before we formally introduce relevant ancestors, we 
recall some important
abbreviations concerning configurations and runs.
Abusing notations we apply functions defined on stacks to
configurations. For example if $c=(q,s)$ we write $\lvert c \rvert$
for $\lvert s \rvert$ or $\Pop{2}(c)$ for $\Pop{2}(s)$. 

We further abuse this notation by application
of functions defined on stacks to some run $\rho$, meaning that we
apply the function to the last stack occurring in $\rho$. For example,
we write 
$\TOP{2}(\rho)$ for $\TOP{2}(s)$ and $\lvert \rho \rvert$ for $\lvert s
\rvert$  if $\rho(\length(\rho))=(q,s)$.

In the same sense one has to understand equations like
$\rho(i)=\Pop{1}(s)$. This equation says that
$\rho(i)=(q,\Pop{1}(s))$ for some $q\in 
Q$. Keep in mind that $\lvert \rho \rvert$ denotes the width of  
the last stack of $\rho$ and not the length $\length(\rho)$ of the run
$\rho$. 
Recall also that we write $\rho \preceq \rho'$ if the run $\rho$ is an
initial segment of the run $\rho'$. 

\begin{definition}
  Let $\mathfrak{N}$ be some $2$-\HONPT.
  Define the relation  $\overset{+1}{\hookrightarrow} \subseteq
  \mathfrak{N}\times \mathfrak{N}$ by
  \begin{align*}
  \rho \overset{+1}{\hookrightarrow} \rho'  \text{ if }
  \rho\prec \rho', \lvert \rho \rvert = \lvert \rho' \rvert -1, \text{
    and }
    \lvert \pi \rvert > \lvert \rho \rvert \text{ for all } \rho\prec
    \pi \prec \rho'. 
  \end{align*}
  We define the 
  \emph{relevant $l$-ancestors of $\rho$} by induction on $l$.
  The relevant $0$-ancestors of $\rho$ are the elements of the set
  $\RelAnc{0}{\rho}:=\{\rho\}$. 
 Inductively, we set 
  \begin{align*}
    \RelAnc{l+1}{\rho}&:= \RelAnc{l}{\rho} 
    \cup \left\{\pi\in \mathfrak{N}: \exists \pi'\in \RelAnc{l}{\rho}\
      \pi\trans{} \pi' \text{ or } \pi\hookrightarrow \pi' 
    \text{ or } \pi\overset{+1}{\hookrightarrow} \pi' \right\}.
  \end{align*}

  If $\bar\rho=(\rho_1, \rho_2, \dots, \rho_n)$ then we write
  $\RelAnc{l}{\bar\rho}:=\bigcup\limits_{i=1}^n \RelAnc{l}{\rho_i}$. 
\end{definition}

\begin{remark}
  Note that for each $\rho'$ there is at most one $\rho$ such that
  $\rho\overset{+1}{\hookrightarrow}\rho'$ while $\rho$ may have arbitrary
  many $\overset{+1}{\hookrightarrow}$ successors along each branch. 

  The relation $\overset{+1}{\hookrightarrow}$ can be characterised as
  follows: For runs $\rho, \rho'$, it holds that
  $\rho\overset{+1}{\hookrightarrow} \rho'$ if and only if 
  $\rho'=\rho\circ \pi$ for  some run $\pi$ starting at some stack $s_\rho$ and
  ending in some stack $s_\rho:w$, the first operation of $\pi$ is a clone
  and $\pi$ visits $s_\rho$ only in its initial configuration. 

  The motivation for these definitions is the following. 
  If there are elements $\rho,\rho'\in \mathfrak{N}$ such that 
  $\rho'\preceq \rho$
  and there 
  is a path in $\mathfrak{N}$ of length at most $l$ that witnesses
  that $\rho'$ is an 
  ancestor of $\rho$, then we want that $\rho'\in \RelAnc{l}{\rho}$. 
  The relation $\overset{+1}{\hookrightarrow}$ is  tailored
  towards this idea. Assume that there are runs $\rho_1 \prec \rho_2
  \trans{\Pop{2}} \rho_3$ such that $\rho_2 \trans{\Pop{2}} \rho_3
  \hookleftarrow \rho_1$. This path of length $2$ witnesses that
  $\rho_1$ is a predecessor of $\rho_2$. By definition, one  sees
  immediately that $\rho_1 \overset{+1}{\hookrightarrow} \rho_2$
  whence
  $\rho_1\in\RelAnc{1}{\rho_2}$. In this sense,
  $\overset{+1}{\hookrightarrow}$ relates the ancestor $\rho_1$ of
  $\rho_2$ with $\rho_2$ if  $\rho_1$ may be reachable from $\rho_2$
  via a short path passing a descendant of $\rho_2$. 

  In the following, it may be helpful to think of a relevant
  $l$-ancestor $\rho'$
  of a run $\rho$ as  an ancestor of $\rho$ that may have a path of
  length up to $l$  witnessing that $\rho'$ is an ancestor of $\rho$. 
  We  do  not state this idea more precisely, but it may be helpful
  to keep this picture in mind. 
\end{remark}
From the definitions, we obtain immediately the following lemmas. 

\begin{lemma}
  Let $\rho$ and $\rho'$ be runs such that $\rho \hookrightarrow
  \rho'$. Let $\hat\rho$ be the predecessor of $\rho'$, i.e.,
  $\hat\rho$ is the unique element such that $\hat\rho \trans{}
  \rho'$. Then $\rho\overset{+1}{\hookrightarrow}\hat\rho$. 
\end{lemma}

\begin{lemma} \label{LemmaRelAnclocalIso}
  If $\rho,\rho'\in \mathfrak{N}$ are connected by a single edge $\trans{}$ or
  $\hookrightarrow$ then 
  either $\rho\in \RelAnc{1}{\rho'}$ or $\rho'\in\RelAnc{1}{\rho}$. 
\end{lemma}

\begin{lemma}
  For all $l\in\N$ and $\rho\in \mathfrak{N}$, 
  $\lvert \RelAnc{l}{\rho} \rvert \leq 4^l$.
\end{lemma}

\begin{lemma}
  $\RelAnc{l}{\rho}$ is linearly ordered by $\preceq$. 
\end{lemma}
\begin{proof}
  By induction, one obtains easily that $\RelAnc{l}{\rho}$ only contains
  initial segments of the run $\rho$. These are obviously ordered
  linearly by $\preceq$.
\end{proof}

In the following we investigate the relationship between relevant
ancestors of different runs.  
First, we characterise the minimal element of $\RelAnc{l}{\rho}$. 
\begin{lemma}\label{LemmaMinRelAnc}
  Let $\rho_l\in \RelAnc{l}{\rho}$ be minimal with respect to $\preceq$. 
  \begin{align*}
    \text{Either }&\lvert \rho_l \rvert = 1 \text{ and }
    \lvert \rho \rvert \leq l, \\  
    \text{or } &\rho_l = \Pop{2}^l(\rho) \text{ and } 
    \lvert \rho_l \rvert < \lvert \rho' \rvert \text{ for all }
    \rho'\in \RelAnc{l}{\rho}\setminus\{\rho_l\}.      
  \end{align*}
\end{lemma}
\begin{remark}
  Recall that $\lvert \rho \rvert \leq l$ implies that
  $\Pop{2}^{l}(\rho)$ is undefined. 
\end{remark}
\begin{proof}
  The proof is by induction on $l$. For $l=0$, there is nothing to
  show because \mbox{$\rho_0 = \rho = \Pop{2}^0(\rho)$.}
  Now assume that the statement is true for some $l$. 

  Assume that $\lvert \rho \rvert \leq l+1$. Then $\rho_l$ satisfies $\lvert
  \rho_l \rvert = 1$. If $\rho_l$ has no
  predecessor it is also the minimal element of $\RelAnc{l+1}{\rho}$ and we are
  done. Otherwise, there is a maximal ancestor  $\hat \rho\prec \rho_l$ such
  that $\lvert \hat\rho \rvert = 1$. Either $\hat \rho\trans{} \rho_l$ or
  $\hat \rho\hookrightarrow \rho_l$ whence 
  $\hat \rho\in \RelAnc{l+1}{\rho}$. Furthermore, no ancestor of $\hat
  \rho$ can 
  be contained in $\RelAnc{l+1}{\rho}$. We prove this claim by
  contradiction. 

  Assume that there is some element $\tilde \rho\prec\hat\rho$ such
  that $\tilde\rho\in\RelAnc{l+1}{\rho}$. Then there is some
  $\tilde\rho'\in\RelAnc{l}{\rho}$ such that $\tilde\rho$ and
  $\tilde\rho'$ are connected by some edge. Due to the definition of
  $\hat\rho$, we have $\tilde\rho \prec \hat\rho \prec
  \tilde\rho'$. Thus, the edge between $\tilde\rho$ and $\tilde\rho'$
  has to be $\hookrightarrow$ or
  $\overset{+1}{\hookrightarrow}$. Thus, 
  $\tilde\rho$ must have width less than $\hat \rho$, i.e., width
  $0$. 
  Since there are no stacks of width $0$, this is a contradiction. 
  
  Thus, the minimal element of $\RelAnc{l+1}{\rho}$ is
  $\rho_{l+1}=\hat \rho$. This completes the case $\lvert \rho \rvert
  \leq l+1$. 

  Now assume that $\lvert \rho \rvert > l+1$. 
  Let $\hat \rho$ be the maximal ancestor of $\rho_l$ such that $\lvert
  \hat \rho \rvert +1 = \lvert \rho_l \rvert$. Then 
  $\hat \rho \overset{+1}\hookrightarrow \rho_l$
  or $\hat\rho \trans{\gamma} \rho_1$, whence 
  $\hat \rho\in \RelAnc{l+1}{\rho}$.  
  We have to show that $\hat \rho$ is the minimal element of
  $\RelAnc{l+1}{\rho}$ and
  that there is no other element of width $\lvert \hat\rho \rvert$ in
  $\RelAnc{l+1}{\rho}$.  
  For the second part, assume that there is some 
  $\rho'\in \RelAnc{l+1}{\rho}$ with
  $\lvert \rho' \rvert = \lvert \hat \rho \rvert$. Then $\rho'$ has to
  be connected via 
  $\trans{}, \overset{+1}{\hookrightarrow}$, or
  $\hookrightarrow$  to some element 
  $\rho''\in \RelAnc{l}{\rho}$. By
  definition of these relations $\lvert \rho'' \rvert \leq \lvert \rho' \rvert
  +1$. By induction hypothesis, this implies $\rho''=\rho_l$. But then it is
  immediately clear that $ \rho' = \hat \rho$ by definition. 

  Similar to the previous case, the minimality of $\hat\rho$ in
  $\RelAnc{l+1}{\rho}$ is proved by 
  contradiction.
  Assume that there is some $\rho'\prec \hat \rho$
  such that $\rho'\in \RelAnc{l+1}{\rho}$. Then there is some 
  \mbox{$\hat \rho\prec \rho_l \preceq \rho''\in 
  \RelAnc{l}{\rho}$} such that 
  \mbox{$\rho'\overset{+1}{\hookrightarrow} \rho''$} or 
  $\rho' \hookrightarrow \rho''$. By the definition of  
  $\hookrightarrow$
  and $\overset{+1}{\hookrightarrow}$, we obtain
  $\lvert  \rho'' \rvert \leq \lvert \hat \rho \rvert$. But 
  this contradicts  
  $\lvert \rho'' \rvert \geq \lvert \rho_l \rvert > \lvert \hat \rho \rvert$. 
  Thus, we conclude that $\hat\rho$ is the minimal element of
  $\RelAnc{l+1}{\rho}$, i.e., $\hat \rho=\rho_{l+1}$.
\end{proof}

The previous lemma shows that the width of stacks among the relevant
ancestors cannot decrease too much. Furthermore, the width cannot grow
too much. This is shown in the following corollary.

\begin{corollary} \label{CorDistRelAnc}
  Let $\pi,\rho\in \mathfrak{N}$ such that
  $\pi\in\RelAnc{l}{\rho}$. Then $\big\lvert 
  \lvert \rho \rvert - \lvert \pi \rvert \big\lvert \leq l$.
\end{corollary}
\begin{proof}
  From the previous lemma, we know that the minimal width of the last
  stack of an element in $\RelAnc{l}{\rho}$ is $\lvert \rho \rvert -l$. 
  We prove by induction that the maximal width is 
  $\lvert \rho \rvert +l$. The case $l=0$ is trivially true. 
  Assume that $\lvert \pi \vert \leq \lvert \rho \rvert +l-1$ for all
  $\pi\in\RelAnc{l-1}{\rho}$.
  Let $\hat\pi\in \RelAnc{l}{\rho}\setminus\RelAnc{l-1}{\rho}$. Then
  there is a 
  $\pi\in\RelAnc{l-1}{\rho}$ such that $\hat\pi\trans{} \pi$,
  $\hat\pi\hookrightarrow \pi$, or $\hat\pi \overset{+1}{\hookrightarrow} \pi$.
  In the last two cases the width of $\hat\pi$ is smaller than the width
  of $\pi$ whence $\lvert \hat\pi \rvert \leq \lvert \rho \rvert +l-1$. 
  In the first case, recall that all stack operations of an level $2$
  higher order pushdown system alter the width of the stack by at most
  $1$. Thus, $\lvert \hat\pi \rvert \leq \lvert \pi \rvert +1 \leq
  \lvert \rho \rvert + l$. 
\end{proof}

The next lemma shows a kind of triangle inequality of the relevant
ancestor relation. If $\rho_2$ is a relevant ancestor of $\rho_1$ then
all relevant ancestors of $\rho_1$ that are prefixes of $\rho_2$ are
relevant ancestors of $\rho_2$.

\begin{lemma} \label{LemmaRelAncComposition}
  Let $\rho_1, \rho_2 \in \mathfrak{N}$ and let $l_1,l_2\in\N$. 
  If $\rho_1\in \RelAnc{l_1}{\rho_2}$, then
  \begin{align*}
    &\RelAnc{l_2}{\rho_1}\subseteq \RelAnc{l_1+l_2}{\rho_2} \text{ and}\\
    &\RelAnc{l_2}{\rho_2}\cap \{\pi: \pi\preceq \rho_1\} \subseteq
    \RelAnc{l_1+l_2}{\rho_1}.
  \end{align*}
\end{lemma}
\begin{proof}
  The first relation holds directly because of the inductive
  definition of relevant ancestors. 

  For the second claim, we proceed by induction on $l_2$. For
  $l_2=0$ the claim holds because  $\RelAnc{0}{\rho_2}=\{\rho_2\}$ and
  $\rho_1\preceq\rho_2$ imply that $\RelAnc{0}{\rho_2}\cap\{\pi:
  \pi\preceq \rho_1\}\neq\emptyset$ if and only if $\rho_1=\rho_2$ and
  $\{\rho_2\}\in\RelAnc{0}{\rho_1}$. 
 
  For the induction step assume that 
  \begin{align*}
    \RelAnc{l_2-1}{\rho_2}\cap \{\pi: \pi\preceq \rho_1\} \subseteq
    \RelAnc{l_1+l_2-1}{\rho_1}.
  \end{align*}
  Furthermore, assume that 
  $\pi\in \RelAnc{l_2}{\rho_2} \cap \{\pi:\pi \preceq \rho_1\}$. 
  We show that $\pi\in\RelAnc{l_1+l_2}{\rho_1}$.
  By
  definition there is 
  some $\pi\prec \hat\pi$ such that $\hat\pi\in
  \RelAnc{l_2-1}{\rho_2}$ and $\pi \in  \RelAnc{1}{\hat\pi}$.
  We distinguish the following cases. 
  \begin{itemize}
  \item   Consider the case $\hat\pi \preceq \rho_1$. Due to the
    induction hypothesis,  
    \mbox{$\hat\pi\in \RelAnc{l_1+l_2-1}{\rho_1}$}. Thus,
    $\pi\in \RelAnc{l_1+l_2}{\rho_1}$.
  \item  Consider the case $\hat\pi = \rho_1$. Then 
    $\pi\in \RelAnc{1}{\rho_1}\subseteq \RelAnc{l_1+l_2}{\rho_1}$.
  \item  Finally, consider the case
    $\pi\prec \rho_1\prec \hat\pi \prec \rho_2$. 
    This implies that $\pi\hookrightarrow\hat\pi$ or
    $\pi\overset{+1}{\hookrightarrow} \hat\pi$ whence 
    $\lvert \pi\rvert = \lvert \hat\pi \rvert -j<
    \lvert \rho_1 \rvert$ for some
    $j\in \{0,1\}$. 
    From Corollary \ref{CorDistRelAnc}, we know that
    \begin{align*}
      \left\rvert \lvert \hat\pi \rvert - \lvert \rho_2 \rvert \right\rvert
      \leq l_2-1 \text{ and}
      \left\rvert \lvert \rho_1 \rvert - \lvert \rho_2 \rvert \right\rvert
      \leq l_1.
    \end{align*}
    This implies that $\lvert \rho_1 \rvert - \lvert \pi \rvert \leq
    l_1+l_2$. 
    By definition of $\hookrightarrow$ and
    $\overset{+1}{\hookrightarrow}$, there cannot be any 
    element \mbox{$\pi\prec \pi' \prec \hat\pi$} with 
    $\lvert \pi' \rvert = \lvert \pi \rvert$.
    Thus, $\pi$ is the maximal predecessor of $\rho_1$ with
    $\pi=\Pop{2}^{\lvert \rho_1 \rvert - \lvert \pi \rvert}(\rho_1)$. 
    Application of Lemma \ref{LemmaMinRelAnc} shows that 
    $\pi$ is the minimal
    element of 
    $\RelAnc{\lvert \rho_1 \rvert - \lvert \pi \rvert}{\rho_1}$. 
    Hence, 
    \begin{align*}
      \pi \in \RelAnc{ \lvert \rho_1 \rvert - \lvert \pi \rvert}{\rho_1}
      \subseteq \RelAnc{l_1+l_2}{\rho_1}.      
    \end{align*}
    \qedhere
  \end{itemize}
\end{proof}

\begin{corollary} \label{CorRelAncDistBound}
  For $\rho\in \RelAnc{l}{\rho_1} \cap \RelAnc{l}{\rho_2}$, we have
  \mbox{$\RelAnc{l}{\rho_1} \cap \{\pi: \pi\preceq \rho\} \subseteq
  \RelAnc{3l}{\rho_2}$.}  
\end{corollary}
\begin{proof}
  By the previous lemma,  $\rho\in \RelAnc{l}{\rho_1}$ implies
  $\RelAnc{l}{\rho_1} \cap \{\pi: \pi\preceq \rho\} \subseteq
  \RelAnc{2l}{\rho}$.  
  Using the lemma again, $\rho\in \RelAnc{l}{\rho_2}$ implies 
  $\RelAnc{2l}{\rho}\subseteq \RelAnc{3l}{\rho_2}$.
\end{proof}

The previous corollary shows that
if the relevant $l$-ancestors of two elements $\rho_1$ and $\rho_2$
intersect at some point $\rho$, then all relevant $l$-ancestors of
$\rho_1$ that are ancestors of $\rho$ are contained in the relevant
$3l$-ancestors of $\rho_2$.  Later, we will use the contraposition of
this result in order to
prove that relevant ancestors of certain runs are disjoint sets. 

The following proposition describes how $\RelAnc{l}{\rho}$
embeds into the full $2$-\HONPT $\mathfrak{N}$. 
Successive relevant ancestors of some run $\rho$ are either connected
by a single edge or by a $\overset{+1}{\hookrightarrow}$-edge.
Later, we will see that this proposition allows to explicitly construct
for any run $\rho$ an isomorphic relevant ancestor set that consists of
small runs.  

\begin{proposition} \label{Prop:NextRelAnc}
  Let $\rho_1\prec \rho_2 \prec \rho$ such that 
  $\rho_1,\rho_2\in  \RelAnc{l}{\rho}$.  
  If $\pi\notin \RelAnc{l}{\rho}$ for all $\rho_1\prec \pi\prec \rho_2$, then
  either $\rho_1\trans{} \rho_2$ or $\rho_1\overset{+1}{\hookrightarrow} \rho_2$.
\end{proposition}
\begin{proof}
  Assume that $\rho_1\not\trans{} \rho_2$. 
  Consider the set
  \begin{align*}
  M:=\{\pi\in\RelAnc{l}{\rho}:\rho_1\overset{+1}{\hookrightarrow} \pi\}.
  \end{align*}
  $M$ is nonempty because there is some $\pi\in\RelAnc{l-1}{\rho}$ such that
  either $\rho_1\overset{+1}{\hookrightarrow} \pi$ (whence $\pi\in M)$ or
  $\rho_1\hookrightarrow \pi$ (whence the predecessor $\hat\pi$ of
  $\pi$ satisfies $\hat\pi\in M$). 
  Let $\hat\rho\in M$ be minimal. It suffices to show that
  $\hat\rho=\rho_2$. For this purpose, we show that  
  $\pi\notin\RelAnc{l}{\rho}$ for all $\rho_1\prec \pi \prec
  \hat\rho$. Since $\hat\rho\in\RelAnc{l}{\rho}$, this implies that
  $\hat\rho=\rho_2$.   

  We start with two general observations.
  \begin{enumerate}
  \item \label{ImmerGroesser}
    For all $\rho_1\prec \pi\prec \hat\rho$,  
    $\lvert \pi \rvert \geq \lvert \hat\rho \rvert$
    due to the definition of $\rho_1\overset{+1}{\hookrightarrow}
    \hat\rho$. Furthermore, due to the minimality of $\hat\rho$ in
    $M$, for all  
    $\rho_1\prec \pi \prec \hat\rho $ with $\pi\in\RelAnc{l}{\rho}$, $\lvert
    \pi\rvert>\lvert \hat\rho \rvert$ (otherwise we have $\pi\in M$
    contradicting the minimality of $\hat\rho$). 
  \item \label{KeineHaken}
    Note that there cannot exist $\rho_1\prec \pi \prec \hat\rho \prec \hat\pi$
    with $\pi \hookrightarrow \hat\pi$ or $\pi \overset{+1}{\hookrightarrow}
    \hat\pi$ because $\lvert \pi \rvert \geq \lvert \hat\rho \rvert$.
  \end{enumerate}
  Heading for a contradiction, assume that there is some 
  $\rho_1\prec \pi \prec \hat\rho$ such that \mbox{$\pi\in\RelAnc{l}{\rho}$.}

  Due to observation \ref{KeineHaken}, there is
  a chain $\pi_0:=\pi, \pi_1, \dots, \pi_{n-1}, \pi_n:=\hat\rho$ such
  that 
  for each $0\leq i<n$ there is $*\in\{\trans{}, \hookrightarrow,
  \overset{+1}{\hookrightarrow}\}$ such that $\pi_i\mathrel{*} \pi_{i+1}$
  and $\pi_i\in \RelAnc{l-i}{\rho}$. By assumption, 
  $n\neq 0$, whence $\hat\rho\in\RelAnc{l-1}{\rho}$. 
  Due to observation \ref{ImmerGroesser}, we have 
  $\lvert \rho_1 \rvert < \lvert \hat\rho \rvert < \lvert \pi
  \rvert$. Since each 
  stack operation alters the width 
  of the stack by at most $1$, we conclude that the set 
  \begin{align*}
    M':=\left\{\pi': \rho_1\prec \pi' \prec \hat\rho, \lvert \hat\rho \rvert =
      \lvert \pi'  \rvert\right\}     
  \end{align*}
  is nonempty because on the path from $\rho_1$ to $\pi$ there
  occurs at least one run with final stack of width $\lvert \hat\rho
  \rvert$. But the maximal element $\pi'\in M'$ satisfies  
  $\rho_1 \overset{+1}{\hookrightarrow} \pi'\trans{} \hat\rho$ or 
  $\rho_1 \overset{+1}{\hookrightarrow} \pi' \hookrightarrow
  \hat\rho$. 
  Since
  $\hat\rho\in\RelAnc{l-1}{\rho}$, this would imply $\pi' \in M$ which
  contradicts the minimality of $\hat\rho$ in $M$. Thus, no 
  $\rho_1\prec \pi 
  \prec \hat\rho$ 
  with $\pi\in\RelAnc{l}{\rho}$ can exist. 

  Thus,
  $\pi\notin\RelAnc{l}{\rho}$ for all
  $\rho_1\prec\pi\prec \hat\rho$ and
  $\rho_1\overset{+1}{\hookrightarrow} \hat\rho=\rho_2$. 
\end{proof}

In the final part of this section, we consider relevant ancestors of
two different runs $\rho$ and $\rho'$. Since we aim at a construction
of small runs $\hat\rho$ and $\hat\rho'$ such that the relevant
ancestors of $\rho$ and $\rho'$ are isomorphic to the relevant
ancestors of $\hat\rho$ and $\hat\rho'$, we need to know how sets of
relevant ancestors touch each other. Every isomorphism from the
relevant ancestors of $\rho$ and $\rho'$ to those of $\hat\rho$ and
$\hat\rho'$ has to preserve edges between a relevant ancestor of
$\rho$ and another one of $\rho'$. 

The positions where the relevant $l$-ancestors of $\rho$ and
$\hat\rho$ touch can  be identified by looking at the intersection of
their relevant \mbox{$(l+1)$-ancestors}. This is shown in the following
Lemma. 
For $A$ and $B$ subsets of some $2$-\HONPT $\mathfrak{N}$ and $\rho$ some
run of $\mathfrak{N}$, we 
say $A$ and $B$ \emph{touch after $\rho$} if there are runs $\rho\prec
\rho_A,\rho\prec\rho_B$ such  
that $\rho_A\in A$, $\rho_B\in B$ and either $\rho_A=\rho_B$ or
$\rho_A * \rho_B$ for some $*\in \{\trans{}, \invtrans{}, \hookrightarrow,
\hookleftarrow\}$. In this case we say $A$ and $B$ touch at 
$(\rho_A,\rho_B)$. 
In the 
following, we reduce the question whether $l$-ancestors of two
elements touch after 
some $\rho$ to the
question whether the $(l+1)$-ancestors of these elements intersect after
$\rho$. 

\begin{lemma}
  If $\rho_1,\rho_2$ are runs such that $\RelAnc{l_1}{\rho_1}$ and 
  $\RelAnc{l_2}{\rho_2}$
  touch after some $\rho_0$, then  
  $\RelAnc{l_1+1}{\rho_1} \cap \RelAnc{l_2+1}{\rho_2}
  \cap \{\pi: \rho_0\preceq \pi\} \neq \emptyset$.
\end{lemma}
\begin{proof}
  Let $\rho_0$ be some run, 
  $\rho_0\prec \hat\rho_1\in \RelAnc{l_1}{\rho_1}$, and $\rho_0\prec
  \hat\rho_2\in 
  \RelAnc{l_2}{\rho_2}$ such that the 
  pair $(\hat\rho_1,\hat\rho_2)$ is minimal and $\RelAnc{l_1}{\rho_1}$
  and $\RelAnc{l_2}{\rho_2}$ 
  touch at $(\hat\rho_1,\hat\rho_2)$. 
  Then one of the following holds.
  \begin{enumerate}
  \item $\hat\rho_1=\hat\rho_2$: there is nothing to prove because
    $\hat\rho_1 \in
    \RelAnc{l_1}{\rho_1}\cap\RelAnc{l_2}{\rho_2}\cap\{\pi:\rho_0\preceq\pi\}$.
  \item $\hat\rho_1\rightarrow \hat\rho_2$ or
    $\hat\rho_1\hookrightarrow \hat\rho_2$ 
    or $\hat\rho_1 \overset{+1}{\hookrightarrow} \hat\rho_2:$ this implies that
    $\hat\rho_1\in \RelAnc{l_2+1}{\rho_2}\cap\RelAnc{l_1}{\rho_1}$.
  \item 
    $\hat\rho_2\rightarrow \hat\rho_1$ or $\hat\rho_2\hookrightarrow \hat\rho_1$
    or $\hat\rho_2 \overset{+1}{\hookrightarrow} \hat\rho_1:$ this
    implies that that 
    $\hat\rho_2\in \RelAnc{l_1+1}{\rho_1}\cap\RelAnc{l_2}{\rho_2}$.     \qedhere
  \end{enumerate}
\end{proof}

\begin{corollary} \label{CorTouch}
  If $\rho$ and $\rho'$ are runs such that $\RelAnc{l_1}{\rho}$ and
  $\RelAnc{l_2}{\rho'}$ touch after some run $\rho_0$ then there 
  exists some $\rho_0\prec \rho_1 \in \RelAnc{l_1+1}{\rho} \cap
  \RelAnc{l_2+1}{\rho'}$ such that
  \begin{align*}
    \RelAnc{l_1+1}{\rho}\cap \{x: x \preceq \rho_1\} \subseteq 
    \RelAnc{l_2 + 2 l_1 + 3}{\rho'}.    
  \end{align*}
\end{corollary}
\begin{proof}
  Use the previous lemma and Lemma \ref{LemmaRelAncComposition}.
\end{proof}

\subsection{A Family of Equivalence Relations on Words and Stacks}
\label{subsec:EquivalenceonWords}

In this section we introduce a family of equivalence relations on
words. The
basic idea is to classify 
words according to the  $\FO{k}$-type of the word model associated
to the word $w$ enriched by information about certain runs between
prefixes of $w$. This additional information 
describes
\begin{enumerate}
\item  the number of possible loops and returns with certain initial
  and final state of each prefix $v\leq
  w$, and 
\item the number of runs from $(q,w)$ to $(q',v)$ for each  
  prefix $v\leq w$ and all pairs $q,q'$ of states.
\end{enumerate}
It turns out that this equivalence has the following property:
if $w$ and $w'$ are equivalent and $\rho$ is a run starting in $(q,w)$
and ending in $(q',w:v)$, 
then there is a run from $(q,w')$
to $(q,w':v')$ such that the loops and returns of $v$ and $v'$
agree. This is important because runs of this 
kind connect consecutive elements of relevant ancestor sets
(cf. Proposition \ref{Prop:NextRelAnc}). 

In order to copy relevant ancestors, we want to apply this kind of
transfer 
property iteratively, e.g., we want to take a run from
$(q_1, w_1)$ via $(q_2,w_1:w_2)$ to $(q_3,w_1:w_2:w_3)$ and translate
it into some run from $(q_1,w_1')$ via $(q_2, w_1':w_2')$ to 
$(q_3,w_1':w_2':w_3')$ such that the loops and returns of $w_3$ and $w_3'$
agree. 
Analogously, we want to take a run creating $n$ new words and
transfer it to a new run starting in another word and creating $n$
words such that the last words agree on their loops and returns. 
If we can do this, then we can transfer the whole set of relevant
ancestors from some run to another one. Using the results of Section
\ref{ChapterLoops}, this allows us to construct isomorphic
relevant ancestors that consist only of short runs. 

The family of equivalence
relations that we define have the following transfer property. Words 
that are equivalent with respect to the
$n$-th relation allow a transfer of runs creating $n$ new words. 
The idea of the definition is as follows. Assume that we have already
defined the 
$(i-1)$-st equivalence relation. 
We  take the
word model of some word $w$ and annotate each prefix of the
word by its equivalence class with respect to the $(i-1)$-st
relation. Then we define two words to be equivalent with respect to
the $i$-th relation if the $\FO{k}$-types of their enriched word models
agree.  

These equivalence relations and the transfer properties that they induce
are an important tool in the next section. There we 
apply them to an arbitrary  set of relevant ancestors $S$ in order to
obtain isomorphic 
copies of the substructure induced by $S$. 
For the next definition, recall that $w_{-n}$ is an abbreviation for
$\Pop{1}^n(w)$.  

\begin{definition}
  Fix a level $2$ pushdown system $\mathcal{N}$. 
  Let $w\in\Sigma^*$ be some word. 
  We are going to define expanded word models
  $\Lin{n}{k}{z}(w)$ by induction on $n$. Note that for $n=0$ the
  structure will be 
  independent of the parameter $k$ but for greater $n$ this parameter
  influences with which kind of information the structure is
  enriched.  
  Let
  $\Lin{0}{k}{z}(w)$  be the expanded word model 
  \begin{align*}
    \Lin{0}{k}{z}(w):=(\{0,1,\dots, \lvert w \rvert -1\},
    \mathrm{succ}, 
    (P_\sigma)_{\sigma\in\Sigma}, (S^j_{q,q'})_{(q,q')\in Q^2, j\leq z}, 
    (R_j)_{j\in J}, (L_j)_{j\in J}, 
    (H_j)_{j\in J})    
  \end{align*}
  such
  that for $0\leq i <\lvert w \rvert$ the following holds.
  \begin{itemize}
  \item $\mathrm{succ}$ and $P_\sigma$ form the standard word model of
    $w$ in reversed order, i.e., $\mathrm{succ}$ is the successor
    relation on the domain and 
    $i\in P_{\sigma}$ if and only if
    $\TOP{1}(w_{-i})=\sigma$, 
  \item $i\in S^j_{q,q'}$, if there are $j$ pairwise distinct runs $\rho_1,
    \dots, \rho_j$ starting in $(q,w)$ and ending in 
    $(q',w_{-i})$ such that for all $1\leq k \leq j$ and
    $0\leq l < \length(\rho_k)$ the stack 
    at $\rho_k(l)$ is not $w_{-i}$.  
  \item The predicates $R_j$ encode at every position $i$ the function
    $\ReturnFunc{z}(w_{-i})$ (cf. Definition \ref{Def:ReturnFunc}). 
  \item The predicates $L_j$ encode at every position $i$ the function
    $\LoopFunc{z}(w_{-i})$  (cf. Definition \ref{Def:LoopFunc}). 
  \item The predicates $H_j$ encode at every position $i$ the function
    $\HighLoopFunc{z}(w_{-i})$. 
  \end{itemize}

  Now, set $\Typ{0}{k}{z}(w):=\FO{k}[\Lin{0}{k}{z}(w)]$, the quantifier
  rank $k$ 
  theory  of $\Lin{0}{k}{z}(w)$. We call it the \emph{$(0,k,z)$-type of
  $w$}.   
  Note that there are only finitely many 
  $(0,k,z)$-types (cf. example \ref{Ex:colouredWordStructures}).

  Inductively, we define $\Lin{n+1}{k}{z}(w)$ to be the expansion of
  $\Lin{n}{k}{z}(w)$ by predicates describing $\Typ{n}{k}{z}(v)$ for
  each prefix $v\leq w$. 
  More formally,  fix a maximal list $\theta_1, \theta_2, \dots,
  \theta_m$ of pairwise distinct $\FO{k}$-types that are realised 
  by some $\Lin{n}{k}{z}(w)$. 
  We define predicates $T_1, T_2, \dots, T_m$ such that 
  $i\in T_j$ if
  $\Typ{n}{k}{z}\left(w_{-i}\right)=\theta_j$ 
  for all $0\leq i\leq n$. Now, let $\Lin{n+1}{k}{z}(w)$ be the
  expansion of $\Lin{n}{k}{z}(w)$ by the 
  predicates $T_1, T_2, \dots, T_m$. 
  We conclude the inductive definition by setting
  $\Typ{n+1}{k}{z}(w):= \FO{k}[\Lin{n+1}{k}{z}(w)]$. 
\end{definition}
\begin{remark}
  Each element of $\Lin{n}{k}{z}(w)$ corresponds to a prefix of
  $w$. In this sense, we write $v\in S^j_{q,q'}$ for
  some prefix $v\leq w$ if $v=w_{-l}$ and 
  $\Lin{n}{k}{z}(w) \models l\in S^j_{q,q'}$.  

  It is an important observation that $\Lin{n}{k}{z}(w)$ is  a finite 
  successor structure with finitely many colours. Thus, there are only
  finitely many $(n,k,z)$-types for each $n,k,z\in\N$
  (cf. Example \ref{Ex:colouredWordStructures}). 
\end{remark}
For our application, $k$ and $z$ can be chosen to be
some fixed large numbers, depending on the number of rounds we are
going to play in the Ehrenfeucht-\Fraisse game. 
Furthermore, it will turn out that the conditions on $k$ and $z$
coincide whence we will 
assume that $k=z$. This is due to the fact that both parameters are 
counting thresholds in some sense: $z$ is the threshold for counting
the existence of loops and returns, while $k$ can be seen as the
threshold for distinguishing different prefixes of $w$ which 
have the same atomic type. Thus, we
identify $k$ and $z$ in  
the following definition of the equivalence relation induced by
$\Typ{n}{k}{z}$. 
\begin{definition}
  For words $w,w'\in\Sigma^*$, we write $ w\wordequiv{n}{z} w'$ if
  $\Typ{n}{z}{z}(w)=\Typ{n}{z}{z}(w')$. 
\end{definition}

As a first step, we want to show that $\wordequiv{n}{z}$ is a right
congruence. We prepare the proof of this fact in the following lemma. 
\begin{lemma}
  Let $n\in\N$, $z\geq 2$ and $\mathcal{N}$ be some pushdown system
  of level $2$. 
  Let $w$ be some word and $\sigma\in\Sigma$ some letter. 
  For each $0\leq i <\lvert w \rvert$, the atomic type of $i$ and of
  $0$ in
  $\Lin{n}{z}{z}(w)$ determines the atomic type of $i+1$ in
  $\Lin{n}{z}{z}(w\sigma)$.
\end{lemma}
\begin{proof}
  Recall that $i\in\Lin{n}{z}{z}(w)$ represents 
  $w_{-i}$ and $i+1\in\Lin{n}{z}{z}(w\sigma)$ represents
  $w\sigma_{-(i+1)}$. 
  Since $w_{-i}=w\sigma_{-(i+1)}$, it follows directly that the two
  elements agree on  
  $(P_\sigma)_{\sigma\in\Sigma}$, $(R_j)_{j\in J}$, $(L_j)_{j\in J}$, and
  $(H_j)_{j\in J}$ and that 
  $w_{-i}\wordequiv{n-1}{z} w\sigma_{-(i+1)}$ (recall that the
  elements in $\Lin{n}{z}{z}(w)$ are coloured by
  $\wordequiv{n-1}{z}$-types).

  We claim that the function $\ReturnFunc{z}(w)$ and the set
  \begin{align*}
    \{(j,q,q')\in\N\times Q\times Q: j\leq z, \Lin{n}{z}{z}(w)\models i\in
    S^j_{q,q'}\}
  \end{align*}
  determine whether $\Lin{n}{z}{z}(w\sigma)\models (i+1)\in
  S^j_{q,q'}$.  Recall that the predicates $S^j_{q,q'}$ in
  $\Lin{n}{z}{z}(w)$ encode at each position $l$ the number of runs $\rho$
  from $(q,w)$ to $(q' , w_{-l})$ that do not pass $w_{-l}$ before
  $\length(\rho)$. 
  We now want to determine the number of runs $\rho$ from $(q,w\sigma)$ to
  $(q', w\sigma_{-(i+1)}) = (q', w_{-i})$ that do not pass $w_{-i}$
  before $\length(\rho)$. 

  It is clear that such a run starts with a high loop from
  $(q,w\sigma)$ to some $(\hat q,w\sigma)$. Then it performs some
  transition of the form $(\hat q, \sigma, \gamma, \hat q', \Pop{1})$
  and then it continues with a run from $(\hat q', w)$ to 
  $(q', w_{-i})$ that do not pass $w_{-i}$ before its last configuration. 
  
  In order to determine whether 
  $\Lin{n}{z}{z}(w\sigma)\models (i+1)\in S^j_{q,q'}$, we have to
  count whether $j$ runs of this form exist. 
  To this end,
  we define the  numbers
  \begin{align*}
    &k_{(\hat q, \hat q')}:=\HighLoopFunc{z}(w\sigma)(q,\hat q),\\
    &j_{(\hat q, \hat q')}:=\lvert \{ (\hat q, \sigma, 
    \gamma, \hat q', \Pop{1})\in \Delta\}\rvert,
    \text{ and}  \\
    &i_{(\hat q, \hat q')}:=\max\{k: \Lin{n}{z}{z}(w) \models
    w_{-i}\in S^{k}_{(\hat q',q')}\}
  \end{align*}
  for each pair $\bar q=(\hat q, \hat q') \in Q^2$. 
  It follows directly that there are
  $\sum\limits_{\bar q\in Q^2} i_{\bar q}j_{\bar q} k_{\bar q}$
  many such runs up to threshold $z$. 
  Note that $j_{\bar q}$ only depends on the pushdown system. 
  Due to Corollary
  \ref{Cor:InductiveComputabilityHighLoops},
  $\HighLoopFunc{z}(w\sigma)$ is determined by $\sigma$ and
  $\ReturnFunc{z}(w)$.
  Thus, $k_{\bar q}$ is determined by the atomic type of $0$ in
  $\Lin{n}{z}{z}(w)$. 
  $i_{\bar q}$ only depends on the atomic type of $i$ in 
  $\Lin{n}{z}{z}(w)$. These observations  complete the proof. 
\end{proof}

\begin{corollary}
  Let $n,z\in\N$ such that $z\geq 2$.  
  Let $w_1$ and $w_2$ be words such that 
  \mbox{$w_1 \wordequiv{n}{z} w_2$.}
  Any strategy of Duplicator in the $z$ round Ehrenfeucht-\Fraisse game on
  $\Lin{n}{z}{z}(w_1)$ and $\Lin{n}{z}{z}(w_2)$ translates directly
  into a strategy of Duplicator
  in the $z$ round Ehrenfeucht-\Fraisse game on
  $\Lin{n}{z}{z}(w_1\sigma){\restriction}_{[1,\lvert w_1\sigma\rvert]}$ and
  $\Lin{n}{z}{z}(w_2\sigma){\restriction}_{[1,\lvert w_2\sigma\rvert]}$.
\end{corollary}
\begin{proof}
  It suffices to note that the existence of Duplicators strategy
  implies that the atomic types of $0$ in 
  $\Lin{n}{z}{z}(w_1)$ and $\Lin{n}{z}{z}(w_2)$ agree. Hence, the
  previous lemma applies. 
  Thus, if the atomic type of $i\in\Lin{n}{z}{z}(w_1)$ and
  $j\in\Lin{n}{z}{z}(w_2)$ agree, then the atomic types of
  $i+1\in\Lin{n}{z}{z}(w_1\sigma)$ and 
  $j+1\in\Lin{n}{z}{z}(w_2\sigma)$ agree. 
  Hence, we can obviously
  translate Duplicator's strategy on 
  $\Lin{n}{z}{z}(w_1)$ and $\Lin{n}{z}{z}(w_2)$ into a strategy on 
  $\Lin{n}{z}{z}(w_1\sigma){\restriction}_{[1,\lvert w_1\sigma\rvert]}$ and
  $\Lin{n}{z}{z}(w_2\sigma){\restriction}_{[1,\lvert w_2\sigma\rvert]}$.
\end{proof}

The previous corollary is the main ingredient for the following
lemma. It states that $\wordequiv{n}{z}$ is a right congruence. 

\begin{lemma} \label{LemmaTypeConcatenation}
  For $z\geq 2$, $\wordequiv{n}{z}$ is a right congruence, i.e., 
  if $\Typ{n}{z}{z}(w_1) = \Typ{n}{z}{z}(w_2)$ for some $z\geq 2$, then
  $\Typ{n}{z}{z}(w_1w) = \Typ{n}{z}{z}(w_2w)$ for all
  $w\in\Sigma^*$. 
\end{lemma}
\begin{proof}
  It is sufficient to prove the claim for $w=\sigma\in\Sigma$.
  The lemma then follows by induction on $\lvert w \rvert$. 
  First observe that 
  \begin{align*}
    &\LoopFunc{z}(w_1\sigma) = \LoopFunc{z}(w_2\sigma),\\ 
    &\HighLoopFunc{z}(w_1\sigma) = \HighLoopFunc{z}(w_2\sigma),\text{ and}\\
    &\ReturnFunc{z}(w_1\sigma) = \ReturnFunc{z}(w_2\sigma),     
  \end{align*}
  because these values are determined by the values of the corresponding
  functions at  $w_1$ and $w_2$ (cf. Propositions
  \ref{Prop:AutomatonForReturns} and
  \ref{Prop:AutomatonForLoops}). These functions agree on $w_1$ and
  $w_2$ because the first elements 
  of $\Lin{n}{z}{z}(w_1)$ and $\Lin{n}{z}{z}(w_2)$ are
  $\FO{2}\subseteq\FO{z}$ definable.

  For $i\in\{1,2\}$, $\Lin{n}{z}{z}(w_i\sigma) \models 0\in S^j_{(q,q')}$ if
  and only if $j=1$  and $q=q'$ because $S^j_{(q,q')}$ counts at position $0$ 
  the runs $\rho$ from $(q,w_i\sigma)$ to $(q',w_i\sigma)$ that do not pass
  $w_i\sigma$ before $\length(\rho)$ and, apparently, this implies
  $\length(\rho)=0$. Since
  $\HighLoopFunc{z}(w_1\sigma)=\HighLoopFunc{z}(w_2\sigma)$, we
  conclude that the atomic types of the first elements of 
  $\Lin{0}{z}{z}(w_1\sigma)$ and of 
  $\Lin{0}{z}{z}(w_2\sigma)$ coincide. 
  
  Due to the previous corollary, we know that Duplicator has a
  strategy in the $z$ round Ehrenfeucht-\Fraisse game on
  $\Lin{n}{z}{z}(w_1\sigma){\restriction}_{[1,\lvert w_1\rvert]}$ and
  $\Lin{n}{z}{z}(w_2\sigma){\restriction}_{[1,\lvert w_2\rvert]}$.

  Standard composition arguments for Ehrenfeucht-\Fraisse games on
  word structures directly imply that
  $\Lin{0}{z}{z}(w_1\sigma) \simeq_z \Lin{0}{z}{z}(w_2\sigma)$. But
  this directly implies that the atomic types of the first elements of 
  $\Lin{1}{z}{z}(w_1\sigma)$ and of 
  $\Lin{1}{z}{z}(w_2\sigma)$ coincide. If $n\geq 1$, we can apply the
  same standard argument and obtain that
  $\Lin{1}{z}{z}(w_1\sigma) \simeq_z \Lin{1}{z}{z}(w_2\sigma)$.
  By induction one concludes that
  $\Lin{n}{z}{z}(w_1\sigma) \simeq_z \Lin{n}{z}{z}(w_2\sigma)$.
  But this is the definition of 
  $w_1\sigma \wordequiv{n}{z} w_2\sigma$. 
\end{proof}

The next lemma can be seen as the inverse direction of the previous
lemma. Instead of appending a word, we want to remove the topmost
symbols from the word. For this operation,  we cannot preserve the
equivalence at the same level but at one level below.

\begin{lemma} \label{Lem:EquivalencesDownwardCompatible}
  Let $m< 2^{z-1}-1$ and $w, w'\in\Sigma^*$. 
  If $w\wordequiv{n}{z} w'$ then 
  $w_{-m} \wordequiv{n-1}{z} w'_{-m}$.
\end{lemma}
\begin{proof}
  Quantifier rank $z$ suffices to define the $m$-th element of a
  word structure. Hence, $w\wordequiv{n}{z} w'$ implies that 
  $\Typ{n-1}{z}{z}(w_{-m}) =   \Typ{n-1}{z}{z}(w'_{-m})$. But this is
  equivalent to \mbox{$w_{-m} \wordequiv{n-1}{z} w'_{-m}$.}
\end{proof}

The previous lemmas can be seen as statements concerning the
compatibility of the stack operations $\Push{\sigma}$ and $\Pop{1}$ 
with the equivalences $\wordequiv{n}{z}$. 
Later, we need a  compatibility result of the equivalences with all
level $2$ stack operations. For this purpose, we
first lift these
equivalences to equivalences on level $2$ stacks. We compare the
stacks word-wise beginning with the topmost word, then the word below
the topmost one, etc. up to some threshold $m$. 
The following
definition introduces the precise notion of these equivalence
relations on stacks. 

\begin{definition}
  Let $s,s'$ be stacks. We write $s \stackequiv{n}{z}{m} s'$ if
  for all $0\leq i \leq m$
  \begin{align*}
    \TOP{2}\left(\Pop{2}^i(s)\right)\wordequiv{n}{z}
    \TOP{2}\left(\Pop{2}^i(s')\right).
  \end{align*}
\end{definition}
\begin{remark}
  If $\width(s) \leq m$ or $\width(s')\leq m$ then $\Pop{2}^m(s)$ or
  $\Pop{2}^m(s')$ is undefined. In this case we write $s
  \stackequiv{n}{z}{m} s'$ iff $\width(s)=\width(s')$ and
  $s \stackequiv{n}{z}{m'} s'$ for $m':=\width(s)-1$. 
\end{remark}

Next, we prove that these equivalence relations on stacks are
compatible with all stack operations. 

\begin{proposition} \label{Prop:CompatibilityStackOpTypeq}
  Let $z\geq 2$ and let $s_1, s_2, s_1', s_2'$ be stacks such that
  $s_1' = \op(s_1)$ and \mbox{$s_2' = \op(s_2)$} for some stack
  operation $\op$.   
  If $s_1 \stackequiv{n}{z}{m} s_2$ then the following hold:
  \begin{itemize}
  \item for $\op=\Push{\sigma}$, $s_1'\stackequiv{n}{z}{m} s_2'$,
  \item for $\op=\Pop{1}$, $s_1' \stackequiv{n-1}{z}{m} s_2'$,
  \item for $\op=\Clone{2}$, $s_1' \stackequiv{n}{z}{m+1} s_2'$,
    and
  \item for $\op=\Pop{2}$, $s_1' \stackequiv{n}{z}{m-1} s_2'$. 
  \end{itemize}
\end{proposition}
\begin{proof}
  For $\op=\Pop{1}$, we use Lemma
  \ref{Lem:EquivalencesDownwardCompatible}. For 
  $\op=\Push{\sigma}$ we use Lemma \ref{LemmaTypeConcatenation}. 
  For $\Clone{2}$ and $\Pop{2}$, the claim follows directly from the
  definitions.  
\end{proof}

The previous proposition shows that the equivalence relations on
stacks are compatible with the stack operations. 
Recall that successive relevant ancestors of a given run $\rho$ are
runs $\rho_1 \prec \rho_2 \preceq \rho$ such that $\rho_2$ extends
$\rho_1$ by either a single transition or by some run that creates
some new word on top of the last stack of $\rho_1$ 
(cf. Proposition \ref{Prop:NextRelAnc}). 
In the next section, we are concerned with the construction of a
short run $\hat \rho$ such that its relevant ancestors are isomorphic
to those of $\rho$. A necessary condition for a run $\hat\rho$ to be
short is that it only passes small stacks. We
construct $\hat\rho$ using the following construction. Let 
$\rho_0 \prec \rho_1 \prec \rho_2 \ldots \prec \rho$ be the set of
relevant ancestors of $\rho$. 
We then first define a run $\hat\rho_0$ that ends in some small stack
that is equivalent to the last stack of $\rho_0$. Then, we iterate the
following construction. If $\rho_{i+1}$ extends $\rho_i$ by a single
transition, then we define $\hat\rho_{i+1}$ to be the extension of
$\hat\rho_{i}$ by the same transition. Due to the previous proposition
this preserves equivalence of the topmost stacks of $\rho_{i}$ and
$\hat\rho_{i}$. Otherwise, $\rho_{i+1}$ extends $\rho_i$ by some run
that creates a new word $w_{i+1}$ on top of the last stack of
$\rho_i$. Then we want to construct a short run that creates a new
word $w_{i+1}'$ on top of the last stack of $\hat\rho_i$ such that
$w_{i+1}$ and $w_{i+1}'$ are equivalent and $w_{i+1}'$ is small. Then
we define
$\hat\rho_{i+1}$ to be $\hat\rho_i$ extended by this run. 

Finally, this procedure defines a run $\hat\rho$ that corresponds
to $\rho$ in the sense that the relevant ancestors of the two runs are
isomorphic but $\hat\rho$ is a short run.

In the following, we prepare this construction. We 
show that for any run $\rho_0$ there is a run $\hat\rho_0$ that ends
in some small stack that is equivalent to the last stack of $\rho_0$. 
This is done in Corollary \ref{Cor:RunsendInShorStacks}. 
Furthermore, we show that for runs $\rho_i$ and $\hat\rho_i$ that end
in equivalent stacks, any run that extends the last stack of $\rho_i$
by some word $w$ can be transferred into a run that extends $\hat\rho_i$ by
some small word that is equivalent to $w$. This is shown in
Proposition \ref{Prop:ConstructOneStep}. 

The proofs of Corollary \ref{Cor:RunsendInShorStacks} and 
Proposition \ref{Prop:ConstructOneStep}
are based on the property that 
prefixes of equivalent stacks share the same number of loops and
returns for each pair of initial and final states. 
Recall that our analysis of generalised milestones showed that the
existence of loops with certain initial and final states has a crucial
influence on the question whether runs between certain stacks exist. 

In the following, we first state three main lemmas 
concerning the reachability of small stacks that are equivalent to
some given stack. Together, these lemmas directly imply the
Corollary \ref{Cor:RunsendInShorStacks}. 
Afterwards, we present the Proposition \ref{Prop:ConstructOneStep}.  
In the end of this section, we provide the technical details for the
proofs of the main lemmas and the proposition. 

The first lemma allows to translate an arbitrary run $\rho$ into
another run $\rho'$ that ends in a stack with a small topmost word
such that the topmost words of $\rho$ and $\rho'$ are equivalent. 
We first define a function that is used to define what small means in
this context. 
\begin{definition}
  Let $\mathcal{N}=(Q,\Sigma,\Gamma,\Delta,q_0)$ be a pushdown system
  of level $2$. 
  Set
  \begin{align*} 
    \FuncBoundTopWord:\N^2 &\rightarrow \N\\
    \FuncBoundTopWord(n,z) &= \lvert Q \rvert \cdot 
    \lvert \nicefrac{\Sigma^*}{\wordequiv{n}{z}} \rvert+1,
  \end{align*}  
  where $\lvert \nicefrac{\Sigma^*}{\wordequiv{n}{z}} \rvert$ is the
  number of equivalence classes of $\wordequiv{n}{z}$. 
\end{definition}

\begin{lemma}  \label{lemmaPumpingTOP}
  For all $z,n\in\N$ with $z\geq 2$ and for each run $\rho$ with
  $\lvert \TOP{2}(\rho) \rvert > \FuncBoundTopWord(n,z)$ 
  there is some run $\hat\rho$ with 
  \begin{align*}
    &\lvert\TOP{2}(\rho)\rvert - \FuncBoundTopWord(n,z) \leq \lvert
    \TOP{2}(\hat\rho) \rvert < 
    \lvert \TOP{2}(\rho) \rvert\text{ and}\\
    &\TOP{2}(\rho) \wordequiv{n}{z} \TOP{2}(\hat\rho).     
  \end{align*}
\end{lemma}

The previous lemma gives the possibility to replace a given run by
some run that ends in an equivalent but small topmost word. 
After bounding the topmost word, we want to bound the height of all
the words occurring in the last configuration of some run $\rho$. 
This is done with the next lemma. 

\begin{definition}
  Let $\mathcal{N}=(Q,\Sigma,\Gamma,\Delta,q_0)$ be a pushdown system
  of level $2$. 
  Set
  \begin{align*}
    \ConstBoundHeightWord :=\lvert
    \nicefrac{\Sigma^*}{\wordequiv{0}{2}} \rvert 
    \cdot \lvert Q \rvert^2.    
  \end{align*}
\end{definition}

\begin{lemma}  \label{lemmaPumpingHeight}
  If $\rho$ is some run with $\height(\rho) > \lvert \TOP{2}(\rho)\rvert +
  \ConstBoundHeightWord $, then there is a run $\hat\rho$ with 
  \begin{align*}
    &\height(\rho)-\ConstBoundHeightWord \leq 
    \height(\hat\rho) < \height(\rho) 
     \text{ and}\\
    &\TOP{2}(\rho) = \TOP{2}(\hat\rho).
  \end{align*}
\end{lemma}

Finally, we want to bound the width of the last stack of some run in
terms of its height while preserving the topmost word. 
This is done in the following lemma. 

\begin{definition}
  Set
  \begin{align*}
    \FuncBoundWidthWord:\N&\rightarrow \N \\
    n &\mapsto \lvert Q \rvert \cdot (\lvert\Sigma\rvert+1)^n.
  \end{align*}  
\end{definition}
\begin{remark}
  $\FuncBoundWidthWord(n)$ is an upper bound for the number of pairs
  of states and words of length up to $n$. 
\end{remark}
\begin{lemma} \label{lemmaPumpingWidth}
  For every run  $\rho$ with $\width(\rho) >
  \FuncBoundWidthWord(\height(\rho))$ there is 
  a run $\hat\rho$ with 
  \begin{align*}
    &\width(\rho)-
    \FuncBoundWidthWord\left(\height(\rho)\right) \leq
    \width(\hat\rho) < \width(\rho) \text{ and}\\
    &\TOP{2}(\hat\rho)=\TOP{2}(\rho).
  \end{align*}
\end{lemma}

The previous three lemmas are summarised in the
following corollary. It asserts that for every run there is a run ending
in a small stack with equivalent topmost word.  

\begin{corollary} \label{Cor:RunsendInShorStacks}
  For each run $\rho$ starting in the initial configuration, there is
  a run $\rho'$ starting in the initial configuration such that 
  \begin{align*}
    &\lvert \TOP{2}(\rho') \rvert \leq \FuncBoundTopWord(n,z), \\
    &\height(\rho') \leq \lvert \TOP{2}(\rho')\rvert + 
    \ConstBoundHeightWord,\\
    &\width(\rho') \leq \FuncBoundWidthWord(\height(\rho')) \text{
      and}\\
    & \TOP{2}(\rho) \wordequiv{n}{z} \TOP{2}(\rho'). 
  \end{align*}
\end{corollary}

The previous corollary deals with the reachability of some stack from the
initial configuration. The following proposition is concerned with the
extension of a given stack by just one word. We first define the
function that is used to bound the size of the new word. 

\begin{definition}
  Let $\mathcal{N}$ be a level $2$ pushdown system with state set
  $Q$. 
  Set
  \begin{align*}
    \BoundHeightOnestepConstructionSimultanious:\N^4 &\rightarrow \N\\
      (a,b,c,d)&\mapsto b+ a( \lvert Q \rvert \lvert
      \nicefrac{\Sigma^*}{\wordequiv{c}{d}} \rvert).
  \end{align*}
\end{definition}

Before we state the proposition, we explain its meaning. The
proposition says that given two equivalent words $w$ and $\hat w$ and a run
$\rho$ from $(q, s:w)$ to $(q', s:w:w')$ that does not pass a substack
of $s:w$, then, for each stack $\hat s:\hat w$, we find a run $\hat
\rho$ from  $(q, \hat s: \hat w)$ to $(q', \hat s: \hat w: \hat w')$
for some short word $\hat w'$ that is equivalent to $w'$. Furthermore, 
this transfer of runs works simultaneously on a tuple of such runs,
i.e., given $m$ runs starting at $s:w$ of the form described above, we
find $m$ corresponding runs starting at $\hat s :\hat w$. 
This simultaneous transfer becomes important when we search 
an isomorphic copy of the relevant ancestors of 
several runs. In this case the simultaneous transfer allows to copy
the relevant ancestors of a certain run while avoiding an intersection
with relevant ancestors of other given runs.

\begin{proposition} \label{Prop:ConstructOneStep}
  Let $\mathcal{N}$ be a level $2$ pushdown system and
  $n,z,m\in \N$ such that $n\geq 1$, $z> m$, and $z\geq
  2$. 
  Let $c=(q,s:w),\hat c=(q,\hat s:\hat w)$ be configurations such that
  $w \wordequiv{n}{z} \hat w$. 
  Let $\rho_1, \dots, \rho_m$ be pairwise distinct runs such that for
  each $i$, $\lvert \rho_i(j) \rvert > \lvert s:w \rvert$ for all
  $j\geq 1$ and such that $\rho_i$
  starts at $c$ and ends in $(q_i, s:w:w_i)$. 
  Analogously, let $\hat \rho_1, \dots, \hat \rho_{m-1}$ be pairwise
  distinct runs such 
  that each $\hat \rho_i$ starts at $\hat c$ and ends in $(q_i, \hat s:\hat
  w:\hat w_i)$ and 
  $\lvert \hat \rho_i(j) \rvert > \lvert \hat s:\hat w \rvert$ for all
  $j\geq 1$. 
  If 
  \begin{align*}
  &w_i\wordequiv{n-1}{z} \hat w_i \text{ for all }1\leq i
  \leq m-1,   
  \end{align*}
  then there is some run $\hat \rho_m$ from $\hat c$ to $(q_0,\hat
  s:\hat w:\hat w_m)$ such
  that 
  \begin{align*}
    & w_m \wordequiv{n-1}{z} \hat w_m, \\ 
    &\hat \rho_m \text{ is distinct from each } \hat\rho_i \text{ for
    } 1\leq i < 
    m \text{, and}\\
    &\lvert \hat w_m\rvert \leq 
    \BoundHeightOnestepConstructionSimultanious(m,\lvert \hat w \rvert, n, z).
  \end{align*}
\end{proposition}

The rest of this section is concerned with the proofs of Lemmas 
\ref{lemmaPumpingWidth}, \ref{lemmaPumpingHeight}, and
\ref{lemmaPumpingTOP} and with the proof of Proposition
\ref{Prop:ConstructOneStep}. The reader who is not interested in the
technical details of these proofs may  skip the rest of this
section and continue reading Section
\ref{SectionEquivalenceonTuples}. In that section show how the
results of this section can be used to construct isomorphic relevant
ancestors that consist of runs ending in small stacks.

\paragraph{Prefix Replacement Revisited}
Recall that we defined the prefix replacement for  runs that are
prefixed by a certain stack
(cf. Lemma \ref{Lem:BlumensathHOLevel2}).
We want to extend the notion of prefix replacement to
runs that are only prefixed at the beginning and at the end by
some stack $s$ and that never visit the substack $\Pop{2}(s)$. 
We apply this new form of prefix replacement in the
 proofs of Lemmas \ref{lemmaPumpingWidth},
\ref{lemmaPumpingHeight} and \ref{lemmaPumpingTOP}.
The following lemma prepares this new kind of prefix replacement. 

\begin{lemma}
  Let $\mathcal{N}$ be some level $2$ pushdown system and let
  $\rho$ be a run of $\mathcal{N}$ of length $n$. 
  Let $s$ be a stack with topmost word $w:=\TOP{2}(s)$ such that
  \begin{enumerate}
  \item $s\prefixeq\rho(0)$,
  \item $s\prefixeq\rho(n)$, and
  \item $\lvert s \rvert \leq \lvert \rho(i) \rvert$ for all $0\leq i
    \leq n$. 
  \end{enumerate}
  There is a unique sequence $0=i_0 \leq j_0 < i_1 \leq j_1 < \dots <
  i_{m-1} \leq  j_{m-1} < i_m \leq j_m=n$ such that
  \begin{enumerate}
  \item $s\prefixeq \rho{\restriction}_{[i_k,j_k]}$ for all $0\leq k
    \leq m$ and 
  \item $\TOP{2}(\rho(j_k+1))=\Pop{1}(w)$,
    $\rho{\restriction}_{[j_k,i_{k+1}]}$ is either a loop or a return,
    and $\rho{\restriction}_{[j_k,i_{k+1}]}$ does not visit the stack of
    $\rho(j_k)$ between its initial configuration and its final
    configuration for all $0\leq k < m$.
  \end{enumerate}
\end{lemma}
\begin{proof}
  If $s\prefixeq\rho$, then we set $m:=0$ and we are done. 
  Otherwise, we proceed by induction on the length of $\rho$. 

  There is a minimal position $j_0+1$ such that
  $s\notprefixeq\rho(j_0+1)$. By assumption on $s$,
  \mbox{$\rho(j_0+1)\neq\Pop{2}(s)$}. Thus,
  $\TOP{2}(\rho(j_0))=w$ and  
  $\TOP{2}(\rho(j_0+1))=\Pop{1}(w)$. 
  Now, let $i_1>j_0$ be minimal such that $s\prefixeq\rho(i_1)$.
  Concerning the stack at $i_1$ there are the following possibilities. 
  \begin{enumerate}
  \item 
    If $\rho(i_1)=\Pop{2}(\rho(j_0))$
    then $s\prefixeq\rho(i_1)$ (cf. Lemma
    \ref{Lem:PrefixPop1Pop2wieder}). Furthermore, 
    $\rho{\restriction}_{[j_0,i_1]}$ is a return.
  \item 
    Otherwise, 
    the stacks of $\rho(j_0)$ and $\rho(i_1)$ coincide whence
    $\rho{\restriction}_{[j_0,i_1]}$ is a loop (note that between $j_0$
    and $i_1$ the stack $\Pop{2}(\rho(j_0))$ is never visited due to the
    minimality of $i_1$ and due to assumption 3). 
  \end{enumerate}
  $\rho{\restriction}_{[i_1,n]}$ is shorter than $\rho$. Thus, it
  decomposes by induction hypothesis and 
  the lemma follows immediately. 
\end{proof}

This lemma gives rise to the following extension of the prefix
replacement. 
\begin{definition}\label{Def:StackreplacementinMilestones}
  Let $s$ be some stack and $\rho$ be a run of some pushdown system
  $\mathcal{N}$ of level $2$ such that 
  $s\prefixeq \rho(0)$, $s\prefixeq\rho(\length(\rho))$ and
  $\lvert s \rvert \leq \lvert \rho(i) \rvert$ for all
  $i\in\domain(\rho)$. 
  Let $u$ be 
  some stack such that $\TOP{1}(u)=\TOP{1}(s)$, 
  $\LoopFunc{1}(u) = \LoopFunc{1}(s)$ and
  $\ReturnFunc{1}(u)=\ReturnFunc{1}(s)$.

  Let $0=i_0 \leq j_0 < i_1 \leq j_1 < \dots <
  i_{m-1} \leq  j_{m-1} < i_m \leq j_m=\length(\rho)$ be the sequence
  corresponding to $\rho$
  in the sense of the previous lemma. 
  We set $(q_k,s_k):=\rho(j_k)$ and $(q'_k,s'_k):=\rho(i_{k+1})$. 
  By definition,
  $\rho{\restriction}_{[j_k,i_{k+1}]}$ is a loop or a return 
  from $(q_k,s_k)$ to $(q_k', s_k')$ and
  $\TOP{2}(s_k)=\TOP{2}(s)$ and $s\prefixeq s_k$. 
  Thus, $\TOP{2}(s_k[s/u])=\TOP{2}(u)$. Since
  \mbox{$\ReturnFunc{1}(u)=\ReturnFunc{1}(s)$} and
  \mbox{$\LoopFunc{1}(u) = \LoopFunc{1}(s)$}, there is a run from
  $(q_k,s_k[s/u])$ to $(q'_{k+1}, s'_{k+1}[s/u])$. 
  We set $\rho_k$ to be the lexicographically shortest run from
  $(q_k,s_k[s/u])$ to $(q'_{k+1}, s'_{k+1}[s/u])$. 

  Then we define the run 
  \begin{align*}
    \rho[s/u]:=\rho{\restriction}_{[i_0,j_0]}[s/u] \circ \rho_0
    \circ \rho{\restriction}_{[i_1,j_1]}[s/u] \circ \rho_1 \circ \dots
    \circ \rho_{m-1} \circ \rho{\restriction}_{[i_{m},j_{m}]}[s/u].  
  \end{align*}
  Note that $\rho[s/u]$ is a well-defined run from $\rho(0)[s/u]$ to
  $\rho(\length(\rho))[s/u]$. 
\end{definition}

\paragraph{Proof of Lemma \ref{lemmaPumpingTOP}}

Recall that Lemma \ref{lemmaPumpingTOP} asserts for every run $\rho$
the existence of a
run $\rho'$ that ends in a stack with small topmost word that is
equivalent to the topmost word of the last stack of $\rho$. 
The proof of this lemma is as follows. 

\begin{figure}
  \centering
  $
  \begin{xy}
    \xymatrix@=0.2mm{
      &      &      &  f   &   e\\
      &  b   &   d  &  e   &   e \\
      &  a   &   a  &  a   &   a \\
      s= & \bot & \bot & \bot & \bot
    }
  \end{xy}
  $\hskip 3mm
  $
  \begin{xy}
    \xymatrix@=0.2mm{
      &\textcolor{white}{ f}\\
      &\textcolor{white}{ d}\\
      &\textcolor{white}{ a}\\
      m_1= & \bot &
    }
  \end{xy}
  $\hskip 3mm
  $
  \begin{xy}
    \xymatrix@=0.2mm{
      \textcolor{white}{ f}\\
      &  b   &   d  &   \\
      &  a   &   a  &  a    \\
      m_2= & \bot & \bot & \bot 
    }
  \end{xy}
  $\hskip 3mm
  $
  \begin{xy}
    \xymatrix@=0.2mm{
      &      &      &  \textcolor{white}{ f}   \\
      &  b   &   d  &  e    \\
      &  a   &   a  &  a    \\
      m_3= & \bot & \bot & \bot 
    }
  \end{xy}
  $
  \caption{Illustration for the construction in the  proof of Lemma
    \ref{lemmaPumpingTOP}.} 
  \label{fig:LemmaPumpingTopExample}
\end{figure}
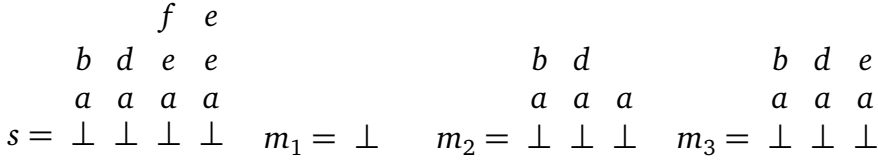

\begin{proof}[Proof of Lemma \ref{lemmaPumpingTOP}.]
  Let $\rho$ be some run with \mbox{$\lvert \TOP{2}(\rho)\rvert >
  \FuncBoundTopWord(n,z)$.} 
  \mbox{$(q,s):=\rho(\length(\rho))$} denotes the final configuration of
  $\rho$. For each $k\leq \FuncBoundTopWord(n,z)$, there is a
  maximal milestone $m_k\in\Milestones(s)$ with
  $\lvert\TOP{2}(m_k)\rvert = k$. 
  Figure \ref{fig:LemmaPumpingTopExample} illustrates this
  definition. 
  Let $w_k:=\TOP{2}(m_k)$ and let $\rho_k
  \preceq \rho$ be the largest initial segment of $\rho$ that ends in
  $m_k$. Note that $m_k\prefixeq m_{k'} \prefixeq s$ for all $k\leq k'
  \leq \FuncBoundTopWord(n,z)$  by the maximality of $m_k$ and
  $m_{k'}$.  
  
  Then there are $i < j \leq \FuncBoundTopWord(n,z)$ such that
  $\TOP{2}(\rho_i) \wordequiv{n}{z} \TOP{2}(\rho_j)$ and the
  final states of $\rho_i$ and $\rho_j$ agree. 
  
  Due to the maximality of $\rho_j$, no substack of $\Pop{2}(m_j)$ is
  visited by $\rho$ after $k:=\length(\rho_j)$. 
  Thus, the run
  $\pi:=(\rho{\restriction}_{[k,\length(\rho)]})[m_j/m_i]$ is
  well-defined   (cf. Definition
  \ref{Def:StackreplacementinMilestones}). Note that
  $\pi$ starts by definition in $(q',m_i)$ for $q'\in Q$ the final
  state of $\rho_i$. Thus, we can set $\hat\rho:=\rho_i\circ \pi$. 
  Due to $w_i\wordequiv{n}{z} w_k$ and the right congruence of
  $\wordequiv{n}{z}$ (cf.   Lemma \ref{LemmaTypeConcatenation}.), it
  is clear that 
  $\TOP{2}(\hat\rho)\wordequiv{n}{z}\TOP{2}(\rho)$.  
  Since $0< \lvert w_j \rvert - \lvert w_i \rvert <
  \FuncBoundTopWord(n,z)$, it also follows directly that
  \begin{align*}
    \lvert\TOP{2}(\rho)\rvert - \FuncBoundTopWord(n,z) \leq \lvert
    \TOP{2}(\hat\rho) \rvert <  \lvert \TOP{2}(\rho)\rvert. 
  \end{align*}
\end{proof}

\paragraph{Proof of Lemma \ref{lemmaPumpingHeight}}

Recall that Lemma \ref{lemmaPumpingHeight} asserts that for each run
$\rho$ there is a run $\rho'$ such that
\mbox{$\TOP{2}(\rho)=\TOP{2}(\rho')$} and such that the height of the last
stack of $\rho'$ is bounded in terms of $\lvert \TOP{2}(\rho)\rvert$. 

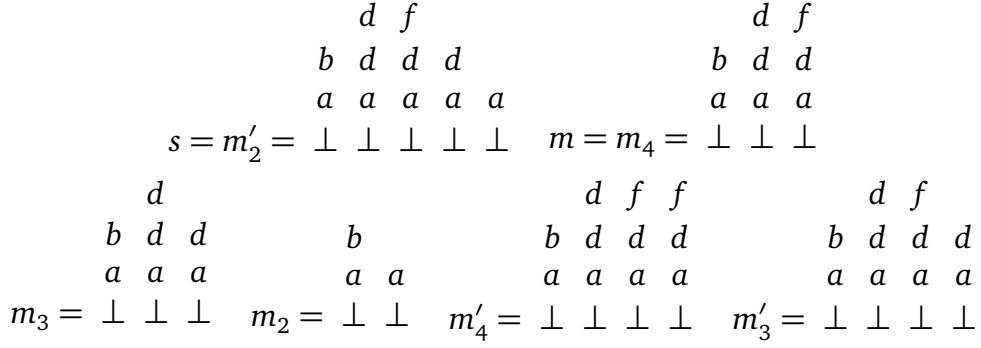
\begin{figure}
  \centering
  $
  \begin{xy}
    \xymatrix@=0.2mm{
               &      &   d  &  f   &   \\
               &  b   &   d  &  d   &   d \\
               &  a   &   a  &  a   &   a &   a\\
      s=m_2' = & \bot & \bot & \bot & \bot & \bot
    }
  \end{xy}
  $\hskip 3mm
  $
  \begin{xy}
    \xymatrix@=0.2mm{
             &         &   d  &  f   \\
             &  b      &   d  &  d    \\
             &  a      &   a  &  a   \\
      m=m_4= & \bot & \bot & \bot
    }
  \end{xy}
  $\\
  $
  \begin{xy}
    \xymatrix@=0.2mm{
           &      &   d  &        \\
           &  b   &   d  &  d       \\
           &  a   &   a  &  a     \\
      m_3= & \bot & \bot & \bot 
    }
  \end{xy}
  $\hskip 3mm
  $
  \begin{xy}
    \xymatrix@=0.2mm{
         & \textcolor{white}{f}   &      \\
         &  b   &    \\
         &  a   &   a  \\
      m_2= & \bot & \bot 
    }
  \end{xy}
  $\hskip 3mm
  $
  \begin{xy}
    \xymatrix@=0.2mm{
            &      &   d  &  f   &  f \\
            &  b   &   d  &  d   &   d \\
            &  a   &   a  &  a   &   a \\
      m'_4= & \bot & \bot & \bot & \bot
    }
  \end{xy}
  $\hskip 3mm
  $
  \begin{xy}
    \xymatrix@=0.2mm{
            &      &   d  &  f   &   \\
            &  b   &   d  &  d   &   d \\
            &  a   &   a  &  a   &   a \\
      m_3'= & \bot & \bot & \bot & \bot
    }
  \end{xy}
  $\hskip 3mm
  \caption{Illustration for the construction in the  proof of Lemma
    \ref{lemmaPumpingHeight}.} 
  \label{fig:LemmaPumpingHeightExample}
\end{figure}

\begin{proof}[Proof of Lemma \ref{lemmaPumpingHeight}]
  The proof is by induction on the
  number of words in 
  the last stack of $\rho$ that have length $h:=\height(\rho)$. 
  Assume that $\rho$ is some run such that
  \begin{align*}
    \height(\rho) > \lvert \TOP{2}(\rho)\rvert +
    \ConstBoundHeightWord.    
  \end{align*}
  In the following, we define several generalised milestones of the final
  stack $s$ of 
  $\rho$. An illustration of  these definitions can be found in 
  Figure \ref{fig:LemmaPumpingHeightExample}.

  Let $m\in \Milestones(s)$ be a milestone of the last stack of
  $\rho$ such that
  $\lvert \TOP{2}(m)\rvert = h$. 
  For each $\lvert \TOP{2}(\rho) \rvert \leq i \leq h$ let
  $m_i\in\Milestones(m)$ be the maximal milestone of $m$ with $\lvert
  \TOP{2}(m_i)\rvert = i$. 
  Let $n_i$ be maximal such that $\rho(n_i)= (q', m_i)$ for some $q'\in
  Q$. 
  Let $m'_i\in \genMilestones(s)\setminus \genMilestones(m)$ be the
  minimal generalised milestone after $m$ such that
  \mbox{$\TOP{2}(m'_i) = \TOP{2}(m_i)$.} 
  Let $n'_i$ be maximal with $\rho(n'_i) = (q', m'_i)$ for some $q'\in
  Q$. 

  There are 
  $\lvert \TOP{2}(\rho) \rvert \leq k < l \leq \height(\rho)$
  satisfying the following conditions.
  \begin{enumerate}
  \item There is a $q\in Q$ such that $\rho(n_k)=(q,m_k)$ and
    $\rho(n_l) = (q, m_l)$.
  \item There is a $q'\in Q$ such that $\rho(n'_k)= (q', m'_k)$
    and $\rho(n'_l) = (q', m'_l)$.
  \item $\TOP{2}(m_k)  \wordequiv{0}{2} \TOP{2}(m_l)$ (this assumption
    implies that
     \mbox{$\LoopFunc{1}(m_k)=\LoopFunc{1}(m_l)$} and
    \mbox{$\ReturnFunc{1}(m_k)= \ReturnFunc{1}(m_l)$}). 
  \end{enumerate}
  By definition, we have $m_l\prefixeq m_l'$. 
  Thus, the run $\pi_1:= (\rho{\restriction}_{[n_l,n'_l]})[m_l/m_k]$ is
  well defined (cf. Definition
  \ref{Def:StackreplacementinMilestones}).
  Note that $\pi_1$ starts in 
  $(q,m_k)$ and ends in $(q',\hat s)$ for $\hat s:=m'_l[m_l/m_k])$.
  Moreover, 
  $\TOP{2}(\pi_1) = \TOP{2}(m'_k) = \TOP{2}(\rho(n'_k))$. 
  Furthermore,
  $\rho{\restriction}_{[n'_k,\length(\rho)]}$ never looks below the
  topmost word of $m'_k$ because $n'_k$ is the maximal node where
  the generalised milestone $m'_k$ is visited. 
  Thus, $\Pop{2}(m'_k):\bot \prefixeq
  \rho{\restriction}_{[n'_k,\length(\rho)]}$. Due to Lemma 
  \ref{Lem:BlumensathHOLevel2}, 
  \begin{align*}
    \pi_2:=\rho{\restriction}_{[n'_k, \length(\rho)]}[\Pop{2}(m'_k):\bot/
    \Pop{2}(\hat s):\bot]    
  \end{align*}
  is well defined. It starts in the last stack of $\pi_1$. 
  
  Now, we define the run 
  \begin{align*}
    \hat\rho:= \rho{\restriction}_{[0,n_k]} \circ \pi_1 \circ \pi_2.     
  \end{align*}
  Either $\height(\hat\rho) < \height(\rho)$ and we are done
  or there are less words of height
  $\height(\rho)$ in the last stack of $\hat\rho$ than in the last
  stack of $\rho$ and we conclude by induction. 
\end{proof}

\paragraph{Proof of Lemma \ref{lemmaPumpingWidth}}

Recall that Lemma \ref{lemmaPumpingWidth} asserts that for any run
$\rho$ there is another run $\rho'$ such that $\rho$ and $\rho'$ end
in stacks with equal topmost word but the width of the final stack of
$\rho'$ is bounded in terms of its height. 

\begin{proof}[Proof of Lemma \ref{lemmaPumpingWidth}]
  Assume that $\rho$ is a run with 
  $\length(\rho)> \FuncBoundWidthWord(\height(\rho))$.
  We denote by $n_i$ the maximal position in $\rho$ such that the stack
  at $\rho(n_i)$ is $\Pop{2}^i(\rho)$ for each 
  \mbox{$0\leq i \leq \width(\rho)$.}
  There are less than
  $\frac{\FuncBoundWidthWord(\height(\rho))}{\lvert Q 
    \rvert}$ many words of length up to $\height(\rho)$. 
  Thus, there are $i < j$ such that
  \begin{enumerate}
  \item $\TOP{2}\left(\Pop{2}^i(\rho)\right) =
      \TOP{2}\left(\Pop{2}^j(\rho)\right)$, and 
  \item $\rho(n_i)=\left(q,\Pop{2}^i(\rho)\right)$ and
    $\rho(n_j)=\left(q,\Pop{2}^j(\rho)\right)$ for some $q\in Q$. 
  \end{enumerate}
  Now, let $s_i:=\Pop{2}^{i+1}(\rho)$ and $s_j:=\Pop{2}^{j+1}(\rho)$. 
  There is a unique stack $s$ such that
  $\rho(\length(\rho)) = (\hat q, s_i:s)$. 
  $\rho{\restriction}_{[n_i, \length(\rho)]}$ is a run from
  $\Pop{2}^i(\rho)$ to $s_i:s$ that never visits $s_i$. 
  Thus,
  \begin{align*}
    \hat\rho_1:=\rho{\restriction}_{[n_i, \length(\rho)]}[s_i:\bot/s_j:\bot]    
  \end{align*}
  is a well defined run. 
  The composition $\hat\rho:=\rho{\restriction}_{[0,n_j]} \circ \hat\rho_1$
  satisfies the claim. 
\end{proof}

\paragraph{Proof of Proposition \ref{Prop:ConstructOneStep}}

We decompose the proof of Proposition \ref{Prop:ConstructOneStep} in
several lemmas. Recall that this Proposition is about
a run $\rho$ from some stack $s:w$ to some stack $s:w:w'$ that does not visit
substacks of $s:w$. Such a run decomposes into three parts. First, it
performs a $\Clone{2}$ operation. Then there is a run from $w:w$ to
$w:(w\sqcap w')$, i.e., a run that removes letters from $w$ until it
reaches the greatest common prefix of $w$ and $w'$. Finally, the run
constructs $w'$ from the prefix $w\sqcap w'$. In the following we 
first treat the second and the third part separately. We prove  lemmas
that allow to transfer each of these parts from one starting stack to
another one. Afterwards, we compose these arguments in order to obtain
the proof of Proposition \ref{Prop:ConstructOneStep}. 

The first lemma is concerned with a transfer of several runs starting
with the same stack $s:w$ and ending in the same stack $s:w'$ for 
some prefix $w'$ of $w$. 

\begin{lemma} 
  Let $z,m,n\in\N$ such that $z\geq 2$ and $z > m$. 
  Let $v,w,w'$ be words with $v\leq w$, and $q,\hat q\in Q$
  states. Let there be 
  pairwise distinct runs $\rho_1, \dots, \rho_m$ from
  $(q,w)$ to $(\hat q, v)$ such that each $\rho_i$  does not visit $v$
  before $\length(\rho_i)$. 
  If $ w \wordequiv{n+1}{z} w'$, then there exist a word $v'\leq
  w'$ and pairwise 
  distinct runs $\rho_1', \dots,  \rho_m'$ from
  $(q,w')$  to $(\hat q, v')$ such that
  \begin{align*}
    &v \wordequiv{n}{z} v'\text{, }\TOP{1}(v)=\TOP{1}(v'),\\
    &v=w\text{  iff }  v'=w'     \text{and }\\
    &\rho_i' \text{ does not visit } v' \text{ before } \length(\rho_i').    
  \end{align*}
\end{lemma}
\begin{proof}
  Let $i\in\N$ be such that $v= w_{-i}$. Then $i$ is labelled by
  $S^m_{q,\hat q}$, $P_{\TOP{1}(v)}$, and $\Typ{n}{z}{z}(v)$ in
  $\Lin{n+1}{z}{z}(w)$. Since 
  $\Lin{n+1}{z}{z}(w) \simeq_z \Lin{n+1}{z}{z}(w')$, there is some 
  $i'\in \N$ such that $i'$ in $\Lin{n+1}{z}{z}(w')$ is labelled
  by the same relations as $i$ in $\Lin{n+1}{z}{z}(w)$. Due to $z\geq 2$,
  we can choose this $i'$ such that $i'\neq 0$ if and only if $i\neq 0$.  
  Note that for $v':=w'_{-i'}$, $\TOP{1}(v')=\TOP{1}(v)$ because $i'$
  and $i$ agree on the label $P_{\TOP{1}(v)}$. Since $i'\in
  S^m_{q,\hat q}$, there are $m$ pairwise distinct runs from $(q,w')$
  to $(\hat q, v')$ that visit $v'$ only in the final
  configuration. Finally, $v \wordequiv{n}{z} v'$ due to the 
  fact that $i'$ and $i$ agree on the labels characterising 
  $\Typ{n}{z}{z}(v')$ and $\Typ{n}{z}{z}(v)$, respectively. 
\end{proof}

The next lemma is in some sense the ``otherwise'' to the previous
one. This lemma allows to transfer runs starting in the same word $w$ but
ending in different prefixes of $w$. 

\begin{lemma} 
  Let $z,m,n\in\N$ such that $z \geq 2$ and $z > m$. 
  Let $w,w'$ be words, \mbox{$v_1, \dots, v_m\leq w$} pairwise
  distinct prefixes.  Let there be
  pairwise distinct runs $\rho_1, \dots, \rho_m$ such that $\rho_i$ is
  from $(q,w)$ to $(q_i, v_i)$ and $\rho_i$ does not visit $v_i$
  before $\length(\rho_i)$. 
  If $w \wordequiv{n+1}{z} w'$, then there are pairwise distinct
  prefixes  $v_1', \dots, v_m' \leq w'$ and pairwise 
  distinct runs $\rho_1', \dots,  \rho_m'$ such that each $\rho_i'$ starts in
  $(q,w')$  and ends in $(q_i, v_i')$, $v_i' \leq w'$, and
  $v_i' \wordequiv{n}{z} v_i$ such that $\rho_i'$ does not visit
  $v_i'$ before $\length(\rho_i')$. 
\end{lemma}
\begin{proof}
  Any winning strategy in the  $z$-round Ehrenfeucht-\Fraisse game on
  $\Lin{n+1}{z}{z}(w)$ and $\Lin{n+1}{z}{z}(w')$ chooses responses
  $v_1', \dots, v_m'$ for $v_1, \dots, v_m$ such that the labels of
  the nodes associated to $v_i'$ in $\Lin{n+1}{z}{z}(w')$ and the
  nodes associated to $v_i$ in $\Lin{n+1}{z}{z}(w)$ agree and the
  $v_i'$ are pairwise distinct. 
  Hence, $v_i \wordequiv{n}{z} v_i'$. Now, there are runs $\rho_i'$ as
  desired due to the fact that the node associated with $v_i'$ in
  $\Lin{n+1}{z}{z}(w')$ is labelled by $S^1_{q,q_i}$. 
\end{proof}
\begin{remark} \label{Rem:CommonprefixProblem}
  Since $z>m$, the strategy preserves also the labels of the left and
  right neighbour. Thus, if $v_i=w_{-k}$ and $v'_i=w'_{-k'}$, then
  $k=0$ if and only if $k'=0$, and if $k\neq 0$, then 
  $\TOP{1}(w_{-k+1}) = \TOP{1}(w'_{-k'+1})$. 
\end{remark}
The combination of the previous two lemmas yields the following corollary.
\begin{corollary} \label{CorPrefixReachabilityTransfer}
  Let $z,m,n\in\N$ such that $z\geq 2$ and $z>m$ and 
  let $w,w'$ be words.
  Let there be
  pairwise distinct runs $\rho_1, \dots, \rho_m$ such that $\rho_i$ is
  a run from $(q_i,w)$ to 
  $(\hat q_i, v_i)$ for $v_i\leq w$ and $\rho_i$ does not visit
  $v_i$  before $\length(\rho_i)$. 
  If $w \wordequiv{n+1}{z} w'$,
  then there are prefixes $v_1', \dots, v_m'\leq w'$ and pairwise
  distinct runs $\rho_1', \dots,  \rho_m'$ such that $\rho_i$ starts in
  $(q_i,w')$  and ends in $(\hat q_i, v_i')$, 
  $v_i' \wordequiv{n}{z} v_i$, and $\rho_i'$ does not visit $v_i'$
  before $\length(\rho_i')$. 
\end{corollary}

This corollary provides the transfer of runs from some stack  $s:w$ to
stacks $s:w_i$ with $w_i\leq w$ to another starting stack $s':w'$ if 
$w$ and $w'$ are equivalent words. 

Now, we start the investigation of the other direction. 
We analyse runs from some word $w$ to some extension $wv$. 
If $w'$ is equivalent to $w$, then we first transfer 
the run from $w$ to $wv$ to a run from $w'$ to $w'v$. Afterwards, we
even provide a lemma that allows to shrink $v$ during this transfer
process.  

\begin{lemma} \label{LemmaBuildword}
  Let $n,m,z\in \N$ with $z\geq 2$ and $z>m$. 
  Let $\rho_0, \dots, \rho_m$ be pairwise distinct runs such that $\rho_i$
  starts in $(q,w)$ and ends in $(\hat q, ww_i)$. If $w'$ is a word
  such that $w\wordequiv{n}{z} w'$, then there are
  pairwise distinct runs $\rho_0', \dots, \rho_m'$ such that $\rho_i'$
  starts in $(q,w')$ and ends in $(\hat q, w'w_i)$. 
\end{lemma}
\begin{proof}
  This follows directly from the fact that
  \begin{align*}
  \ReturnFunc{z}(w)=\ReturnFunc{z}(w')\text{ and }
  \LoopFunc{z}(w)=\LoopFunc{z}(w').
  \end{align*}
  Each run from $w$ to $ww_i$ is a 
  sequence of loops and push operations. But the existence of the 
  corresponding loops when starting in $w'$ are guaranteed by the
  inductive computability of the number of loops and returns (cf.
  Proposition \ref{Prop:AutomatonForLoops}).  
\end{proof}

As already indicated, we will now improve this lemma in the sense that
we shrink the word $w_i$ to a short word $w_i'$. 

We fix some pushdown system $\mathcal{N}$ of level $2$ with state set
$Q$ and stack alphabet $\Sigma$. 
Let 
\begin{align*}
  &\BoundHeightOnestepConstruction:
  \N^2 \rightarrow \N \\
  &(n,z)\mapsto 1+\lvert Q \rvert \cdot 
  \left\lvert \nicefrac{\Sigma^*}{\wordequiv{n}{z}}\right\rvert
\end{align*}
for $\left\lvert \nicefrac{\Sigma^*}{\wordequiv{n}{z}}\right\rvert$ the
index of $\wordequiv{n}{z}$. 

\begin{lemma}
  For $q,\bar q\in Q$, $w$ and $v$ words, 
  let $\rho$ be a run from $(q,w)$ to $(\bar q, wv)$. 
  If  $\lvert v \rvert > \BoundHeightOnestepConstruction(n,z)$ then
  there is a run 
  $\rho'$ from $(q,w)$ to $(\bar q,wv')$ such that
  \begin{align*}
    &\lvert v \rvert - \BoundHeightOnestepConstruction(n,z) \leq
    \lvert v' \rvert < \lvert v \rvert, \\
    &wv \wordequiv{n}{z} wv' \text{ and}\\
    &\text{the first letters of } v \text{ and } v' \text{ coincide.}    
  \end{align*}
\end{lemma}

\begin{proof}
  Assume that $\lvert v \rvert >
  \BoundHeightOnestepConstruction(n,z)$. 
  Then
  there are two distinct  prefixes $\varepsilon < u_1 < u_2 \leq v$ such that 
  \begin{enumerate}
  \item 
    $\rho$ passes $wu_1$ and $wu_2$ in the same state
    $\hat q\in Q$,
  \item $wu_1 \wordequiv{n}{z} wu_2$, and
  \item $1\leq \lvert u_1 \rvert < \lvert u_2 \rvert \leq f(n,z)+1$.
  \end{enumerate}  
  Set $v_i$ to be the unique word such that $v=u_iv_i$ for $i\in\{1,2\}$. 
  Since $\wordequiv{n}{z}$ is a right congruence,  $wu_1v_2
  \wordequiv{n}{z} wu_2v_2 =wv$. 
  $wu_1 \wordequiv{n}{z} wu_2$ implies that 
  $\ReturnFunc{z}(wu_1)=\ReturnFunc{z}(wu_2)$ and
  $\LoopFunc{z}(wu_1)=\LoopFunc{z}(wu_2)$. 
  Since there is a run from $(\hat q, wu_2)$ to $(\bar q,
  wu_2v_2)=(\bar q, wv)$,  we can use the prefix replacement 
  $[wu_2 / wu_1]$ we obtain a 
  run
  $\hat\rho$ from $(\hat q,wu_1)$ to 
  $(\bar q, wu_1v_2)$. 
  Composition of the initial part of $\rho$ up to 
  $(\hat q, wu_1)$ with $\hat\rho$ yields a run $\rho'$.
  By construction $\rho'$ satisfies the claim. 
\end{proof}

The following corollary uses the previous lemma in such a way that we
can transfer some run to a new run that does not coincide with certain
given runs 

\begin{corollary}  \label{BuildWordsShort}
  Let $\rho$,$v$, and $w$ be as in the previous lemma and let
  $\rho_1, \dots, \rho_m$ be runs distinct from $\rho$, then we can find a run
  $\rho'$ distinct from 
  all $\rho_i$ for $1\leq i \leq m$ from $(q,w)$ to
  $(\bar q, wv')$ such that 
  $\lvert v' \rvert \leq m \cdot \BoundHeightOnestepConstruction(n,z)$,
  $wv \wordequiv{n}{z} wv'$, and $v$ and $v'$ start with the same letter.
\end{corollary}
\begin{proof}
  If $v$ is long enough, we find $m+2$ words that are visited
  in the same state and which are of the same type. There is
  one pair among these words which can be used as $u_1$ and $u_2$ in
  the previous 
  lemma such that the final configuration does not agree with that of
  any of the other runs $\rho_1, \dots, \rho_m$. 
\end{proof}

For the proof of Proposition \ref{Prop:ConstructOneStep}, we now
compose the previous lemmas. 
Recall that the proposition says the following: given
$m$ runs $\rho_1, \dots, \rho_m$ that only add one word to a given stack
and given $m$ words that are equivalent to the words on top of the
initial stacks of the $\rho_i$, we can transfer the runs $\rho_1, \dots,
\rho_m$ to runs $\rho_1', \dots, \rho_m'$ that start at the given $m$
words and extend these by one word each such that the resulting new
words are equivalent to the words originally created by $\rho_1,
\dots, \rho_m$.  

\begin{proof}[Proof of Proposition \ref{Prop:ConstructOneStep}.]
  Let $\rho_1, \rho_2, \dots, \rho_m$ and $\hat\rho_1, \hat\rho_2,
  \dots, \hat\rho_{m-1}$ be runs as required in the proposition. 
  
  First, assume that $w_i \wordequiv{n-1}{z} w_j$ 
  and that all runs $\rho_i$ end in the same state, i.e., $q_i=q_j$,
  for all $1\leq i \leq j \leq m$. At the end of the proof we discuss
  the case that this assumption is not true. 
 
  Each run $\rho_i$ decomposes as $\rho_i=\rho^0_i\circ \rho^1_i\circ
  \rho^2_i$ where 
  $\rho^0_i$ performs only one clone operation, and $\rho^1_i$ is the 
  run from $s:w:w$ to the first occurrence of $s:w:(w\sqcap w_i)$. 

  Due to $\TOP{1}(w)=\TOP{1}(\hat w)$, there are runs ${\hat\rho^0}_i$
  from $\hat c$ to
  $\hat s:\hat w:\hat w$ performing only one clone operation and ending in the
  same state as $\rho^0_i$. 

  By Corollary \ref{CorPrefixReachabilityTransfer}, we can transfer
  the $\rho^1_i$ to runs ${\hat\rho^1}_i$ starting at $(q,\hat s:\hat
  w:\hat w)$ and ending
  at $\hat s:\hat w:\hat u_i$ with $\hat u_i\leq \hat w$  and with
  $w\sqcap w_i \wordequiv{n-1}{z} \hat u_i$. 
  The lemma guarantees that ${\hat \rho^1}_i={\hat\rho^1}_j$ iff
  ${\rho^1}_i={\rho^1}_j$. 
  
  For $v_i$ the word such that $w_i= (w\sqcap w_i) \circ v_i$,
  we use Lemma
  \ref{LemmaBuildword} to construct runs  from $\hat s:\hat w:\hat u_i$ to
  $(q_i,\hat s:\hat w:\hat u_iv_i)$ such that $\hat u_iv_i \wordequiv{n-1}{z}
  w_i$. 
  Applying Corollary \ref{BuildWordsShort}, we find words 
  $\hat v_1, \dots, \hat v_m$, and runs ${\hat\rho^2}_1, \dots,
  {\hat\rho^2}_m$ such that 
  ${\hat\rho^2}_i$ is a run from \mbox{$\hat s:\hat w:\hat u_i$} to
  \mbox{$(q_i,\hat  
  s:\hat w:\hat u_i\hat v_i)$} such that $\hat u_i \hat v_i
  \wordequiv{n-1}{z} w_i$ and such that
  $\hat u_i\hat v_i$ has length bounded by 
  \begin{align*}
    \BoundHeightOnestepConstructionSimultanious(m,\lvert \hat w\rvert,
    n, z) = \lvert \hat w \rvert +
    m\cdot\BoundHeightOnestepConstruction(n,z).    
  \end{align*}
  Furthermore, Corollary \ref{BuildWordsShort} assures that 
  $\hat\rho^2_i$ and $\hat\rho^2_j$ coincide if and only
  if $\rho^2_i$ and $\rho^2_j$ coincide.
  
  We conclude that 
  the runs
  \begin{align*}
    {\hat\rho^0}_1 \circ {\hat\rho^1}_1 \circ {\hat\rho^2}_1,
    {\hat\rho^0}_2 \circ {\hat\rho^1}_2 \circ {\hat\rho^2}_2, \dots, 
    {\hat\rho^0}_m \circ {\hat\rho^1}_m \circ {\hat \rho^2}_m    
  \end{align*}
  are
  pairwise distinct as follows.
  Heading for a contradiction assume that  there are $i\neq j$ such that
  ${\hat\rho^0}_i \circ {\hat\rho^1}_i \circ {\hat\rho^2}_i =
  {\hat\rho^0}_j \circ {\hat\rho^1}_j \circ {\hat\rho^2}_j$. 
  Since $\length({\hat\rho}^0_i) = 1 = \length({\hat\rho}^0_j)$, it
  follows that ${\hat\rho}^0_i={\hat\rho}^0_j$ and
  ${\hat\rho^1}_i \circ {\hat\rho^2}_i =
  {\hat\rho^1}_j \circ {\hat\rho^2}_j$. Since the runs coincide, we
  have $\hat u_i \hat v_i = \hat u_j \hat v_j$. 
  Using Remark \ref{Rem:CommonprefixProblem} and the fact that the
  first letters of $v_i$ and $\hat v_i$ agree, one concludes that
  \mbox{$\hat u_i = \hat w \sqcap \hat u_i \hat v_i = \hat w \sqcap
    \hat u_j  \hat v_j = \hat u_j$}.  But by definition,
  $\length({\hat\rho}^1_i)$ and 
  $\length({\hat\rho}^1_j)$, respectively, is the first occurrence
  of $\hat u_i =  \hat u_j$ in ${\hat\rho}^1_i \circ {\hat\rho}^2_i$ and
  ${\hat\rho}^1_j\circ {\hat\rho}^2_j$, respectively. 
  It follows that ${\hat\rho}^1_i={\hat\rho}^1_j$ and
  ${\hat\rho}^2_i={\hat\rho}^2_j$. 
  By definition, this implies that $\rho^0_i = \rho^0_j$,
  $\rho^1_i=\rho^1_j$, and $\rho^2_i=\rho^2_j$. Thus, we have
  $\rho_i=\rho_j$ contradicting the assumptions of the proposition.  

  Thus, the runs are pairwise distinct and one of them does not
  coincide with any 
  of the $\hat\rho_i$ for \mbox{$1\leq i \leq m-1$}. Without loss of generality,
  assume  that 
  \begin{align*}
    \hat\rho_m:={\hat\rho^0}_m \circ {\hat\rho^1}_m \circ {\hat\rho^2}_m \neq
    \hat\rho_i\text{ for all }1\leq i < m.    
  \end{align*}
  Note that $\hat\rho_m$ satisfies the claim of the proposition
  by construction.

  In the case that the runs $\rho_1, \dots, \rho_m$ end in
  configurations with 
  different states or different $\wordequiv{n-1}{z}$-types of their
  topmost words. 
  In this case, we just concentrate on those $\rho_i$ which end in the
  same state as 
  $\rho_0$ and with a topmost word of the same
  type as $w_0$. This is sufficient because some run $\rho$ can only
  coincide  with $\hat\rho_i$ if both runs end up in  stacks whose
  topmost words have the same type.  
\end{proof}

\subsection{Small-Witness Property via Isomorphisms of 
  Relevant Ancestors} 
\label{SectionEquivalenceonTuples}
In this section, we want to define a family of  equivalence relations
on  tuples of runs of a level $2$ nested pushdown tree. 
The equivalence class of a tuple $\rho_1, \dots, \rho_m$ with respect to
one of these relations is the isomorphism type of the 
substructure induced by the relevant $l$-ancestors of $\rho_1, \dots,
\rho_m$ extended by some information for
preserving this isomorphism during an Ehrenfeucht-\Fraisse game (while
decreasing $l$ in every round of the game). Recall that such a game
ends in a winning position for Duplicator if
the relevant $1$-ancestors of the elements that were chosen in the two
structures are isomorphic (cf. Lemma \ref{LemmaRelAnclocalIso}).

An important property of these equivalence relations is that they
have finite index because 
the sets of $l$-ancestors are finite and the information we add to the
structure can be encoded by a bounded number of unary predicates. 

Finally, we show how to construct small representatives for each
equivalence class. 
As explained in Section \ref{Sec:EFGame}, this property can be turned
into an \FO{} model checking algorithm on the class of $2$-\HONPT. 

\begin{definition}
  Let $\bar \rho=(\rho_1, \rho_2, \dots, \rho_m)$ be runs of a level
  $2$ pushdown 
  system $\mathcal{N}$ and let 
  $\mathfrak{N}:=\NPT(\mathcal{N})$ be the $2$-\HONPT generated by
  $\mathcal{N}$. 
  Let $l, n_1,n_2,z\in \N$. We define the following relations
  on $\RelAnc{l}{\bar\rho}$. 
  \begin{enumerate}
  \item For $k\leq l$ and $\rho\in \bar\rho$, let $P^k_\rho:=\{\pi \in
    \RelAnc{l}{\bar\rho}: \pi \in \RelAnc{k}{\rho}\}$.
  \item 
    Let $\stackequivTyp{n_2}{z}{n_1}$ denote the function 
    that maps a run $\pi$ to the equivalence
    class of the last stack of $\pi$ with respect to
    $\stackequiv{n_2}{z}{n_1}$. 
  \end{enumerate}
  We write $\AncestorClass{l}{n_1}{n_2}{z}(\bar\rho)$ for the following
  expansion of the relevant ancestors of $\bar\rho$:
  \begin{align*}
    \AncestorClass{l}{n_1}{n_2}{z}(\bar\rho):=
    (\mathfrak{N}{\restriction}_{\RelAnc{l}{\bar\rho}},  
    \trans{},
    \hookrightarrow,
    \overset{+1}{\hookrightarrow},
    \stackequivTyp{n_2}{z}{n_1},
    (P^k_{\rho_j})_{k\leq l, 1\leq j \leq m}).
  \end{align*}
  For tuples of runs $\bar\rho=(\rho_1, \dots, \rho_m)$ and 
  $\bar\rho'=(\rho_1', \dots, \rho_m')$ we set 
  $\bar\rho \RelAncequiv{l}{n_2}{z}{n_1} \bar\rho'$ if
  \begin{align*}
    \AncestorClass{l}{n_1}{n_2}{z}(\bar\rho) \simeq
    \AncestorClass{l}{n_1}{n_2}{z}(\bar\rho').
  \end{align*}
\end{definition}
\begin{remark}
  \begin{itemize}
  \item[]
  \item   If $\bar\rho \RelAncequiv{l}{n_2}{z}{n_1} \bar\rho'$ then there is a
    unique isomorphism $\varphi:
    \AncestorClass{l}{n_1}{n_2}{z}(\bar\rho) \simeq 
    \AncestorClass{l}{n_1}{n_2}{z}(\bar\rho')$ witnessing this
    equivalence.
    Note that due to the predicate $P^0_j$, $\rho_j$ is mapped to
    $\rho_j'$ for all $1\leq j \leq m$. Due to the predicate $P^l_j$, the
    relevant ancestors of $\rho_j$ are mapped to the relevant ancestors
    of $\rho'_j$. Finally, $\varphi$ must preserve the order of the
    relevant ancestors of $\rho_j$ because they form a chain with
    respect to $\trans{} \cup \overset{+1}{\hookrightarrow}$ 
    (cf. Proposition \ref{Prop:NextRelAnc}).
  \item 
    Due to Lemma \ref{LemmaRelAnclocalIso}, it is clear that
    $\bar\rho \RelAncequiv{l}{n_2}{z}{n_1} \bar\rho'$ implies that there
    is a partial 
    isomorphism mapping $\rho_i\mapsto \rho_i'$ for all $1\leq i \leq m$. 
  \end{itemize}
\end{remark}

Since equivalent relevant ancestors induce partial isomorphisms, a
strategy that preserves the equivalence between relevant ancestors is
winning for Duplicator in the Ehrenfeucht-\Fraisse-game.

Given a level $2$ pushdown system $\mathcal{N}$ 
we are going to show that there is a strategy in the
Ehrenfeucht-\Fraisse game on $\NPT(\mathcal{N})=:\mathfrak{N},
\bar\rho$ and $\mathfrak{N},\bar\rho'$ 
in which Duplicator can always
choose small elements compared to the size of the elements chosen so
far in the structure where he has to choose. Furthermore, this
strategy will 
preserve equivalence of the relevant ancestors in the following
sense. 
Let 
$\bar\rho , \bar\rho'\subseteq 
\mathfrak{N}$ be the  $n$-tuples chosen in the previous
rounds of the game. Assume that Duplicator managed to maintain the
relevant ancestors of these tuples equivalent, i.e., it holds that 
$\bar\rho\RelAncequiv{l}{n}{z}{k} \bar\rho'$.
Now, Duplicators strategy enforces that these tuples are extended by
runs $\pi$ and $\pi'$ satisfying the following.
There are numbers $k_i,l_i, n_i$ such that
$\bar\rho,\pi \RelAncequiv{l_i}{n_i}{z}{k_i} \bar\rho', \pi'$ and
furthermore, the size of the run chosen by Duplicator is small
compared to the elements chosen so far.

The exact claim is given in the following proposition. 

\begin{proposition} \label{Prop:2NPT-Strategy}
  Let $\mathcal{N}$ be a level $2$ pushdown system defining the higher
  order nested 
  pushdown tree 
  $\mathfrak{N}:=\NPT(\mathcal{N})$.  
  Given $\mathcal{N}$, we can compute functions
  \begin{align*}
    &\BoundHeight:\N^5\rightarrow \N,\\
    &\BoundWidth:\N^5\rightarrow \N, \text{and }\\
    &\BoundRunLength:\N^5\rightarrow \N
  \end{align*}
  with the following property.

  Let $n,z,n_1',n_2',l'\in\N$,  
  $l := 4l'+5, n_1 := n_1'+2(l'+2)+1,$ and
  $n_2 := n_2' + 4^{l'+1}+1$ such that $z\geq 2$ and $z> n\cdot 4^l$.
  Furthermore, let $\bar\rho$ and $\bar\rho'$ be $n$-tuples of runs of
  $\mathfrak{N}$ such that  
  \begin{enumerate}
  \item $\bar\rho \RelAncequiv{l}{n_2}{z}{n_1} \bar\rho'$, and
  \item $\length(\pi)\leq\BoundRunLength(n,l,n_1,n_2,z)$
    for all $\pi\in\RelAnc{l}{\bar\rho'}$,
  \item $\height(\pi) \leq  \BoundHeight(n,z,l,n_1,n_2)$ for all
    $\pi\in\RelAnc{l}{\bar\rho'}$, and
  \item $\width(\pi) \leq \BoundWidth(n,z,l,n_1,n_2)$ for all
    $\pi\in\RelAnc{l}{\bar\rho'}$. 
  \end{enumerate}
  For each $\rho\in \mathfrak{N}$ there is some $\rho'\in
  \mathfrak{N}$ such that 
  \begin{enumerate}
  \item $\bar\rho, \rho\RelAncequiv{l'}{n_2'}{z}{n_1'} \bar\rho', \rho'$,
  \item $\length(\pi)\leq\BoundRunLength(n+1,l',n_1',n_2',z)$
    for all $\pi\in\RelAnc{l'}{\bar\rho',\rho'}$,
  \item $\height(\pi) \leq  \BoundHeight(n+1,z,l',n_1',n_2')$ for all
    $\pi\in\RelAnc{l}{\bar\rho',\rho'}$, and
  \item $\width(\pi) \leq \BoundWidth(n+1,z,l',n_1',n_2')$ for all
    $\pi\in\RelAnc{l}{\bar\rho',\rho'}$. 
  \end{enumerate}
\end{proposition}

In the next section we show how this proposition can be used to
define an \FO{} model checking algorithm on nested pushdown trees of
level $2$. The rest of this section proves the main proposition. 
For this purpose we split the claim into several pieces. 
The proposition asserts bounds on the length of the runs and on the
sizes of the final stacks of the relevant ancestors. As the first step
we prove that Duplicator has a strategy that chooses runs with
small final stacks. This result relies mainly on the 
Proposition \ref{Prop:CompatibilityStackOpTypeq} and 
Proposition \ref{Prop:ConstructOneStep}. These results allow to
construct equivalent relevant 
ancestor sets that contain runs ending in small stacks.\footnote{In
  the following we sometimes say ``Duplicator can choose small
  stacks''. This expression always means that ``Duplicator can
  choose a run such that all its relevant ancestors end in small
  stacks''.}
Afterwards, we apply the general bounds on short loops
(cf. Proposition \ref{Prop:FuncBoundLoopLengthLemma}) 
in order to shrink the length of the runs involved.

The reader who is not interested in the details of the proof of
Proposition \ref{Prop:2NPT-Strategy}, may skip this part and continue
reading Section \ref{SubsectionRunBounds}.

\paragraph{Construction of Isomorphic Relevant Ancestors}

Before we prove that Duplicator can choose short runs, 
we state some 
auxiliary lemmas concerning the construction of isomorphic relevant
ancestors. The following lemma gives a sufficient criterion that
allows to check that the relevant ancestors of two runs are
equivalent. Afterwards, we show that for each run $\rho$ we can
construct a second run $\rho'$ satisfying this criterion.

\begin{lemma} \label{LemmaConstructedRelAncEquiv}
  Let $\rho_0 \prec \rho_1 \prec \dots \prec \rho_m=\rho$ be runs such that 
  $\RelAnc{l}{\rho}=\{\rho_i:0\leq i \leq m\}$.
  If 
  $\hat\rho_0\prec \hat\rho_1 \prec \dots \prec \hat\rho_m $ are
  runs such that 
  \begin{itemize}
  \item the final states of $\rho_i$ and $\hat\rho_i$ coincide,
  \item $\rho_0=\Pop{2}^l(\rho_m)$ or $\lvert \rho_0 \rvert = \lvert
    \hat \rho_0 \rvert=1$, 
  \item $\rho_0 \stackequiv{n_2}{z}{n_1} \hat\rho_0$, and
  \item $\rho_i \mathrel{*} \rho_{i+1}$ iff $\hat\rho_i \mathrel{*}
    \hat\rho_{i+1}$ for 
    all $1\leq i < m$ and $*\in\{\overset{+1}{\hookrightarrow}\}\cup\{
    \trans{\gamma}: \gamma\in\Gamma\}$,
  \end{itemize}
  then
  \begin{align*}
    \RelAnc{l}{\hat\rho_m}=\{\hat\rho_i:0\leq i \leq m \}.     
  \end{align*}
  If additionally $\TOP{2}(\rho_i) \wordequiv{n_2'}{z}
  \TOP{2}(\hat\rho_i)$, then 
  \begin{align*}
    \hat\rho_m \RelAncequiv{l}{n_2'}{z}{n_1} \rho_m    
  \end{align*}
  for $n_2':= n_2- 4^l$. 
\end{lemma}
\begin{proof}
  First, we show that for all $0 \leq i < j \leq m$, 
  the following statements are true:
  \begin{align}
    \rho_i \trans{\gamma} \rho_j &\text{ iff } 
    \hat\rho_i \trans{\gamma} \hat\rho_j,  \label{nextbyop} \\
    \rho_i \hookrightarrow  \rho_j &\text{ iff } 
    \hat\rho_i \hookrightarrow \hat\rho_j, \text{ and}\\
    \rho_i \overset{+1}{\hookrightarrow} \rho_j &\text{ iff } 
    \hat\rho_i \overset{+1}{\hookrightarrow} \hat\rho_j. \label{nextbyplusone}
  \end{align}
  Note that $\rho_i\trans{\gamma} \rho_j$ implies $j=i+1$.
  Analogously, $\hat\rho_i\trans{\gamma} \hat\rho_j$ implies $j=i+1$.
  Thus, \ref{nextbyop} is true by definition of the sequences. 

  For the other parts, 
  it is straightforward to see that 
  $\lvert \rho_k \rvert - \lvert \rho_j \rvert = 
  \lvert  \hat\rho_k \rvert - \lvert \hat\rho_j \rvert$ for all 
  $0\leq j\leq k  \leq m$: for $k=j$ the claim holds trivially. For
  the induction step from $j$ to $j+1$, the claim follows from the
  assumption that  
  $\rho_j * \rho_{j+1}$ if and only if $\hat\rho_j * \hat\rho_{j+1}$
  for all
  $*\in\{\overset{+1}{\hookrightarrow}\}\cup
  \{\trans{\gamma}:\gamma\in\Gamma\}$. 

  Furthermore, assume that there is some $\hat\pi$ such that
  $\hat\rho_k \prec \hat\pi \prec \hat\rho_{k+1}$. Then it cannot be
  the case that  
  $\hat\rho_k \trans{\gamma} \hat\rho_{k+1}$. 
  This implies that $\rho_k \overset{+1}{\hookrightarrow}
  \rho_{k+1}$. Due to our assumptions, it follows that
  $\hat\rho_k \overset{+1}{\hookrightarrow} \hat\rho_{k+1}$. 
  We conclude directly that $\lvert\hat\pi \rvert \geq \lvert
  \hat\rho_{k+1} \rvert > \lvert \hat \rho_k \rvert$. 
  Thus, 
  \begin{align*}
    &\rho_j \hookrightarrow \rho_k\text{ iff}\\
    &\lvert \rho_j \rvert = \lvert \rho_k \rvert\text{ and }
    \lvert \pi \rvert > \lvert \rho_j \rvert\text{ for all }\rho_j
    \prec \pi \prec \rho_k \text{ iff}\\    
    &  \lvert \hat\rho_j \rvert = \lvert \hat\rho_k \rvert\text{ and }
    \lvert \hat\pi \rvert > \lvert \hat\rho_j \rvert\text{ for all }\hat\rho_j
    \prec \hat\pi \prec \hat\rho_k \text{ iff}\\
    &\hat\rho_j \hookrightarrow \hat\rho_k.
  \end{align*}
  Analogously, one concludes that \ref{nextbyplusone} holds. 

  We now show that $\RelAnc{l}{\hat\rho_m}=\{\hat\rho_i:0\leq i \leq
  m\}$. 
  Note that 
  \begin{align*}
    \RelAnc{l}{\hat\rho_m} \cap \{\pi: \hat\rho_{m} \preceq \pi\} = 
    \{\hat\rho_m\}.    
  \end{align*}
  Now assume that  there is some $0 \leq m_0 \leq m$ such that    
  \begin{align*}
    &\RelAnc{l}{\hat\rho_m} \cap \{\pi: \hat\rho_{m_0} \preceq \pi\} = 
    \{\hat\rho_i: m_0\leq i \leq m\} \text{ and}\\
    &\rho_i\in\RelAnc{k}{\rho}\text{ iff }
    \hat\rho_i\in\RelAnc{k}{\hat\rho_m}\text{ for all }  k\leq l\text{
      and } i\geq m_0.
  \end{align*}
  Now, we distinguish the following cases. 
  \begin{itemize}
  \item 
    If $\rho_{m_0-1} \trans{\op}
    \rho_{m_0}$ for some stack-operation $\op$ then $\hat\rho_{m_0-1}
    \trans{\op} \hat\rho_{m_0}$ due to \ref{nextbyop}. 
    Thus, there are no runs $\rho_{m_0-1} \prec \pi \prec
    \rho_{m_0}$. Hence, we only have to show that
    \mbox{$\rho_{m_0-1}\in\RelAnc{k}{\rho_m}$} if and only if
    $\hat\rho_{m_0-1}\in\RelAnc{k}{\hat\rho_m}$ for all $k\leq l$. 

    If $\rho_{m_0-1}\in\RelAnc{k}{\rho_m}$, then  there
    is some $j\geq m_0$ such that
    \mbox{$\rho_{j}\in\RelAnc{k-1}{\rho_m}$} and  
    $\rho_{m_0-1}$ is connected to $\rho_j$ via some edge. 
    But then \mbox{$\hat\rho_j\in\RelAnc{k-1}{\hat\rho_m}$} and 
    $\hat\rho_{m_0-1}$ is connected with $\hat\rho_j$ via the same
    sort of edge. Thus,
    \mbox{$\hat\rho_{m_0-1}\in\RelAnc{k}{\hat\rho_m}$}. 

    The other direction is completely analogous. 
  \item 
    Now, consider the case that there is some $\rho_{m_0-1} \prec \pi
    \prec \rho_{m_0}$. Since its direct predecessor is not in
    $\RelAnc{l}{\rho_m}$, $\rho_{m_0}\notin\RelAnc{l-1}{\rho}$. 
    Thus, $\hat\rho_{m_0}\notin\RelAnc{l-1}{\hat\rho}$. 
    By
    construction of the $\hat\rho_i$, 
    $\hat\rho_{m_0-1} \overset{+1}{\hookrightarrow} \hat\rho_{m_0}$.
    Thus, $\lvert \hat \pi \rvert \geq \lvert \hat\rho_{m_0} \rvert$ for
    all $\hat\rho_{m_0-1} \prec \hat \pi \prec \hat\rho_{m_0}$.
    This implies that 
    $\pi \not\hookrightarrow \hat\rho_i$ and
    $\pi\not\overset{+1}{\hookrightarrow} \hat\rho_i$
    for all $m_0 < i \leq m$. This shows that
    $\pi\notin\RelAnc{l}{\hat\rho_m}$. 
    
    Now, for all $k\leq l$ we conclude completely analogous to the
    previous case that 
    $\hat\rho_{m_0-1}\in\RelAnc{k}{\hat\rho_m}$ iff 
    $\rho_{m_0-1}\in\RelAnc{k}{\rho_m}$. 
  \end{itemize}
  Up to now, we have shown that 
  $\RelAnc{l}{\hat\rho_m} \cap \{\pi: \hat\rho_0 \preceq \pi \} = 
  \{\hat\rho_i: 0 \leq i \leq m\}$.  
  In order to prove 
  $\RelAnc{l}{\hat\rho_m} = 
  \{\hat\rho_i: 0 \leq i \leq m\}$, we have to show that $\hat\rho_0$
  is the minimal element of $\RelAnc{l}{\hat\rho_m}$. 
  
  There are the following cases
  \begin{enumerate}
  \item $\rho_0=\Pop{2}^l(\rho_m)$. In this case, we conclude that
    $\hat\rho_0 = \Pop{2}^l(\hat\rho_m)$ by construction. But 
    Lemma \ref{LemmaMinRelAnc} then implies that $\hat\rho_0$ is the
    minimal element of $\RelAnc{l}{\hat\rho_m}$.
  \item $\lvert \rho_0 \rvert = \lvert \hat\rho_0 \rvert = 1$.
    Note that $\rho_0\notin\RelAnc{l-1}{\rho_m}$
    because $\rho_0$ is minimal in $\RelAnc{l}{\rho_m}$. Thus, we know
    that $\hat\rho_0\notin\RelAnc{l-1}{\hat\rho_m}$. 
    
    Heading for a contradiction, assume that there is some
    $\hat\pi\in\RelAnc{l}{\hat\rho_m}$ with 
    $\hat\pi\prec\hat\rho_0$. 
    We conclude that $\hat\pi\hookrightarrow
    \hat\rho_k$ or $\hat\pi\overset{+1}{\hookrightarrow}\hat\rho_k$ for
    some $\hat\rho_k\in\RelAnc{l-1}{\hat\rho_m}$. But this implies
    that $\lvert \hat\pi \rvert < \lvert \hat\rho_0 \rvert =1$. Since
    there are no stacks of width $0$, this is a contradiction.
    
    Thus, there is no $\hat\pi\in\RelAnc{l}{\hat\rho_m}$ that is a
    proper prefix of $\hat\rho_0$.
  \end{enumerate}
  We conclude that $\RelAnc{l}{\hat\rho_m}=\{\hat\rho_i: 0\leq i \leq
  m\}$. 
  
  Now, we prove the second part of the lemma. Assume that
  $\TOP{2}(\rho_i)\wordequiv{n_2'}{z} \TOP{2}(\hat\rho_i)$ for all 
  \mbox{$0\leq i \leq m$}. Since $\hat\rho_i$ and $\hat\rho_{i+1}$
  differ in at most one 
  word, a straightforward induction shows that $\rho_i
  \stackequiv{n_2'}{z}{n_1-\lvert\rho_0\rvert + \lvert \rho_i \rvert}
  \hat\rho_i$ (cf. Proposition
  \ref{Prop:CompatibilityStackOpTypeq}). But this implies $\hat\rho_m 
  \RelAncequiv{l}{n_2'}{z}{n_1} \rho_m$ because 
  \mbox{$\lvert \rho_0 \rvert \leq \lvert \rho_i \rvert$}  as we have
  seen in Lemma \ref{LemmaMinRelAnc}. 
\end{proof}

The previous lemma gives us a sufficient condition for the equivalence
of relevant ancestors of two elements. Now, we show how to
construct such a chain of relevant ancestors. 

\begin{lemma} \label{LemmaAncestorConstruction}
  Let $l,n_1, n_2, m,z\in\N$ such that $n_2\geq 4^l$ and $z\geq 2$. 
  Let 
  \begin{align*}
    &\rho_0 \prec \rho_1 \prec \dots \prec \rho_m=\rho\text{ be runs
      such that}\\ 
    &\RelAnc{l}{\rho}\cap \{\pi: \rho_0\preceq \pi\preceq \rho\} =
    \{\rho_i: 0\leq i \leq m\}.      
  \end{align*}
  Let $\hat\rho_0$ be a run such that
  $\rho_0 \stackequiv{n_2}{z}{n_1} \hat\rho_0$.
  Then we can effectively construct runs 
  \begin{align*}
    \hat\rho_0\prec \hat\rho_1\prec \dots \prec \hat\rho_m=:\hat\rho    
  \end{align*}
  such that 
  \begin{itemize}
  \item the final states of $\rho_i$ and $\hat\rho_i$ coincide for all $0\leq
  i \leq m$,
  \item $\rho_i\trans{\gamma} \rho_{i+1}$ iff
    $\hat\rho_i\trans{\gamma} \hat\rho_{i+1}$ and
    $\rho_i \overset{+1}{\hookrightarrow} \rho_{i+1}$ iff
    $\hat\rho_i\overset{+1}{\hookrightarrow} \hat\rho_{i+1}$ 
    for all $0\leq
    i < m$, and
  \item $\TOP{2}(\rho_i) \wordequiv{n_2- 4^l}{z} \TOP{2}(\hat\rho_i)$
    for all $0\leq
    i \leq m$.
  \end{itemize}
\end{lemma}
\begin{proof}
 Assume we have constructed 
 \begin{align*}
   \hat\rho_0 \prec \hat\rho_1\prec \dots \prec
   \hat\rho_{m_0},   
 \end{align*}
 for some $m_0<m$ such that for all $0 \leq i\leq
  m_0$
  \begin{enumerate}
  \item the final states of $\rho_i$ and $\hat\rho_i$
    coincide, \label{condition_States} 
  \item $\rho_i\trans{\gamma} \rho_{i+1}$ iff \label{TEST_three}
    $\hat\rho_i\trans{\gamma} \hat\rho_{i+1}$ and
    $\rho_i \overset{+1}{\hookrightarrow} \rho_{i+1}$ iff
    $\hat\rho_i\overset{+1}{\hookrightarrow} \hat\rho_{i+1}$ 
    (note that either $\rho_i \trans{} \rho_{i+1}$ or $\rho_i
    \overset{+1}{\hookrightarrow} \rho_{i+1}$ hold 
    due to Proposition \ref{Prop:NextRelAnc}), and
  \item \label{COndition_Wordequivalence}
    $\TOP{2}(\rho_i) \wordequiv{n_2-i}{z} \hat\rho_i$.
  \end{enumerate}
  We extend this chain by a new element $\rho'_{m_0+1}$ such that
  all these conditions are again satisfied. We distinguish two cases.
  
  First, assume that $\rho_{m_0} \trans{\gamma} \rho_{m_0+1}$. 
  Since $\rho_{m_0} \wordequiv{n_2-m_0}{z} \hat\rho_{m_0}$, 
  $\TOP{1}(\rho_{m_0}) = \TOP{1}(\hat\rho_{m_0})$. Due to Condition 
  \ref{condition_States}, their final states also coincide. 
  Hence, we can define
  $\hat\rho_{m_0+1}$ such that
  $\hat\rho_{m_0}\trans{\gamma} \hat\rho_{m_0+1}$. 
  Due to Proposition  
  \ref{Prop:CompatibilityStackOpTypeq},
  $\hat\rho_{m_0+1}$ satisfies Condition 
  \ref{COndition_Wordequivalence}.
  
  Now, consider the case $\rho_{m_0} \overset{+1}\hookrightarrow
  \rho_{m_0+1}$.
  The run from $\rho_{m_0}$ to $\rho_{m_0+1}$ starts from
  some stack $s$ and ends in some stack $s:w$ for $w$ some word, the
  first operation is a clone and then $s$ is never reached
  again. Hence, we can use Proposition \ref{Prop:ConstructOneStep}
  in order to find some appropriate $\hat\rho_{m_0+1}$ that satisfies
  Condition \ref{COndition_Wordequivalence}.
\end{proof}

The previous lemmas give us the possibility to construct an isomorphic
copy of the 
relevant ancestors of a single run $\rho$. In our proofs, we want to
construct such a copy while avoiding relevant ancestors of certain
other runs. Using the full power of Proposition \ref{Prop:ConstructOneStep}
we obtain the following stronger version of the lemma. 


\begin{corollary} \label{Cor:AncestorConstruction}
  Let $l,n_1, n_2, m,z\in\N$ be numbers such that $z> m \cdot 4^l$ and
  $n_2\geq 4^l$.
  As before, let 
  \mbox{$\RelAnc{l}{\rho_m}=\{\rho_i: 0 \leq i \leq m\}$} and $\hat\rho_0$
  some run such that
  $\rho_0 \stackequiv{n_2}{z}{n_1} \hat\rho_0$.
  Let $\bar\rho$ and $\bar\rho'$ be $m$-tuples such that 
  $\bar\rho \RelAncequiv{l}{n_2}{z}{n_1} \bar\rho'$ and $\varphi_l$ is an
  isomorphism witnessing this equivalence. 
  
  If  $\rho_0\in\RelAnc{l}{\bar \rho}$,
  $\varphi_l(\rho_0)=\hat\rho_0$, and if $\rho_1\notin
  \RelAnc{l}{\bar\rho}$   
  then we can construct $\hat\rho_1, \hat\rho_2, \dots,
  \hat\rho_m$ satisfying the conditions from the previous 
  lemma but additionally with the property that
  $\hat\rho_1\notin\RelAnc{l}{\bar\rho'}$. 
\end{corollary}
\begin{proof}
  We distinguish two cases. 
  \begin{enumerate}
  \item Assume that $\rho_0 \trans{} \rho_1$. 
    Due to the equivalence of $\rho_0$ and $\hat\rho_0$, we can apply
    the transition connecting $\rho_0$ with $\rho_1$ to $\hat\rho_0$
    and obtain a run $\hat\rho_1$. We have to prove that
    $\hat\rho_1\notin\RelAnc{l}{\bar\rho'}$.
    
    Heading for a contradiction assume that
    $\hat\rho_1\in\RelAnc{l}{\bar\rho'}$. Then $\varphi^{-1}_l$
    preserves the edge between $\hat\rho_0$ and $\hat\rho_1$, i.e., 
    $\rho_0=\varphi^{-1}_l(\hat\rho_0) \trans{}
    \varphi^{-1}_l(\hat\rho_1)$. But this implies that
    $\varphi^{-1}_l(\hat\rho_1)=\rho_1$ which contradicts the
    assumption that $\rho_1\notin\RelAnc{l}{\bar\rho}$. 
  \item 
    Assume that $\rho_0 \overset{+1}{\hookrightarrow} \rho_1$. 
    Up to threshold $z$, for each $\hat\pi$ such that
    $\hat\rho_0\overset{+1}{\hookrightarrow} \hat\pi$ and
    \mbox{$\hat\pi\in\RelAnc{l}{\bar\rho'}$} there is  a run
    $\rho_0 \overset{+1}{\hookrightarrow} \varphi_l^{-1}(\hat\pi)$. Since 
    $\rho_1\notin \RelAnc{l}{\bar\rho}$, we find another run
    $\hat\rho_1$ that satisfies 
    the conditions of 
    the previous lemma and
    $\hat\rho_1\notin \RelAnc{l}{\bar\rho'}$. This is due to the fact that
    Proposition \ref{Prop:ConstructOneStep} allows to transfer up to
    $z>\lvert \RelAnc{l}{\bar\rho'} \rvert$ many runs simultaneously.\qedhere
  \end{enumerate}
\end{proof}

\paragraph{Strategy for Choosing Small Stacks}

By now we are prepared to prove that Duplicator has a strategy that
preserves the isomorphism type of the relevant ancestors but chooses
short runs. 
First, we  prove the existence of a strategy choosing runs with small
final stacks. Afterwards, we show how to bound the length of such
runs. 
The analysis of this strategy decomposes into the local and the
global case. 
We say Spoiler makes a local move if he chooses a new element such
that the relevant ancestors of this element intersect with the
relevant ancestors of elements chosen so far. 
In this case Duplicator has to extend the other tuple by an element whose
relevant ancestors intersect with the relevant ancestors of this
tuple. 

We say Spoiler makes a global move if he  chooses an element such that
the relevant ancestors of this new element do not intersect with the
relevant ancestors of the
elements chosen so far. 
In this case Duplicator has to extend the other tuple by an element
whose relevant ancestors do not intersect with the relevant ancestors
of this tuple. 

We first head for the result that Duplicator can manage the
local case 
in such a way that he chooses an element such that all its relevant
ancestors end in small stacks. Then we show that Duplicator can manage
the global case analogously. 
Finally,  we show that Duplicator
can choose a short run ending in small stacks. 

\begin{lemma}   \label{LemmaLocalStep}
  Let $n,z,l,l',n_1,n_1',n_2,n_2'\in\N$ be numbers such that 
  $l = 4l'+4$,  $z > n \cdot 4^l$, $z\geq 2$, $n_1 = n_1' + 2 (l'+1)+1$,
  $n_1'>0$, $n_2 = n_2' + 4^{l'+1}+1$, and $n_2'>0$.

  Let $\bar\rho, \bar\rho'$ be $n$-tuples of runs such that
  $\varphi_l: \AncestorClass{l}{n_1}{n_2}{z}(\bar\rho) \simeq 
  \AncestorClass{l}{n_1}{n_2}{z}(\bar\rho')$ witnesses
  \mbox{$\bar\rho \RelAncequiv{l}{n_2}{z}{n_1} \bar\rho'$}. 
  Furthermore, let $\rho$ be some run such that 
  $\RelAnc{l'+1}{\rho}\cap \RelAnc{l'+1}{\bar\rho}\neq\emptyset$.
  Then there is some run $\rho'$ such that
  $(\bar\rho,\rho) \RelAncequiv{l'}{n_2'}{z}{n_1'} (\bar\rho', \rho')$. 
\end{lemma}
\begin{proof}
  Let $\rho_0\in\RelAnc{l'+1}{\rho}$ be maximal such that 
  \begin{align*}
    \RelAnc{l'+1}{\rho}\cap\{\pi : \pi \preceq \rho_0\}\subseteq
    \RelAnc{4l'+3}{\bar\rho}\subseteq\RelAnc{l}{\bar\rho}.    
  \end{align*}
  We choose numbers  $m_0 \leq 0 \leq m_1$ and runs 
  \begin{align*}
    \rho_{m_0} \prec \rho_{m_0+1} \prec
    \dots \prec \rho_0 \prec \rho_1 \prec \dots \prec \rho_{m_1}
  \end{align*}
  such that  
  $\RelAnc{l'+1}{\rho} = \{\rho_i: m_0 \leq i \leq m_1\}$. 
  We
  set $\rho_i':=\varphi_l(\rho_i)$ for all $m_0 \leq i\leq 0$. Next, we
  construct $\rho'_1, \dots, \rho'_{m_1}$  such that
  $\rho':=\rho'_{m_1}$ has relevant ancestors isomorphic to those of $\rho$.

  Analogously to the previous corollary, we can construct runs
  $\rho_1', \rho_2', \dots, \rho'_{m_1}$ such that 
  \begin{enumerate}
  \item $\rho_1'\notin \RelAnc{4l'+3}{\bar\rho'}$,
  \item the final states of $\rho_i$ and $\rho'_i$ coincide 
    for all \mbox{$0\leq i \leq m_1$}, and
  \item $\rho_i\trans{\gamma} \rho_{i+1}$ iff 
    $\rho'_i\trans{\gamma} \rho'_{i+1}$ and
    $\rho_i \overset{+1}{\hookrightarrow} \rho_{i+1}$ iff
    $\rho'_i\overset{+1}{\hookrightarrow} \rho'_{i+1}$ for all
    \mbox{$0\leq i < m_1$}.
  \end{enumerate}
  By definition, it is clear that condition 2 and 3 hold also for all
  $m_0\leq i < 0$. Using Lemma
  \ref{LemmaConstructedRelAncEquiv}, we obtain that 
  $\rho \RelAncequiv{l'+1}{n_2'}{z}{n_1'} \rho'$. 
  
  As a next step, we have to show that the isomorphism between
  $\RelAnc{l}{\bar\rho}$ and $\RelAnc{l}{\bar\rho'}$ and the
  isomorphism between 
  $\RelAnc{l'}{\rho}$ and $\RelAnc{l'}{\rho'}$ are compatible in the sense
  that they may be composed to an isomorphism between
  $\RelAnc{l'}{\bar\rho,\rho}$ and $\RelAnc{l'}{\bar\rho',\rho'}$. 
  
  The only possible candidate for such a combined isomorphism is of
  the form
  \begin{align*}
    \varphi_{l'}:\RelAnc{l'}{\bar\rho,\rho} &\to
    \RelAnc{l'}{\bar\rho',\rho'}\\ 
    \pi &\mapsto
    \begin{cases}
      \rho'_i &\text{for }  \pi=\rho_i, m_0\leq i\leq m_1\\
      \varphi_l(\pi) & \text{for } \pi\in\RelAnc{l'+1}{\bar\rho}.
    \end{cases}
  \end{align*}
  In order to see that this is a well-defined function, we have to
  show that if $\rho_i\in\RelAnc{l'+1}{\bar\rho}$
  then $\rho_i'=\varphi_l(\rho_i)$ for each $m_0\leq i \leq m_1$.
  Note that
  $\rho_i\in\RelAnc{l'+1}{\bar\rho}\cap\RelAnc{l'+1}{\rho}$ implies
  (using Corollary 
  \ref{CorRelAncDistBound}) that $\pi\in\RelAnc{3l'+3}{\bar\rho}$ for all
  $\pi\in\RelAnc{l'+1}{\rho}$ with $\pi\preceq \rho_i$. But then by definition
  $i\leq 0$ and by construction
  $\rho'_i = \varphi_l(\rho_i)$. 

  We claim that $\varphi_{l'}$ is an isomorphism. 
  Since we composed $\varphi_{l'}$ of existing isomorphisms
  $\RelAnc{l'}{\bar\rho}\simeq \RelAnc{l'}{\bar\rho'}$ and
  $\RelAnc{l'}{\rho}\simeq \RelAnc{l'}{\rho'}$, respectively, we only have
  to consider the following question: 
  let $\pi\in \RelAnc{l'}{\bar\rho}$ and
  $\hat\pi\in\RelAnc{l'}{\rho}$. Does $\varphi_{l'}$ preserve the
  existence and nonexistence of edges between $\pi$ and $\hat\pi$?

  In other words, we have to show that for each
  \begin{align*}
    *\in\{\hookrightarrow, \hookleftarrow,
    \overset{+1}{\hookrightarrow}, \overset{+1}{\hookleftarrow}\}&\cup\{
    \trans{\gamma}:
    \gamma\in\Gamma\}\cup\{\invtrans{\gamma}:\gamma\in\Gamma\},\\
    \pi * \hat\pi &\text{ iff } \varphi_{l'}(\pi) *
    \varphi_{l'}(\hat\pi).    
  \end{align*}
   The following case distinction treats all these cases. 
  \begin{itemize}
  \item   Assume that there is some 
    $*\in\{\hookrightarrow, 
    \overset{+1}{\hookrightarrow}\}\cup
    \{\trans{\gamma}:\gamma\in\Gamma\}$ such that  
    $\pi * \hat\pi$. Then 
    \mbox{$\pi\in\RelAnc{l'+1}{\rho}$}. Thus, there are $m_0 \leq i < j \leq
    m_1$ such that $\pi=\rho_i$ and $\hat\pi=\rho_j$. We have already seen
    that then $\varphi_{l'}(\pi) = \rho'_i$ and
    $\varphi_{l'}(\hat\pi)=\rho'_j$ and 
    these elements are connected by an edge of the same type due to
    the construction of $\rho'_i$ and $\rho'_j$. 
  \item Assume that there is some 
    $*\in\{\hookleftarrow, 
    \overset{+1}{\hookleftarrow}\} \cup
    \{\invtrans{\gamma}:\gamma\in\Gamma\}$ such that
    $\pi * \hat\pi$. Then 
    $\hat\pi\in\RelAnc{l'+1}{\bar\rho}$ whence $\varphi_{l'}$
    coincides with the 
    isomorphism $\varphi_l$ on $\pi$ and $\hat\pi$. But $\varphi_l$
    preserves edges whence $\pi*\hat\pi$  implies
    $\varphi_{l'}(\pi) *  \varphi_{l'}(\hat\pi)$.  
  \item Assume that there is some
    $*\in\{\hookrightarrow, 
    \overset{+1}{\hookrightarrow}\}\cup
    \{\trans{\gamma}:\gamma\in\Gamma\}$ such that  
    $\varphi_{l'}(\pi) * \varphi_{l'}(\hat\pi)$. 
    By definition,
    $\varphi_{l'}(\hat\pi)\in\RelAnc{l'}{\rho'}$ 
    whence \mbox{$\varphi_{l'}(\hat\pi)=\rho'_j$} for some $m_0\leq j \leq
    m_1$. Thus,  \mbox{$\varphi_{l'}(\pi)\in\RelAnc{l'+1}{\rho'}$} whence
    \mbox{$\varphi_{l'}(\pi)=\rho'_i$} for some $m_0\leq i <j$. 
    We claim that \mbox{$\pi = \rho_i$}. Note that due to Corollary
    \ref{CorRelAncDistBound} for all $m_0\leq k \leq i$ we have
    $\rho'_k\in\RelAnc{3l'+3}{\varphi_{l'}(\pi)}$. 
    Since $\varphi_{l'}(\pi)=\varphi_l(\pi)\in\RelAnc{l'}{\bar\rho'}$, 
    we conclude that $\rho'_k\in\RelAnc{4l'+3}{\bar\rho'}$ for all
    \mbox{$m_0\leq k \leq i$}. By
    construction,  
    this implies $i\leq 0$ and
    $\varphi_{l'}(\pi) = \rho_i'=\varphi_l(\rho_i)$.
    Furthermore, since $\pi\in\RelAnc{l'}{\bar\rho}$, 
    $\varphi_{l'}(\pi)=\varphi_l(\pi)$. 
    Since $\varphi_l$ is an isomorphism, it follows that 
    $\pi=\rho_i$.  
    But this implies that there is an edge from $\pi=\rho_i$ to
    $\hat\pi=\rho_j$.
  \item 
    Assume that there is some 
    $*\in\{\hookleftarrow, 
    \overset{+1}{\hookleftarrow}\} \cup
    \{\invtrans{\gamma}:\gamma\in\Gamma\}$ such that
    $\varphi_{l'}(\pi)*\varphi_{l'}(\hat\pi)$. This implies
    \begin{align}
     \label{Ancestorinstersection}
     \varphi_{l'}(\hat\pi)\in\RelAnc{l'+1}{\bar\rho'}\cap\RelAnc{l'+1}{\rho}.  
    \end{align}
    By definition, 
    $\hat\pi=\rho_j$ and 
    $\varphi_{l'}(\hat\pi)=\rho'_j$ for some $m_0\leq j \leq m_1$. Due
    to \ref{Ancestorinstersection}, 
    \mbox{$\rho'_i\in\RelAnc{4l'+3}{\bar\rho'}$} for all $m_0 \leq i \leq j$. 
    Since $\rho_1'\notin\RelAnc{4l'+3}{\bar\rho'}$,  $j\leq 0$. Thus, 
    $\rho_j\in\RelAnc{4l'+3}{\bar\rho}$ and 
    $\varphi_{l'}(\hat\pi)=\varphi_l(\hat\pi)$. Since
    $\varphi_l$ preserves the relevant ancestors of $\bar\rho$ level
    by level, we 
    obtain that $\hat\pi\in\RelAnc{l'+1}{\bar\rho}$. Since
    $\pi\in\RelAnc{l'+1}{\bar\rho}$, we obtain that 
    $\varphi_{l'}(\pi) = \varphi_l(\pi)$ and $\varphi_{l'}(\hat\pi) =
    \varphi_l(\hat\pi)$. Since
    $\varphi_l$ is an isomorphism, we conclude that $\pi*\hat\pi$ 
  \end{itemize}
  Thus, we have shown that $\varphi_{l'}$ is an isomorphism witnessing 
  $\bar\rho,\rho \RelAncequiv{l'}{n_2'}{z}{n_1'} \bar\rho', \rho'$. 
\end{proof}
Due to the iterated use of Proposition \ref{Prop:ConstructOneStep} in
the construction of  
$\rho'$ we can require the following further properties for $\rho'$:
\begin{corollary} \label{CorGrowthRateLocal}
  Let $n,z,l,n_1,n_2,l',n_1',n_2', \bar\rho, \bar\rho', \rho$ be as in
  the previous lemma. 
  Let 
  \begin{align*}
    &H:=\max\{\height(\pi): \pi\in \RelAnc{l}{\bar\rho'}\}\text{ and}\\
    &W:=\max\{\width(\pi): \pi\in \RelAnc{l}{\bar\rho'}\}.    
  \end{align*}
  Set 
  \begin{align*}
    g^{n,z}_{l',n_2}(x):= 
    \BoundHeightOnestepConstructionSimultanious(n \cdot 4^{l'+1}, x,
    n_2,z)    
  \end{align*}
  for 
  $\BoundHeightOnestepConstructionSimultanious$ the monotone function 
  defined in Proposition \ref{Prop:ConstructOneStep}. 
  
  We can construct a run $\rho'$ such that 
  \begin{align*}
    &\height(\rho')\leq
    (g^{n,z}_{l',n_2})^{(4^{l'+1})}(H),\\
    &\width(\rho') \leq W + 2l'+2
    \text{ and }\\
    &(\bar\rho,\rho) \RelAncequiv{l'}{n_2'}{z}{n_1'} (\bar\rho', \rho').
  \end{align*}
\end{corollary}
\begin{proof}
  By construction of $\rho_0'$, it is a relevant ancestor of
  $\bar\rho'$. 
  It follows that \mbox{$\height(\rho_0')\leq H$.}
  Then an easy induction shows that
  $\height(\rho_i')\leq (g^{n,z}_{l',n_2})^i(H)$ for all $i\leq m_1$:
  in each step, we either apply Proposition \ref{Prop:ConstructOneStep} or
  $\rho'_{i+1}$ is generated from $\rho'_i$ by applying a single
  stack operation $\op$. In the latter case we conclude by noting
  that 
  \mbox{$\height(\rho'_{i+1}) \leq \height(\rho'_i) + 1$.} 
  Note that $g$ is a monotone function. Since $m_1\leq
  4^{l'+1} \leq \lvert \RelAnc{l'+1}{\rho} \rvert$, we conclude that 
  \mbox{$\height(\rho')\leq (g^{n,z}_{l',n_2})^{(4^{l'+1})}(H)$.} 
  
  Now, consider the width of the $\rho_i'$. By assumption we know that
  $\width(\rho_i')\leq W$ for $m_0\leq i \leq 0$. Furthermore,
  as all $\rho_i'$ are  relevant $l'+1$-ancestors  of $\rho'$, their width
  differ in at most $2l'+2$. Therefore,
  $\width(\rho_i') \leq W + 2l' +2$  for
  all $m_0\leq i\leq m_1$. 
\end{proof}
\begin{remark}
  Note that the monotonicity of
  $\BoundHeightOnestepConstructionSimultanious$ carries over to the
  monotonicity
  of $g$ (in all parameters, i.e., in $n,z,l',n_2$, and $x$). 
\end{remark}


The previous lemmas showed that Duplicator can respond to local moves 
in such a way that she preserves 
isomorphisms of  relevant  ancestors while choosing small stacks.

Now, we deal with global moves of Spoiler. 
We present a strategy for Duplicator that allows to answer a global move
by choosing a run with the following property. 
Duplicator chooses a run such that the isomorphism of relevant
ancestors is preserved and such that all relevant ancestors of
Duplicator's choice end in small stacks. 
We split this proof into several lemmas. 
First, we address
the problem that Spoiler may choose an element  far away from
$\bar\rho$ but close to 
$\bar\rho'$. Then Duplicator has to find a run that has isomorphic
relevant ancestors but that is far away from $\bar\rho'$.

\begin{lemma} \label{LemmaGlobalStep1}
  Let $l, l', n, z, n_1, n_1', n_2, n_2'\in\N$ be numbers such that
  $l>3l'+3$, $z > n \cdot 4^l$,
  $z\geq 2$, 
  $n_1 > n_1' + 2 (l'+1)$, and $n_2 > n_2' + 4^{l'+1}$.
  Let $\bar\rho$ and $\bar\rho'$ be $n$-tuples of runs such that 
  \mbox{$\bar\rho \RelAncequiv{l}{n_2}{z}{n_1}
    \bar\rho'$}. Furthermore, let $\rho$ be a run such 
  that $\RelAnc{l'+1}{\bar\rho} \cap \RelAnc{l'+1}{\rho} = \emptyset$. Then
  there is some run 
  $\rho'$ such that $\bar\rho,\rho \RelAncequiv{l'}{n_2'}{z}{n_1'}
  \bar\rho',\rho'$. 
\end{lemma}
\begin{proof}
  We write $\varphi_l$ for the isomorphism witnessing
  $\RelAnc{l}{\bar\rho} \RelAncequiv{l}{n_2}{z}{n_1}
  \RelAnc{l}{\bar\rho'}$.  
  If 
  \begin{align*}
    \RelAnc{l'+1}{\bar\rho'} \cap \RelAnc{l'+1}{\rho} = \emptyset,    
  \end{align*}
  we can
  set $\rho':=\rho$ and we are done. 

  Otherwise, let $\pi^0_0 \prec \pi^0_1 \prec \dots \prec \pi^0_{n_0}$
  be an enumeration of all elements of 
  \mbox{$\RelAnc{l'+1}{\bar\rho'} \cap \RelAnc{l'+1}{\rho}$}. 
  Due to Corollary \ref{CorTouch}, 
  $\RelAnc{l'+1}{\rho} \cap\{ \pi:\pi\preceq \pi^0_{n_0}\} \subseteq 
  \RelAnc{3l'+3}{\bar\rho'}$. Since $l> 3l'+3$, we can set
  $\pi^1_i:=\varphi_{l}^{-1}(\pi^0_i)$ for all $0\leq i \leq n_0$. 
  Due to Lemma \ref{LemmaAncestorConstruction}, there
  is an 
  extension $\pi^1_{n_0}\prec \rho^1$ such that 
  $\RelAnc{l'+1}{\rho} \RelAncequiv{l'+1}{n_2'}{z}{n_1'}
  \RelAnc{l'+1}{\rho^1}$ and $\pi^1_i\in\RelAnc{l'+1}{\rho^1}$ for all $0\leq
  i \leq 
  n_0$. If $\RelAnc{l'+1}{\rho^1} \cap \RelAnc{l'+1}{\bar\rho'}= 
  \emptyset$ we set 
  $\rho':=\rho^1$ and we are done. 

  Otherwise we can repeat this process,
  defining $\pi^2_i:=\varphi_{l}^{-1}(\pi^1_i)$ for the maximal $n_1\leq
  n_0$ such that  
  $\pi^1_{i}\in \RelAnc{3l'+3}{\bar\rho'}$
  for all $0\leq i\leq n_0$. Then we extend this run to some run
  $\rho^2$. 
  If this process 
  terminates with the construction of some run $\rho^i$ such that 
  $\RelAnc{l'+1}{\rho^i} \cap \RelAnc{l'+1}{\bar\rho'} = \emptyset$,
  we
  set $\rho':=\rho^i$ and we are done. 
  If this is not the case, recall that $\RelAnc{3l'+3}{\bar\rho'}$ is
  finite. 
  Thus, we eventually reach the
  step were we have defined $\pi^0_0, \pi^1_0, \dots, \pi^m_0$ for some $m\in\N$
  such that for the first time $\pi^m_0 = \pi^i_0$ for some
  $i< m$.  But if $i>0$, then 
  \begin{align*}
    \pi^{m-1}_0 =\varphi_{l}(\pi^m_0)= \varphi_{l}(\pi^i_0)=\pi^{i-1}_0.  
  \end{align*}
  But this
  contradicts the minimality of 
  $m$. We conclude that $\pi^m_0=\pi^0_0$ which implies that 
  \mbox{$\pi^0_0\in \RelAnc{3l'+3}{\bar\rho}$.} 
  Furthermore, by definition 
  we have  $\pi^0_0\in \RelAnc{l'+1}{\rho}$ and
  there is a maximal $i$ such that  $\pi^0_i\in\RelAnc{3l'+3}{\bar\rho}$. 
  Since $z>\lvert \RelAnc{l}{\bar\rho} \rvert$, 
  we can apply Lemma
  \ref{LemmaAncestorConstruction} 
  and construct a chain $\varphi_l(\pi^0_i)\prec \rho_{i+1}' \prec
  \rho_{i+2}'\prec  \dots
  \prec \rho'$ such that
  $\bar\rho, \rho \RelAncequiv{l'}{n_2'}{z}{n_1'} \bar\rho',
  \rho'$. 
\end{proof}

The previous lemma showed that there is an answer to every global
challenge of Spoiler. 
In the following, we use the pumping constructions from Lemmas
\ref{lemmaPumpingTOP} -
\ref{lemmaPumpingWidth} in order to show that Duplicator may answer
global moves with  a $\rho'$ such that
\begin{align*}
  \RelAnc{l}{\bar\rho,\rho} \RelAncequiv{l'+1}{n_2'}{z}{n_1'}
  \RelAnc{l}{\bar\rho',\rho'}  
\end{align*}
and such that
$\RelAnc{l}{\bar\rho',\rho'}$ only contains runs that end in
small stacks.

Before we state this lemma, we have to give a precise notion of small
stacks. For this purpose, we introduce the following functions. 

\begin{definition}
  Let $l,l',n,n_1,n_2,z\in\N$ such that $l\geq 3l'+3$. 
  Set 
  \begin{align*}
    &\beta: \{-n_1, -n_1+1, \dots , 4^{(l'+1)}\}\times\N
    \rightarrow \N  \\ 
    &\beta(i,H):=
    \begin{cases}
      H + \ConstBoundHeightWord  + \FuncBoundTopWord(n_1+n_2+4^l,k,z)
      &\text{for } i= -n_1\\  
      \beta(i-1,H) + 
      \BoundHeightOnestepConstructionSimultanious(0,\beta(i-1,H),
      n_2+4^{l'+1}-i,k,z) &\text{otherwise,}  
    \end{cases}
  \end{align*} 
  and
  \begin{align*}
    \alpha(n_1,H, W, l'):=
    \max\left\{W,\FuncBoundWidthWord(\beta(-n_1,H))+n_1+2(l'+1)\right\}.
  \end{align*}
  where
  $\BoundHeightOnestepConstructionSimultanious, \FuncBoundTopWord,
  \FuncBoundWidthWord$ 
  the monotone functions from Proposition \ref{Prop:ConstructOneStep} and
  Lemmas 
  \ref{lemmaPumpingTOP} and \ref{lemmaPumpingWidth}, and
  $\ConstBoundHeightWord $ the constant from Lemma
  \ref{lemmaPumpingHeight}. 
\end{definition}

\begin{lemma} \label{LemmaGlobalStep2}
  Let $l,l',n,n_1,n_2,z\in\N$ such that $l\geq 3l'+3$. 
  Furthermore,   let $\bar\rho$ be an $n$-tuple of runs and $\rho$ a run such
  that $\RelAnc{l'+1}{\bar\rho} \cap \RelAnc{l'+1}{\rho}=\emptyset$. 
  Let $H,W\in\N$ be bounds such that 
  $\height(\pi)\leq H$ and
  $\width(\pi) \leq W$
  for all
  $\pi\in \RelAnc{l}{\bar\rho}$. 
  There is
  some run $\rho'$ such that
  \begin{align*}
    &\bar\rho,\rho \RelAncequiv{l'}{n_2}{z}{n_1} \bar\rho, \rho',\\
    &\height(\pi)\leq \beta(4^{l'+1},H), \text{ and}\\
    &\width(\pi)\leq \alpha(n_1,H,W, l')    
  \end{align*}
  for 
  all $\pi\in \RelAnc{l'}{\bar\rho,\rho'}$.
\end{lemma}
\begin{proof}
  Let $\rho_0 \prec \rho_1 \prec \dots \prec \rho_m := \rho$ be runs such that
  $\RelAnc{l'+1}{\rho}=\{\rho_i:0\leq i \leq m\}$. We have to find an
  isomorphic copy of $\RelAnc{l'+1}{\rho}$ consisting of small words
  but not intersecting with $\RelAnc{l'+1}{\bar\rho}$. 
  Using Lemmas \ref{LemmaConstructedRelAncEquiv},
  \ref{LemmaAncestorConstruction} and 
  \ref{Prop:ConstructOneStep}, we can construct such an isomorphic copy
  as soon as we find some small $\rho'_0$ with 
  $\rho'_0 \stackequiv{n_2 + 4^l}{z}{n_1}\rho_0$. Thus, as a first
  step we construct such a run $\rho'_0$. 

  Let $m_0 \geq -n_1$ be minimal such that there are runs
  $\rho_{m_0} \overset{+1}{\hookrightarrow}  \rho_{{m_0}+1}
  \overset{+1}{\hookrightarrow} \dots \overset{+1}{\hookrightarrow}
  \rho_0$. Note that by Lemma \ref{LemmaMinRelAnc} either $m_0=-n_1$
  or $\width(\rho_0)\leq n_1$.
  
  If $\height(\rho_{m_0}) \leq H$ and 
  $\width(\rho_{m_0}) \leq W$, then we
  choose $m_1$ 
  maximal such that 
  \mbox{$\height(\rho_i) \leq \beta(i,H)$} and
  $\width(\rho_i) \leq W$ 
  for all $m_0 \leq i \leq m_1$. 
  In this case, we set $\rho'_{m_1} := \rho_{m_1}$. 
  Otherwise, we set $m_1:=m_0$ and by Lemmas 
  \ref{lemmaPumpingTOP},\ref{lemmaPumpingHeight}, and
  \ref{lemmaPumpingWidth}, there is a run $\rho'_{m_1}$ such that
  \begin{align*}
    &\height(\rho'_{m_1})\leq \beta(-n_1,H) \leq \beta(m_1,H),\\
    &\width(\rho'_{m_1}) \leq \FuncBoundWidthWord(\beta(-n_1,H)) 
    \text{, and}\\
    &\rho_{m_1} \stackequiv{n_2 + n_1 +
      4^l}{z}{1} \rho'_{m_1}. 
  \end{align*}
  The last condition just says that $\TOP{2}(\rho_{m_1})
  \wordequiv{n_2+n_1+4^l}{z} \TOP{2}(\rho'_{m_1})$. 
  Furthermore, we construct $\rho'_{m_1}$ in such a way that either 
  $\height(\rho'_{m_1})>H$ or $\width(\rho'_{m_1})>W$. 
  
  Having constructed $\rho'_{m_1}$ according to one of the two cases, in both
  cases we continue  with the following construction. 
  Note that $\rho'_i=\rho_i$ for all $m_0\leq i \leq m_1$ or 
  $H < \height(\rho'_{m_1})$ or  
  \mbox{$W < \width (\rho'_{m_1})$.}

  By Proposition \ref{Prop:ConstructOneStep}, we can construct 
  $\rho'_{m_1} \prec \rho'_{m_1+1} \prec \dots \prec
  \rho'_m=:\rho'$ such that the following holds.
  \begin{enumerate}
  \item For $*\in\{\overset{+1}{\hookrightarrow},
    (\trans{\op})_{\op\in\Op}\}$ and for all $m_1 \leq i <
    m$, $\rho'_i \mathrel{*} \rho'_{i+1}$ iff $\rho_i \mathrel{*}
    \rho_{i+1}$.
  \item $\TOP{2}(\rho_{i}) \wordequiv {n_2}{z}\TOP{2}(\rho'_i)$. 
  \item $\height(\rho'_i) \leq \beta(i,H)$ and 
    $\width(\rho'_i) \leq \width(\rho'_{m_1})
    + n_1 + 2l' + 2$ for all $m_1
    \leq i \leq m$.
  \item $\rho'_i = \rho_i$ iff for all $m_1 \leq j \leq i$ we have
    $\height(\rho_j)\leq H$ and $\width(\rho_j) \leq W$ (this
    just requires to construct $\rho'_{m_0+1}$ such that 
    $\height(\rho'_{m_0+1})> H$ or $\width(\rho'_{m_0+1}) > W$. 
  \end{enumerate}
  From Lemma \ref{LemmaConstructedRelAncEquiv}, it follows that
  $\rho' \RelAncequiv{l'+1}{n_2}{z}{n_1} \rho$. 
  
  Furthermore, $\RelAnc{l'+1}{\rho'} \cap
  \RelAnc{l'+1}{\bar\rho}=\emptyset$:
  heading for a contradiction, assume that 
  \begin{align*}
    \rho'_i\in
    \RelAnc{l'+1}{\rho'}\cap\RelAnc{l'+1}{\bar\rho}    
  \end{align*}
 for some $0\leq i \leq m$. 
  Then $\rho'_j\in\RelAnc{3l'+3}{\bar\rho} \subseteq
  \RelAnc{l}{\bar\rho}$ for all 
  $0\leq j \leq i$. Thus, 
  \begin{align*}
    &\height(\rho'_j)\leq H \leq \beta(j,H)\text{ and}\\ 
    &\width(\rho'_j) \leq  W.  
  \end{align*}
   This implies that $\rho'_j=\rho_j$ 
  for all $m_0 \leq i \leq j$. But then
  $\rho_j=\rho'_j\in\RelAnc{l'+1}{\rho}\cap\RelAnc{l'+1}{\bar\rho}$
  which contradicts our 
  assumptions on $\bar\rho$ and $\rho$. 

  Hence, $\RelAnc{l'}{\bar\rho}$ and $\RelAnc{l'}{\rho'}$ do not touch whence 
  $\bar\rho,\rho \RelAncequiv{l'}{n_2}{z}{n_1} \bar\rho, \rho'$. 
\end{proof}

Combining the previous lemmas, we obtain a proof that for each
$n$-tuple in $\mathfrak{N}(\mathcal{S})$ there is an
$\FO{k}$-equivalent one such that the relevant 
ancestors of the second tuple only contain runs that end in small
stacks. This result is summarised in the following corollary. 

\begin{corollary}
  Let $\mathfrak{N}$ be a $2$-\HONPT. There are
  monotone functions 
  \begin{align*}
    &\BoundHeight:\N^5\to \N\text{ and}\\
    &\BoundWidth:\N^5 \to \N    
  \end{align*}
  such that the following holds. 
  Let $n,n_1',n_2',l'\in\N$. We set
  \begin{align*}
    &l := 4l' +5,\\    
    &n_1 := n_1'+2(l'+2)+1\text{ and}\\
    &n_2 := n_2' + 4^{l'+1}+1. 
  \end{align*}
  Let $z\in\N$ such that  $z\geq 2$ and $z> n\cdot 4^l$.  

  For all pairs of $n$-tuples $\bar\rho=\rho_1, \dots,
  \rho_n\in\mathfrak{N}$, $\bar\rho'=\rho'_1, \dots, \rho'_n\in
  \mathfrak{N}$ such that
  \begin{align*}
    &\height(\rho'_i) \leq \BoundHeight(n,z,l,n_1,n_2),\\
    &\width(\rho'_i) \leq \BoundWidth(n,z,l,n_1,n_2),\text{ and}\\
    &\bar\rho \RelAncequiv{l}{n_2}{z}{n_1} \bar\rho',
  \end{align*}
  and  for all runs $\rho\in \mathfrak{N}(\mathcal{S})$,
  there is a run $\rho'$ such that
  \begin{align*}
    &\bar\rho,\rho \RelAncequiv{l'}{n_2'}{z}{n_1'} \bar\rho', \rho', \\
    &\height(\rho')\leq \BoundHeight(n+1,z,l',n_1',n_2'),\text{ and}\\
    &\width(\rho')\leq \BoundWidth(n+1,z,l',n_1',n_2').    
  \end{align*}
\end{corollary}
\begin{proof}
  The proof is by induction on $n$. Assume that we have defined 
  $\BoundHeight(x_1,x_2,x_3,x_4,x_5)$ and
  $\BoundWidth(x_1,x_2,x_3,x_4,x_5)$ for all $x_2, \dots, x_5\in
  \N$ and $x_1\leq n$ such that for all tuples where  $x_1\leq n$ the
  claim holds. 
  For $\bar x_n:=(n,z,l,n_1,n_2)$ and
  $\bar x_{n+1}:=(n+1,z,l',n_1',n_2')$ , we set
  \begin{align*}
   &\BoundHeight(\bar x_{n+1}):=\max\left\{
  \beta(4^{l'+1},\BoundHeight(\bar x_n)), 
  (g^{n,z}_{l',n_2})^{4^{l'+1}}(\BoundHeight(\bar x_n))
  \right\} \text{ and}\\
  &\BoundWidth(\bar x_{n+1}) :=\max\left\{
  \FuncBoundWidthWord\left(\beta(4^{(l'+1)},\BoundHeight(\bar
  x))+n_1'+2l'+2\right),  
  \BoundWidth(\bar x_n)+2l'+2\right\}
  \end{align*}
  where 
  $g^{n,z}_{l',n_2}$ is the function from
  Corollary \ref{CorGrowthRateLocal} and 
  $\beta$ the function from Lemma \ref{LemmaGlobalStep2}.
   The following
  case distinction proves that this definition satisfies the claim. 
  \begin{enumerate}
  \item 
    First assume that $\RelAnc{l'+1}{\rho}\cap
    \RelAnc{l'+1}{\bar\rho}\neq\emptyset$. 
    Then we can apply Lemma \ref{LemmaLocalStep} and obtain an element
    $\rho'\in \mathfrak{N}(\mathcal{S})$ such that
    $\bar\rho,\rho \RelAncequiv{l'}{n_2'}{z}{n_1'} \bar\rho',
    \rho'$. Furthermore, by 
    Corollary \ref{CorGrowthRateLocal}, $\rho'$ can be chosen such that
    \begin{align*}
      &\height(\rho')\leq
      (g^{n,z}_{l',n_2})^{4^{l'+1}}(\BoundHeight(n,z,l,n_1,n_2)) =
      (g^{n,z}_{l',n_2})^{4^{l'+1}}(\BoundHeight(\bar x_n)) \text{ and}\\
      &\width(\rho')\leq
      \BoundWidth(\bar x_n)+2l'+2.      
    \end{align*}
  \item 
    Otherwise, 
    $\RelAnc{l'+1}{\rho}\cap \RelAnc{l'+1}{\bar\rho}=\emptyset$.
    Since $l> 4l' +4$,  we can apply Lemmas 
    \ref{LemmaGlobalStep1} and  \ref{LemmaGlobalStep2} and obtain
    $\rho'\in \mathfrak{N}(\mathcal{S})$ such that
    \begin{align*}
      &\bar\rho,\rho \RelAncequiv{l'}{n_2'}{z}{n_1'} \bar\rho',
      \rho', \\
      &\height(\rho') \leq \beta\left(4^{l'+1},
      \BoundHeight(n,z,l,n_1,n_2)\right)=
      \beta\left(4^{l'+1},
      \BoundHeight(\bar x_n)\right), \text{ and}\\
      &\width(\rho') \leq \FuncBoundWidthWord\left(\beta(4^{(l'+1)},
      \BoundHeight(\bar x_n)) + n_1' + 2(l'+2)\right).
    \end{align*}
  \end{enumerate}    
  By induction, our definition satisfies the claim. Note that the
  monotonicity of $\BoundHeight$ and $\BoundWidth$ follows from the
  monotonicity of all the components involved in the definition.  
\end{proof}

\paragraph{Strategy for Bounding the Length of Runs}
\label{SubsectionRunBounds}

For each relevant ancestor set, there is an equivalent one which 
only contains runs that end in small stacks. But the runs leading to
these stacks can still be arbitrary long. In the next lemmas, we show
that we can also bound the length. For this proof, the
Corollaries \ref{Cor:GlobalBoundRun} and \ref{Cor:GlobalBoundRun2}
are important tools because they allow to replace long runs by shorter
ones. 

\begin{lemma} \label{Lemma:BoundingFunctions}
  Let $\mathcal{N}$ be a level $2$ pushdown system defining the higher
  order nested 
  pushdown tree 
  $\mathfrak{N}:=\NPT(\mathcal{N})$.  
  We can compute a function $\BoundRunLength:\N^5\rightarrow
  \N$ such that the following hold:

  Let $n,z,n_1',n_2',l'\in\N$,  
  $l := 4l' +5, n_1 := n_1'+2(l'+2)+1,$ and
  $n_2 := n_2' + 4^{l'+1}+1$ such that $z\geq 2$ and $z> n\cdot 4^l$.
  Furthermore, let $\bar\rho$ and $\bar\rho'$ be $n$-tuples of runs of
  $\mathfrak{N}$ such that  
  \begin{enumerate}
  \item $\bar\rho \RelAncequiv{l}{n_2}{z}{n_1} \bar\rho'$, and
  \item $\length(\pi)\leq\BoundRunLength(n,l,n_1,n_2,z)$
    for all $\pi\in\RelAnc{l}{\bar\rho'}$,
  \item $\height(\pi) \leq  \BoundHeight(n,z,l,n_1,n_2)$ for all
    $\pi\in\RelAnc{l}{\bar\rho'}$, and
  \item $\width(\pi) \leq \BoundWidth(n,z,l,n_1,n_2)$ for all
    $\pi\in\RelAnc{l}{\bar\rho'}$. 
  \end{enumerate}
  For each $\rho\in \mathfrak{N}$ there is some $\rho'\in
  \mathfrak{N}$ such that 
  \begin{enumerate}
  \item $\bar\rho, \rho\RelAncequiv{l'}{n_2'}{z}{n_1'} \bar\rho', \rho'$,
  \item $\length(\pi)\leq\BoundRunLength(n+1,l',n_1',n_2',z)$
    for all $\pi\in\RelAnc{l'}{\bar\rho',\rho'}$,
  \item $\height(\pi) \leq  \BoundHeight(n+1,z,l',n_1',n_2')$ for all
    $\pi\in\RelAnc{l}{\bar\rho',\rho'}$, and
  \item $\width(\pi) \leq \BoundWidth(n+1,z,l',n_1',n_2')$ for all
    $\pi\in\RelAnc{l}{\bar\rho',\rho'}$. 
  \end{enumerate}
\end{lemma}
\begin{proof}
  Using the Lemmas  \ref{LemmaLocalStep} -- \ref{LemmaGlobalStep2}, we
  find some candidate $\rho'$ such that
  \begin{align*}
    \bar\rho,\rho \RelAncequiv{l'}{n_2'}{z}{n_1'} \bar\rho', \hat \rho'    
  \end{align*}
  and the
  height and width 
  of the last stacks of all $\pi\in \RelAnc{l'+1}{\hat \rho'}$ are bounded by
  \mbox{$\BoundHeight(n+1,z,l',n_1',n_2')$} and 
  $\BoundWidth(n+1,z,l',n_1',n_2')$, respectively. 

  Recall that there is a chain
  $\hat \rho_0'\prec 
  \hat \rho_1' \prec \dots \prec \hat \rho_m'=\hat \rho'$ for some 
  $0\leq m\leq 4^{(l'+1)}$ 
  with  
  $\hat \rho_i'\trans{\gamma} \hat \rho_{i+1}'$ or 
  $\hat \rho_i' \overset{+1}{\hookrightarrow} \hat \rho_{i+1}'$ for all
  $0\leq i < m$ 
  such that $\RelAnc{l'+1}{\hat \rho'} 
  = \{\hat \rho_i': 0 \leq i \leq m\}$.

  If $\hat \rho_0'\notin\RelAnc{3l'+3}{\bar\rho'}$, then
  we can use Corollary \ref{Cor:GlobalBoundRun} and choose some $\rho_0'$  
  that ends in the same configuration as $\hat \rho_0'$ such
  that $\rho_0'\notin \RelAnc{3l'+3}{\bar\rho'}$ and 
  \begin{align*}
    \length(\rho_0')\leq 1+ &2\cdot \BoundHeight(n+1,z,l',n_1',n_2') \cdot
    \BoundWidth(n+1, z, l', n_1', n_2'))\\ &\cdot 
    (1+ \FuncBoundLoopLength{\mathcal{N}}{z}(\BoundHeight(n+1, z, l',
    n_1', n_2'))).     
  \end{align*}

  If $\hat \rho_0'\in\RelAnc{3l'+3}{\bar\rho'}$ 
  let $0 \leq i \leq m$ be maximal such that 
  $\hat \rho_i'\in\RelAnc{l}{\bar\rho'}$.
  In this case let $\rho_j' :=\hat \rho_j'$ for all $0\leq j \leq i$. 
  
  By now, we have obtained a chain $\rho_0' \prec \rho_1' \prec \dots \prec
  \rho_i'$ for some $0\leq i \leq m$.   
  Using the previous lemma, we can extend this chain to a chain
  $\{\rho_i': 0\leq i \leq m\}$ such that 
  \begin{enumerate}
  \item $\rho'_i(\length(\rho'_i)) = \hat\rho_i(\length(\hat\rho_i))$,
  \item   $\rho_i'\trans{\gamma} \rho_{i+1}'$  iff  
    $\hat \rho_i'\trans{\gamma} \hat \rho_{i+1}'$   for all $0\leq i < m$,
  \item $\rho_i' \overset{+1}{\hookrightarrow} \rho_{i+1}'$  iff 
    $\hat \rho_i' \overset{+1}{\hookrightarrow} \hat \rho_{i+1}'$ for all
    $0\leq i < m$,
  \item   $\length(\rho_{i+1}') \leq \length(\rho_i') + 
    2 \cdot\BoundHeight(n+1,z,l',n_1',n_2') \cdot 
    (1+ \FuncBoundLoopLength{\mathcal{N}}{z}(
    \BoundHeight(n+1, z, l', n_1',
    n_2')))$, and
  \item  $\hat \rho_j'\in\RelAnc{3l'+3}{\bar\rho'}$  for
    all $0\leq j \leq i$ implies  $\rho_i'=\hat \rho_i'$.
  \end{enumerate}
  Assume that we have constructed the chain up to 
  $\rho_0' \prec \dots \prec \rho_{m_0}'$ for some $m_0 < m$. 
  Note that  $\hat \rho_{m_0+1}\notin
  \RelAnc{3l'+3}{\bar\rho'}$ by definition of the initial segment of
  the $\rho_i'$.  
  We can use Corollary \ref{Cor:GlobalBoundRun2} in order to construct
  $\rho_{m_0+1}'$ as 
  required. In this construction, we can enforce
  $\rho_{m_0+1}'\notin\RelAnc{3l'+3}{\bar\rho'}$ if $\rho_{m_0}'=\hat
  \rho_{m_0}'$.   

  Using Lemma \ref{LemmaConstructedRelAncEquiv}, we conclude that $\rho'
  \RelAncequiv{l'+1}{n_2'}{z}{n_1'} 
  \hat \rho'$ for $\rho':=\rho_m'$. Furthermore, we claim that 
  $\RelAnc{l'+1}{\bar\rho'} \cap \RelAnc{l'+1}{\hat \rho'} =   
  \RelAnc{l'+1}{\bar\rho'} \cap \RelAnc{l'+1}{\rho'}$. 
  By definition the inclusion from left to right is clear. For the
  other direction, assume that
  there is some element 
  \mbox{$\rho_i'\in \RelAnc{l'+1}{\bar\rho'} \cap
  \RelAnc{l'+1}{\rho'}$.}
  By Lemma \ref{CorRelAncDistBound}, this implies that 
  $\rho_j'\in \RelAnc{3l'+3}{\bar\rho'}$ 
  for all $0\leq j \leq i$. Thus, $\rho_i'=\hat \rho_i'$, which implies that
  $\rho_i'\in \RelAnc{l'+1}{\bar\rho'} \cap \RelAnc{l'+1}{\rho'}$. 
  
  We conclude that $\bar\rho, \rho 
  \RelAncequiv{l'}{n_2'}{z}{n_1'} \bar\rho', \hat \rho'
  \RelAncequiv{l'}{n_2'}{z}{n_1'} \bar\rho', \rho'$ because 
  \mbox{$\RelAnc{l'+1}{\bar\rho'}\cap \RelAnc{l'+1}{\hat
    \rho'}$} is isomorphic to   
  $\RelAnc{l'+1}{\bar\rho'}\cap\RelAnc{l'+1}{\hat \rho'}$. 
  By definition, the length of $\rho'$ is bounded by a polynomial in 
  \begin{align*}
    &\BoundHeight(n+1, z, l', n_1', n_2'),\\ 
    &\BoundWidth(n+1, z, l', n_1', n_2'),\\
    &\FuncBoundLoopLength{\mathcal{N}}{z}(\BoundHeight(n+1, z, l',
    n_1', n_2')),\text{ and}\\
    &\BoundRunLength(n,l,n_1,n_2,z).    
  \end{align*} 
  This polynomial can be used to inductively define 
  $\BoundRunLength(n+1, l', n_1', n_2', z)$.  
\end{proof}
Note that the previous lemma completes the proof of Proposition 
\ref{Prop:2NPT-Strategy}.

\subsection{FO Model Checking Algorithm for Level 2  Nested Pushdown
  Trees} 
\label{SectionFODecidability}

In the previous section, we have shown that each existential
quantification on a nested pushdown tree
$\mathfrak{N}:=\NPT(\mathcal{N})$ can be 
witnessed by a run 
$\rho$ of small length. Even when we add parameters $\rho_1, \dots,
\rho_n$
this result still holds, in the sense that there is a witness $\rho$ of
small length compared to the length of the parameters. Hence, we can decide
first-order logic on level $2$ nested pushdown trees with the
following algorithm.

\begin{enumerate}
\item Given the pushdown system $\mathcal{N}$ and a
  first-order formula $\varphi$, the algorithm first computes the
  quantifier rank $q$ of $\varphi$.
\item Then it computes numbers $z,l^1,l^2,l^3,\dots, l^q, n_1^1, n_1^2,
  n_1^3 \dots, n_1^q, n_2^1, n_2^2, n_2^3, \dots, n_2^q\in\N$ such
  that for each $i< q$ the numbers $z,l^i,l^{i+1},n_1^i, n_1^{i+1},
  n_2^i, n_2^{i+1}$ can be used as parameters in 
  Proposition \ref{Prop:2NPT-Strategy}.
%
\item These numbers define a constraint $S=(S^{\mathfrak{N}}(i))_{i\leq
    q}$ for Duplicator's strategy in the $q$-round game on
  $\mathfrak{N}$ and $\mathfrak{N}$ as follows. 
  We set $(\rho_1, \rho_2, \dots, \rho_m)\in S^{\mathfrak{N}}_m$
  if for each $i\leq m$ and $\pi\in\RelAnc{l_i}{\rho_i}$
  \begin{align*}
    &\length(\pi)\leq \BoundRunLength(i,l^i,n_1^i,n_2^i,z), \\
    &\height(\pi)\leq \BoundHeight(i,z,l^i,n_1^i,n_2^i), \text{ and}\\
    &\width(\pi) \leq \BoundWidth(i,z,l^i,n_1^i,n_2^i). 
  \end{align*}
\item Due to Lemma \ref{Lemma:BoundingFunctions}, Duplicator has an
  $S$-preserving strategy in the $q$-round game on $\mathfrak{N}$
  and $\mathfrak{N}$. Thus, applying the algorithm SModelCheck
  (cf. Algorithm \ref{AlgoSPReservingModelCheck} in Section
  \ref{Sec:EFGame}) decides whether $\mathfrak{N}\models \varphi$. 
\end{enumerate}

\paragraph{Complexity of the Algorithm}
For the case of nested pushdown trees (of level $1$) our approach
resulted in an 
$2$-EXPSPACE $\FO{}$ model checking algorithm. In the case of level $2$
nested pushdown trees, we
cannot prove such a nice result. 
At the moment, we cannot prove an elementary complexity bound for the
$\FO{}$ model checking on $2$-\HONPT because we cannot
determine the length of short loops.
Our algorithm can only be efficient if we have a good bound on the
length of the $k$ shortest loops of any stack because we use loops as a
main ingredient in the construction of equivalent relevant ancestors. 
But such a good bound is not known to exist. 
We do not know any elementary algorithm that, given a level $2$
pushdown system $\mathcal{N}$ and a 
number $k$, calculates the shortest $k$ loops
from $(q_0, \bot)$ to $(q_1, \bot)$ of $\mathcal{N}$. 
The underlying problem is that we cannot derive an elementary bound on the
length of such loops. The best bound we know can be derived as follows.

From Hayashi's pumping lemma for indexed grammars \cite{Hayashi73}, we
can derive that  
the shortest loop of $\mathcal{N}$ has size 
$\exp(\exp(\exp(p(\lvert \mathcal{N}\rvert))))$ for some polynomial
$p$. Unfortunately, Hayashi's pumping 
lemma does not yield any bound on the second shortest loop. Thus, the
only known way of calculating the second shortest loop is to design a
copy of the pushdown system which simulates the first one but avoids
this first 
loop. This involves increasing the number of states by the length
of the shortest loop, i.e., we design a system $\mathcal{N}'$ with $\lvert
\mathcal{N}'\rvert \approx \exp(\exp(\exp(p(\lvert \mathcal{N} \rvert))))$ many
states. Using this system 
we obtain a $6$-fold exponential bound in $\lvert \mathcal{N} \rvert$ for the
second shortest loop of 
$\mathcal{N}$ (which is the shortest one of $\mathcal{N}'$) the same
way as we obtained 
the bound for the first loop. Thus, 
the best bound known for the $k$ shortest loops is an exponential
tower of height $3k$ in the size of the pushdown system. 
But it is quite clear that there are level $2$ nested pushdown trees
where we can define the existence of  $k$ 
loops from $(q_0, \bot)$ to $(q_1,\bot)$ by a first-order formula of
quantifier rank linear in $k$. Thus, our model checking algorithm would have
to choose $k$ short loops. 
Given the bounds on short loops, we expect that our algorithm then
needs space up to a tower of exponentials of height $3k$ in order to verify
this formula. Since $k$ is arbitrary, 
the algorithm has nonelementary space
consumption in the quantifier rank of the formula. 

It remains open to determine the exact complexity of our algorithm. 
We neither know whether our algorithm has elementary complexity nor do
we know a good lower bound on the complexity of model checking on
nested pushdown trees of level $2$. These questions require further
study. 



\section{Decidability of Ramsey Quantifiers on Tree-Automatic
  Structures} 

\label{SectionRamseyQuantifier} 
Recently, Rubin \cite{Rubin2008} proved the decidability of 
Ramsey quantifiers on 
string-automatic structures using the concept of word-combs. In this
section we will lift his 
techniques to the tree-case, i.e., we prove the decidability of 
Ramsey quantifiers on (tree-)automatic structures. 
Actually, our proof can also be seen as an adaption of
To's and Libkin's proof \cite{ToL08}
of the decidability of the recurrent reachability problem on
automatic structures. Nevertheless, our result was developed
independently 
from To's and Libkin's work. 

Let us briefly recall Rubin's ideas. 
His main tool is the concept of a word-comb. A 
word-comb is an infinite sequence of finite $\Sigma$-words 
such that
there is a sequence of natural numbers $g_1 < g_2 < g_3 <
\dots$ 
such that all but the shortest $n$ words of
the word-comb agree on the first $g_n$ letters.
A word-comb can be represented using infinite words as
follows.  
Let $w_1\in\Sigma^\omega$, $w_2\in(\Sigma\cup \{\Box\})^\omega$ be
infinite words and $G\subseteq \N$ an infinite set. 
A finite word $w$ belongs to the word-comb represented by $(w_1,
w_2, G)$ if the following holds: 
$w$ decomposes as $w=v_1\circ v_2$ where $v_1$ is a prefix of $w_1$
and $v_2$ is a subword of $w_2$ such that 
\begin{enumerate}
\item $\lvert v_1 \rvert \in G$,
\item there is some $k\in \N$ such that $v_2\Box^k$ is the subword of
  $w_2$ induced by the $(\lvert v_1 \rvert+1)$-st to the 
  $(\lvert w \rvert+k)$-th letter of $w_2$ such that $\lvert w \rvert
  +k$ is the successor of $\lvert v_1\rvert$ in $G$. 
\end{enumerate}
Figure \ref{fig:Wordcomb} illustrates such a representation of a
word-comb. 

Now, we explain how the notion of a word-comb can be used to decide
Ramsey quantifiers on string-automatic structures. 
The first important observation is that every infinite set of finite
words contains a subset which is a word-comb, i.e., a subset that can
be represented by some triple $(w_1, w_2, G)$ as explained above. 
Secondly, $\omega$-string-automata can be used to 
extract the words of the word-comb from the representation. 

Recall that the Ramsey quantifier asserts the existence of an infinite
subset that is homogeneous with respect to a certain formula $\varphi$, i.e.,
all pairwise distinct $n$-tuples from this set satisfy $\varphi$. 
Now, for each string-automatic structure $\mathfrak{A}$,  
this can be translated into the
assertion that there is a representation of a word-comb such that
each pairwise distinct $n$-tuple from the comb satisfies
$\varphi$. This assertion 
can be formulated in a first-order formula $\varphi'$ on a certain
$\omega$-string-automatic extension $\mathfrak{A'}$ of $\mathfrak{A}$.
This extension $\mathfrak{A'}$ enriches $\mathfrak{A}$ by those
infinite strings that occur in the representation of word-combs.  
The classical correspondence between first-order logic on
$\omega$-string-automatic structures and $\omega$-string-automata
yields an $\omega$-string-automaton that represents $\varphi'$ on
$\mathfrak{A'}$. Finally, this $\omega$-string-automaton can be turned into a
string-automaton that represents $\varphi$ on $\mathfrak{A}$. 

This idea carries even further. 
Kuske \cite{Kuske2009} introduced
a logic which he  calls $\mathrm{FSO}$. The Formulas of $\mathrm{FSO}$
are formed according to the formation rules of first-order logic
and the following two rules. 
First, one may use variables for $n$-ary relations, i.e., for $X$ an
$n$-ary relation variable and $x_1, x_2, \dots, x_n$ element
variables, $Xx_1x_2\dots x_n$ is an atomic formula of $\mathrm{FSO}$. 
Second, for $X$ a relation variable that only occurs negatively in
some $\varphi\in\mathrm{FSO}$, $\exists X$ inf. $\varphi$ is in
$\mathrm{FSO}$. This formula is satisfied if there is an infinite
interpretation for $X$ that satisfies $\varphi$. 
$\mathrm{FSO}$ is a generalisation of $\FO{}((\RamQ{n}){n\in\N})$ as
follows.
$\RamQ{n}\bar x \varphi$ is equivalent to $\exists X (\forall x_1,
\dots x_n ( \bigwedge_{1\leq i \leq n} x_i\in X) \rightarrow
\varphi)$. 
On string-automatic structures, 
Rubin's technique generalises to $\mathrm{FSO}$: analogously
to the decidability of the Ramsey quantifier, one obtains the
decidability of $\mathrm{FSO}$ on string-automatic structures. 
The reason why this result extends to $\mathrm{FSO}$ is a 
closure under subsets of witnesses for $\mathrm{FSO}$ formulas: 
if an $\mathrm{FSO}$ formula $\varphi$ asserts the existence of some
infinite set $X$ then $X$ appears only negatively in the subformulas
of $\varphi$. Without loss of generality this means that the only
occurrences of $X$ in subformulas of $\varphi$ are of the form $x\notin
X$. If $A$ is an infinite set witnessing the assertion of $\varphi$,
then any infinite subset $A'\subseteq A$ also witnesses the
statement $x\notin A'$ if $A$ witnesses $x\notin A$. Thus, taking an
infinite subset 
of some witness of a formula in  $\mathrm{FSO}$ is again a witness
of this formula. Since every infinite subset of a set of words
contains a word-comb, it suffices to look for witnesses of
$\mathrm{FSO}$ formulas among the word-combs. Hence, Rubin's reduction
works also for $\mathrm{FSO}$.  

\begin{figure}
  \begin{xy}
    \xymatrix{
      w_1 = && a & a & a & a & a & a & a & a & \dots\\
      w_2 = && b & \Box & b & \Box & b & \Box & b & \Box& \dots\\
      G =\{& 0, && 2, && 4, && 6, &&  \dots & \}\\ 
      \\
      w_1 = && a & a & a & a & a & a & a & a & \dots\\
      w_2 = && b & \Box& b & \Box & b & \Box & b & \Box & \dots\\
      \save "6,3"."6,4"*+[F]\frm{}
      &&b\\
      \\
      w_1 = && a & a & a & a & a & a & a & a & \dots\\
      w_2 = && b & \Box & b & \Box & b & \Box & b & \Box & \dots\\
      \save "9,3"."9,4"*+[F]\frm{}
      \save "10,5"."10,6"*+[F]\frm{}
      && a & a &  b\\
      \\
      w_1 = && a & a & a & a & a & a & a & a & \dots\\
      w_2 = && b & \Box & b & \Box & b & \Box & b & \Box & \dots\\
      \save "13,3"."13,6"*+[F]\frm{}
      \save "14,7"."14,8"*+[F]\frm{}
      && a & a &  a & a & b\\
      \vdots
    }
  \end{xy}
  \caption{Word-comb $(w_1, w_2, G)$ encoding the set $\{a^{2n}b:
    n\in\N\}$.} 
  \label{fig:Wordcomb}
\end{figure}

Our goal is to lift the concept of a comb from strings to finite trees.
We use
three infinite trees for representing an infinite set of finite
trees. Unfortunately, the correspondence we obtain is not as tight as
in the string case:  each infinite set of
finite words contains a word-comb that is represented by some triple $(w_1,
w_2, G)$. Furthermore,  
there is an $\omega$-string-automaton that decides, on input some finite word
$w$ and the representation $(w_1, w_2, G)$ whether $w$ is contained in
the word-comb. 
The notion of word-combs smoothly generalises to the notion of
tree-combs. Unfortunately, tree-combs do not form
$\omega$-tree-regular sets. 
This makes the tree case  more involved. 

The outline of our proof is as follows. Given an
infinite set of finite trees,  there is an infinite subset
called a \emph{tree-comb}. 
A tree-comb is an infinite set that allows a unique representation as
a triple 
$(T_1, T_2, G)$ where $T_1$ and $T_2$ are infinite trees and
$G\subseteq\{0,1\}^*$. 
We then define an $\omega$-automaton that extracts finite trees
from the representation of a tree-comb. The set of all these trees 
is called the \emph{closure} of the tree-comb. 
The connection between a tree-comb and its closure is as follows. 
Firstly, every tree-comb is
contained in its closure. Secondly, each tree $T$ contained in the closure
is locally equal to the trees in the tree-comb:
given an arbitrary infinite branch, there
is a tree $t'$ in the tree-comb such that $t$ and $t'$ coincide along
this infinite branch.  

In order to decide Ramsey quantifiers on automatic structures, 
we first prove that each Ramsey quantifier is witnessed by the
closure of some tree-comb. 
In order to explain the single steps of this proof, 
we fix a formula $\varphi\in\FO{}$ and consider 
the formula
\begin{align*}
  \RamQ{n}\bar x(\varphi).   
\end{align*}
We fix an automatic structure $\mathfrak{A}$. On $\mathfrak{A}$,
$\varphi$ corresponds to some automaton $\mathcal{A}_\varphi$. 
We prove the following.
\begin{enumerate}
\item A straightforward generalisation of the string case shows that
  $\mathfrak{A}\models\RamQ{n}\bar x(\varphi)$ if and only if there is
  a tree-comb $C$ witnessing this Ramsey quantifier on
  $\mathfrak{A}$.
\item We show that $C$ can be chosen to be \emph{homogeneous} with respect
  to $\mathcal{A}_\varphi$. Roughly speaking, homogeneity means
  that the runs of $\mathcal{A}_\varphi$ on all pairwise distinct
  $n$-tuples from $C$ look similar.
\item For a homogeneous $C$, we show that $\mathcal{A}_\varphi$
  accepts all $n$-tuples from the closure of $C$. The proof idea of
  this step is as follows. Since each $n$-tuple $\bar c$ from the
  closure is locally equal to $n$-tuples from $C$, we can locally copy
  the accepting runs of $\mathcal{A}_\varphi$ on
  the latter 
  tuples and obtain  a function defined on the domain of $\bar
  c$. Since all runs that we locally copy are similar, this function
  turns out to be a run of $\mathcal{A}_\varphi$ on $\bar c$. 
  Since it is composed from accepting runs, it is also accepting.
\end{enumerate}

Putting these steps together, we obtain that $\mathfrak{A}$
satisfies some Ramsey quantifier if and only if there is a
closure of some tree-comb witnessing this quantifier. 
The proof of the decidability of
Ramsey quantifiers on $\mathfrak{A}$ continues
analogously to the string case.
We obtain an $\omega$-automatic
extension of $\mathfrak{A}$. On this extension, the existence of the
closure of a tree-comb that witnesses the Ramsey quantifier is
expressible in first-order
logic. The resulting first-order formula is turned
into an $\omega$-automaton using standard techniques. 
This $\omega$-automaton can then be turned into an automaton
corresponding to $\RamQ{n}\bar x(\varphi)$ on 
$\mathfrak{A}$.   
Thus, for any formula in
$\FO{}(\exists^\mathrm{mod},(\RamQ{n})_{n\in\N})$
and any automatic structure $\mathfrak{A}$, there is an automaton
$\mathcal{A}_\varphi$ that corresponds to $\varphi$ on
$\mathfrak{A}$. 
Unfortunately, this approach does not extend directly to Kuske's logic
$\mathrm{FSO}$. Thus, it remains an open problem whether
$\mathrm{FSO}$ is decidable on all tree-automatic structures.

\begin{remark}
  As already indicated, 
  we deal with finite and infinite trees in this section. 
  Because of this, we deviate from our notational conventions in the
  following way. 
  Throughout Section \ref{SectionRamseyQuantifier}, where we have
  to distinguish between infinite and finite trees, we write
  ``tree'' for an object that is either a finite or an infinite tree,
  i.e., a $\Sigma$-tree is an element of 
  $\allTrees{\Sigma} =\Trees{\Sigma}\cup\infTrees{\Sigma}$. 
  Thus, whenever we want to consider an element of $\Trees{\Sigma}$,
  we will explicitly write \emph{finite tree}. 
\end{remark}

\subsection{Tree-Combs}
Recall that  Ramsey quantifiers allow a restricted
form of second-order quantification. In order to translate these
quantifiers over an automatic structure into
first-order quantifiers over an $\omega$-automatic structure, we want
to represent infinite sets of finite trees by 
a tuple of infinite trees. 

In Definition \ref{def:tree-comb}, we formally introduce tree-combs. 
Before, we develop some machinery that allows to extract finite trees
from a tuple of infinite trees. 
This machinery is not necessary for understanding the definition of
tree-combs, but it is used to define the closure of a tree-comb. 
Since our interest is in the relationship of tree-combs and their
closures, we postpone the definition of tree-combs. 

In the following, we write
$\Sigma_\Box$ for $\Sigma\cup\{\Box\}$ where
$\Box\notin\Sigma$ is some new symbol. 

Recall that we defined the  following notation. If $t$ is a
$\Sigma$-labelled tree, 
then we denote by $t^\Box$ the full binary tree which consists of $t$
padded by $\Box$-labels. 
We define a kind of inverse to this operation which returns the
maximal $\Sigma$-labelled tree contained in a given
$\Sigma_\Box$-labelled tree. 

\begin{definition}
  Let $T\in \allTrees{\Sigma_{\Box}}$ be an arbitrary tree. Then
  $\prune(T)$ 
  denotes the maximal initial segment of $T$ that is in $\allTrees{\Sigma}$. 
\end{definition}
\begin{remark}
  We stress that $\prune$ yields a $\Sigma$-labelled tree from a
  $\Sigma_{\Box}$-labelled tree. This is done by extracting the
  initial segment up to the first occurrence of $\Box$ along each
  branch. In this sense, $\Box$-labelled positions in $T$ mark
  undefined positions in the domain of $\prune(T)$. 
\end{remark}
Recall that we extract an element of a word-comb from its
representation $(w_1, w_2, G)$ by taking the prefix $v_1$ of $w_1$ of length
$g_1$ for some $g_1\in G$ and appending a subword $v_2$ of $w_2$. 
$v_2$ consists of  the $(g_1+1)$-st to the $g_2$-th letter of $w_2$ where
$g_2$ is the direct successor of $g_1$ in $G$. 
The function $\prune$ will be used to extract an analogue of
$v_1$ in the tree-case. Now, we define another
function, called $\extract$,  that is the analogue to the extraction
of $v_2$. 
It extracts the $\Sigma$-labelled subtree of an infinite tree from a
given position up to the first occurrence of an element from $G$ along
each branch.  

In the string case, we obtain an element encoded in $(w_1, w_2, G)$ by
composition of $v_1$ and $v_2$. Analogously, after defining $\extract$
we need a kind of composition of $\prune$ and $\extract$ which
extracts a tree from a triple $(T_1, T_2, G)$. This composition
is a function called $\pretree$.

\begin{definition}
  Let $T:\{0,1\}^*\to\Sigma_\Box$, $G\subseteq\{0,1\}^*$ and
  $e\in\{0,1\}^*$. Then 
  $\extract(e,T,G)$ is the maximal initial segment of 
  $\inducedTreeof{e}{T}$ (the subtree of $T$ rooted at $e$)
  such that the following two conditions are satisfied. 
  \begin{itemize}
  \item $\extract(e,T,G)$ is a $\Sigma$-labelled tree, i.e.,
    it does not contain $\Box$-labelled nodes.
  \item For all $d\in \domain(\extract(e,T,G))$, $ed\in G$ implies
    $d=\varepsilon$. 
  \end{itemize}
\end{definition}
\begin{remark}
  Note that $\Box$-labelled nodes in $T$ mark again positions that are
  undefined in the domain of $\extract(e,T,G)$. 
\end{remark}
Note that  $\extract(T,G,e)$ is the empty tree if and only if
$T(e)=\Box$. 
Furthermore, it is a finite tree if every branch starting at $e$
contains a node 
$e\leq e'$ with $T(e')=\Box$ or $e'\in G$. If it is a finite tree,
then it is a $\Sigma$-labelled finite tree by the very definition.

Next, we define  $\pretree$. 
In general, $\pretree$ may extract infinite
trees from a triple $(T_1, T_2, G)$. But later we use it only on
inputs where it extracts finite trees.  

Recall that we write $H^+$ for the  border of a tree-domain $H$, 
i.e., $H^+$ is the set of minimal elements of $\{0,1\}^*\setminus H$. 
\begin{definition}
  Let $T_1,T_2:\{0,1\}^*\to\Sigma_\Box$ be trees and 
  $G\subseteq\{0,1\}^*$ some set. Assume that  
  $H\subseteq\{0,1\}^*$ is a finite tree-domain.
  Set
  \begin{align*}
    &P:=\domain(\prune(T_1{\restriction}_H))
    \text{ and}\\
    &D:= P \cup \bigcup_{e\in H^+\cap P^+} \domain(\extract(e,T_2,G)).  
  \end{align*}
  Let $t:=\pretree(H,T_1,T_2,G)$ denote the tree with
  domain $D$ that is defined by
  \begin{align*}
    t(q):=
    \begin{cases}
      T_1(q) & q\in \domain( \prune(T_1{\restriction}_H)),\\
      T_2(q) &\text{otherwise}.
    \end{cases}
  \end{align*}
\end{definition}
\begin{remark}
  $\pretree(H, T_1, T_2, G)$ extracts a tree from $(T_1, T_2, G)$ that
  coincides with $T_1$ on domain $H$ (where positions that are
  labelled by $\Box$
  in $T_1$ count as undefined positions). For each of those branches
  that are  defined up to the border of $H$, we append the
  corresponding subtree of $T_2$. That is, for
  $d\in H^+$ such that $T_1$ is 
  defined on all ancestors of $d$, we append $\extract(d,T_2,G)$. 
\end{remark}

Let us illustrate these definitions in an example.
\begin{example}
  Consider the following infinite trees $T_1$ and $T_2$:\\
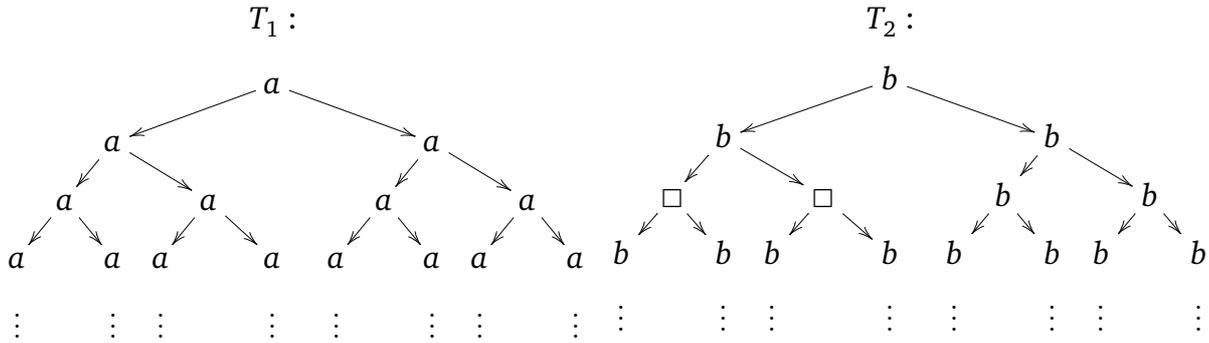
  \begin{xy}
    \xymatrix@R=10pt@C=5pt{
      &&&&&T_1:\\
      &&&&&a\ar[dlll] \ar[drrr]\\
      &&a\ar[drr]\ar[dl]&&&&&&a\ar[dl]\ar[drr]&&\\
      &a\ar[dr]\ar[dl]&&&a\ar[dr]\ar[dl]&&&a\ar[dr]\ar[dl]&&&a\ar[dr]\ar[dl]\\
      a&&a&a&&a&a&&a&a&&a\\
      \vdots&&\vdots&\vdots&&\vdots&\vdots&&\vdots&\vdots&&\vdots
    }
  \end{xy}   
  \begin{xy}
    \xymatrix@R=7pt@C=5pt{
      &&&&&T_2:\\
      &&&&&b\ar[dlll] \ar[drrr]\\
      &&b\ar[drr]\ar[dl]&&&&&&b\ar[dl]\ar[drr]&&\\
      &\Box\ar[dr]\ar[dl]&&&\Box\ar[dr]\ar[dl]&&&
      b\ar[dr]\ar[dl]&&&b\ar[dr]\ar[dl]\\  
      b&&b&b&&b&b&&b&b&&b\\
      \vdots&&\vdots&\vdots&&\vdots&\vdots&&\vdots&\vdots&&\vdots
    }
  \end{xy}
  \vskip 1cm
  Consider $G:=\{w\in\{0,1\}^*: \lvert w
  \rvert$ is odd$\}$ and $H:=\{\varepsilon\}$. Then
  $H^+=\{0,1\}\subseteq G$
  and we obtain the following trees using $\extract$ on $H^+$:\\
  \begin{xy}
    \xymatrix@R=10pt@C=10pt{
      \extract(0, T_2, G):\\
      &b
    }
  \end{xy}  \hskip 5cm 
  \begin{xy}
    \xymatrix@R=10pt@C=10pt{
      \extract(1, T_2, G):\\
      &&b\ar[dl] \ar[dr]&\\
      &b&&b
    }
  \end{xy}\\
  Note that $\prune(T_1{\restriction}_H) = a$. 
  Hence, we conclude that 
  $\pretree(H, T_1, T_2, G)$ is the following tree:
  
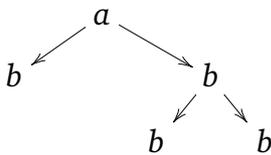
\begin{xy}
    \xymatrix@R=10pt@C=7pt{
      &&a\ar[dll] \ar[drr]\\
      b &&&&b\ar[dl]\ar[dr]\\
      &&&b&&b
    }
  \end{xy}   
\end{example}

Now, we use the function $\pretree$ to define the set of infinite
trees that is encoded by a triple $(T_1, T_2, G)$. 

\begin{definition}
  Let $G\subseteq \{0,1\}^*$. A finite tree-domain $H$ is called
  a \emph{$G$-tree} if $H^+\subseteq G$.
  We set
  \begin{align*}
     \Set(T_1,T_2,G):=\left\{\pretree(H,T_1,T_2,G): H \text{ is a
       }G\text{-tree}\right\}.
  \end{align*}
\end{definition}
\begin{remark} 
  Note that $\Set(T_1, T_2, G)$ may contain infinite trees. Moreover,
  this set may be finite, e.g., if  $G=\emptyset$. 
  In our applications, we 
  always ensure that this definition
  yields an infinite set of finite trees. 
\end{remark}

As the next step, we define the notions of a tree-comb and of the
closure of a tree-comb.
These definitions aim at the following: we look for an infinite
sequence $C=(T_i)_{i\geq 1}$ of finite 
trees that can be encoded by a tuple $(T^C_1, T^C_2, G^C)$ of infinite trees
such that $\Set(T^C_1, T^C_2, G^C)$ contains $C$. 
More precisely, $T_i=\pretree(H,T^C_1,T^C_2,G^C)$ 
for the $G^C$-tree $H$ induced by the $i$-th layer of $G$ in the
following sense. Let $B_1$ be the set of all infinite branches $b$
such that $\lvert G\cap b \rvert \geq i$, i.e., those branches along
which an element of $G$ occurs at least $i$ times. 
Along every infinite branch $b\in B_1$, $H\cap b$ is the
finite branch up to the predecessor of the $i$-th element of $b\cap G$
(if $b\cap G$ contains at least $i$ elements). 
Along every infinite branch $b$ in the complement of $B_1$, $H\cap b$
is maximal in the sense that $H$ contains all predecessors of the
maximal element of $b\cap G$. 

Before we state the precise definition of a tree-comb, let us explain
how this notion generalises the notion of a word-comb. A word-comb is
an infinite 
set of strings such that there is a sequence $g_1 < g_2 < g_3 <
\dots$ of natural numbers such that all but the shortest $n$ words of
the word-comb agree on the first $g_n$ letters. Furthermore, the
length of the words
forming the word-comb grows unbounded. We transfer this principle
to the tree case as follows: we replace the notion of ``length of a
string'' by the notion of ``depth of a 
tree''. Thus, we want a tree-comb to be an infinite sequence of trees
of growing depth such that all but the first $n$ trees coincide on a
certain initial part $D_n$ of their domain. 

Before we state the definition, recall that $\domain(T)^\oplus$
denotes the union of $\domain(T)$ with its border $\domain(T)^+$. 
\begin{definition}\label{def:tree-comb}
   An infinite sequence of finite trees $C=(T_i)_{i\geq 1}$ is called
   a \emph{tree-comb} if 
   $T_j{\restriction}_{\domain(T_i)^\oplus} =
   T_k{\restriction}_{\domain(T_i)^\oplus}$
   for all natural numbers $1\leq i < j<k$.
\end{definition}
\begin{remark}
  A tree-comb $C=(T_i)_{i\geq 1}$ is an
  infinite sequence of finite trees. Abusing notation, we will
  identify $C$ with the set $\{T_i: i\geq 1\}$ if no confusion
  arises. In this sense, we write $D\subseteq C$ for the fact that
  $D$ is an infinite subsequence of $C$. In this case, $D$ is also a
  tree-comb. 
\end{remark}

We will soon see that any infinite set of trees contains a subset
which forms a tree-comb. Before we come to this, let us define the
notion of a representation of a tree-comb by a triple of infinite
trees.

\begin{definition}
  Let $C=(T_i)_{i\geq 1}$ be a tree-comb. We define the 
  trees $T_1^C, T_2^C, G^C$ as follows:
  \begin{align*}
   &T^C_1:\{0,1\}^*\rightarrow \Sigma_\Box \text{ with}\\
   &T^C_1(d):=
   \begin{cases}
     T_{i}(d) &\text{for } d\in \domain(T_i)\cap \bigcup_{j<i} \domain(T_j),\\
     \Box &\text{otherwise,}
   \end{cases}\\
   &T^C_2:\{0,1\}^*\rightarrow \Sigma_\Box \text{ with}\\
   &T^C_2(d):=
   \begin{cases}
     T_{i}(d) &\text{for }d\in \domain(T_i) \setminus \bigcup_{j<i}
     \domain(T_{j}),\\
     \Box &\text{ otherwise,}
   \end{cases}\\
   &G^C:= \{\varepsilon\} \cup \bigcup_{i\geq 1}\left(  \domain(T_i)^+
   \setminus (\domain(T_{i-1})^\oplus)\right).
  \end{align*}
  We call the triple $(T^C_1, T^C_2, G^C)$ the \emph{representation}
  of $C$.  
\end{definition}
\begin{remark}  \label{EncodingTreeCombWellDefined}
  Note that $T^C_1$ is well-defined:
  if there are $i'>i>j$ such that
  \begin{align*}
    d\in\domain(T_{i'})\cap\domain(T_{i})\cap\domain(T_j),    
  \end{align*}
  then
  $T_{i'}(d) = T_i(d)$ by the tree-comb property. 
  
  Furthermore, if there is some node $d$ such that $T_1^C(d)=\Box$
  then $T_1^C(de)=\Box$ for all $e\in\{0,1\}^*$. This is due to the
  fact that $T_1^C(d)=\Box$ if $d$ is in the domain of at most one of
  the $T_i$. But then all descendents of $e$ satisfy this condition,
  too. 

  Note that 
  $G^C = \{\varepsilon\}\cup \bigcup_{i\geq 1} \left(\domain(T_i)^+ \setminus
  \bigcup_{j<i}\domain(T_j)^\oplus\right)$: $(\supseteq)$ is trivially
  true. 
  For $(\subseteq)$ assume that 
  $\varepsilon\neq d\in \domain(T_i)^+ \setminus
  (\domain(T_{i-1})^\oplus)$. Heading for a contradiction, assume that  
  $d\in \domain(T_j)^\oplus$ for some $j<i-1$. By definition of a
  tree-comb,  this implies that $T_i$ and $T_{i-1}$ agree on $d$ which
  contradicts the assumption $d\notin\domain(T_{i-1})^\oplus$ and
  $d\in\domain(T_i)^\oplus$. Thus, $d\notin\domain(T_j)^\oplus$ for
  all $j<i-1$ whence 
  $d\in \bigcup_{i\geq 1} \left (
    \domain(T_i)^+\setminus \bigcup_{j<i} \domain(T_j)^\oplus
  \right)$. 
\end{remark}

\begin{definition}
  Let $C$ be a tree-comb. We call
  $\Closure(C):=\Set(T^C_1, T^C_2, G^C)$ the \emph{closure of $C$}. 
\end{definition}
\begin{remark}
Calling $\Set(T^C_1, T^C_2. G^C)$ a closure of $C$ requires some
justification: we postpone this justification for a while. But 
in Lemma \ref{TreeCombContainedInClosure}, we will see that
$\Closure(C)$ contains each element of $C$.  
\end{remark}

In the following, we study tree-combs,
their representations and their closures. 
First, we show that any infinite set of trees contains a tree-comb. 
Then we show that every tree-comb is contained in its closure. 
Furthermore, we show that the closure of every tree-comb is an
infinite set of finite trees. 
Finally, we introduce a partial order on the closure of every
tree-comb. This order  plays
a crucial technical role in our reduction of the Ramsey quantifier. 
Each Ramsey quantifier that asserts a certain property of all pairwise
distinct $n$-tuples of some infinite set will be reduced to 
the assertion that all pairwise comparable $n$-tuples of the closure
of some tree-comb have this property.  

We apply Ramsey's Theorem  in many of the following proofs. 
Thus, we recall this theorem briefly.

\begin{theorem}[\cite{Ramsey30}] \label{Thm:RamseyTheorem}
  Let $S$ be an infinite set, $C$ a finite set of colours. We write
  $P_n(S)$ for the set of $n$-element subsets of $S$. 
  For each colouring $f:P_n(S)\to C$ of the $n$-element subsets of
  $S$ there is an infinite subset $S'\subseteq S$ such that $f$ is
  constant on $P_n(S')$.
\end{theorem}

We are now prepared to prove that every infinite set of finite trees
contains a subset that induces a tree-comb. 

\begin{lemma}\label{L combs in sets}
  Let $S$ be an infinite set of finite trees. Then there is 
  a tree-comb $C$ such that each element of $C$ is contained in $S$.
\end{lemma}
\begin{proof}
  We define $C=(T_i)_{i\geq 1}$ by induction. 

  Choose $T_1\in S$ arbitrarily. 
  Since $\domain(T_1)^\oplus$ is finite
  and due to Ramsey's Theorem, 
  there is an infinite set $S^1\subseteq S$ such that 
  for all  $T,T'\in S^1$, we have 
  $T{\restriction}_{\domain(T_1)^\oplus} = T'{\restriction}_{\domain(T_1)^\oplus}$.  

  Choose $T_2\in S^1$ arbitrarily. Again, $\domain(T_2)$ is finite
  whence there
  is some infinite $S^2\subseteq S^1$ such that 
  $T{\restriction}_{\domain(T_2)^\oplus} =
  T'{\restriction}_{\domain(T_2)^\oplus}$ for all $T,T'\in S^2$.  

  Continuing this construction, we obtain infinitely many finite trees
  $T_1, T_2, T_3, \dots$ Because of the definition of these trees,
  $C:=(T_i)_{i\geq 1}$ is a tree-comb.
\end{proof}
Note that by definition of the
Ramsey quantifier, the witnesses for Ramsey quantifiers are closed
under taking infinite subsets: if $S$ is an infinite set of
finite trees witnessing some Ramsey quantifier, then every infinite
subset of $S$ also witnesses this Ramsey quantifier.
From this point of view, the previous lemma says that the search space
for
witnesses for Ramsey quantifiers on automatic structures can be
restricted  to tree-combs. 

The next lemma collects some technical facts about the representation
$(T^C_1, T^C_2, G^C)$ of a tree-comb $C$. 

\begin{lemma} \label{T1BasicFacts}
  Let $C=(T_j)_{j\geq 1}$ be a tree-comb and let $i\geq 1$ be some
  natural number.
  \begin{enumerate}
  \item For all $d\in\bigcup_{j<i} \domain(T_j)$ 
    we have $T^C_1(d)\neq \Box$ iff $d\in \domain(T_i)$.
  \item For all $D\subseteq\bigcup_{j<i} \domain(T_j)$ 
    and for $E:=\domain(\prune(T^C_1{\restriction}_D))$, 
    we have
    \begin{align*}
      \prune(T^C_1{\restriction}_D)=T_i{\restriction}_E.      
    \end{align*}
  \item For all $d\in\domain(T_{i-1})^+\setminus \bigcup_{j<i-1}
    \domain(T_j)^\oplus$, we have 
    $\extract(d,T^C_2,G^C)=\inducedTreeof{d}{T_i}$. 
  \end{enumerate}
\end{lemma}
\begin{proof}
  \begin{enumerate}
  \item Let $d\in \bigcup_{j<i} \domain(T_j)$. 
    By definition, $T^C_1(d)\neq \Box$ if and only if there is some
    $k\in \N$ such that 
    $d\in \domain(T_k)\cap \bigcup_{j<k} \domain(T_j)$.
    By assumption on $d$, this is the case if and only if
    there is some $k$ such that 
    $d\in \domain(T_k)\cap \bigcup_{j<k} \domain(T_j) \cap
    \bigcup_{j<i} \domain(T_j)$.  
    We have to show that this is the case if and only if 
    $d\in \domain(T_i)\cap \bigcup_{j<i} \domain(T_j)$.
    
    Assume that $d\in \domain(T_i) \cap \bigcup_{j<j}
    \domain(T_j)$. Setting 
    $k:=i$, we obtain directly that 
    $d\in \domain(T_k)\cap \bigcup_{j<k} \domain(T_j) \cap
    \bigcup_{j<i} \domain(T_j)$. 

    For the other direction, assume that there is some $k\in\N$ such
    that 
    \begin{align*}
      d\in \domain(T_k)\cap \bigcup_{j<k} \domain(T_j) \cap
      \bigcup_{j<i} \domain(T_j).
    \end{align*}
    Due to the definition of a tree-comb, for all $i,k\in\N$ the trees
    $T_i$ and $T_k$ agree on $\bigcup_{j<\min(k,i)}
    \domain(T_j)^\oplus$. It follows immediately that
    \begin{align*}
      \domain(T_k)\cap \bigcup_{j<\min(k,i)} \domain(T_j) =
      \domain(T_i)\cap \bigcup_{j<\min(k,i)} \domain(T_j).
    \end{align*}
    From this we derive directly that 
    \begin{align*}
      &\domain(T_k) \cap \bigcup_{j<k} \domain(T_j)
        \cap \bigcup_{j<i} \domain(T_j) \\
      =
      &\left(\domain(T_k) \cap \bigcup_{j<\min(i,k)}
        \domain(T_j)\right) \cap \bigcup_{j<i} \domain(T_j)\\
      =
      & \left(\domain(T_i) \cap \bigcup_{j<\min(i,k)}
        \domain(T_j)\right) \cap \bigcup_{j<i} \domain(T_j)\\
      \subseteq &\domain(T_i).
    \end{align*}
    Thus, we conclude that $d\in \domain(T_i) \cap \bigcup_{j<i}
    \domain(T_i)$. 
  \item Let $D\subseteq\bigcup_{j<i} \domain(T_j)$. 
    The previous part of this Lemma showed that 
    \begin{align*}
      D\cap\domain(T_i) =
      D\cap \{d: T_1^c(d)\neq \Box\}.      
    \end{align*}
    By definition of the function $\prune$, it follows that
    \begin{align*}
      E:=\domain(\prune(T^C_1{\restriction}_D))= \domain(T_i) \cap D.      
    \end{align*}
    Together with the definition of $T^1_C$, this implies that 
    $T^1_C(d) = T_i(d)$ for all $d\in E$.
    Thus, $\prune(T_1^C{\restriction}_D)= T_i{\restriction}_E$. 
  \item 
    Let $d\in\domain(T_{i-1})^+\setminus \bigcup_{j<i-1}
    (\domain(T_j)^\oplus)$.
    We have to show that 
    \begin{align*}
      \extract(d, T_2^C, G^C) = \inducedTreeof{d}{T_i}.      
    \end{align*}
    There are the following cases.
    \begin{enumerate}
    \item $d\notin\domain(T_k)$ for all $k\in\N$: by definition of
      $T^C_2$ this implies $T^C_2(d)=\Box$ whence
      $\extract(d,T^C_2,G^C)=\inducedTreeof{d}{T_i}=\emptyset$.
    \item Otherwise, 
      there is some $k\in\N$ such that $d\in\domain(T_k)$: in this
      case, $k\geq i$ 
      because $d\notin\bigcup_{j<i-1}\domain(T_j)$. 
      But then $T_k$ and $T_i$ agree on $d$ because
      $d\in\domain(T_{i-1})^\oplus$ and due to the  definition of a
      tree-comb. Hence, $d\in\domain(T_i)$. 
      Furthermore, $T^C_2(e) = T_i(e)$
      for all \mbox{$d\leq e \in\domain(T_i)$} due to the definition
      of $T^C_2$.  
      Moreover,   $e\notin
      \domain(T_j)^+$ for all \mbox{$d<e\in\domain(T_i)$} and all $j\leq i$. 
      Remark
      \ref{EncodingTreeCombWellDefined} then implies that $e\notin
      G^C$ for all $d<e\in\domain(T_i)$. 
      Finally, due to
      $d\in\domain(T_i)\setminus\domain(T_{i-1})$, 
      \begin{align*}
        \domain(T_i)^+\cap\{e: d\leq e\} \subseteq
        \domain(T_i)^+\setminus (\domain(T_{i-1})^\oplus) \subseteq G^C        
      \end{align*}
      Thus, we conclude that
      $\extract(d,T^C_2,G^C)=\inducedTreeof{d}{T_i}$. \qedhere
    \end{enumerate}
  \end{enumerate}
\end{proof}

In the next lemma we show that the representation $(T^C_1,T^C_2,
G^C)$  of a tree-comb $C$ is a correct representation in the following
sense: 
all elements of the tree-comb can be extracted from this
representation, i.e., for each tree $T\in C$,
it holds that $T\in\Closure(C)$.

\begin{lemma} \label{TreeCombContainedInClosure}
  Let $C=(T_i)_{i\geq 1}$ be a tree-comb. For each $i\geq 1$,
  $T_i\in\Closure(C)$, i.e., \mbox{$T_i \in \Set(T^C_1,T^C_2,G^C)$.} 
\end{lemma}
\begin{proof} 
  For each $i\geq 1$, we construct a $G^C$-tree $H$ such that
  $T_i=\pretree(H,T^C_1,T^C_2,G^C)$.
  Set $H:=\bigcup_{j<i} \domain(T_j)$. 
  First, we show that $H$ is a $G^C$-tree, then we show that
  \mbox{$T_i=\pretree(H, T^C_1, T^C_2, G^C)$.}
  \begin{enumerate}
  \item  We have
    to show that $H^+\subseteq G$. 
    Let $x^-$ be the predecessor of some $x\in H^+$
    and let $j$ be minimal such that $x^-\in\domain(T_j)$. By definition $x\in
    \domain(T_j)^+ \setminus (\domain(T_{j-1})^\oplus)$, whence
    $x\in G^C$. 
  \item 
    Let us first consider the restriction of this tree to $H$.
    Set 
    \begin{align*}
      P:= \domain(\prune(T^C_1{\restriction}_H)) =
      \domain(\prune(T^C_1{\restriction}_{\bigcup_{j<i} \domain(T_j)})).
    \end{align*}
    By Lemma  \ref{T1BasicFacts} $T^C_1$ agrees with $T_i$ on $P$. 

    Due to the definition of a $G^C$-tree, for each $d\in H^+$ there
    is some $k<i$ with 
    \mbox{$d\in\domain(T_k)^+\setminus\bigcup_{j<k} \domain(T_j)^\oplus$.} 

    If $k=i-1$, then the third item of Lemma \ref{T1BasicFacts} implies
    $\extract(d,T^C_2,G^C) = \inducedTreeof{d}{T_{i}}$. 

    Otherwise, $k < i-1$. 
    By the definition of a tree-comb, we know that
    \mbox{$d\notin\domain(T_{k+1})$}  iff 
    \mbox{$d\notin\domain(T_{j})$} for all $j>k$.        

    By $d\in H^+$ we know that $d\notin\domain(T_{j})$ for $j<i$. 
    Since $k+1 < i$, $d\notin\domain(T_{k+1})$ whence
    $d\notin\domain(T_j)$ for all $j>k$. 

    Due to $k<i$, we conclude that $d\notin\domain(T_j)$ for
    all $j\geq 1$. 
    Thus, $T^C_2(d) = \Box$ whence 
    $\extract(d,T^C_2,G^C) = \emptyset = \inducedTreeof{d}{T_{i}}$. \qedhere
  \end{enumerate}
\end{proof}

The proof of the previous lemma implies the following corollary.
\begin{corollary} \label{Cor:Min1}
  Let $C$ be a tree-comb and $(T^C_1, T^C_2, G^C)$ be its
  representation. For each $g\in G^C$, there is some $G^C$-tree $H$
  such that $g\in H^+$. 
\end{corollary}
\begin{proof}
  By definition of $G^C$, there is some  $i\in\N$ such
  that $g\in \domain(T_i)^+ \setminus (\domain(T_{i-1})^\oplus)$. By
  Remark \ref{EncodingTreeCombWellDefined}, we know that
  $g\in \domain(T_i)^+ \setminus \bigcup_{j<i} \domain(T_j)^\oplus$. 
  In the proof of Lemma \ref{TreeCombContainedInClosure} we have
  already seen that $H:=\bigcup_{j\leq i} \domain(T_j)$ forms a 
  $G$-tree. The claim follows from $g\in H^+$. 
\end{proof}

The next lemma shows that for the representation $(T^C_1, T^C_2, G^C)$ 
of an arbitrary tree-comb $C$, the set $\Set(T^C_1, T^C_2, G^C)$ is an 
infinite set of finite trees. 
Since we aim at representing infinite sets
of finite trees, we will call any  triple $(T_1, T_2, G)$
\emph{coherent} if it induces
an infinite set of finite trees via the operator $\Set$. 

\begin{lemma}\label{Lemma:CombCoherent}
 Let $C=(T_i)_{i\geq 1}$ be an arbitrary tree-comb. 
 Its closure $\Closure(C)$ is coherent, i.e., 
 $\Set(T^C_1, T^C_2, G^C)$ is an infinite set of finite trees.  
\end{lemma}
\begin{proof}
  By Lemma \ref{TreeCombContainedInClosure}, we have already seen that
  all trees from $C$ are contained in $\Closure(C)$. Hence, 
  $\Closure(C)$ contains an infinite set of finite trees.  
  Thus, it is only left to show
  that each $G^C$-tree 
  $H$ induces a finite tree. 

  Since a $G^C$-tree is a finite tree-domain by definition, it
  suffices to show the 
  finiteness of $\extract(d,T^C_2, G^C)$ for all 
  $d\in G^C$. 

  For this purpose, let $d\in G^C$. Then there is some
  $i\in\N$ such that 
  $d\in T_i^+\setminus \bigcup_{j<i} \domain(T_j)^\oplus$.
  Due to the last item of Lemma \ref{T1BasicFacts}, 
  $\extract(d,T^C_2,G^C)=\inducedTreeof{d}{T_{i+1}}$. Since $T_{i+1}$
  is a finite tree, its subtree rooted at $d$ is also finite. 
\end{proof}

In order to reduce Ramsey quantifiers on automatic structures to
first-order logic on $\omega$-automatic structures, we need to
introduce one further concept concerning tree-combs: for 
$(T^C_1, T^C_2, G^C)$ a representation of some tree-comb $C$, we
define a partial order $<_{G^C}$ on $\Closure(C)$. 
The purpose of this order is the following: the Ramsey quantifier
asserts that there is an infinite set such that its pairwise distinct
$n$-tuples satisfy a certain formula. This assertion will be reduced to
the assertion that there is a closure of some tree-comb such that all
pairwise $<_{G^C}$ comparable $n$-tuples satisfy the formula. 
We are going to define $<_{G^C}$ in such a way that the tree-comb $C$
is ordered linearly. Thus, if there is a tree-comb $C$ such that its
closure 
$\Closure(C)$ witnesses the 
reduced assertion, then $C$ witnesses the original assertion:
with respect to $C$, the notions of ``pairwise distinct'' and ``pairwise
comparable'' coincide whence $C$ witnesses the Ramsey quantifier.  

The order $<_{G^C}$ is defined on trees from $\Closure(C)$ by
comparing the underlying \mbox{$G^C$-trees} with respect to $\subsetneq$. 
We call a $G^C$-tree $H$ the underlying tree for $T\in\Closure(C)$, 
if 
\mbox{$T=\pretree(H, T_1^C, T_2^C, G^C)$} and $H$ is maximal with this
property. Unfortunately, for an arbitrary representation $(T_1, T_2, G)$
this notion is not well-defined. For an extremely degenerated example,
take 
$T_1$ to be the constant $\Box$-labelled tree and $G=\{0,1\}^*$. Any
finite tree domain $H$ forms a $G$-tree and $\pretree(H, T_1, T_2, G)$
is the empty tree for all $H$. Thus, 
there is no maximal $G$-tree underlying the empty tree in this
representation. 
In order to obtain a well-defined notion of  underlying $G$-tree, 
we first define the notion of a
\emph{\SMALL} representation. Afterwards, we show that 
there is an underlying
$G$-tree for every tree $T$ contained in a \SMALL
representation. Furthermore, we prove that the representation of every
tree-comb is \SMALL. Finally, we formally define the order
$<_G$ for each \SMALL representation $(T_1, T_2, G)$. 

\begin{definition}
  We call a representation $(T_1, T_2, G)$ \emph{\SMALL} if the
  following two conditions hold.
  \begin{enumerate}
  \item For all $g\in G$ there is some $G$-tree $H$ such that $g\in
    H^+$.
  \item If there are $d < e\in\{0,1\}^*$ with $d,e\in G$, then for all
    $c<d$ we have $T_1(c)\neq \Box$. 
  \end{enumerate}
\end{definition}
\begin{remark}
  It does not depend on $T_2$ whether $(T_1, T_2, G)$ is \SMALL. Thus,
  we will also say $(T_1, G)$ is \SMALL meaning that $(T_1, T_2, G)$ is
  small.  
  
  Note that the representation of every tree-comb satisfies the first
  condition due to Corollary \ref{Cor:Min1}. 
\end{remark}

\begin{lemma} \label{LemmaCombSmall}
  Let
  $T_1$ and $T_2$ be $\Sigma_{\Box}$-labelled infinite binary trees
  and $G\subseteq\{0,1\}^*$. Assume that $(T_1, T_2, G)$ is \SMALL. 
  For each $T\in\Set(T_1, T_2, G)$ there is a unique maximal $G$-tree
  $H_T$ such that $T=\pretree(H_T, T_1, T_2, G)$. 
\end{lemma}
\begin{proof}
  Fix a $T\in\Set(T_1, T_2, G)$. Let 
  \begin{align*}
    S_T:=\{H\subseteq\{0,1\}^*: H\text{ a }
    G\text{-tree and }T=\pretree(H, T_1, T_2, G)\}.     
  \end{align*}
  Furthermore,  let $H_T:=\bigcup S_T$ be the union of all these
  $G$-trees. 
  
  First, we show that $H_T$
  is a finite tree-domain. By definition of $H_T$ this implies that
  $H_T$ is a $G$-tree. 
  Afterwards, we show that it generates $T$. 

  $H_T$ is infinite if and only if there is an infinite chain $d_0 <
  d_1 < d_2 < \dots \in \{0,1\}^*$ such that for each $i\in\N$ there
  is some $H_i\in S_T$ with $d_i\in H_i^+$, i.e.,  the trees-domains
  in $S_T$ grow unbounded along some infinite branch 
  $d_0 < d_1 < d_2 < \dots < b\in\{0,1\}^\omega$. 

  Heading for a contradiction, assume that such a chain $d_0< d_1< d_2
  <\dots$ exists. 
  
  Since $H_i$ is a $G$-tree, each $d_i\in G$. Because of 
  $d_1\in G$ and $d_1\notin\domain(H_0)^\oplus$, 
  the definition of
  $\pretree(H_0, T_1, T_2, G)$ implies that $d_1\notin
  \domain(\pretree(H_0, T_1, T_2, G))$. Due to 
  \mbox{$\domain(T)=\pretree(H_0, T_1, T_2, G))$}, we conclude that
  $d_1\notin\domain(T)$. 
  
  On the other hand, $(T_1, T_2, G)$ is \SMALL whence  $d_1<d_2 < d_3 \in
  G$ implies that  
  $T_1(c)\neq \Box$ for all $c\leq d_1$. Since $d_1\in H_2$, 
  $\pretree(H_2, T_1, T_2, G)$ and $T_1$ coincide up to $d_1$. 
  We conclude that 
  $d_1\in \domain(T)=\domain(\pretree(H_2, T_1,  T_2, G)$. 

  This
  contradicts $d_1\notin\domain(T)$. Thus, the tree-domains in $S_T$
  cannot grow 
  unbounded along any infinite branch and we conclude that $H_T$ is a
  well-defined finite tree-domain.

  We  come to our second claim: $H_T\in S_T$, or equivalently
  $T=\pretree(H_T, T_1, T_2, G)$.  
  In order to prove this claim, let $b\in\{0,1,\}^\omega$ be an
  arbitrary infinite branch. 
  There is a unique element $d$ that is on the border of $H_T$ in the
  branch $b$, i.e., there is a unique element $d\in b\cap H_T^+$. Let
  $d^-$ be the 
  direct predecessor of $d$. 
  
  By definition of $H_T$, $d^-\in H$ for some $H\in S_T$ and there is
  no $H'\in S_T$ with $d\in H'$. Thus, $H_T$ and $H$ agree
  along the branch $b$ whence $T=\pretree(H, T_1, T_2, G)$ and
  $\pretree(H_T, T_1, T_2, G)$ coincide along $b$. 

  For each infinite branch there is such a $H\in S_T$ whence we
  conclude that 
  $T$ and $\pretree(H_T, T_1, T_2, G)$ coincide along each infinite
  branch. Hence, \mbox{$T=\pretree(H_T, T_1, T_2, G)$.}
\end{proof}

\begin{lemma}
  Let $C=(T_i)_{i\geq 1}$ be some tree-comb and $(T_1^C, T^C_2, G^C)$ its
  representation. Then $(T_1^C, T^C_2, G^C)$ is \SMALL. 
\end{lemma}
\begin{proof}
  We have to show the following two claims:
  \begin{enumerate}
  \item 
    For all $g\in G^C$ there is some $G^C$-tree $H$ such that $g\in
    H^+$.
  \item If there are $d < e\in\{0,1\}^*$ with $d,e\in G^C$, then for all
    $c<d$ we have $T^C_1(c)\neq \Box$. 
  \end{enumerate}
  The first claim holds due to Corollary \ref{Cor:Min1}. 
  The second part is an easy consequence of the definition of $T_1^C$:
  by definition of $G^C$, $d,e\in G^C$ with $d\leq e$ implies that  
  there are numbers $i$ and $j$ such that $i\neq j$,
  $d\in\domain(T_i)^+$ and $e\in\domain(T_j)^+$. But this implies that 
  for all $c<d$, $c\in\domain(T_i)\cap\domain(T_j)$ whence by
  definition of $T^C_1$, $T^C_1(c)= T_k(c)$ for $k=\max(i,j)$. 
\end{proof}

We conclude the section on tree-combs by defining the order $<_G$
for all \SMALL representations $(T_1, T_2, G)$ and by showing that
each tree-comb $C$ is linearly ordered by the induced order
$<_{G^C}$. 

\begin{definition}
  Let $G\subseteq\{0,1\}^*$. Furthermore, let $H$ and $H'$ be
  $G$-trees. We define $H <_G H'$ if the following
  two conditions hold:
  \begin{enumerate}
  \item $H\subsetneq H'$ and
  \item for each infinite branch $b\in\{0,1\}^\omega$,  
    $H\cap b = H'\cap b$ implies $(b\setminus (H^\oplus)) \cap G = \emptyset$. 
  \end{enumerate}
\end{definition}
This means that $H<_G H'$ holds if $H'$ extends $H$ properly along
each branch where this is possible for a $G$-tree. In other words, 
if there is a descendent of some $d\in H^+$ which is in $G$, then
$H'$ must contain $d$. Thus, $H'$ extends $H$ properly along this
branch.  

We extend this order to $S,T\in\Set(T_1,T_2,G)$ for \SMALL  representations
$(T_1, T_2, G)$ as follows.
\begin{definition}
  Let $(T_1, T_2, G)$ be a \SMALL representation. 
  Let $H_S$ ($H_T$) denote the maximal $G$-tree such that
  $S=\pretree(H_S,T_1,T_2,G)$ ($T=\pretree(H_T,T_1,T_2,G)$,
  respectively), i.e., $H_S$ and $H_T$ are the underlying trees for $S$
  and $T$, respectively. 
  We set 
  \begin{align*}
    S<_G T\text{ iff } H_S <_G H_T.    
  \end{align*}
\end{definition}

This order formalises the idea that the underlying $G$-tree $H'$
extends $H$ in each possible direction. Since a $G$-tree ends along each
path just in 
front of a node from $G$, the branches where a $G$-tree cannot be
extended are those where no further elements from $G$ follow after
$H^\oplus$. 

We conclude this section by showing that any tree-comb $C$ is linearly
ordered by the induced order $<_{G^C}$. 

\begin{lemma} \label{Lemma:TreeCombLinear}
  Let  $C=(T_i)_{i\geq 1}$ be a tree-comb. 
  Then $T_i <_{G^C} T_{k}$ for all $1\leq i \leq k$. 
\end{lemma}
\begin{proof}
  Let $H_i:=\bigcup_{j<i} \domain(T_j)$
  and let $\hat H_i$ be the maximal $G^C$-tree generating $T_i$. 
  From the proof of lemma \ref{TreeCombContainedInClosure} we know
  that $H_i\subseteq \hat H_i$ because $H_i$ also generates $T_i$. 
  
  By definition of $T_i<_{G^C} T_{i+1}$, it suffices to show that
  $\hat H_i  <_{G^C} \hat H_{i+1}$. We prove this 
  claim in two steps.
  First we show that $\hat H_i \subseteq H_{i+1}$. This
  implies $\hat H_i \subseteq \hat H_{i+1}$ and furthermore, these
  two trees cannot coincide because they generate two different trees,
  namely, $T_i$ and $T_{i+1}$.
  Afterwards,  we show that for each infinite branch $b$ the following
  holds. If \mbox{$\hat H_i \cap b = \hat H_{i+1} \cap b$} then
  $(b\setminus (\hat H_i^\oplus)) \cap G^C = \emptyset$. 
  \begin{enumerate}
  \item Since $\hat H_i$ and $H_{i+1}$ are tree-domains, $\hat H_i
    \not\subseteq H_{i+1}$ would imply that $\hat H_i \cap
    (H_{i+1}^+)\neq\emptyset$. 
    Heading for a contradiction, assume that there is some $d\in
    \hat H_i \cap  (H_{i+1}^+)$. 

    By definition of
    $H_{i+1}$, $d\in\left(\bigcup_{j\leq i} \domain(T_i)\right)^+$. 
    Let $D:=\domain(\prune(T^C_1{\restriction}_{\hat H_i}))$. By
    definition of $\hat H_i$, 
    $\prune(T^C_1{\restriction}_{\hat H_i}) = T_i{\restriction}_D$. Since
    $d\notin\domain(T_i)$, this implies $T^C_1(d)=\Box$. 
    Thus, by definition of $T^C_1$ it is not possible that there are
    two numbers $j_1\neq j_2 \in\N$ such that
    $d\in\domain(T_{j_1})\cap\domain(T_{j_2})$. 

    We claim that then 
    $d\notin\domain(T_j)$ for all $j\geq 1$.

    For $j\leq i$, $d\notin\domain(T_j)$ due to 
    $d\in \left(\bigcup_{j\leq i} \domain(T_i)\right)^+$.
    
    Nevertheless, for the same reason, 
    $d\in\domain(T_k)^\oplus$ for some $k\leq i$. 
    Due to the tree-comb property, this implies that 
    $T_{j_1}$ and $T_{j_2}$ agree at $d$ for all 
    $j_1 > j_2 > k$. 
    Since we have already seen that there cannot be two different trees
    $T_{j_1}$ and $T_{j_2}$ defined at $d$, we conclude that there is
    no $j>k$ such that $d\in\domain(T_j)$. 
    
    Thus, we conclude that $d\notin\domain(T_j)$ for all $j\geq 1$. 
    This implies that    
    all $d < e$ satisfy $e\notin\domain(T_j)^+$ for all $j\geq
    1$. Due to the definition of $G^C$, it follows that $e\notin G^C$
    for all $d<e$. 
    Thus, $d$ cannot be contained in any $G^C$-tree. 

    But this contradicts the assumption that $d\in\hat H_i$. 
    
    We conclude that $\hat H_i \cap (H_{i+1}^+)=\emptyset$ which
    implies $\hat H_i \subseteq H_{i+1} \subseteq \hat
    H_{i+1}$.
  \item Fix some infinite branch $b$ such that 
      $\hat H_i \cap b = \hat H_{i+1} \cap b $. Due to
      $\hat H_i \subseteq H_{i+1} \subseteq \hat H_{i+1}$ this implies 
    \begin{align}\label{HatHi}
      \hat H_i \cap b = \hat H_{i+1} \cap b =
      H_{i+1} \cap b = \bigcup_{j\leq i} \domain(T_j) \cap b.        
    \end{align}
    As a direct consequence of the coincidence of $\hat H_i$ and $\hat
    H_{i+1}$ along $b$, we obtain that
    \begin{align*}
      T_i{\restriction}_b = 
      \pretree(\hat H_i,T^C_1,T^C_2,G^C){\restriction}_b=
      \pretree(\hat H_{i+1},T^C_1,T^C_2,G^C){\restriction}_b=
      T_{i+1}{\restriction}_b.         
    \end{align*}
    Thus,
    \begin{align}\label{TIMOreEqual}
      b\cap \bigcup_{j\leq i} \domain(T_j) = b\cap \bigcup_{j\leq
        i+1} \domain(T_j).
    \end{align}
    Now, let $d$ be the unique element of $b\cap (\hat H_{i+1}^+)$. 
    
    \ref{HatHi} implies that there is some $k\leq i$ such that 
    $d\in \domain(T_k)^\oplus$ while $d\notin\domain(T_j)$ for all $j\leq i$.
    Due to \ref{TIMOreEqual}, this implies $d\notin\domain(T_{i+1})$. 
    Since $i+1> k$,  it follows from the tree-comb property that
    $d\notin\domain(T_j)$ for all $j\geq i+1 > k$. 

    We conclude that $d\notin\domain(T_j)$ for all $j\geq 1$. 
    But this implies that no proper descendant of $d$ is contained in
    $\domain(T_j)^+$ for any $j\geq 1$. Hence, no proper descendant of $d$ 
    is contained in $G^C$. Since $d\in b\cap(\hat H_{i+1}^+)$, it follows
    that $(b\setminus (\hat H_{i+1}^\oplus)) \cap G =
    \emptyset$, which concludes the proof.\qedhere
  \end{enumerate}
\end{proof}

\subsection{Reduction of the Ramsey Quantifier}

We now reduce $\FO{}((\RamQ{n})_{n\in\N})$ on an automatic structure
$\mathfrak{A}$ to $\FO{}$ on an $\omega$-automatic structure 
$\OmegaExp(\mathfrak{A})$.

\paragraph{Adding Tree-Comb Representations to an Automatic Structure}
From now up to the end of Section \ref{EndOfRamsey}, we fix an
automatic structure $\mathfrak{A}$. We 
assume that, without loss of generality,
the identity $\Id$ is a tree presentation of $\mathfrak{A}$. This means that
the universe of $\mathfrak{A}$ is a regular
subset $A\subseteq \Trees{\Sigma}$ and all relations of $\mathfrak{A}$
are automatic.

We next define a structure $\OmegaExp(\mathfrak{A})$ corresponding to
$\mathfrak{A}$ in the following sense. 
$\OmegaExp(\mathfrak{A})$ is the disjoint union of $\mathfrak{A}$ with
a structure that allows to reason about tree-combs in the following sense:
each
$\FO{}((\RamQ{n})_{n\in\N})$ formula 
over $\mathfrak{A}$ can be reduced to an $\FO{}$ formula over
$\OmegaExp(\mathfrak{A})$.  
Furthermore, $\OmegaExp(\mathfrak{A})$ turns out to be
$\omega$-automatic whence 
this reduction proves the decidability of $\FO{}((\RamQ{n})_{n\in\N})$
over $\mathfrak{A}$. Later, we use the reduction of
an $\FO{}((\RamQ{n})_{n\in\N})$ formula $\varphi$ on $\mathfrak{A}$ to
an $\FO{}$ formula on
$\OmegaExp(\mathfrak{A})$ in order to design 
an $\omega$-automaton that represents $\varphi$ on
$\OmegaExp(\mathfrak{A})$ and that can be turned into 
an automaton
$\mathcal{A}_\varphi$ that corresponds to $\varphi$ on
$\mathfrak{A}$.

\begin{definition}
  Let $\OmegaExp(\mathfrak{A})$ be the following structure. 
  \begin{itemize}
  \item The universe is
    $A':=\Trees{\Sigma}\cup\fullTrees{\Sigma_\Box}\cup\fullTrees{\{0,1\}}$
    where
    \begin{align*}
      \fullTrees{\Sigma}:=\{T\in\infTrees{\Sigma}:\domain(T)=\{0,1\}^*\}    
    \end{align*}
    is the set of all full
    infinite binary $\Sigma$-trees.
    We identify a subset $G\subseteq \{0,1\}^*$ with its characteristic
    map in $\fullTrees{\{0,1\}}$. 
  \item The basic relations are those of $\mathfrak{A}$ including the unary
    relation~$A$ which denotes the universe of the structure~$\mathfrak{A}$.
  \item We add the following new relations:
    \begin{enumerate}
    \item $\Trees{\Sigma}$, $\fullTrees{\Sigma_\Box}$, and
      $\fullTrees{\{0,1\}}$,
    \item $\In:=\left\{(T,T_1,T_2,G)\in\Trees{\Sigma} \times
        (\fullTrees{\Sigma_\Box})^2  
        \times \fullTrees{\{0,1\}}: T\in\Set(T_1,T_2,G)\right\}$,
    \item $\mathrm{Coherent}:=\left\{(T_1,T_2,G)\in (\fullTrees{\Sigma_\Box})^2  
        \times \fullTrees{\{0,1\}}: \Set(T_1,T_2,G)\text{ is
          coherent}\right\}$
          (recall that coherent means that $\Set(T_1, T_2, G)$ is an
          infinite set of finite trees),
    \item $\SMALLRel:=\left\{ (T_1, G)\in
        \fullTrees{\Sigma_\Box}\times \fullTrees{\{0,1\}}:
        (T_1, G) \text{ is \SMALL}\right\}$,
    \item $\mathrm{Comp}:=$
      \begin{align*}
        \left\{(S,T,T_1,T_2,G)\in \mathcal{T}:
          S,T\in\Set(T_1,T_2,G)\text{ and either }  
          S <_G T\text{ or }T <_G S\right\}
      \end{align*}
      for
      $\mathcal{T}:=(\Trees{\Sigma})^2
      \times (\fullTrees{\Sigma_\Box})^2
      \times \fullTrees{\{0,1\}}$.
    \end{enumerate}
  \end{itemize}
\end{definition}
Now, we  construct an $\omega$-presentation 
of $\OmegaExp(\mathfrak{A})$ over the 
alphabet \mbox{$\Gamma=\Sigma\cup\{\bot,\Box,0,1\}$.}
Recall that we write $T^\bot$ for the lifting of a tree $T$ to the
full domain $\{0,1\}^*$ where we use $\bot$ as a padding symbol.
We define the domain of the presentation to be the set 
\begin{align*}
  L:= \{T^\bot: T\in\Trees{\Sigma}\}
  \cup \fullTrees{\Sigma_\Box} \cup \fullTrees{\{0,1\}}.
\end{align*}
This set is obviously $\omega$-regular. 
Furthermore, it is easy to describe a bijection
$h:L\to A'$ by stating its inverse
$h^{-1}:A'\to L$. For $T\in\Trees{\Sigma}$, we set
$h^{-1}(T) := T^\bot$, for 
all other elements $T$ of $A'$, we set
$h^{-1}(T) := T$.  It 
remains to show that the ($h$-preimages of the) relations of 
$\OmegaExp(\mathfrak{A})$ are 
$\omega$-automatic. This is trivial for
$h^{-1}(\Trees{\Sigma})=\{T^\bot: T\in\Trees{\Sigma}\}$,
$h^{-1}(\fullTrees{\Sigma_\Box})=\fullTrees{\Sigma_\Box}$ and
$h^{-1}(\fullTrees{\{0,1\}})=\fullTrees{\{0,1\}}$. For the
relation $A$ (the universe of $\mathfrak{A}$), we have
$h^{-1}(A)=\{T^\bot: T\in A\}$. Recall that 
$\mathrm{id}: A\to A$ is
a tree presentation for $\mathfrak{A}$, so $A\subseteq\Trees{\Sigma}$
can be accepted by an automaton. This automaton can be
transformed into an $\omega$-automaton accepting~$h^{-1}(A)$
(cf. Lemma \ref{Lemma:FinToOmegaAutomaton}). A
similar argument applies to the basic relations of
$\mathfrak{A}$. Thus, it remains to consider the relations $\In$,
$\mathrm{Coherent}$, $\SMALLRel$ and $\mathrm{Comp}$.

\begin{lemma}
  The relations $h^{-1}(\In)$, $h^{-1}(\mathrm{Coherent})$,
  $h^{-1}(\SMALLRel)$, and
  $h^{-1}(\mathrm{Comp})$ are $\omega$-automatic.
\end{lemma}

\begin{proof}
  \begin{enumerate}
  \item $h^{-1}(\In)$: The property ``$H$ is a $G$-tree'' is an
    MSO-definable property of the infinite tree $H \otimes G$ (where
    we consider $H$ 
    and $G$ as characteristic maps). Similarly,
    ``$T=\pretree(H,T_1,T_2,G)$'' is an MSO-definable property of the
    infinite tree $\bigotimes(T^\bot,T_1,T_2,H, G)$. Thus,  also
    ``$T\in\Set(T_1,T_2,G)$'' is an MSO-definable property of the infinite tree
    $\bigotimes(T^\bot,T_1,T_2,G)$. Hence, $\omega$-automaticity of
    $h^{-1}(\In)$ follows from Theorem \ref{Rab69AutomaticEqualsMSO}.
  \item $h^{-1}(\mathrm{Coherent})$: the condition ``for any $G$-tree
    $H$, $\pretree(H,T_1,T_2,G)$ is 
    actually a finite tree'' is an MSO-definable property of the tree
    $\bigotimes(T_1,T_2,G)$ and can therefore be checked by an
    $\omega$-automaton by
    Theorem \ref{Rab69AutomaticEqualsMSO}. 
    Now assume that 
    $\Set(T_1,T_2,G)$ is a set of finite trees. It is infinite if and
    only if the union of the domains of its elements is infinite,
    i.e., if this union contains an infinite branch. 
    But the property ``there is an infinite branch $b$ such that for
    each element $d$ in $b$ there is a $G$-tree $H_d$ such that 
    $d\in\pretree(H_d,T_1,T_2,G)$'' is an \MSO-definable property
    of the tree $\bigotimes(T_1,T_2,G)$. 
  \item $h^{-1}(\SMALLRel)$: the property ``for each
    $d\in G$, there is  $G$-tree $H$ with 
    $d\in H^+$'' is an \MSO-definable property.
    Furthermore, ``for all $c<d<e$ with $e\in G$ and $d\in G$, it
    holds that
    $T_1(c)\neq\Box$'' is first-order definable on $(G\otimes T_1, <)$
    and the prefix order $<$ on $G\cap\{0,1\}^*$ is \MSO definable on $G$. 
  \item $h^{-1}(\mathrm{Comp})$: 
    the assertion ``$H <_G H'$'' is an
    MSO-definable property of the infinite 
    tree $\bigotimes(H, H', G)$: there is a formula that
    checks for each branch $b$ 
    that either 
    \begin{align*}
      H\cap b \subsetneq H' \cap b\text{ or }H\cap b = H' \cap b      
    \end{align*}
    and there is 
    no $d\in (b\cap G) \setminus (H^\oplus)$ with $G(d)=1$. Furthermore,
    the maximal $G$-trees generating $S$ and $T$ are MSO definable in 
    $\bigotimes(S^\bot,T^\bot,T_1,T_2,G)$: we have already seen that
    $T=\pretree(H,T_1,T_2,G)$ is \MSO definable. 
    Thus, the set of $G$-trees generating $T$ (or $S$) are
    definable. The maximal of these trees is the $G$-tree underlying
    $T$ (or $S$). But maximality of a tree among an \MSO-definable set
    of trees is clearly \MSO definable. \qedhere
  \end{enumerate}
\end{proof}

Thus, summarising these results, we obtain the following corollary.

\begin{corollary}
  For each automatic structure $\mathfrak{A}$, the corresponding
  structure $\OmegaExp(\mathfrak{A})$ is $\omega$-automatic.
\end{corollary}

\paragraph{Reduction of the Ramsey quantifier}
We now inductively translate an $\FO{}((\RamQ{n})_{n\in\N})$ formula in
the language of 
$\mathfrak{A}$ into an $\FO{}$ formula in the language
of $\OmegaExp(\mathfrak{A})$. The idea is to replace the
occurrence of a Ramsey quantifier like $\RamQ{n} \bar x (\varphi)$ by
the assertion that there is a \SMALL and
coherent representation $(T_1, T_2, G)$ such that all pairwise
$<_G$-comparable $n$-tuples from $\Set(T_1, T_2, G)$ satisfy
$\varphi$. Our intention is to consider $(T_1, T_2, G)$ as the
representation of some tree-comb. We will first define this reduction
in detail. Then we prove its soundness. Finally, we show that this
reduction is correct.  

\begin{definition}
  For each $\FO{}(\exists^\mathrm{mod},(\RamQ{n})_{n\in\N})$ formula
  $\varphi$ in the 
  language of $\mathfrak{A}$, we define its reduction $\reduction(\varphi)$
  to the  $\FO{}(\exists^\mathrm{mod})$ language of $\OmegaExp(\mathfrak{A})$ by
  \begin{align*}
    \reduction(\varphi) :=& \varphi \text{ for $\varphi$ an atomic formula},\\
    \reduction(\varphi\lor\psi) :=& \reduction(\varphi)\lor\reduction(\psi),\\
    \reduction(\neg\varphi) :=& \neg\reduction(\varphi),\\
    \reduction(\exists x \varphi) :=& \exists x   (x \in A \land
    \reduction(\varphi)),\\
    \reduction(\exists^{k,l} x \varphi) :=& \exists^{k,l}  x (x \in A \land
    \reduction(\varphi)),\\ 
    \reduction(\RamQ{n} \bar x(\varphi)) :=&
    \exists T_1,T_2\in \fullTrees{\Gamma_\Box},G\in\fullTrees{\{0,1\}}
    \psi_{\mathrm{CoSm}}(T_1,T_2,G) \land  \psi_{\mathrm{Ram}}(T_1,T_2,G),
  \end{align*}
  where
  \begin{align*}
    \psi_{\mathrm{CoSm}}(T_1, T_2, G) :=&
    \mathrm{Coherent}(T_1,T_2,G)  \land  \SMALLRel(T_1, G)
    \land \forall x \left(\In(x,T_1,T_2,G)\to x\in A\right) 
  \end{align*}
  and 
  \begin{align*}
    \psi_{\mathrm{Ram}} &:=
    \forall x_1,\dots,x_n\in \Trees{\Sigma}\\
    &\left(\left(
    \bigwedge_{1\leq i \leq n}\In(x_i,T_1,T_2,G)
    \land  \bigwedge_{1 \leq i< j \leq
      n}\mathrm{Comp}(x_i,x_j,T_1,T_2,G) \right)
    \to \reduction(\varphi)\right).
  \end{align*}
\end{definition}
\begin{remark}
  This reduction of a Ramsey quantifier asserts that there is a
  representation of an infinite set such
  that each pairwise comparable $n$-tuple from this set satisfies
  $\varphi$. At first, this seems to be a weaker condition than the
  assertion of the Ramsey quantifier because there are
  tuples of pairwise distinct elements that are not tuples of
  pairwise $<_G$ comparable elements.  
  But it turns out that this condition is
  sufficient: there is an infinite linearly $<_G$-ordered subset $S'$ for each
  $\Set(T_1, T_2, G)$ where $(T_1, T_2, G)$ is a
  coherent and \SMALL representation.
\end{remark}

In the following we first show 
that this translation is sound, i.e.,
for any formula $\varphi$, if
$\OmegaExp(\mathfrak{A})\models\reduction(\varphi)$, then 
$\mathfrak{A}\models\varphi$. Afterwards, we  prove the correctness,
i.e., for any formula $\varphi$, if $\mathfrak{A}\models\varphi$, then 
$\OmegaExp(\mathfrak{A})\models\reduction(\varphi)$.

\subsection{Soundness of the Reduction}

In order to prove the soundness of our reduction, we start with a
technical lemma. It asserts that for $(T_1, T_2, G)$ some
representation of an infinite set of finite trees, there is at least
one branch with infinitely many nodes in $G$. We use this fact in
order to prove the existence of an infinite linear $<_G$-ordered
subset of every small representation. 

\begin{lemma}
  Let $T_1,T_2:\{0,1\}^*\to\Sigma_\Box$ and $G\subseteq \{0,1\}^*$
  such that $(T_1, T_2, G)$ is coherent, i.e., $S=\Set(T_1,T_2,G)$ is
  an infinite set of finite $\Sigma$-trees.  
  If $b$ is an infinite branch in $\bigcup_{T\in S} \domain(T)$, then
  $\lvert b\cap G \rvert=\infty$. 
\end{lemma}
\begin{proof}
  Let $b$ be an infinite branch in $\bigcup_{T\in S}\domain(T)$. 
  Assume that 
  $\lvert b\cap G \rvert < \infty$. Then there are 
  $d_1 < d_2 < \dots <  d_n\in\{0,1\}^*$ such that 
  $G\cap b =\{d_1, d_2, \dots, d_n\}$.  
  Under this assumption, $b\subseteq\bigcup_{T\in S} \domain(T)$
  implies that  there is some $i\leq n$  and
  some 
  $G$-tree $H$ with $d_i\in H^+$ such that 
  $\domain(\extract(d_i,T_2,G))\cap B$ is infinite. But this implies that
  $S$ contains an infinite tree which is a contradiction to the
  assumption that $S\subseteq\Trees{\Sigma}$. 
\end{proof}

\begin{lemma}\label{L-linear}
  Let $T_1,T_2:\{0,1\}^*\to\Sigma_\Box$ and $G \subseteq \{0,1\}^*$
  such that $S:=\Set(T_1,T_2,G)$ is coherent and $(T_1, G)$ is
  \SMALL. Then there is an  infinite subset
  $S'\subseteq S$ which is  linearly ordered by $<_G$. 
\end{lemma}
\begin{proof}
  We show that for every tree $T\in S$
  there is a tree $T'\in S$  
  with $T <_G T'$. Let $T\in S$ and 
  $H$ be the maximal $G$-tree such that $T=\pretree(H,T_1,T_2,G)$. 

  Let 
  \begin{align*}
    D:=\{d\in H^+ : \text{ there is an infinite branch } b \text{
      such that }d<b \text{ and }(b\setminus (H^\oplus)) \cap G \neq
    \emptyset\}.     
  \end{align*}

  Since $\lvert S \rvert = \infty$, there is an infinite branch 
  in $\bigcup_{T\in S} \domain(T)$. Together with the previous lemma,
  this implies that $D$ is nonempty. 
  
  Since $(T_1, G)$ is \SMALL, for each $d\in D$ there exists a 
  $G$-tree $H_d$ with $d\in H_d$. Set 
  $H':=H\cup \bigcup_{d\in D} H_d$. We claim that $H'$ is a $G$-tree
  with $H <_G H'$. 

  Since $D$ is finite, $H'$ is a finite tree. Furthermore for each $e\in
  {H'}^+$ either $e\in H^+$ or $e\in H_d^+$ for some $d\in D$. Thus,
  $H'$ is a $G$-tree. 
  We claim that $H <_G H'$. It  is clear that 
  $d\in H'\setminus H$ for all $d\in D\neq\emptyset$ and that
  $H\subseteq H'$ whence $H\subsetneq H'$. 
  Now assume that $b$ is an infinite branch such that 
  $H\cap b = H'\cap b$. We have to show that $(b\setminus (H'^\oplus))
  \cap G =\emptyset$. 

  Heading for a contradiction assume 
  that  there is some element $e$ contained in this set. 
  Let $d$ be the unique element in $b\cap H^+$.  
  By definition of $D$, $d\in D$. Thus, $d\in H'\setminus H$ which 
  contradicts $H\cap b = H'\cap b$.  

  Hence, for $T':=\pretree(H',T_1,T_2,G)$ we have $T<_G T'$. Repeating
  this construction ad infinitum we obtain an infinite, linearly
  $<_G$-ordered subset of $S$.
\end{proof}

\begin{lemma}
  Let $\varphi\in\FO{}(\exists^\mathrm{mod}, (\RamQ{n})_{n\in\N})$ be
  a sentence. If   
  $\OmegaExp(\mathfrak{A}) \models\reduction(\varphi)$, then
  $\mathfrak{A}\models\varphi$. 
\end{lemma}

\begin{proof}
  Since we want to prove the proposition by induction on the
  construction of $\varphi$, we also have to consider 
  formulas with free variables. Hence, the statement we actually prove
  is the following: 
  \begin{claim}
    Let $\varphi\in\FO{}(\exists^{\mathrm{mod}}, (\RamQ{n})_{n\in\N})$
    be a formula with free 
    variables among 
    $x_1,\dots,x_n$ and let $a_1,a_2, \dots,a_n\in A$. If
    $\OmegaExp(\mathfrak{A}),(a_1,a_2,
    \dots,a_n)\models\reduction(\varphi)$, then  
    $\mathfrak{A},(a_1,a_2, \dots,a_n)\models\varphi$.
  \end{claim}

  The inductive proof of this claim is rather clear except for the
  case $\varphi=\RamQ{n}\bar x(\psi)$. So let
  $a_1, a_2, \dots, a_m\in\mathfrak{A}$, and let
  $\OmegaExp(\mathfrak{A}), (a_1, a_2, \dots, a_m)
  \models\reduction(\varphi)(y_1, y_2, \dots, y_m)$. 
  Then there are infinite trees 
  $T_1,T_2$ and $G$ with the properties given by $\reduction(\varphi)$. In
  particular, 
  $S=\Set(T_1,T_2,G)\subseteq A$ is an infinite set of finite trees and
  $(T_1, G)$
  is \SMALL. By
  Lemma \ref{L-linear}, there is an infinite $S'\subseteq S$
  that is linearly ordered by $<_G$. Hence, from the
  properties of $T_1,T_2,G$, we obtain that
  \begin{align*}
    \OmegaExp(\mathfrak{A}),(t_1,\dots,t_n, a_1, a_2, \dots, a_m,)
    \models\reduction(\psi)(x_1, x_2, \dots, x_n, y_1, y_2, \dots,
    y_m)
  \end{align*}
  for every tuple  
  $(t_1,\dots,t_n)\in (S')^n$ such that the $t_i$ are pairwise
  $<_G$-comparable. Since the pairwise $<_G$-comparable tuples are
  exactly the  pairwise distinct tuples in $S'$, 
  $S'\subseteq A$ witnesses 
  $\mathfrak{A}, (a_1, a_2, \dots, a_m) 
  \models\RamQ{n}\bar x(\psi)(y_1, y_2, \dots, y_m)$. 
\end{proof}

\subsection{Correctness of the Reduction}

The outline of the correctness proof is as follows. 
We fix some formula $\varphi:=\RamQ{n}\bar x(\psi)$. 
We have already seen that every witness for the Ramsey quantifier in
$\varphi$ contains a subset that forms a tree-comb. 
In the following we show that this tree-comb 
contains a certain tree-comb $D\subseteq C$
such that the trees  $T^D_1, T^D_2$, and $G^D$ witness the reduction
$\reduction(\varphi)$.
This means that the  pairwise $<_G$-comparable tuples from the closure
$\Closure(D)$ witness $\reduction(\psi)$. 
In order to prove this, we introduce the notion
of \emph{homogeneity} of some tree-comb with respect to an automaton
$\mathcal{A}$. We will show that any tree-comb witnessing a
Ramsey quantifier contains a homogeneous tree-comb. Furthermore, an
automaton accepts all pairwise distinct $n$-tuples from a homogeneous
tree-comb if and only if it accepts all pairwise comparable $n$-tuples
from the closure of this tree-comb. 
This completes the proof because the existence of such a set is
exactly what the reduction of $\varphi$ asserts. 

Before we give formal definitions, let us informally explain what the
concept of homogeneity is. 
Consider a formula $\RamQ{n} \bar x (\psi)$ with
$\psi\in\FO{}$. 
Assume that there is some
tree-comb $C=(T_i)_{i\geq 1} \subseteq \mathfrak{A}$ witnessing this
quantifier on $\mathfrak{A}$. Let $\mathcal{A}$ denote the automaton
corresponding to $\psi$, i.e., 
$\mathfrak{A}, \bar a \models \psi$ if and only if $\mathcal{A}$
accepts $\bigotimes \bar a$.   
Thus, any pairwise distinct $n$-tuple from $C$ is accepted by
$\mathcal{A}$. Recall that for any finite tree-domain $D$ most of the
elements from $C$ agree on $D$. More precisely, if $D=\domain(T_i)$
then $T_{i+1}, T_{i+2}, T_{i+3} ,\dots$ agree on $D$. We call the
tree-comb homogeneous with respect to $\mathcal{A}$, if all $n$-tuples
from the closure of the tree-comb that agree on some finite domain $D$ are
accepted by runs that coincide on $D$. 

The purpose of this concept is the following: 
First of all, note that every tree $T$ from the closure $\Closure(C)$ 
locally coincides with a tree from $C$ in the following sense. 
Let \mbox{$D\subseteq
\domain(T_i)^\oplus\setminus\domain(T_{i-1})$} be some
tree-domain, i.e., there is a unique minimal element $d\in D$ and for
every $d'\in D\setminus\{d\}$ the predecessor of $d'$ is contained in
$D$. Then, $T{\restriction}_D$ coincides with either
$T_{i-1}{\restriction}_D$ or 
$T_i{\restriction}_D$ or $T_{i+1}{\restriction}_D$. We denote by $H$ the tree
underlying $T$. The three cases 
correspond to the following three conditions on $H$.
\begin{enumerate}
\item If there are $e_1 < e_2 \leq d$ such that $e_1\in H^+$ and
  $e_2\in G^C$, then 
  \begin{align*}
    D\cap\domain(T)=\emptyset = D\cap
    \domain(T_{i-1}).    
  \end{align*}
\item If there is some $e_1\leq d$ such that $e_1\in H^+$ and $G$ does
  not contain any element between $e_1$ and $d$, then 
  $T{\restriction}_D$ coincides with $T_i{\restriction}_D$ 
  (cf. Lemma \ref{T1BasicFacts}).
\item If $d\in H$,
  $T{\restriction}_D$ coincides with 
  $T_{i+1}{\restriction}_D$. Moreover, the definition of a tree-comb
  implies that
  $T{\restriction}_D$ then coincides with 
  $T_k{\restriction}_D$ for all $k>i$. 
\end{enumerate}
Now, given a pairwise
$<_{G^C}$-comparable $n$-tuple $\bar T$ from $\Closure(C)$, $\bar T$
coincides locally with pairwise distinct $n$-tuples from $C$. We can
then define a function $\rho$ on $\bar T$ by
locally copying the accepting runs on the $n$-tuples from $C$. If $C$ is
homogeneous with respect to $\mathcal{A}$, $\rho$ is
an accepting run on $\bar T$ due to the following fact. 
Let  
\begin{align*}
  D_1 \subseteq \domain(T_i)^\oplus\setminus\domain(T_{i-1})
  \text{ and }
  D_2\subseteq \domain(T_{i+1})^\oplus\setminus\domain(T_{i})    
\end{align*}
be maximal tree domains such
that
$D_1$ and $D_2$ are touching. Then there are tuples 
\mbox{$\bar C_1, \bar C_2\in C$} such that $\rho$
coincides on $D_1$ with 
the accepting run on $\bar C_1$ and  $\rho$ coincides on $D_2$ with the
accepting run on $\bar C_2$. Due to homogeneity, the accepting run on
$\bar C_2$ coincides with the accepting run on $\bar C_1$  on the
path from the minimal element of $D_1$ to the minimal element of
$D_2$. Thus, $\rho$ respects the transition relation at the
border between $D_1$ and $D_2$. 
Since this argument applies at all borders where $\rho$ consists of
copies of different accepting runs, $\rho$ respects the transition
relation whence it is a run on
$\bar T$. Moreover,  the 
function
copies the behaviour of an accepting run on
each branch. Hence, the run is an accepting
run on $\bar T$. 

The precise definition of homogeneity is more complicated than
indicated above because we
have to deal with different permutations of fixed $n$-tuples. 
When we investigate pairwise $<_{G^C}$-comparable tuples, we can
order these tuples in various ways. But the accepting run for each
permutation of a tuple may differ from all the accepting runs on
the other permutations. Thus, our 
definition of homogeneity asserts that $n$-tuples from the tree-comb
share similar accepting runs if their elements are ordered by
$<_{G^C}$ in the same manner.
Let us first define some auxiliary notation. Afterwards, we state the
exact definition of homogeneity. 

\begin{definition}
  \begin{itemize}
  \item[]
  \item Let $C$ be some tree-comb and $D\subseteq\Closure(C)$. 
    Then we write $\left[\frac{D}{n}\right]$ for the set of 
    \mbox{$<_{G^C}$-increasing} $n$-tuples from $D$. 
  \item   For $\sigma:\{1, 2, \dots, n\}\to\{1, 2, \dots, n\}$ a
    permutation and $\bar x=(x_1, x_2, \dots, x_n)$ 
    some \mbox{$n$-tuple}, we write 
    $\sigma(\bar
    x):=(x_{\sigma(1)},x_{\sigma(2)},\dots,x_{\sigma(n)})$.  
  \item Let $D$ and $E$ be subsets  of some tree-comb $C$ such that
    $d <_{G^C} e$ for all $d\in D$ and $e\in E$. Furthermore, let
    $\mathcal{A}$ be a deterministic automaton 
    recognising an $n$-ary relation of $\Sigma$-trees. 
    We write  $\rho_{\sigma(\bar T)}$ for the run of $\mathcal{A}$ on
    $\bigotimes \sigma(\bar T)$ for all $n$-tuples $\bar T$ of trees and
    all permutations $\sigma$. 
    
    Set $F:=\bigcup_{T\in D} \domain(T)^\oplus$.
    We say  
    $\mathcal{A}$ \emph{runs homogeneously} on $E$ with respect to $D$, 
    if for each permutation $\sigma$ of $n$ elements and each number
    $0\leq m\leq n$ the following holds: 
    
    For all  $\bar T\in \left[\frac{D}{m}\right]$, 
    all $\bar U\in \left[\frac{E}{n-m}\right]$ and
    all $\bar V\in \left[\frac{E}{n-m}\right]$,
    \begin{align*}
      \rho_{\sigma(\bar T\bar U)}{\restriction}_F =
      \rho_{\sigma(\bar T\bar V)}{\restriction}_F, 
    \end{align*}
    i.e., the runs on $\bigotimes\sigma(\bar T\bar U)$ and
    $\bigotimes\sigma(\bar T\bar V)$ coincide on the domain $F$. 
    
    This means that different tuples from $E$ that are in the same order
    with respect to the tree-comb order $<_G$ have identical runs with each
    fixed tuple from $D$ on domain $F$, where $F$ may be seen as the
    ``maximal domain'' of $D$. 
    Note that this assertion is symmetric in the order in which we mix the
    tuple from $D$ with the tuples from $E$. 
  \end{itemize}
\end{definition}  
\begin{definition}
  Let $\mathcal{A}$ be some deterministic automaton and
  $C=(T_i)_{i\geq 1}$ be a  tree-comb. 
  Set 
  \begin{align*}
    D_n:=\{T_i: 1 \leq i < n\}\text{ and }E_n:=\{T_i: n\leq i\}\text{ for
      all } n\geq 1.    
  \end{align*}
  We say $C$ is
  \emph{homogeneous with respect to $\mathcal{A}$}, if, for all $n\in\N$,
  $\mathcal{A}$ runs  homogeneously on   $E_n$ with respect to $D_n$. 
\end{definition}

The crucial observation for the correctness proof is the following. Any
tree-comb whose pairwise distinct $n$-tuples are all accepted by an
automaton $\mathcal{A}$ contains a subcomb that is homogeneous
with respect to $\mathcal{A}$. Every set $M$ that witnesses the Ramsey
quantifier $\RamQ{n}\bar x (\psi)$ contains a tree-comb $C$
which also witnesses the Ramsey quantifier. We are going to show that $C$
contains a subcomb $C'$ which is homogeneous with respect to $\mathcal{A}_\psi$ 
(where $\mathcal{A}_\psi$ corresponds to $\psi$). 
Because of this homogeneity, we can then construct an accepting run of
$\mathcal{A}_\psi$ on
each pairwise comparable $n$-tuple from $\Closure(C')$. Since
$\mathcal{A}_\psi$ corresponds to $\psi$, this implies that every
pairwise comparable $n$-tuple from $\Closure(C')$ satisfies $\psi$. 
Thus, the representation of such a closure is  a witness for
$\reduction(\RamQ{n}\bar x(\psi))$.

\begin{lemma}
  Let $C$ be a tree-comb and $\mathcal{A}$ some 
  deterministic automaton such that $\mathcal{A}$ accepts 
  $\bigotimes \sigma(\bar  T)$ for all 
  $\bar T\in \left[ \frac{C}{n} \right]$ and all 
  permutations $\sigma$. Then there is a 
  subcomb $C^\mathcal{A}\subseteq C$ which is homogeneous
  with respect to $\mathcal{A}$. 
\end{lemma}
\begin{proof}
  We generate 
  $C^{\mathcal{A}}$ by the use of
  Ramsey's Theorem (Theorem \ref{Thm:RamseyTheorem}).

  For each $\bar T \in  \left[\frac{C}{n}\right]$ and each permutation
  $\sigma$, we denote the accepting 
  run  of $\mathcal{A}$ on $\bar T$ by $\rho_{\sigma(\bar T)}$.  

  We are going to define two infinite chains 
  \begin{align*}
    &D_0\subsetneq D_1\subsetneq D_2 \subsetneq\dots \text{ and}\\
    &C_0\supsetneq C_1\supsetneq C_2 \supsetneq\dots
  \end{align*}
 such that $\mathcal{A}$ runs homogeneously on $C_i$ with respect to
 $D_i$. $D_{i+1}$ will extend $D_i$ by exactly one finite tree
 $T_{i+1}$. The
 sequence of these trees then forms a tree-comb that is homogeneous
 with respect to $\mathcal{A}$. 

 At the beginning we set $D_0 := \emptyset$. Since
 $\emptyset^\oplus=\{\varepsilon\}$, we have to provide an infinite set
 \mbox{$C_0\subseteq C$} such that for each permutation $\sigma$ all $\bar
 T, \bar T' \in\left[\frac{C_0}{n}\right]$ satisfy 
 $\rho_{\sigma(\bar T)}(\varepsilon) = \rho_{\sigma(\bar
   T')}(\varepsilon)$. 
 This set $C_0$ can be obtained by applying Ramsey's Theorem as
 follows:
 For each \mbox{$\bar T\in\left[\frac{C_0}{n}\right]$} and for each
 permutation $\sigma$ the 
 function $\rho_{\sigma(\bar T)}{\restriction}_{\{\varepsilon\}}$ has
 finite domain and range. 
 Let $\sigma_1, \sigma_2, \dots, \sigma_m$ be a fixed enumeration of 
 all permutations of $n$ elements.
 Assigning
 \begin{align*}
   \bar T \mapsto 
   (\rho_{\sigma_1(\bar T)}{\restriction}_{\{\varepsilon\}}, 
   \rho_{\sigma_2(\bar T)}{\restriction}_{\{\varepsilon\}}, \dots, 
   \rho_{\sigma_m(\bar T)}{\restriction}_{\{\varepsilon\}})   
 \end{align*}
 for each $\bar
 T\in \left[\frac{C}{n}\right]$ induces a
 finite colouring of all $n$-element subsets of $C$:  
 since $C$ is linearly ordered by $<_{G^C}$ (see Lemma
 \ref{Lemma:TreeCombLinear}), each pairwise distinct
 $n$-tuple has a unique representative among the $<_{G^C}$ increasing
 sequences of  length $n$. Furthermore, the range of this map is
 finite.

 By Ramsey's
 theorem, there is an infinite subset $C_0\subseteq C$ that is
 homogeneous with respect to this colouring, i.e., if the
 $<_{C^G}$-order of 
two tuples from $C$ coincides, then their accepting runs
 coincide on the state at the root. 

 We now construct $D_{i+1}$ and $C_{i+1}$ from $D_i$ and $C_i$ by
 generalising this 
 process. 
 Assume that $D_i, C_i
 \subseteq C$ are disjoint sets
 such that $D_i$ is finite, $C_i$ is infinite,  and $T<_{G^C} T'$ for each
 $T\in D_i$ and $T'\in C_i$. 
 Furthermore,  assume that $\mathcal{A}$ runs   homogeneously on 
 $C_i$ with respect to $D_i$. 
 
 Let $T_i$ be the minimal element of $C_i$ with respect to
 $<_{G^C}$. We set 
 $D_{i+1}:=D_i\cup\{T_i\}$. Applying
 Ramsey's Theorem iteratively for each $0\leq k<n$ and each
 $\bar T\in \left[\frac{D_{i+1}}{k}\right]$, we can choose an
 infinite $C_{i+1}\subseteq C_i\setminus\{T_i\}$ such that $\mathcal{A}$
 runs  homogeneously on $C_{i+1}$ with respect to $D_{i+1}$. We
 explain one of these applications of 
 Ramsey's Theorem in  detail: 
 
 Fix a number $k\leq n$ and some $\bar T\in \left[\frac{D_{i+1}}{k}\right]$.
 In this step we consider the colouring that maps each 
 $\bar U\in \left[\frac{C_{i}}{n-k}\right]$ to
 \begin{align*}
   &(\rho_{\sigma_1(\bar T\bar U)}{\restriction}_F, \rho_{\sigma_2(\bar
     T \bar U)}{\restriction}_F, \dots, 
   \rho_{\sigma_m(\bar T \bar U)}{\restriction}_F)\text{ where}\\
   &F:=\bigcup_{T\in D_{i+1}} \domain(T)^\oplus.   
 \end{align*}
 Since $F$ is finite, 
 this induces a colouring of finite range on the $k$-tuples of $C_i$. 
 Applying Ramsey's Theorem, there is a
 homogeneous infinite subset $C'\subseteq C_i$ with respect to this
 colouring.  

 Iterating this
 construction for each $k\leq n$ and each $\bar T\in
 \left[\frac{D_{i+1}}{k}\right]$, we obtain a subset 
 $C_{i+1}\subseteq C'\subseteq C_i$ such that $\mathcal{A}$ runs
 homogeneously on $C_{i+1}$ with 
 respect to $D_{i+1}$. 
   
 Furthermore, it is clear that $T <_{G^C} T'$ for all 
 $T\in D_{i+1}$ and $T' \in C_{i+1}$ because of the following facts. 
 The same claim is true
 for $D_i$ and $C_i$. Furthermore, $D_{i+1}$ is $D_i$ extended by the
 minimal element of $C_i$ and $C_{i+1}$ does not contain the minimal
 element of $C_i$.  
 
 Repeating this construction, we obtain a sequence of trees $T_1, T_2,
 T_3, \dots$ such that \mbox{$D_i=\bigcup_{j\leq i} \{T_j\}$.}
 The sequence  $(T_i)_{i \geq 1} $ is a
 tree-comb because it is a subsequence of $C$. We set
 $C^{\mathcal{A}}:=(T_i)_{i\geq 1}$. Note that $C^{\mathcal{A}}$ is
 homogeneous with respect to $\mathcal{A}$ by construction. 
\end{proof}

Next, we show the following. 
Let $C$ be some tree-comb such that all pairwise distinct $n$-tuples
from $C$ are accepted by some automaton $\mathcal{A}$. 
If $C$ is homogeneous with respect to
$\mathcal{A}$, then  $\mathcal{A}$ accepts all
 pairwise comparable $n$-tuples of
$\Closure(C)$. The proof of this claim
relies on the fact that every tree in $\Closure(C)$ is locally
similar to one of the trees in $C$ and accepting runs for different
trees coincide on equal prefixes.

In the next lemma we use the following notation. 
Let $T$ be some tree and $D\subseteq \domain(T)^\oplus$ an initial segment. 
We call a map $\rho:D \to Q$ a \emph{partial run} of $\mathcal{A}$
on $T$ if $\rho$ respects the transition relation of $\mathcal{A}$ on
$T$. We call a partial run $\rho$ accepting if $\rho(\varepsilon)$ is
a final state. 
Let $\rho$ be some partial run on some tree $T$. 
For some  $d\in\{0,1\}^*$ we call $\rho$ total and correctly
initialised (\tci) towards $d$ if 
there is some $e\leq d$ such that $e\in \domain(T)^+$
and $\rho(e) = q_I$, i.e., the domain of $T$ ends at some ancestor
of $d$ and the border of its 
domain along this branch is labelled by the initial state. 
\begin{remark}
  Note that a partial run $\rho$ on some tree $T$ is an accepting run
  if and only if  it is accepting and \tci towards all $d\in\domain(T)^+$. 
\end{remark}

\begin{lemma}
  Let $\mathcal{A}=(Q,\Sigma,q_I,F,\Delta)$ be some deterministic
  automaton.
  Let $C$ be a tree-comb such that $\mathcal{A}$ accepts all pairwise
  disjoint $n$-tuples of $C$. 
  If $C$ is homogeneous with respect to $\mathcal{A}$, 
  then $\mathcal{A}$ accepts all pairwise comparable $n$-tuples from
  $\Closure(C)$, i.e., $\mathcal{A}$ accepts $\sigma(T)$ for all
  permutations $\sigma$ 
  and all $\bar T=(T_1, T_2, \dots, T_n)\in
  \left[\frac{\Closure(C)}{n}\right]$.
\end{lemma}
\begin{proof}
  We write $(T_1^C, T_2^C, G^C)$ for  the representation of
  $C=(C_i)_{i\geq 1}$.  
  Assume that 
  \begin{align*}
    \bar T= (T_1, T_2, \dots, T_n) \in
    \left[\frac{\Closure(C)}{n}\right].    
  \end{align*}
   Furthermore, assume that  $H_1, H_2, \dots, H_n$ 
  are the $G^C$-trees that underlie the trees $T_1, T_2, \dots, T_n$
  (i.e.,  for each $i$, $T_i=\pretree(H_i, T_1^C, T_2^C, G^C)$ and
  $H_i$ is maximal with this property). 

  We assume 
  that $\sigma=\Id$ (due to the symmetric definition of homogeneity,
  the proof is completely analogous for any other permutation). 

  For each $\bar C\in\left[\frac{C}{n}\right]$, we write 
  $\rho_{\bar C}$ for the accepting run of  
  $\mathcal{A}$ on $\bar C$.   

  Set $F_k:=\bigcup_{i<k} \domain(C_i)^\oplus$. We will
  define an accepting run $\rho_{\bar T}$ of $\mathcal{A}$ on
  $\bigotimes \bar T$ as the union of accepting  partial runs 
  $\rho_{\bar T}{\restriction}_{F_k}$. 

  We start with the definition of 
  $\rho_{\bar T}{\restriction}_{F_1}$. Note that $F_1=\{\varepsilon\}$. 
  We set $\rho_{\bar T}(\varepsilon) := \rho_{\bar C}(\varepsilon)$ for an
  arbitrary $\bar C \in \left[\frac{C}{n}\right]$. Recall that by
  homogeneity of $\bar C$, this definition is independent of the
  concrete choice of $\bar C$. Furthermore, since $\rho_{\bar C}$ is
  accepting $\rho_{\bar T}{\restriction}_{F_1}$ is an accepting
  partial run. 

  For $d=\varepsilon, k=1$, and $m=1$, 
  $\rho_{\bar T}{\restriction}_{F_1}$ satisfies the
  following properties. 
  \begin{enumerate}
  \item $d\in \domain(C_{k-1})^+ \setminus \bigcup_{j<k-1}
    (\domain(C_j)^\oplus)$ (where we define $C_0:=\emptyset$),
  \item $d\in\domain(T_m)^\oplus$ (just by definition of $^\oplus$),
  \item $d\in H_j^\oplus$ for $m\leq j \leq n$, and
  \item $\rho_{\bar T}(d) = \rho_{\bar C} (d)$ for any $\bar
    C=C_{i_1}, C_{i_2}, \dots, C_{i_n}$ with $k\leq i_1 < i_2 <
    \dots < i_n$. 
  \end{enumerate}
  For each $k\geq 1$,  we inductively extend the accepting partial run
  $\rho_{\bar T}{\restriction}_{F_{k-1}}$ to an accepting partial run on domain 
  \mbox{$F_k\cap\domain(\bigotimes \bar T)^\oplus$} . 
  In each step of this construction, we preserve the
  property that  
  for each maximal element $d\in F_k$, at least one of the following
  conditions hold.  
  \begin{enumerate}
  \item $\rho_{\bar T}$ is a \tci accepting partial run on $\bigotimes
    \bar T$ towards $d$. 
  \item There is some $1\leq m \leq n$ such that the following
    conditions are  satisfied:
    \begin{enumerate}
    \item $d\in\domain(C_{k-1})^+\setminus F_{k-1}$, i.e., 
      $d\in \domain(C_{k-1})^+)\setminus\bigcup_{j<k-1}(\domain(C_j)^\oplus)$,
    \item $d\notin \domain(T_j)$ for all $1\leq j < m$,
    \item $d\notin H_j^\oplus$ for $1\leq j < m$,
    \item $d\in\domain(T_m)^\oplus$,
    \item $d\in H_j^\oplus$ for all $m\leq j \leq n$, \label{ConditiondInHgrosj}
    \item there are natural numbers 
      \mbox{$1\leq i_1 < i_2 < \dots < i_{m-1} < k\leq i_m<i_{m+1}<
        \dots <i_n$} 
      such that \mbox{$\rho_{\bar T}(d) = \rho_{\bar C}(d)$} where $\bar
      C:=(C_{i_1}, C_{i_2}, \dots, C_{i_n})$. We stress that
      $\rho_{\bar C}(d)$ does not depend 
      on the concrete choice of $i_m, i_{m+1}, \dots, i_n$, i.e., for all 
      \mbox{$\bar C' :=(C_{i_1}, C_{i_2}, \dots, C_{i_{m-1}}, C_{i_m'},
      C_{i_{m+1}'}, \dots C_{i_n'})$} with $k\leq i_m' < i_{m+1}' <
      \dots < i_n'$, we have $\rho_{\bar T}(d) = \rho_{\bar C}(d)$
      due to the homogeneity of  $\mathcal{A}$ on $(C)_{i\geq k}$ with
      respect to $(C)_{i<k}$. \label{BarCCondition}
    \end{enumerate}
  \end{enumerate}
  Note that these conditions imply that $\rho_{\bar
    T}{\restriction}_{F_k}$ is defined on 
  $F_k \cap \domain (\bigotimes \bar T)$ and that 
  $\rho_{\bar  T}{\restriction}_{F_k}$ may
  be extendable to an accepting run of $\mathcal{A}$ on $\bar
  T$. Especially, if the first condition applies to all
  maximal $d\in F_k$, then $\rho_{\bar T}$ is an
  accepting run on $\bigotimes \bar T$. 
  
  We now extend $\rho_{\bar T}$ to the maximal possible
  segment of  \mbox{$F_{k+1}=\bigcup_{i<k+1} \domain(C_i)^\oplus$}, i.e., to
  $F_{k+1} \cap (\domain(\bigotimes\bar T)^\oplus)$. 

  By assumption, we only have to extend $\rho_{\bar T}$ at the maximal
  positions $d\in F_k$ where the second condition holds. 
  We will distinguish the following three cases. 
  \begin{enumerate}
  \item $d\notin\domain(C_j)$
    for all $j\geq 1$,
  \item $d\in \domain(C_k)$ and $d\in H_m$, and
  \item   $d\in\domain(C_k)$ but $d\notin H_m$. 
  \end{enumerate}
  Let us first explain why this case distinction is complete:
  Assume that there is some $j\geq 1$ such that $d\in \domain(C_j)$.
  Then $j\geq 
  k$ because $d\in\domain(C_{k-1})^+\setminus F_{k-1}$. Since
  $d\in\domain(C_{k-1})^\oplus$ and due to the tree-comb property,
  $d\in\domain(C_j)$ for some $j\geq k$ if and only if
  $d\in\domain(C_j)$ for all $j\geq k$. Thus, we conclude that
  $d\in\domain(C_k)$.  

  For the case distinction, let us fix a tuple $\bar C=(C_{i_1},
  C_{i_2}, \dots, C_{i_n})$ witnessing condition \ref{BarCCondition}. 
  \begin{enumerate}
  \item $d\notin\domain(C_j)$ for all $j\geq1$: first of all, note
    that in this case either $d\in\domain(C_{i_j})^+$ for some $1\leq j
    \leq n$ and $\rho_{\bar C}(d) = q_I$ or
    $d\notin\domain(C_{i_j})^\oplus$ for all $1\leq j \leq n$ and 
    $d\notin\domain(\rho_{\bar C})$. 
    
    Secondly, by definition of $G^C$, we have
    $G^C\cap\{e:d<e\} = \emptyset$ whence $d\notin H$ for all $G$-trees
    $H$. Especially, $d\notin H_i$ for $1\leq i \leq n$. 
    Furthermore, by Lemma \ref{T1BasicFacts}
    $\extract(d,T^C_2,G^C)=\emptyset$.
    Hence, $d\notin\domain(T_i)$ for all $i$. Recall that
    $d\in\domain(T_m)^\oplus$ by assumption, whence 
    $d\in\domain(\bigotimes \bar T)^+$ and furthermore, 
    $d\in F_k\cap (\domain(\bigotimes \bar T)^\oplus)$. Thus,
    $\rho_{\bar T}{\restriction}_{F_k}$ is defined at $d$. 
    
    Putting these two facts together, it is only possible that
    $\rho_{\bar C}$ and $\rho_{\bar T}$ agree on $d$ if 
    $\rho_{\bar T}(d) = \rho_{\bar C}(d) = q_I$ whence $\rho_{\bar T}$
    is an accepting partial run on $\bigotimes \bar T$ that is \tci
    towards $d$.  
    Thus,  the first condition is satisfied
    for all  
    $d'\in \{e: d\leq e\}$. 
  \item $d\in\domain(C_k)$ and $d\in H_m$:
    first of all, we claim that  $\domain(C_k)\cap \{e: d\leq e\}
    \subseteq H_m$. 

    By definition of $H_m$ and  $G^C$,
    $e\in H_m^+$ implies that $e\in
    \domain(C_{l})^+\setminus\bigcup_{l'<l} (\domain(C_{l'})^\oplus)$
    for some $l\in\N$. Due to 
    $d\in\domain(C_{k-1})^+\setminus F_{k-1}$, no proper successor $e$
    of $d$ is contained in $\domain(C_{l})^+$ for $l<k$. Thus, the
    first descendants of $d$ that are contained in $G$ are contained
    in $\domain(C_k)^+$. Thus, all elements between $d$ and
    $\domain(C_k)^+$ are contained in $H_m$. 
    
    Since  $H_m\subsetneq H_j$, the same holds for all $H_j$ with
    $j\geq m$.  
    Thus,  $T_m, T_{m+1}, \dots, T_{n}$ agree with 
    $C_{k+1}, C_{k+2},\dots, C_{k+n-m}$ on
    $D:=\domain(C_k)\cap \{e:d\leq e\}$ (cf. Lemma \ref{T1BasicFacts}). 
    
    Furthermore, $D\cap \domain(T_j)=\emptyset$ for $j < m$ by the
    assumption that $d\notin\domain(T_j)$ for all $j < m$. 
    Similarly, $d\notin \domain(C_{i_j})$ for $j<m$ due to $i_j<k$ and
    $d\notin\bigcup_{l<k} \domain(C_l)$ (recall that $d$ is maximal in
    $F_k$. 
    
    Thus,  $T_j$ agrees with $C_{i_j}$ for all $j<m$ on the subtree
    rooted at $d$.
    
    It follows that, for  
    $\bar C':=(C_{i_1}, C_{i_2}, \dots, C_{i_{m-1}},
    C_{k+1}, C_{k+2}, \dots C_{k+n-m})$,
    the trees $\bigotimes \bar C'$ and $\bigotimes \bar T$ agree on
    $D$. By condition \ref{BarCCondition},  
    $\rho_{\bar T}(d) = \rho_{\bar C'}(d) = \rho_{\bar C}(d)$.  
    Thus, setting 
    $\rho_{\bar T}(e):=\rho_{\bar C'}(e)$ for all $e\in D^\oplus\cap
    (\domain(\bigotimes \bar C')^\oplus)$
    extends $\rho_{\bar T}$ in such a way that it still is a
    partial run on $\bigotimes\bar T$. Note that the maximal
    elements of $D^\oplus$ are by definition maximal elements of $F_{k+1}$. 
    We claim that for any such element $d'$ condition $1$ or condition
    $2$ holds. There are the following cases. 
    
    For the first case, assume that 
    $d'\notin\domain(\rho_{\bar T})$. This implies that 
    $d'\notin\domain(\bigotimes\bar T)^\oplus$ whence by
    coincidence of $\bigotimes \bar T$ and $\bigotimes \bar C'$ on
    $D$, it follows that  $d'\notin
    \domain(\bigotimes \bar C')^\oplus$. Thus,
    $\bigotimes \bar C'$ and 
    $\bigotimes \bar T$ agree on the path from $d$ to $d'$. 
    $\rho_{\bar C'}$
    is \tci on $\bigotimes \bar C'$ towards $d'$ because 
    $d'\notin\domain(\bigotimes \bar C')$.
    Since $\rho_{\bar C'}$ and $\rho_{\bar T}$ agree on this path, 
    $\rho_{\bar T}$ is \tci on $\bigotimes \bar T$ towards $d'$.
    
    For the second case, assume that $d'\in\domain(\rho_{\bar T})$. In
    this case, we show that the second condition holds for $k+1$ and $m$.
    \begin{enumerate}
    \item we have to show that $d'\in\domain(C_k)^+\setminus
      F_k$. $d'\in\domain(C_k)^+$ follows from its definition while
      $d'\notin F_k$ follows from the facts that $d<d'$, $d\notin
      F_{k-1}$ and \mbox{$d\in\domain(C_{k-1})^+$}: note that
      $F_k= F_{k-1}\cup  \domain(C_{k-1})^\oplus$ and $d$ is by
      definition a maximal element of this set.
    \item $d'\notin\domain(T_j)^\oplus$ for all $1\leq j < m$ because
      $d<d'$ and $d\notin\domain(T_j)$ for all $1\leq j < m$ by
      assumption.
    \item $d'\notin H_j^\oplus$ for $1\leq j < m$ because $d<d'$ and 
      $d\notin H_j^\oplus$ for $1\leq j < m$ by assumption.
    \item Since $d'\in \domain(\rho_{\bar T})$, 
      $d'\in \bigcup_{j=1}^n \domain(T_j)^\oplus$. Thus, there is some 
      $1\leq j \leq n$ such that $d'\in\domain(T_j)^\oplus$. b)
      implies that $j\geq m$. Furthermore, due to 
      \begin{align*}
        \domain(C_k)\cap\{e:d\leq e\} \subseteq H_m \subseteq H_{m+1}
        \subseteq \dots \subseteq H_n,  
      \end{align*}
      the 
      trees $T_m, T_{m+1}, \dots, T_n$ agree on 
      $\domain(C_k)\cap\{e:d\leq e\}$. But the predecessor of
      $d'$ is contained in $\domain(C_k)$. Thus, 
      we conclude that 
      $d'\in\domain(T_j)^\oplus$
      for all  $j\geq m$.  
    \item $d'\in H_j^\oplus$ for all $j\geq m$ follows directly from $d<d'$,
      $d'\in\domain(C_k)^+$ and 
      \mbox{$\domain(C_k)\cap\{e:d\leq e\}\subseteq H_j$.}
    \item According to the definition of $\rho_{\bar T}$ on $D^\oplus$, 
      $\rho_{\bar T}(d') = \rho_{\bar C'}(d')$ where 
      \begin{align*}
        \bar C' =(C_{i_1}, C_{i_2}, \dots, C_{i_{m-1}}, C_{k+1},
        C_{k+2}, \dots, C_{k+m-n}).        
      \end{align*}
      Hence, this tuple witnesses 
      condition \ref{BarCCondition}. 
    \end{enumerate}
  \item $d\in\domain(C_k)$ and $d\notin H_m$:
    due to \ref{ConditiondInHgrosj}, $d\in H_m^\oplus$ whence we know
    that $d\in H_m^+$. 
    Furthermore, $d\notin F_{k-1} = \bigcup_{j<k-1}
    (\domain(C_j)^\oplus)$. 
    Since $H_m$ is a $G^C$-tree, we conclude that
    $d\in\domain(C_{k-1})^+$. Due to $d\in\domain(C_k)$, it follows
    immediately that 
    \begin{align*}
      \emptyset\neq \domain(C_k)^+ \cap \{e: d < e\}\subseteq G^C.      
    \end{align*}
    Thus,
    for all $j>m$, $H_j$ extends $H_m$ along the subtree rooted at $d$
    because  $H_m <_{G^C} H_j$. 

    For $m<j\leq n$, this implies $d\in H_j$  whence
    $\domain(C_k)\cap\{e:d\leq e \} \subseteq H_j$.
    Hence, $T_{m+1}, T_{m+2}, \dots, T_n$ agree with $C_{k+1}, C_{k+2},
    \dots, C_{k+n-m}$ on $\domain(C_k)\cap\{e:d\leq e\}$.

    Furthermore, Lemma \ref{T1BasicFacts} implies that $C_k$ and $T_m$
    agree on $\{e:d\leq e\}$. 

    Condition 2b implies $d\notin\domain(T_j)$
    for $j<m$. By condition 2a, $d\notin\domain(C_{i_j})$ for $j<m$
    (recall that the $i_j$ are defined as in condition
    \ref{BarCCondition}).      

    We conclude that $\bar T$ and $\bar C':=(C_{i_1}, C_{i_2}, \dots,
    C_{i_{m-1}}, C_k, C_{k+1}, \dots, C_{k+n-m})$ agree on
    $\domain(C_k)\cap\{e: d\leq e\}$. 

    Due to condition \ref{BarCCondition}, $\rho_{\bar T}(d)=\rho_{\bar
      C}(d) = \rho_{\bar C'}(d)$. 
    Thus, setting $\rho_{\bar T}(e):= \rho_{\bar C}(e)$ for all 
    $e\in D:=\domain(C_k)^\oplus \cap\{e: d< e\}$ extends
    $\rho_{\bar T}$ in such a way that it is still a partial run.
    
    Note that the maximal elements of $\domain(C_k)^\oplus \cap
    \{e:d\leq e\}$ are the maximal elements of $F_{k+1}\cap\{e:d\leq
    e\}$. 
    We claim that for each maximal element $d'$ in
    \mbox{$\domain(C_k)^\oplus\cap \{e:d\leq e\}$} condition $1$ or condition
    $2$ with $k+1$ and $m+1$ are satisfied. Again, we prove this claim
    by case distinction. 

    For the first case, assume that $d'\notin\domain(\rho_{\bar
      T})$. This implies that 
    $d'\notin\domain(\bigotimes\bar T)^\oplus$ whence by
    coincidence of $\bigotimes \bar T$ and $\bigotimes \bar C'$ on
    all $d\leq e < d'$, it follows that  \mbox{$d'\notin
    \domain(\bigotimes \bar C')^\oplus$}. Thus,
    $\bigotimes \bar C'$ and 
    $\bigotimes \bar T$ agree on the path from $d$ to $d'$
    and $\rho_{\bar C'}$
    is \tci on $\bigotimes \bar C'$ towards $d'$ because 
    $d'\notin\domain(\bigotimes \bar C')$ and $\rho_{\bar C'}$ is an
    accepting run on  $\bigotimes\bar C'$. But for the coincidence of
    $\rho_{\bar C'}$ and $\rho_{\bar T}$ on this path, 
    $\rho_{\bar T}$ is then \tci on $\bigotimes \bar T$ towards $d'$.
    
    For the second case, assume that $d'\in\domain(\rho_{\bar T})$. 
    We show that condition 1 or condition 2  is
    satisfied; 
    if $m=n$ or $m<n$ and $d'\notin\domain(T_{m+1})^\oplus$, condition 1
    is satisfied, i.e., the partial run $\rho_{\bar T}$ is \tci towards $d'$. 
    Otherwise, we show that condition 2 is satisfied. 

    Independent of the case we are in, we first show conditions 2a --
    2c. These
    are also  helpful when discussing the cases $m=n$ or
    $d'\notin\domain(T_{m+1})^\oplus$.
    \begin{enumerate}
    \item We show that $d'\in\domain(C_k)^+\setminus
      F_k$. $d'\in\domain(C_k)^+$ follows from its definition while
      $d'\notin F_k$ follows from $d<d'$, $d\notin F_{k-1}$ and
      $d\in\domain(C_{k-1})^+$: note that $F_k= F_{k-1}\cup
      \domain(C_{k-1})^\oplus$ and $d$ is by definition a maximal
      element of this set.
    \item $d'\notin\domain(T_j)^\oplus$ for all $1\leq j < m$ because
      $d<d'$ and $d\notin\domain(T_j)$ for all $1\leq j < m$ by
      induction hypothesis.
      Furthermore, $d'\in\domain(C_k)^+$ and $T_m$ agrees with
      $C_k$ on  $\{e:d\leq e\}$. Thus,
      $d'\in\domain(T_m)^+$ whence $d'\notin\domain(T_m)$. 
    \item $d'\notin H_j^\oplus$ for $1\leq j < m$ because $d\notin
      H_j^\oplus$ for $1\leq j < m$ and $d<d'$. Since we
      are in the case $d\notin H_m$, we also have $d'\notin H_m^\oplus$
      because $d < d'$. 
    \item In order to satisfy condition 2d, we would have to show
      that $m+1\leq n$ and $d'\in\domain(T_{m+1})^\oplus$. 
      Instead, we show the
      following: if this is not the case, then condition 1 is
      satisfied. 

      First assume that $m=n$. We have seen that
      $d'\notin\domain(T_j)$ for $j\leq m=n$. Since $\bar T$
      and $\bar C'$ coincide on $d\leq e$,
      $d'\notin\domain(\bigotimes \bar C')$ and $\rho_{\bar C'}$ is
      \tci on $\bigotimes\bar C'$ towards $d'$. But then
      $\rho_{\bar T}$ is \tci on $\bigotimes \bar T$ towards $d'$
      because $\rho_{\bar T}$ and $\rho_{\bar C'}$ agree on $d\leq
      e \leq d'$. 

      Now assume that $m<n$ and $d'\notin\domain(T_{m+1})^\oplus$. 
      We have already seen that 
      $\domain(C_k)\cap\{e:d\leq e\}\subseteq H_{m+1}$.  Hence,
      the predecessor of $d'$ is in $H_{m+1}$. Thus,
      $d'\notin\domain(T_{m+1})^\oplus$ implies that either $d'\in
      H_{m+1}$ and $T_1^C(d')=\Box$ or $d'\notin H_{m+1}$ and
      $T_2^C(d')=\Box$. 
      Recalling the  definitions of $T_1^C$ and $T_2^C$, we 
      conclude that in the first case there is at most one $j\geq 1$
      such that $C_j(d')$ is defined while in the second case there is
      no $j\geq 1$ such that $C_j(d')$ is defined. 

      Heading for a contradiction, assume that there is
      a $j\in\N$ such that $C_j(d')$ is defined. By a), we know $j >
      k$. Due to  the tree-comb 
      property and because of $d'\in \domain(C_k)^\oplus$, we know
      that $C_{i}$ and 
      $C_j$ agree on $d'$ for all $k<i<j$. Thus, we arrive at the
      contradiction that there are infinitely many $j$ such that
      $C_j(d')$ is defined. 
      
      Thus, $d'\notin\domain(C_j)$ for all $j\geq 1$ whence
      $d'\notin\domain(T_j)$ for all $1\leq j \leq n$. But then
      $\rho_{\bar C'}$ is \tci on $\bigotimes\bar C'$ towards
      $d'$. 
      For the coincidence 
      of $\rho_{\bar C'}$ and $\rho_{\bar T}$ on the path between
      $d$ and $d'$, $\rho_{\bar T}$ is \tci on
      $\bigotimes\bar T$ towards $d'$.  

      We conclude that either
      $\rho_{\bar T}$ is \tci towards $d'$ whence $d'$ satisfies
      condition 1 or $m<n$ and 
      $d'\in\domain(T_{m+1})^\oplus$ whence it satisfies condition 2d
      for $m+1$. In case that $m<n$ and
      $d'\in\domain(T_{m+1})^\oplus$, we continue by showing that
      conditions 2e and 
      2f are also satisfied.
    \item $d'\in H_j^\oplus$ for $m+1\leq j \leq n$ follows directly from
      $d'\in\domain(C_k)^+$ and \mbox{$\domain(C_k)\cap\{e:d\leq e\}
      \subseteq H_j$}.
    \item By the very definition, $\bar C'$ is a witness for the claim
      $\rho_{\bar T}(d') = \rho_{\bar C'}(d')$ whence condition
      \ref{BarCCondition} is satisfied. 
    \end{enumerate}
    This completes the third case.
    We have shown that one of the following holds:
    \begin{itemize}
    \item $m=n$ and $d'$ satisfies condition 1, i.e., $\rho_{\bar T}$
      is \tci towards $d'$.
    \item $m<n$, $d'\notin\domain(T_{m+1})^\oplus$ and $d'$ satisfies
      condition 1. 
    \item $m<n$, $d'\in\domain(T_{m+1})^\oplus$
    and condition $2$ is satisfied for $k$ replaced by $k+1$ and $m$
    replaced by $m+1$. 
    \end{itemize}
  \end{enumerate}
  Repeating this inductive definition for all $k\in\N$, we define
  a partial run $\rho_{\bar T}$ on 
  \mbox{$\domain(\bigotimes \bar T) \cap F_k$} for all $k\in \N$. 
  Note that this inductive process terminates at some step because
  $\bigotimes \bar T$ is a finite tree with 
  $\domain(\bar T)\subseteq \bigcup_{k\in\N} F_k$. 
  Due to the finiteness of $\bigotimes \bar T$, there is some $i\in\N$
  such that $\domain(\bar T) \subseteq \bigcup_{k=0}^n F_k$.

  Note that this process stops if and only if all maximal elements of
  $\domain(\rho_{\bar T})$ satisfy condition $1$,
  i.e., $\domain(\rho_{\bar T})$ is \tci on $\bigotimes \bar T$
  towards all $d\in\domain(\bigotimes \bar T)^+$. This is equivalent
  to the fact that 
  $\rho_{\bar T}$ is an  run on $\bigotimes \bar T$. Since its root is
  labelled by an accepting state,
  we have constructed an accepting run of $\mathcal{A}$ on $\bar T$ as
  required by the lemma. 
\end{proof}

By now, we have obtained the following result. 
For each tree-comb $C$ whose pairwise distinct
$n$-tuples are accepted by some automaton
$\mathcal{A}$, there is a subcomb $C'$ that is homogeneous with
respect to $\mathcal{A}$. Due to homogeneity, all pairwise comparable
$n$-tuples 
from the closure of $C'$ are accepted by $\mathcal{A}$. 

We apply this result in order to prove the correctness of our
reduction of the Ramsey quantifier and to prove decidability of the 
$\FO{}(\exists^\mathrm{mod}, (\RamQ{n})_{n\in\N})$-theory of
automatic structures. We prove these two facts
simultaneously. The correctness of the reduction relies on the fact
that every formula induces an  automaton that corresponds to this
formula. On the other hand, the correctness of the reduction allows
the construction of an automaton corresponding to a formula. Thus, we
prove both facts by parallel induction. Let us start with an
auxiliary lemma that allows to extend the correctness proof for  one
construction step.  

\begin{lemma}
  Let $\varphi\in\FO{}(\exists^\mathrm{mod}, (\RamQ{n})_{n\in\N})$ be
  a formula such that for each proper subformula $\psi$ of $\varphi$,
  there is an automaton $\mathcal{A}_\psi$ that corresponds to $\psi$
  on $\mathfrak{A}$. 
  For each $\bar a \in \mathfrak{A}$, 
  $\mathfrak{A}, \bar a \models\varphi$ implies
  $\OmegaExp(\mathfrak{A}), \bar a \models\reduction(\varphi)$. 
\end{lemma}
\begin{proof}
  Except for the
  case $\varphi =\RamQ{n} \bar x (\psi)$,   the inductive proof of
  this claim is straightforward (for these cases we even 
  do not need the fact that there is a corresponding automaton for
  each proper subformula).

  Assume that $\varphi =\RamQ{n} \bar x (\psi)$ and assume that there
  is some $\bar a=a_1, a_2, \dots, a_m  \in\mathfrak{A}$ such that
  $\mathfrak{A}, (a_1, a_2, \dots, a_m)\models\varphi(y_1, y_2,
  \dots, y_m)$. 
  By assumption, we know that there is an automaton $\mathcal{A}_\psi$
  corresponding to $\psi$, i.e., 
  for all $\bar b=b_1, b_2, \dots, b_n \in \mathfrak{A}$, 
  \begin{align*}
    \mathfrak{A}, (b_1, b_2, \dots, b_n, a_1, a_2, \dots a_m) \models
    \psi(\bar x, \bar y)    \text{  if and only if }
    \mathcal{A}_\psi\text{ accepts }\bigotimes \bar b \otimes\bigotimes \bar
    a.
  \end{align*}
  Due to $\mathfrak{A}, \bar a \models \varphi$, 
  there is an infinite set $S\subseteq
  \mathfrak{A}$ such that 
  \begin{align*}
    \mathfrak{A},(b_1, b_2, \dots, b_n, a_1, a_2, \dots a_m) \models
    \psi(x_1, x_2, \dots, x_n, y_1, y_2, \dots, y_m)      
  \end{align*}
  for all pairwise distinct $n$-tuples $\bar b= b_1, b_2, \dots, b_n$
  from $S$.  
  By Lemma \ref{L combs in sets}, there is a tree-comb $C'$ contained
  in $S$. 
  This tree-comb contains a subcomb $C$ that is homogeneous 
  with
  respect to the automata $\mathcal{A}_A$
  and  $\mathcal{A}_\psi$ where $\mathcal{A}_A$  
  recognises the domain of $\mathfrak{A}$. 
  By the previous lemma, it follows that $(T^C_1, T^C_2, G^C)$
  witnesses
  \begin{align*}
    \OmegaExp(\mathfrak{A}), (a_1, a_2, \dots,
    a_m)\models\reduction(\varphi)    
  \end{align*}
  due to the 
  following facts:
  \begin{enumerate}
  \item $(T_1^C, T_2^C, G^C)$ is coherent due to Lemma
    \ref{Lemma:CombCoherent},
  \item $(T_1^C, G^C)$ is \SMALL due to Lemma \ref{LemmaCombSmall}
  \item Since $C$ is homogeneous with respect to $\mathcal{A}_A$, the
    previous lemma shows that $x\in A$ for all $x\in\Closure(C)$.
  \item Since $C$ is homogeneous with respect to $\mathcal{A}_\psi$, the
    previous lemma shows that for each pairwise comparable $n$-tuple
    $\bar b$ from
    $\Closure(C)$, 
    \begin{align*}
      \mathfrak{A}, \bar b, \bar a
      \models \psi(x_1, x_2, \dots, x_n, y_1, y_2, \dots, y_m).      
    \end{align*}
    By
    induction hypothesis, this implies that 
    \begin{align*}
      \OmegaExp(\mathfrak{A}), \bar b, \bar a
      \models \reduction(\psi)(x_1, x_2, \dots, x_n, y_1, y_2,
      \dots, y_m).
    \end{align*}
    Thus, $(T_1^C, T_2^C, G^C)$ witnesses $\reduction(\varphi)$ which
    concludes the proof. \qedhere
  \end{enumerate}
\end{proof}

Using the previous lemma, we can now prove that
there is an automata construction corresponding to the Ramsey
quantifier. 

\begin{lemma}
  Let $\varphi(\bar y) \in\FO{}(\exists^\mathrm{mod}, (\RamQ{n})_{n\in\N})$ be
  a formula with free variables among $\bar y$. Then there is an automaton
  $\mathcal{A}_\varphi$ such that for all $\bar a \in \mathfrak{A}$
  \begin{align*}
    \mathfrak{A}, \bar a \models \varphi(\bar y)\text{ iff }
    \mathcal{A}_\varphi \text{ accepts } \bigotimes \bar a. 
  \end{align*}
\end{lemma}
\begin{proof}
  Except for the
  case $\varphi =\RamQ{n} \bar x (\psi)$, 
  the inductive proof 
  of this claim is a straightforward adaption of the proof of Lemma
  \ref{Thm:FOTreeAutomaticDecidable}. 
  
  Now, consider the case $\varphi(\bar y) = \RamQ{n} \bar x (\psi)$
  where $\bar y = y_1, y_2, \dots, y_m$. 
  By induction hypothesis there is an automaton
  $\mathcal{A}_\psi$ that corresponds to $\psi$ on $\mathfrak{A}$, 
  i.e., 
  \begin{align*}
    &\text{for all }\bar a= a_1, a_2, \dots, a_m\in \mathfrak{A}\text{ and
    }\bar b=b_1, b_2, \dots, b_n \in \mathfrak{A},\\  
    &\mathfrak{A}, \bar a, \bar b \models \psi(\bar x, \bar y)\text{
      iff}\\
    &\mathcal{A}_\psi\text{ accepts }\bigotimes \bar a\otimes \bigotimes \bar b. 
  \end{align*}
  Due to the soundness of the reduction and due to the previous lemma,
  this implies that 
  \begin{align} \label{SoundandCorrect}
    \mathfrak{A}, \bar a \models \varphi\text{ if and only if 
    }\OmegaExp(\mathfrak{A}), \bar a \models \reduction(\varphi). 
  \end{align}
  By Lemma \ref{Lemma:FinToOmegaAutomaton},
  the $\omega$-automaton $\mathcal{A}_\psi^\infty$ corresponds to
  $\reduction (\psi)$ on $\OmegaExp(\mathfrak{A})$ in the sense
  that
  for all $\bar a, \bar b \in \mathfrak{A}$, 
  \begin{align*}
    &\OmegaExp(\mathfrak{A}), \bar a, \bar b \models \reduction (\psi) \\
    \text{iff }   
    &\mathfrak{A}, \bar a, \bar b \models \psi  \\
    \text{iff }
    & \mathcal{A}_\psi \text{ accepts } \bigotimes \bar a
    \otimes \bigotimes \bar b \\ 
    \text{iff } 
    &\mathcal{A}_\psi^\infty \text{ accepts} \bigotimes \bar a \otimes
    \bigotimes \bar b. 
  \end{align*}
  Recall that the construction of $\reduction(\varphi) =
  \reduction(\RamQ{n} \bar x ( \psi))$ is first-order
  except for the construction of $\reduction(\psi)$. Thus, we can use 
  standard techniques in order to construct an $\omega$-automaton
  $\mathcal{A}^\infty_\varphi$ from $\mathcal{A}^\infty_\psi$ which corresponds to 
  $\varphi$ on $\OmegaExp(\mathfrak{A})$ in the sense that for all
  $\bar a \in \mathfrak{A}$, $\OmegaExp(\mathfrak{A}), \bar a \models
  \reduction(\varphi)$ if and only if $\mathcal{A}^\infty_\varphi$
  accepts $\bigotimes \bar a$. 

  But now, Lemma \ref{Lemma:OmegaToFinAutomaton} provides the automaton
  $\mathcal{A}_\varphi:=(\mathcal{A}^\infty_\varphi)^{\mathrm{fin}}$. 
  Due to \ref{SoundandCorrect}, we obtain that for all $\bar a \in
  \mathfrak{A}$, 
  \begin{align*}
    &\mathfrak{A}, \bar a \models \varphi \\
    \text{iff } &\OmegaExp(\mathfrak{A}),
    \bar a \models \reduction(\varphi)\\
    \text{iff }
    &\mathcal{A}^\infty_\varphi \text{ accepts } \bigotimes\bar a\\
    \text{iff } &\mathcal{A}_\varphi \text{ accepts } \bigotimes \bar
    a. 
  \end{align*}
  Thus, $\mathcal{A}_\varphi$ corresponds to $\varphi$ on
  $\mathfrak{A}$. This concludes our proof. 
\end{proof}
\begin{remark}
  Theorem \ref{Thm:RamseyQuantifierDecidable} is a direct corollary of
  the previous lemma. Every $\FO{}(\exists^\mathrm{mod},
  (\RamQ{n})_{n\in\N})$ formula can be effectively translated into a
  corresponding finite automaton on every given automatic structure. 
  This reduces the model checking problem to the membership problem of
  regular languages. The latter problem is decidable. 
\end{remark}

\subsection{Recurrent Reachability on Automatic Structures}

In this section we review To's and Libkin's result \cite{ToL08} on
the recurrent reachability problem for automatic structures. 
Unaware of the concept of word- or tree-combs, they constructed an
automaton for the recurrent reachability problem by hand. 
In fact, they designed an
automaton that looks for a tree-comb witnessing  recurrent
reachability. 

We first describe the recurrent reachability problem. Then we
show how our method can be adapted to solve this problem. In fact,
we only have to replace the role of pairwise comparable tuples by
increasing chains. 
We conclude this section with an application of our results to the
decision problem whether a definable partial ordering is a 
quasi-well-ordering. 
The \emph{recurrent reachability problem} is defined as follows.

\begin{definition}
  Given a starting point $p$, a relation $R$ and 
  a subset $S$, decide whether there is an infinite $R$-path starting
  at $p$ and reaching $S$ infinitely often. 
\end{definition}

To and Libkin proved that the recurrent reachability problem is
decidable on automatic structures with a regular set $S$ and a
transitive, regular relation $R$. They solve the problem
globally, i.e., they construct an automaton that accepts those starting
points for which the set $S$ is recurrently reachable. 

\begin{theorem}[\cite{ToL08}]
  Let $\mathfrak{A}$ be an automatic structure with an automatic
  transitive relation $R$ and let $S$ be a regular subset of its
  domain $A$. Then
  the recurrent reachability problem  for $R$ and $S$ is
  decidable. Moreover, one can effectively construct an automaton
  $\mathcal{R}(R,S)$
  that
  accepts those nodes $p$ of $\mathfrak{A}$ such that there starts an
  infinite $R$ path at $p$ that passes $S$ infinitely often. 
  The size of $\mathcal{R}(R,S)$ is polynomially bounded in the size of the
  automata for $R$ and $S$. 
\end{theorem}

The proof of To and Libkin gives an explicit construction of
$\mathcal{R}(R,s)$. Roughly speaking, this construction yields an
automaton  corresponding to  an existential quantification over a
tree-comb whose elements form an $R$-chain in $S$.

Our proof can be adapted to reprove the decidability of the recurrent
reachability problem in To's and Libkin's setting. 
If we consider transitive relations, the recurrent reachability problem has
solutions of two different types. Either there is an element $p R q$
such that $q R q$ and $q\in S$ or there is an infinite chain $p R q_1
R q_2 R q_3 \dots$ of pairwise distinct elements 
$q_1, q_2, q_3, \dots \in S$. The first case is 
first-order definable whence it is decidable on automatic
structures. Hence, we only have to provide a decidability 
result for the other case. 
In order to obtain this result,  we modify our reduction of the Ramsey
quantifier to a 
reduction of a kind of chain quantifier.
Recall that the reduction of a Ramsey quantifier is 
of the form
\begin{align*}
  \reduction(\RamQ{n} \bar x(\varphi)) :=&
  \exists T_1,T_2\in \fullTrees{\Gamma_\Box},G\in\fullTrees{\{0,1\}}
  \psi_{\mathrm{CoSm}}(T_1,T_2,G) \land
  \psi_{\mathrm{Ram}}(T_1,T_2,G),
\end{align*}
where
\begin{align*}
    \psi_{\mathrm{CoSm}}(T_1, T_2, G) :=&
  \mathrm{Coherent}(T_1,T_2,G)  \land  \SMALLRel(T_1, G)
  \land \forall x (\In(x,T_1,T_2,G)\to x\in A) 
\end{align*}
and 
\begin{align*}
    \psi_{\mathrm{Ram}} &:=
    \forall x_1,\dots,x_n\in \Trees{\Sigma}\\
    &\left(\left(
    \bigwedge_{1\leq i \leq n}\In(x_i,T_1,T_2,G)
    \land  \bigwedge_{1 \leq i< j \leq
      n}\mathrm{Comp}(x_i,x_j,T_1,T_2,G) \right)
    \to \reduction(\varphi)\right).\\
\end{align*}
In order to solve the recurrent reachability problem, we propose to 
replace $\psi_{\mathrm{Ram}}$ by the formula 
\begin{align*}
  \forall x_1,x_2 \in \Trees{\Sigma}\left(\left(
    \bigwedge_{i\in\{1,2\}} \In(x_i,T_1,T_2,G) 
    \land  x_1 <_G x_2 \right) \to  (pRx_1 \land x_1Rx_2 \land Sx_2)\right)
\end{align*}
Note that $pRx_1 \land x_1 R x_2 \land Sx_2$ is
represented by some automaton due to the regularity of $S$ and $R$. 
Analogously to the 
automaton recognising $\mathrm{Comp}$, there is an $\omega$-automaton
for ``$x_1 <_G x_2$'' on input $(x_1, x_2, T_1, T_2, G)$ for all
\SMALL and coherent triples $(T_1, T_2, G)$. 

The formula asserts that there is a
closure of some coherent and \SMALL tree-comb such that each element
$a$ of this closure 
satisfies the following conditions:
\begin{enumerate}
\item $a$ is an $R$ successor of $p$,
\item $a$ is in $S$ and
\item if there is some $a'$ with $a <_G a'$ then $a R a'$ holds. 
\end{enumerate}
A witness $(T_1, T_2, G)$ for this assertion
induces an infinite
increasing $<_G$ chain in $S$. Hence, the soundness of this reduction is
obvious. For the completeness, we use  a
tree-comb that is homogeneous with respect to the automaton
corresponding to $pRx_1 \land x_1 R x_2 \land Sx_2$. Recall 
that any witness of the recurrent reachability of $S$ is
an infinite set that is linearly ordered by $R$. Analogously to the
fact that any infinite set 
contains a tree-comb, one proves that any ascending infinite chain
contains a tree-comb whose induced order coincides with the order of
the chain. Once we have obtained this result, the decidability proof
for the recurrent reachability problem is analogous to the
decidability proof of the Ramsey quantifier.

Let us conclude this section with an application of our result to
partial orderings.

\begin{example}
  Consider a formula $\varphi(x,y)$ that defines a partial order $\leq$
  on some automatic structure $\mathfrak{A}$. Assume that
  $\varphi$ is represented 
  by some automaton, e.g., assume that $\varphi\in\FO{}$. 
  Now, the Ramsey quantifier can
  be 
  used to formalise the existence of an infinite antichain. Let
  $\psi(x,y):=\RamQ{2} x,y (\neg \varphi(x,y) \land
  \neg\varphi(y,x))$. $\psi$ asserts that there is an infinite set of
  pairwise $\leq$-incomparable sets, i.e., an infinite antichain. 
  Thus, there is an automaton corresponding to the assertion that
  $\leq$ does not contain an infinite antichain.  
  
  We can also construct effectively an automaton that decides whether 
  $<$ contains an infinite descending chain. This is the same as
  deciding whether there is some point $p$ for which $>$ satisfies the
  recurrent reachability problem with respect to the full domain of
  the structure. 
  
  If there is neither an infinite antichain nor an infinite descending
  chain,  $\leq$ is a well-quasi-ordering. Thus, 
  if $\varphi$ defines a partial order, 
  the statement $\mathrm{WQO}(\varphi):=$``$\varphi$ induces a
  well-quasi-ordering'' is 
  decidable on automatic structures. Furthermore, one can
  effectively construct an
  automaton that corresponds to $\mathrm{WQO}(\varphi)$. 
\end{example}

\label{EndOfRamsey}

\chapter{Conclusions}
\label{Chapter_Conclusions}
In the following we summarise the main results of this thesis and
relate these results to open problems. 

We have shown that \FO{}(Reg) model checking on level $2$ collapsible
pushdown graphs is
decidable and Broadbent showed that first-order model
checking on level $3$ is undecidable (even with fixed formula or fixed
graph). The positive result on level $2$ is in fact even stronger: the
extensions by regular reachability, Ramsey quantifiers and
$L\mu$-definable predicates is still decidable. Hence, the structures
in level $3$ of the hierarchy  are much  more complicated than
those structures in level $2$. But it is still an open question what
the reason for this difference is. Broadbent's results point out that
even a very weak 
use of collapse operations already turns the first-order model
checking undecidable on level $3$. 
It would be nice to clarify which structural difference between the
graphs of level $2$ and those of level $3$ provokes the 
rather big difference in the model checking results. 
Another direction of further research is the
question for extensions of our results. What is the largest fragment
of \MSO that is decidable on collapsible pushdown graphs of level $2$?

We introduced the new hierarchy of higher-order nested pushdown
trees and provided first-order model checking algorithms for the
first two levels of this hierarchy. Due to its similarity to a
subclass of collapsible pushdown graphs, we also conjecture that the
$L\mu$ model checking is decidable on this hierarchy. 
But the proof of this conjecture is still open. Another open question
concerns 
first-order model checking on levels $3, 4, 5, \dots$ in this new
hierarchy. Our approach on level $2$, i.e., the use of the analysis of
strategies in the Ehrenfeucht-\Fraisse game via the notion of relevant
ancestors is extendable to higher levels. But on higher levels, we miss
an analysis of ``higher-order loops'' in analogy
to our results for loops of  collapsible pushdown systems of level $2$. 
Further research is necessary in order to clarify whether \FO{}
model checking on higher-order nested pushdown trees is decidable. 

Focusing on the second level of the nested pushdown tree hierarchy,
we are still 
lacking a characterisation of the complexity of first-order model
checking on nested pushdown trees of level $2$. Is there another
approach that yields an 
elementary complexity? Can we derive any reasonable lower bound for
the first-order model checking on nested pushdown trees?
We already know that first-order with reachability model checking
has nonelementary complexity on nested pushdown trees. 

Another more general question concerning first-order model checking
and collapsible pushdown graphs is the classification of those
graphs in the hierarchy that have decidable first-order theories. 
What kind of restrictions can one impose on the transition relation of 
a collapsible pushdown graph in order to obtain decidability of its
first-order theory?

Another open question concerns the  characterisation of the
differences between collapsible pushdown graphs and higher-order
pushdown graphs. We propose the further study of higher-order nested
pushdown trees in order to approach this question. 
Higher-order nested pushdown trees 
can be seen as collapsible pushdown graphs with a rather tame
application of collapse. 

We now turn to a more general direction of research. 
In this thesis, we have focused on what is called the local model
checking. Global model checking on the other hand asks for identifying
all elements in a given structure that satisfy some formula. Recently,
Broadbent et al.\ \cite{BroadbentCOS10} showed the global
$L\mu$ model checking on collapsible pushdown graphs to be decidable. 
Furthermore, they showed that collapsible pushdown graphs themselves
are sufficient to describe the result of the global $L\mu$ modal
checking in the following sense: for each formula and each pushdown
system 
there is another one that generates the same graph but marks each
element that satisfies the given formula.
The analogous questions for first-order model checking on
(higher-order) nested pushdown trees or collapsible pushdown graphs
have not been investigated yet and their investigation may reveal
interesting insights into these classes.

We have also shown that Ramsey quantifiers on tree-automatic
structures are decidable. This 
extends the corresponding result and proof techniques for the
string-automatic case. But 
in fact, on string automatic structure a far stronger logic is
decidable\cite{Kuske2009}. Kuske called this logic FSO. 
It is the extension of $\FO{}$ by existential quantification over infinite
relations that only occur negatively, i.e., under the scope of an odd
number of negations. Ramsey quantifiers can be rewritten in terms of
FSO. Thus, our result shows the decidability of a fragment of FSO on
tree-automatic structures. 
It is an open problem whether FSO is decidable on
all tree-automatic structures or whether some undecidable problem
may be encoded into FSO on some tree-automatic structure. 
For most results on string-automatic structures there has been found
an analogous one for tree-automatic structures. If this were not the
case for FSO model checking this may point to new insights into the
difference between tree-automata and string-automata.


\appendix
\chapter{Undecidability of ${\mathbf{L\mu}}$ on the Bidirectional Half-Grid}
\label{Appendix_Lmuundecidabilty}
In this appendix, we  show the undecidability of 
$L\mu$ on the bidirectional $\N\times\N$-grid, i.e., the grid
$\N\times\N$ with modalities
``left'', ``right'', ``up'', and ``down'' (denoted by $\leftarrow$,
$\rightarrow$, $\uparrow$, $\downarrow$). Note that we do not allow
atomic propositions apart from $\True$ and $\False$. 
The proof is by reduction to the halting problem of Turing machines. 
At the end, we will see how this proof generalises to
the case of the bidirectional half-grid $\BiHalfgrid$ (cf. Figure
\ref{HALFGRID}).
Before we start the proof, we will shortly recall our notation
concerning Turing machines and recall the definition of the halting
problem. 

\section{Turing Machines}

In order to fix notation, we briefly
recall the notion of a Turing machine.

\begin{definition}
  The tuple $\mathcal{M}=(Q, \Sigma, q_I, q_F, \Delta)$ with $Q$ the
  finite sets of states, $\Sigma$ the finite tape alphabet, 
  $q_I,q_F\in Q$ the initial, respectively, final state, and
  \begin{align*}
    \Delta: Q \times \Sigma \rightarrow
    \Sigma \times \{l,r\} \times Q    
  \end{align*}
  a transition function is called a
  \emph{Turing machine}.  
  The \emph{set of configurations of $\mathcal{M}$} is
  \begin{align*}
    \mathrm{CONF}:=\Sigma^\omega\times\N\times Q.      
  \end{align*}
\end{definition}

\begin{remark}
  The elements in $\{l,r\}$ are called 
  head instructions where $l$ denotes ``move to the left'' and $r$
  denotes ``move to the right''.

  In the literature this definition is normally called a
  deterministic Turing machine, while in the general nondeterministic
  case $\Delta$ is assumed to be a relation instead of a
  function. We restrict ourselves to
  deterministic Turing machines because they have the same computational
  power as nondeterministic ones (cf. \cite{HopcroftU79}).
\end{remark}

The notion of the computation of a Turing machine is captured by the
notion of a run of the machine. A run is a list of configurations
where the $(n+1)$-st configuration evolves from the $n$-th by applying
$\Delta$. 

\begin{definition}
  Let $M=(Q, \Sigma, q_I, q_F, \Delta)$ be a Turing machine. 
  For $w\in\sigma^\omega$ we write $w(i)$ for the $i$-th letter in
  $w$. 
  $M$ induces a function $\trans{}: \mathrm{CONF}\to \mathrm{CONF}$ as
  follows. Let $(w,i,q)$ and $(w',i',q')$ in $\mathrm{CONF}$. 
  Assume that $\Delta(q,w(i))=(\sigma',o,q')$. 
  It holds that 
  $(w,i,q) \trans{} (w',i',q')$ if 
  \begin{enumerate}
  \item $w'(j) = w(j)$ for all $i\neq j \in\N$, 
  \item $w'(i) = \sigma'$, and
  \item $i'=i-1$ if $o=l$ and $i'=i+1$ if $o=r$. 
  \end{enumerate}
  A \emph{run of $M$ on input
  $w\in\Sigma^*\{0\}^\omega$} is a function $\rho: \N \to \mathrm{CONF}$ such
  that $\rho(0)=(w,1,q_i)$ and $\rho(i)\trans{}\rho(i+1)$ for all
  $i\in\N$. 

  We say that the computation of $M$ on $w$
  \emph{terminates} if there is some $i\in\N$ such that 
  \mbox{$\rho(i) = (w',i',q_F)$}  for $\rho$ the run of $M$ on input $w$. 
\end{definition}

For a detailed introduction into the theory of Turing machines
we recommend \cite{HopcroftU79}. For the purpose of this proof, we
only need one of the cornerstones of computability theory: the
undecidability of the halting problem. 
The halting problem is  the problem whether the computation of a given
Turing machine terminates on input $0^\omega$. 
This is one of the classical examples of undecidable problems. 

\begin{theorem}[\cite{turing1936a}]
  The halting problem is undecidable. 
\end{theorem}

\section{Reduction to the Halting Problem}

For simplicity, we only consider Turing machines with tape
alphabet $\Sigma=\{0,1\}$. We assume that the state set is
$Q=\{q_1, q_2, \ldots, q_{\lvert Q \rvert}\}$. 
We will represent the run of an arbitrary fixed
Turing machine on input $0^\omega$ 
as an $L\mu$-definable colouring of the bidirectional grid. Then using
an $L\mu$ 
definable reachability query for a final state of the Turing machine,
the halting problem is reduced to $L\mu$ model-checking on the
bidirectional grid. 
The idea is as follows. 

Each configuration $c=(w,p,q)$ of a Turing machine $M$ can be encoded as
an infinite bitstring 
$v\in\{0,1\}^\omega$. We use the letters $v((\lvert Q \rvert +3)\cdot
i), \ldots, 
v((\lvert Q \rvert +3)\cdot (i+1)-1)$ to encode the information
concerning the $i$-th cell of $c$. We use the last two bits of the
representation of a cell to indicate whether the corresponding cell
contains the letter $0$ or $1$. Thus, we define $v((\lvert Q \rvert
+3)\cdot (i+1)-1) :=1$ if and only if $w(i)=1$ and
$v((\lvert Q \rvert+3)\cdot (i+1)-2) :=1$ if and only if $w(i)=0$. 
The other $\lvert Q \rvert +1$ bits are used to encode the information
about the position of the head $p$ and the state of the machine
$q$. We set the first bit of a cell to $1$ if the head 
is not above this cell and we set the $j$-th bit of a cell, if the
head is above this cell and $q=q_{j-1}$. Formally, we set 
$v(k):=1$ if $k= (\lvert Q \rvert+3)\cdot i$ for some $i\neq p$ or
$k = (\lvert Q \rvert+3)\cdot p +j$ for $q=q_j$. 
All other positions in $v$ are set to $0$. We denote as
$W:\mathrm{CONF}\rightarrow \{0,1\}^\omega$ the function that
translates each configuration into the corresponding encoding.

Now, it is easy to encode the run $\rho$ of $M$ on $0^\omega$ in the infinite
$\N\times\N$-grid using a set $X_M\subseteq \N\times\N$. 
We write $c_n:=\rho(n)$ and $v_n:=W(c_n)$. 
We set $(i,j)\in X_M$ if and only if $v_i(j) = 1$. 

As a next step, we show that $X_M$ is definable in $L\mu$. Using this
fact, we later define an $L\mu$ formula that is true on the infinite
grid if and only if $\rho$ terminates.  

\begin{lemma}
  There is an effective translation from a given Turing machine  $M$
  into an $L\mu$ formula 
  $\varphi_M$ such that $\varphi_M$ defines $X_M$ on the infinite grid. 
\end{lemma}
\begin{proof}
  As a preliminary step, we want to define those $(i,j)\in\N\times\N$
  where the encoding of one of the cells start, i.e., the set
  $\{i,j:\exists k\in\N\ j = k \cdot (\lvert Q \rvert +3)\}$. 

  This is done by the formula 
  \begin{align*}
    \TMCells:= \mu X. \left(([\uparrow] \False \land [l]
    \False) \lor \langle \uparrow \rangle^{\lvert Q \rvert +3} X \lor
    \langle \leftarrow \rangle X \right).     
  \end{align*}
  The first part of this formula defines the position $(0,0)$, the
  second part adds the position $(i,j +\lvert q \rvert +3)$ to
  $\TMCells$ for each $(i,j)\in\TMCells$ and the last part adds $(i+1,j)$
  to $\TMCells$ for $(i,j)\in\TMCells$. 
  We call a node in $\TMCells$ initial position of an encoding of
  a cell, or simply initial position of a cell. In the following we
  identify the initial position of a cell with the cell
  itself. 

  In the following, we use some auxiliary formulas:
  \begin{itemize}
  \item We set $\TMPos_i:= \langle \uparrow \rangle^i \TMCells$ for $0\leq
    \lvert Q \rvert +2$. $\TMPos_i$ holds on $(j,k)$ if $(j,k)$ is the
    $i+1$-st bit of the encoding of one of the cells. 
  \item For $\varphi\in L\mu$, we write $\TMFrom_0(\varphi)$ for the
    formula
    $\bigvee_{i=0}^{\lvert Q \rvert +2} \left(\TMPos_i \wedge \langle \uparrow
    \rangle^i \varphi\right)$. Some node $(i,j)$ satisfies $\TMFrom_0(\varphi)$
    if and only if the initial position corresponding to the same cell
     as $(i,j)$ satisfies $\varphi$. 
  \item We need formulas for navigation on the cells in encoded
    form. We set 
    \begin{align*}
      &\TMLeft(\varphi):= \TMFrom_0( \langle \uparrow^{\lvert Q\rvert +3}\rangle
      \varphi)\text{ and}\\
      &\TMRight(\varphi):= \TMFrom_0( \langle \downarrow \rangle^{\lvert Q\rvert +3}
      \varphi).      
    \end{align*}
    These formulas are satisfied at $(i,j)$ if
    ``the cell to the left of the one corresponding to $(i,j)$
    satisfies $\varphi$'', respectively ``the cell to the right
    \ldots''. 
    Similarly we use $\TMBefore(\varphi):= \TMFrom_0 ( \langle \leftarrow
    \rangle \varphi)$. This formula is satisfied by node whose cell
    satisfied $\varphi$ in the preceding configuration. 
  \end{itemize}

  As a next step we introduce some formulas that recover information
  about the configuration from its encoding in the grid. For this we
  assume that $X$ is a colouring that colours each column of the grid
  with the encoding of some configuration. 
  \begin{itemize}
  \item For $q_i\in Q$ set 
    $\TMState_{q_i}:= \TMFrom_0( \langle \downarrow \rangle^i X)$. 
    This formula is satisfied at $(i,j)$ if the head of $M$ is above
    the corresponding cell and the state of $M$ is $q_i$. 
    Analogously, we write 
    $\TMState_0:= \TMFrom_0(X)$ for the formula specifying the head
    is not at the cell corresponding to this position. 
  \item For $\sigma\in\{0,1\}$ we set
    $\TMTape_\sigma:= \TMFrom_0( \langle \downarrow \rangle^{\lvert Q\rvert
      +1+\sigma} X)$. $\TMTape_\sigma$ is satisfied if the corresponding
    cell contains the letter $\sigma$. 
  \end{itemize}
  Since we want to generate the encoding of the run of $M$ on input
  $0^\omega$ by a fixpoint formula in $L\mu$, we have to define
  ``update-formulas'' which, given the encoding of a valid
  configuration in the $i$-th column, return the encoding of the
  following configuration in the $i+1$-st column. 
  
  Let us first consider the information concerning the update of the
  tape. There are two possibilities for each cell: either the head is
  not at this cell, then its value is preserved by $\trans{}$, or the
  head is at this cell, then its value depends on $\Delta(q,\sigma)$
  where $q$ is the state of the machine and $\sigma$ is the symbol of
  this cell. 
  For $\sigma\in\Sigma$, set
  \begin{align*}
    Set_\sigma:=\{(\tau,q_k)\in\Sigma\times Q:
    \Delta(q_k,\tau)=(\sigma, d, q), d\in\{l,r\}, q\in Q\}.
  \end{align*}
  Note that $Set_1$ contains those combinations of a letter $\sigma$ and
  a state $q$  such that the head of $M$ will write $1$ onto the tape
  if it is in state $q$ and reads $\sigma$. The analogous claim holds
  for $Set_0$. 
  Thus, 
  \begin{align*}
    \TMUpdate_\sigma:= \TMPos_{\lvert Q \rvert +1+\sigma} \land
    \big( \TMBefore(\TMState_0\land\TMTape_\sigma) \lor
    \bigvee_{(i,q)\in Set_\sigma} \TMBefore(\TMState_q \wedge
    \TMTape_i)\big)
  \end{align*}
  are the correct update formulas for the information concerning the
  tape. 
  
  The update for the information concerning the head of the Turing
  machine are slightly more complicated because the head can reach a
  certain cell either from the cell to the left or from the cell to
  the right. Furthermore, we also have to update the information on
  those cells where the head is not positioned. There, we have to set
  the first bit encoding the cell which represents the absence of the
  head from this cell. 
  We start by collecting those combinations of states and letters that
  induce the head to move to the left or to the right,  respectively.
  For $q_j\in Q$, set 
  \begin{align*}
    &\mathrm{Left}_{q_j}:=\{(\sigma,q)\in \Sigma\times Q: \exists \tau\in\Sigma
   \ \Delta(q,\sigma) = (\tau, l, q_j)\}, \text{ and} \\
    &\mathrm{Right}_{q_j}:=\{(\sigma,q)\in \Sigma\times Q: \exists \tau\in\Sigma
   \  \Delta(q,\sigma) = (\tau, r, q_j)\}. \\
  \end{align*}
  Now, we can use these sets to define the position of the head in the
  next configuration. We set
  \begin{align*}
    &\TMUpdate_{q_j}:=\TMPos_j \land \\
    &(\bigvee_{(\sigma,q)\in \mathrm{Left}_{q_j}}
    \TMBefore(\TMRight( \TMState_q \land \TMTape_\sigma)) \lor
    \bigvee_{(\sigma,q)\in \mathrm{Right}_{q_j}}
    \TMBefore(\TMLeft( \TMState_q \land \TMTape_\sigma))).
  \end{align*}
  These formulas update correctly the information on the head of the
  tape, i.e., if $X$ encodes some configuration in the $i$-th column
  of the grid, then $\TMUpdate_{q_k}$ will hold at some $(i+1,k)$ if
  and only if the next configuration has state $q_k$, $(i+1,k)$
  corresponds to the encoding of state $q_k$ in some cell, and the
  head is at this cell in the next configuration. 
  Of course, we also have to propagate the ``no-head'' information
  along all other cells. This is the case if either the head is
  neither to the left nor to the right in the previous configuration
  or it is positioned one step to the left, respectively right,  but
  will move further left, respectively right. Note that by the
  definition of a Turing machine, the sets $(\mathrm{Left}_{q})_{q\in Q}$ and
  $(\mathrm{Right}_{q})_{q\in Q}$ form a partition of $\Sigma\times Q$. Hence,
  the following formula updates the ``no-head'' information:
  \begin{align*}
    \TMUpdate_0:=\big( \TMPos_0 \land 
    &\TMBefore(\TMLeft(\TMState_0))\land
    \TMBefore(\TMRight(\TMState_0))\big) \\
    &\lor
    \bigvee_{(\sigma,q)\in \mathrm{Left}_{q_j}}
    \TMBefore(\TMLeft( \TMState_q \land \TMTape_\sigma)) \\
    &\lor
    \bigvee_{(\sigma,q)\in \mathrm{Right}_{q_j}}
    \TMBefore(\TMRight( \TMState_q \land \TMTape_\sigma)). 
  \end{align*}
  The first part of this formula deals with positions where the head
  in the previous configuration is not one step to the left or to the
  right and the second and third part update the ``no-head'' information
  if the head is close but is going to move further away. 

  Having defined the necessary formulas for the update from one
  configuration to the next, we have to define the colouring of
  the initial configuration of the run of $M$ on input $0^\omega$ in
  order to obtain the encoding of the full run by a least fixpoint
  induction. For this purpose, we assume that $q_I=q_1$ and we set
  \begin{align*}
    \TMInit:= \langle \leftarrow \rangle \False \land 
    (\TMPos_{\lvert Q \rvert +1} \lor (\TMLeft(\True) \land \TMPos_0)
    \lor (\TMLeft(\False) \land \TMPos_1)).
  \end{align*}
  $\TMInit$ holds only at positions on the leftmost column of the
  grid which means that it only initialises the encoding of the first
  configuration. The first part sets in each cell the position $\lvert Q
  \rvert +1$ which corresponds to setting all cells to $0$. Secondly,
  we set in all but the first cell the ``no-head''
  information. Finally, we set the state $q_I=q_1$ in the first
  cell. Hence, $\TMInit$ is the definition of the encoding of the first
  configuration of the run of $M$ on input $0^\omega$.

  We claim that the formula
  \begin{align*}
    \varphi_M:=\mu X. (\TMInit \lor \bigvee_{q\in Q} \TMUpdate_q \lor
    \TMUpdate_0 \lor \bigvee_{\sigma\in\Sigma} \TMUpdate_\sigma)
  \end{align*}
  defines the encoding of the run of $M$ on input $0^\omega$ on the
  bidirectional grid. In fact, an easy but technical induction shows
  that the $n$-th stage of the fixpoint of $\varphi_M$ defines exactly
  the encoding of the first $n$ configurations of the run of $M$ on
  input $0^\omega$. 
\end{proof}

From the translation of runs of Turing machines into $L\mu$ formulas
on the grid, the undecidability of the $L\mu$ model checking on the
grid follows immediately.

\begin{lemma}
  $L\mu$ model checking is undecidable on the bidirectional grid.  
\end{lemma}
\begin{proof}
  By reduction to the halting problem:
  Deciding the halting problem for $M$ is the same as deciding
  whether $\varphi_M$ defines some cell where
  $\TMState_{q_F}$ holds. But this is the same as deciding whether
  \begin{align*}
    \TMHalting:=\mu Y. \TMState_{q_F}(\varphi_M) \lor \langle \downarrow \rangle
    Y \lor \langle \rightarrow \rangle Y
  \end{align*}
  is satisfied in the position $(0,0)$ of the bidirectional
  grid. Thus, a model checking algorithm of  $L\mu$ on the
  bidirectional grid would lead to a decision procedure for the
  halting problem. This proves the undecidability of
  $L\mu$ on the grid. 
\end{proof}

\begin{remark}
  In the presence of a universal modality, the proof can be adapted to
  show the undecidability of the $\N\times\N$ grid only with the
  modalities left and up. The search for a cell with state $q_F$ can
  be done by using the universal modality. Hence, we only have to
  remove the down modalities from the update formulas. This can be
  achieved by shifting the beginning of the encoding of the $i$-th
  column by $2i (\lvert Q \rvert +3)$. 
\end{remark}

Having obtained the undecidability of $L\mu$ on the grid, we want to
refine the result such that it  applies to the half-grid. But this is
easy by 
noting that the head of the Turing machine in the $i$-th configuration
of a run can only have visited the first $i$-cells of the tape. 

\begin{corollary}
  $L\mu$ on the bidirectional half-grid is undecidable.
\end{corollary}
\begin{proof}
  Instead of encoding the $i$-th configuration of the run of $M$ in
  the $i$-th row, we can use the $i (\lvert Q \rvert +\lvert \Sigma
  \rvert +1)$-st row instead. Doing this, the head of the Turing
  machine is in all configurations at some cell which is encoded by
  elements in the grid of the form $(i,j)$ where $i>j$. Thus, we can
  treat the missing nodes in the half-grid $\{(i,j)\in \N\times\N:
  i>j\}$ as the encodings of cells which contain the symbol $0$ and
  which contain the ``no-head'' marker. 
\end{proof}


\bibliography{Standard}
\bibliographystyle{plain}

\chapter*{Wissenschaftlicher Werdegang von Alexander Kartzow,\\
  geboren  am 12. Januar 1983 in Gie\ss en}

\begin{tabular}{lc}
{\bf Schulabschluss}\\
 Juni 2002  & Abitur\\
 \\
 {\bf Studium}\\
 2002-2007 &Studium der Mathematik mit Nebenfach Informatik\\
 & an der TU
 Darmstadt \\
 2005-2006 & Studium der Mathematik und der Informatik\\
 & an der Universidad de Salamanca,
 Spanien\\
 November 2007 & Diplom in Mathematik\\
\\
{\bf Promotion}\\
2007 - 2011 & Promotionsstudium \\
&an der TU Darmstadt
\end{tabular}

\end{document}